\documentclass[11pt]{article}
\usepackage[font=small,labelfont=bf]{caption}
\usepackage{latexsym,amsfonts,amsthm,amsmath,amssymb}
\usepackage{hyperref}
\usepackage{xcolor}
\usepackage{verbatim}
\usepackage{enumerate}
\usepackage{tabularx,graphicx}
\usepackage[ngerman]{babel}
\usepackage{epigraph}

\newtheorem{Thm}{Theorem}

\newtheorem{Def}{Definition}

\newcolumntype{L}{>{\displaystyle}l}
\newcolumntype{R}{>{\displaystyle}r}

\begin{document}
\title{Eulers Horizonte -- \\
Möglichkeiten und Grenzen seiner Arbeitsweise in der Mathematik}
\author{Alexander Aycock}  \vspace{15mm}
\date{\textit{``Repetitio cognitorum est fortuna ignorantium."{}}}
\maketitle

\newpage 
\tableofcontents

\newpage

\listoffigures

\newpage

\section{Geleitwort}

 \epigraph{Lest EULER - er ist unser aller Meister!}
 {Pierre--Simon Laplace}%

Leonhard Euler (1707--1783) war unumstritten einer der produktivsten und vielseitigsten Mathematiker aller Zeiten. Seine gesammelten Werke umfassen mittlerweile 4 Reihen bestehend aus insgesamt 84 Büchern von je 300 bis 700 Seiten, sodass sich sein  Opus in Gänze auf ungefähr 30000 beläuft. Reihe 1 ist dabei seinen mathematischen Untersuchungen gewidmet und zählt 30 Bände, Reihe 2 mit 31 Bänden enthält Eulers Untersuchungen zur Physik, die 12 Bände der 3. Reihe haben allgemeine Untersuchungen zum Inhalt. Schließlich ist in der 4. Reihe Eulers wissenschaftlicher Briefwechsel zu finden. Eulers Opus ist demnach nicht nur am reinem Umfang gemessen,  sondern auch der Qualität nach beurteilt immens. Die Leistung der Herausgeber seiner gesammelten Werke als \textit{Opera Omnia Leonhardi Euleri} kann demnach kaum überbewertet werden.\\

Angesichts des gewaltigen Umfangs Euler'schen Werkes mag es jedoch verwundern, dass Klein (1849--1925) auf den Seiten 4--5 über die Darstellung mathematischer Werke im 19. Jahrhundert in seinem Buch \textit{``Vorlesungen über die Entwicklung der Mathematik"} (\cite{Kl56}, 1956) schreibt :\\

\textit{``Freilich  werden die großen Schönheiten dieser kristallklaren, abgeschlossenen, klassischen Darstellungsweise nicht ohne Einbuße erkauft. Es ist\\ nämlich diesen Meisterwerken kaum mehr ihre Werdegeschichte zu entnehmen. Dadurch ist dem Leser die eigentümliche und für einen selbständigen Geist größte Freude versagt, unter der Führung des Meisters die gefundenen Resultate selbstständig gleichsam noch einmal zu entdecken. In diesem Sinne mangelt es den Werken der klassischen Zeit das eigentliche erzieherische Moment. Der Gedanke, den Leser nicht nur zu erfreuen und zu belehren, sondern in ihm über das Werk hinausgehende Kräfte zu wecken, zur eigenen Tätigkeit anzuregen -- eine Wirkung, wie sie etwa von Monges, von Jacobis oder auch von Faradays Schriften ausgeht -- gehört durchaus dem 19. Jahrhundert an."}\\

Das Wundersame ist, dass Leonhard Euler als einer der, wenn nicht gar der führende Mathematiker, was die Erklärungen seiner Entdeckungen betrifft, in obigem Zitat keine Erwähnung findet\footnote{Eine Aufklärung dieses Umstandes findet man in einer Fußnote der Herausgeber von Kleins Buch: Nämlich,  dass sich die mangelnde Würdigung Eulers von Kleins Seite in der Nichtbeschäftigung des letzteren mit den Werken von erstem begründet liegt.}. Die nicht oder kaum stattfindende  Betrachtung des Euler'schen Opus in Büchern zur Mathematikgeschichte ist jedoch kein Einzelfall.  Entweder sind  die entsprechenden Werke sämtlich Euler gewidmet oder der Euler'sche Beitrag wird in gedrungener Form wiedergeben, um schnellstmöglich zum eigentlichen Gegenstand überzugehen. Nicht zuletzt daran macht man den Paradigmenwechsel in der Mathematik von Euler, der die Mathematik im 18. Jahrhundert entscheidend geprägt hat, zur Moderne aus.\\

 Eulers \textit{Opera Omnia} war die Grundlage der Arbeit des \textit{Euler--Kreis Mainz}$^{\copyright}$,  einer Gruppe historisch interessierter Mathematiker und Mathematikerinnen, die sich mit den Schriften und dem Werk Leonhard Eulers beschäftigen. Die Vorträge, Abhandlungen sowie Übersetzungen ausgewählter Euler'scher Arbeiten aus dem Lateinischen und Französischen in die deutsche und englische Sprache sind online\footnote{Der Link zum Euler--Kreis Mainz lautet \url{https://www.agtz.mathematik.uni-mainz.de/algebraische-geometrie/van-straten/euler-kreis-mainz/}.} zugänglich\footnote{ Basierend auf den Euler'schen Arbeiten fanden ebenfalls studentische Seminare statt.}. Ein Teil derselben Translationen sind auch beim \textit{Euler Archive} unter \url{http://eulerarchive.maa.org/} von dessen Betreibern veröffentlicht worden. Letzteres hat ebenfalls gescannte Versionen von allen Originalarbeiten Eulers bereitgestellt\footnote{Alle Scans von Eulers Arbeiten, die in die vorliegende Abhandlung eingefügt wurden, stammen aus dieser Quelle.}. Bei der Übersetzungsarbeit im \textit{Euler--Kreis} sind überdies kleinere Entdeckungen gemacht worden, welche anschließend den Inhalt eigener Abhandlungen gebildet haben. Selbige haben auch Eingang in das \textit{``Euleriana Archive"}\footnote{Der Link lautet: \url{https://scholarlycommons.pacific.edu/euleriana/}.}  gefunden. Dabei handelt es sich um ein wissenschaftliches Journal, das Arbeiten von und über Euler gewidmet ist.\\

Warum jedoch gerade die Beschäftigung mit Euler? Zum einen, weil sich dem auf Seite 235 in der \textit{``Encyclopedia of Mathematics Education"} (\cite{Gr01}, 2001)  Gauß (1777--1855)  zugeschriebenen Ausspruch\\ 

\textit{``The study of Euler's works will remain the best school for the different fields of mathematics and nothing else can replace it."}\\

auch heute noch beipflichten lässt, was den reinen Inhalt betrifft. Zum anderen gestatten Eulers Ausarbeitungen, die Hürden zu erkennen, welche sich beim mathematischen Schaffensprozess auftun, welche Euler in ihn auszeichnender Art auch stets eingeräumt hat. Obschon sich die Mathematik durch ihre geforderte Beweisstrenge auszeichnet, lässt dies doch die Frage offen, wie diese Schlüsse, genauer die Prämissen und die Folgerungen, ursprünglich ausgemacht worden sind. Die Virtuosität Eulers in dieser Tätigkeit in Kombination mit seinen ausführlichen Erläuterungen sucht  wohl bis zum heutigen Tage ihresgleichen. Ein Zeugnis von Eulers meisterhaftem Umgang mit Formeln und  der Erstellung derselben ist, dass sich viele Ausdrücke, welche die exakte Auswertung von Reihen und Integrale betreffen, bereits bei Euler finden, obschon sie den Namen seiner Nachfolger tragen\footnote{Ein ganzer Abschnitt ist eigens diesem Sachverhalt gewidmet.}. \\

Angesichts dessen mag es schwer zu glauben sein, dass Euler überhaupt an irgendwelche Grenzen stoßen konnte, welche nicht durch seine endliche Lebenszeit begründet waren. Dennoch trägt sich das natürlich zu, dass selbst Euler manche Probleme nicht auflösen konnte. Die Ursachenforschung dessen führt, je nach behandeltem Themenbereich, zu unterschiedlichen Antworten. Allem geht aber die schlichte, wenngleich rückblickend offenkundig anmutende, Feststellung voraus, dass Euler  den entsprechenden Lehrsatz, die etwaige Formel oder dergleichen nicht zwingend in derselben Form zu erhärten versucht hat, wie es heutzutage geschähe. Stets ist nämlich zu beachten, dass etwaige heute wohl vertraute Konzepte  im 18. Jahrhundert noch gar nicht existierten oder gar eine andere Bedeutung hatten\footnote{Ein prägnantes Beispiel dessen ist das des Begriffs der Funktion. Von Euler in seinem berühmten Lehrbuch \textit{``Introductio in analysin infinitorum, volumen primum"} (\cite{E101}, 1748, ges. 1745) (E101: ``Einleitung in die Analysis des Unendlichen, erster Band") gleichsam in die Mathematik eingeführt, ist die Euler'sche Definition einer Funktion  nicht mit der verträglich, die man heute in einem Einführungsbuch zur Analysis vorfindet.}. \\

Ein Vergleich mit der heutigen Mathematik zeichnet hierneben auch von Euler weniger bis gar nicht berührte Gebiete aus\footnote{Bezüglich der Euler'schen Arbeiten ist hier insbesondere die komplexe Analysis zu nennen, zu welcher  Euler freilich wichtige Beiträge geleistet hat, aber nie eine umfassende Theorie vorgestellt hat, wie dies Cauchy (1789--1857) und Riemann (1826--1866) später getan haben.}. Was andererseits Abhandlungen \textit{über} Euler betrifft, findet man insbesondere das Thema, was er -- aus der modernen Sicht betrachtet -- bereits alles vorweg genommen hat und was er selbst aus seinen Entdeckungen heraus noch zu leisten vermocht hätte, kaum diskutiert. Die Aufhebung dieses Umstandes ist eines der Kernanliegen der gegenwärtigen Arbeit. \\

Unabhängig vom Inhalt,  tritt weiter die paradigmatische\footnote{In diesem Kontext referiert dies auf die Definition von Thomas Kuhn (1922--1996) in seinem Buch \textit{``The Structure of Scientific Revolutions"} (\cite{Ku12}, 2012).} Struktur der Mathematik hervor, welche sich in unterschiedlicher Gewichtung einiger mathematischer Grundbegriffe zwischen dem 18. Jahrhundert und der Moderne manifestiert.  So mutet Eulers zur Auffassung dessen, was einen Beweis darstellt, heute eigentümlich an und  seine Argumentationen sind  mancherorts zu ergänzen, um den Ansprüchen des heutigen Rigors zu entsprechen\footnote{Hiervon werden im Haupttext einige Beispiele auftreten, welche auch korrigiert werden, sofern dies möglich ist.}. Manche seiner Beweisführungen sind indes nicht zu reparieren, auch  wenn sie zu seiner Zeit Akzeptanz gefunden haben.\\

Die vorliegende Ausarbeitung intendiert, die zuvor erwähnten teils auch anderenorts publizierten  und an dieser Stelle erstmals mitgeteilten Funde des \textit{Euler-Kreis Mainz} bezüglich der Euler'schen Werke, gesammelt und möglichst kohärent darzustellen\footnote{Folgende bereits veröffentlichte Arbeiten des Verfassers haben Eingang in die vorliegende Arbeit gefunden und bilden wesentliche Bestandteile der entsprechenden Abschnitte: \textit{``On proving some of Ramanujan’s formulas for $\frac{1}{\pi}$ with an elementary method"} (\cite{Ay13}, 2013) sowie die mit weiteren Autoren veröffentlichte Arbeit \textit{``Proof of some conjectured formulas for $\frac{1}{\pi}$ by Z.W.Sun"} (\cite{Al11}, 2011) (Abschnitt \ref{subsubsec: Ramanujans Formeln zur Kreisquadratur}), \textit{``Euler and the Legendre Polynomials"} (\cite{Ay23}, 2023) (Abschnitt \ref{subsubsec: Den Kontext betreffend: Die Legendre Polynome} und \ref{subsubsec: Weitere Untersuchungen zu den Legendre-Polynomen}), \textit{``Euler and the Gammafunction"} (\cite{Ay21a}, 2021) (Abschnitt \ref{subsubsec: Ein anderes Vorhaben: Das Weierstraß-Produkt} und \ref{subsubsec: Die Darstellung betreffend -- Die hypergeometrische Reihe}), \textit{``Answer to a question concerning Euler’s paper Variae considerationes circa series hypergeometricas"} (\cite{Ay22}, 2022) und \textit{``Euler and the Duplication Formula for the Gamma-Function"} (\cite{Ay23c}, 2021) (Abschnitt \ref{subsubsec: Durch Anwendung einer Methode: Seine Konstante A}), \textit{``Euler and Homogeneous Difference Equations with Linear Coefficients"} (\cite{Ay24a}, 2024) (Abschnitt \ref{subsubsec: Die Mellin--Transformierte bei Euler}), \textit{``On Euler's Solution of the Simple Difference Equation"} (\cite{Ay23b}, 2023) (Abschnitt \ref{subsec: Ein Prototyp -- Differenzengleichungen} und \ref{subsubsec: Herleitung Formel für die Potenzsummen der Reziproken}),  \textit{``Euler and the Gaussian Summation Formula for the Hypergeometric Series"} (\cite{Ay24b}, 2024) (Abschnitt \ref{subsubsec: Die Darstellung betreffend -- Die hypergeometrische Reihe}), \textit{``Euler and a Proof of the Functional Equation for the Riemann Zeta-Function He Could Have Given"} (\cite{Ay24c}, 2024) (Abschnitt \ref{subsubsec: Durch Kombinieren von Ergebnissen: Die zeta-Funktion}).\\
Hier erstmalig mitgeteilte Entdeckungen sind unter anderem der Vergleich der Beiträge zur Theorie der hypergeometrischen Funktion von Euler und Gauß (Abschnitt \ref{subsubsec: Die Darstellung betreffend -- Die hypergeometrische Reihe}), die heuristische Herleitung des Primzahlsatzes aus den Euler'schen Formeln (Abschnitt \ref{subsubsec: Weg fehlender Formulierung: Der Primzahlsatz}) sowie die Diskussion der Gründe für die Nichtentdeckung der komplexen Analysis durch Euler (Abschnitt \ref{subsubsec: Komplexe Analysis}). Weiterhin scheint der Zugang zu Formeln für gewisse Orthogonalpolynome mit Euler'schen Methoden (Abschnitt \ref{subsubsec: Explizite Formeln aus Eulers Theorie zu Differenzengleichungen}) bisher in der Literatur nicht behandelt worden zu sein.}. Dies ließe sich freilich auf vielerlei Art leisten. Jedoch scheint die Orientierung am Euler'schen Gedankengang selbst dem gesteckten Ziel am ehesten zuzukommen, zumal sie unter anderem  zum  Verständnis von Eulers Kreativität und seiner Schaffenskraft beizutragen vermag. Darum  nimmt die Präsentation ihren Ausgangspunkt bei Eulers Ideen zur Mathematik selbst, geht dann zu exemplarischen Auszügen seiner Arbeitsweise über und beleuchtet die Schranken, an die er gestoßen ist (und aus moderner Sicht gesprochen, teilweise stoßen musste), um abschließend Anwendungen seiner Ergebnisse vorzustellen. Eine gewisse Redundanz der behandelten Themen ließ sich dabei nicht gänzlich vermeiden, da -- wie auch aufgezeigt werden wird -- Euler des Öfteren zu einem behandelten Gegenstand zurückkehrte, um ihn dann mit mittlerweile (meist von ihm selbst bei Behandlung scheinbar unverwandter Fragestellungen) geschaffenen Techniken neu anzugehen. Diese Herangehensweise gestattet des Weiteren ein Urteil darüber, wie weit Eulers Forschungen gereicht haben und wohin sie nicht mehr zu reichen vermochten; eine Frage, die man in der Literatur kaum diskutiert findet. \\

Um dem Kontext einer historischen motivierten Abhandlung gerecht zu werden, beabsichtigt der Stil der vorliegenden Abhandlung, die Euler'sche Schaffensweise möglichst getreu nachzuzeichnen und folgt entsprechend weitestgehend nicht der von Euklid vorgegebenen synthetisch--deduktiven Struktur zur Präsentation mathematischer Sachverhalte. Weiterhin sind insbesondere bei Eulers Arbeiten die Titel seiner Abhandlungen beim ersten Auftreten im Text in Originalsprache und Übersetzung desselben\footnote{Die Übersetzung wird nur bei den Titeln in lateinischer und französischer Sprache gegeben. Datum der Veröffentlichung und ggf. das Erstverfassung sind stets mit angegeben, eine Praxis, welche von den Richtlinien des \textit{Euleriana Archive} herrührt. Zusätzlich wird bei Zitierungen von Eulers Arbeiten stets die Eneström--Zahl mit aufgeführt. Eneström (1852--1923) hat eigens einen Katalog \textit{``Die Schriften Eulers chronologisch nach den Jahren geordnet, in denen sie verfasst worden sind"} (\cite{En10}, 1910) mit kurzen Übersichten zu jeder einzelnen Euler'schen Abhandlung verfasst und ihnen nach ihrer Erscheinungschronologie eine Zahl, eben die heute so genannte Eneström-Zahl, zugewiesen.} angegeben. Wie aus dem Quellenverzeichnis zu entnehmen, bilden die Euler'schen Originalabhandlungen das Fundament für die vorliegende Ausarbeitung,  was sich in einem gleichsam ``Euler'sch--solipsistischen"{} Schreibstil niederschlägt\footnote{Man wird etwa deswegen in der vorliegenden Arbeit kaum eine Definition finden, weil, wie auch herausgearbeitet werden wird, die explizite Definition von Begriffen bei Euler vergleichsweise in den Hintergrund trat.}. Das Hauptaugenmerk liegt stets bei Euler und andere Beiträge werden ergänzend genannt. In den Fußnoten haben größtenteils Erwähnungen ihren Platz gefunden, welche Antworten auf Fragen aus Diskussionen des Verfassers mit verschiedensten Interessenten der Mathematikhistorie geben, jedoch zu weit vom Haupttraktat wegführen würden.

\newpage


\section{Einleitung}
\label{sec: Einleitung}

\epigraph{History is a kind of introduction to more interesting people than we can possibly meet in our restricted lives; let us not neglect the opportunity.}{Dexter Perkins}

In diesem einleitendem Abschnitt so ein rudimentärer Überblick über bereits über Euler verfasste Arbeiten (Abschnitt \ref{subsec: Was wurde schon über Euler geschrieben?}) gegeben werden. Weiterhin wird zur Erleichterung des Verständnisses die Struktur der vorliegenden Abhandlung vorgestellt werden  (Abschnitt \ref{subsec: Aufbau dieser Arbeit}), wonach noch die Zielsetzung ihre Erwähnung finden wird (Abschnitt \ref{subsec: Zielsetzung}).

\subsection{Abfassungen über Euler und sein Opus}
\label{subsec: Was wurde schon über Euler geschrieben?}

\epigraph{Euler published so much and in so many different fields that [...] no one person will know enough to span all of his work.}{Fernando Q. Gouvêa}

Es wird zuträglich sein, dem Hauptgegenstand der gegenwärtigen Abhandlung eine kurze (und notwendig unvollständige) Übersicht einiger bereits über Euler verfasster Arbeiten, sowohl seine Mathematik und Physik als auch seine Biografie betreffend, voranzustellen.

\subsubsection{Übersicht über Arbeiten über Euler}
\label{subsubsec: Übersicht über Arbeiten über Euler}


\epigraph{Much of our knowledge is due to a comparatively few great mathematicians such as Newton, Euler, Gauss, or Riemann; few careers can have been more satisfying than theirs. They have contributed something to human thought even more lasting than great literature, since it is independent of language.}{Edward Charles Titchmarsh}

Neben der bekannten  Biographie von E. Fellman (1927--2012) mit dem Titel \textit{``Leonhard Euler"} (\cite{Fe95}, 1995) und einzelnen Übersichtsartikeln wie etwa \textit{``The Mathematics and Science of Leonhard Euler (1707–1783)"} (\cite{Th06}, 2006) von Thiele (1943--) ist das Buch \textit{``Leonhard Euler: Mathematical Genius in the Enlightenment"} (\cite{Ca15}, 2015) von Calinger dermaßen detailliert, dass sich allerhöchstens noch Fußnoten zur Ergänzung beifügen ließen. Überdies sind die Euler'schen Werke allesamt mit größter Sorgfalt editiert leicht zugänglicher Form als \textit{Opera Omnia} nach Themengebieten geordnet herausgegeben worden. Die jeweiligen Vorworte enthalten  sonst keinerorts nieder geschriebene Informationen über den Inhalt von Eulers Arbeiten. Während letztgenannte  angesichts des gewaltigen Umfangs von Eulers Werk den Charakter einer Übersicht behalten mussten, sind an vielen Orten synoptische Einzelbesprechungen von ausgesuchten Euler'schen Arbeiten publiziert worden. Hervorzuheben sind in diesem Kontext die Kolumne \textit{``How Euler did it"} von E. Sandifer (1951--2022), welche auch als Buch mit demselben Titel \textit{``How Euler Did It"}  (\cite{Sa07b}, 2007) erschienen ist, und sein Buch \textit{``The Early Mathematics of Leonhard Euler"} (\cite{Sa07a}, 2007), welches nahezu 50 Arbeiten aus der frühen Euler'schen Schaffensperiode bespricht. \\

Weiterhin sind einzelne, umfassende Werke, zu Teilgebieten erschienen, welche Euler mit seinen Beiträgen bereichert hat. Genannt sei an dieser Stelle stellvertretend das Buch von Verdun \textit{``Leonhard Eulers Arbeiten zur Himmelsmechanik"} (\cite{Ve14}, 2014) über Eulers Untersuchungen zur Himmelsmechanik.  Besonders kreative Beweise von Euler sind  ebenfalls Inhalt von Abhandlungen gewesen und bilden außerdem den Inhalt von Büchern wie dem von Dunham (1947--) \textit{``Euler: The Master of Us All"}(\cite{Du99}, 1999). \\

Angesichts der Breite des Euler'schen Opus  überrascht es wenig, dass sich darin auch Entdeckungen finden, welche seinen Nachfolgern zugeschrieben worden sind. Arbeiten, welche die Euler'sche Prioritätsansprüche geltend machen, sind eine logische Konsequenz dessen und ebenfalls zahlreich verfasst worden. Hier sei die Arbeit von Ayoub (1923--2013) \textit{``Euler and the Zeta Function"} (\cite{Ay74}, 1974) hervorgehoben, welche den Euler'schen Beitrag zur  Riemann'schen $\zeta$--Funktion zum Inhalt hat\footnote{Euler entdeckt bereits alle Eigenschaften der $\zeta$--Funktion, welche man zutage födern kann, wenn sie als Funktion einer reellen Variable  behandelt wird. Riemann (1826--1866) reicht in seiner berühmten Arbeit (\textit{``Über die Anzahl der Primzahlen unter einer gegebenen Größe"} \cite{Ri60}, 1860, ges. 1859) strenge Beweise für die Euler'schen Entdeckungen nach. Dies wird alles in Abschnitt (\ref{subsubsec: Durch Kombinieren von Ergebnissen: Die zeta-Funktion}) detailliert ausgeführt werden.}. Weiterhin haben sich Forscher wie Polya (1887--1985) der Präsentation mathematischen Arbeitsweise angenommen und dabei Eulers Ausführungen zur Illustration herangezogen\footnote{Die Methoden findet man unter anderem in seinem Buch \textit{``Mathematics and Plausible Reasoning"} (\cite{Po14}, 2014).}. Das Buch  von Suisky \textit{``Euler as Physicist"} (\cite{Su10}, 2010)  diskutiert diverse Aspekte des Euler'schen Beitrags zur Physik, wobei insbesondere Eulers Arbeit \textit{``Anleitung zur Naturlehre"} (\cite{E842}, 1862, ges. ca. 1750) im Vordergrund steht.\\

\subsection{Aufbau dieser Arbeit}
\label{subsec: Aufbau dieser Arbeit}


\epigraph{Information is a source of learning. But unless it is organized, processed, and available to the right people in a format for decision making, it is a burden, not a benefit.}{William Pollard}

Die vorliegende Ausarbeitung  gliedert sich in drei Hauptteile. Der erste Teil ist ganz den Euler'schen Entdeckungen gewidmet. Allem voran wird die Euler'sche Ansicht auf die Mathematik betrachtet (Abschnitt \ref{sec: Eulers Ansicht zur Mathematik}); mit einen besonderen Fokus auf seine Auffassung eines Beweises (Abschnitt \ref{subsubsec: Eulers Auffassung eines Beweises}). Die Rolle der Definition und anderer Konzepte werden in den historischen Kontext eingeordnet (Abschnitt \ref{subsubsec: Die Rolle der Definition in der Mathematik}).\\

Nach Präsentation der Euler'schen Mathematik--Philosophie folgt die Darstellung seines Schaffensprozesses anhand ausgewählter Beispiele (Abschnitt \ref{sec: Eulers Arbeitsweise anhand ausgewählter Beispiele}). Das bekannte Beispiel der Lösung des Baseler Problems (Abschnitt \ref{subsec: Die Lösung des Baseler Problems}) und  Eulers in der Literatur weniger diskutierte allgemeine Lösung der einfachen Differenzengleichung (Abschnitt \ref{subsec: Ein Prototyp -- Differenzengleichungen}) werden es gestatten, anhand einzelner Abhandlungen unmittelbar zahlreiche charakteristische Facetten seiner Arbeitsweise zu illustrieren.  So werden  sowohl sein nahezu unermesslicher Ideenreichtum als auch seine -- mathematisch gesprochen -- Waghalsigkeit, die ihn bisweilen auch zu falschen Ergebnissen leitet, zutage treten (Abschnitt \ref{subsec: Ein falsches Ergebnis}). \\

Der nächste Abschnitt der vorliegenden Abhandlung nimmt sich der Klärung von Prioritätsfragen an (Abschnitt \ref{sec: Von Euler vorweggenommene Entdeckungen}); genauer werden einige  Ergebnisse ihre Erwähnung finden, die Euler bereits zutage gefördert hat, jedoch seinen Nachfolgern zugeschrieben worden sind. Unterteilt werden die Resultate entsprechend der verschiedenen Ursachen für diese jeweiligen Missattributionen: Einmal ereignete es sich, dass Euler einen Schatz zwar als erster gehoben hatte, selbigen mit einem Beweis untermauert hat, und die Entdeckung gar in heute üblicher Form niedergeschrieben hat, seine Nachfolger jedoch von seinen Arbeiten keine Kenntnis hatten. Ein Exempel dessen stellt die Formel die Fourierkoeffizienten (Abschnitt \ref{subsubsec: Die Fourierkoeffizienten}) dar. Bei anderen Begebenheiten präsentierte Euler ebenfalls einen Beweis für eine Entdeckung, allerdings findet sich das Ergebnis nicht in der modernen Gestalt in seinen Arbeiten und bedarf erst einer Einkleidung in das Gewand der modernen Sprache der Mathematik, um ersichtlich zu werden. Die Unterschiede der Euler'schen Darstellungen reichen über die reine Form (Abschnitt \ref{subsubsec: Die Multiplikationsformel für die Gammafunktion}), zu einem anderen Kontext (Abschnitt \ref{subsubsec: Den Kontext betreffend: Die Legendre Polynome}) und  anderen Fragestellungen (Abschnitte \ref{subsubsec: Die Mellin--Transformierte bei Euler} und \ref{subsubsec: Ein anderes Vorhaben: Das Weierstraß-Produkt}) bis hin zu einer anderen Präsentation (Abschnitt \ref{subsubsec: Die Darstellung betreffend -- Die hypergeometrische Reihe}).\\

Weiterhin hat Euler Theoreme formuliert, die er, wie er an entsprechender Stelle auch einräumt, nicht mit einem Beweis unterlegt hat (Abschnitt \ref{sec: Von Euler Entdecktes, aber nicht Bewiesenes}). Hier lassen sich seine Resultate weiter  zwei Unterkategorien zuordnen. Die eine bilden diejenigen Funde, in welchen es Euler freilich selbst vermocht hätte, noch zur vermissten rigorosen Bestätigung zu gelangen. Die nachzureichenden Begründungen erstrecken sich dabei von einer einfachen Anwendung einer von ihm anderenorts bewiesenen Methode (Abschnitt \ref{subsubsec: Durch Anwendung einer Methode: Seine Konstante A}) weiter über eine leichte Kombination seiner anderenorts, meist in anderen Kontexten, erzielten Ergebnissen (Abschnitt \ref{subsubsec: Durch Kombinieren von Ergebnissen: Die zeta-Funktion}) bis hin zu einem einfachen Gedankengang (Abschnitt \ref{subsubsec: Durch einen neuen Gedanken: Die Theta-Funktion}), welchen er zweifelsohne leicht hätte haben können. Die andere Unterkategorie bilden indes Entdeckungen, welche Euler zwar gemacht hatte und nachzuweisen suchte, jedoch das Fehlen der Mittel  (Abschnitt \ref{subsubsec: Wegen fehlender Mittel: Das Reziprozitätsgesetz}) oder seine Formulierung des Sachverhaltes (Abschnitt \ref{subsubsec: Weg fehlender Formulierung: Der Primzahlsatz}) oder eine übersehene, jedoch für Euler schwer zu erkennende, Unvollständigkeit in der Argumentation (Abschnitt \ref{subsubsec: Wegen übersehener Unvollständigkeit: Der große Satz von Fermat für n=3}) dies nicht zuließ.\\

Diese letzte Unterkategorie bildet gleichzeitig den Übergang zum nächsten Teil, welcher Grenzen Eulers zum Gegenstand hat (Abschnitt \ref{sec: Grenzen Eulers gezogen durch das Paradigma}). Einmal wurde diese Schranke von den  Begriffen selbst (wie dem der Funktion (Abschnitt \ref{subsubsec: Der Begriff der Funktion})), des Grenzwertes  (Abschnitt \ref{subsubsec: Der Begriff des Grenzwerts}), der Summe (von divergenten) Reihen (Abschnitt \ref{subsubsec: Der Begriff der Summe einer Reihe})) gesetzt, welche bei Euler noch nicht die moderne Gestalt hatten. Abseits dieser mathematischen Hinderungsgründe  vereitelte bei anderen Begebenheiten  seine eigene Herangehensweise an Problemstellungen ein weiteres Fortschreiten. Denn letztere veranlasste Euler bisweilen zu einer eingeschränkten Sichtweise  oder zur Ausbildung einer irreführenden Frage  (Abschnitt \ref{subsec: Durch eine falsche Frage}). Einmal wählte Euler eine nicht zielführende Kategorisierung für die Objekte seines Interesses (Abschnitt \ref{subsubsec: Die Normalform von elliptischen Integralen}), ein anderes Mal versuchte er etwas zu beweisen, von dem man erst nach ihm wusste, dass es unmöglich ist (Abschnitt \ref{subsubsec: Wegen Unbeweisbarkeit: Wurzeln von Polynomen}). Weshalb Euler die Entwicklung mancher Techniken seinen Nachfolgern überlassen hat, wird in der folgenden Sektion (Abschnitt \ref{subsec: Grenzen durch den eigenen Arbeitsethos}) im Kontext verschiedener Themenkomplexe erläutert werden. \\

Im letzten Teil (Abschnitt \ref{sec: Herleitungen aus Eulers Formeln und Ideen}) werden die  vorgestellten Euler'schen Ergebnisse ihre Anwendung finden.  So wird die Anzahl der Beweise für der Werte  der Potenzsummen der Reziproken der natürlichen Zahlen  um eine Einheit erhöht (Abschnitt \ref{subsubsec: Herleitung Formel für die Potenzsummen der Reziproken})  sowie weitere Spezialfälle seiner Ergebnisse zu Legendre--Polynomen präsentiert  (Abschnitt \ref{subsubsec: Weitere Untersuchungen zu den Legendre-Polynomen}) und  Formeln für weitere orthogonale Polynome mitgeteilt (Abschnitt \ref{subsubsec: Explizite Formeln aus Eulers Theorie zu Differenzengleichungen}).  Abschließend soll eine Verbindung zwischen Euler und Ramanujan (1887--1920), welche oft wegen ihrer ausgeprägten mathematischen Intuition verglichen worden sind, hergestellt werden. Zum einen geschieht dies über die Herleitung von zwei bekannteren Integralidentitäten Ramanujans  (Abschnitt \ref{subsubsec: Ein bestimmtes Integral von Ramanujan} und \ref{subsubsec: Ramanujans Mastertheorem}), zum anderen über die berühmten Formeln von Ramanujan zur Berechnung der Kreiszahl $\pi$, welche sich in Euler'scher Manier finden lassen (Abschnitt \ref{subsubsec: Ramanujans Formeln zur Kreisquadratur}). Die Zusammenfassung (Abschnitt \ref{sec: Zusammenfassung}) wird den Schlusspunkt dieser Arbeit bilden.\\

\subsection{Zielsetzung dieser Arbeit}
\label{subsec: Zielsetzung}

\epigraph{No matter what you do, you first have to have a vision [...]. If you don't have a goal or a vision, then you have nothing.}{Arnold Schwarzenegger}

Das ausgeschriebene Ziel der hiesigen Ausführungen besteht im Beitrag zum umfassenderen Verständnis des Euler'schen Werks. Angestrebt wird dies über die Vermittlung  seiner Auffassungen zur Mathematik und die sich daraus ableitenden Konsequenzen für seinen Arbeitsprozess. Überdies wird durch die Mitteilung weniger diskutierter Ergebnisse Eulers ein Beitrag zur Geschichte der Mathematik beabsichtigt. Als Nebenprodukt bleibt die Hoffnung, dass sich aus der  Beschäftigung mit dem Euler'schen Opus der Beantwortung von Fragen aus anderen Wissenschaften näher kommen lässt. Dies zielt insbesondere auf Fragen aus der Psychologie wie derjenigen, welche Bedingungen für kreatives Schaffen notwendig sind\footnote{Hierzu wurde trotz intensiver Bemühungen auf verschiedensten Wegen keine zufriedenstellende bzw. allgemein in der Psychologie akzeptierte Antwort gefunden; man siehe etwa die Ausführungen im Buch \textit{``Thought and Knowledge: An Introduction to Critical Thinking"}  (\cite{Ha22}, 2022). Da Euler wie kaum ein anderer seines Leistungsvermögens seinen Entdeckungsprozess samt seiner Misserfolge  in seinen Arbeiten kund getan hat, ist diese Hoffnung auf die Aufhebung dieses Umstandes nicht allzu unbegründet.}.

\newpage

\section{Eulers Ansicht zur Mathematik}
\label{sec: Eulers Ansicht zur Mathematik}

\epigraph{If there is a God, he's a great mathematician.}{Paul Dirac}

Für die weitere Progression der hiesigen Ausführungen ist eine Diskussion dessen, was für Euler die Essenz der Mathematik bildet, unabdingbar. Dieser Besprechung ist der erste Teil dieser Sektion gewidmet (Abschnitt \ref{subsec: Eulers Mathematikphilosophie}), wohingegen der zweite Teil anhand kleinerer Beispiele Besonderheiten seiner Darstellungen präsentiert (Abschnitt \ref{subsec: Abriss von Eulers Vorgehensweise}).

\subsection{Eulers Mathematikphilosophie}
\label{subsec: Eulers Mathematikphilosophie}

\epigraph{Mathematics is not just a language. Mathematics is a language plus reasoning.}{Richard P. Feynman}

Unter dem Begriff der \textit{Mathematikphilosophie} soll im folgenden stets das gesammelte Bündel an expliziten und impliziten Auffassungen und Konzepten dessen verstanden werden, was dem Begriff \textit{Mathematik} zugehörig ist.  Es wird im folgenden zunächst das diskutiert, was Mathematik für Euler bedeutet (\ref{subsubsec: Sein Essay zur Nützlichkeit der höheren Mathematik}) und die Vorrangstellung  der Analysis bei Euler herausgestellt (Abschnitt \ref{subsubsec: Euler Auffassung der Analysis}), bevor dann auf das fundamentale Konzept der Definition (Abschnitt \ref{subsubsec: Die Rolle der Definition in der Mathematik}), das Verständnis eines Beweises (Abschnitt \ref{subsubsec: Eulers Auffassung eines Beweises}) und die Methode der Induktion bei Euler (\ref{subsubsec: Beweistechnik der Induktion bei Euler}) eingegangen wird. Den Abschluss bildet die Diskussion von Eulers Ansicht auf die Physik (Abschnitt \ref{subsubsec: Eulers Ansicht zur Physik}).

\subsubsection{Eulers Auffassung zum Wesen der Mathematik}
\label{subsubsec: Sein Essay zur Nützlichkeit der höheren Mathematik}


\epigraph{God used beautiful mathematics in creating the world.}{Paul Dirac}

Seine Auffassung vom Wesen der Mathematik setzt Euler in seinem Text \textit{``Commentatio de matheseos sublimioris utilitate"} (\cite{E790}, 1847, ges. 1741) (E790: ``Kommentar über die Nützlichkeit der höheren Mathematik") auseinander. Zunächst betont er den praktischen Nutzen der Mathematik. Die diskutieren Gebiete umfassen die klassische Mechanik, die Hydromechanik, Artillerie, Schifffahrt und die Physik, der gesonderten Besprechung welcher Anwendungsbereiche Euler folgende Passage voranstellt\footnote{Es folgt die Übersetzung des ursprünglich von Euler in Latein verfassten Textes, welche Burckhardt (1903--2006) für die Serie 3, Band 2 der \textit{Opera Omnia} (\cite{OO32}, 1942) ins Deutsche übertragen hat. Man findet seine ganze Übersetzung auf den Seiten 408--416 des Buches.}:\\

\textit{``Weil aber die ganze Mathematik im Auffinden unbekannter Größen besteht, und zu diesem Zweck entweder Methoden oder gleichsam zur Wahrheit führende Wege eröffnet, oder äußerst tief verborgene Wahrheiten ausfindig macht und klar beleuchtet, von welchen zum einen die Kraft unseres Geistes geschärft wird, zum anderen unsere Erkenntnis erweitert wird, kann gewiss auf keine der beiden zu wenig Mühe aufgewandt werden. Weil nämlich die Wahrheit nicht nur an sich sehr lobenswert ist, sondern auch wegen des großen Zusammenhangs, in welchem alle Wahrheiten zusammenhängen, nicht ohne Nutzen sein kann, auch wenn selbiger nicht sofort erkannt wird, ist jener Einwand, dass die höhere Mathematik, allzu sehr von der Wahrheitsfindung eingenommen ist, in einen Lobesauspruch umgekehrt, als ein Tadel zu sein."}\\

Demnach sieht  Euler die Mathematik als \textit{notwendig und hinreichend}  für die Wahrheitsfindung, wobei die Wahrheit das oberste Gut ist, welches es zu finden gilt. Erst  hieraus ergeben sich letztendlich die ganzen praktischen Anwendungen in den anderen Wissenschaften. Von Interesse ist auch die Auffassung, dass letztendlich alle Wahrheiten\footnote{Weiter unten (Abschnitt \ref{subsubsec: Eulers Auffassung eines Beweises}) werden noch die verschiedenen Arten von Wahrheiten, welche Euler hierunter versteht, ihre Erwähnung finden.}, welche mithilfe der Mathematik der Kenntnis zugeführt werden, zusammenhängen. Dies erklärt zugleich Eulers Motivation, für eine Behauptung in der Mathematik, ob der eigentlichen Redundanz, meist mehrere Begründungen zu liefern, sofern ihm dies freilich möglich war. \\

\subsubsection{Eulers Auffassung der Analysis}
\label{subsubsec: Euler Auffassung der Analysis}

\epigraph{The person who did most to give to analysis the generality and symmetry which are now its pride, was also the person who made Mechanics analytical; I mean Euler.}{William Whewell}

Kaum jemand hat die Analysis dermaßen geprägt und gleichsam personifiziert wie Euler\footnote{So begründet sich es wohl auch, dass Eulers Portrait die Vorderseite des Buches \textit{``3000 Jahre Analysis"} (\cite{So16}, 2016) ziert.}, was sich im folgenden Zitat von Hermann Hankel (1839--1873) auf Seite 14 in seinem Buch \textit{``Die Entwickelung der Mathematik in den letzten Jahrhunderten"} (\cite{Ha69}, 1869) widerspiegelt:\\

\textit{``Es ist das unschätzbare Verdienst des grossen Baseler Mathematikers Leonhard Euler, den analytischen Calcul von allen geometrischen Fesseln befreit zu haben, und damit die Analysis, als selbstständige Wissenschaft begründet zu haben."}\\

\paragraph{Zwei kleinere Erfolgsbeispiele seiner Auffassung}
\label{para: Zwei kleine Erfolgsbeispiele seiner Auffassung }

Euler selbst äußert des Öfteren die Ansicht, die Analysis finde in jedem Feld der Mathematik, wenn auch nicht unmittelbar, jedoch eventuell ihre Anwendung.  In § 2 seiner Abhandlung \textit{``De numeris amicibilicus"} (\cite{E152}, 1752, ges. 1747) (E152: ``Über befreundete Zahlen")  formuliert er diesbezüglich die Aussage:\\

\textit{``Weil aber außer Zweifel steht, dass die Analysis auch aus diesem Zweig [der Mathematik] nicht zu verachtende Zuwächse nehmen wird, wenn eine Methode erschlossen wird, mit welcher sich um vieles mehr Paare von dieser Art ausfindig machen lassen, halte ich es nicht für unpassend, gewisse sich hierauf beziehende Methoden, auf welche ich zufällig gestoßen bin, hier mitzuteilen."}\\

Die Paare dieser Art meinen dabei Paare von sogenannten befreundeten Zahlen; dies sind solche Zahlenpaare, von welchen die Summe der echten Teiler der einen Zahl der jeweils anderen Zahl gleich ist. Das kleinste Beispiel befreundeter Zahlen besteht aus den Zahlen $220$ und $284$. In seiner Abhandlung führt Euler insbesondere die Teilerfunktion $\sigma$ ein, welche einer natürlichen Zahl die Summe ihrer echten Teiler zuordnet\footnote{Euler benutzt das Symbol $\int$, womit er  explizit machen will, dass die Funktion eine Summe ist.}. Weiterhin weist er die grundlegenden Eigenschaften

\begin{equation*}
    \sigma(p^k)= \dfrac{p^{k+1}-1}{p-1} \quad \text{für} \quad p \in \mathbb{P}, k\in \mathbb{N},
\end{equation*}
wobei $\mathbb{P}$ die Menge der Primzahlen anzeigt, sowie

\begin{equation*}
    \sigma(m \cdot n)= \sigma(m) \cdot \sigma(n)
\end{equation*}
für teilerfremde Zahlen $m$, $n$ nach. Euler zeigt demnach, dass $\sigma$ eine zahlentheoretische Funktion ist, wie Funktionen mit den just mitgeteilten Beschaffenheiten heute genannt werden. Diese verwendet er im Folgenden, um durch geschickte Ansätze und Probieren weitere befreundete Zahlen zu finden. Am Ende seiner Arbeit \cite{E152} hat er die Menge der Paare von den 3 bekannten auf 61 erhöht\footnote{Man vergleiche dazu Sandifers Artikel \textit{``Amicable numbers"} (\cite{Sa05nov}, 2005).}.\\

Das eingangs mitgeteilte Euler'sche Zitat ist bemerkenswert, weil die referierte Arbeit sich einem Problem aus der Zahlentheorie annimmt, einem Gebiet dessen systematische Entwicklung  Euler wesentliche Impulse verdankt. Im weiteren Verlauf der vorliegenden Ausführungen wird sich zeigen, dass diese ``Analytisierung"{} Euler auf der einen Seite weit vordringen lässt, ihm aber andererseits auch Grenzen setzt, die er mangels der dafür notwendigen Konzepte nicht überwinden kann.\\

Auch in der Geometrie lässt Euler die analytische Fassung eines Problems Entdeckungen machen, die seinen Vorgängern verschlossen geblieben waren.  Obwohl etwa umfassende Untersuchungen über  Dreiecke viele ihrer Eigenschaften enthüllten, bedurfte es zur Entdeckung der Euler'schen Geraden seines in der Arbeit \textit{``Solutio facilis problematum quorundam geometricorum difficillimorum"} (\cite{E325}, 1767, ges. 1763) (E325: ``Eine leichte Lösung gewisser sehr schwieriger geometrischer Probleme") gewählten analytischen Ansatzes. Einen rein geometrischen Beweis der Tatsache, dass sich Höhenschnittpunkt, Schnittpunkt der Mittelsenkrechten und der Seitenhalbierenden stets auf einer Geraden --  der heute so genannten Euler'schen Geraden -- befinden, ist vermutlich mit noch höheren Schwierigkeiten verbunden als die analytischen Rechnungen Eulers, welche man in Sandifers Artikel \textit{``The Euler line"} (\cite{Sa09jan}, 2009) und auch in Dunhams Buch \textit{``Euler: The Master of Us All"} (\cite{Du99}, 1999) nachgezeichnet findet.

\paragraph{Was für Euler nicht zur Mathematik gehörte}
\label{para: Was für Euler nicht zur Mathematik gehört}

Nach Eulers Auffassung,  Mathematik bestehe im Auffinden von unbekannten Größen (siehe das Zitat aus Abschnitt \ref{subsubsec: Sein Essay zur Nützlichkeit der höheren Mathematik}), erscheint folgende Äußerung aus Eulers Arbeit \textit{``Solutio problematis ad geometriam situs pertinentis} (\cite{E53}, 1741, ges. 1735) (E53: ``Lösung eines sich auf die Geometriam Situs beziehenden Problems") zum berühmten Königsberger Brückenproblem selbiges von der Mathematik auszuschließen. Er schreibt direkt in § 1 zu dem damals von Leibniz als \textit{Geometriam Situs} bezeichneten Teilbereich der Mathematik:\\

\textit{``So wird dieser Teilbereich [der Mathematik] von selbigem [Leibniz] allein  von der Bestimmung der Lage und dem Eruieren der Beschaffenheiten des Ortes eingenommen festgelegt, bei welchem Unterfangen weder auf die Größen zu achten noch ein Kalkül mit Größen zu gebrauchen ist. [...] Deswegen, nachdem neulich die Erwähnung eines gewissen Problems gemacht worden ist, das sich auf die Geometrie zu beziehen schien, aber so beschaffen war, dass es weder die Bestimmung von Größen erfordert, noch eine Lösung mithilfe des Kalküls von Größen zuließ, habe ich nicht gezweifelt, es zur  Geometriam Situs zu zählen."}\\

\begin{figure}
    \centering
   \includegraphics[scale=0.7]{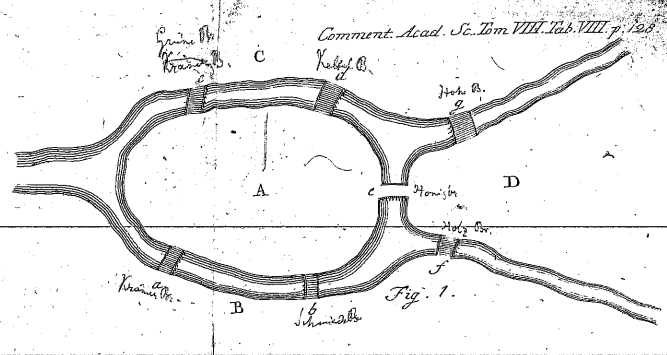}
    \caption{Eulers Figur zum Königsberger Brückenproblem aus \cite{E53}. Die Aufgabe besteht im Überschreiten jeder der sieben Brücken (mit Kleinbuchstaben gekennzeichnet), ohne eine von ihnen doppelt zu passieren. Euler zeigt im Verlauf die Unmöglichkeit dieser Aufgabe auf.}
    \label{E53Brücken}
\end{figure}

Das ``gewisse Problem"{} meint bereits das zuvor angesprochene Problem der Brücken von Königsberg, welches  häufig angesehen wird, die Topologie als Teildisziplin der Mathematik begründet zu haben. Euler hingegen hätte die Frage wohl nicht als eine mathematische eingeordnet. Jedenfalls ist die Arbeit \cite{E53} das einzige Euler'sche, was sich aus moderner Sicht der Graphentheorie zuordnen lässt.

\subsubsection{Die Rolle der Definition in der Mathematik bei Euler}
\label{subsubsec: Die Rolle der Definition in der Mathematik}

\epigraph{What is a good definition? For the philosopher or the scientist, it is a definition which applies to all the objects to be defined, and applies only to them; it is that which satisfies the rules of logic. But in education it is not that; it is one that can be understood by the pupils.}{Henri Poincaré}

Heute sind Definitionen in der Mathematik ein wesentlicher Bestandteil und überdies wurden eigens Bücher zum Gegenstand zur Natur des Begriffs und Gebrauchs der Definition verfasst. Man vergleiche etwa das Lehrbruch von Dubislav\footnote{Dubislav hat insbesondere Mathematik studiert, einer seiner Lehrer war David Hilbert (1862--1943).} (1895--1937) \textit{``Die Definition"}  (\cite{Du94}, 1994). Zu Eulers Lebenszeiten lag das Hauptaugenmerk vielmehr auf den Anwendungen und das implizite Verständnis eines Begriffs ließ die explizite Ausformulierung von Begriffen entsprechend in den Hintergrund treten.  Diese Situation wird treffend von Hardy (1877--1947)  in seinem Buch \textit{``Divergent Series"}  (\cite{Ha49}, 1949)  zusammengefasst. Man findet auf Seite 15 folgende Ausführung:\\

\textit{``It is a mistake to think of Euler as a `loose' mathematician, though his language sometimes might seem loose to modern ears; and even his language sometimes suggests a point view far in advance of the general ideas of his time. [...] It is impossible to state Euler's principle accurately without clear ideas about functions of a complex variable and analytic continuation."}\\

Euler erfasst mit seinen Ideen, von Hardy beschrieben am Konzept der divergenten Reihen, zwar oft den Kern der Sache, jedoch bedurfte es der genauen Ausarbeitung entsprechender Konzepte, hier unter anderem dem der analytischen Fortsetzung von Funktionen einer komplexen Variablen, um sie auf solidere Füße zu stellen.\\

Indes hat Euler auch selbst den Umgang mit Begriffsbildung und Definitionen in seinem Artikel \textit{``De usu functionum discontinuarum in analysi} (\cite{E322}, 1767, ges. 1762) (E322: ``Über den Gebrauch von unstetigen Funktionen in der Analysis) thematisiert und bei dieser Gelegenheit bemängelt. Hier schreibt er am Ende von § 6:\\

\textit{``Daher entbehren die sehr häufigen Beschwerden, dass die Idee der Analysis des Unendlichen niemals genau beschrieben und fundamental begründet gefunden wird, nicht jedweder Grundlage, [...]."}\\

Die zitierte Textpassage bezieht sich vornehmlich auf die Natur von Differentialen, welche damals auch von Euler selbst in seinem Lehrbuch \textit{``Institutiones calculi differentialis cum eius usu in analysi finitorum ac doctrina serierum"} (\cite{E212}, 1755, ges. 1748) (E212: ``Grundlagen des Differentialkalküls mit seinem Gebrauch in der Analysis des Endlichen und der Reihenlehre") noch als unendlich kleine Größen betrachtet werden, was teilweise zu paradoxen Schlussfolgerungen führt\footnote{Dies wird weiter unten Abschnitt \ref{subsubsec: Der Begriff des Grenzwerts} noch explizit besprochen und mit entsprechenden Beispielen unterlegt werden.}. Zusammengefasst wird dies treffend von Hardy (in einen anderen Zusammenhang) auf S. 5--6 von \cite{Ha49}:\\

\textit{``[I]t does not occur to a modern mathematician that a collection of mathematical symbols should have a 'meaning' until one has been assigned to it by definition. It was not a triviality even to the greatest mathematicians of the eighteenth century. They had not the habit of definition: it was not natural to them to say, in so many words, by $X$ we mean $Y$."}

\subsubsection{Eulers Auffassung eines Beweises}
\label{subsubsec: Eulers Auffassung eines Beweises}

\epigraph{Mathematics is not a deductive science - that's a cliché. When you try to prove a theorem, you don't just list the hypotheses, and then start to reason. What you do is trial and error, experimentation, guesswork.}{Paul Halmos}

Die Konzeption eines Beweises nimmt das Primat der Begriffe in der Mathematik ein. Wohingegen sich die moderne Mathematik heute auf die Beweistheorie\footnote{Man vergleiche etwa das Buch \textit{``Handbook of mathematical logic"} (\cite{Ba77}, 1977) zum Gebiet der mathematischen Logik und das Buch \textit{``Basic Proof Theory"} (\cite{Po08}, 2008), welches eigens der Beweistheorie in der Mathematik gewidmet ist.} zurückblicken kann, lag selbige zu Eulers Lebzeiten noch nicht vor. Nichtsdestoweniger besitzt Euler ebenfalls eine Auffassung eines Beweises, welche er jedoch selten explizit erläutert zu haben scheint, aber natürlich implizit in seinen Untersuchungen gebraucht\footnote{Ab welchem Punkt ein Beweis als solcher anzuerkennen ist, bildet immer noch den Gegenstand vieler Diskussionen über mehrere Wissenschaftsgebiete hinweg.}. Da die Euler'schen Ansichten mancherorts von ihren modernen Gegenstücken abweichen, wird eine Analyse seiner Äußerungen zu diesem Gegenstand an dieser Stelle eingeschaltet.

\paragraph{Allgemeine Ansichten}
\label{para: Allgemeine Ansichten}

In seinen \textit{``Lettres à une princesse d'Allemagne sur divers sujets de physique et de philosophie -- Tome second"} (\cite{E343}, 1768, ges. 1760--1761)  (E343: ``Briefe an eine deutsche Prinzessin über verschiedene Gegenstände der Physik und Philosophie -- Buch 2") setzt Euler neben vielen anderen Themen seine Theorie des Wissenserwerbs und die damit verknüpften Wahrheitsbegriffe sowie Beweismethodiken auseinander. Seine Erläuterungen diesbezüglich finden sich in den Briefen 115--120. In Brief 115 unterscheidet Euler nämlich drei Wahrheitsarten, welche für ihn gleichzeitig auch die einzigen drei sind: Die Erfahrungswahrheiten (\textit{``Verités des sens"}), die Vernunftswahrheiten (\textit{``Verités de l'entendement"}) und Glaubenswahrheiten (\textit{``Veritès de la foi"}). Dies mit der gängigen Definition des Terminus \textit{Wissen} als gerechtfertigte, wahre Meinung\footnote{Man siehe zum Beispiel das einführende Lehrbuch \textit{``Einführung in die theoretische Philosophie"} (\cite{Hü15}, 2015) für umfassendere Erläuterungen.} vergleichend, sieht man die Teilfacetten der modernen Auffassung des Wissensbegriffs von der Euler'schen Aufteilung abgedeckt. Während jedoch heute Wissen aus diesen Teilen konstituiert verstanden wird\footnote{Dass die Begriffe Überzeugung, Wahrheit und Meinung lediglich notwendig, jedoch keinesfalls hinreichend für Wissen sind, ist von Gettier (1927--2021) in seiner Arbeit \textit{``Is Justified True Belief Knowledge?"} (\cite{Ge63}, 1963) aufgezeigt worden. Gettier demonstriert, dass  nebst der drei obigen Begriffe für Wissen der Zufall ausgeschlossen werden müsste. Das daraus resultierende heute sogenannte Gettier--Problem, wie dies zu bewerkstelligen ist, ist eine offene Frage.}, betont Euler explizit die strikte Verschiedenheit seiner drei Wahrheitsgattungen, welche dementsprechend jeweils eigener Beweismethoden bedürfen. In Brief 119 erläutert Euler, die Erfahrungswahrheiten verlangen physische Gewissheit (\textit{``certitude physique"}), die Vernunftswahrheiten logische oder demonstrative Gewissheit (\textit{``certitude logique ou démonstrative"}) und die Glaubenswahrheiten schließlich moralische Gewissheit (\textit{``certitude morale"}). Die Mathematik\footnote{In seinen Briefen schreibt er den damals noch gleichwertig zu verstehenden Begriff \textit{``Géometrie"}.} ordnet Euler in die zweite Klasse ein, sagt demnach, man habe sich von der Gültigkeit mathematischer Wahrheiten mithilfe der Logik zu überzeugen. Die Art des logischen Schließens hatte er zuvor in seinen Briefen 102--109 auseinander gesetzt, wo er unter anderem auch die heute nach ihm benannten Diagramme zur Illustration von Syllogismen heranzieht\footnote{Euler--Diagramme werden heute verwendet, um mengentheoretische Beziehungen zu veranschaulichen, wohingegen Venn--Diagramme alle möglichen Mengenrelationen in die Betrachtung mit einbeziehen. Letztere erlauben im Gegensatz zu ersteren ebenfalls eine Inferenz auf gegebenenfalls fehlende Zusammenhänge. Der Leser ist diesbezüglich auch auf Sandifers Artikel \textit{``Venn Diagrams"} (\cite{Sa04jan}, 2004) aus seiner Kolumne verwiesen.}. Dazu passend schreibt Euler schreibt im Vorwort (Seite  XII) zu seinen \textit{Calculi Differentialis}  \cite{E212} zur Mathematik:\\

\textit{``Obschon nämlich demjenigen, der die Bestimmung der Größe der gesamten Erdkugel mit dem Kalkül in Angriff genommen hat, ein Fehler, nicht nur von einem sondern vielleicht gar mehreren Steinchen, leicht nachgegeben zu werden pflegt, scheut indes die mathematische Strenge auch vor dem kleinsten Fehler zurück, und allzu schwer wäre dieser Einwand, sollte jedweden Bestand behalten."}\\

\paragraph{Ein explizites Exempel}

Die Mathematik ist demnach für Euler durch ihre Exaktheit geprägt, während  folgende Äußerung aus seiner Arbeit \textit{``Remarques sur un beau rapport entre les series des puissances tant directes que reciproques"} (\cite{E352}, 1768, ges. 1749) (E352: ``Bemerkung zur schönen Beziehungen über die Reihen der direkten und reziproken Potenzen") dazu widersprüchlich anmutet, wo er in § 12  formuliert:\\

\textit{``Siehe nun diesen neuen Beweis, welcher mit dem vorherigen zusammengenommen schon als vollständiger Beweis unserer Vermutung angesehen werden kann. Nichtsdestotrotz ist es überaus gerechtfertigt, noch einen direkten Beweis zu fordern, welcher alle möglichen Fälle auf einmal abhandelt."}\\

Der Zusammenhang ist hierbei folgender: In besagter Arbeit gelangt Euler durch findiges Hantieren mit divergenten Reihen\footnote{Seine Auffassung divergenter Reihen wird in Abschnitt (\ref{subsubsec: Der Begriff der Summe einer Reihe}) diskutiert.} zu folgender Gleichung:

\begin{equation}
\label{eq: Zeta-Euler}
    \dfrac{1-2^{n-1}+3^{n-1}-4^{n-1}+5^{n-1}-6^{n-1}+\text{etc.}}{1-2^{-n}+3^{-n}+4^{-n}+5^{-n}-6^{-n}+\text{etc.}} 
\end{equation}
\begin{equation*}
    = \dfrac{-1\cdot 2 \cdot 3 \cdots (n-1)\cdot (2^n-1)}{(2^{n-1}-1)\cdot \pi^n}\cdot \cos \dfrac{n\pi}{2},
\end{equation*}
welche für ihn zunächst einmal für die natürlichen Zahlen Geltung haben soll\footnote{Die Euler'sche Vermutung stellt die Funktionalgleichung der Riemann'schen $\zeta$-Funktion dar. Weiter unten (Abschnitt \ref{para: Beweis der Funktionalgleichung aus Eulers Formeln}) wird ein Beweis nachgereicht, der von Euler selbst hätte geführt werden können.}. Mit den ``Beweisen"{} drückt Euler in diesem Zusammenhang aus, dass er die Richtigkeit der Formel überdies explizit für die speziellen Fälle $n=0$ und $n=1$ erhärtet hat, in welchen die Summen auf der linken Seite anderswoher bekannt sind\footnote{Euler weiß, dass $1-\frac{1}{2}+\frac{1}{3}-\frac{1}{4}+\text{etc.}=\log(2)$ und bedient sich des Wertes $1-1+1-1+1-1+\text{etc.}=\frac{1}{2}$. Wie Euler zu letzterem Wert gelangt ist, wird unten (Abschnitt \ref{subsubsec: Der Begriff der Summe einer Reihe}) erläutert werden.}. Unter Verwendung der l'Hospital'schen Regel kann Euler jeweils rechnerisch verifizieren, dass die Formel auf der rechten Seite mit dem Quotienten der Reihenwerte auf der linken übereinstimmt. Mit seiner Auffassung der divergenten Reihen kann Euler seine Vermutung insgesamt für die ganzen Zahlen validieren. An anderer Stelle (§ 16), nachdem er auch noch den Wert $n=\frac{1}{2}$ numerisch überprüft hat,  in selbiger Arbeit äußert er des Weiteren, dass eine richtige Formel wohl kaum zufällig zu so vielen richtigen Ergebnissen führen könnte, was also umgekehrt für ihre Wahrheit spricht. Genauer schreibt er:\\

\textit{``Weil unsere Vermutung nun zum höchsten Grad an Gewissheit erhoben worden ist und nicht einmal mehr ein Zweifel in den Fällen bestehen bleibt, in denen man für den Exponenten $n$ Brüche einsetzt, wollen wir [...]"}\\

Euler demonstriert noch für $n=\dfrac{3}{2}, \dfrac{5}{2}, \dfrac{7}{2}$, $\dfrac{9}{2}$ und $\dfrac{11}{2}$ die Übereinstimmung und führt demnach insgesamt einen regressiven Beweis für die Richtigkeit der Gleichung (\ref{eq: Zeta-Euler}), welcher heute in der Mathematik nicht mehr als Beweis herangezogen werden kann. Euler verleiht indes seiner Argumentation im letzten Paragraphen der Arbeit (§ 20) noch einmal Nachdruck. Bezugnehmend auf seine Untersuchungen zu Gleichung (\ref{eq: Zeta-Euler}) schreibt er:\\

\textit{``In gleicher Manier kann man die Summen dieser zwei Reihen}

\begin{equation*}
    1-3^{n-1}+5^{n-1}-7^{n-1}+\text{etc.} \quad \textit{und} \quad 1-\dfrac{1}{3^n}+\dfrac{1}{5^n}-\dfrac{1}{7^n}-\dfrac{1}{9^n}+\text{etc.}
\end{equation*}

\textit{miteinander vergleichen und eine ähnliche Vermutung liefert dieses Theorem}\footnote{Das Wort \textit{théoreme} ist in Eulers Arbeit kursiv gedruckt, womit Euler vermutlich die Gewissheit über die Richtigkeit der Formel unterstreichen will.}

\begin{equation*}
    \dfrac{1-3^{n-1}+5^{n-1}-7^{n-1}+\text{etc.}}{1-3^{-n}+5^{-n}-7^{-n}+\text{etc.}}= \dfrac{1 \cdot 2 \cdot 3 \cdots (n-1) \cdot 2^n}{\pi^n} \cdot \sin \left(\dfrac{n\pi}{2}\right)."
\end{equation*}

Kontrastiert man dies mit seinen Untersuchungen in der Arbeit \textit{``Consideratio quarumdam serierum, quae singularibus proprietatibus sunt praeditae"} (\cite{E190}, 1753, ges. 1743) (E190: ``Betrachtung gewisser Reihen, die einzigartige Eigenschaften aufweisen") zu der Reihe\footnote{Diese bildet auch den Untersuchungsgestand von Sandifers Arbeit \textit{``A false logarithm series"} (\cite{Sa07dec}, 2007).}:

\begin{equation*}
    s(a,x):= \dfrac{1-x}{1-a}+\dfrac{(1-x)(a-x)}{a-a^3}+\dfrac{(1-x)(a-x)(a^2-x)}{a^3-a^6} + \textit{etc.}
\end{equation*}
mit allgemeinen Term 

\begin{equation*}
    \dfrac{(1-x)(a-x)(a^2-x)(a^3-x)\cdots (a^{n-1}-x)}{a^{\frac{n(n-1)}{2}}-a^{\frac{n(n+1)}{2}}},
\end{equation*}
wird die Paradoxie der geschilderten Euler'schen Beweisaufassungen noch eindrücklicher.
Wie Euler in letztgenannter Abhandlung richtig bemerkt, hat man für die Wahl $x=a^n$ allgemein $s(a, a^n)=n$, eine Eigenschaft, die der Funktion $\log_a (x)$ zukommt. Jedoch gilt, wie Euler durch ein numerisches Beispiel nachweist, für nicht natürliche Zahlen $n$ bereits nicht mehr $s(a,a^n) = \log_a(x)$, sodass man die Gültigkeit eines für die natürlichen Zahlen gültigen Ausdrucks nicht ohne Hinzutreten eines weiteren Arguments auf die reellen Zahlen überführen kann\footnote{Der Nicht--Eindeutigkeit des Interpolationsproblem ist sich Euler selbstredend bewusst gewesen. Dies wird in Abschnitt (\ref{subsubsec: Methodus Inveniendi über Methodus Demonstrandi}) eine ausführlichere Verwendung finden.}.\\

Der Vollständigkeit wegen sei in diesem Zusammenhang auch die Euler'sche Arbeit \textit{``Exercitationes analyticae"} (\cite{E432}, 1773, ges. 1772) (E432: ``Analytische Übungen") erwähnt; sie hat ebenfalls obige Funktionalgleichung (\ref{eq: Zeta-Euler}) zum Inhalt. Gleich in § 2 findet man die Formel

\begin{small}
\begin{equation*}
    1-{2^{n-1}}+3^{n-1}-4^{n-1}+\cdots = \dfrac{2 \cdot 1 \cdot 2 \cdot 3 \cdots (n-1)}{\pi^n} \cdot N \cdot \left(1+\frac{1}{3^n}+\dfrac{1}{5^n}+\dfrac{1}{7^n}+\cdots\right).
\end{equation*}
\end{small}Obgleich  Euler bezüglich ihrer Allgemeingültigkeit seine Worte behutsamer wählt als in \cite{E352}, scheint dieselbe Meinung  in § 3  hervor; so schreibt er:\\

\textit{``Aber es lässt sich nicht bezweifeln, dass die einfachste und natürlichste [Formel] hier Geltung hat, das heißt $N=\cos \frac{n-2}{2}\pi$, während $\pi$ hier den zwei rechten gleichen Winkel bezeichnet, weil ja der ganze Sinus $=1$ angenommen wird, dass $\pi$ der Halbumfang des Kreises ist."}\\

Seine Aussage bezieht sich  auf den Term $\cos \frac{(n-2)\pi}{2}$ in der Formel, welchen er als Interpolation für die Folge von Werten $-1,0,+1,0,-1,0,+1 \cdots$ annimmt. Euler sieht diesen Term wegen der Einfachheit und Schönheit als richtige Wahl an, welches Prinzip  insbesondere von Paul Dirac (1902--1984)  in der theoretischen Physik vertreten wurde; man konsultiere diesbezüglich seine Arbeit \textit{``The Relation between Mathematics and Physics"} (\cite{Di39}, 1939).\\

Als allgemeines Resümee bleibt  die Differenz bezüglich der Stellung eines Beweises in der heutigen Mathematik und der von Euler.  Bezogen auf das vorgestellte Beispiel führt Euler gleichsam das \textit{Wunderargument zur Existenz unbeobachtbarer Realität}\footnote{Dieses Argument wird oft Hilary Putnam (1926--2016) in Bezug auf den wissenschaftlichen Realismus zugeschrieben. Er formuliert in seiner Arbeit \textit{``What is mathematical truth?"} (\cite{Pu75}, 1975) :
\textit{``The positive argument for realism is that it is the only philosophy that doesn't make the success of science a miracle."}} als Grund für die Wahrheit der Formel  (\ref{eq: Zeta-Euler}) ins Feld. Daraus ließe weiter folgern, dass Euler als Vertreter der Meinung, Mathematik werde \textit{entdeckt} und nicht \textit{kreiert}\footnote{Das sind auch die von Klein ausgemachten Archetypen: Er schreibt auf Seite 72 seines Buches \cite{Kl56}: \textit{``Die eine Gruppe von Mathematikern hält sich für unbeschränkte Selbstherrscher in ihrem Gebiet, das sie nach eigener Willkür logisch deduzierend aus sich heraus schaffen; die andere geht von der Auffassung aus, daß die Wissenschaft in ideeller Vollendung vorexistiere, und daß es uns nur gegeben ist, in glücklichen Augenblicken ein begrenztes Neuland davon zu entdecken."}.}, anzusehen ist\footnote{Eine gegenteilige Position wird von Vertretern der Kohärenztheorie der Wahrheit wie etwa Rescher vertreten. Letzterer erläutert seine Auffassungen in seinem Buch \textit{``The Coherence Theory of Truth"} (\cite{Re73}, 1973). Als einen Vertreter dieser Meinung in den mathematischen Wissenschaften ließe sich hier Poincaré (1854--1912) nennen. Man konsultiere dazu seine Essays zur Mathematikphilosophie, welche sich als Nachdruck im Buch \textit{``La science selon Henri Poincaré: La science et l'hypothèse - La valeur de la science - Science et méthode"} (\cite{Po24}, 2024) (``Die Wissenschaft nach Henri Poincaré: Wissenschaft und Hypothese -- Der Wert der Wissenschaft -- Wissenschaft und Methode") als Nachdruck finden. Man betrachte seine Beschreibung aus \textit{``Wissenschaft und Hypothese"} (Kapitel 2): \textit{``Mathematiker untersuchen nicht die Objekte, sondern die Beziehungen zwischen den Objekten; es ist daher für sie nebensächlich, die einen Objekte durch andere zu ersetzen, solange sich nur die Beziehungen nicht ändern."}}\\

Eng verknüpft mit dem just Gesagten ist, dass Euler des Öfteren ein und dieselbe Tatsache auf mehreren und unabhängigen Wegen beweist, wofür seine Beiträge zur Auswertung der Summe

\begin{equation*}
    1+\dfrac{1}{4}+\dfrac{1}{9}+\dfrac{1}{16}+\dfrac{1}{25}+\text{etc.}
\end{equation*}
ein konkretes Beispiel liefern (siehe Abschnitt \ref{subsec: Die Lösung des Baseler Problems} für eine eingehendere Diskussion). Dieses Vorgehen ist dabei für Euler jedoch nicht redundant, sondern gar nötig, um sich der Wahrheit einer Sache sicherer zu werden. Selbige Schlussweise führt Euler zumeist zu richtigen Resultaten; in der Tat haben die von Euler falsch vermuteten Sätze eine geringe Zahl an Argumenten für ihre Richtigkeit gemein\footnote{Neben der gleich zu diskutierenden Potenzvermutung hat sich die Euler'sche Vermutung zu den sogenannten griechisch--lateinischen Quadraten als falsch herausgestellt, welche Euler in seinen Arbeiten \textit{``Recherches sur un nouvelle  espèce de quarrés"} (\cite{E530}, 1782) (E530: ``Untersuchungen über eine neue Gattung von magischen Quadraten") und \textit{``De quadratis magicis"} (\cite{E795}, 1849, ges. 1776) (E795: ``Über magische Quadrate") untersucht. Die Euler'schen Ausführungen sowie die Unrichtigkeit seiner Vermutung sind umfassend in der Arbeit \textit{``Graeco-Latin Squares and a Mistaken Conjecture of Euler"} (\cite{Kl06}, 2006) dargestellt.}.

\paragraph{Beispiel einer unrichtigen Vermutung Eulers}

Als Beispiel für eine unrichtig von Euler vermutete Behauptung sei die Euler'sche Potenzenvermutung angeführt. Euler spricht sie in seinen Arbeiten \textit{``Resolutio formulae Diophanteae $ab(maa+nbb) = cd(mcc+ndd)$ per numeros rationales"} (\cite{E716}, 1802, ges. 1778) (E716: ``Auflösung der Diophant'schen Formel $ab(maa+bb)=cd(mcc+ndd)$ mit rationalen Zahlen") und \textit{``Dilucidationes circa binas summas duorum biquadratorum inter se aequales"} (\cite{E776}, 1830, ges. 1780) (E776: ``Erläuterungen zu den Summen zweier einander gleicher Biquadrate") aus. In § 3 erstgenannter Arbeit schreibt er:\\

\textit{``Wie nämlich keine zwei Kuben dargeboten werden können, deren Summe oder Differenz ein Kubus ist}\footnote{Den Beweis, hierfür hat Euler in § 243 des zweiten seines Buches \textit{``Vollständige Anleitung zur Algebra"}, bestehend aus den zwei Teilen  (\cite{E387}, 1770, ges. 1767) und (\cite{E388}, 1770, ges. 1776) gegeben, was auch unten (Abschnitt \ref{subsubsec: Wegen übersehener Unvollständigkeit: Der große Satz von Fermat für n=3}) noch diskutiert werden wird. Die Lösbarkeit der Gleichung $x^3+y^3+z^3=u^3$ beweist Euler in § 16 seiner Arbeit \textit{``Solutio generalis quorundam problematum Diophanteorum, quae vulgo nonnisi solutiones speciales admittere videntur"} (\cite{E255}, 1761, ges. 1754) (E255: ``Die allgemeine Lösung gewisse Diophant'scher Probleme, welche für gewöhnlich nur spezielle Lösungen zuzulassen scheinen").},\textit{ so ist auch gewiss, dass nicht einmal drei Biquadrate dargeboten werden können, deren Summe gleichermaßen ein Biquadrat ist, sondern dass mindestens vier Biquadrate erforderlich sind, damit deren Summe ein Biquadrat ergibt, obgleich noch niemand vier solche Biquadrate angeben konnte. In gleicher Weise scheint bekräftigt werden zu können, dass keine vier fünften Potenzen dargeboten werden können, deren Summe auch eine fünfte Potenz ist; ebenso wird sich die Angelegenheit bei den höheren Potenzen verhalten; [...]"}\\

In § 1 der zweiten Arbeit formuliert Euler seinen Verdacht etwas behutsamer; er schreibt:\\

\textit{``Weil solche Formeln $A^2\pm B^2=0$, $A^3\pm B^3 \pm C^3=0$ als unmöglich nachgewiesen worden sind, wenn freilich gleiche und verschwindende Zahlen ausgeschlossen werden, kann es den Anschein haben, dass auch diese Form $A^4 \pm B^4 \pm C^4 \pm D^4 =0$, und entsprechende für die höheren Potenzen $A^5\pm B^5 \pm C^5 \pm D^5 \pm E^5=0$ und $A^6 \pm B^6 \pm C^6 \pm D^6 \pm E^6 \pm F^6=0$ unmöglich sind."}\\

Modern ließe sich die Euler'sche Vermutung wie folgt fassen: Ist die Gleichung

\begin{equation}
\label{eq: Euler Power Conj}
    \sum_{k=1}^{p}x_k^n =0
\end{equation}
über den rationalen Zahlen lösbar,  gilt $p >n$. In dieser Form ist die Euler'sche Vermutung als unrichtig nachgewiesen worden. Der Fall $n=4$ wurde von  Elkies (1966--)  abgehandelt, welcher in seiner Arbeit \textit{``On $A^4+B^4+C^4=D^4$"} (\cite{El88}, 1988) das Beispiel 

\begin{equation}
\label{eq: Elkies}
    20615673^4=2682440^4+15365639^4+18796760^4
\end{equation}
angibt, nachdem bereits zuvor von Landir und Parkin in der Arbeit \textit{``Counterexample to Euler’s conjecture on sums of like powers"} (\cite{La66}, 1966) für den Fall $n=5$ mit der Lösung 

\begin{equation*}
    144^5=27^5+84^5+110^5+133^5
\end{equation*}
ein Gegenbeispiel angegeben hatten. Ob tatsächlich stets $n=p$ auch das Vorhandensein einer rationalen Lösung von (\ref{eq: Euler Power Conj}) bedeutet, ist eine offene Frage. \\

In diesem Zusammenhang könnte als weitere indirekte Bestätigung des Gesagten die Widerlegung der Fermat'schen Vermutung, bei allen Zahlen der Form

\begin{equation*}
    F_n = 2^{2^n}+1
\end{equation*}
handele es sich um Primzahlen, durch Euler, angeführt werden. Während die ersten fünf Zahlen:

\begin{equation*}
   F_0=3,~~ F_1=5,~~ F_2=17,~~ F_3=257,~~ F_4=65537 
\end{equation*}
als Primzahlen Fermat (1601--1665) zu seiner Vermutung veranlasst haben, widerlegt Euler in seiner Arbeit \textit{``Observationes de theoremate quodam Fermatiano aliisque ad numeros primos spectantibus"} (\cite{E26}, 1738, ges. 1732) (E26: ``Beobachtungen zu einem gewissen Fermat'schen Lehrsatz und anderen, welche sich auf die Primzahlen beziehen") die Fermat'sche Behauptung, indem er $641$ als Teiler von $2^{2^5}+1$ angibt\footnote{Für Eulers genaues Vorgehen konsultiere man zum Beispiel Sandifer's Artikel \textit{``Factoring $F_5$"} (\cite{Sa07mar}, 2007).}. Trotz ähnlicher Sachlage bezüglich der Evidenz wie bei der Potenzvermutung (Richtigkeit für die Zahlen $1$, $2$ und $3$), lässt sich Euler bei den Fermat'schen Primzahlen (Richtigkeit für die ersten fünf Zahlen) nicht von der Richtigkeit der Vermutung überzeugen. Dies liegt natürlich an der expliziten Angabe eines Gegenbeispiels, wohingegen zu Eulers Lebzeiten das Auffinden des Beispiels von Elkies aus(\ref{eq: Elkies}) wegen fehlender Methoden und der Größe der Zahlen gleichsam eine Unmöglichkeit darstellten.  \\ 

Abschließend sei bezüglich der Euler'schen Herangehensweise an Beweise  die Bewertung von G. Polya (1887--1985) aus seinem Buch \textit{``Mathematics and Plausible Reasoning"} (\cite{Po14}, 2014)\footnote{Die ursprüngliche Version des Buchs wurde von Polya 1954 verfasst.} mitgeteilt:\\

\textit{``Euler's reasons are not demonstrative. Euler does not reexamine the grounds for his conjecture [...] only its consequences. [...] He examines also the consequences of closely related analogous conjectures [...] Euler's reasons are, in fact, inductive."}

\subsubsection{Beweistechnik der Induktion bei Euler}
\label{subsubsec: Beweistechnik der Induktion bei Euler}

\epigraph{It is by logic that we prove, but by intuition that we discover. To know how to criticize is good, to know how to create is better.}{Henri Poincaré}

Das Polya'sche Zitat schafft einen unmittelbaren Übergang zum nächsten zu besprechenden Thema: Die Rolle der Induktion bei Euler.

\paragraph{Eulers Auffassung vom Prinzip der Induktion}
\label{para: Eulers Auffassung vom Prinzip der Induktion}

Euler beschreibt seine Ansicht  zur Anwendung des Induktionsprinzips  unumwunden in Arbeiten zur Zahlentheorie. In pädagogischer Weise mit expliziten Hinweisen für die Nützlichkeit der Methode beim Forschungsprozess tut er dies etwa in der Arbeit \textit{``Specimen de usu observationum in mathesi pura"} (\cite{E256}, 1761, ges. 1754) (E256: ``Ein Beispiel über den Nutzen von Beobachtungen in der reinen Mathematik")\footnote{Sandifer hat die Arbeit \textit{``$2aa+bb$"} (\cite{Sa06jan}, 2006) seiner Kolumne den Euler'schen Ausführungen diesbezüglich gewidmet.}. Euler schreibt im Vorwort zu besagter Abhandlung:\\

\textit{``Aus diesen Erläuterungen schließen wir mit Recht, dass in der Erforschung der Natur der Zahlen der Beobachtung und der Induktion [...] sehr große Wichtigkeit einzuräumen ist. [...] Denn mit dieser Methode sind wir zur Kenntnis von Eigenschaften dieser Art gelangt, welche uns andernfalls verborgen geblieben wären. [...] Obwohl aber Beschaffenheiten von dieser Art durch scharfsinnige Beobachtungen entdeckt worden sind, [...] können wir dennoch, wenn nicht ein strenger Beweis hinzutritt, bezüglich ihrer Richtigkeit nicht hinreichend sicher sein."}\\

 Die Euler'schen Ausführungen zum Gegenstand der Induktion im Teilbereich der Zahlentheorie sind aus moderner Sicht umso interessanter, weil das Fehlen der für einen Fortschritt notwendigen Konzepte Euler zwingt seine Beobachtungen anders zu untermauern\footnote{So entgeht Euler etwa selbst das quadratische Reziprozitätsgesetz nicht. Es findet sich bereits in der Arbeit \textit{``Theoremata circa divisores numerorum in hac forma $paa \pm qbb$ contentorum"} (\cite{E164}, 1751, ges. 1747) (E164: ``Theoreme über die in der Form $paa \pm qbb$ enthaltenen Zahlen"), wie zuerst von Kronecker (1823--1891) in seiner Arbeit \textit{``Bemerkungen zur Geschichte der Reciprocitätsgesetzes"}  (\cite{Kr76}, 1876) bemerkt worden ist. Besagte Euler'sche Arbeit enthält, wie Euler auch selbst zugibt, keinen einzigen Beweis, aber nahezu 60 Theoreme. Ein Nachweis des Reziprozitätsgesetzes ist Euler verwehrt geblieben und wurde erst von Gauß in seinem Buch \textit{``Disquisitiones arithemeticae"} (\cite{Ga01}, 1801, ges. 1798) (``Untersuchungen zur Zahlentheorie") geliefert.}.

\paragraph{Ein Beispiel falscher Induktion bei Euler}
\label{para: Ein Beispiel falscher Induktion bei Euler}

Da Euler Induktion  als den Schluss von Einzelfällen auf ein allgemeines Ergebnis versteht, ist ihm die Fehleranfälligkeit dieser Methode nicht entgangen. Dies hat Euler klar in §§ 1--3 der Veröffentlichung \textit{``Varia artificia in serierum indolem inquirendi"} (\cite{E551}, 1783, ges. 1772) (E551: ``Verschiedene Kunstgriffe, die Natur von Reihen ausfindig zu machen") erläutert. Hier betrachtet er die Reihe der größten Koeffizienten in der Entwicklung von $(1+x+x^2)^n$ nach Potenzen von $x$. Die ersten Terme dieser Reihe sind:

\begin{equation*}
    1,~~3,~~7,~~19,~~51,~~141,~~393 ~~\text{etc.}
\end{equation*}
Die  wesentlichen Euler'schen Gedanken bezüglich der Induktion werden im Folgenden kompakt nachgezeichnet. In § 1 schreibt Euler:\\

\textit{``Demjenigen der diese Reihe betrachtet, gelangt schnell in den Sinn, dass jeder beliebige Term trefflich mit dem Dreifachen des vorhergehenden verglichen werden kann, weil es aus ihrem Ursprung ersichtlich ist, dass diese Reihe ins Unendliche fortgesetzt mit der verdreifachten geometrischen Progression zusammenfallen muss."}\\

Euler versucht nun durch weitere Operationen die allgemeinen Terme der Reihe solange umzuformen, bis er auf eine Progression stößt, für welche er einen expliziten Term anzugeben imstande ist. Aus seinen Überlegungen konstruiert er folgende Tabelle:

\begin{equation*}
    \renewcommand{\arraystretch}{2,0}
\begin{array}{l|ccccccccccc}
   \text{Index} & 0  & 1 & 2 & 3 & 4 & 5 & 6 & 7 & 8 & 9  \\  
    A(n) & 1 & 1 & 3 &7 & 19 & 51 & 141 & 393 & 1107 & 3139 \\
    B(n+1)(=3A(n)) &   & 3 & 3 & 9 & 21 & 57 & 153 & 423 & 1179 & 3321 \\ \hline
    C(n)(=B(n)-A(n)) & & 2 & 0 & 2 &  2 & 6 & 12 & 30 & 72 & 182  \\
    D(n)(=\frac{1}{2}C(n)) & & 1 & 0 & 1 & 1 & 3 & 6 & 15 & 36 &  91 \\
    E(n)(=D^{-1}(n))  &  & 1  & 0 & 1&  1 & 2 & 3 & 5 & 8 & 13 
\end{array}
\end{equation*}

In seiner Notation unterdrückt Euler die Abhängigkeit von $n$ und auch die in Klammern geschriebenen Ausdrücke. Wie Euler jedoch erklärt, ergibt sich $B(n)$ aus Verdreifachung von $A(n)$. $C(n)$ ist gerade die Differenz von $B(n)$ und $A(n)$. $D(n)$ ist die Hälfte von $C(n)$. Entscheidend ist nun, dass Euler die Folge $D(n)$ als Dreieckszahlen erkennt, woraus er die Folge $E(n)$ als deren Indizes\footnote{Euler selbst schreibt das lateinische Wort für Wurzel.} definiert. Die Folge $E(n)$ sind aber gerade die Fibonacchi--Zahlen, für welche Euler mithilfe der Theorie der rekurrenten Reihen eine explizite Formel anzugeben weiß\footnote{In seiner Arbeit \cite{E551} gibt er unmittelbar die finale Formel für $A(n)$ an. Seine Methode, die $n$--te Fibonacchi--Zahl zu finden, erläutert er unter anderem in seiner \textit{Introductio}. Man betrachte insbesondere seine allgemeinen Ausführungen in §§ 212--214. Aus seinem Beispiel 4 von § 216 ergibt sich die explizite Formel für die Fibonacchi--Zahlen als Spezialfall.}. Daraus kann man durch Umkehren aller erläuterten Schritte aus der expliziten Formel für $E(n)$ die für $A(n)$ finden.  Über die Reihe $E(n)$, welche Euler schlicht $E$ nennt, schreibt er in § 2:\\

\textit{
    ``In dieser Reihe $E$ scheint die Struktur der Terme so beschaffen zu sein, dass jeder der Summe der zwei vorhergehenden gleich ist, und diese Schlussfolgerung, allein durch Beobachtung gebildet, scheint dermaßen gewiss, dass [...]."}\\
    
Daraus schließt Euler dann in erläuterter Weise auf die ursprüngliche Progression, welcher er mit $A$ bezeichnet, zurück und fasst in § 3 zusammen:\\

\textit{
    ``Wenn diese Induktion mit der Wahrheit einherginge, wäre sie als eine Entdeckung von größter Bedeutung anzusehen, weil dann der allgemeine Term der vorgelegten Reihe $A$ angegeben werden könnte."}\\
    
Diesen gibt Euler auch explizit an, aber sagt dann in § 4:\\

\textit{``[...], weil ja die obige Induktion, so wahr sie auch erschien, dennoch von der Wahrheit abkommt."}\\

Die Abweichung stellt man beim Vergleich des $11.$ Terms der durch Induktion gefundenen Formel mit dem tatsächlichen Wert fest, womit das Euler'sche Beispiel ``falscher Induktion"{} treffend gewählt ist\footnote{Sandifer hat über dieses  Beispiel von Euler den Artikel \textit{``A memorable example of false induction"} (\cite{Sa05aug}, 2005) verfasst.}.

\paragraph{Ein Beispiel vollständiger Induktion bei Euler: Der kleine Satz von Fermat}
\label{para: Kleiner Satz von Fermat}

Jedoch findet man auch Beispiele von Induktion in heutigem Sinne bei Euler, ohne dass Euler selbst dies auch als Induktion bezeichnet hat. Dazu betrachte man den nachstehenden bekannten Lehrsatz. \\

\begin{figure}
    \centering
   \includegraphics[scale=0.8]{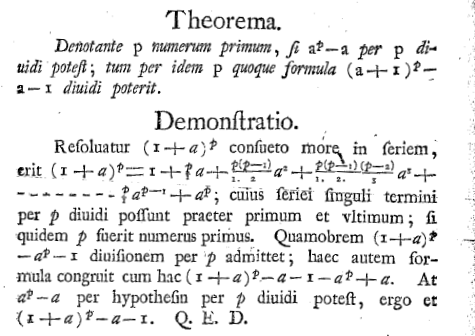}
    \caption{Eulers Beginn des Beweises kleinen Satzes von Fermat mittels vollständiger Induktion nach  natürlichen Zahl $a$ aus seiner Arbeit \cite{E54}.}
    \label{E56Fermat}
\end{figure}

\begin{Thm}[Kleiner Satz von Fermat]
\label{Theorem: Kleiner Satz von Fermat}
Für eine Primzahl $p$ und eine zu $p$ teilerfremde Zahl $a$ gilt: 
\begin{equation*}
    a^{p-1} \equiv 1 \quad \operatorname{mod} p
\end{equation*}
\end{Thm}
 Obschon Euler in der später veröffentlichten Arbeit \textit{``Theoremata arithmetica nova methodo demonstrata"} (\cite{E271}, 1763, ges. 1758) (E271: ``Zahlentheoretische Theoreme -- bewiesen mit einer neuen Methode") einen allgemeinen Beweis, der sehr an moderne Versionen erinnert\footnote{In besagter Arbeit führt Euler die Euler'sche $\varphi$-Funktion ein, welche die Anzahl der zur natürlichen Zahl $n$  teilerfremden Zahlen, welche $<n$ sind, angibt, und leitet den kleinen Satz von Fermat als Spezialfall seiner allgemeineren Untersuchungen ab. Man vergleiche insbesondere seine Ausführungen aus §61. In besagter Arbeit benutzt Euler allerdings noch kein eigenes Symbol für die $\varphi$--Funktion, sondern erläutert ihren Inhalt verbal. In seiner später verfassten Abhandlung \textit{``Speculationes circa quasdam insignes proprietates numerorum"} (\cite{E564}, 1784, ges. 1775) (E564: ``Betrachtungen zu gewissen außerordentlichen Eigenschaften von Zahlen") nutzt er den Formelbuchstaben $\pi$.}, gegeben hat, war sein erster Beweis des obigen Lehrsatzes ein sehr untypischer. In seiner Arbeit  \textit{``Theorematum quorundam ad numeros primos spectantium demonstratio"} (\cite{E54}, 1741, ges. 1736) (E54: ``Beweis gewisser sich auf Primzahlen beziehender Theoreme") beweist er selbigen  ab § 3  mit Induktion -- wie man sie heute in der Mathematik versteht -- nach $a$.\\

 Trotz der Vorteile der mathematischen Induktion (oder auch vollständigen Induktion) als ein generisches Beweisverfahren, muss sie für Euler unbefriedigend gewesen, weil sie die Kenntnis des nachzuweisenden Ergebnisses voraussetzt. Seinem Arbeitsethos entsprechend hat Euler im Verlauf seiner Karriere noch  weitere Beweise des kleinen Satzes von Fermat geliefert, welche alle in der Arbeit \textit{``Euler and Number Theory: A Study in Mathematical Invention"}  (\cite{Su07}, 2007) detailliert besprochen werden. Einen Teil dieser Nachweise bespricht auch Sandifer in seinem Artikel \textit{``Fermat’s Little Theorem"} (\cite{Sa03nov}, 2003).

\subsubsection{Eulers Ansicht zur Physik}
\label{subsubsec: Eulers Ansicht zur Physik} 

\epigraph{The miracle of the appropriateness of the language of mathematics for the formulation of the laws of physics is a wonderful gift which we neither understand nor deserve.}{Eugene Wigner}


Zu Eulers Zeiten herrschte zwischen Mathematik und Physik noch nicht eine solch scharfe Trennung wie es heute der Fall ist, weswegen  die Euler'sche Ansicht zum Wesen der Physik bündig referiert wird. In seinem Essay \cite{E790} findet man diesbezüglich die Passage\footnote{Die Übersetzung ist wie oben aus der von Burckhardt aus Serie 3, Band 2 der \textit{Opera Omnia} (\cite{OO32}, 1942) entnommen.}:\\

\textit{``Mag auch die Physik, welche die Ursachen aller Vorgänge in der Natur untersucht, eines offensichtlichen Nutzens entbehren, so muss doch jeder, der die Wahrheit liebt, die Würde und Erhabenheit des Ziels anerkennen. [...] Denn alle Veränderungen, die wir in der Natur beobachten, stammen von Bewegungen her; es ist also klar, dass die Mechanik, das heißt die Wissenschaft von der Bewegung notwendig ist, um selbst die kleinste Bewegung im Universum zu erklären."} \\

Hier kristallisiert sich  vor allem die Ansicht heraus, die sei Physik für die Wahrheitsfindung notwendig. Insbesondere bezogen auf die Bewegungslehre ist sie womöglich überdies hinreichend. Die Einschränkung geschieht hier, weil es unter Umständen noch andere Erklärungen für einen Sachverhalt geben kann. Moderne Terminologie gebrauchend, ist es wohl zutreffend, dass die Mechanik eine Abduktion, wie Peirce (1839--1914) den Begriff in seinen \textit{``Vorlesungen über Pragmatismus"} (\cite{Pe03}, 1903) erläutert\footnote{Peirce formuliert diesen Begriff wie folgt: \textit{``Abduktion ist der Vorgang, in dem eine erklärende Hypothese gebildet wird."} (CP 5.171).}, für die Bewegung ist. \\

Euler äußert insbesondere implizit die Ansicht,  aus den Axiomen der Mechanik abgeleitete Sätze manifestieren sich zwangsläufig in der Natur, worin man -- sich Kant'scher Terminologie bedienend -- eine physiko--teleologische Weltsicht erkennen kann\footnote{Kant (1724--1804) hat diesen Begriff in seinem Werk \textit{``Kritik der reinen Vernunft"} (\cite{Ka81}, 2009, urspr. 1781) eingeführt und meinte damit ein Weltbild, welches jede Wirkung in der Natur als einem höheren Ziel zweckdienlich sieht.}.  Er schreibt dazu hingegen explizit in seiner Arbeit \textit{``Meditationes super problemate nautico, quod illustrissima regia Parisiensis Academia Scientiarum proposuit"} (\cite{E4}, 1738, ges. 1729) (E4: ``Betrachtungen über das Problem aus der Nautik, welches die illustre Pariser Akademie der Wissenschaften gestellt hat") in § 50:\\

\textit{``Ich habe es nicht als nötig erachtet, diese meine Theorie durch das Experiment zu bestätigen, weil sie gänzlich aus sichersten und unanfechtbarsten Prinzipien abgeleitet worden ist, und daher kein Zweifel  aufkommen kann, ob sie wahr ist und auch in der Praxis ihre Gültigkeit behalten kann."}\\

Folgender Ausspruch aus seiner Arbeit \textit{``Réflexions sur l'espace et le tems"} (\cite{E149}, 1750, ges. 1748) (E149: ``Gedanken zur Natur von Raum und Zeit") unterstreicht dies; dort schreibt er unmittelbar zu Beginn:\\

\textit{``Die Prinzipien der Mechanik sind schon dermaßen etabliert, dass es höchst töricht wäre, wenn man immer noch an ihrer Wahrheit zweifeln wollte. Obgleich wir nicht in der Lage sind, sie durch allgemeine Prinzipien der Metaphysik zu demonstrieren, sollte die wunderbare Übereinstimmung aller Folgerungen, welche man aus ihnen vermöge des Kalküls zieht, mit den Bewegungen der Festkörper sowie der Fluide und gar der Himmelskörper, hinreichen, ihre Gültigkeit zweifelsfrei zu belegen."}\footnote{Seine Auffassungen zu den Grundlagen der Mechanik hat Euler in seinem Opus \textit{``Theoria motus corporum solidorum seu rigidorum"} (\cite{E289}, 1765, ges. 1760) (E289: ``Die Theorie der Bewegung von soliden und starren Körper") gebündelt aufbereitet. Die Hydrostatik bespricht Euler in der Abhandlung \textit{``Sectio prima de statu aequilibrii fluidorum"} (\cite{E375}, 1769, ges. 1766) (E375: ``Erster Abschnitt über den Gleichgewichtszustand von Fluiden"), seine Ideen zur Hydrodynamik sind in der Arbeit \textit{``Sectio secunda de principiis motus fluidorum"} (\cite{E396}, 1770, ges. 1766) (E396: ``Zweiter Abschnitt über die Prinzipien der Bewegung von Fluiden") zusammengefasst zu finden. Für Eulers Erkenntnisse zur Himmelsmechanik ist auf das Buch \textit{``Leonhard Eulers Arbeiten zur Himmelsmechanik"} (\cite{Ve14}, 2014) verwiesen.}\\

In § V. schreibt er weiter zur Natur des Raums und der Zeit in diesem Zusammenhang: \\

\textit{``Man muss vielmehr schließen, dass so die absolute Zeit wie der absolute Raum, wie Mathematiker sie darstellen, reelle Dinge sind, welche über die Vorstellung hinaus ebenfalls Bestand haben."}\\

Heute wird die bis mindestens zu Leibniz (1646--1716) und Newton (1642--1726) zurückreichende Diskussion, ob Zeit und Raum relativ oder absolut sind\footnote{Newton vertritt in seinen \textit{``Philosophiae Naturalis Principia Mathematica"} (\cite{Ne87}, 1687) (``Mathematische Prinzipien der Naturphilosophie") die These, dass Raum und Zeit absolut sind.}, als unentschieden betrachtet\footnote{Die ursprüngliche Debatte zwischen Leibniz und Samuel Clarke (1675--1729), Newtons Fürsprecher, besteht aus 5 Abhandlungen von Leibniz und 5 jeweiligen Repliken von Clarke. Man findet sie im Buch \textit{``The--Leibniz Correspondence"} (\cite{Al98}, 1998) abgedruckt. Die wesentlichen Argumente werden beispielsweise in dem Büchlein \textit{``Philosophy of Science -- A Very Short Introduction"} (\cite{Ok16}, 2016) vorgestellt und diskutiert.}, wohingegen Euler eindeutig Stellung bezieht. Somit ist folgende sehr apodiktisch anmutende Aussage aus § XV, aus welcher sich auch ein deterministisches Weltbild herauslesen lässt, fast schon folgerichtig:\\

\textit{``Ich gehe damit konform, dass alle Dinge, welche existieren, vollständig bestimmt sind, [...]"}\\

Geht man nämlich von ein deterministischen Weltbild aus, so kann man auch die Physik bzw. die Mechanik als wahrheitskonform ansehen. Die Umkehrung dieser Aussage ist hingegen unrichtig.\\

 Folgendes Zitat von Fellmann in seiner Euler--Biografie \cite{Fe95} (S. 24) fasst alles Vorgetragene zusammen:\\

\textit{``Dieses fast blinde Vertrauen in die Stringenz  der Prinzipien und in die apriorischer Deduktionen begleitete Euler bis in sein hohes Alter und kennzeichnet ein Paradigma seines Schaffens."}\\

Umgekehrt, so würde Euler begründen, schreibt auch die Physik vor, welche Objekte in der Mathematik überhaupt einer Betrachtung wert sind. An kaum einer Stelle wird diese wechselseitige Befruchtung deutlicher als am Beispiel der schwingenden Saite, welche Euler in \cite{E322} gleichsam dazu zwang, seinen Funktionsbegriff dahingehend abzuwandeln, auch das mit einzuschließen, was heute als Distribution bezeichnen wird\footnote{Man konsultiere diesbezüglich insbesondere den Übersichtsartikel \textit{``Euler’s Vision of a General Differential Calculus for a Generalized Kind of Function"} (\cite{Lu83}, 1983).}. \\

Ein Vergleich zu dem modernen Artikel \cite{Di39}, welcher auch die Relation von Mathematik und Physik zum Inhalt hat, bestätigt das Erläuterte auf noch einmal mehr. Diracs Position wird prägnant durch folgenden Satz zusammengefasst:\\

\textit{``One may describe the situation by saying that the mathematician plays a game in which he himself invents the rules while the physicist plays a game in which the rules are provided by Nature, but as times goes on it becomes increasingly evident that the rules which the mathematician finds interesting are the same as those which Nature has chosen."}\\

\subsection{Abriss von Eulers Vorgehensweise}
\label{subsec: Abriss von Eulers Vorgehensweise} 

\epigraph{Continuous improvement is better than delayed perfection.}{Mark Twain}


Die vorausgeschickten Abschnitte zur Euler'schen Auffassung der Essenz der Mathematik und Physik gestattet nun die Betrachtung der sich daraus ergebenen Konsequenzen für die konkrete Arbeitsweise Eulers. Es werden die Leitfäden seines Schaffens (Abschnitt \ref{subsubsec: Leitfäden der Arbeitsweise}) und Alleinstellungsmerkmale seines Präsentationsstils (Abschnitt \ref{subsubsec: Besonderheiten der Präsentationsweise}) anhand Diskussion ihrer Teilfacetten auseinander gesetzt.

\subsubsection{Leitfäden der Arbeitsweise}
\label{subsubsec: Leitfäden der Arbeitsweise}

\epigraph{Immer mit den einfachsten Beispielen anfangen.}{David Hilbert}

Aus den grundlegenden sich aus dem Essay \cite{E790} herauskristallisierenden Euler'schen Auffassungen zum Wesen der Mathematik manifestieren sich insbesondere in seiner Arbeits-- und seines Präsentationsstils grundlegende Unterschiede im Vergleich zum heutigen Leitbild. Während heutzutage erneut das Euklid'sche Vorbild des Dreischritts aus Definition, Satz und Beweis die Darstellung mathematischer Abhandlungen prägt, arbeitet sich Euler, nahezu diametral dazu, von Beispielen zur Allgemeinheit vor. Man vergleiche hierzu die Worte von Hankel (1839--1873) auf Seite 16 seines Buchs \textit{``Die Entwickelung der Mathematik in den letzten Jahrhunderten"} (\cite{Ha69}, 1869):\\

\textit{``Euler's Art zu arbeiten, war in der Hauptsache die, dass er zunächst seine Kräfte auf ein specielles Problem concentrierte, und so zu einer speciellen Auflösungsmethode gelangte. Daran schließt sich dann in einer folgenden Abhandlung häufig ein zweites, jenem verwandtes, ein drittes, viertes Problem, das er wiederum mit einer speciellen, jener ersten verwandten, aber dem neuen Probleme angepassten Form behandelt. Hierin ist er unübertrefflich; kein zweiter Mathematiker kommt ihm gleich an Fülle analytischer Gedanken und Geschick, die Methoden der speciellen Problem[SIC] zu accomodieren. [...]"}. \\

Dies kann man mit David Hilberts (1862--1943) Beschreibung aus seinem Buch \textit{``Geometry and the Imagination"} (\cite{Hi52}, 1952) ins Verhältnis setzen, wo er über den Forschungsprozess in der Mathematik schreibt:\\

\textit{``In mathematics, as in any scientific research, we find two tendencies [...] [T]he tendency toward abstraction seeks to crystallize the logical relations inherent in the maze of material in a systematic and orderly manner. On the other hand, the tendency toward intuitive understanding fosters a more immediate grasp of the objects [...] a live rapport with them [...] which stresses the concrete meaning of their relations."}\\

Hilbert lässt demnach Euler als Vertreter der zweiten Kategorie erkennen.
\paragraph{Methodus Inveniendi und praktische Fragestellungen als Leitfaden}
\label{para: Methodus inveniendi}

Angesichts des Hankel'schen Zitats bildet sich die richtige Erwartung heraus, dass eine Großzahl von Eulers Arbeiten von der Entwicklung von Methoden zur Lösung eines spezifischen Problems handeln wird\footnote{Eine hohe zweistellige Anzahl (ca. 70) seiner Arbeiten wie etwa sein wegweisendes Buch zur Variationsrechnung \textit{``Methodus inveniendi lineas curvas maximi minimive proprietate gaudentes, sive solutio problematis isoperimetrici lattissimo sensu accepti"} (\cite{E65}, 1744, ges. 1743) (E65: ``Eine Methode sich der Eigenschaft des Maximums oder Minimums erfreuende Kurven zu finden, oder die Lösung des im weitesten Sinne aufgefassten isoperimetrischen Problems") tragen schon im Titel die Worte \textit{``Methodus inveniendi"} oder ähnliche Formulierungen.}. Während Euler demnach eigens spezielle Methoden zur Auflösung eines gestellten Problems entwickelt, besteht der moderne Ansatz zumeist zur Problemlösung in der Verallgemeinerung der vorhandenen Begriffe zu einem höheren Grad an Abstraktion, mit welchem dann Theoreme entwickelt werden, welche das ursprüngliche Problem als Spezialfall der abstrakten Theorie begreifen. \\

\begin{figure}
    \centering
   \includegraphics[scale=0.7]{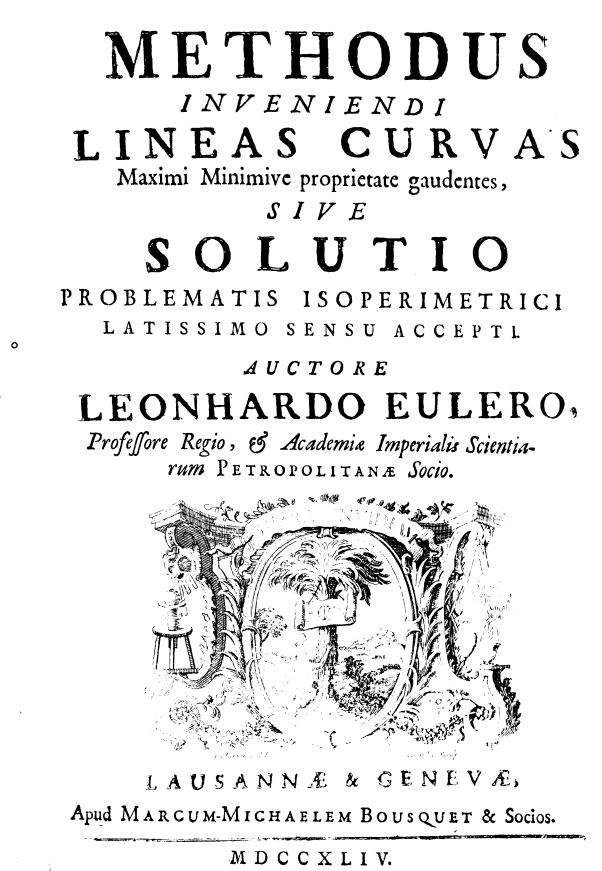}
    \caption{Das Titelblatt von Eulers Buch zur Varitionsrechnung \cite{E65} als Pars Pro Toto für die von ihm bevorzugte \textit{Methodus inveniendi}.}
    \label{fig:E65Titelblatt}
\end{figure}

 Während im Haupttext die Entdeckung der ``Numeri Idonei`` (siehe Abschnitt \ref{para: Quadratische Formen}) als Beispiel herangezogen werden wird, dass der Praxisbezug Euler bisweilen einen weiteren Fortschritt versperrte, soll hier ein Exempel gemacht werden,  dass Euler selbst abstraktere Fragen zumeist mit praktischen Anwendungen im Hinterkopf behandelt hat. Die Arbeit \textit{``Investigatio functionum ex data differentialium conditione"} (\cite{E285}, 1764, ges. 1759) (E285: ``Das Finden von Funktionen aus einer gegebenen Bedingung an die Differentiale") nimmt sich der Frage an, wie partielle Differentialgleichungen mit allgemeiner gegebener Struktur gelöst werden können. Die Motivation für diese Untersuchung gibt Euler in der Mitte von  § 5 dieser Abhandlung:\\

\textit{``Obwohl aber Fragen von dieser Art  fast gänzlich neu erscheinen, besteht indes kein Zweifel, dass die Methoden,  sie entsprechend aufzulösen für die gesamte Mathematik einen Nutzen haben wird. Denn beim ganzen Problem der schwingenden Saite ist die Essenz der Lösung zu dieser Gattung zu rechnen, weil sie in einer gewissen Beziehung zwischen den Funktionen $P$ und $Q$ verborgen liegt. Überdies habe ich ebenfalls diese gesamte Wissenschaft von der Bewegung von Fluiden in Differentialformen von dieser Art erfasst, wo eine bestimmte Relation zwischen den Anteilen der jeweiligen Differentiale vorgeschrieben wird, aus welcher sich aber wegen des Mangels einer solchen Lösungsmethode kaum etwas ableiten lässt."}\footnote{Euler referiert dabei vermutlich auf seine Arbeiten \textit{``De vibratione chordarum exercitatio"} (\cite{E119}, 1749, ges. 1748) (E119: ``Über die Vibration einer zum Schwingen angeregten Saite") bezüglich der schwingenden Saite und \textit{``Principia motus fluidorum"} (\cite{E258}, 1761, ges. 1753) (E258: ``Prinzipien der Bewegung von Fluiden") bezüglich der Hydrodynamik. In letztgenannter Arbeit leitet Euler die nach ihm benannten Euler'schen Gleichungen für die Dynamik inkompressibler Fluide her.}\\

Als weiteres Beispiel für Euler'sche Praxisnähe  ließen sich im Allgemeinen  seine Methoden nennen, die einer effizienteren numerischen Berechnung versprechen. Zu jedem Themenfeld, welchem Euler seine Aufmerksamkeit gewidmet hat, finden sich Arbeiten, die die explizite Berechnung zum Gegenstand haben. Allen voran geht dabei die Euler--Maclaurin'schen Summationsformel. Hergeleitet hat er sie zu ersten Mal vollständig in der Arbeit \textit{``Inventio summae cuiusque seriei ex dato termino generali"} (\cite{E47}, 1741, ges. 1735) (E47: ``Das Finden jedweder Summe aus einem gegebenen allgemeinen Term") aus der Einsicht, dass der Ausdruck

\begin{equation*}
    f(x):= \sum_{k=1}^{x-1} g(k)
\end{equation*}
 der Gleichung $f(x+1)-f(x)=g(x)$ Genüge leistet. In just erwähnter Arbeit wählt Euler einen geschickten Ansatz und gelangt so zur bekannten Summenformel. Voll ausgeschrieben findet man die Summenformel überdies  in  der Arbeit \textit{``De seriebus quibusdam considerationes"} (\cite{E130}, 1750, ges. 1739) (E130: ``Betrachtungen über gewisse Reihen") (§§ 26--27)\footnote{Dies wird auch von Sandifer in seinem Artikel \textit{``Bernoulli numbers"} (\cite{Sa05sep}, 2005) im Kontext der Bernoulli'schen Zahlen angemerkt.}.

\begin{figure}
    \centering
     \includegraphics[scale=0.8]{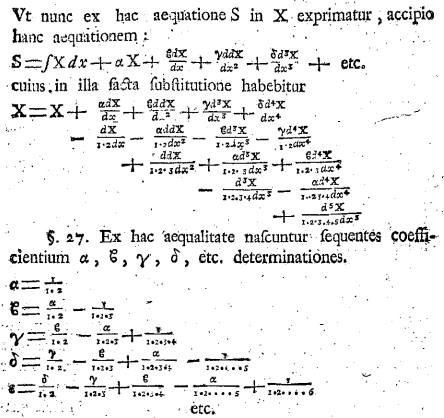}
    \caption{Euler gibt die Euler-Maclaurin'sche Summenformel in \cite{E130} an. Zu sehen ist ganz oben die allgemeine Formel für eine Summe $S$ bis hin zu $n$ Termen summiert über den allgemeinen Term $X$.  Am unteren Rand sind noch die ersten Entwicklungskoeffizienten $\alpha$, $\beta$, $\gamma$, $\delta$, $\varepsilon$ zu erkennen.}
    \label{fig:E130EulerMaclaurin}
\end{figure}

\paragraph{Dialektischer Ansatz}
\label{para: Dialektischer Ansatz}

In seinen \textit{``Briefen an eine deutsche Prinzessin"} \cite{E343} hat Euler auf den dialektischen Prozess bei epistemologischen Fragen hingewiesen und nimmt damit insbesondere die Ausführungen Hegels (1770--1831) vorweg, welcher den Fortschritt in dem Wort \textit{``Aufhebung"} subsummiert hat. Eulers Ausführungen sind insofern erhellender, zumal er unterstreicht, dass sich gewisse Paradoxa durch entsprechende Einführung von neuen Konzepten auflösen und die Einzelteile sich synthetisieren lassen. Ein  Beispiel dessen findet sich in Eulers Auflösung des Cramer'schen Paradoxons. Dieses besteht darin, dass eine Kubik, definiert über

\begin{equation*}
     a_1 x^3 +a_2y^3 +a_3x^2y+a_4x y^2+ a_5x^2+a_6y^2+a_7xy +a_8x+y_9y+a_{10}=0
\end{equation*}
mit den variablen $x$ und $y$ und den Koeffizienten $a_1$ bis $a_{10}$ im Allgemeinen durch $9$ Punkte eindeutig festgelegt wird, sich jedoch zwei verschiedene Kubiken ebenfalls in $9$ Punkten schneiden. Euler löst dieses in seiner Arbeit \textit{``Démonstration sur le nombre des points, où deux lignes des ordres quelconques peuvent se couper"} (\cite{E148}, 1748) (E148: ``Beweis zur Anzahl der Punkte, in welchen sich zwei Kurven beliebiger Ordnung schneiden können") auf. Dies wird in Sandifers Arbeit \textit{``Cramer’s Paradox"} (\cite{Sa04aug}, 2004) diskutiert.\\

Jedoch ist Eulers Theorie des komplexen Logarithmus zur Demonstration dieses Prozesses noch illustrativer. Denn sie zeigt, wie sehr er über seine Vorgänger -- insbesondere Leibniz (1646--1716) und Johann Bernoulli (1667--1748), aber auch seinen Zeitgenossen d’Alembert (1717--1783) hinausgeht\footnote{Eulers gleich zu besprechende Betrachtung dieses Gegenstandes ist auch ein Musterbeispiel für eine gute Argumentation, wenn man die Maßstäbe der in der Philosophie viel beachteten Bücher zum Gegenstand der Argumentationstheorie \textit{``Attacking Faulty Reasoning: A Practical Guide to Fallacy-Free Arguements"} (\cite{Da13}, 2013) und \textit{``A Practical Study of Argument"} (\cite{Go13}, 2013)  heranzieht.}. Eine Darstellung von Eulers weiteren Leistungen zu diesem Gegenstand findet man auch im Artikel \textit{``Euler, D’Alembert and the Logarithm Function"} (\cite{Br07}, 2007) beschrieben. \\

Die Auflösung der Schwierigkeiten in der Theorie des Logarithmus gelingt Euler durch Aufdeckung von dessen Unendlichwertigkeit. Euler stellt die Berechnung des Logarithmus einer allgemeinen komplexen Zahl in Problem 3 seiner Abhandlung \textit{``De la controverse entre Mrs. Leibnitz et Bernoulli sur les logarithmes des nombres négatifs et imaginaires"} (\cite{E168}, 1751, ges. 1747) (E168: ``Über die Kontroverse der Herren Leibniz und Bernoulli über die Logarithmen von negativen und imaginären Zahlen")\footnote{Eine Übersetzung ins Deutsche findet sich auch in dem Buch \textit{``Ostwalds Klassiker der exakten Wissenschaften Band 261: Zur Theorie komplexer Funktionen
von Leonhard Euler"} (\cite{OK07}, 2007).} vor. Das vorausgehende Problem 1 mit der Behandlung einer beliebigen positiven Zahl und Problem 2 mit einer beliebigen negativen Zahl sind als Vorbereitungen auf den allgemeinen Fall zu verstehen. Die grundlegenden Gedanken Eulers werden an dieser Stelle unter Verwendung moderner Terminologie nachgezeichnet. Der folgende Grenzwert\footnote{Diesen  nutzt Euler unter anderem in seinem berühmten Lehrbuch \textit{``Introductio in analysin infinitorum, volumen primum"} \cite{E101}, um -- unter Verwendung des binomischen Lehrsatzes -- die Taylorreihenentwicklung von $\log (1+x)$ abzuleiten. } bildet für ihn den Ausgangspunkt:

\begin{equation}
\label{eq: Grenzwert ln}
    \log (z) =\lim_{n\rightarrow \infty} n \left(z^{\frac{1}{n}}-1\right),
\end{equation}
 Für den Ausdruck $z^{\frac{1}{n}}$ gilt:

\begin{equation*}
    z^{\frac{1}{n}}=|z|^{\frac{1}{n}}\left(\cos \left(\dfrac{\varphi+2k \pi}{n}\right)+i(\sin \left(\dfrac{\varphi+2k \pi}{n}\right)\right)
\end{equation*}
mit  $\varphi$ als Argument der komplexen Zahl\footnote{Euler benutzt in diesem Zusammenhang meist das Wort ``Winkel".},  $k$ ist eine beliebige ganze Zahl. Für endliches $n$ hat man nur $n$ verschiedene Werte,  für unendliches $n$, so Euler, wird  $z^{\frac{1}{n}}$ folglich unendlich viele Werte aufweisen, was sich demnach auf den Logarithmus übertragt, denn:

\begin{equation*}
    \log (z) =\lim_{n\rightarrow \infty} n \left(|z|^{\frac{1}{n}}\left(\cos \left(\dfrac{\varphi+2k \pi}{n}\right)+i(\sin \left(\dfrac{\varphi+2k \pi}{n}\right)\right) -1\right),
\end{equation*}
Die Funktionalgleichung $\log  (x\cdot y)$ erlaubt die Fokussierung auf den Fall $|z|=1$ ohne die Allgemeinheit einzuschränken:

\begin{equation*}
    \log (z) =\lim_{n\rightarrow \infty}\left( n \left(\cos \left(\dfrac{\varphi+2k \pi}{n}\right)+i(\sin \left(\dfrac{\varphi+2k \pi}{n}\right)\right) -1\right),
\end{equation*}
was zu diesem gleichwertig ist

\begin{equation*}
    \log  (z)= \lim_{n\rightarrow 0} \dfrac{1}{n}\left(\cos(n(\varphi+2k\pi))+i\sin(n(\varphi+2k\pi))-1\right).
\end{equation*}
Eine Aufteilung in

\begin{equation*}
    \log (z)= \lim_{n\rightarrow 0}\dfrac{\cos (n(\varphi+2k\pi))-1}{n}+i \cdot \lim_{n\rightarrow 0} \dfrac{\sin (n(\varphi+2k\pi))}{n}. 
\end{equation*}
führt unter Anwendung bekannten Potenzreihen 

\begin{equation*}
    \cos(x)=1- \dfrac{x^2}{2!}+\dfrac{x^4}{4!}+\cdots \quad \text{und} \quad \sin (x) = x- \dfrac{x^3}{3!}+\dfrac{x^5}{5!}-\cdots
\end{equation*}
zu

\begin{equation*}
    \log  (z)= 0+i \varphi + 2k\pi \cdot i,
\end{equation*}
woraus man für den Fall $|z| \neq 1$ schließt:

\begin{equation}
    \label{eq: log komplexe Zahl}
    \log  (z)= \log  |z| + i \varphi + 2k \pi \cdot i.
\end{equation}
Nimmt man nun noch hinzu, dass sich das Argument $\varphi$ mithilfe von inversen trigonometrischen Funktionen berechnen lässt, drückt Euler den komplexen Logarithmus somit über elementare Ausdrücke aus. Ganz konkret gibt Euler in § 100 seiner Arbeit  \textit{``Recherches sur les racines imaginaires des equations"} (\cite{E170}, 1751, ges. 1746) (E170: ``Untersuchung über die imaginären Wurzeln von Gleichungen") den Ausdruck

\begin{equation}
\label{eq: log complex}
    \log(x+iy) = \log \sqrt{x^2+y^2} +i  \arctan \left(\dfrac{y}{x}\right),
\end{equation}
an, wobei er jedoch nicht nur $\arctan \left(\frac{y}{x}\right)$, sondern auch die gleichwertigen Ausdrücke $\arcsin \left(\frac{y}{\sqrt{x^2+y^2}} \right)$ und $\arccos \left(\frac{x}{\sqrt{x^2+y^2}}\right)$ für das Argument angibt. Er gelangt hier über den Weg der Integration der Differentialform

\begin{equation*}
    d\arctan \left(\dfrac{y}{x}\right)= \dfrac{-ydx+xdy}{x^2+y^2} \quad \text{und} \quad d \log \sqrt{x^2+y^2}= \dfrac{xdx+ydx}{x^2+y^2}
\end{equation*}
zu den besagten Ausdrücken, also auf anderen Wege als in \cite{E168}. Auch hier stellt Euler die Vieldeutigkeit des Logarithmus fest, zieht allerdings diesmal die aus geometrischen Betrachtungen evidente Vieldeutigkeit der inversen trigonometrischen Funktionen als Begründung heran, welche sich  auf den Logarithmus überträgt.\\

\begin{figure}
    \centering
    \includegraphics[scale=0.7]{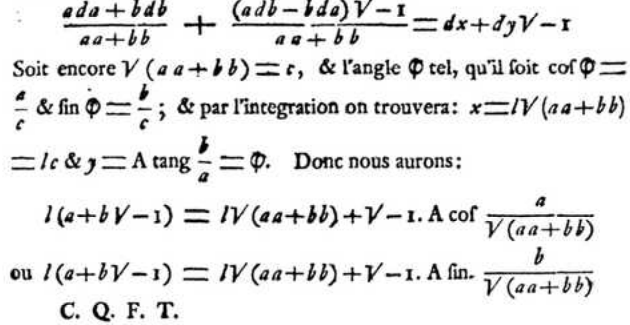}
    \caption{Eulers Herleitung von Ausdrücken für Logarithmen einer komplexen Zahl mithilfe von Differentialformen aus seiner Arbeit \cite{E170}. Euler schreibt hier, wie zumeist, für den natürlichen Logarithmus nur den Formelbuchstaben $l$.}
    \label{fig:E170komplexerln}
\end{figure}

Euler hat indes noch ein weiteres Argument für die Vieldeutigkeit des Logarithmus in seiner Abhandlung \textit{``Sur les logarithmes des nombres negativs et imaginaires"} (\cite{E807}, 1862, ges. 1747) (E807: ``Über die Logarithmen von negativen wie imaginären Zahlen") vorgestellt, was diese Zusammenhänge nicht verwendet. Da dieses ein Paradebeispiel Euler'scher Ingeniosität ist (siehe auch \cite{Du99}), soll es hier nicht vorenthalten werden:\\

Euler beginnt mit dem Integral
\begin{equation*}
    \int\dfrac{dx}{\sqrt{1+x^2}}=\log(x+\sqrt{x^2+1})+C,
\end{equation*}
wobei $C$ eine Integrationskonstante anzeigt. Nun setzt Euler $x=iu$, sodass

\begin{equation*}
    \int \dfrac{idu}{\sqrt{1-u^2}}=\log(iu+\sqrt{1-u^2})+C.
\end{equation*}
Das Integral linker Hand ist  bekanntermaßen eine Arkusfunktion, sodass:

\begin{equation*}
    i \arcsin(u)=\log (iu+\sqrt{1-u^2})+C.
\end{equation*}
Nun setzt Euler $u=\sin \varphi$, sodass

\begin{equation*}
    i \arcsin (\sin (\varphi))=i \varphi = \log(i \sin (\varphi)+ \sqrt{1-\sin^2 (\varphi}))+C
\end{equation*}
\begin{equation*}
    = \log(\cos (\varphi) + i\sin \varphi)+C.
\end{equation*}
Den Wert der Konstante $C$ kann man aus dem Spezialfall $\varphi=0$ als $=0$ ermitteln, woraus die Formel

\begin{equation*}
   \log( \cos (\varphi)+i\sin (\varphi))= i \varphi
\end{equation*}
entspringt. Durch Exponenzieren gelangt man zur berühmten Euler'schen Formel und Eulers Beweis ist abgeschlossen. In der Tat ist dieser auch im Vergleich zum obigen aus moderner Sicht vollkommen streng. Umso mehr mag es verwundern, dass Euler  das erste  der Arbeit \cite{E168} nachempfundene Argument gegenüber dem  just vorgestellten für eine Veröffentlichung den Vorzug gegeben hat\footnote{Die Abhandlung \cite{E807} wurde zwar erst posthum veröffentlicht, aber zeitlich leicht vor \cite{E168} verfasst. Im Vorwort zur Band 19 der \textit{Opera Omnia} \cite{OO19} wird \cite{E807} gar als die sorgsame Ausarbeitung von \cite{E168} gesehen.}. Seine Wahl mag in der womöglich von ihm antizipierten Präferenz seiner Zeitgenossen für auf der Exponentialfunktion fußende Argumentationen begründet liegen, zumal die Wesensverwandschaft dieser zu Logarithmen zunächst enger ist als die zu den Kreisfunktionen.  Bekräftigt wird diese Ansicht dadurch, dass  trotz der Verbindung von Kreis-- und Exponentialfunktionen über die Euler'sche Identität

\begin{equation*}
    e^{i \varphi}= \cos(\varphi)+i \sin (\varphi)
\end{equation*}
 eine  Unterscheidung zwischen ihren Umkehrfunktionen, den Logarithmen und Arkusfunktionen, üblich gewesen zu sein scheint.  Euler selbst spricht insbesondere im Zusammenhang von Integration in späteren Arbeiten oft von ``Logarithmen und Kreisbögen"{} und trennt sie somit ungeachtet ihrer mathematischen Äquivalenz sprachlich voneinander. \\
 
 Explizit erwähnt hat er diese mathematische Verbindung der beiden Funktionen   unter anderem in § 6 der in der Arbeit \textit{``Speculationes analyticae"} (\cite{E475}, 1776, ges. 1774) (E475: ``Analytische Betrachtungen") durch Angabe dieser Formel
\begin{equation*}
    \log \left(\dfrac{1+(p+q)\sqrt{-1}}{{1-(p+q)\sqrt{-1}}} \right)= 2\sqrt{-1}\arctan(p+q),
\end{equation*}
wobei sich $i=\sqrt{-1}$ zu denken ist. $p$ und $q$ sind für Euler in diesem Kontext beliebige reelle Zahlen. Eine Erwähnung des Nexus von Logarithmen komplexer Zahlen zu Kreiszahlen findet sich bereits in  § 17 seiner Ausarbeitung \textit{``De progressionibus harmonicis observationes"} (\cite{E43}, 1740, ges. 1734) (E43: ``Beobachtungen zu harmonischen Progressionen"). Die Formel erläutert Euler in Worten:\\

\textit{``Es ist nämlich $\frac{\sqrt{-3}}{2}\log \frac{3+\sqrt{-3}}{3-\sqrt{-3}}$ die Peripherie des Kreises dividiert durch $\sqrt{3}$, sofern der Durchmesser $=1$ gesetzt worden ist, und $\frac{\sqrt{-3}}{2}\log \frac{1-\sqrt{-3}}{1+\sqrt{-3}}$ die Hälfte dessen."}\\

Euler liegt grundlegend richtig, ihm scheint indes ein Rechenfehler unterlaufen zu sein. Denn man hat, sofern man die Hauptwerte des $\log$ nimmt:

\begin{equation*}
 \dfrac{\sqrt{3}i}{2}\log \left(\dfrac{1-\sqrt{3}i}{1+\sqrt{3}i}\right)= \dfrac{\pi}{\sqrt{3}} \quad \text{sowie} \quad   -\dfrac{\pi}{2\sqrt{3}}= \dfrac{\sqrt{3}i}{2}\log \left(\dfrac{3+\sqrt{3}i}{3-\sqrt{3}i}\right).
\end{equation*}

\subsubsection{Besonderheiten des Präsentationstils}
\label{subsubsec: Besonderheiten der Präsentationsweise}

\epigraph{The mediocre teacher tells. The good teacher explains. The superior teacher demonstrates. The great teacher inspires.}{William Arthur Ward}


Nach Beleuchtung der Facetten von Eulers  Arbeitsweise soll die gleiche Aufmerksamkeit nun seiner Darstellungsweise zukommen. Jedes Merkmal wird dabei exemplarisch illustriert.

\paragraph{Ausführliche Darstellung des Gedankengangs und Präsentation von Rechnungen}
\label{para: Ausführliche Darstellung des Gedankengangs und Präsentation von Rechnungen}

Euler hat selten an der ausführlichen Erläuterungen seiner Gedanken gespart, was ungeachtet der daraus resultierenden Länge seiner Arbeiten das Verständnis des Lernenden fördert. Hankel fasst dies auf Seite 16--17 seines Buches \textit{``Die Entwickelung der Mathematik in den letzten Jahrhunderten"} (\cite{Ha69}, 1869) wie folgt zusammen: \\

\textit{``Er [Euler] war eine wesentlich concrete Natur, die sich mit wirklicher Liebe und Begeisterung dem Stoffe hingab und sich von ihm treiben liess. [...]. Daher geht durch all seine Schriften ein warmer, lebendiger Zug; man ließt zwischen den Zeilen überall begeisterte Freude über die Schönheit und wunderbare Tiefe, die ihm der Gegenstand offenbart. Mit behaglicher Breite, die nicht jedes einzelne Wort ängstlich abwägt, erzählt er, was ihn seine Untersuchungen gelehrt haben -- und so lesen sich -- wie man gesagt hat, seine Schriften wie "Novellen"."}\\

Hankels Beschreibung bestätigt sich in Eulers ausführlichen verbalen Erklärungen, sehr umfassend vorgetragenen Rechnungen und expliziter Exempel, welche zum Ende in einer allgemeinen Formel münden. Obschon diese Merkmale vielen Euler'schen Arbeiten zukommen, möge im Hinblick auf den späteren Teil der vorliegenden Abfassung die folgende Reihe an Euler'schen Abhandlungen als pars pro toto herhalten: \textit{``Theorema maxime memorabile circa formulam integralem $\int \frac{\partial \phi \cos \lambda \phi}{(1+aa-2a\cos \phi)^{n+1}}$"} (\cite{E672}, 1794, ges. 1778) (E672: ``Ein höchst bemerkenswertes Theorem über das Integral $\int \frac{\partial \phi \cos \lambda \phi}{(1+aa-2a\cos \phi)^{n+1}}$"), \textit{``Disquitio coniecturalis super formula integrali $\int \frac{\partial \phi \cos i \phi}{(\alpha +\beta \cos \phi)^n}$"} (\cite{E673}, 1794, ges. 1778) (E673: ``Eine auf Vermutung basierende Untersuchung über die Integralformel $\int \frac{\partial \phi \cos i \phi}{(\alpha +\beta \cos \phi)^n}$"), \textit{``Demonstratio theorematis insignis per coniecturam eruti circa intagrationem formulae $\int \frac{\partial \phi \cos i \phi}{(1+aa-2a \cos \phi)^{n+1}}$"} (\cite{E674}, 1794, ges. 1778)\footnote{Die Arbeiten \cite{E672}, \cite{E673}, \cite{E674} waren ein Teil eines einzigen Abschnitts innerhalb des Buches  \textit{``Institutionum calculi integralis -- Volumen 4"}.} (E674: ``Beweis des außergewöhnlichen Theorems, welches nur auf einer Vermutung basierend über die Integralformel $\int \frac{\partial \phi \cos i \phi}{(1+aa-2a \cos \phi)^{n+1}}$ entdeckt worden ist") und  \textit{``Specimen transformationis singularis serierum"} (\cite{E710}, 1801, ges. 1778) (E710: ``Ein Beispiel einer einzigartigen Transformation von Reihen"). In der Tat scheinen sie nahezu ohne Nachbearbeitung in die Publikation übergegangen zu sein, da sie -- selbst für Euler'sche Verhältnisse -- sehr viele minutiös evaluierte Spezialfälle beinhalten. Den Inhalt bilden dabei die Integrale der Form

\begin{equation}
\label{eq: Familie E672}
    \int\limits_{0}^{\pi} \dfrac{d\varphi \cos(\lambda\varphi)}{(1+2\cos(\varphi)a+a^2)^n}
\end{equation}
für beliebige Zahlen $a$, und natürliche Zahlen $\lambda$ und $n$, deren Geltungsbereich Euler im Verlauf seiner Investigationen  auf die ganzen Zahlen ausweitet. Dies geschieht durch systematische Betrachtung verschiedener Tupel $(\lambda,n)$ für die kleinsten natürlichen Zahlen. Seine Anstrengungen kulminieren in einer Differenzengleichung in den Buchstaben $\lambda$ und $n$ sowie einer Funktionalgleichung in der Variablen $n$, welche die Berechnung der Integrale (\ref{eq: Familie E672}) von ganzzahligen negativen $n$ auf die positiven zurückführt\footnote{Sowohl die Funktionalgleichung als auch die Differenzgleichung werden unten in Abschnitt (\ref{subsubsec: Den Kontext betreffend: Die Legendre Polynome}) Gegenstand der Diskussion bilden, weshalb sie an dieser Stelle noch nicht aufgeführt sind. Es handelt sich respektive um die Gleichungen (\ref{eq: Euler Functional Equation}) und (\ref{eq: Euler's Difference Equation})}. Zunächst als Vermutungen in \cite{E672} und \cite{E673} vorgestellt, enthält das Papier \cite{E674} schließlich die Beweise, welche allerdings auf Ergebnisse von \cite{E710} zurückgreifen. Unten (Abschnitt \ref{subsubsec: Den Kontext betreffend: Die Legendre Polynome}) werden die Integrale (\ref{eq: Familie E672}) demonstriert, die nach Legendre (1752--1833) benannten Polynome auszudrücken, und die Euler'schen Untersuchungen als eine Eruierung einiger ihrer Eigenschaften umgedeutet. Letztgenannte Abhandlung enthält überdies wohl die erste Definition der hypergeometrischen Reihe, welche den Inhalt von Abschnitt (\ref{subsubsec: Die Darstellung betreffend -- Die hypergeometrische Reihe}) bildet. 

\begin{figure}
    \centering
  \includegraphics[scale=0.9]{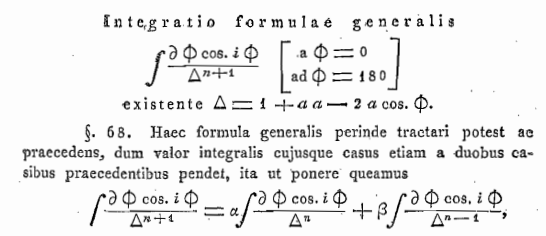} 
    \caption{Die Integralfamilie von Interesse aus Eulers Arbeit \cite{E673}. Die angenommene Differenzengleichung mit zu bestimmenden Koeffizienten $\alpha$ und $\beta$ ist unten zu sehen.}
    \label{fig:E673Diff}
\end{figure}

\paragraph{Mitteilung von kuriosen Ergebnissen}
\label{para: Mitteilung von kuriosen Ergebnissen}

Euler hat selten die Gelegenheit verstreichen lassen,  Ergebnisse mitzuteilen, welche er selbst als kurios empfand, obgleich sie nach seiner Einschätzung nur eine lose Verbindung zum Gestand der eigentlichen Untersuchung haben. Zur Bekräftigung dieser Aussage sei ein Euler'sches Fundstück beigefügt, was in den sogenannten gebrochenen Ableitungen besteht. Diese findet er bereits in seiner Arbeit \textit{``De progressionibus transcendentibus seu quarum termini generales algebraice dari nequeunt"} (\cite{E19}, 1738, ges. 1729) (E19: ``Über transzendente Progressionen oder solche, deren allgemeine Terme algebraisch nicht angegeben werden können"). In Paragraph § 28 entdeckt man, sofern in moderner Sprache präsentiert, den Ausdruck:

\begin{equation}
    \label{eq: Euler gebrochene Ableitung}
    \dfrac{d^n z^{\alpha}}{dz^n}= \dfrac{\Gamma(\alpha +1)}{\Gamma(\alpha+1-n)}z^{\alpha-n},
\end{equation}
wobei $\alpha$ und $n$  nicht ganzzahlig zu sein brauchen. Dabei ist hier

\begin{equation*}
    \Gamma(x):= \int\limits_{0}^{\infty} t^{x-1}e^{-t}dt \quad \text{mit} \quad \operatorname{Re}(x)>0
\end{equation*}
die $\Gamma$-Funktion. Euler leitet den Leser dabei in § 27 der erwähnten Arbeit wie folgt zu seiner Entdeckung:\\

\textit{``Anstelle des Schlussschnörkels möchte ich noch etwas eher Kurioses als Nützliches mitteilen. Es ist freilich bekannt, dass unter $d^nx$  das Differential $n$-ter Ordnung verstanden wird und $dp^n$, wenn $p$ eine beliebige Funktion von $x$ bedeutet und $dx$ als konstant festgelegt wird, in einem endlichen Verhältnis zu $dx^n$ stehen wird; aber immer wenn $n$ eine ganze positive Zahl ist, wird dieses Verhältnis, welches $d^np$ zu $dx^n$ hat, algebraisch ausgedrückt werden können; wie etwa wenn $n=2$ und $p=x^3$ ist, so wird $d^2(x^3)$ sich zu $dx^2$ verhalten wie $6x$ zu $1$. Nun stellt sich aber die Frage, wie das Verhältnis sein wird, wenn $n$ eine gebrochene Zahl ist. Die Schwierigkeit in diesen Fällen sieht man leicht ein, denn wenn $n$ eine ganze positive Zahl ist, wird $d^n$ durch wiederholte Differentiation gefunden, aber ein solcher Weg  steht einem nicht offen, wenn $n$ eine gebrochene Zahl ist. Aber mithilfe der Interpolation von Progressionen, welche ich hier erläutert habe, wird sich die Angelegenheit klären lassen."}\footnote{Euler schreibt in der Originalarbeit beim ersten Auftreten von $dx^n$ noch $d^nx$, geht dann aber im weiteren Text zur heute üblichen Schreibweise $d^nx$ für das Differential im Nenner über."}\\

Was die Kuriosität betrifft, hat Euler Recht behalten, was die geringe Nützlichkeit betrifft dahingegen nicht. Anwendungen werden beispielsweise im Buch \textit{``Special Functions"} (\cite{An11}, 2011) gezeigt. Diese Entdeckung von Euler ist demnach nicht unbeachtet geblieben -- wohl wegen der Wichtigkeit der Arbeit \cite{E19}, welche die erste Definition der $\Gamma$-Funktion im Druck beinhaltet\footnote{Siehe dazu auch Sandiders Artikel \textit{``Gamma the Function"} (\cite{Sa07sep}, 2007).} --, wohingegen die Euler'sche Abhandlung \textit{``Specimen aequationum differentialium indefiniti gradus earumque integrationis"} (\cite{E681}, 1794, ges. 1781) (E681: ``Ein Beispiel einer Differentialgleichung unbestimmten Grades und deren Integration") einen geringeren Bekanntheitsgrad zu besitzen scheint.\\

\begin{figure}
    \centering
     \includegraphics[scale=1.1]{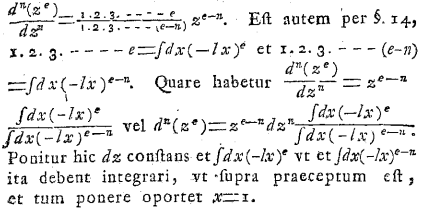}
    \caption{Euler nutzt in seiner Arbeit \cite{E19} die zuvor eruierte Interpolation der Fakultät zur Definition von gebrochenen Ableitungen. Die Integrale sind alle von $x=0$ bis zu $x=1$ erstreckt. }
    \label{fig:E19GebrocheneAbleitung}
\end{figure}

In besagter Abhandlung betrachtet Euler die Differentialgleichung

\begin{equation*}
    \dfrac{d^qf}{dx^q}=g(x),
\end{equation*}
wobei nun $q$ nicht notwendig eine ganze Zahl zu sein braucht. Und vermöge Interpolationstechniken findet er in § 40 -- in moderner Notation geschrieben -- den Ausdruck:

\begin{equation*}
    f(1)= \dfrac{1}{\Gamma(q)}\int\limits_{0}^1 g(x)(1-x)^{q-1}dx. 
\end{equation*}
Aus der Euler'schen Formel findet man leicht die Verallgemeinerung:

\begin{equation}
    \label{eq: Euler gebrochenes Integral}
    f(x)= \dfrac{1}{\Gamma(q)}\int\limits_{0}^x g(t)(x-t)^{q-1}dt,
\end{equation}
welche für ganzzahliges $n$  oft Cauchy zugeschrieben wird, welcher sie in seinen \textit{``Cours d’analyse de l’ école royale polytechnique} (\cite{Ca21}, 1821) (``Vorlesungen zur Analysis an der Ecole royale polytechnique") in der 25. Vorlesung beweist.

\paragraph{Mitteilung bizarrer sowie unrichtiger Ergebnisse}
\label{para: Auch Falsches bzw. Bizarres wird mitgeteilt}

Ein weiteres Charakteristikum von Eulers Präsentation mathematischer Sachverhalte besteht darin,  wissentlich Formeln mitzuteilen, die in einigen Fällen  freilich das richtige Ergebnis liefern, jedoch in anderen unrichtig sind. Dies quittiert Euler an entsprechender Stelle anschließend meist mit einer Bemerkung, die Mitteilung gemacht zu haben, um den Forschergeist seiner Kollegen zu wecken oder Ähnliches.  Ein konkretes Beispiel dessen macht man in seiner Arbeit \textit{``Consideratio progressionis cuiusdam ad circuli quadraturam inveniendam idoneae"} (\cite{E125}, 1750, ges. 1739) (E125: ``Die Betrachtung einer gewissen Progression, welche sich zur Kreisquadratur eignet") aus, wo das mitzuteilende Kuriosum die letzten drei Paragraphen einnimmt (§§ 17--20). Euler beginnt mit der Reihe

\begin{equation*}
    a_1+a_2+a_3+\cdots
\end{equation*}
bis ins Unendliche\footnote{Euler benutzt und schreibt keine Indizes $a+b+c+d+e+f+g+h+\text{etc.}$. Zur Erleichterung der Darstellung wird jedoch weiter die Indexschreibweise verwendet.}. Aus der Differenzenrechnung lässt sich nach der Newton'schen Interpolationsformel\footnote{Man findet die entsprechende Formel bereits in Newtons \textit{``Philosophiae Naturalis Principia Mathematica"} (\cite{Ne87}, 1687) (``Mathematische Prinzipien der Naturphilosophie").} die Summe bis hin zu $x$ Termen wie folgt darstellen:

\begin{equation*}
    a_x = a_1+ \dfrac{x-1}{1}(a_2-a_1)+\dfrac{(x-1)(x-2)}{1\cdot 2}(a_3-2a_2+a_1)
\end{equation*}
\begin{equation*}
    +\dfrac{(x-1)(x-2)(x-3)}{1\cdot 2 \cdot 3}(a_4-3a_3+3a_2-a_1)+\text{etc.},
\end{equation*}
was sich mit Differenzenoperator $\Delta^n$ für die $n$--te Differenz  wie folgt schreiben lässt:

\begin{equation*}
    a_x = \sum_{n=0}^{x-1} \binom{x-1}{n}\Delta^n a_1. 
\end{equation*}
Diese Formel benutzt Euler, um die Terme für negative ganze Zahlen und $0$ als Index zu bestimmen. Durch Summierung all dieser Terme findet auf formalen Wege:

\begin{equation*}
    \sum_{n=0}^{\infty} f(n)\Delta^n a_1
\end{equation*}
mit

\begin{equation*}
    f(n) = \lim_{x \rightarrow 1} \dfrac{1}{(1-x)^{n+1}},
\end{equation*}
sodass die einzelnen Terme der Reihe divergieren. Bei Euler findet man die entsprechende Identität wie folgt niedergeschrieben:

\begin{equation*}
    \dfrac{a_1}{1-1}+\dfrac{a_2-a_1}{(1-1)^2}+\dfrac{a_3-2a_2+a_1}{(1-1)^3}+\dfrac{a_4-3a_2+3a_2-a_1}{(1-1)^3}+\text{etc.}
\end{equation*}
Diesen formalen Ausdruck ordnet Euler um und gelangt zu:

\begin{equation*}
\renewcommand{\arraystretch}{2,0}
\setlength{\arraycolsep}{0.5mm}
\begin{array}{rlclclclclclclclcl}
     +a_1 &\bigg(\dfrac{1}{(1-1)^1}&+&\dfrac{1}{(1-1)^2}&+&\dfrac{1}{(1-1)^3}&+&\dfrac{1}{(1-1)^4}&+&\text{etc.}\bigg)\\
      -a_2 &\bigg(\dfrac{1}{(1-1)^2}&+&\dfrac{2}{(1-1)^3}&+&\dfrac{3}{(1-1)^4}&+&\dfrac{4}{(1-1)^5}&+&\text{etc.}\bigg)\\
       +a_3 &\bigg(\dfrac{1}{(1-1)^3}&+&\dfrac{3}{(1-1)^4}&+&\dfrac{6}{(1-1)^4}&+&\dfrac{10}{(1-1)^{6}}&+&\text{etc.}\bigg)\\
        -a_4 &\bigg(\dfrac{1}{(1-1)^4}&+&\dfrac{4}{(1-1)^5}&+&\dfrac{10}{(1-1)^6}&+&\dfrac{20}{(1-1)^7}&+&\text{etc.}\bigg)\\
\end{array}
\end{equation*}
Euler identifiziert  die Ausdrücke in Klammern jeweils als entsprechende Ableitung der geometrischen Reihe, sodass er für die ersten Zeile findet:

\begin{small}
\begin{equation*}
    +a_1 \bigg(\dfrac{1}{(1-1)^1}+\dfrac{1}{(1-1)^2}+\dfrac{1}{(1-1)^3}+\dfrac{1}{(1-1)^4}+\text{etc.}\bigg)= a_1 \dfrac{1}{(1-1)-1}=-a_1
\end{equation*}
\end{small}sowie für die zweite
\begin{small}
\begin{equation*}
    -a_2 \bigg(\dfrac{1}{(1-1)^2}+\dfrac{2}{(1-1)^3}+\dfrac{3}{(1-1)^4}+\dfrac{4}{(1-1)^5}+\text{etc.}\bigg)=-a_2 \dfrac{1}{(1-1)^2-1}=+a_2
\end{equation*}
\end{small}und die dritte

\begin{small}
   \begin{equation*}
    +a_3 \bigg(\dfrac{1}{(1-1)^3}+\dfrac{3}{(1-1)^4}+\dfrac{6}{(1-1)^4}+\dfrac{10}{(1-1)^{6}}+\text{etc.}\bigg)=+a_3 \dfrac{1}{(1-1)^3-1}=-a_3
\end{equation*} 
\end{small}und schließlich die vierte

\begin{small}
    \begin{equation*}
    -a_4 \bigg(\dfrac{1}{(1-1)^4}+\dfrac{4}{(1-1)^5}+\dfrac{10}{(1-1)^6}+\dfrac{20}{(1-1)^7}+\text{etc.}\bigg)=-a_4 \dfrac{1}{(1-1)^4-1}=a_4,
\end{equation*}
\end{small}sodass das Muster schnell sichtbar wird und die Divergenzen nicht mehr zu finden sind. Schlussendlich gelangt Euler zur Identität

\begin{equation*}
    \sum_{n=0}^{\infty} a_{-n}= - \sum_{n=1}^{\infty}a_n,
\end{equation*}
was ihn zu folgender Konklusion bewegt:

\begin{equation*}
    \sum_{n=-\infty}^{\infty}a_n=0.
\end{equation*}
Euler drückt dies in Worten am Ende von  § 18 wie folgt aus:\\

\textit{``Wenn also irgendeine unendliche Reihe}

\begin{equation*}
    a+b+c+d+e+\text{etc.}
\end{equation*}
\textit{auch nach links ins Unendliche fortgesetzt würde, wäre die Summe  nach beiden Seiten hin ins Unendliche fortgesetzt immer $=0$; wenn diese Begründung freilich richtig wäre."}\\

Trotz der offenkundigen Invalidität  der Formel führt sie, wie Euler auch zeigt, für die Wahl $a_n=k^n$ zum richten Ergebnis, da unter völliger Missachtung von Konvergenzbetrachtungen gilt

\begin{equation*}
    \sum_{k=1}^{\infty} k^n = \dfrac{k}{1-k} \quad \text{und} \quad \sum_{k=0}^{\infty} k^{-n}= \dfrac{k}{k-1}
\end{equation*}
und

\begin{equation*}
    \dfrac{k}{1-k}+\dfrac{k}{k-1}=0.
\end{equation*}
Jedoch führt die Wahl $a_n= \frac{1}{(2n-1)^2}$ zum falschen Ergebnis, da

\begin{equation*}
    \sum_{k=1}^{\infty} \dfrac{1}{(2k-1)^2}= \dfrac{\pi^2}{4}
\end{equation*}
Euler schließt seine Arbeit mit den Worten:\\

\textit{``Ich glaube, all dies vorgestellt zu haben, wird nicht weniger Nutzen haben als mit größter Strenge bewiesene Wahrheiten."}\\

In gewisser Weise sollte Euler Recht behalten, weil der Durchgang durch das Divergente ihn und seine Nachfolger bisweilen zu interessanten Entdeckungen geführt hat, was in den Abschnitten über die $\zeta$-Funktion (Abschnitt \ref{subsubsec: Durch Kombinieren von Ergebnissen: Die zeta-Funktion}) und über seine Theorie divergenter Reihen (Abschnitt \ref{para: Eulers Definition  von divergenten Reihen}) evident werden wird. 

\paragraph{Mithilfe formaler Rechnungen zu richtigen Erkenntnissen}
\label{para: Mithilfe formaler Rechnungen zu richtigen Erkenntnissen}

Eulers Arbeiten zum Gebiet der Partitionen sind ein treffliches Beispiel für die Illustration dessen, wie sich aus rein formalen Rechnungen ebenfalls Theoreme ableiten lassen\footnote{Solche formalen Operationen bespricht Sandifer auch kurz in seinem Artikel \textit{``Formal sums and products"} (\cite{Sa06jul}, 2006).}. Er erläutert seine Ideen diesbezüglich in Kapitel 16 seiner \textit{Introductio} \cite{E101}, man konsultiere aber auch seine Einzelarbeiten \textit{``De partitione numerorum"} \cite{E191} (E191: ``Über die Partition von Zahlen") sowie \textit{``De partitione numerorum in partes tam numero quam specie datas"} \cite{E394} (E394: ``Über die Partition von Zahlen in so von der Anzahl wie der Gattung her vorgegebene Zahlen"), in welchen Arbeiten man überdies Grundsteine  der Theorie der erzeugenden Funktionen erkennen kann. Eulers Einsicht und das sich daraus ableitende Vorgehen besteht darin, gewünschte zahlentheoretische Eigenschaften aus den Koeffizienten von gewissen Potenzreihen abzuleiten. Hier soll Eulers Nachweis für die Eindeutigkeit der Darstellbarkeit   einer jeden Zahl mit den Potenzen von $3$ in §§ 330--331 der \textit{Introductio} \cite{E101} als Beispiel dienen. Euler betrachtet den Ausdruck:

\begin{equation*}
    f(x)= (x^{-1}+1+x^1) (x^{-3}+1+x^3) (x^{-9}+1+x^9) (x^{-27}+1+x^{27})\cdots
\end{equation*}
und setzt an, dass gilt

\begin{equation}
\label{eq: Ansatz Partitionen 3}
    f(x)=\sum_{n=-\infty}^{\infty} a_n x^n,
\end{equation}
was, moderne Terminologie gebrauchend, eine Laurententwicklung um den Ursprung mit zu bestimmenden Koeffizienten $a_n$  ist. Euler bemerkt nun aus der Produktdarstellung, dass 

\begin{equation*}
    f(x)= (x^{-1}+1+x^{1})f(x^3)
\end{equation*}
gilt, sodass für den Ausdruck rechter Hand folgt:

\begin{equation*}
    (x^{-1}+1+x^{1})\sum_{n=-\infty}^{\infty} a_n x^{3n}=  \sum_{n=-\infty}^{\infty} a_n x^{3n-1}+\sum_{n=-\infty}^{\infty} a_n x^{3n}+\sum_{n=-\infty}^{\infty} a_n x^{3n+1}.
\end{equation*}
Zumal dieser Ansatz nun mit (\ref{eq: Ansatz Partitionen 3}) gleichwertig sein muss, entspringen nachstehende Gleichheiten für die Koeffizienten:

\begin{equation*}
    a_n=a_{3n-1} \quad \text{und} \quad a_n=a_{3n} \quad \text{und ebenso} \quad a_n=a_{3n+1},
\end{equation*}
welche für alle $n \in \mathbb{Z}$ ihre Gültigkeit beibehalten müssen. Mit dem offenkundigen Wert $a_0=1$ findet man, dass $a_n=1$ für jedwede ganze Zahl $n$ gelten muss. Damit ist der Euler'sche Beweis vollendet. Euler hat im Verlauf seiner Karriere viele andere zahlentheoretische Tatsachen mit dieser Methode bewiesen\footnote{Ein Beispiel: Die Anzahl der Partitionen einer Zahl $n$ in  ungerade Zahlen ist gleich der Anzahl an Partitionen derselben Zahl $n$ in lauter verschiedene Zahlen.}.  Für ausführlichere moderne Darstellung konsultiere man etwa die Bücher \textit{``The Theory of Partitions"} (\cite{An08}, 2008) und auch \textit{``An Introduction To The Theory Of Numbers"} (\cite{Ha09}, 2009). Für eine kleine Diskussion der Ätiologie der Partitionstheorie bei Euler konsultiere man auch den Artikel \textit{``Philip Naudé’s problem"} (\cite{Sa05oct}, 2005). \\

Es sei an dieser Stelle nochmals auf die Nebensächlichkeit der Konvergenz der resultierenden Reihe hingewiesen, denn die Reihe

\begin{equation*}
     f(x)=\sum_{n=-\infty}^{\infty}  x^n
\end{equation*}
konvergiert tatsächlich für keinen einzigen Punkt in der gesamten komplexen Ebene. Dennoch ergibt sich allein aus den Koeffizienten bereits das zu Beweisende\footnote{Die Summe hatte Euler in \cite{E125} in einem anderen Zusammenhang zu $=0$ summiert. Man vergleiche die Ausführung des letzten Paragraphen über die von Euler mitgeteilten Kuriositäten.}.

\paragraph{Mitteilung von Misserfolgen}
\label{para: Mitteilung von Misserfolgen}

Die letzte in dieser Ausarbeitung vorgestellte Eigentümlichkeit von Eulers Arbeitsweise soll die vielleicht heutzutage unüblichste aller erläuterten sein: Die vollständige Darstellung von Rechnungen und Gedankengängen, \textit{obwohl} sie schließlich nicht zum gewünschten Ziel führen. Ein epitomes Exempel dessen sind die Euler'schen Versuche der Evaluation der Summe

\begin{equation*}
    \zeta(3):= 1+\dfrac{1}{2^3}+\dfrac{1}{3^3}+\dfrac{1}{4^3}+\text{etc.}
\end{equation*}
oder der Apéry'schen Konstante, wie sie seit Apérys (1916--1994) Nachweis ihrer Irrationalität in seiner Abhandlung \textit{``Irrationalité de $\zeta(2)$ et $\zeta(3)$"} (\cite{Ap79}, 1979) (``Die Irrationalität von $\zeta(2)$ und $\zeta(3)$") bezeichnet zu werden pflegt. Bis zum heutigen Tage sind alle Bestrebungen, diese Zahl mithilfe bereits bekannter Konstanten auszudrücken, fehlgeschlagen. Dieses Problem hat auch Euler zeitlebens nicht losgelassen. Wie im folgenden Abschnitt (\ref{subsec: Die Lösung des Baseler Problems}) aufgezeigt werden wird, gelang Euler die Summation im Fall der Quadrate auf mannigfaltige Weise; ein Erfolg, welcher zu seinen größten Triumphen gerechnet werden darf. Bezüglich der Kuben konstatiert Euler nach einem weiteren Fehlschlag in § 17 von \cite{E352}:\\

\textit{``Bezüglich der Reihen der Reziproken der Potenzen}

\begin{equation*}
    1-\dfrac{1}{2^n}+\dfrac{1}{3^n}-\dfrac{1}{4^n}+\dfrac{1}{5^n}-\frac{1}{6^n}+\text{etc.}
\end{equation*}

\textit{habe ich schon angemerkt, dass ihre Summen lediglich in den Fällen angegeben werden können}\footnote{Euler hat seine Arbeit \cite{E352} in französischer Sprache verfasst und benutzt an dieser Stelle stelle das Wort ``savoir", welches mit ``können"{} übersetzt werden kann, dann allerdings auf ein durch fehlendes Wissen begründetes nicht bestehendes Vermögen verweist.}\textit{, in denen der Exponent $n$ eine ganze gerade Zahl ist, jedoch für die Fälle, wo $n$ eine ungerade Zahl ist, all meine Bestrebungen diesbezüglich ertraglos geblieben sind.``}\\

Die nachfolgenden Untersuchungen in derselben Arbeit, obwohl sie für Euler selbst zunächst vielversprechend erschienen, sind ebenfalls nicht von Erfolg gekrönt. In § 19 resümiert er:\\

\textit{``Es mag freilich möglich sein, die Summen dieser zusammengesetzten Reihen der Form}

\begin{equation*}
    1^{2\lambda} \log(1)-2^{2\lambda} \log(2)+3^{2\lambda} \log(3)-4^{2\lambda} \log(4)+\text{etc.}
\end{equation*}
\textit{zu finden. Allerdings ist dieses Unterfangen womöglich schwieriger als das gegenwärtige und ich kann keine Anzeichen einer Methode erkennen, die uns zum gesteckten Ziel führt."}\footnote{Das gesteckte Ziel besteht hier nach wie vor in der Summation der Reihen für ungerades $n$.}\\

 Euler hat in seinen  Investigationen zu diesem Gegenstand jedoch auch Formeln entdeckt, von denen ausgewählte hier Erwähnung finden sollen. Nennt man

\begin{equation*}
    \chi(s):= \sum_{n=0}^{\infty} \dfrac{1}{(2n+1)^s},
\end{equation*}
teilt Euler in § 38 seiner Arbeit \cite{E130} die Identität

\begin{equation*}
    1-\dfrac{1}{2^3}+\dfrac{1}{3^3}-\dfrac{1}{4^3}+\text{etc.}
\end{equation*}
\begin{equation*}
    = \dfrac{\pi^2 \log(2)}{1\cdot 2 \cdot 3}- 2\pi^2 \left(\dfrac{\chi(2)}{2\cdot 3 \cdot 4 \cdot 5}+\dfrac{\chi(4)}{4\cdot 5 \cdot 6 \cdot 7}+\dfrac{\chi(6)\pi^8}{6 \cdot 7 \cdot 8 \cdot 9}+\text{etc.}\right),
\end{equation*}
mit, welche sich von Euler mit anderen Formelbuchstaben ausgedrückt ebenfalls in § 20 von \cite{E432} findet. In letztgenannter Arbeit gelangt Euler überdies in § 21 zur kompakten Formel:

\begin{equation*}
   \chi(3)= \dfrac{1}{1^3}+\dfrac{1}{3^3}+\dfrac{1}{5^3}+\dfrac{1}{7^3}+\text{etc.}= \dfrac{\pi^2}{4}\log(2)+2 \int\limits_{0}^{\frac{\pi}{2}} \phi \log (\sin (\phi))d\phi.
\end{equation*}
Nicht zuletzt diese Formel mag Euler in § 16 seiner Arbeit \textit{``De relatione inter ternas pluresve quantitates instituenda"} (\cite{E591}, 1785, ges. 1775) (E591: ``Über das Herstellen einer Beziehung zwischen drei oder mehr Größen") zur Vermutung verlasst haben, es könnte

\begin{equation*}
    a \chi(3) + b (\log(2))^3+c \log (2) \dfrac{\pi^2}{6} = 0
\end{equation*}
mit ganzen Zahlen $a$, $b$, $c$ gelten. Jedoch führen seine Untersuchungen ihn nicht zu gewünschten Ziel, sodass er sie nicht weiterverfolgt.

\newpage

\section{Eine Illustration von Eulers Arbeitsweise anhand ausgewählter Beispiele}
\label{sec: Eulers Arbeitsweise anhand ausgewählter Beispiele}

\epigraph{Es ist wunderbar, dass man noch heute jede der Abhandlungen Eulers nicht bloss mit Belehrung, sondern mit Vergnügen liest.}{Carl Gustav Jacob Jacobi}

Wohingegen im vorherigen Abschnitt die Arbeitsweise allgemeiner diskutiert wurde, soll das Besprochene in diesem Abschnitt anhand zweier umfassender diskutierter Exempel noch stärker bekräftigt werden.  Das erste Beispiel ist die Euler'sche Lösung des Baseler Problems (Abschnitt \ref{subsec: Die Lösung des Baseler Problems}), das zweite widmet sich der Euler'schen Betrachtung von Differenzengleichungen mit konstanten Koeffizienten (Abschnitt \ref{subsec: Ein Prototyp -- Differenzengleichungen}). Zwischen den beiden genannten Exempeln wird die Betrachtung eines  richtigen, jedoch falsch abgeleiteten, Resultats behandelt (Abschnitt \ref{subsec: Ein falsches Ergebnis}).

\subsection{Die Lösung des Baseler Problems}
\label{subsec: Die Lösung des Baseler Problems}

\epigraph{Every mathematical discipline goes through three periods of development: the naive, the formal, and the critical.}{David Hilbert}

Ein sehr bekanntes Beispiel für den Erfolg von Eulers Arbeitsweise ist die Lösung des sogenannten Baseler Problems, dies ist die Summierung der Reihe

\begin{equation}
    \label{eq: Baseler Problem}
    1+\dfrac{1}{4}+\dfrac{1}{9}+\dfrac{1}{16}+\cdots = \sum_{k=1}^{\infty} \dfrac{1}{k^2}= \dfrac{\pi^2}{6},
\end{equation}
welches die gesamte Bandbreite von Eulers Einfallsreichtum und zugleich Hartnäckigkeit demonstriert.
Lediglich auf seine publizierten Abhandlungen blickend sind bis zum ersten lückenlosen Beweis für die Richtigkeit der obigen Summe die Arbeiten \textit{``De summatione innumerabilium progressionum} (\cite{E20}, 1738, ges. 1731) (E20: ``Über die Summation von unzähligen Progressionen"), \textit{``De summis serierum reciprocarum"} (\cite{E41}, 1740, ges. 1735) (E41: ``Über die Summe der Reihen von Reziproken), \textit{``De summis serierum reciprocarum ex potestatibus numerorum naturalium ortarum dissertatio altera, in qua eaedem summationes ex fonte maxime diverso derivantur"} (\cite{E61}, 1743, ges. 1742) (E61: ``Ein andere Dissertation über die Summen der Reihen der Reziproken, welche aus den Potenzen der natürlichen Zahlen entspringen, in welcher dieselben Summationen aus einer gänzlich anderen Quellen abgeleitet werden"), und \textit{``Demonstration de la somme de cette suite $1+\frac{1}{4}+\frac{1}{9}+\frac{1}{16}+\cdots$"} (\cite{E63}, 1743, ges. 1741) (E63: ``Ein Beweis der Summe der Reihe $1+\frac{1}{4}+\frac{1}{9}+\frac{1}{16}+\cdots$") zu nennen, welche einen Zeitraum von 10 Jahren umspannen. Die Teilschritte bis zum einem modernen Ansprüchen genügenden Nachweis sind:

\begin{itemize}
    \item[1.] Transformation der Reihe aus (\ref{eq: Baseler Problem}) in § 22 von \cite{E20}
    \item[2.] Entdeckung des exakten Werts in \cite{E41}
    \item[3.] Ein alternativer Beweis in \cite{E61}
    \item[4.] Ein ``moderner"{} Beweis in \cite{E63}
\end{itemize}
Da die Geschichte der Lösung des Baseler Problem schon mehrfach dargestellt worden ist, -- man siehe z.B. das Buch von Dunham \textit{``Euler: The Master of Us All"}  (\cite{Du99}, 1999) und die Artikel aus Sandifers Kolumne  \textit{``Estimating the Basel Problem"} (\cite{Sa03dec}, 2003)\footnote{Diese Arbeit zeichnet den ersten Schritt der obigen Liste nach.} und \textit{``The Basel Problem with Integrals"} (\cite{Sa04mar}, 2004)\footnote{Diese Arbeit betrachtet den 4. Schritt aus obiger Liste.} --, erfolgt in dieser Abhandlung eine  komprimierte  Darstellung.

\subsubsection{Transformation der Reihe}
\label{subsubsec: Transformation der Reihe}

\epigraph{I think this is the start of something really big. Sometimes that first step is the hardest one, and we've just taken it.}{Steve Jobs}

Die langsame Konvergenz der Reihe aus (\ref{eq: Baseler Problem}) disqualifiziert sie in dieser Form  für numerische Untersuchungen ihres Wertes, weswegen sich Euler in § 22 von \cite{E20} einer konvergenzbeschleunigenden Transformation bedient. Ausgehend vom Integral

\begin{equation*}
   I= -\int\limits_{0}^{\frac{1}{2}} \dfrac{\log(1-x)}{x}dx.
\end{equation*}
gelangt er über die Potenzreihenentwicklung des Logarithmus zunächst zu

\begin{equation*}
    I= \int\limits_{0}^{\frac{1}{2}} \sum_{n=1}^{\infty} \dfrac{x^{n-1}}{n}dx= \sum_{n=1}^{\infty} \int\limits_{0}^{\frac{1}{2}} \dfrac{x^{n-1}}{n}dx.
\end{equation*}
Termweises Integrieren liefert:

\begin{equation}
\label{eq: I1}
    I= \sum_{n=1}^{\infty} \dfrac{1}{n^2\cdot 2^{n}}.
\end{equation}
Andererseits gibt die Substitution $1-x=u$ und somit $-dx=du$, wenn anschließend wieder $x$ statt $u$ im Integral geschrieben wird:

\begin{equation*}
    I= -\int\limits_{1}^{\frac{1}{2}} \dfrac{\log(x)}{1-x}dx.
\end{equation*}
Mit der geometrischen Reihe ist dies:

\begin{equation*}
I= - \sum_{n=0}^{\infty} \int\limits_{1}^{\frac{1}{2}} x^n \log(x)dx.
\end{equation*}
Die Berechnung des Integrals führt Euler nicht explizit aus, was an dieser Stelle nachgereicht wird. Es gilt:

\begin{equation*}
    \int\limits_{\frac{1}{2}}^{1} x^n dx = \dfrac{1-\left(\frac{1}{2}\right)^{n+1}}{n+1}.
\end{equation*}
Differenziert man nun beide Seiten nach $n$, erhält man auch:

\begin{equation*}
    \int\limits_{\frac{1}{2}}^{1} x^n \log(x) dx = \dfrac{\log(2)}{2^{n+1}\cdot (n+1)}+\dfrac{1}{2^{n+1}(n+1)^2}-\dfrac{1}{(n+1)^2},
\end{equation*}
somit

\begin{equation*}
    \renewcommand{\arraystretch}{2,5}
\begin{array}{rcL}
   I  & = & \sum_{n=0}^{\infty} \left(-\dfrac{\log(2)}{2^{n+1}\cdot (n+1)}-\dfrac{1}{2^{n+1}(n+1)^2}+\dfrac{1}{(n+1)^2}\right)   \\
     & = & -\sum_{n=0}^{\infty} \dfrac{\log(2)}{2^{n+1}\cdot (n+1)} - \sum_{n=0}^{\infty} \dfrac{1}{2^{n+1}(n+1)^2} + \sum_{n=0}^{\infty}  \dfrac{1}{(n+1)^2} \\
     & = & -\log(2) \sum_{n=0}^{\infty} \dfrac{1}{2^{n+1}\cdot (n+1)} - \sum_{n=0}^{\infty} \dfrac{1}{2^{n+1}(n+1)^2} + \sum_{n=0}^{\infty}  \dfrac{1}{(n+1)^2}.
\end{array}
\end{equation*}
Die erneute Verwendung der Potenzreihe des Logarithmus leitet einen zu:

\begin{equation}
    \label{eq: I2}
    I= -\log^2(2)- \sum_{n=0}^{\infty} \dfrac{1}{2^{n+1}(n+1)^2}+ \sum_{n=0}^{\infty}  \dfrac{1}{(n+1)^2}.
\end{equation}
Setzt man  (\ref{eq: I1}) mit diesem Ausdruck gleich, erreicht man nach einer Umformung -- wie Euler -- die Formel:

\begin{equation*}
    \sum_{n=1}^{\infty} \dfrac{1}{n^2}= \sum_{n=1}^{\infty} \dfrac{1}{n\cdot 2^{n-1}}+\dfrac{\log^2(2)}{2}.
\end{equation*}
$\log^2(2)$ bzw. $\log(2)$ lässt sich Tabellen entnehmen. Überdies konvergiert die Reihe auf der rechten Seite  schneller als die auf der linken, sodass sie Euler befähigt, in § 22 von \cite{E20} den numerischen Wert

\begin{equation*}
    1+\dfrac{1}{4}+\dfrac{1}{9}+\dfrac{1}{16}+\cdots \approx 1,644934
\end{equation*}
anzugeben\footnote{Extensive numerische Berechnungen sind, wie zuvor schon einmal angesprochen, ein Charakteristikum für Eulers Arbeitsweise. Zahlreiche seiner Artikel enthalten entweder Passagen mit ausgiebigen Berechnungen oder sind ganz Methoden gewidmet, wie man selbige noch beschleunigen kann. Speziell im Kontext von $\pi$--Approximationen kann man die Arbeiten  \cite{E125}  sowie \textit{``Annotationes in locum quendam Cartesii ad circuli quadraturam spectantem"} (\cite{E275}, 1763, ges. 1758) (E275: ``Anmerkungen zur Descartes'schen Quadratur des Kreises") und weitere nennen.}. 

\subsubsection{Ermittlung des exakten Wertes -- mit formalen Methoden}
\label{subsubsec: Ermittlung des exakten Wertes}

\epigraph{No great discovery was ever made without a bold guess.}{Isaac Newton}

Die Arbeit \cite{E41} enthält bereits den Wert aus (\ref{eq: Baseler Problem}). Jedoch dürfte selbst zu Eulers Zeiten sein Argument nicht als Beweis dienen, welcher berechtigten Kritik Euler  auch in \cite{E61} entsprechend entgegnet ist. Genauer schreibt er in § 6 des letztgenannten Papiers auf die Einwände bezugnehmend:\\

\textit{``Diese fast vollkommen in Vergessenheit geratene Sorge hat neulich ein vom hoch geehrten Herrn Daniel Bernoulli empfangener Brief in mir erneuert, in welchem er mir Gründe darlegt, meine Methode zu bezweifeln, und zugleich deutet er an, dass der hoch geehrte Herr Cramer dieselben Zweifel hegt, ob er meine Methode billigen kann."}\\

Geteilt wird diese Kritik auch von Polya in seinem Buch \cite{Po14} bezüglich dieses Beweises:\\

\textit{``Euler's step was daring. In strict logic, it was an outright fallacy [...] Yet it was justified by analogy, by the analogy of the most successful achievements of a rising science that he called [...] ``Analysis of the Infinite."{} Other mathematicians, before Euler, passed from finite differences to infinitely small differences, from sums with a finite number of terms to sums with an infinity of terms, from finite products to infinite products. And so Euler passed from equations of a finite degree (algebraic equations) to equations of infinite degree, applying the rules made for the finite [...]."}\\

Eulers Ansatz in der Arbeit \cite{E41} besteht   darin, die für Polynome gültigen Relationen zwischen deren Koeffizienten und Nullstellen ohne weitere Rechtfertigung auf unendliche Reihen auszudehnen. Zur Lösung des Baseler Problems bedient Euler sich der Taylorreihe des Sinus

\begin{equation*}
    \sin(x)= x-\dfrac{x^3}{3!}+\dfrac{x^5}{5!}+\cdots
\end{equation*}
und \textit{behauptet}  an dieser Stelle die Gültigkeit der alternativen Darstellung 

\begin{equation}
    \label{eq: sine-Product}
    \sin (x) =  x \prod_{k=1}^{\infty} \left(1-\dfrac{x^2}{(k\pi)^2}\right)
\end{equation}
und gewinnt durch Koffizientenvergleich der beiden Ausdrücke die Summe (\ref{eq: Baseler Problem}).
Denn die formale Multiplikation der führenden Terme des Sinus--Produktes gibt:

\begin{equation*}
    \sin (x) = x - \left(\dfrac{1}{1^2\pi^2}+\dfrac{1}{2^2\pi^2}+\dfrac{1}{3^2\pi^2}+\cdots\right)x^3+\text{etc.}.
\end{equation*}
Euler verwendet demnach  den Satz von Vieta (1540--1603) für den Fall eines  Polynoms mit ``unendlichem" Grad\footnote{Euler hat sich weit später in seiner Laufbahn die umgekehrte Frage gestellt, wie sich von der rechten Seite von (\ref{eq: sine-Product}) die Gleichheit mit dem Sinus auf direktem Weg demonstrieren lässt. Die Antwort gibt er in der Arbeit \textit{``Exercitatio analytica"} (\cite{E664}, 1794, ges. 1776) (E664: ``Eine analytische Übung"), aus Gründen der Kürze der Darstellung, jedoch für den Kosinus und nicht den Sinus.}.  Bei dieser Begebenheit führt ihn sein Vorgehen zu einem richtigen Ergebnis, an anderer Stelle sollte ihn in ein analoger Übergang ins Unendliche indes zu Fehlern verleiten, was in Abschnitt (\ref{subsec: Ein falsches Ergebnis}) noch detaillierter besprochen werden wird.

\subsubsection{Ein strengerer Beweis}
\label{subsubsec: Ein strengerer Beweis}

\epigraph{Ever tried. Ever failed. No matter. Try again. Fail again. Fail better.}{Samuel Beckett}

Die  aus der Kritik von Cramer und D. Bernoulli (1700--1782) bereits zu erahnende mangelnde Validität der Herleitung des Sinus--Produktes (\ref{eq: sine-Product}) betrifft den fehlenden Nachweis der Nichtexistenz von weiteren Nullstellen der $\sin$--Funktion, was Euler implizit voraussetzt\footnote{Einen ähnlichen Fehler begeht Euler auch anderenorts, an welchem er versuchte aus den Nullstellen des charakteristischen ``Polynoms"{} die Lösung einer Differentialgleichung unendlicher Ordnung zu konstruieren. Dies wird Abschnitt (\ref{subsec: Ein falsches Ergebnis}) diskutiert werden.}. Daneben zeigt das Beispiel $e^x \sin x$ die nichteindeutige Festlegung einer transzendenten Funktion mithilfe ihrer Nullstellen allein auf. \\

Diesen Einwänden begegnet Euler nun in seiner Arbeit \cite{E61} durch Nachweis des Sinusproduktes (\ref{eq: sine-Product}) mit den Mitteln der Infinitesimalrechnung. Die grundlegende Idee besteht in der allgemeinen Zerlegung des Ausdrucks $a^n-b^n$ in trinomiale Faktoren der Form $a^2-2ab\cos \zeta +b^2$. Anschließend hat  man lediglich noch $a=1+\frac{ix}{n}$ und $b=1+\frac{-ix}{n}$ festzulegen und  $n$ gegen $\infty$ streben zu lassen, wonach vermöge der Euler'schen Formel $e^x=\cos x +i \sin x$ das Sinusprodukt (\ref{eq: sine-Product}) entspringt. Eulers Beweis lässt sich etwas ausführlicher wie folgt darstellen:
Man betrachte das Polynom

\begin{equation}
\label{eq: Def Q_n für Sin}
    Q_n(x):= \dfrac{\left(1+\frac{ix}{n}\right)^n-\left(1-\frac{ix}{n}\right)^n}{2ix}
\end{equation}
Mithilfe der Euler'schen Zerlegungsformel für die Differenz zweier $n$--ter Potenzen findet man:

\begin{equation*}
    Q_n(x)= \dfrac{1}{2ix} \prod_{k=1}^n \left(\left(1+\frac{ix}{n}\right)^2 -2\left(1+\frac{ix}{n}\right)\left(1-\frac{ix}{n}\right)\cos \zeta_k +\left(1-\frac{ix}{n}\right)^2\right)
\end{equation*}
oder nach Vereinfachen

\begin{equation*}
    Q_n(x)= \dfrac{1}{2ix} \prod_{k=1}^n \left( 2\left(1-\frac{x^2}{n^2}\right)-2\left(1+\frac{x^2}{n^2}\right) \cos \zeta_k\right),
\end{equation*}
wobei $\zeta_k$ hier $=\frac{2k\pi}{n}$ ist. Da für unendliches $n$ das Argument des Kosinus gegen Null strebt, setzt Euler die Näherung

\begin{equation*}
    \cos \zeta_k \approx 1- \frac{2k^2\pi^2}{n^2}
\end{equation*}
in der Darstellung von $Q_n$ ein. Für, wie Euler sagt, unendlich großes $n$ kann demnach das Produkt als eines folgender Gestalt präsentiert werden:

\begin{equation*}
    Q_n(x)=  \dfrac{C_n}{2ix}\prod_{k=1}^n  \left(1- \frac{x^2}{k^2\pi^2}\right),
\end{equation*}
wobei $C_n$ eine noch zu bestimmende Konstante ist\footnote{Euler formuliert das nicht mithilfe einer Konstanten, sondern sagt, was auf dasselbe zurückführt, dass der allgemeine Faktor des Produktes im Fall von unendlichem $n$ die Form $1-\frac{ss}{kk\pi \pi}$ hat.}. Weiter ist bekannt, dass die rechte Seite von (\ref{eq: Def Q_n für Sin}) (vermöge der Euler'schen Identität) gerade $\frac{\sin x}{x}$ ist. Die Konstante $C_n$ kann man anschließend aus dem Fall $x=0$ zu $C_n=2i$ bestimmen, sodass man das Sinus--Produkt (\ref{eq: sine-Product}) gewinnt. Auch dieser Beweis erfüllt noch nicht die modernen Ansprüche an den mathematischen Rigor; ein solider Nachweis entlang der Linien des eben skizzierten bedarf des Begriffs der normalen Konvergenz, welchen erst Weierstraß (1815--1897) später in die Analysis eingebracht hat\footnote{Genau genommen wird der Begriff ``normale Konvergenz"{} erst von Baire (1874--1932) in seinem Buch \textit{``Le\c{c}ons sur les théories générales de l'analyse"} (\cite{Ba08}, 1908) (``Verlesungen zur allgemeinen Theorie der Analysis") verwendet, dort aber Weierstraß zugeschrieben. Baire referiert auf das, was heute gleichmäßige absolute Konvergenz genannt wird und der Inhalt des Weierstraß'schen Majorantentest ist.}. Einen solchen ausführlichen Beweis nach dem Euler'schen Vorbild mit ausführlichen Erläuterungen findet man zum Beispiel im Buch \textit{``Euler Through Time: A New Look at Old Themes"}  (\cite{Va06}, 2006).\\

\begin{figure}
    \centering
   \includegraphics[scale=0.6]{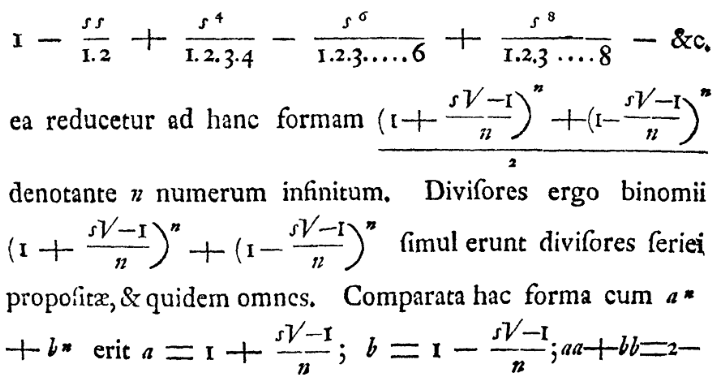}
    \caption{Der Beginn von Eulers Herleitung des Sinusproduktes aus der Zerlegung des Binoms $a^n-b^n$ in lauter trinomische Faktoren in seiner Arbeit \cite{E61}.}
    \label{fig:E61Faktor}
\end{figure}

Jedoch belässt es Euler  in \cite{E61} nicht mit einem Beweis, sondern enthüllt mit einer weiteren nicht über (\ref{eq: sine-Product}) verlaufenden Argumentation zum Nachweis von (\ref{eq: Baseler Problem}) seine Motivation des Wortes ``altera"{} im Titel seiner Arbeit. Sein neuer Ansatz fußt auf der Identität

\begin{equation*}
    \int\limits_{0}^1 \dfrac{x^{p-1}-x^{q-p-1}}{1-x^q}dx= \dfrac{\pi \cos \frac{p \pi}{q}}{q\sin \frac{p \pi}{q}},
\end{equation*}
welche für all jene Zahlen $p$ und $q$ gilt, für welche das Integral konvergiert. Entnommen hat Euler diese Formel  aus seiner Abhandlung \textit{``De inventione integralium, si post integrationem variabili quantitati determinatus valor tribuatur"} (\cite{E60}, 1743, ges. 1742) (E60: ``Über das Finden von Integralen, wenn nach der Integration der variablen Größe ein bestimmter Wert zugeteilt wird"), wo er eigens solche Integrale untersucht hat. Aus seinen Forschungen ergibt sich die Gültigkeit des obigen Ausdrucks zunächst nur für entsprechende natürliche Zahlen $p$ und $q$, Euler benötigt allerdings den Ausdruck auch für beliebige reelle Werte.  Er greift daher in § 16 von \cite{E61} auf das \textit{Kontinuitätsprinzip}\footnote{Euler beruft sich hierbei vermutlich auf Leibniz (1646--1716), welcher in seiner Schrift \textit{``Nouveaux Essais sur l'entendement humain"} (\cite{Le65}, 1765, ges. 1704) (``Neue Essays über den menschlichen Verstand") den Ausspruch \textit{``Natura non saltum facit"} -- Die Natur macht keinen Sprung -- verwendet.} zurück, welchem ihm die benötigte Inferenz gestattet\footnote{Da der Begriff der Stetigkeit noch nicht streng eingeführt worden ist, kann Euler natürlich nicht das gebräuchliche Argument Cauchys verwenden, mit dessen Hilfe man die den Gültigkeitsbereich zunächst auf die ganzen, dann die rationalen und schließlich auf die reellen Zahlen ausdehnen kann.}. Er findet somit die für reelle $s$ gültige Formel:

\begin{equation}
    \label{eq: Partialbruchzerlegung cot}
    \cot (\pi s) = \dfrac{1}{s}-\dfrac{1}{1-s}+\dfrac{1}{1+s}-\dfrac{1}{2-s}+\dfrac{1}{2+s}-\cdots,
\end{equation}
welche natürlich  auch unmittelbar durch logarithmisches Differenzieren aus dem Sinus--Produkt (\ref{eq: sine-Product}) folgt, welches umgekehrt aus der Partialbruchzerlegung des $\cot(x)$ entspringt. Damit schließt Euler also die Lücke seines Arguments aus \cite{E41}  und bestätigt überdies ein Charakteristikum seiner Arbeitsweise erneut: Die mehrfache Bestätigung eines Theorems auf verschiedenen Wegen.

\subsubsection{Ein moderner Beweis}
\label{subsubsec: Ein moderner Beweis}

\epigraph{Inspirational thunderbolts do not appear out of the blue. [...] Coming up with  a grand idea without ever having been closely involved with an area of study is not impossible, but is very improbable.}{Ulrich Kraft}

Die letzte Aussage, Euler beweise ein und denselben Lehrsatz auf mehrere Weise, findet auch Bestätigung durch  einen weiteren Beweis von (\ref{eq: Baseler Problem}) in \cite{E63}, der sich nochmal von den beiden vorherigen Methoden unterscheidet und außerdem  nicht heute obsolete Begriffe wie das Kontinuitätsprinzip verwendet. In der Tat ist er auch nach modernen Ansprüchen als vollkommen streng anzusehen und kann s  wörtlich ohne Einwände übernommen werden. \\

Für diesen Beweis bemerkt Euler zunächst, dass sich der Wert der gesuchten Reihe aus
\begin{equation*}
    \int\limits_{0}^{1} \dfrac{dx}{\sqrt{1-x^2}}\int\limits_{0}^{x}\dfrac{dy}{\sqrt{1-y^2}}
\end{equation*}
ergibt. Man bedarf überdies der Reduktionsformel

\begin{equation*}
  \int\limits_{0}^{1} \dfrac{x^{n+2}dx}{\sqrt{1-x^2}}= \dfrac{n+2}{n+1} \int\limits_{0}^{1} \dfrac{x^ndx}{\sqrt{1-x^2}},
\end{equation*}
welche Euler mit partieller Integration nachweist. Zudem ergibt sich unter Anwendung des binomischen Lehrsatzes\footnote{Ein strenger Nachweis des binomischen Lehrsatzes für nicht natürliche Zahlen existierte zum Zeitpunkt der Verfassung der Arbeit \cite{E63} noch nicht. Einen solchen hat   Euler selbst erst in der Abhandlung \textit{``Demonstratio theorematis Neutoniani de evolutione potestatum binomii pro casibus, quibus exponentes non sunt numeri integri"} (\cite{E465}, 1775, ges. 1773) (E465:``Beweis des Newton'schen Theorems über die Entwicklung der Potenzen des Binoms für die Fälle, in denen die Exponenten keine ganzen Zahlen sind") gegeben.} auf den Integranden und anschließender gliedweiser Integration

\begin{equation*}
    \int\limits_{0}^{x}\dfrac{dy}{\sqrt{1-y^2}}= x+\dfrac{1}{2\cdot 3}x^3+\dfrac{1\cdot 3}{2\cdot 4 \cdot 5}x^5+\dfrac{1\cdot 3 \cdot 5}{2\cdot 4 \cdot 6 \cdot 7}x^7+\cdots.
\end{equation*}
sowie

\begin{equation*}
     \int\limits_{0}^{1} \dfrac{dx}{\sqrt{1-x^2}}\int\limits_{0}^{x}\dfrac{dy}{\sqrt{1-y^2}}
\end{equation*}
\begin{equation*}
    =\int\limits_{0}^{1} \dfrac{xdx}{\sqrt{1-x^2}}+\dfrac{1}{2\cdot 3}\int\limits_{0}^{1} \dfrac{x^3dx}{\sqrt{1-x^2}}+\dfrac{1\cdot 3}{2\cdot 4 \cdot 5}\int\limits_{0}^{1} \dfrac{x^5dx}{\sqrt{1-x^2}}+\cdots
\end{equation*}
Mit der Reduktion vereinfacht sich das aber gerade zu:

\begin{equation*}
    \int\limits_{0}^{1} \dfrac{dx}{\sqrt{1-x^2}}+\dfrac{1}{3\cdot 3}\int\limits_{0}^{1} \dfrac{dx}{\sqrt{1-x^2}}+\dfrac{1}{5\cdot 5}\int\limits_{0}^{1} \dfrac{dx}{\sqrt{1-x^2}}+\text{etc.}
\end{equation*}
Schnell berechnet man

\begin{equation*}
    \int\limits_{0}^{1} \dfrac{dx}{\sqrt{1-x^2}} =1,
\end{equation*}
wonach insgesamt gilt

\begin{equation*}
    \int\limits_{0}^{1} \dfrac{dx}{\sqrt{1-x^2}}\int\limits_{0}^{x}\dfrac{dy}{\sqrt{1-y^2}}= 1+\dfrac{1}{3^2}+\dfrac{1}{5^2}+\dfrac{1}{7^2}+\text{etc.}
\end{equation*}
Jedoch hat man allgemein

\begin{equation*}
    \int f(x)\cdot f'(x)dx = \dfrac{f^2(x)}{2}+C,
\end{equation*}
somit wegen $(\arcsin(x))'=\frac{1}{\sqrt{1-x^2}}$ an die Fragestellung angepasst

\begin{equation*}
    \int\limits_{0}^{1} \dfrac{dx}{\sqrt{1-x^2}}\int\limits_{0}^{x}\dfrac{dy}{\sqrt{1-y^2}}=\dfrac{1}{2}\left((\arcsin(1))^2-(\arcsin(0))^2\right)= \dfrac{\pi^2}{8}.
\end{equation*}
Bringt man alles zusammen, lautet das Endergebnis:

\begin{equation*}
    \dfrac{\pi^2}{8}=1+\dfrac{1}{3\cdot 3}+\dfrac{1}{5\cdot 5}+\cdots
\end{equation*}
Die Bemerkung

\begin{equation*}
    1+\dfrac{1}{3^2}+\dfrac{1}{5^2}+\cdots = 1+\dfrac{1}{2^2}+\dfrac{1}{3^2}+\cdots -\left(\dfrac{1}{2^2\cdot 1^2}+\dfrac{1}{2^2\cdot 3^2}+\dfrac{1}{2^2\cdot 5^2}+\cdots\right)
\end{equation*}
\begin{equation*}
    =\dfrac{3}{4}\left(1+\dfrac{1}{2^2}+\dfrac{1}{3^2}+\cdots\right).
\end{equation*}
lässt Euler schließlich zur gewünschten Formel gelangen:

\begin{equation*}
   \dfrac{4}{3} \cdot \dfrac{\pi^2}{8}=\dfrac{\pi^2}{6}=1+\dfrac{1}{2^2}+\dfrac{1}{3^2}+\cdots,
\end{equation*}
was  wieder (\ref{eq: Baseler Problem}) ist.\\

Neben diesem Beweis gibt Euler in selbigem Aufsatz \cite{E63} noch einen weiteren Beweis im selben Gewand. Statt $\arcsin(x)$ gebraucht er die Potenzreihe für $\arcsin^2(x)$, welche Euler in Analogie zum vorgestellten  Beweis mithilfe des Ansatzes von unbestimmten Koeffizienten aus der von $\arcsin^2(x)$ erfüllten Differentialgleichung ableitet. Höhere Potenzsummen kann Euler, wie er auch selbst einräumt, indes mit den in \cite{E63} vorgestellten Methoden nicht behandeln.

\begin{figure}
    \centering
     \includegraphics[scale=1.1]{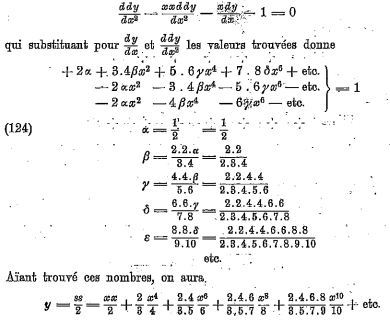}
    \caption{Euler leitet in seiner Arbeit \cite{E63} aus einem Potenzreihenansatz die Potenzreihe für $\arcsin^2(x)$ her.}
    \label{fig:E63arcsin}
\end{figure}

\subsubsection{Zusammenfassung zum Baselproblem}

\epigraph{There is eternal simplicity to a solution once it has been discovered!}{Aleksandr Solzhenitsyn}
\label{subsubsec: Zusammenfassung zum Baselproblem}

Probleme haben Euler selten losgelassen, bevor er sie nicht selbst gelöst hatte. Beim Baseler Problem umfasst der Zeitraum mindestens 10 Jahre, wobei er seine Teilfortschritte wie die Reihentransformation aus \cite{E20} ebenfalls mit der Öffentlichkeit geteilt hat\footnote{Dies bedeutet im Zusammenhang des Baselerproblems wie gesehen etwa die Transformation der Reihe in eine schneller konvergente in \cite{E20}, obwohl auch diese ihn dem eigentlichen Problem nicht näher zu bringen scheint. Dieser Schein trügt aber, da der durch die Reihe ermittelte Näherungswert Euler wohl auf die Verbindung zur Kreiszahl gebracht hat. Bekanntermaßen ist es unlängst leichter, sofern man bereits in Kenntnis dessen ist, was man zu zeigen versucht. Dies ist auch ein in der  Psychologie   -- siehe etwa die Arbeit \footnote{``Analysis of minimal complex systems and complex problem solving require different forms of causal cognition"} (\cite{Fu14}, 2014) -- bekannter Sachverhalt. Für die definitorischen Grundlagen in diesem Bereich konsultiere man etwa das Buch \textit{``Thinking, Problem Solving, Cognition"} (\cite{Ma92}, 1992) und gängige Einteilung von Problem in einfache und komplexe Probleme, wie sie unter Anderem im Buch \textit{``Complex Problem Solving: The Europeen Perspective"} (\cite{Fr14}, 2014) erläutert werden.  Man findet überdies in Fellmanns Buch \cite{Fe95}  Daniel Bernoulli mit dem Ausspruch zitiert, dass er unter Kenntnis des exakten Wertes auch zu einer Lösung des Baseler Problems gelangt wäre.}. Zu dem Beweis in \cite{E41}, welcher zu diesem Zeitpunkt noch nicht als einer gelten konnte, was Euler auch bewusst war, ist hervorzuheben, dass Euler ihn publiziert hat, \textit{obschon} er um seine Unvollständigkeit wusste. Er resümiert in \cite{E63} über seine Methode in \cite{E41}:\\

\textit{``Sie ist also ebenso sicher und fundiert wie jede andere Methode, die man gewöhnlich bei der Summierung von unendlichen Reihen angewendet: Dies habe ich durch die perfekte Übereinstimmung einiger bereits bekannter Fälle und durch die Annäherungen gezeigt, die uns eine einfache Möglichkeit bieten, die Gültigkeit in der Praxis zu prüfen."}\\

Die oben (Abschnitt \ref{subsubsec: Eulers Auffassung eines Beweises}) aufgestellte Behauptung, Euler vertrete die Meinung, Mathematik werde entdeckt und nicht geschaffen, findet demnach eine weitere Bestätigung. Gleiches gilt für die enge Verzahnung von Praxis und Theorie in seinem Schaffen.  Diese Einstellung bezüglich Publikationen auch unvollständig behandelter Gegenstände unterscheidet ihn insbesondere von Gauß, welcher nichts vor seiner Vollendung publiziert hat\footnote{Der Ausspruch \textit{``Pauca, sed matura"} (``Weniges, aber dafür Reifes") zierte überdies seinen Petschaft.}. Die Veröffentlichung des unfertigen Arguments aus Gründen der Sicherung der Priorität darf als unwahrscheinlich angesehen werden, zumal Euler bei anderen Gelegenheiten seine Entdeckungen entweder zurückgehalten hat\footnote{Gegenüber seinen Ideen zur Hydromechanik gab er Daniel Bernoullis Buch {``Hydrodynamica -- Sive de viribus et motibus fluidorum commentarii} (\cite{Be38}, 1738) (``Hydrodynamik -- Oder Kommentare über die Kräfte und Bewegungen von Fluiden") den Vorrang.} -- oder Ideen trotz Prioritätsansprüchen nicht geltend gemacht hat, wie etwa bei der Summenformel von Euler und Maclaurin (siehe den entsprechenden Abschnitt in \cite{Va06}.). Zusätzlich hat Euler die Gewohnheit der Kommunikation von Teilfortschritten im Verlauf seiner Karriere nicht abgelegt. Außerdem ist die Aussage, dass er etwas schon jetzt für ein sicheres Theorem halte, es aber noch nicht beweisen kann,  typisch für Euler, wovon oben (Abschnitt \ref{subsubsec: Eulers Auffassung eines Beweises}) ein anderes Exempel gegeben worden ist. \\

Es soll aber noch einmal kurz auf die Argumentationsstruktur im obigen Zitat eingegangen werden.  So zeigt die Arbeit \cite{E61}  die Anwendung des Prinzip der konvergenten Evidenz zur Erhärtung von Thesen. Konkret nimmt Euler Bezug auf die Reihen

\begin{equation*}
    \dfrac{\pi}{4}=1-\dfrac{1}{3}+\dfrac{1}{5}-\text{etc.},
\end{equation*}
welche Leibniz bereits vorher anders gezeigt hatte, was Euler in § 10 auch erwähnt, sowie

\begin{equation*}
    \dfrac{\pi}{2\sqrt{2}}=1+\dfrac{1}{3}-\dfrac{1}{5}-\dfrac{1}{7}+\text{etc.},
\end{equation*}
welche er aus dem Studium von Newtons Werken kennt. Diese Reihe gibt Euler in § 14 an. All dies muss ihn hinreichend überzeugt haben, dass insbesondere (\ref{eq: sine-Product}) trotz der mehrfach erwähnten lückenhaften Beweisführung richtig ist. \\

Nach diesen Beweisen sollte Euler im Verlauf seiner Karriere noch mehrmals auf die allgemeinen Summen
\begin{equation*}
    \sum_{k=1}^{\infty} \dfrac{1}{k^{2n}}.
\end{equation*}
zurückkommen. In der Tat formuliert Euler in § 24 von \cite{E130} einen allgemeinen Ausdruck, welche sich in moderner Weise so darstellt

\begin{equation}
      \sum_{k=1}^{\infty} \dfrac{1}{k^{2n}} = \dfrac{(-1)^n\cdot (2\pi)^{2n}B_{2n}!}{2(2n)!},
\end{equation}
wobei $B_n$ die Bernoulli-Zahlen sind, definiert über die Potenzreihenrentwicklung

\begin{equation*}
    \dfrac{x}{e^x-1}= \sum_{n=0}^{\infty} \dfrac{B_nx^n}{n!}.
\end{equation*}
Diese Definition geht ebenfalls auf Euler zurück, man entnimmt sie zwar bereits aus seiner Abhandlung \textit{``Methodus generalis summandi progressiones"} (\cite{E25}, 1738, ges. 1732) (E25: ``Eine allgemeine Methode Progressionen zu summieren") in § 2, deutlicher erkennbar ist sie allerdings der Arbeit \textit{``Inventio summae cuiusque seriei ex dato termino generali"} (\cite{E47}, 1741, ges. 1735) (E47: ``Das Finden einer Summe aus dem gegebenen allgemeinen Term"). Dort findet man sie in § 15.\\

 Die Lösung des Baseler Problems markiert zweifelsohne eine von Eulers größten Erfolgen. Sein Stolz über selbige  ist ihm nicht zu verdenken, welcher sich auch darin ausdrückt, dass  (\ref{eq: Baseler Problem}) eine der wenigen Wahrheiten ist, bei welchen Euler in späteren Arbeiten explizit seine Erstentdeckerschaft erwähnt\footnote{Man vergleiche etwa  § 60 von \textit{``De serierum determinatione seu nova methodus inveniendi terminos generales serierum"} (\cite{E189}, 1753, ges. 1749) (E189: ``Über die Bestimmung von Reihen oder eine neue Methode die allgemeinen Termen von Reihen zu finden"), das Ende von § 9 seiner Arbeit \textit{``De valore formulae integralis $\int \frac{z^{\lambda-\omega}\pm z^{\lambda +\omega}}{1\pm z^{2\lambda}}\frac{dz}{z}(\log z)^{\mu}$ im Fall, in welchem nach der Integration $z = 1$ gesetzt wird"} (\cite{E463}, 1775, ges. 1774) (E463: ``Über den Wert der Integralformel $\int \frac{z^{\lambda-\omega}\pm z^{\lambda +\omega}}{1\pm z^{2\lambda}}\dfrac{dz}{z}(\log z)^{\mu}$ im Fall, in welchem nach der Integration $z=1$ gesetzt wird"), das Ende von § 74 der Abhandlung \textit{``De formulis integralibus implicatis earumque evolutione et transformatione"} (\cite{E679}, 1794, ges. 1778) (E679: ``Über verschachtelte Integrale und deren Entwicklung sowie Transformation") sowie § 1 von \textit{``De summatione serierum in hac forma contentarum $\frac{a}{1}+\frac{a^2}{4}+\frac{a^3}{9}+\frac{a^4}{16}+\frac{a^5}{25}+\frac{a^6}{36}+\text{etc.}$"} (\cite{E736}, 1811, ges. 1779) (E736: ``Über die Summe der in der Form $\frac{a}{1}+\frac{a^2}{4}+\frac{a^3}{9}+\frac{a^4}{16}+\frac{a^5}{25}+\frac{a^6}{36}+\text{etc.}$ enthaltenen Reihen"), in welcher Euler die heute als Dilogarithmus bezeichnete Funktion einführt und untersucht.}. Die Anzahl seiner Beweise für diese Summe zeigt zugleich auch sein riesiges Ausmaß an Kreativität. Dies betrifft nicht nur die reine Anzahl an verschiedenen Beweisen -- insgesamt hat Euler in seiner Karriere mindestens 6 in ihrer Natur verschiedene Beweise von (\ref{eq: Baseler Problem}) vorgelegt\footnote{Den ersten Beweis findet man in \cite{E41}, sofern man das Sinusprodukt (\ref{eq: sine-Product}) zulässt, einen in \cite{E61}, zwei in \cite{E63}, einen weiteren in \textit{``De eximio usu methodi interpolationum in serierum doctrina"} (\cite{E555}, 1783, ges. 1772) (E555: ``Über den außerordentlichen Nutzen von Interpolationen in der Reihenlehre") über einen Spezialfall der Lagrange-Interpolation, und einen in \textit{``De resolutione fractionum transcendentium in infinitas fractiones simplices"} (\cite{E592}, 1785, ges. 1775) (E592:``Über die Auflösung von transzendenten Brüchen in unendlich viele einfache Brüche") mittels Partialbruchzerlegung von transzendenten Funktionen.} -- sondern auch die Natur der Beweise selbst. Der Beweis über (\ref{eq: Baseler Problem}) zeigt die Unerschrockenheit beim Ausdehnen von Konzepten über die eigentlichen vom Standpunkt der Logik erlaubten Grenzen hinaus, manifestiert in der Verwendung für den endlichen Fall gültigen Formeln im den Fall des Unendlichen. Insbesondere der Beweis aus \cite{E63} zeigt, dass Euler ebenfalls sehr gut scheinbar lose Dinge verknüpfen kann, um einen Beweis zu konstruieren. Er verlangt nur elementare Konzepte, welche richtig zusammengebracht, zu der wohl elementarsten Lösungen des Baseler Problems führen, welche Euler jemals angegeben hat.

\subsection{Ein richtiges Ergebnis mit unzulässiger Herleitung}
\label{subsec: Ein falsches Ergebnis}

\epigraph{All you need is ignorance and confidence and the success is sure.}{Mark Twain}

Anhand seiner Resultate zu Differentiagleichungen unendlicher Ordnung mit konstanten Koeffizienten wird ein Beispiel vorgestellt, in welchem Euler mit einer unzulässigen Herleitung dennoch ein richtiges Ergebnis deriviert. Im Kontrast zu im in \cite{E41} (Abschnitt \ref{subsubsec: Ermittlung des exakten Wertes}) gegebenen Argument, sollte es Euler nicht gelingen, einen strengen Beweis nachzureichen. Die Präsentation folgt der Chronologie der Euler'schen Publikationen: Erst werden die Beiträge zu homogenen Differentiagleichungen (Abschnitt \ref{subsubsec: Differentialgleichungen unendlicher Ordnung: Der homogene Fall}) betrachtet, gefolgt von den inhomogenen (Abschnitt (\ref{subsubsec: Der inhomogene Fall}). 

\subsubsection{Differentialgleichungen unendlicher Ordnung: Der homogene Fall}
\label{subsubsec: Differentialgleichungen unendlicher Ordnung: Der homogene Fall}

\epigraph{[W]ho takes a chance, who walks the line between the known and unknown, who is unafraid of failure, will succeed.}{Gordon Parks}

 Mit glücklichem Erfolg hatte Euler die allgemeine Differentialgleichung
    \begin{equation}
    \label{eq: homogeneous Finite}
    \left(a_0+a_1 \dfrac{d}{dx}+a_2 \dfrac{d^2}{dx^2} +\cdots+ a_n \dfrac{d^n}{dx^n} \right)f(x)= 0,
\end{equation}
mit komplexen Koeffizienten $a_1$, $a_2$, $\cdots$, $a_n$, welche Euler (in unnötiger Weise) durchgehend alle als reell annimmt\footnote{Die konsequente Vermeidung Eulers von komplexen Zahlen in der Lösung solcher und verwandter Fragen wird in Abschnitt \ref{subsubsec: Komplexe Analysis} noch einmal aufgenommen werden.}, behandelt und gelöst. Die Lösungsmethode, welche er in der Arbeit \textit{``De integratione aequationum differentialium altiorum graduum"} (\cite{E62}, 1743, ges. 1742) (E62: ``Über die Integration von Differentialgleichungen höherer Grade") vorstellt, stimmt immer noch mit der heutigen überein. Sie findet sich in §§ 16--27 besagter Arbeit, woran sich Anwendungen auf spezielle Fälle anschließen, geordnet gemäß der Natur der Nullstellen des heute so bezeichneten charakteristischen Polynoms: 

\begin{equation*}
    P(z)= a_0+a_1 z+a_2 z^2 +\cdots+ a_n z^n=0.
\end{equation*}
Selbiges erhält man, wie Euler auch beschreibt, durch Ersetzung von $\frac{d^n}{dx^n}$ mit $z^n$ in (\ref{eq: homogeneous Finite}). Von Interesse ist hier zunächst, dass Euler in seinen finalen Formeln keine komplexen Zahlen mehr auftreten lassen will\footnote{Die Annahme ausschließlicher reeller Koeffizienten von $P(z)$ impliziert das Auftreten von komplexen Lösungen in konjugierten Paaren, sodass sich je zwei Linearfaktoren $(z-\alpha)$ und $(z-\overline{\alpha})$ zu einem reellen quadratischen zusammenfassen lassen.}. Dies sei an Eulers Beispiel 4 aus § 32 illustriert. Hier möchte er folgende Gleichung lösen:

\begin{equation}
\label{eq: Euler Beispiel 4}
    0=y-a^4\dfrac{d^4y}{dx^4}.
\end{equation}
$a$ ist hierbei reell und man findet von $P(z)=1-a^4z^4$ ohne Mühe die Faktorisierung

\begin{equation*}
    P(z)=(1-az)(1+az)(1+a^2z^2).
\end{equation*}
Aus den ersten beiden Faktoren deriviert Euler unmittelbar beiden Lösungen

\begin{equation*}
    y_1(x)= \mathfrak{A}e^{\frac{x}{a}} \quad \text{und} \quad   y_2(x)= \mathfrak{B}e^{-\frac{x}{a}},
\end{equation*}
wobei $\mathfrak{A}$ und $\mathfrak{B}$ beliebige Integrationskonstanten sind. Für die restlichen Lösungen zerlegt Euler nun nicht, wie heute üblich, $1+a^2z^2$ weiter in $(1-iaz)(1+iaz)$ und  konstruiert die weiteren Lösungen

\begin{equation*}
    y_3(x)= \mathfrak{C}e^{-i\frac{x}{a}} \quad \text{und} \quad   y_4(x)= \mathfrak{D}e^{-i\frac{x}{a}},
\end{equation*}
wobei $\mathfrak{C}$ und $\mathfrak{D}$ wieder beliebige Integrationskonstanten sind, sondern vermeidet das Komplexe gänzlich und setzt an, dass sich für gewisses $f$ und $\varphi$ schreiben lässt:

\begin{equation*}
    \dfrac{1}{a^2}+z^2= f^2-2f z\cos \varphi + z^2.
\end{equation*}
Dies ist wegen der bekannter  Eigenschaften der komplexen Lösungen von reellen Polynomen stets möglich und Euler findet durch direkten Vergleich entsprechende Werte und daraus die Lösungsfunktion:

\begin{equation*}
    y_{E}(x)= \gamma \sin \left( \dfrac{x}{a}-\beta\right),
\end{equation*}
mit beliebigen Konstanten $\beta$ und $\gamma$. Dieses bewusste Umschiffen des Komplexen in der finalen Lösung ist typisch für Euler und zieht sich durch seinen gesamten Opus und wird unten (Abschnitt \ref{subsubsec: Komplexe Analysis}) Gestand der Ausführungen werden.\\

Nach diesen einleitenden Worten soll nun der Übergang zum Fall einer Differentialgleichung unendlicher Ordnung geschehen. Im letzten Problem dieser Arbeit \cite{E62} (Problem XI § 50) behandelt Euler eine nämliche, welche zugleich das erste Beispiel einer solchen Fragestellung überhaupt zu sein scheint.  Euler möchte hier folgende Gleichung lösen:

\begin{equation}
    \label{eq: Euler erste DGL unendlich}
    0= y -\dfrac{1}{2!}\dfrac{d^2y}{dx^2}+\dfrac{1}{4!}\dfrac{d^4y}{dx^4}-\dfrac{1}{6!}\dfrac{d^6y}{dx^6}+\text{etc.}.
\end{equation}
Ein wohlüberlegtes Beispiel, führt sie gemäß seiner Methode zum charakteristischen ``Polynom":

\begin{equation*}
    P(z)=1- \dfrac{z^2}{2!}+\dfrac{z^4}{4!}-\dfrac{z^6}{6!}+\text{etc.}, 
\end{equation*}
was gerade die Potenzreihe für $\cos z$ ist. Dessen unendlich viele Nullstellen sind -- wie Euler dann kurz zuvor in \cite{E61} bewiesen hatte --  gegeben als

\begin{equation*}
    \pm \dfrac{\pi}{2}, \pm \dfrac{3\pi}{2},  \pm \dfrac{5\pi}{2},  \pm \dfrac{7 \pi}{2}. \cdots
\end{equation*}
Dies führt somit zur Lösung:

\begin{equation}
\label{eq: Euler Diff Cos}
    y(x)= \alpha e^{\frac{\pi}{2}x}+a e^{-\frac{\pi}{2}x}+\beta e^{\frac{3\pi}{2}x}+b e^{-\frac{3\pi}{2}x}+\gamma e^{\frac{5\pi}{2}x}+c e^{-\frac{5\pi}{2}x}+\text{etc.}
\end{equation}
mit beliebigen Konstanten $\alpha$, $a$ etc. Diese Lösung erfüllt per Konstruktion  (\ref{eq: Euler erste DGL unendlich}). Jedoch ist an dieser Stelle noch nicht sichergestellt, dass dies auch die \textit{vollständige} Lösung darstellt. Im Fall endlicher Ordnung gibt die Euler'sche Vorgehensweise entsprechend eine Lösung mit $n$ frei wählbaren Konstanten, woraus  -- wie Euler auch bemerkt -- sich die Vollständigkeit der Lösung ergibt. Im Fall unendlicher Ordnung ist dies weniger offenkundig, beinhaltet etwa

\begin{equation*}
   y_{\text{teil}}(x)= \alpha e^{\frac{\pi}{2}x} +\beta \alpha e^{\frac{3\pi}{2}x} +\gamma e^{\frac{5\pi}{2}x}+\text{etc.}
\end{equation*}
ebenfalls unendlich viele frei wählbare Konstanten, kann jedoch als Teillösung der zuvor hergeleiteten (\ref{eq: Euler Diff Cos}) nicht die vollständige sein. Gleichermaßen könnte man einwenden, dass auch die von Euler gefundene Lösung nur eine Teillösung einer noch allgemeineren ist.  Die Kernfrage, ob es auch im Unendlichen zulässig ist\footnote{Auf die folgenden beiden Abschnitte vorgreifend, sei angemerkt, dass dies im Allgemeinen nicht richtig ist, was auch anhand Beispielen illustriert werden wird.}, die Lösung allein aus den Nullstellen des charakteristischen ``Polynoms"{} zu konstruieren, lässt Euler somit unbehandelt zurück.\\

Ergänzend und in Vorbereitung auf den kommenden Abschnitt sei noch bemerkt, dass sich (\ref{eq: Euler erste DGL unendlich}) als Differenzengleichung auffassen lässt:

\begin{equation*}
    0= y(x+i)+y(x-i),
\end{equation*}
wobei Euler selbst wohl $\sqrt{-1}$ statt $i$ geschrieben hätte. Schnell wird verifiziert, dass Eulers Lösung (\ref{eq: Euler Diff Cos}) in dieser enthalten ist, überdies ist

\begin{equation*}
    y(x)= A (i)^{ix}+B(-i)^{ix}
\end{equation*}
mit beliebigen $A$ und $B$ die allgemeinste Lösung der Differenzengleichung. Vermöge der Euler'schen Formel

\begin{equation*}
    e^{ix}= \cos(x)+i \sin(x)
\end{equation*}
lässt sich die Euler'sche unendliche viele Terme umfassende Lösung schnell begreiflich machen. Dazu hat man

\begin{equation*}
    i = e^{i\frac{\pi}{2}} \quad \text{oder} \quad i = e^{5i\frac{\pi}{2}} \quad \text{oder auch} \quad i = e^{9i\frac{\pi}{2}} \quad \text{etc.}
\end{equation*}
zu schreiben. Für $(-i)$ gehe man analog vor. Dann gibt eine unmittelbare Anwendung der Potenzgesetze:

\begin{equation*}
    i^{ix}= (e^{(4k+1)\frac{i\pi}{2}})^{ix}= e^{-(4k+1)\frac{\pi}{2}x},
\end{equation*}
analog 

\begin{equation*}
    (-i)^{ix}= (e^{(4l+3)\frac{i\pi}{2}})^{ix}= e^{-(4l+3)\frac{\pi}{2}x}
\end{equation*}
mit ganzen Zahlen $l$ und $k$. Dieser Interpretation einer Differentialgleichung unendlicher Ordnung als Differenzengleichung hat sich Euler in der Arbeit  \textit{``De serierum determinatione seu nova methodus inveniendi terminos generales serierum"} (\cite{E189}, 1753, ges. 1749) (E189: ``Über die Bestimmung von Reihen oder eine neue Methode die allgemeinen Termen von Reihen zu finden") in umgekehrter Weise bedient und dies bildet überdies den Inhalt von Abschnitt (\ref{subsec: Ein Prototyp -- Differenzengleichungen}). 

\subsubsection{Der inhomogene Fall}
\label{subsubsec: Der inhomogene Fall}

\epigraph{The only real mistake is the one from which we learn nothing.}{Henry Ford}

An der Behandlung homogener Differentialgleichungen unendlicher Ordnung kommt wie beim Baseler Problem (Abschnitt \ref{subsubsec: Ermittlung des exakten Wertes}) die heutzutage unbedarft wirkende Haltung gegenüber dem Unendlichen zum Vorschein, welches für Euler bei dieser Begebenheit just als weitere Zahl gesehen wird. Während sich an dieser Stelle noch berechtigt vorbringen ließe, die Unvollständigkeit seiner Herangehensweise wäre aufgrund der Richtigkeit der Lösung und konvergierender Evidenz für Euler schwierig zu erkennen, kann dieser Einwand bei seiner Behandlung des inhomogenen Falles nicht mehr geltenden gemacht werden. Hier, so wird im Verlaufe gezeigt werden, gelangt er zu einer unrichtigen Lösung, welche er außerdem als solche hätte erkennen können.\\

Um den genauen Ursprung von Eulers Fehler zu eruieren, ist es hilfreich, zunächst seine Ausführungen zum Gegenstand der inhomogenen Differentialgleichungen mit konstanten Koeffizienten aus der Arbeit \textit{``Methodus aequationes differentiales altiorum graduum integrandi ulterius promota"} (\cite{E188}, 1753, ges. 1750) (E188: ``Die Methode Differentialgleichungen höherer Grade zu integrieren weiter entwickelt") nachzuverfolgen. Diese beginnen ab § 7 und sind in diesem Sinne typisch Euler'sch, dass er sich von den einfachsten Fällen zum allgemeinen vorarbeitet. Er beginnt mit dem Fall erster Ordnung

\begin{equation}
\label{eq: Inhomogen n=1}
    X(x) = Ay(x)+By'(x).
\end{equation}
Euler bevorzugt in seinen Arbeiten zuhauf die Darstellung mit Differentialen und schreibt letztgenannte Gleichung als:

\begin{equation*}
    X dx =A y dx +Bdy
\end{equation*}
und multipliziert dies mit $e^{\alpha x}$:

\begin{equation*}
    e^{\alpha x} X dx = A e^{\alpha x} y dx + e^{\alpha x}B dy.
\end{equation*}
$\alpha$ ist nun so zu wählen, dass die rechte Seite wie die linke ein vollständiges Differential ist\footnote{Die Idee der Multiplikation mit einer Exponentialfunktion motiviert Euler hier nicht, allerdings hat er  solche Integrale in seiner Arbeit \textit{``De infinitis curvis eiusdem generis seu methodus inveniendi aequationes pro infinitis curvis eiusdem generis"} (\cite{E44}, 1740, ges. 1734) (E44: ``Über unendlich viele Kurven derselben Art oder eine Methode die Gleichungen für unendlich viele Kurven derselben Art zu finden") betrachtet und ist sich seitdem ihrer Eigenschaften bewusst gewesen.}.  Nach Euler ist dies notwendigerweise als $d(Be^{\alpha x}y)$ gegeben; damit hat man

\begin{equation*}
     d(Be^{\alpha x}y)= B\alpha e^{\alpha x}ydx +B e^{\alpha x}dy.
\end{equation*}
Dies führt  zur Gleichung

\begin{equation*}
   B\alpha e^{\alpha x}ydx +B e^{\alpha x}dy =  A e^{\alpha x} y dx + e^{\alpha x}B dy,
\end{equation*}
welche notwendig gelten muss. Ein Vergleich der entsprechenden Terme impliziert: $A=\alpha B$ oder $\alpha=\frac{A}{B}$, woraus mit

\begin{equation*}
   \int e^{\alpha x} X dx =  \int  d(Be^{\alpha x}y) = Be^{\alpha x}y
\end{equation*}
nach Einsetzen des Wertes für $\alpha$ zur Lösung

\begin{equation*}
    y= \dfrac{\alpha}{A} e^{-\alpha x}\int e^{\alpha x}Xdx.
\end{equation*}
von (\ref{eq: Inhomogen n=1}) gelangt wird. Die beliebige Integrationskonstante geht über das Integral ein und weist die Lösung als die vollständige aus\footnote{Eulers Ansatz, das gesuchte Differential $d(Be^{\alpha x}y)$ zu betrachten, wird nachvollziehbar, wenn ein Vergleich mit der heute üblichen Methode der Variation der Konstanten zur Lösung von (\ref{eq: Inhomogen n=1}) angestellt wird. Selbiger offenbart, dass Euler in seinen Erklärungen die beiden Schritte dieser Methode (erst die Lösung der homogenen Gleichung zur Bestimmung von $\alpha$ und die anschließende Variation der Konstanten zum Finden des vollständigen Differentials $ d(Be^{\alpha x}y) $) in seinem Ansatz  gleichzeitig vollzieht.}.\\

Im nächsten Schritt seiner allgemeinen Lösung von (\ref{eq: Inhomogeneous Finite}) betrachtet Euler die Gleichung von zweiter Ordnung, schreitet also im Grad einen Schritt voran. Er setzt an:

\begin{equation*}
    X=A +B\dfrac{dy}{dx}+C\dfrac{dy^2}{dx^2}. 
\end{equation*}
Multipliziert mit $e^{\alpha x}dx$ führt dies zu

\begin{equation*}
    e^{\alpha x}X dx= e^{\alpha x}Adx+ Be^{\alpha x} dy+ C e^{\alpha x} \dfrac{dy}{dx}, 
\end{equation*}
was er nun auf den obigen Fall zurückführt, indem er

\begin{equation*}
    \int   e^{\alpha x}X dx = e^{\alpha x}\left(A'y+B'\dfrac{dy}{dx}\right)
\end{equation*}
mit zu bestimmenden Größen $A'$, $B'$ festlegt. Durch Differenzieren dieser Gleichung und einen Vergleich mit der vorgelegten gelingt es ihm, $A'$ und $B'$ über $A$ und $B$ sowie $\alpha$ auszudrücken, sodass er im letzten Schritt erneut eine Gleichung für $\alpha$ ableiten kann; in diesem Fall lautet diese

\begin{equation*}
    0=A-B\alpha+ C \alpha^2.
\end{equation*}
Da die vorherige Integralgleichung  dem vorherigen Fall erster Ordnung entspricht, ist
das vorgelegte Problem mit dem Auffinden von $\alpha$ abgehandelt. \\

Aus diesen Darstellungen scheint Euler das Muster offenkundig geworden zu sein: Der Fall $n-$ter Ordnung kann stets auf den $n-1$-ter Ordnung zurückgeführt werden, welche Einsicht er in § 9 -- 13 mit einem Beweis unterlegt. Da er auf diesem Wege zu einer Lösungsformel gelangt, die $n$ Integrationskonstanten involviert, ist die Auflösung der allgemeinen inhomogenen Differentialgleichung  von $n$--ter Ordnung grundlegend abgehandelt; ergänzt wird Eulers Ausführung noch durch eine Diskussion der Sonderfälle mehrfacher und komplexer Nullstellen und die Angabe expliziter Formeln.\\

\begin{figure}
    \centering
   \includegraphics[scale=1.0]{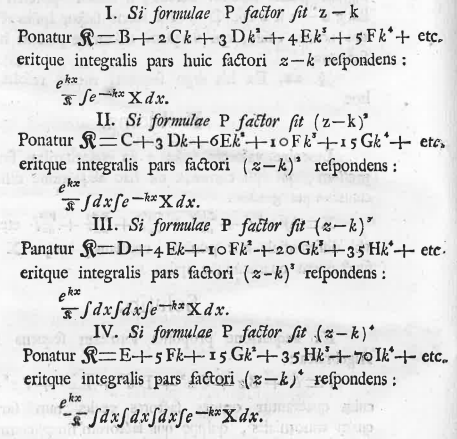}
    \caption{Eulers Auflistung der Lösung der allgemeinen inhomogenen gewöhnlichen Differentialgleichung  aus seiner Arbeit \cite{E188} geordnet nach Vielfachheit der Nullstellen des charakteristischen Polynoms.}
    \label{fig:E188LösungDGL}
\end{figure}

Aus moderner Sicht leitet Euler eine Rekursionsvorschrift zur Lösung von (\ref{eq: Inhomogeneous Finite}) ab, weshalb seine gefundene Lösung (\ref{eq: Inhomogeneous Finite Solution}) als richtig und  vollständig ausweisen lässt. Nichtsdestotrotz liegt er nicht richtig, wenn er hieraus schließt, seine Methode behielte für $n=\infty$ unverändert Geltung, welche Unrichtigkeit  sich außerdem unmittelbar dem Euler'schen Gedanken selbst demonstrieren lässt: Die essentielle Größe $\alpha$  bestimmt sich nach Euler aus der Nullstellensuche für das charakteristische Polynom:

\begin{equation*}
    P(z)=A+Bz+Cz^2+\cdots +Nz^n=0.
\end{equation*}
Diese besitzt -- vermöge des Fundamentalsatzes der Algebra\footnote{Von der Gültigkeit des Fundamentalsatzes war Euler überzeugt, er hatte zuvor in \cite{E170} einen Beweis vorgelegt, welcher sich allerdings -- nicht mehr zu seinen Lebzeiten -- Kritik, insbesondere von Gauß, ausgesetzt sah. Dies wird weiter unten in Abschnitt \ref{subsubsec: Methodus Inveniendi über Methodus Demonstrandi} eingehender erläutert werden.} -- auch $n$ Lösungen, die man zur vollständigen Lösung der Gleichung (\ref{eq: Inhomogeneous Finite}) benötigt. Im Unendlichen Fall hingegen kann $P(z)$ Funktionen ohne Nullstellen bedeuten, wie etwa:

\begin{equation*}
    P(z)=1+z+z^2+z^3+ \cdots = \dfrac{1}{1-z}.
\end{equation*}
Demnach kann man mit der Euler'schen Methode keine Lösung der Gleichung

\begin{equation*}
    y+ \dfrac{dy}{dx}+\dfrac{dy^2}{dx^2}+\dfrac{dy^3}{dx^3}+\cdots = X
\end{equation*}
finden, obwohl eine solche mit $y=X-X'$ gegeben ist. Die Existenz einer Lösung von solchen Differentialgleichungen ist später überdies durch das Theorem von Malgrange (1928--2024) und Ehrenpreis (1930--2010) gesichert worden\footnote{Die Namensgaber beweisen den nach Ihnen benannten Satz in den Arbeiten \textit{``Existence et approximation des solutions des équations aux dérivées partielles et des équations de convolution"} (\cite{Ma55}, 1955), (``Existenz und Approximation der Lösung von partiellen Differentialgleichungen und Faltungsgleichungen") (Malgrange (1928--2024)) und \textit{``Solution of some problems of division. I. Division by a polynomial of derivation"} (\cite{Eh54}, 1955) und \textit{``Solution of some problems of division. II. Division by a punctual distribution"} (\cite{Eh55}, 1955) (Ehrenpreis (1930--2010). Wie aber auch unten (Abschnitt \ref{subsubsec: Methodus Inveniendi über Methodus Demonstrandi}) erläutert werden wird, haben Existenzfragen Euler selten berührt.}. Auch wenn letzterer für inhomogene lineare partielle Differentialgleichungen mit konstanten Koeffizienten formuliert ist, kann er auf den gewöhnlichen Fall appliziert werden. Dieser Umstand des tatsächlichen Vorliegens einer Lösung scheint Euler indes nicht aufgefallen zu sein, obwohl er  in seinem Lehrbuch zur Integralrechnung \textit{``Institutionum calculi integralis volumen secundum"} (\cite{E366}, 1769, ges. 1763) (E366: ``Grundlagen des Integralkalküls -- Zweiter Band") den endlichen Fall des Beispiels betrachtet (Problem 158, § 1194) und berechnet. Allerdings zeigt seine Lösung  auch auf, dass sie nicht ins Unendliche überführbar ist. Selbiges betrifft sein Beispiel aus § 1209, welches ihn zum selben (unrichtigen) Ausdruck für die Lösung der einfachen Differenzengleichung wie in \cite{E189} führt, zur Besprechung welcher Arbeit im Folgenden übergangen wird.

\subsection{Euler und inhomogene Differenzengleichungen mit konstanten Koeffizienten}
\label{subsec: Ein Prototyp -- Differenzengleichungen}

\epigraph{Man muss immer generalisieren.}{Carl Gustav Jacob Jacobi}

Wäre man gezwungen, Eulers Arbeitsweise anhand einer einzigen Abhandlung zu beschreiben, so wäre  die Arbeit \textit{``De serierum determinatione seu nova methodus inveniendi terminos generales serierum"} (\cite{E189}, 1753, ges. 1749) (E189: ``Über die Bestimmung von Reihen oder eine neue Methode die allgemeinen Termen von Reihen zu finden")  eine treffliche Wahl. Dies begründet sich wie folgt: Zum einen erlaubt sie einen Einblick, wie Euler ein spezielles Gebiet -- hier das der Differenzengleichungen -- mit neuen Ideen expandiert. Präsentiert wird es vom Autor von bekannten Beispielen ausgehend. Alte Formeln werden mit der neuen Methode abgeleitet und erhalten so eine weitere Bestätigung. Zudem werden andere Themengebiete von Euler gestreift -- hier  die Darstellung einer periodischen Funktion als eine Reihe von trigonometrischen Funktionen.  Ein unendliches Produkt, was dieser Arbeit zu eigen sein scheint, erwähnt Euler ebenfalls. Auf alles Aufgelistete wird in der vorliegenden Ausarbeitung entsprechend eingegangen werden. Vorangestellt wird diesen Einzelentdeckungen jedoch die Euler'sche Auflösungsmethode der einfachen Differenzengleichung (Abschnitt \ref{subsubsec: Eulers Lösung der einfachen Differenzengleichung}), welche in einer Herleitung der Stirling'schen Formel für die Fakultät ihre Anwendung finden wird (Abschnitt \ref{subsubsec: Stirling'sche Formel nach Euler}). Die Abweichung von Eulers Ergebnis wird ebenfalls diskutiert und entsprechend korrigiert (Abschnitt \ref{subsubsec: Ein Fehlschluss von Euler}).

\subsubsection{Eulers Lösung der einfachen Differenzengleichung}
\label{subsubsec: Eulers Lösung der einfachen Differenzengleichung}

\epigraph{Solving a problem simply means representing it so as to make the solution transparent.}{Herbert Simon}

Das zentrale Ergebnis von \cite{E189} ist  die Lösung der einfachen Differenzengleichung

\begin{equation}
    \label{eq: Simple Difference Equation}
    f(x+1)-f(x)=g(x).
\end{equation}
Bereits in seiner Arbeit \cite{E47} hatte Euler bereits mit der heute als nach ihm und Maclaurin (1698--1746)benannten Summenformel  eine \textit{partikuläre} Lösung dieser Gleichung angegeben, in dieser Abhandlung reicht er die allgemeine Lösung nach. Diese soll in ihrer Herleitung nachvollzogen werden. 

\begin{figure}
    \centering
    \includegraphics[width=0.9\linewidth]{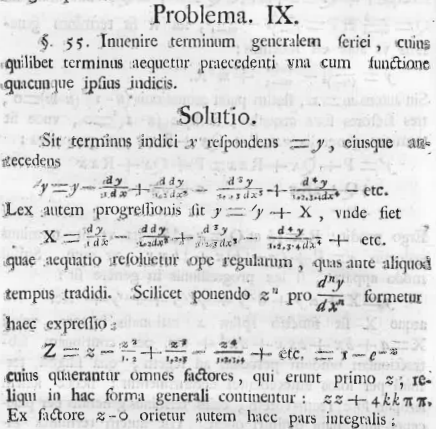}
    \caption{Euler formuliert das Problem zur Lösung der einfachen Differenzengleichung in \cite{E189}.}
    \label{fig:E189AnsatzDifferenzengleichung}
\end{figure}

\paragraph{Gewöhnliche Differentialgleichungen mit konstanten Koeffizienten}
\label{para: Gewöhnliche Differentialgleichungen mit konstanten Koeffizienten}

 Wie oben (Abschnitt \ref{subsec: Ein falsches Ergebnis}) auseinander gesetzt, leitet Euler in \cite{E188} die Lösung der Differentialgleichung

\begin{equation}
    \label{eq: Inhomogeneous Finite}
    \left(a_0+a_1 \dfrac{d}{dx}+a_2 \dfrac{d^2}{dx^2} +\cdots+ a_n \dfrac{d^n}{dx^n} \right)f(x)= g(x),
\end{equation}
mit komplexen Koeffizienten $a_1$, $a_2$, $\cdots$, $a_n$ her\footnote{Für Euler waren die Koeffizienten zwar stets reell, seine Methode funktioniert aber gleichermaßen für den komplexen Fall.}. Aus den Nullstellen des Polynoms

\begin{equation*}
    P(z)=a_0+a_1 z+a_2 z^2 +\cdots+ a_n z^n
\end{equation*}
wird nun explizit die Lösung konstruiert. Sei nun $z=k$ eine solche Lösung bzw. $P(k)=0$. Überdies sei $k$ eine einfache Nullstelle\footnote{Im Folgenden wird ausschließlich  der Fall einer einfachen Nullstelle benötigt und dementsprechend die Mitteilung der Euler'schen Ergebnisse darauf reduziert. Euler hat in \cite{E188} die Fälle der Ordnungen $1$ bis $4$ explizit angegeben, woraus sich leicht der allgemeine Fall erschließen lässt.} von $P(z)$; dann ist die  Lösung zu (\ref{eq: Inhomogeneous Finite}) gegeben als Summe über alle Ausdrücke der Form

\begin{equation}
    \label{eq: Inhomogeneous Finite Solution}
    f_k(x)= \dfrac{e^{kx}}{P'(k)} \int e^{-kx}g(x)dx.
\end{equation}
Da jeder der $n$ Bestandteile jeweils eine Integrationskonstante mitbringt, ist dies auch die vollständige Lösung.\\

Diese Erkenntnis wendet Euler folgend auf die einfache Differenzengleichung (\ref{eq: Simple Difference Equation}) an, indem er zunächst bemerkt, dass  nach dem Satz von Taylor (1685--1731) gilt:

\begin{equation*}
    f(x+1)=f(x)+\dfrac{d}{dx}f(x)+\dfrac{1}{2}\dfrac{d^2}{dx^2}f(x)+\dfrac{1}{3!}\dfrac{d^3}{dx^3}f(x)+\cdots,
\end{equation*}
sodass sich die einfache Differenzengleichung als

\begin{equation}
    \label{eq: Differential Simple Difference}
    \left(\dfrac{d}{dx}+\dfrac{1}{2}\dfrac{d^2}{dx^2}+\dfrac{1}{3!}\dfrac{d^3}{dx^3}+\cdots\right)f(x) =g(x)
\end{equation}
schreiben lässt. Gemäß Eulers  Lösungsmethode sind die Nullstellen des Ausdrucks

\begin{equation}
\label{eq: P(z) infinite}
    P(z) = \dfrac{z}{1!}+\dfrac{z^2}{2!}+\dfrac{z^3}{3!}+\dfrac{z^4}{4!}+\cdots = e^z-1.
\end{equation}
zu ermitteln. Die allgemeine Lösung dieser Gleichung ist mit $z=\log  1$ gegeben.
Kurz zuvor hatte Euler aber auch die Mehr-- bzw. Unendlichwertigkeit der Logarithmusfunktion in seiner Abhandlung \cite{E168} nachgewiesen (siehe auch oben in Abschnitt (\ref{subsubsec: Leitfäden der Arbeitsweise})), sodass er entsprechend eine unendliche Menge von Wurzeln von  (\ref{eq: P(z) infinite}) folgert, nämlich neben der offenkundigen reellen Lösung $z=0$ die weiteren:

\begin{equation*}
    \pm 2 \pi i, \pm 4 \pi i, \pm 6 \pi i, \pm 8 \pi i, \cdots.
\end{equation*}
Demnach liefert die Formel (\ref{eq: Inhomogeneous Finite Solution}) angewandt auf (\ref{eq: Differential Simple Difference}) als Lösung von (\ref{eq: Simple Difference Equation}) den Ausdruck:

\begin{equation}
    \label{eq: Euler Solution Simple Difference}
    f(x) = \int g(x)dx + e^{-2 \pi ix}\int g(x) e^{2 \pi ix}dx +  e^{2 \pi ix}\int g(x) e^{-2 \pi ix}dx
\end{equation}
\begin{equation*}
   +  e^{-4 \pi ix}\int g(x) e^{4 \pi ix}dx +  e^{4 \pi ix}\int g(x) e^{-4 \pi ix}dx+ \cdots
\end{equation*}
Diese Darstellung findet sich entsprechend bei Euler mit $\sin$ und $\cos$ ausgedrückt. Jedoch ist diese in dieser Form nicht ganz korrekt, kommt aber nur wenig von der richtigen Lösung ab, wie weiter unten (Abschnitt \ref{subsubsec: Ein Fehlschluss von Euler}) erkennbar werden wird. Zuvor soll aber noch ein Nebenerkenntnis Erwähnung finden, welches Euler ebenfalls mitteilt.

\subsubsection{Ein Nebenergebnis: Ein unendliches Produkt}
\label{subsubsec: Unendliches Produkt}

\epigraph{One more thing.}{Steve Jobs}

In § 20 von \cite{E189} gibt Euler, in moderner Terminologie, folgendes Produkt an:
\begin{equation}
    \label{eq: Product exp-1}
    e^x -1= \lim_{n \rightarrow \infty} x\prod_{k=1}^{\frac{n}{2}} \left(1+\dfrac{x}{n}+\dfrac{x^2}{4k^2\pi^2}\right).
\end{equation}

\begin{figure}
    \centering
    \includegraphics[scale=1.0]{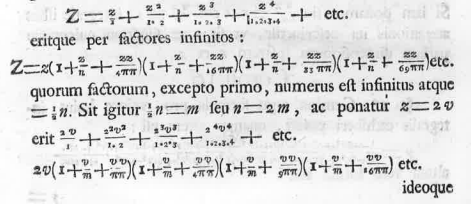}
    \caption{Eulers erklärt die Produktentwicklung von $e^x-1$ aus § 20 von seiner Arbeit \cite{E189}. $n$ und $m$ versteht Euler hier als unendlich große Zahlen.}
    \label{fig:E189Produktex-1}
\end{figure}

Dieses Produkt ist bemerkenswert, weil es den Term $\frac{x}{n}$ enthält, der aber trotz des unendlich werdenden $n$ nicht fortgelassen werden darf, was Euler entsprechend anmerkt. Insbesondere ergibt sich daraus, dass sich das Produkt -- eben nicht wie etwa bei $\sin x$ -- bereits  aus den Nullstellen ableiten lässt. Wie Euler zu dieser Erkenntnis gelangt ist, wird aus seinen Ausführungen in seiner Arbeit nicht deutlich\footnote{Auch die Ausführung in seiner \textit{Introductio}, wo sich das Produkt in § 155 findet, geben keinen Aufschluss.}. Es ist aber anzunehmen, dass er sie auf indirektem Wege erlangt hat. Bereits in seiner  \textit{Introductio} \cite{E101} hatte er zuvor die Formel

\begin{equation*}
    \dfrac{e^x-e^{-x}}{2}= x \prod_{n=1}^{\infty}   \left(1+\dfrac{x^2}{(n\pi)^2}\right)
\end{equation*}
angegeben, welche sich durch Fortlassen des Terms $\frac{x}{n}$ in  (\ref{eq: Product exp-1}) ergäbe. Demnach würde es zu einem Widerspruch führen, selbigen tatsächlich auszulassen.
 Man findet gleichermaßen durch direktes Ausmultiplizieren von (\ref{eq: Product exp-1}), dass das Weglassen vom richtigen Ergebnis wegführt, zumal

\begin{equation*}
    \prod_{k=1}^{\frac{n}{2}} \left(1+\dfrac{x}{n}+\dfrac{x^2}{4k^2\pi^2}\right)=1+ \dfrac{n}{2}\cdot \dfrac{x}{n}+\cdots,
\end{equation*}
woraus sich direkt der richtige Term der Potenzreihenentwicklung ergibt, wenn man den Grenzwert $n \rightarrow \infty$ nimmt. Weiter gibt Euler in § 21 von \cite{E189} auch die Partialbruchzerlegung 

\begin{equation}
    \label{eq: Partial Fraction Decomposition ex-1}
    \dfrac{2e^{2v}}{e^{2v}-1}-\dfrac{1}{v}=\lim_{m \rightarrow \infty} \sum_{k=1}^{m} \dfrac{\frac{1}{m}+\frac{2v}{4\pi \pi}}{1+\frac{v}{m}+\frac{v^2}{k^2\pi^2}},
\end{equation}
an, aus welcher er die Werte der Summen $\sum_{k=1}^{\infty} \frac{1}{k^{2n}}$ ableitet\footnote{In Abschnitt (\ref{subsubsec: Herleitung Formel für die Potenzsummen der Reziproken}) wird gar ein weiterer Beweis vorgestellt, diese Summen zu evaluieren, welcher alleinig auf den Erkenntnissen von Eulers Arbeit \cite{E189} basiert.}.\\

 Obschon das Euler'sche Resultat richtig ist,  scheint er keine Methode vorgelegt zu haben, mit welcher sich  a priori zu (\ref{eq: Product exp-1}) gelangen ließe. Für Beispiele wie das Sinus--Produkt (\ref{eq: sine-Product}) lassen sich die Erläuterungen aus seiner Arbeit \textit{``Dilucidationes in capita postrema calculi mei differentalis de functionibus inexplicabilibus"} (\cite{E613}, 1813, ges. 1780) (E613: ``Erläuterungen zu den letzten Kapiteln meines \textit{Caculi Differentialis} über unerklärbare Funktionen") dahingehend auslegen, dass sie eine solche Methode darbieten. Wie etwa in der Arbeit \textit{``Euler and the Gammafunction"} (\cite{Ay21}, 2021) argumentiert und weiter unten (Abschnitt \ref{subsubsec: Ein anderes Vorhaben: Das Weierstraß-Produkt}) noch einmal aufgegriffen nimmt Euler in besagter Arbeit die Weierstraß'sche Produktentwicklung für eine Funktion mit vorgegebenen Nullstellen vorweg, welche indes nicht auf $e^x-1$ anwendbar ist.

\subsubsection{Stirling'sche Formel nach Euler}
\label{subsubsec: Stirling'sche Formel nach Euler}

\epigraph{The most exciting phrase to hear in science, the one that heralds new discoveries, is not 'Eureka!' but 'That's funny...'}{Isaac Asimov}

Im letzten Abschnitt seiner Abhandlung \cite{E189} (ab §§ 56--60) möchte Euler sein allgemeines Ergebnis (\ref{eq: Euler Solution Simple Difference}) auf die Fakultät anwenden\footnote{Euler betrachtet dabei nicht die Differenzengleichung für $x!$, sondern die für $(x-1)!$, also eine, die von der $\Gamma$--Funktion erfüllt wird.}; selbige erfüllt die Gleichung

\begin{equation*}
    y(x+1)=xy(x),
\end{equation*}
die durch Logarithmieren in eine Differenzengleichung der benötigten Form überführt wird

\begin{equation*}
    \log  y(x+1)-\log  y(x) =\log  (x).
\end{equation*}
Am Ende intensiver Rechnungen steht die Stirling'sche Formel

\begin{equation}
\label{eq: Stirling}
    n! \sim \sqrt{2\pi n}\left(\dfrac{n}{e}\right)^n \quad \text{für} \quad n \rightarrow \infty,
\end{equation}
jedenfalls beinahe. Euler unterläuft wie weiter oben schon angedeutet  auf der konzeptuellen Ebene ein Fehler, welcher aus der Vorstellung des Euler'schen Gedankenganges erkennbar werden wird. \\

Die Applikation der allgemeinen Formel (\ref{eq: Euler Solution Simple Difference}) führt in einem ersten Schritt zu

\begin{equation}
\label{eq: General Solution Factorial}
    f(x) = x\log x-x +C + e^{2 \pi i x}\int \log(x) e^{-2 \pi i x}dx+  e^{-2 \pi i x}\int\log(x) e^{+2 \pi i x}dx
\end{equation}
\begin{equation*}
    +  e^{4 \pi i x}\int \log(x) e^{-4 \pi i x}dx + e^{-4 \pi i x}\int \log(x) e^{+4 \pi i x}dx + \cdots
\end{equation*}
wo $\int \log(x)dx$ bereits zu $x\log(x) -x$ ausgewertet worden ist und $C$ eine beliebige Integrationskonstante bedeutet{\footnote{Dies ist die Lösung, welche Euler in \cite{E189} \S 59 angibt, er verwendet lediglich $\cos$ und $\sin$-Funktion statt $e$--Funktionen. Überdies ist diese Darstellung als eine asymptotische und nicht wirkliche Gleichheit zu verstehen.}. Um nun zur Stirling'schen Formel gelangen, sucht Euler zunächst den Ausdruck

\begin{equation*}
    e^{2 k \pi i x} \int e^{-2k \pi ix}\log (x)dx.
\end{equation*}
Dafür verwendet er partielle Integration unendlich viele Male hintereinander mit $e^{-2 k \pi ix}$ als zu integrierender Funktion. In moderner und kompakter Notation ist das Ergebnis\footnote{Weil Euler $\sin(2k \pi x)$ und $\cos(2k \pi x)$ anstatt $e^{-2 k \pi i x}$ verwendet, weicht sein endgültiges Ergebnis in seiner Gestalt von dem ab, zu welchem hier gelangt werden wird. Die Herleitung ist allerdings in beiden Fällen dieselbe.}:

\begin{equation*}
    e^{2 k \pi i x} \int e^{-2k \pi ix}\log (x)dx= -\dfrac{\log(x)}{2 k \pi i}+ \sum_{n=1}^{\infty} \dfrac{(-1)^n  (n-1)!}{(2k\pi i)^{n+1} x^n}+ C_k e^{2 k \pi i x}.
\end{equation*}
$C_k$ ist eine beliebige durch Integration eingehende Konstante. Indem man so für alle unendlich vielen Ausdrücke verfährt, gelangt man zum Ausdruck:

\begin{equation*}
    \log y(x)= x\log x- x + C +\sum_{k\in \mathbb{Z} \setminus \lbrace 0 \rbrace} \left(C_k e^{2k \pi i x}-\dfrac{\log(x)}{2 k \pi i}+ \sum_{n=1}^{\infty} \dfrac{(-1)^n  (n-1)!}{(2k\pi i)^{n+1}x^n} \right)
\end{equation*}
Zwecks Vereinfachung dieser Summe setze man zunächst

\begin{equation*}
    C+ \sum_{k\in \mathbb{Z} \setminus \lbrace 0 \rbrace} C_k e^{2k \pi i x}  =: h(x);
\end{equation*}
$h$ ist demnach eine periodische Funktion, welche der Gleichung $h(x)=h(x+1)$ Genüge leistet. Ferner gilt

\begin{equation*}
    \sum_{k\in \mathbb{Z} \setminus \lbrace 0 \rbrace} \dfrac{\log (x)}{2 k \pi i} =0,
\end{equation*}
da sich die Terme gegenseitig streichen. Demnach verbleibt die Berechnung der Doppelsumme. Mit einer formalen Rechnung eruiert man:

\begin{equation}
\label{eq: Double Sum}
\sum_{k\in \mathbb{Z} \setminus \lbrace 0 \rbrace}  \sum_{n=1}^{\infty} \dfrac{(-1)^n  (n-1)!}{(2k\pi i)^{n+1}x^n}   =  \sum_{n=0}^{\infty} \sum_{k=1}^{\infty} \dfrac{2}{k^{2n+2}} \cdot \dfrac{ (-1)^{n} (2n)!}{(2\pi)^{2n+2}\cdot x^{2n+1}}.
\end{equation}
Die Summe rechter Hand hatte  Euler, wie oben in Form der Gleichung (\ref{eq: zeta(2n)}) erwähnt, zuvor berechnet;  man entnimmt sie in \cite{E130} als:

\begin{equation}
    \sum_{k=1}^{\infty} \dfrac{1}{k^{2n}} = \dfrac{(-1)^{n-1} (2\pi)^{2n}B_{2n}}{2(2n)!},
\end{equation}
wo $B_n$ die Bernoulli Zahlen sind. Setzt man dies nun in (\ref{eq: Double Sum}) ein, wird dies zu

\begin{equation*}
    \sum_{k\in \mathbb{Z} \setminus \lbrace 0 \rbrace}  \sum_{n=1}^{\infty} \dfrac{(-1)^n  (n-1)!}{(2k\pi i)^{n+1}x^n}= \sum_{n=0}^{\infty} 2 \cdot \dfrac{(-1)^{n} (2\pi)^{2n+2}B_{2n+2}}{2(2n+2)!}\cdot \dfrac{ (-1)^{n} (2n)!}{(2\pi)^{2n+2}\cdot x^{2n+1}}.
\end{equation*}
Viele Terme heben sich auf, sodass diese Gleichheit verbleibt:

\begin{equation*}
      \sum_{k\in \mathbb{Z} \setminus \lbrace 0 \rbrace}  \sum_{n=1}^{\infty} \dfrac{(-1)^n  (n-1)!}{(2k\pi i)^{n+1}x^n}= \sum_{n=1}^{\infty} \dfrac{B_{2n}}{(2n-1)2n x^{2n-1}}.
\end{equation*}
Alles in eingesetzt (\ref{eq: General Solution Factorial}) gibt:

\begin{equation}
\label{eq: Euler Factorial}
    \log y(x) = x\log x -x +h(x)+ \sum_{n=1}^{\infty} \dfrac{B_{2n}}{(2n-1)2n x^{2n-1}}.
\end{equation}
Diese Gleichung ist  in moderner Terminologie als asymptotische Reihe aufzufassen und es ist die Formel zu welcher Euler selbst in § 60 von \cite{E189} gelangt ist, Euler hat lediglich die expliziten Zahlenwerte für die Bernoulli Zahlen substituiert. (\ref{eq: Euler Factorial}) mit (\ref{eq: Stirling}) vergleichend, fehlt der Term $\log(\sqrt{\frac{2\pi}{x}})$, was Euler auch bemerkt. Er argumentiert diesbezüglich, dass besagter aus einem Spezialfall folge, also etwa $x=1$\footnote{Genauer argumentiert Euler, dass in diesem Fall $h(x)$ als konstante Größe zu verstehen ist und der Wert dieser Konstante hier dann  der Summe $1-\sum_{n=1}^{\infty} \frac{B_{2n}}{(2n-1)2n}$ gleich ist, welche Euler ohne Beweis behauptet $\frac{1}{2}\log (2\pi)$ zu sein, obwohl die Reihe wegen des raschen Anwachsens der Bernoulli Zahlen nicht konvergiert. Aber Euler war sich indes gewahr, dass man den zuvor erwähnten Wert der Summe zuschreiben kann, weil er der Konstante $\sqrt{2\pi}$ in Stirlings Formel (\ref{eq: Stirling}) entspricht. Euler selbst hat einen korrekte Berechung der Stirling'schen Konstante in seinem Buch \textit{``Calculi Differentialis"} (§§ 157--158) \cite{E212} gegeben.} und der Anfangsbedingung $y(1)=1$ zur Gleichung $y(x+1)=xy(x)$, sodass man bei dieser finalen Formel angelangt: 

\begin{equation}
\label{eq: Euler Final log}
     \log y(x) = x\log x -x +\log (\sqrt{2 \pi})+ \sum_{n=1}^{\infty} \dfrac{B_{2n}}{(2n-1)2n x^{2n-1}},
\end{equation}
sofern $x$ unendlich groß ist. \\

\begin{figure}
    \centering
   \includegraphics[scale=1.0]{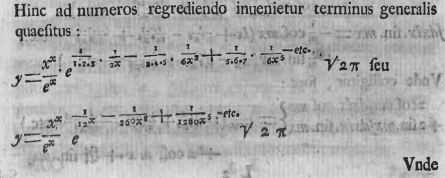}
    \caption{Eulers (unrichtige) Version der Stirling'schen Formel aus seiner Arbeit \cite{E189} abgeleitet aus der Differenzengleichung der Fakultät.}
    \label{fig:E189Stirling}
\end{figure}

Eulers Version der Formel liest sich dann wie folgt:

\begin{equation}
    \label{eq: Euler Factorial Final}
    y(x) = \dfrac{x^x}{e^x}\left(1+\dfrac{1}{12x}+\dfrac{1}{288x^2}-\dfrac{139}{51840x^3}+\cdots\right)\sqrt{2\pi},  
\end{equation}
welche folgt, indem man die Zahlenwerte für die Bernoulli-Zahlen einsetzt, die Exponentialfunktion nimmt und dann ihre Taylorentwicklung heranzieht.\\

Wie Faber (1877--1966) in einer Fußnote der \textit{Opera Omnia} Version der Arbeit \cite{E189} anmerkt, ist die von Euler angegebene Lösung (\ref{eq: Euler Factorial Final}) inkorrekt, und ist wie folgt zu korrigieren

\begin{equation}
\label{eq: Log Stirling correct}
     \log y(x) = x\log x -x +\log (\sqrt{\dfrac{2 \pi}{x}})+ \sum_{n=1}^{\infty} \dfrac{B_{2n}}{(2n-1)2n x^{2n-1}},
\end{equation}
In der Euler'schen Lösung ist demnach der Term $-\log  \sqrt{x}$ abhanden\footnote{Da sich $y(x)$ als $\Gamma(x)$ deuten lässt, gilt entsprechend für die Fakultät

\begin{equation*}
    \log(x!) = x \log x -x + \log \sqrt{2\pi x},
\end{equation*}
sofern auch hier das Gleichheitszeichen hier wieder als asymptotisch äquivalent gelesen wird.} . 

\subsubsection{Ein Fehlschluss von Euler}
\label{subsubsec: Ein Fehlschluss von Euler}

\epigraph{The one who insists on never uttering an error must remain silent.}{Werner Heisenberg}

Zweifelsohne ist Eulers Antwort (\ref{eq: Euler Solution Simple Difference}) zur Auflösung der Gleichung (\ref{eq: Simple Difference Equation}) unrichtig, weicht  jedoch  wenig von der vollständigen Lösung

\begin{equation}
    \label{eq: Simple Difference Correct Solution}
    f(x) ={\color{red}-\dfrac{1}{2}g(x)}+ \int g(x)dx + e^{-2 \pi ix}\int g(x) e^{2 \pi ix}dx +  e^{2 \pi ix}\int g(x) e^{-2 \pi ix}dx
\end{equation}
\begin{equation*}
     e^{-4k \pi ix}\int g(x) e^{4 \pi ix}dx +  e^{4 \pi ix}\int g(x) e^{-4 \pi ix}dx+ \cdots,
\end{equation*}
ab. Genauer kommt Eulers Ergebnis lediglich um den Term $-\frac{1}{2}g(x)$ von der Wahrheit ab. Euler  hat diesen Term, wie bereits angedeutet, nicht aus Flüchtigkeit übersehen, sondern der Lapsus ist der Konstruktion der Lösung aus den Nullstellen des Ausdrucks $P(z)=e^z-1$ geschuldet. Obschon nämlich die Lösungsmethode für den Fall endlicher Ordnung zum richtigen Ergebnis führt, gelingt der Übergang zum Unendlichen nicht reibungslos. Für den erfolgreichen Übergang zum Fall unendlicher Ordnung wäre  die Lösung  aus dem Kehrwert des charakteristischen Polynoms zu konstruieren. Mit der Festlegung $z=\frac{d}{dx}$ schreibe man (\ref{eq: Differential Simple Difference}) wie folgt:

\begin{equation}
\label{eq: Solution formal}
    f(x) = \dfrac{1}{P(z)}g(x).
\end{equation}
Zur Anwendung des Operators $\frac{1}{P(z)}$ auf $g(x)$ ist selbiger in ganzzahligen Potenzen von $z$ anzugeben, wofür mehrere mehrere Möglichkeiten offen stehen. Für den Nachweis von
 (\ref{eq: Simple Difference Correct Solution}) ist die folgende Partialbruchzerlegung zuträglich, welche sich leicht mit Methoden aus der komplexen Analysis zeigen lässt\footnote{Euler selbst hätte zu dieser Formel gelangen können, allerdings nicht mit einer allgemeinen Methode.}.

 \begin{equation}
    \label{eq: Partial Fraction Decomposition}
    \dfrac{1}{e^z-1}=-\frac{1}{2}+ \sum_{k \in \mathbb{Z}\setminus \lbrace 0\rbrace}^{\infty} \dfrac{1}{z-2k\pi i}.
\end{equation}
Die nur bedingte Konvergenz der letzten Formel ist für das Folgende nicht weiter wichtig. ist Man bedarf des Ausdrucks

\begin{equation}
    \label{eq: Single Term}
    \dfrac{1}{z-2 k \pi i}g(x),
\end{equation}
welcher sich mit $2k\pi =\alpha$ als

\begin{equation*}
    \dfrac{1}{z-\alpha}g(x) = \dfrac{1}{z \left(1- \frac{\alpha}{z}\right)}g(x) = \sum_{n=0}^{\infty}\dfrac{\alpha^n}{z^{n+1}}g(x)
\end{equation*}
erweist. Wegen $\frac{d}{dx}=z$  ist $\frac{1}{z}$ als $\int$--Operator aufzufassen, als Konsequenz $\frac{1}{z^n}$ als $n$--fach iteriertes Integral\footnote{Für iterierte Integrale hat Euler ebenfalls Formeln angegeben, und zwar in \cite{E681}, allerdings erst nach der Veröffentlichung der Arbeit \cite{E189}. Man vergleiche \ref{para: Mitteilung von kuriosen Ergebnissen}.}. Mit der Notation $\int^n$ für das $n$--fach iterierte Integral gelangt man zur folgenden Formel:

\begin{equation}
    \label{eq: Iterated integral}
    \int^n g(x)dx = \int\limits_{}^{x}\dfrac{(x-t)^{n-1}}{(n-1)!}g(t)dt.
\end{equation}
Eingesetzt in (\ref{eq: Single Term}) gibt dies

\begin{equation*}
     \dfrac{1}{z-\alpha}g(x)= \sum_{n=0}^{\infty} \alpha^n \int\limits_{}^{x}\dfrac{(x-t)^{n}}{n!}g(t)dt = \int\limits_{}^{x} e^{\alpha (x-t)}g(t)dt.
\end{equation*}
 Vermöge (\ref{eq: Partial Fraction Decomposition}) und (\ref{eq: Solution formal}) zeigt dies:
\begin{equation*}
f(x) = -\dfrac{1}{2}g(x)+\sum_{k\in \mathbb{Z}} \int\limits_{}^{x} e^{2 k\pi i (x-t)}g(t)dt =  -\dfrac{1}{2}g(x)+\sum_{k\in \mathbb{Z}}e^{2 k \pi i x} \int\limits_{}^{x} e^{-2 k\pi i t}g(t)dt,
\end{equation*}
welche  (\ref{eq: Simple Difference Correct Solution}) ist und mit der Lösung in  der Arbeit 
 \textit{``Note on the simple Difference Equation} (\cite{We14}, 1914) übereinstimmt, wo sie mithilfe der komplexen Analysis abgeleitet wird. Auf die Art, wie im vorgestellten Argument mit Differentialoperatoren hantiert worden ist, ist auch Bourlet (1866--1913) in seiner Arbeit \textit{``Sur certaines équations analogues aux équations différentielles"} (\cite{Bo99}, 1899) (``Über gewisse den differentiellen analoge Gleichungen") verfahren.\\

\paragraph{Verbindung zur Euler--Maclaurin'schen Summenformel}

Gründe der Vollständigkeit gebieten es, den Euler'schen Ansatz zur Lösung von Differentialgleichungen unendlicher Ordnung zur Herleitung  der Euler-Maclaurin'schen Summenformel anzuwenden, welche wie eingangs angemerkt eine partikuläre Lösung der einfachen Differenzengleichung (\ref{eq: Simple Difference Equation}) gibt. 
 Man hat demnach immer noch die Differentialgleichung (\ref{eq: Differential Simple Difference}) zu lösen, welche sich für $z=\frac{d}{dx}$ als

\begin{equation}
\label{eq: Formelle Lösung der Differenzengleichung}
    f(x) = \dfrac{1}{e^z-1}g(x)
\end{equation}
schreiben lässt. Statt nun die Funktion $\frac{1}{e^z-1}$ in Partialbrüche zu umzuwandeln, suche man stattdessen die Laurent--Entwicklung um den Ursprung. Mit der schon oben eingeführten Definition der Bernoulli--Zahlen

\begin{equation*}
    \dfrac{z}{e^z-1}=\sum_{n=0}^{\infty} \dfrac{B_n}{n!}z^n= B_0+\dfrac{B_1}{1!}+\dfrac{B_2}{2!}z^2+\cdots
\end{equation*}
findet man:

\begin{equation*}
    \dfrac{1}{e^z-1}=\dfrac{B_0}{z}+\dfrac{B_1}{1!}+\dfrac{B_2}{2!}z+\cdots,
\end{equation*}
sodass mit den Euler'schen Ideen $\int$ als $z^{-1}$  und $z^n$ als $\frac{d^n}{dx^n}$ zu sehen aus (\ref{eq: Formelle Lösung der Differenzengleichung}) nachstehende Formel entspringt:

\begin{equation*}
    f(x)= B_0\int g(x)dx + \dfrac{B_1}{1!}g'(x)+\dfrac{B_2}{2!}g''(x)+\cdots,
\end{equation*}
was gerade die  Summenformel in einer von Euler bevorzugten Form ist. Er selbst scheint sie allerdings nicht auf diesem Wege in irgendeiner seiner Arbeiten hergeleitet zu haben.

\paragraph{Mögliche Ursachen für den Euler'schen Fehlschluss}
\label{para: Warum ihm dieser Fehler wohl unterlaufen ist}

Abschließend soll die   Euler'sche Lösung (\ref{eq: Euler Solution Simple Difference}) zu (\ref{eq: Simple Difference Equation}) auf die Ursachen des Fehlen des Terms $-\frac{1}{2}g(x)$ hin besprochen werden. Ein erster Grund für das Übersehen dieses Terms von Eulers Seite mag folgender sein: In seiner Arbeit \cite{E189} behandelt Euler ausschließlich Beispiele, deren Lösung  ihm bereits aus anderer Quelle bekannt ist. Daher wird es zur Überzeugung der Korrektheit seiner Methode beigetragen haben, dass sie in all seinen vorgestellten Beispielen zu den richtigen Ergebnissen führt. So wird etwa die Gleichung

\begin{equation*}
    f(x+1)=f(x)
\end{equation*}
natürlich von der allgemeinen periodischen Funktion

\begin{equation*}
    f(x)= \sum_{k \in \mathbb{Z}} c_ke^{2k\pi i x}
\end{equation*}
gelöst. Auch für den Spezialfall 

\begin{equation*}
    f(x+1)-f(x)=K,
\end{equation*}
für eine beliebige Konstante $K$ gelangt er vermöge seiner Formel (\ref{eq: Euler Solution Simple Difference}) zur korrekten Lösung, da das Fehlen des Terms $-\frac{1}{2}K$ hier nicht ins Gewicht fällt. Die einzigen Ausnahmen bilden  besagte einfache Differenzengleichung (\ref{eq: Simple Difference Equation}) und der danach behandelte Spezialfall der Fakultät. Da er  den Vergleich mit der Euler--Maclaurin'schen Summenformel nicht unternimmt, ist das Fehlen des Terms nicht zu bemerken, zumal Euler von der Richtigkeit seiner Methode zweifelsohne sehr überzeugt sein musste. Anders lässt sich es wohl auch nicht erklären, dass er am Beispiel der Stirling'schen Formel mit einem falschen Ergebnis konfrontiert -- er hatte unter Anderem in seinen \textit{Calculi Differentialis} \cite{E212} (§§ 157--158) einen Beweis der korrekten Formel gegeben -- versucht, selbiges mit einer angeblichen Richtigkeit für unendlich große $x$ zu erhärten versucht. Er scheint diesen Fehler in seiner Argumentation in keiner seiner späteren Arbeiten korrigiert zu haben.\\

Obgleich der Fehler ein kleiner ist, hätte er, so sollte plausibel geworden sein, Euler auffallen können. Selbstredend ist  Euler von jedwedem Vorwurf frei zu sprechen, trotz Kenntnis aller seiner Bestandteile nicht den obigen Gedankengang angegeben zu haben\footnote{Zumal dieser erst ruhigen Gewissens mit Kenntnis der Euler noch fremden Fourier-Analyse akzeptiert werden kann, welche die algebraische Handhabung der Differentialoperatoren überhaupt erst rechtfertigt.}.  
Gleichermaßen kann Euler die Unkenntnis der Partialbruchzerlegung (\ref{eq: Partial Fraction Decomposition}) nicht zur Last gelegt werden, welche ebenfalls Methoden aus der komplexen Analysis verlangt, um a priori zweifelsfrei  bestätigt zu werden, obschon sie sich aus derjenigen für den $\cot$ auch herleiten ließe\footnote{Man hat zunächst 
\begin{equation*}
     \cot x = \dfrac{\cos x}{\sin x}=i \cdot \dfrac{e^{ix}+e^{-ix}}{e^{ix}-e^{-ix}}= i \cdot \dfrac{e^{2ix}+1}{e^{2ix}-1}.
\end{equation*}
Mit $y=2ix$ gilt jedoch auch 
\begin{equation*}
    \dfrac{1}{2}\left(\dfrac{e^y+1}{e^y-1}-1\right)= \dfrac{1}{2}\cdot \dfrac{e^y+1-e^y+1}{e^y-1}= \dfrac{1}{e^y-1},
\end{equation*}
sodass sich mit der Partialbruchzerlegung für $\cot (x)$ aus (\ref{eq: Partialbruchzerlegung cot}) diese Formel ergibt
\begin{equation*}
    \dfrac{1}{e^y-1}= \dfrac{1}{y}-\dfrac{1}{2}+\sum_{n=1}^{\infty} \dfrac{2y}{x^2+4n^2\pi^2}.
\end{equation*}}. Zwar hat Euler auch die Partialbruchzerlegung von weiteren transzendenten Funktionen in seiner Arbeit  \textit{``De resolutione fractionum transcendentium in infinitas fractiones simplices"} (\cite{E592}, 1785, ges. 1775) (E592: ``Über die Auflösung von transzendenten Brüchen in unendlich viele einfache Brüche")  zum Gegenstand gemacht und dort viele Beispiele mit der dortigen Methode mit glücklichem Erfolg behandelt, jedoch versagt selbige just beim für die Differenzengleichung benötigten Beispiel $\frac{1}{e^z-1}$. Interessanterweise lässt sich aber die entsprechende Partialbruchzerlegung aus dem unendlichen Produkt (\ref{eq: Product exp-1}) ableiten und Euler selbst hat in seiner Arbeit \cite{E189} (§ 21) eine Darstellung gefunden, die dieser gleichwertig ist. Er definiert nämlich zunächst

\begin{equation*}
    V:= \dfrac{e^{2v}-1}{2v},
\end{equation*}
sodass 

\begin{equation*}
    \dfrac{dV}{Vdv}= \dfrac{2e^{v}}{e^{v}-e^{-v}}-\dfrac{1}{v},
\end{equation*}
aber, wie Euler schreibt, gilt ebenfalls

\begin{equation*}
    \dfrac{dV}{Vdv}= \dfrac{\frac{1}{m}+\frac{2v}{1\pi\pi}}{1+\frac{v}{m}+\frac{vv}{1\pi \pi}}+ \dfrac{\frac{1}{m}+\frac{2v}{4\pi\pi}}{1+\frac{v}{m}+\frac{vv}{4\pi \pi}}+\dfrac{\frac{1}{m}+\frac{2v}{9\pi\pi}}{1+\frac{v}{m}+\frac{vv}{9\pi \pi}}+\cdots,
\end{equation*}
wobei $m$ eine unendlich große Zahl ist. Modern ließe sich diese Formel wie folgt schreiben:

\begin{equation*}
    \dfrac{dV}{Vdv}= \lim_{m \rightarrow \infty} \sum_{n=1}^{m} \dfrac{\frac{1}{m}+\frac{2v}{(n\pi)^2}}{1+\frac{v}{m}+\frac{v^2}{(n \pi)^2}}.
\end{equation*}
Euler hätte jedoch bereits aus dem unendlichen Produkt (\ref{eq: Product exp-1})  das Fehlen des Terms $-\frac{1}{2}g(x)$ in seiner Lösung (\ref{eq: Euler Solution Simple Difference}) auffallen können.

\newpage 

\section{Von Euler vorweggenommene Entdeckungen}
\label{sec: Von Euler vorweggenommene Entdeckungen}

\epigraph{All mathematicians alive are his disciples: there is no one who is not guided and sustained by the genius of Euler.}{Marquis de Condorcet}


Dieser Abschnitt handelt von Schätzen, Funden, Lehrsätzen und dergleichen, die trotz Euler'scher Erstentdeckerschaft den Namen seiner Nachfolger tragen\footnote{Dieses Phänomen, dass eine Entdeckung nicht immer nach ihrem Erstentdecker benannt wird, tritt neben der Mathematik auch in anderen Wissenschaften auf und wird in dem Artikel \textit{``Stigler's law of eponymy"} von S. Stigler (1941--) aus dem Buch \textit{``Science and social structure: A festschrift for Robert K. Merton"} (\cite{Gi80}, 1980) in humoristischer Weise diskutiert und seitdem, nach dem Artikel, als \textit{Stigler's law of eponymy} bezeichnet.}. Diese Entdeckungen umfassen schlicht übersehene (Abschnitt \ref{subsec: Übersehenes}) zum einen, in äquivalenter Gestalt bewiesene und dadurch missachtete (Abschnitt \ref{subsec: Von Euler in anderer Gestalt Bewiesenes}) zum anderen. Wohingegen der erstere durch direkte Angabe der Entdeckungen eine schnelle Behandlung erlauben, ist bei zweiteren der jeweilige Kontext mit einzubeziehen, um das Übergehen von Eulers Beiträgen zu erklären. Es wird sich zeigen, dass Euler, freilich in anderer Form, die Multiplikationsformel für die $\Gamma$--Funktion (Abschnitt \ref{subsubsec: Die Multiplikationsformel für die Gammafunktion}) sowie einige Eigenschaften der Legendre--Polynome (Abschnitt \ref{subsubsec: Den Kontext betreffend: Die Legendre Polynome}) zutage gefördert hat. Andere Fragestellungen als die modernen Entsprechungen haben ihn sogar das Weierstraß--Produkt vorweg nehmen  (Abschnitt \ref{subsubsec: Ein anderes Vorhaben: Das Weierstraß-Produkt})  lassen und  zur  Mellin--Transformatierten geführt (Abschnitt \ref{subsubsec: Die Mellin--Transformierte bei Euler}). Abschließend zeigt das Beispiel der hypergeometrischen Reihe (Abschnitt \ref{subsubsec: Die Darstellung betreffend -- Die hypergeometrische Reihe}), dass die Unkenntnis der Euler'schen Beiträge zu einem Gegenstand deren Verstreuung auf verschiedene zeitlich und inhaltlich weit auseinander liegenden Abhandlungen zur Ursache haben kann.

\subsection{Von seinen Nachfolgern übersehene Entdeckungen}
\label{subsec: Übersehenes}

\epigraph{There are three stages in scientific discovery. First, people deny that it is true, then they deny that it is important; finally they credit the wrong person.}{Bill Bryson}


Bei der immensen Anzahl an Euler'schen Publikationen ist es Eulers Nachfolgern leicht nachzusehen, nicht jede Entdeckung ihres Vorgängers zu kennen. Gleichermaßen ist der Vorwegnahme einiger Ergebnisse durch Euler in exakt der Form, wie sie später von anderen unabhängig gefunden worden sind, so weniger verwunderlich, wovon zwei kleinere Beispiele hier erwähnt werden sollen: Das der Fourier--Koeffizienten (Abschnitt \ref{subsubsec: Die Fourierkoeffizienten}) und das der Produktformel der $\Gamma$--Funktion (Abschnitt \ref{subsubsec: Produktdarstellung für die Gamma-Funktion}), die für gewöhnlich Gauß zugeschrieben wird.

\subsubsection{Die Fourierkoeffizienten}
\label{subsubsec: Die Fourierkoeffizienten}

\epigraph{No scientific discovery is named after its original discoverer.}{Stephen Stigler}

Obschon die Fourier--Analyse von Euler kaum berührt und das Prinzipat ihrer Ausarbeitung J. Fourier (1768--1830) bzw. seinem Werk \textit{``Théorie analytique de la chaleur"} (\cite{Fo22}, 1822) (``Analytische Theorie der Wärme") gebührt, findet sich eine zentrale Entdeckung bereits in den Euler'schen Werken vergraben: Die Fourierkoeffizienten. Zusätzlich soll angeschnitten werden, wie nahe Euler der Darstellung einer beliebigen Funktion durch trigonometrische Reihen kommt.

\paragraph{Euler und Fourierreihen}
\label{para: Euler und Fourrierreihen}

In Problem 1 (§ 12--13) seiner Arbeit \cite{E189} gelangt Euler zur Aussage, dass eine periodische Funktion sich in eine Reihe von $\sin$ und $\cos$ zerlegen lässt. Eulers Gedankengang ist dabei folgender: \\

\begin{figure}
    \centering
   \includegraphics[scale=1.0]{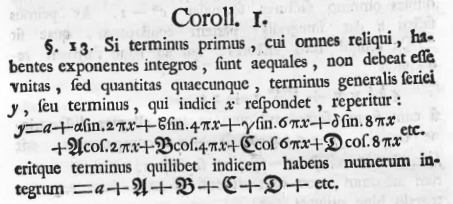}
    \caption{Eulers Ausdruck der Lösung der Gleichung $y(x+1)=y(x)$ aus § 13 seiner Arbeit \cite{E189}, welcher die allgemeine Fourierentwicklung einer periodischen Funktion darstellt.}
    \label{fig:E189Fourierreihe}
\end{figure}

Er unternimmt in \cite{E189} die Auflösung der Gleichung

\begin{equation*}
    y(x+1)=y(x),
\end{equation*}
zu welchem Zweck er gemäß seiner oben (Abschnitt \ref{subsubsec: Differentialgleichungen unendlicher Ordnung: Der homogene Fall}) vorgestellten Methode zur Auflösung von Differentialgleichungen unendlicher Ordnung jedem der unendlich vielen Werte von $\log (1)$ eine entsprechende Lösung zuordnet. Sein Vorgehen oder alternativ seine Formel (\ref{eq: Euler Solution Simple Difference}) führen zu

\begin{equation*}
    y(x)= \sum_{k=-\infty}^{\infty} C_k e^{2k\pi i x},
\end{equation*}
was  gerade die Entwicklung einer periodischen Funktion in eine Fourier--Reihe ist. Euler bevorzugt in \cite{E189} lediglich die Darstellung mit $\sin$- und $\cos$.\\

In dieser Arbeit geht Euler nicht weiter auf diesen Punkt ein. Jedoch nimmt er in den späteren Arbeiten \textit{``Methodus facilis inveniendi series per sinus cosinusve angulorum multiplorum procedentes, quarum usus in universa theoria astronomiae est amplissimus"} (\cite{E703}, 1798, ges. 1777) (E703: ``Eine leichte Methode, nach Sinus und Kosinus vielfacher Winkel fortschreitender Reihen zu summieren, deren Nutzen in der ganzen Astronomie sehr umfassend ist") und \textit{``Disquisitio ulterior super seriebus secundum multipla cuiusdam anguli progredientibus"} (\cite{E704}, 1798, ges. 1777) (E704: ``Weitere Untersuchung über nach Vielfachen eines gewissen Winkels fortschreitende Reihen") diese Frage nach den Koeffizienten wieder auf und gelangt zur Formel für die Koeffizienten. In der ersten Arbeit nähert er sich dem Ausdruck an\footnote{Die Darstellungen und Ausführungen Eulers bezüglich der Entwicklungskoeffizienten lassen sich als Riemann'sche Summe für ein Integral deuten, womit er aus moderner Sicht betrachtet, dem Ziel der Integraldarstellung für die Fourier--Koeffizienten schon sehr nahe kommt.}, wohingegen man in § 4 der zweiten die moderne, oft Fourier zugeschriebene, Formel für die Fourierkoeffizienten entdeckt. In moderner Formulierung lautet sie:

\begin{equation*}
    a_0= \dfrac{1}{\pi}\int\limits_{-\pi}^{\pi}f(x)dx \quad \text{und} \quad a_n= \dfrac{2}{\pi}\int\limits_{-\pi}^{\pi} \cos(nx)f(x)dx \quad \text{für} \quad n\geq 1.
\end{equation*}

\begin{figure}
    \centering
     \includegraphics[scale=0.8]{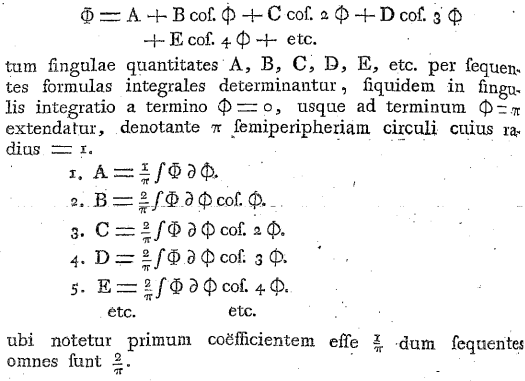}
    \caption{Euler gelangt in seiner Abhandlung \cite{E704} zur heute immer noch gängigen Form der Fourierkoeffizienten.}
    \label{fig:E704Fourierkoeffizienten}
\end{figure}

Unabhängig von seinen Untersuchungen zu Differenzengleichungen findet er aus physikalischen Prinzipien heraus in der Abhandlung \textit{``De vibratione chordarum exercitatio"} (\cite{E119}, 1749, ges. 1748) (E119: ``Über die Vibration einer zum Schwingen angeregten Saite") zur schwingenden Saite erneut eine Fourierreihenentwicklung, widmet sich dort aber nicht der Bestimmung der Koeffizienten, sodass Euler über diesen ersten Berührpunkt mit der Fourier'schen Analysis nicht hinausgeht.


\subsubsection{Produktdarstellung für die Gamma-Funktion}
\label{subsubsec: Produktdarstellung für die Gamma-Funktion}

\epigraph{What is easiest to see is often overlooked.}{Milton Hyland Erickson}


 Ein mittlerweile wohlbekanntes Beispiel für eine Euler'sche Erstentdeckung trotz Benennung nach einem seiner Nachfolger mag die Gauß'sche Produktformel für die Fakultät sein:

\begin{equation}
    \label{eq: Euler Produkt Gamma} 
    \Gamma(x+1)= \prod_{k=1}^{\infty} \dfrac{k^{1-x}\cdot (k+1)^x }{k+x},
\end{equation}
welche sich bei Euler bereits in § 1 seiner Arbeit \cite{E19} findet. Gauß, welcher sie -- wohl unabhängig von Euler -- gefunden hat, gibt sie in seiner Arbeit \textit{``Disquisitiones generales circa Seriem infinitam $1+\frac{\alpha \beta}{1\cdot \gamma}x+\frac{\alpha (\alpha +1)\beta (\beta+1)}{1\cdot 2 \cdot \gamma (\gamma +1)}xx+\frac{\alpha (\alpha +1)(\alpha+2)\beta (\beta +1)(\beta +2)}{1\cdot 2 \cdot 3 \cdot \gamma (\gamma +1)(\gamma +2)}x^3+\text{etc.}$ pars prior"} (\cite{Ga13}, 1813, ges. 1812) (``Allgemeine Untersuchungen über die unendliche Reihe $1+\frac{\alpha \beta}{1\cdot \gamma}x+\frac{\alpha (\alpha +1)\beta (\beta+1)}{1\cdot 2 \cdot \gamma (\gamma +1)}xx+\frac{\alpha (\alpha +1)(\alpha+2)\beta (\beta +1)(\beta +2)}{1\cdot 2 \cdot 3 \cdot \gamma (\gamma +1)(\gamma +2)}x^3+\text{etc.}$ -- erster Teil") mit einem strengen Beweis an, wohingegen Euler sie lediglich in prosaischer Form plausibel macht\footnote{Die Formel (\ref{eq: Euler Produkt Gamma}) wird auch von Weierstraß in seiner Arbeit \textit{``Uber die Theorie der analytischen Facultäten"} (\cite{We54}) studiert. Während Weierstraß Gauß Werk \cite{Ga13} noch erwähnt, lässt er Euler dieses Privileg nicht zukommen.}. Seine Arbeit beginnt nämlich wie folgt:\\

\textit{``Nachdem ich neulich bei Gelegenheit dessen, was  der hochgeehrte Herr Goldbach mit der Sozietät über Reihen geteilt hat, nach einem allgemeinen Ausdruck gesucht habe, welcher alle Terme dieser Progression gäbe}

\begin{equation*}
    1+1\cdot 2+ 1 \cdot 2 \cdot 3+ 1 \cdot 2 \cdot 3 \cdot 4 +\text{etc.},
\end{equation*}

\textit{bin ich, bedenkend, dass sie ins Unendliche fortgesetzt schließlich mit der geometrischen zusammenfällt, auf den folgenden Ausdruck gestoßen}

\begin{equation*}
    \dfrac{1 \cdot 2^n}{1+n} \cdot \dfrac{2^{1-n} \cdot 3^n}{2+n} \cdot \dfrac{3^{1-n} \cdot 4^n}{3+n} \cdot \dfrac{4^{1-n} \cdot 5^n}{4+n}\cdot \text{etc.}
\end{equation*}

\textit{welcher den Term der Ordnung $n$ von besagter Progression darstellt."}\\

Diese Erläuterung findet man erst in der Arbeit \textit{``De termino generali serierum hypergeometricarum"} (\cite{E652}, 1793, ges. 1777) (E652: ``Über den allgemeinen Term von hypergeometrischen Reihen") in einen mathematischen Beweis überführt.\\

\begin{figure}
    \centering
    \includegraphics[scale=1.0]{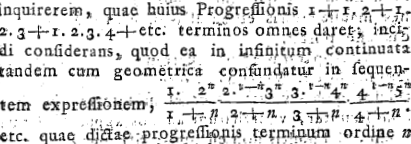}
    \caption{Euler gibt in seiner Arbeit \cite{E19} die Produktdarstellung der $\Gamma$--Funktion bzw. der Fakultät.}
    \label{fig:E19EulerGammaProdukt}
\end{figure}

Weniger bekannt mag die Feststellung sein, dass sich die oft Weierstraß zugeschriebene Formel

\begin{equation}
    \label{eq: Gamma Weierstraß}
    {\Gamma(z)}= \dfrac{e^{-\gamma z}}{z}\prod_{n=1}^{\infty} \left(1+\dfrac{z}{n}\right)^{-1}e^{\frac{z}{n}}
\end{equation}
bereits ohne große Mühen aus den  in dem Artikel  \cite{E613} angestellten Euler'schen Überlegungen ergibt; explizit findet man einen entsprechenden Ausdruck überdies in seinen \textit{Calculi Differentialis} \cite{E212}. In Abschnitt (\ref{subsubsec: Ein anderes Vorhaben: Das Weierstraß-Produkt}) wird dies noch einmal zur Sprache kommen. \\

Solche Prioritätsfragen sollen in der gegenwärtigen Arbeit jedoch lediglich eine Randnotiz bleiben, zumal sie  vergleichsweise leicht aufzuklären sind.

\subsection{Von Euler in anderer Gestalt Bewiesenes}
\label{subsec: Von Euler in anderer Gestalt Bewiesenes}

\epigraph{What does it matter that we take different roads so long as we reach the same goal?}{Mahatma Gandhi}

Nun werden Entdeckungen diskutiert werden, welche Euler gemacht hat, die jedoch in anderer, teils unerwarteter Gestalt in seinen Werken zum Vorschein treten. Das Maß des Nicht--Erwartens korreliert dabei positiv mit den beizufügenden Hintergrundinformationen. Es wird mit der Multiplikationsformel für die $\Gamma$--Funktion begonnen werden (Abschnitt \ref{subsubsec: Die Multiplikationsformel für die Gammafunktion}), gefolgt vom Weierstraß'schen Produktsatz für Funktionen mit vorgegebenen Nullstellen (Abschnitt \ref{subsubsec: Ein anderes Vorhaben: Das Weierstraß-Produkt}) und der Mellin--Transformation (Abschnitt \ref{subsubsec: Die Mellin--Transformierte bei Euler}). Während diese Gegenstände sich noch an einzelnen Euler'schen Ausarbeitungen erklären lassen, sind für die Diskussion von Eulers Beiträgen zu den Legendre'schen Polynomen (Abschnitt \ref{subsubsec: Den Kontext betreffend: Die Legendre Polynome}) und die hypergeometrische Funktion (Abschnitt \ref{subsubsec: Die Darstellung betreffend -- Die hypergeometrische Reihe}) mehrere seiner Abhandlung in die Betrachtung mit einzubeziehen.

\subsubsection{Die Form betreffend: Die Multiplikationsformel für die Gammafunktion}
\label{subsubsec: Die Multiplikationsformel für die Gammafunktion}


\epigraph{ Always there will be someone who can tell it [the story] differently depending on where they are standing [...] this is the way I think the world's stories should be told: from many different perspectives.}{Chinua Achebe}

Die Stelle des ersten Beispiels nimmt die Multiplikationsformel für die $\Gamma$--Funktion ein, welche sich ebenfalls in Eulers Werk findet, jedoch erst erkannt werden kann, wenn der aus der Euler'schen Symbolik gebildete Schleier gelüftet wird. In moderner Sprache ausgedrückt ist die folgende Formel Fokus des Interesses:

\begin{equation}
\label{eq: Multiplikationsformel Gamma}
\Gamma \left(\dfrac{x}{n}\right)\Gamma \left(\dfrac{x+1}{n}\right)\cdots \Gamma \left(\dfrac{x+n-1}{n}\right)= \dfrac{(2\pi)^{\frac{n-1}{2}}}{n^{x-\frac{1}{2}}}\cdot \Gamma(x).
\end{equation}
Sie wurde in dieser Form von Gauß in seiner Abhandlung \cite{Ga13}  in der obigen Form bewiesen\footnote{Gauß verwendet in seiner Abhandlung das Symbol $\Pi$, womit er die Fakultät anzeigt, statt des von Legendre eingeführten und heute gebräuchlichen Symbols $\Gamma$.}.

\paragraph{Herleitung der Multiplikationsformel aus Eulers Formeln}

 Die Multiplikationsformel wird aus Eulers Formeln nach dem Vorbild der Arbeit \textit{``Euler and the multiplication formula for the Gamma Function"} (\cite{Ay21}, 2021) abgeleitet; sie findet sich in seinem Papier \textit{``Evolutio formulae integralis $\int x^{f-1} dx \log ^{\frac{m}{n}}(x)$ integratione a valore $x=0$ ad $x=1$ extensa"} (\cite{E421}, 1772, ges. 1771) (E421: ``Entwicklung der Integralformel $\int x^{f-1} dx \log ^{\frac{m}{n}}(x)$, wobei die Integration vom Wert $x=0$ bis hin zu $x=1$ erstreckt worden ist.")  in der Form

\begin{equation}
\label{eq: Euler Mult Gamma}
\left[ \frac{m}{n} \right] = \frac{m}{n} \sqrt[n]{n^{n-m}\cdot 1 \cdot 2  \cdots (m-1) \left(\frac{1}{m}\right)\left(\frac{2}{m}\right)\cdots \left(\frac{n-1}{m}\right)}.
\end{equation}
In dieser Darstellung ist die Äquivalenz zu (\ref{eq: Multiplikationsformel Gamma}) nicht unmittelbar evident. \\

\begin{figure}
    \centering
     \includegraphics[scale=1.1]{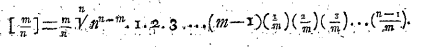}
    \caption{ Eulers Version der Multiplikationsformel für die $\Gamma$-Funktion aus seiner Arbeit \cite{E421} mit von ihm eigens eingeführten Symbolen.}
    \label{fig:E421MultiplicationGamma}
\end{figure}

Eine Erklärung der von Euler gebrauchten Formelzeichen ist vorauszuschicken. Zunächst definiert Euler in § 44 den Ausdruck

\begin{equation*}
\left(\frac{p}{q}\right):= \int\limits_{0}^{1} \dfrac{x^{p-1}dx}{(1-x^n)^{\frac{n-q}{n}}}.
\end{equation*}
Dabei sollen $p$, $q$, $n$  natürliche Zahlen sein. Man beachte, dass der Buchstabe $n$ im Euler'schen Symbol nicht auftritt\footnote{Euler untersucht die durch das Symbol $\left(\frac{p}{q}\right)$ definierten Funktionen ausgiebig in der Abhandlung \textit{``Observationes circa integralia formularum $\int x^{p-1}dx (1-x^n)^{\frac{q}{n}-1}$ posito post integrationem $x = 1$"} (\cite{E321}, 1766, ges. 1765)  (E321: ``Beobachtungen zu den Integralformeln $\int x^{p-1}dx (1-x^n)^{\frac{q}{n}-1}$, wenn nach der Integration $x = 1$ gesetzt wird").}, für welchen man heute den folgenden Ausdruck bevorzugt

\begin{equation}
\label{eq: Def EulerBeta}
\left(\frac{p}{q}\right) = \frac{1}{n} \int\limits_{0}^1 y^{\frac{p}{n}-1}dy(1-y)^{\frac{q}{n}-1}=\frac{1}{n} \cdot B \left(\frac{p}{n},\frac{q}{n}\right),
\end{equation}
wobei 

\begin{equation*}
B(x,y)= \int\limits_{0}^1 t^{x-1}dt(1-t)^{y-1} \quad \text{für} \quad \text{Re}(x),\text{Re}(y) > 0
\end{equation*}
die Beta--Funktion bzw. das Beta--Integral bedeutet. Mit dem Symbol $[\lambda]$ meint Euler $\lambda!$, sodass $[\lambda]=\Gamma(\lambda+1)$. Daher findet man bei Euler auch die Spiegelungsformel in der Gestalt

\[ 
[\lambda]\cdot [-\lambda]=\frac{\pi \lambda}{\sin \pi \lambda}. 
\]
Die Formel findet sich in § 43 von \cite{E421} und wird heute als

\begin{equation*}
    \Gamma(x)\Gamma(1-x)=\dfrac{\pi}{\sin (\pi x)}
\end{equation*}
geschrieben. Unter Verwendung  dieser Formel für $x=\frac{i}{n}$, für $i=1,2,3\cdots, n$, findet man

$$
\begin{array}{rcl}
\Gamma \left(\dfrac{1}{n}\right)  \Gamma \left(\dfrac{n-1}{n}\right)&=&\dfrac{\pi}{\sin \frac{\pi}{n}},\\
\Gamma \left(\dfrac{2}{n}\right)  \Gamma \left(\dfrac{n-2}{n}\right)&=&\dfrac{\pi}{\sin \frac{2\pi}{n}},\\
\Gamma \left(\dfrac{3}{n}\right)  \Gamma \left(\dfrac{n-3}{n}\right)&=&\dfrac{\pi}{\sin \frac{3\pi}{n}},\\
\ldots&=&\ldots\\
\Gamma \left(\dfrac{n-1}{n}\right)  \Gamma \left(\dfrac{1}{n}\right)&=&\dfrac{\pi}{\sin \frac{(n-1)\pi}{n}} .\\
\end{array}
$$
Das Produkt all dieser gibt demnach

\begin{equation}
\label{eq: Prod Gamma (i/n)}
 \prod_{i=1}^{n-1}\Gamma \left(\frac{i}{n}\right)^2= \frac{\pi^{n-1}}{\prod_{i=1}^{n-1} \sin \left(\frac{i \pi}{n}\right)} .
\end{equation}
Vermöge der Formel

\begin{equation*}
\prod_{i=1}^{n-1} \sin \left(\frac{i \pi}{n} \right) = \frac{n}{2^{n-1}},
\end{equation*}
welche sich zum Beispiel aus der allgemeineren aus § 240 der \textit{Introductio} \cite{E101}

\begin{equation*}
\sin n \varphi = 2^{n-1} \sin \varphi \sin \left(\dfrac{\pi}{n}- \varphi\right) \sin \left(\dfrac{\pi}{n}+ \varphi\right)
\end{equation*}
\begin{equation*}
 \sin \left(\dfrac{2\pi}{n}- \varphi\right) \sin \left(\dfrac{2\pi}{n}+ \varphi\right)\cdot\text{etc.}
\end{equation*}
ableitet, lässt sich demnach (\ref{eq: Prod Gamma (i/n)}) geschlossen auswerten:

\begin{equation}
\label{eq: Prod Gamma (i/n) ausgewertet}
\Gamma \left(\frac{1}{n}\right) \Gamma \left( \frac{2}{n}\right) \cdots \Gamma \left(\frac{n-1}{n}\right)=\sqrt{\frac{(2\pi)^{n-1}}{n}}.
\end{equation}

Um nun (\ref{eq: Multiplikationsformel Gamma}) aus (\ref{eq: Euler Mult Gamma}) abzuleiten, bringe man letztere schrittweise in die Form der ersten. Erreicht wird dies durch Ersetzten  der Formeln $\left(\frac{p}{q}\right)$ mit (\ref{eq: Def EulerBeta}) in Eulers Version, Schreiben von $\Gamma(m)$ statt $1\cdot 2 \cdot 3 \cdots (m-1)$ und entsprechend $\Gamma\left(\frac{m}{n}+1\right)$ statt $[\frac{m}{n}]$.  Der Übersichtlichkeit wegen ersetze man noch $m$ durch $x$, sodass gilt

\begin{equation*}
\Gamma \left(\frac{x}{n}\right)= \sqrt[n]{n^{n-x} \Gamma(x) \frac{1}{n^{n-1}}B \left(\frac{1}{n},\frac{x}{n}\right)B \left(\frac{2}{n},\frac{x}{n}\right)\cdots B \left(\frac{n-1}{n},\frac{x}{n}\right)}.
\end{equation*}
Mit der bekannten Formel

\begin{equation}
    \label{eq: Gamma und Beta}
    B(x,y)=\dfrac{\Gamma(x)\Gamma(y)}{\Gamma(x+y)}
\end{equation}
lassen sich  die Beta--Integrale ersetzen, sodass

\begin{equation*}
\Gamma \left(\frac{x}{n}\right)= \sqrt[n]{n^{1-x}\Gamma(x) \dfrac{\Gamma\left(\frac{1}{n}\right)\Gamma \left(\frac{x}{n}\right)}{\Gamma \left(\frac{x+1}{n}\right)} \cdot \dfrac{\Gamma\left(\frac{2}{n}\right)\Gamma \left(\frac{x}{n}\right)}{\Gamma \left(\frac{x+2}{n}\right)} \cdots \dfrac{\Gamma\left(\frac{n-1}{n}\right)\Gamma \left(\frac{x}{n}\right)}{\Gamma \left(\frac{x+n-1}{n}\right)}}.
\end{equation*}
Hinüberschaffen aller $\Gamma$--Funktionen mit gebrochenem Argument auf die linke Seite gibt

\begin{small}
\begin{equation*}
\Gamma \left(\frac{x}{n}\right)\Gamma \left(\frac{x+1}{n}\right) \Gamma \left(\frac{x+2}{n}\right) \cdots \Gamma \left(\frac{x+n-1}{n}\right)= n^{1-x} \Gamma (x) \Gamma \left(\frac{1}{n}\right) \cdots \Gamma \left(\frac{n-1}{n}\right).
\end{equation*}
\end{small}
Mit (\ref{eq: Prod Gamma (i/n) ausgewertet}) kontrahiert dies zu

\begin{equation*}
\Gamma \left(\frac{x}{n}\right)\Gamma \left(\frac{x+1}{n}\right) \Gamma \left(\frac{x+2}{n}\right) \cdots \Gamma \left(\frac{x+n-1}{n}\right)=  n^{1-x} \Gamma(x) \sqrt{\dfrac{(2\pi)^{n-1}}{n}},
\end{equation*}
was  gerade die Gauß'sche Multiplikationsformel (\ref{eq: Multiplikationsformel Gamma}) ist.

\paragraph{Gründe für die Euler'sche Darstellungsweise}

Es bleibt die Frage, was Euler bewogen haben mag, nicht dieselbe Darstellungsweise wie Gauß gewählt zu haben. Für eine Beantwortung dieser ist es ratsam, die leichter zu beantwortende Frage nach Gauß' Gründen für eine Darstellung wie (\ref{eq: Multiplikationsformel Gamma}) zu bevorzugen, voranzustellen.  In seiner Abhandlung \cite{Ga13} widmet er einen eigenen Abschnitt der verallgemeinerten Fakultät und findet verschiedene Ausdrücke zur Interpolation\footnote{Euler hat über verschiedene Werke hinweg all diese Formeln ebenfalls, jedoch teilweise in völlig anderen Zusammenhängen, gefunden. Man vergleiche etwa \cite{Ay21}.}. Dementsprechend war er  auf der Suche nach Eigenschaften der $\Gamma$--Funktion und (\ref{eq: Euler Mult Gamma}) ist eine Funktionalgleichung, welcher sie genügt. Die Euler'sche Version (\ref{eq: Euler Mult Gamma}) hingegen verschleiert dies, sodass Gauß wohl allein schon dieses Grundes wegen die heute nach ihm benannte Darstellung bevorzugt hätte. \\

Euler hingegen intendiert, kompliziertere Funktionen durch einfachere auszudrücken, wobei ``kompliziert"{} in diesem Zusammenhang als Maß der Transzendenz verstanden werden kann\footnote{So sind die algebraischen Funktionen $x^n$ und $\sqrt[n]{x}$ weniger kompliziert als die elementaren Transzendenten $e^x$, $\sin x$ und $\log x$, $\arcsin (x)$, welche wiederum weniger kompliziert sind als elliptische Integrale usw.}. Gleichermaßen verhält es sich in \cite{E421}. Man beachte nämlich, dass in (\ref{eq: Euler Mult Gamma}) das Symbol linker Hand durch das Integral im Titel ausgedrückt wird:

\begin{equation*}
    \left[\dfrac{m}{n}\right]=\int\limits_{0}^1 \left(\log\left(\dfrac{1}{x}\right)\right)^{\frac{m}{n}}dx,
\end{equation*}
welches ein Integral über eine \textit{transzendente} Funktion ist. Alle Ausdrücke auf der rechten Seite von (\ref{eq: Euler Mult Gamma}) sind hingegen \textit{algebraische} Funktionen, sind also insbesondere durch Quadraturen leichter zugänglich\footnote{Es handelt sich gar um Perioden im Sinn der Arbeit \textit{``Periods"} (\cite{Ko01}, 2001) von Kontsevich (1964--) und Zagier (1951--).}. Diesbezüglich schreibt inbesondere Euler in § 1 seiner Arbeit \textit{``De constructione aequationum"} (\cite{E70}, 1744, ges. 1737) (``Über die Konstruktion von Gleichungen") zu einem vergleichbaren Gegenstand:\\

\textit{``Sooft in der Auflösung von Problemen zu Differentialgleichungen gelangt wird, ist vor allem zu untersuchen, ob diese Gleichungen eine Integration zulassen; denn ein Problem ist als vollständigst gelöst anzusehen, welches auf die Konstruktion einer algebraischen Gleichung reduziert worden ist. Aber wenn die Gleichung, was auch sehr oft geschieht, in keiner Weise in eine algebraische Form überführt werden kann, dann müssen entweder Quadraturen oder Rektifikationen von Kurven, deren Konstruktion man gewahr ist, herangezogen werden."}\\

Demnach hätte Euler vermutlich  (\ref{eq: Multiplikationsformel Gamma}) nicht seiner Version vorgezogen. Dass er die Gauß'sche Version zumindest hätte ableiten können, ist aus der vorgestellten Rechnung ersichtlich, zumal der Beweis keine Formeln und Konzepte nutzt, die Euler unbekannt waren.

\subsubsection{Eine andere Intention: Euler und das Weierstraß'sche Produkt}
\label{subsubsec: Ein anderes Vorhaben: Das Weierstraß-Produkt}

\epigraph{One can measure the importance of a scientific work by the number of earlier publications rendered superfluous by it.}{David Hilbert}

Der Weierstraß'sche Produktsatz lehrt,  ob eine  Funktion mit vorgegebenen Nullstellen existiert und gibt im zweiten Schritt -- gegebenenfalls unter Hinzunahme konvergenzerzeugenden Faktoren --   eine Möglichkeit zu ihrer Konstruktion an. Euler gelangt, wie auch im Vorwort zu Band 16,2 von Serie 1 der \textit{Opera Omnia} (\cite{OO162}, 1935) erwähnt wird, ebenfalls zu einem Ausdruck, welcher als Weierstraß--Produkt aufgefasst werden kann. Man findet seine Ideen hierzu in Kapitel 17 des zweiten Teils seiner \textit{Calculi Differentialis} \cite{E212} und in ausführlicherer Darstellung in \cite{E613}. Der Euler'sche Weg ist jedoch vom Weierstraß'schen verschieden.

\paragraph{Eulers Ansatz am Beispiel der harmonischen Reihe}

Euler intendiert die Interpolation der Summe

\begin{equation}
\label{eq: Euler zu Weierstraß}
    F(x):= \sum_{k=1}^x f(k), 
\end{equation}
zu welchem Zweck er folgende trivial anmutende Identität anmerkt\footnote{Dieses geschickte Addieren von $0$ wurde auch von Ramanujan (1887--1920) mit Vorliebe benutzt, siehe etwa das Buch \cite{Be85}.}:

\begin{equation*}
    \renewcommand{\arraystretch}{2,0}
\setlength{\arraycolsep}{0.0mm}
\begin{array}{rclclclclclclcl}
     F(x) &~=~&~ f(1) &~+~&~ f(2)  & ~+~ &~ \cdots &~+~&~  f(x) &~+~&~ f(x+1) &~+~&~ \cdots \\
          &~-~&~ f(x+1) &~-~&~ f(x+2)  & ~-~ &~ \cdots \\
\end{array}
\end{equation*}

Nach Euler'scher Maxime\footnote{Euler verwendet an vielerlei Stellen eine Version des Ausspruchs: \textit{``Man lehrt am besten anhand von Beispielen."}} soll das Gesagte  an einem Beispiel demonstriert werden: In § 17 wählt Euler $f(x)=\frac{1}{x}$\footnote{Dieses Beispiel findet sich auch in Sandifers Artikel \textit{``Inexplicable functions"} (\cite{Sa07nov}, 2007) in anderem Kontext besprochen.}:

\begin{equation*}
    \renewcommand{\arraystretch}{2,0}
\setlength{\arraycolsep}{0.0mm}
\begin{array}{rclclclclclclcl}
     F(x) &~=~&~ \dfrac{1}{1} &~+~&~  \dfrac{1}{2} &~+~&~ \dfrac{1}{3} & ~+~ &~ \cdots &~+~&~  \dfrac{1}{x} &~+~&~ \dfrac{1}{x+1} &~+~& \cdots \\
          &~-~&~ \dfrac{1}{x+1} &~-~&~ \dfrac{1}{x+2} &~-~&~ \dfrac{1}{x+3} & ~-~ &~ \cdots \\
\end{array}
\end{equation*}
Spaltenweises Zusammenfassen der Terme liefert:

\begin{equation*}
    F(x)= \dfrac{x}{x+1}+ \dfrac{x}{2(x+2)}+ \dfrac{x}{3(x+3)}+\text{etc.}
\end{equation*}
Hieraus leitet Euler dann nach einiger Rechnung durch Integration die Potenzreihe für $\int F(x)dx$ her und gelangt zur Reihenentwicklung (§ 23):

\begin{equation*}
    \int F(x)dx = \dfrac{\zeta(2)}{2}x^2- \dfrac{\zeta(3)}{3}x^3+\dfrac{\zeta(4)}{4}x^4-\text{etc.},
\end{equation*}
wobei

\begin{equation*}
    \zeta(s):= \sum_{n=1}^{\infty} \dfrac{1}{n^s}
\end{equation*}
bedeutet. Eulers Wahl $f(x)=\frac{1}{x}$ interpoliert  die harmonische Reihe, sodass die Euler'sche Summation der Funktion

\begin{equation*}
    \psi(x):= \dfrac{d}{dx} \log (\Gamma(x))
\end{equation*}
sehr nahe kommt\footnote{In \cite{E613} drückt Euler mit der letzten Formel nicht ganz diese Funktion aus, da er schlicht termweise integriert und die Integrationskonstanten ignoriert. Indes gelangt er in § 384 von Teil 2 seiner \textit{Calculi Differentialis} \cite{E212} auf gleichem Wege zur Potenzreihenentwicklung von $\psi$.} \\

\begin{figure}
    \centering
     \includegraphics[scale=1.2]{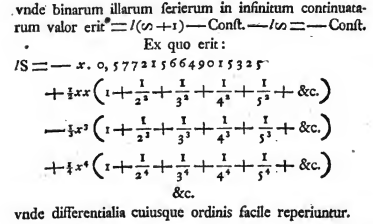}
    \caption{ Euler gelangt in \cite{E212} zur Potenzreihe der $\psi$-Funktion. Die Euler--Mascheroni--Konstante $\gamma$ ist hierbei für die ersten Stellen ausgeschrieben.  }
    \label{fig:E212DigammaFunktion}
\end{figure}

Er findet, modern ausgedrückt:

\begin{equation*}
    \psi(x+1)= -\gamma x +\dfrac{\zeta(2)}{2}x^2-\dfrac{\zeta(3)}{3}x^3+\dfrac{\zeta(4)}{4}x^4-\text{etc.}
\end{equation*}
mit der Euler--Mascheroni--Konstante\footnote{Euler führt diese Konstante zum ersten Male in § 11 seiner Arbeit \textit{``De progressionibus harmonicis observationes"} (\cite{E43}, 1740, ges. 1734) (E43: ``Bemerkungen zu harmonischen Progressionen") ein.}

\begin{equation*}
    \gamma:= \lim_{n \rightarrow \infty} \left(\sum_{k=1}^n \dfrac{1}{k}- \log (n)\right).
\end{equation*}
Nach der Definition von $\psi(x)$ gelangt man durch Integration erneut zum ``Weierstraß'{}schen"{} Produkt (\ref{eq: Gamma Weierstraß}). Diese Euler'sche Vorgehensweise nimmt gleichsam die Idee der konvergenzerzeugenden Faktoren vorweg, zumal das geschickte Addieren von $0$ gestattet, zwei eigentlich divergente Reihen voneinander in solcher Weise zu subtrahieren, dass sie eine konvergente Reihe ergeben. 

\paragraph{Behandlung von Produkten am Beispiel der Fakultät}

Seine Ideen will Euler auch auf kompliziertere Reihen übertragen. Er bemerkt, dass je nach Wahl der Funktion $f(x)$ in (\ref{eq: Euler zu Weierstraß}) dazu gegebenenfalls mehr als eine Differenz gebildet werden muss, um die Konvergenz zu erzwingen. Die Behandlung von Produkten dieser Gestalt lässt sich durch Betrachtung der Logarithmen derselben mit derselben Methoden behandeln und bildet daher ``lediglich"{}  ein Supplement zu \cite{E613}.\\

\begin{figure}
    \centering
    \includegraphics[scale=0.9]{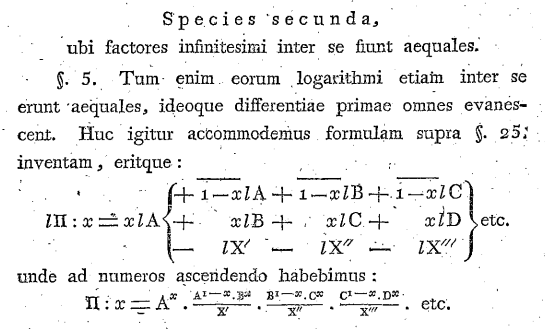}
    \caption{Euler erklärt in seiner Arbeit \cite{E613} die Interpolation von Produktformeln.}
    \label{fig:E613EulerWeierstraß}
\end{figure}

Aber auch die Behandlung der Produkte soll an einem Beispiel illustriert werden. Es sei also die $\Gamma$--Funktion vorgelegt, als ein unendliches Produkt darzustellen. Wegen der Funktionalgleichung $\Gamma(x+1)= x\Gamma(x)$ gilt auch $\log \Gamma(x+1)-\log \Gamma(x)= \log x$, sodass zunächst für natürliches $x$:

\begin{equation*}
    \log \Gamma(x+1) = \sum_{k=1}^{x} \log (k).
\end{equation*}
Zumal $\log(k+1)-\log(k)= \log \left(1+\frac{1}{k}\right)$ für $k \rightarrow \infty$ zu $0$ strebt, gehört diese Reihe zur zweiten Klasse, welche Euler in \cite{E613} definiert\footnote{Die erste Klasse wird dabei von Reihen gebildet, deren Terme mit wachsendem $k$ bereits selbst gegen $0$ streben. Die Reihen der ersten Klasse bedürfen einer Differenzenbildung, welche Differenz dann gegen $0$ strebt, die der zweiten Klasse hingegen zwei usw.}. Eulers Methode lehrt, dass die letzte Summe wie folgt geschrieben werden kann:

\begin{equation*}
      \renewcommand{\arraystretch}{1,5}
\setlength{\arraycolsep}{0.0mm}
\begin{array}{Rccccccccccccccccccccccc}
    \sum_{k=1}^x \log (k) &~=~ & (1-x)\log(1) &~+~&  (1-x)\log (2) &~+~& (1-x)\log(3) &~+~& \text{etc.}     \\
 +   x\log(1)   & ~+~ & x\log(2)  & ~+~ & x\log(3)  & ~+~ & x\log(4)  & ~+~ & \text{etc.} \\
           & ~-~ & \log(x+1)  & ~-~ & \log(x+2)  & ~-~ & \log(x+3)  & ~-~ & \text{etc.}   
\end{array}
\end{equation*}
Spaltenweises Summieren liefert:

\begin{equation*}
    \log (\Gamma(x+1))= \log (1^{1-x}) + \log \left(\dfrac{1^{1-x}2^x}{x+1}\right) + \log \left(\dfrac{2^{1-x}3^x}{x+2}\right) + \log \left(\dfrac{3^{1-x}4^x}{x+3}\right) + \text{etc.}
\end{equation*}
Nehmen der Exponentialfunktionen gibt die Formel für $x!$, welche Euler bereits in der Arbeiz \cite{E19} mitgeteilt hat, nämlich:

\begin{equation*}
    x! = \Gamma(x+1) = 1^{1-x} \cdot \dfrac{1^{1-x}2^x}{x+1} \cdot \dfrac{2^{1-x}3^x}{x+2} \cdot \dfrac{3^{1-x}4^x}{x+3} \cdot \text{etc.} = \prod_{k=1}^{\infty} \dfrac{k^{1-x}(k+1)^x}{x+k}
\end{equation*}

\paragraph{Ein Vergleich zum Weierstraß'schen Produktsatz}

Ein abschließender Vergleich zum Weierstraß'schen Produktsatz aus seinen Vorlesungen (\cite{We78}, 1878) und den Euler'schen Unternehmungen in \cite{E613} wird dem Verständnis förderlich sein, wobei neben den Ausführungen in \cite{Ay21a} die Ausführungen von Faber (1877--1966) im Vorwort von Band 16,2 der Serie 1 der \textit{Opera Omnia} das Gesagte treffend zusammenfassen. Faber schreibt auf Seite XLIII:\\

    \textit{``Tatsächlich hat Euler nicht nur die Produktdarstellung (12)\footnote{Dies ist die Produktdarstellung $ x! = \Gamma(x+1) = 1^{1-x} \cdot \dfrac{1^{1-x}2^x}{x+1} \cdot \dfrac{2^{1-x}3^x}{x+2} \cdot \dfrac{3^{1-x}4^x}{x+3} \cdot \text{etc.} = \prod_{k=1}^{\infty} \dfrac{k^{1-x}(k+1)^x}{x+k}$, welche im Beispiel diskutiert wurde.}, sondern sogar den Gedanken der Konvergenz erzeugenden Faktoren von Weierstra\ss{} vorweggenommen. Denn es bedeutet keinen Unterschied, ob man den Gliedern des divergenten Produktes $\prod_{\nu=1}^{\infty}\left(1+\frac{x}{\nu}\right)$ die Konvergenz erzeugenden Faktoren $e^{-\frac{x}{\nu}}$ oder den Gliedern der divergenten unendlichen Reihe $\sum_{\nu=1}^{\infty}\log \left(1+\frac{x}{\nu}\right)$ die Konvergenz erzeugenden Summanden $-\frac{x}{\nu}$ oder auch $-x \log \left(1+\frac{1}{\nu}\right)$ beifügt. Das tat aber Euler mit voller Absicht in der Abhandlung 613."}

\subsubsection{Eine andere Fragestellung: Die Mellin--Transformierte}
\label{subsubsec: Die Mellin--Transformierte bei Euler}

\epigraph{Over and over again, scientific discoveries have provided answers to problems that had no apparent connection with the phenomena that gave rise to the discovery.}{Isaac Asimov}

In der Arbeit \textit{``Euler and Homogeneous Difference Equations with Linear Coefficients"} (\cite{Ay24a}, 2024) wird darauf hingewiesen, dass einige von  Eulers Abhandlungen, die ihren Ursprung in der Konstruktion und Evaluation von Kettenbrüchen haben, zu einer Methode führen,  zu einer gegebenen Mellin--Transformierten -- ein Begriff, welcher erst weit nach Eulers Tod von Mellin (1854--1933) in seiner Arbeit  \textit{``Über die fundamentale Wichtigkeit des Satzes von Cauchy f\"ur die Theorien der Gamma- und hypergeometrischen Functionen"} (\cite{Me95}, 1896) allgemein eingeführt wurde -- ihre Inverse vermöge eines Systems von Differentiagleichungen, anstatt des heute üblichen Weges über Integration in der komplexen Ebene, zu finden. Dies soll hier noch einmal beleuchtet werden.

\paragraph{Eulers Zugang}
\label{para: Eulers Zugang}

In seinen Ausarbeitungen  \textit{``De fractionibus continuis observationes"} (\cite{E123}, 1750, ges. 1739) (E123: ``Bemerkungen zu Kettenbrüchen")  sowie  \textit{``Methodus inveniendi formulas integrales, quae certis casibus datam inter se teneant rationem, ubi simul methodus traditur fractiones continuas summandi"} (\cite{E594}, 1785, ges. 1775) (E594: ``Eine Methode Integralformeln zu finden, welche in gewissen Fällen ein gegebenes Verhältnis zueinander haben, wo zugleich eine Methode angegeben wird, Kettenbrüche zu summieren") präsentiert Euler eine Methode zur Lösung der Differenzengleichung

\begin{equation}
    \label{eq: Linear Difference Equation}
    (a_0x+b_0)f(x)+(a_1x+b_1)f(x+1)+ \cdots + (a_n x+ b_n) f(x+n) =0
\end{equation}
mit komplexen Koeffizienten $a_0$, $a_1$, $\cdots$, $a_n$ und $b_0$, $b_1$, $\cdots$, $b_n$, wobei mindestens einer der Koeffizienten $a_0$, $a_1$, $\cdots$, $a_n$ nicht verschwinden soll. \\

\begin{figure}
    \centering
    \includegraphics[scale=0.9]{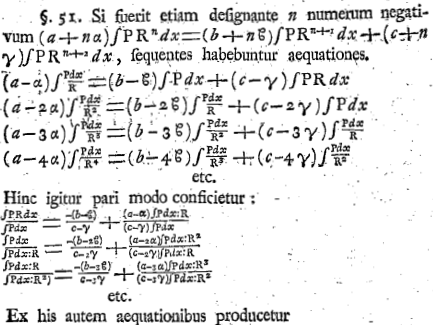}
    \caption{Euler erklärt in seiner Arbeit \cite{E123} die Verwendung einer Differenzengleichung zur Bildung von Kettenbrüchen.}
    \label{fig:E123Differenzengleichung}
\end{figure}

In seinen Arbeiten ist Euler an Kettenbrüchen interessiert\footnote{Der Euler'sche Beitrag zu Kettenbrüchen wird in der Arbeit \textit{``Euler und die analytische Theorie der Kettenbrüche"} (\cite{Ay15}, 2015) umfassend diskutiert.}, weshalb er nur den Fall $n=2$ betrachtet, jedoch erlauben seine Ausführungen eine wortgetreue Übertragung auf den allgemeinen Fall. Seine Idee ist die Lösung von (\ref{eq: Linear Difference Equation}) als 

\begin{equation}
    \label{eq: Ansatz}
    f(x) = \int\limits_{a}^{b} t^xp(t)dt
\end{equation}
anzunehmen und die unbekannte Funktion $p(t)$ sowie die Grenzen $a$, $b$ zu ermitteln. Dazu führt er die Hilfsgleichung

\begin{equation}
    \label{eq: Aux}
     (a_0x+b_0)\int t^xp(t)dt + \cdots + (a_n x+ b_n) \int t^{x+n}p(t)dt +t^x q(t) =0,
\end{equation}
mit einer weiteren unbekannten Funktion $q(t)$ ein. Diese Gleichung werde anschließend nach $t$ differenziert, sodass sie gibt:

\begin{equation*}
    (a_0x+b_0)t^xp(t)+ \cdots + (a_n x+ b_n)  t^{x+n}p(t) +xt^{x-1} q(t)+t^xq'(t) =0.
\end{equation*}
Nach Division durch $t^{x-1}\neq 0$ findet man

\begin{equation*}
    (p(t)(a_0 t+\cdots +a_n t^{n+1})+q(t))\cdot x  +  p(t)(b_0 t+\cdots +b_n t^{n+1})+ tq'(t) = 0. 
\end{equation*}
Nun lehrt Euler, die jeweiligen Koeffizienten der Potenzen von $x$ auf den beiden Seiten der Gleichung zu vergleichen, um ein System von Differentialgleichungen für die Funktionen $q$ und $p$ abzuleiten.\\

Hat man diese Funktionen ausfindig gemacht, verbleibt Bestimmung der Grenzen der Integration. Dazu bedient sich Euler der Gleichung:

\begin{equation}
    \label{eq: Boundaries}
    t^xq(t)=0.
\end{equation}
Je zwei Lösungen dieser Gleichung -- sie seien $t_1$ und $t_2$ -- konstituieren eine Lösung der Gleichung (\ref{eq: Linear Difference Equation})  als

\begin{equation}
    \label{eq: Solution}
    f(x) = \int\limits_{t_1}^{t_2}t^x p(t)dt.
\end{equation}
 Euler bedarf in seiner den Kettenbrüchen gewidmeten Untersuchung nur zweier Lösungen von (\ref{eq: Boundaries}); jedoch geben etwa 3 Lösungen  $t_1$, $t_2$, $t_3$   die Integrale

\begin{equation*}
    f_1(x) = \int\limits_{t_1}^{t_2}t^x p(t)dt, \quad   f_2(x) = \int\limits_{t_1}^{t_3}t^x p(t)dt, \quad   f_3(x) = \int\limits_{t_2}^{t_3}t^x p(t)dt
\end{equation*}
und deren Linearkombinationen Lösungen der zu lösenden Gleichung (\ref{eq: Linear Difference Equation}). Die Verallgemeinerung für $n$ Lösungen von (\ref{eq: Boundaries}) ist leicht ersichtlich. Unabhängig davon eröffnet Euler mit seiner Methode gleichzeitig einen Weg, die inverse Mellin--Transformation aus einer Differentialgleichung zu finden, anstatt sich der Kurvenintegration in der komplexen Ebene zu bedienen.

\paragraph{Das Beispiel der Gammafunktion}
\label{para: Das Beispiel der Gammafunktion}

 Seine Methode wendet Euler in § 13 von \cite{E594} auf die $\Gamma$--Funktion an, dort ist sie sein drittes Beispiel. Folgende Gleichung ist also vorgelegt:

\begin{equation}
    \label{eq: Equation Gamma}
    f(x+1)= xf(x).
\end{equation}
Gemäß der Euler'schen Methode ist diese Hilfsgleichung zu betrachten:

\begin{equation*}
    \int p(t)t^xdt = x\int p(t)t^{x-1}dt+t^{x}q(t).
\end{equation*}
Differenzieren von selbiger gibt

\begin{equation*}
    p(t)t^x= xp(t)t^{x-1}+xt^{x-1}q(t)+t^xq'(t)
\end{equation*}
und daher nach Teilen durch $t^{x}\neq 0$

\begin{equation*}
    p(t)t=xp(t)+xq(t)+tq'(t).
\end{equation*}
Die Koeffizienten der jeweiligen Potenzen von $x$ vergleichend gelangt man zu diesem System von Gleichungen:

\begin{equation}
\label{eq: p, q for Gamma}
    p(t)=q'(t) \quad \text{und} \quad q(t)=-p(t),
\end{equation}
welches von

\begin{equation*}
    p(t)= -q(t)= C\cdot e^{-t} \quad \text{mit} \quad C \neq 0.
\end{equation*}
gelöst wird; somit ist zur Erkenntnis der Integrationsgrenzen die Gleichung

\begin{equation}
    \label{eq: Limiteq Gamma}
    C \cdot t^x e^{-t}=0
\end{equation}
 zu lösen. Diese lässt für $\operatorname{Re}(x)>0$ als Lösungen nur $t_1=0$ und $t=\infty$ zu, sodass  (\ref{eq: Equation Gamma}) gelöst wird von

\begin{equation}
    \label{eq: Solution Gamma}
    f(x) = C \cdot \int\limits_{0}^{\infty} t^{x-1}e^{-t}dt.
\end{equation}
Fordert man nun noch $f(1)=1$ ein, ergibt sich die klassische Darstellung der $\Gamma$-Funktion\footnote{Im Allgemeinen reicht es bekanntermaßen zur \textit{eindeutigen} Definition der $\Gamma$-Funktion nicht aus, lediglich $\Gamma(x+1)=x\Gamma(x)$ und $\Gamma(1)=1$ einzufordern. Denn jede Funktion $f(x)=\Gamma(x)h(x)$ mit $h(x)=h(x+1)$ und $h(1)=f(1)$ genügt der Gleichung ebenfalls. Es muss noch mindestens eine weitere Bedingung hinzutreten, was Euler aber nicht angemerkt hat -- weder in dieser Arbeit \cite{E594} noch in einer anderen. Dennoch gelangt man hier mit der Euler'schen Methode zum bekannten Ergebnis, weil der Ansatz über das Integral das Hinzutreten der periodischen Funktion verhindert.}.

\begin{figure}
    \centering
    \includegraphics[scale=1.3]{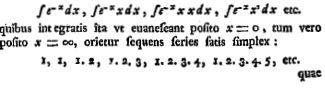}
    \caption{Euler leitet in seiner Arbeit \cite{E594} die bekannte Integraldarstellung der $\Gamma$-Funktion aus seiner Theorie der Differenzengleichungen her.}
    \label{fig:E594Gamma}
\end{figure}

\subsubsection{Den Kontext betreffend: Die Legendre--Polynome}
\label{subsubsec: Den Kontext betreffend: Die Legendre Polynome}

\epigraph{All the truths of mathematics are linked to each other, and all means of discovering them are equally admissible.}{Andrien--Marie Legendre}

Euler hat, wie beispielsweise in der Arbeit \textit{``Euler and the Legendre Polynomials"} (\cite{Ay23}, 2023) dargelegt wird, bereits einige Beschaffenheiten der Legendre--Polynome zutage gefördert. Von Interesse ist, dass er ausgehend von verschiedenen Fragestellungen zu denselben Gleichungen gelangt, dies aber nicht erkannt zu haben scheint.\\

Von Legendre  wurden die nach ihm benannten Polynome in seiner Arbeit  \textit{``Recherches sur l'attraction des sphéroïdes homogènes"} (\cite{Le85}, 1785) (``Untersuchungen zur Attraktion von homogen Sphäroiden") in seinen Studien zum Gravitationspotential

\begin{equation*}
    \dfrac{1}{|\bf{x}-\bf{x}^{\prime}|}= \dfrac{1}{\sqrt{r^2 -2rr^{\prime} \cos (\beta)+(r^{\prime})^2}}= \sum_{n=0}^{\infty} \dfrac{{r^{\prime}}^n}{r^{n+1}}P_n(\cos (\beta)),
\end{equation*}
mit den Längen $r$ und $r^{\prime}$ der Vektoren $\bf{x}$ und $\bf{x}^{\prime}$ und dem Winkel $\beta$ zwischen ihnen als die Polynome $P_n(t)$ eingeführt, womit  das $n$--te Legendre-Polynom als Koeffizient einer Reihenentwicklung entspringt.\\

Vom mathematischen Standpunkt betrachtet sind sie als vollständiges Orthogonalsystem von Polynomen von Interesse und werden heute alternativ über den Ausdruck

\begin{equation}
    \label{eq: Definition Legendre polynomials}
    \dfrac{1}{\sqrt{1-2xt+x^2}}:= \sum_{n=0}^{\infty}P_n(t)x^n
\end{equation}
als Koeffizienten einer Potenzreihenentwicklung definiert. Daneben existieren zahlreiche weitere Möglichkeiten, diese Polynome vorzustellen. In hiesiger Abhandlung soll ihre Definition über die nachstehende Differenzengleichung den Anfang bilden:

\begin{equation}
\label{eq: Difference Equation for Legendre}
(n+1)P_{n+1}(t)= (2n+1)tP_n(t)-nP_{n-1}(t),\quad \text{für} \quad n \in \mathbb{N}_0.
\end{equation}
Die Zusatzforderungen $P_0(t)=1$ sowie $P_1(t)=t$ bestimmen sie eindeutig, woraus sich insgesamt schnell die führenden Legendre--Polynome wie in der unten stehenden Tabelle errechnen. \\

\begin{table}[h]
    \centering
    \begin{tabular}{c}
    \renewcommand{\arraystretch}{2,0}
\setlength{\arraycolsep}{1.5mm}
$\begin{array}{|l|l|} \hline
\quad  n \quad & \quad P_n(t) \\ \hline
\quad 0       & \quad 1 \\
\quad 1       & \quad t \\
\quad 2       & \quad \frac{1}{2}(3t^2-1) \\
\quad 3       & \quad \frac{1}{2}(5t^3-3t) \\
\quad 4       & \quad \frac{1}{8}(35t^4-30t^2+3) \\
\quad 5       & \quad \frac{1}{8}(63t^5-70t^3+15t) \\
\quad 6       & \quad \frac{1}{16}(231t^6-315t^4+105t^4-5) \\
\quad 7       & \quad \frac{1}{16}(429t^7-693t^5+315t^3-35t) \\ \hline
\end{array}$
    \end{tabular}
    \caption{Die ersten Legendre-Polynome für kleine Werte von $n$.}
    \label{Table 1}
\end{table}
Die Gleichung (\ref{eq: Difference Equation for Legendre}) wird sich als Spezialfall allgemeinerer Formeln in einigen Abhandlungen von Euler herauskristallisieren; die Formel  tritt unter anderem in den Euler'schen Publikationen  \textit{``Speculationes super formula integrali $\int \frac{x^ndx}{\sqrt{aa-2bx+cxx}}$, ubi simul egregiae observationes circa fractiones continuas occurrunt"} (\cite{E606}, 1786, ges. 1775) (E606: ``Betrachtungen über die Integralformel $\int \frac{x^ndx}{\sqrt{aa-2bx+cxx}}$, wo zugleich außerordentliche Beobachtungen über Kettenbrüche auftreten") \cite{E672},  \cite{E673}, \cite{E674} sowie \cite{E710} auf. Wohingegen die Darstellung in der ersten Arbeit nicht allzu bekannt zu sein scheint\footnote{Sie werden etwa nicht in der umfassenden Formeltabelle \textit{``Table of Integrals, Series, and Products"} (\cite{Zw14}, 2014) genannt.}, können die in den verbleibenden Werken von Euler mitgeteilten Ausdrücke auf diese bekannte Darstellung zurückgeführt werden:

\begin{equation}
    \label{eq: Laplace Legendre}
    P_n(t) = \dfrac{1}{\pi}\int\limits_{0}^{\pi}\left(t+\sqrt{t^2-1}\cos \varphi\right)^n d\varphi,
\end{equation}
welche von Laplace (1749--1827) in seinem Opus \textit{``Traité de mécanique céleste"} (\cite{La25}, 1825) (``Traktat über die Himmelsmechanik") bewiesen wurde\footnote{In Band 16,2 der Serie 1 der \textit{Opera Omnia} (\cite{OO162}, 1935) wird ebenfalls Laplace als der Entdecker genannt.}. Dieser wird sich nun vor der Darstellung (\ref{eq: Definition Legendre polynomials}) angenommen.

\paragraph{Euler und die Laplace'sche Darstellung}
\label{para : Euler und die Laplace'sche Darstellung}

In den drei unmittelbar aufeinander folgenden Arbeiten \cite{E672},\cite{E673}, \cite{E674} betrachtet Euler die Integrale

\begin{equation}
    \label{eq: Euler Family}
    A_n(a,i):= \int\limits_{0}^{\pi} \dfrac{d \varphi\cos (i \varphi)}{(1+a^2-2a \cos (\varphi))^{n}},
\end{equation}
und zeigt drei Eigenschaften von ihnen\footnote{Hier stellt $i$ eine ganze Zahl und nicht die komplexe Einheit $\sqrt{-1}$ dar.}. Zunächst gibt er eine explizite Formel für die obigen Integrale an. Weiterhin entdeckt er eine Funktionalgleichung, welcher sie Genüge leisten. Zuletzt teilt er eine Differenzengleichung in $n$ für sie mit, welche sich für den Fall $i=0$ auf (\ref{eq: Difference Equation for Legendre}) reduziert.\\

Letztere findet sich in § 74 von \cite{E673} in der Form:
\begin{equation*}
     n(n-1)(1-a^2)^2 \int \dfrac{ d \varphi \cos (i \varphi)}{\Delta^{n+1}}
\end{equation*}
\begin{equation}
\label{eq: Euler's Difference Equation}
    =(n-1)(2n-1)(1+a^2)\int \dfrac{d \varphi \cos(i \varphi)}{\Delta^{n}}+(i^2-(n-1)^2)\int \dfrac{d \varphi \cos(i \varphi)}{\Delta^{n-1}},
\end{equation}
wo Euler die Integrale von $0$ bis hin zu $\pi$ erstreckt versteht. Das Zeichen $\Delta$ entspricht $1-2a\cos (\varphi)+a^2$. Mit der Notation aus (\ref{eq: Euler Family}) drückt sich diese Gleichung folgendermaßen aus:

\begin{equation*}
    n(1-a^2)^2\cdot A_{n+1}(a)=(2n-1)(1+a^2)\cdot A_{n}(a)-(n-1)\cdot A_{n-1}(a).
\end{equation*}
Eine Verschiebung des Index $n \mapsto n+1$ gibt

\begin{equation*}
    (n+1)(1-a^2)^2\cdot A_{n+2}(a)=(2n+1)(1+a^2)\cdot A_{n+1}(a)-n\cdot A_{n}(a).
\end{equation*}
Zur weiteren Vereinfachung mache man die Einführung des Buchstaben $a$ aus \cite{E673} rückgängig. Dort geht Euler nämlich von folgendem Integral aus:
\begin{equation*}
    \int \dfrac{d \varphi \cos (i \varphi)}{(\alpha + \beta \cos (\varphi))^n}.
\end{equation*}
Die Integrale $A_n(a,i)$ aus (\ref{eq: Euler Family}) entspringen aus letzterem Integral durch die Festlegungen  $\alpha = 1+a^2$ und $\beta =-2a$, welche  Euler in § 42 von \cite{E673} ebenfalls macht. Demnach kann die vorausgehende Differenzengleichung  in dieser Manier vorgestellt werden: 

\begin{equation*}
    (n+1)(\alpha^2-\beta^2)\cdot A_{n+2}(a) =(2n+1)\alpha \cdot A_{n+1}(a)-n \cdot A_{n}(a).
\end{equation*}
Setzt man  $\alpha =x$ sowie $\beta= \sqrt{x^2-1}$ und schreibt der Übersicht wegen $A_n(x)$ anstatt $A_n(a)$:

\begin{equation*}
      (n+1)\cdot A_{n+2}(x) =(2n+1)x \cdot A_{n+1}(x)-n \cdot A_{n}(x),
\end{equation*}
was gerade (\ref{eq: Difference Equation for Legendre}) ist. Überdies gilt $A_{n+1}(x)= \pi P_n(x)$, zumal durch direkte Rechnung

\begin{equation*}
    A_1(x) = \int\limits_{0}^{\pi} \dfrac{d \varphi}{x+\sqrt{x^2-1}\cos \varphi}= \pi = \pi \cdot P_0(x),
\end{equation*}
 $x \in [-1,1]$ vorausgesetzt. Mit denselben Restriktionen an $x$, gilt ebenfalls

\begin{equation*}
    A_2(x) = \int\limits_{0}^{\pi} \dfrac{d \varphi}{(x+\sqrt{x^2-1}\cos \varphi)^2}= \pi \cdot x = \pi \cdot P_1(x); 
\end{equation*}
dies impliziert
    
    \begin{equation*}
        P_n(x) = \dfrac{1}{\pi}\int\limits_{0}^{\pi} \dfrac{d \varphi}{(x+\sqrt{x^2-1}\cos(\varphi))^n},
    \end{equation*}
    welche (\ref{eq: Laplace Legendre}) schon sehr nahe kommt. Um gänzlich zu Laplace'schen Darstellung zu gelangen, bedarf es der Funktionalgleichung, welche Euler für die Integrale $A_n(a,i)$ ausfindig gemacht hat; den Beweis reicht er in § 21 von \cite{E710} nach\footnote{Euler teilt dieselbe Formel auch schon in seinen Arbeiten \cite{E672}, \cite{E673}, \cite{E674}  mit. Allerdings ohne Beweis in den beiden erstgenannten, wo er sie mit dem Prädikat Vermutung vorstellt.  Der erste Beweis findet sich in \cite{E674}.}. Besagte Formel lautet:

\begin{equation*}
     \binom{n+i}{i} (1-aa)^{-n} \int \Delta^n d \varphi \cos (i\varphi)
\end{equation*}
\begin{equation}
\label{eq: Euler Functional Equation}
    = \binom{-n-1+i}{i}(1-aa)^{n+1} \int \Delta^{-n-1} d \varphi \cos (i \varphi),
\end{equation}
wo wieder $\Delta = 1-2a \cos (\varphi)+a^2$ ist und die Integrale von $0$ bis hin zu $\pi$ zu nehmen sind.\\ 

 Der Fall $i=0$ ist der wesentliche, in welchem die Binomialkoeffizienten in (\ref{eq: Euler Functional Equation}) alle  $=1$ werden. Nach derselben Substitution wie zuvor, $a$ für $x$ hinauswerfend, liest sich die Euler'sche Gleichung in der eingeführten Notation:

\begin{equation*}
    A_{-n}(a)= A_{n+1}(a),
\end{equation*}
wo $1-a^2=\alpha^2-\beta^2=x^2-(x^2-1)=1$ benutze wurde.  Verbinden der letzten Gleichung mit  $A_{n+1}(a)=\pi \cdot P_n(x)$ offenbart $A_{-n}(x)=\pi \cdot P_n(x)$, was nach der Definition von $A_{n}(a)$ zu dieser Darstellung führt:

\begin{equation*}
    P_n(x) =\dfrac{1}{\pi}\int\limits_{0}^{\pi} (x+\sqrt{x^2-1}\cos (\varphi))^n d \varphi,
\end{equation*}
was gerade (\ref{eq: Laplace Legendre}) ist.\\

Unter  ausschließlicher Applikation von Eulers Formeln ist demnach folgende Identität bewiesen worden:
\begin{equation}
\label{eq: Integral Relation}
   \int\limits_{0}^{\pi} (x+\sqrt{x^2-1}\cos (\varphi))^n d \varphi = \int\limits_{0}^{\pi} (x+\sqrt{x^2-1}\cos (\varphi))^{-n-1} d \varphi, 
\end{equation}
eine Formel, welche im Vorwort zu Band 16,2 der ersten Serie der \textit{Opera Omnia} (\cite{OO162}, 1935) Jacobi (1804 -- 1851)  zugeschrieben wird, der ebenfalls (\ref{eq: Euler Functional Equation}) in seiner Arbeit \textit{ ``Über die Entwickelung des Ausdrucks $(aa-2aa'[\cos \omega \cos \varphi+\sin \omega \sin \varphi \cos(\vartheta-\vartheta')]+a'a')^{-\frac{1}{2}}$"} (\cite{Ja43}, 1843)  bewiesen hat. Es lässt sich demnach behaupten, Euler hätte die Relation ebenfalls zu beweisen vermocht, hätte er sich der Aufgabe angenommen. Seine Funktionalgleichung (\ref{eq: Euler Functional Equation}) ist jedenfalls  allgemeiner als (\ref{eq: Integral Relation}). Überdies ist (\ref{eq: Euler Functional Equation}) ein Spezialfall der Euler'schen Transformation für die hypergeometrische Reihe, sprich der Formel 
\begin{equation*}
    _2F_1(a,b,c;x)= (1-x)^{c-a-b}{}_2F_1(c-a, c-b,c; x),
\end{equation*}
welche Euler in § 9 von \cite{E710} nachweist. Hier steht $_2F_1$ für die Gauß'sche hypergeometrische Funktion, die im nächsten Abschnitt noch eingehender behandelt werden wird (siehe Abschnitt \ref{subsubsec: Die Darstellung betreffend -- Die hypergeometrische Reihe}). Demnach darf (\ref{eq: Integral Relation}) als ein Spezialfall einer gewissen Eigenschaft der hypergeometrischen Funktion ausgelegt werden.\\

Die Vollständigkeit gebietet die Angabe der Gleichung, welche Euler für die Integrale $A_n(a,i)$ deriviert. Diese lautet (siehe, e.g., § 43 in \cite{E673}):

\begin{equation*}
    \int\limits_{0}^{\pi} \dfrac{d \varphi \cos (i \varphi)}{(1+a^2-2a \cos (\varphi))^{n+1}}= \dfrac{\pi a^i}{(1-a^2)^{2n+1}} \cdot V,
\end{equation*}
mit

\begin{equation*}
    V= \binom{n+i}{i}+\binom{n+i}{i+1}\binom{n-i}{1}a^2+\binom{n+i}{i+2}\binom{n-i}{2}a^4+\cdots
\end{equation*}
Euler merkt aber in \cite{E672} wie auch \cite{E673} an, dass diese Formel für $V$ auf einer Vermutung basiert; den Beweis hat er in \cite{E674} nachgeliefert.

\paragraph{Legendres Untersuchung dieser Integrale}
\label{para: Legendres Untersuchung dieser Integrale}

Es ist dem umfassenderen Verständnis zuträglich, an dieser Stelle Legendres Betrachtungen der Integrale (\ref{eq: Euler Family}) einzuschalten, bevor die weiteren Beiträge Eulers zu den Legendre'schen Polynomen ihre Erwähnung finden. \\

Legendre hat eben diese   in seinem Buch \textit{``Exercices du calcul intégral -- Tome premier"} (\cite{Le11}, 1811) (``Übungen im Integralkalkül -- erstes Buch") (p. 372), aber auch im späteren Werk \textit{``Traité des fonctions elliptiques et intégrales Eulériennes -- Tome second"} (\cite{Le26}, 1826) (``Traktat über die elliptischen Funktionen und die Euler'schen Integrale -- zweites Buch"), genauer im ersten Anhang dieses Buches, einer Untersuchung unterworfen.  Zum Zeitpunkt der Verfassung seines Buches \cite{Le26} ist Legendre der Eigenschaften der Integrale bereits gewahr und hat Kenntnis des Euler'schen Beitrags genommen, welchen er an entsprechender Stelle  auch zitiert, wohingegen Euler die ganzen Eigenschaften erst erarbeiteten musste. Die zutage geförderten Relationen unterscheiden sich bei den Autoren nicht grundlegend: So findet Euler, wie eben genauer ausgeführt, die Legendre'sche Gleichung (61), sogar in allgemeinerer Form; zur der in (64) von Legendre mitgeteilten Differentialgleichung gelangt Euler freilich nicht, jedoch gibt er in \cite{E710} (§ 9) gar die Differentialgleichung  für die hypergeometrische Reihe an, welche die Legendre'sche als Spezialfall beinhaltet. Die Differenzengleichung, die sich bei Legendre als (62) findet, nimmt Euler ebenfalls vorweg.  Legendre zueigen ist indes die Gleichung (33), welche sich als die bekannte Partitätsrelation für die Legendre'schen Polynome auffassen lässt. Naturgemäß lässt sie sich leicht aus den Euler'schen Formeln ableiten\footnote{Einen Beweis dieser Behauptung findet man in Abschnitt (\ref{subsubsec: Weitere Untersuchungen zu den Legendre-Polynomen}).}. Legendre übertrifft Euler darin, dass er auch die Integrale 

\begin{equation*}
    \int\limits_{0}^{\pi} \dfrac{ \cos(\lambda\varphi)d\varphi}{(1+2\cos(\varphi)a+a^2)^{\frac{1}{2}}}
\end{equation*}
explizit untersucht. Er deutet sie als elliptische Integrale in einer Weise, welche Euler zeitlebens versperrt geblieben ist. Die Gründe hierfür werden in Abschnitt (\ref{subsubsec: Die Normalform von elliptischen Integralen}) dargelegt werden.

\paragraph{Eulers andere Integraldarstellung}
\label{para: Eulers andere Integraldarstellung}

Wie oben angeschnitten, findet sich in der Arbeit \cite{E606}, welche sich eigentlich Kettenbrüchen widmet, eine weitere Integraldarstellung für die Legendre--Polynome. Nachdem Euler in den Abhandlungen  \cite{E123}  und  \cite{E594}  gezeigt hatte, wie homogene Differenzengleichungen mit linearen Koeffizienten mithilfe von Integralen aufgelöst werden können (siehe auch Abschnitt \ref{subsubsec: Die Mellin--Transformierte bei Euler}), wendet er die dort gewonnenen Einsichten in \cite{E606} nun auf die Familie von Integralen 

\begin{equation*}
    G(n):=\int\limits_{\frac{b-\sqrt{b^2-a^2c}}{c}}^{\frac{b+\sqrt{b^2-a^2c}}{c}} \dfrac{x^ndx}{\sqrt{a^2-2bx+cx^2}}
\end{equation*}
an und weist nach, dass sie dieser Gleichung Genüge leisten:

\begin{equation*}
    na^2\cdot G(n-1)= (2n+1)b\cdot G(n)-(n+1)c\cdot G(n+1).
\end{equation*}

\begin{figure}
    \centering
     \includegraphics[scale=0.8]{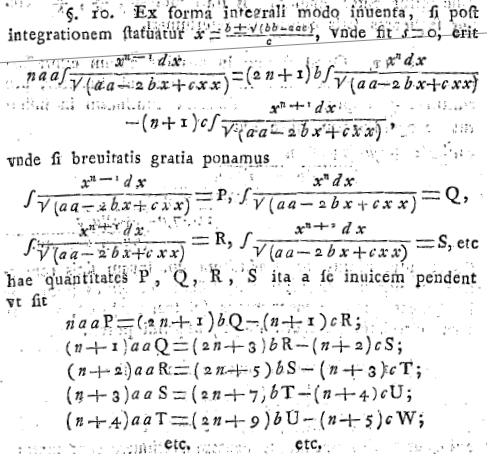}
    \caption{Euler gelangt in seiner Arbeit \cite{E606} zur Darstellung der Legendre--Polynome mithilfe von Integralen. In der ersten Zeile ist die Differenzengleichung für die Legendre--Polynome zu sehen.}
    \label{fig:E606Legendre}
\end{figure}

Beschränkt man sich auf den Fall $a=c=1$  und schreibt $t$ statt $b$, reduziert sich diese Gleichung auf die Differenzengleichung für die Legendre--Polynome (\ref{eq: Difference Equation for Legendre}), was folgenden Schluss gestattet: 

\begin{equation*}
C(t)\cdot    P_n(t) =   \int\limits_{t-\sqrt{t^2-1}}^{t+\sqrt{t^2-1}}\dfrac{x^n dx}{\sqrt{1-2xt+x^2}}.
\end{equation*}
$C(t)$ ist dabei eine von $n$ unabhängige im Weiteren zu eruierende Funktion. Zu diesem Zweck werte man die Integrale für $n=0$ und $n=1$ explizit aus. Eine formale Rechnung liefert:

\begin{equation*}
    \int\limits_{t-\sqrt{t^2-1}}^{t+\sqrt{t^2-1}}\dfrac{x^0 dx}{\sqrt{1-2xt+x^2}} = \log(-1)
\end{equation*}
sowie

\begin{equation*}
    \int\limits_{t-\sqrt{t^2-1}}^{t+\sqrt{t^2-1}}\dfrac{x^1 dx}{\sqrt{1-2xt+x^2}}= \log(-1)\cdot t,
\end{equation*}
was  die Funktion $C(t)$  als Konstante enthüllt. Die Mehrdeutigkeit des Ausdrucks $\log (-1)$ erzwingt die Wahl eines Wertes, jedoch bietet nach entsprechender Festlegung

\begin{equation}
   \label{eq: New Legendre Integral}
    P_n(t)= \dfrac{1}{\log(-1)} \cdot \int\limits_{t-\sqrt{t^2-1}}^{t+\sqrt{t^2-1}}\dfrac{x^n dx}{\sqrt{1-2xt+x^2}},
\end{equation}
diese Gleichung die richtige Formel dar. Für die Wahl des Hauptzweiges des Logarithmus, sodass $\log(-1)=i \pi$, entspringt

\begin{equation}
   \label{eq: New Legendre Speziell}
    P_n(t)= \dfrac{1}{i\pi} \cdot \int\limits_{t-\sqrt{t^2-1}}^{t+\sqrt{t^2-1}}\dfrac{x^n dx}{\sqrt{1-2xt+x^2}},
\end{equation}
in welcher Form Stieltjes (1854--1894) in seiner Arbeit \textit{``Sur les polynômes de Legendre"} (\cite{St90}, 1890) (``Über die Legendre--Polynome") die Legendre--Polynome mit Mitteln der komplexen Analysis ausgedrückt hat. Der Umstand der Mehrdeutigkeit ließe sich überdies durch nachstehendes unbestimmtes Integral in Gänze umgehen:
\begin{equation*}
    \int \dfrac{x^ndx}{\sqrt{1-2xt+x^2}}
\end{equation*}
\begin{equation*}
    = P_n(t)\operatorname{artanh}\left(\dfrac{x-t}{\sqrt{1-2xt+x^2}}\right)+ Q_n(x,t) \sqrt{1-2xt+x^2}+C
\end{equation*}
wobei $Q_n$ ein Polynom in $x$ und $t$ ist und $C$ eine Integrationskonstante ist. Diese Formel kann beispielsweise mit Induktion und (\ref{eq: Difference Equation for Legendre}) nachgewiesen werden, was an dieser Stelle ausgespart wird.

\paragraph{Ein weiterer Ansatz von Euler mit ähnlichen Ergebnissen}
\label{para: Ein weiterer Ansatz mit ähnlichen Ergebnissen}

Hier rückt nun die Definition über die erzeugende Funktion (siehe Gleichung (\ref{eq: Definition Legendre polynomials})) in den Vordergrund, welche in unerwartetem Zusammenhang in Eulers Arbeiten \textit{``Observationes analyticae"} (\cite{E326}, 1767, ges. 1763) (E326: ``Analytische Beobachtungen") und \cite{E551} auftritt. In diesen beiden Werken betrachtet er die allgemeinen Trinome  $(a+bx+cx^2)^n$ für $n \in \mathbb{N}$ und $a,b,c \in \mathbb{C}$ und ist dabei an der Folge der Koeffizienten der Potenz $x^n$ in der Entwicklung dieser Trinome interessiert. Aus der Entwicklung

\begin{equation*}
    (a+bx+cx^2)^n = a^n +\cdots + B_n\cdot x^n+ \cdots +c^nx^{2n},
\end{equation*}
wo die Koeffizienten $B_n$ von den Zahlen $a$, $b$ und $c$ abhängen, definiert Euler die Potenzreihe\footnote{Der Fall $a=b=c=1$ gibt die Koeffizienten der Folge A002426 aus Online Encyclopedia of Integer Sequences (\url{https://oeis.org/}).}:

\begin{equation}
\label{eq: Taylor Series}
    F(x):= \sum_{n=0}^{\infty} B_n x^n.
\end{equation}
Das Problem, welches Euler in \cite{E326} und \cite{E551} eröffnet und löst, besteht im Finden eines geschlossenen Ausdrucks für $F(x)$. Seine Lösung zu diesem Problem (siehe, z.B., § 21 von \cite{E551}) lautet:
\begin{equation*}
    F(x):= \dfrac{1}{\sqrt{1-2bx+(b^2-4ac)x^2}}.
\end{equation*}
Dieser Ausdruck ähnelt der erzeugenden Funktion für die Legendre--Polynome in (\ref{eq: Definition Legendre polynomials}), welche sich  durch Schreiben von $t$ für $b$ und Setzen von $b^2-4ac=1$ auch ergibt. \\

\begin{figure}
    \centering
      \includegraphics[scale=1.1]{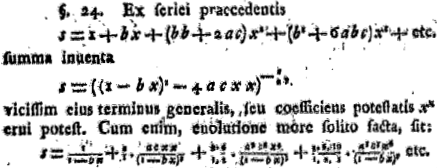}
    \caption{Euler gelangt in seiner Arbeit \cite{E551} zur erzeugenden Funktion der Legendre--Polynome.}
    \label{fig:E551LegendreErzeugende}
\end{figure}

Ein weiterer wichtiger Beitrag aus \cite{E551} besteht in der Entdeckung der Differenzengleichung zwischen drei unmittelbar aufeinander folgenden Legendre--Polynomen, also (\ref{eq: Difference Equation for Legendre}). Euler gibt die folgende Relation zwischen den aufeinander folgenden Koeffizienten der Taylor Reihe für $F$ in (\ref{eq: Taylor Series}) in § 22 von  \cite{E551} an: \\


\begin{equation*}
    r= bq+\dfrac{n-1}{n}(bq-p).
\end{equation*}
In seiner Notation unterdrückt Euler die Abhängigkeit der Buchstaben $r,p,q$ von $n$ sowie dem Buchstaben $b$, welches dasselbe $b$ ist wie im expliziten Ausdruck für die Funktion $F$. Aber Euler lehrt uns, dass $r$ der $n$-te Term in der Folge der Koeffizienten ist, woraus man $q$ aus $r$ durch Schreiben von $n-1$ anstelle von $n$ in $r$ erhält. Gleichermaßen findet man $p$ aus $r$, indem man $n-2$ anstatt von $n$ in $r$ schreibt. Ersetzt man $b$ mit $t$ und stellt alles in moderner Manier mit Indizes dar, setzt demnach $r=B_{n}(t)$, $q=B_{n-1}(t)$ und $p=B_{n-2}(t)$, nimmt Eulers Formel diese Gestalt an:

\begin{equation*}
    B_n(t)= t \cdot B_{n-1}(t)+\dfrac{n-1}{n}(t \cdot B_{n-1}(t)-B_{n-2}(t)).
\end{equation*}
Eine einfache Rechnung zeigt die Gleichwertigkeit jener Gleichung zu dieser:

\begin{equation*}
    n B_n(t)=(2n-1)t\cdot B_{n-1}(t)-(n-1)B_{n-2}(t),
\end{equation*}
welche sich nach der Ersetzung $n \mapsto n+1$ als Differenzengleichung für die Legendre-Polynome (\ref{eq: Difference Equation for Legendre}) hervortut. \\ 

Aus diesen Überlegungen ließe sich aus den Euler'schen Überlegungen heraus folgende Definition der Legendre--Polynome formulieren

\begin{Def}[Legendre--Polynome]
Das $n$--te Legendre--Polynom $P_n(t)$ ist gegeben als Koeffizient der Potenz $x^n$ in der Entwicklung des folgenden Trinoms:

\begin{equation*}
   T(t,x)= \left(\dfrac{t-1}{2}+xt+\dfrac{t+1}{2}x^2\right)^n.
\end{equation*}
\end{Def}
Trotz der mannigfachen Berührungspunkte mit den Legendre'schen Polynomen in diversen Kontexten, scheint Euler von alledem keine Notiz genommen zu haben, was folgende Äußerung in § 7 von \cite{E551} unterstreicht:\\

\textit{``Auch wenn das Ausfindigmachen der Zahl $N$ auf eine Differenzen--Differentialgleichung zurückgeführt werden kann, ist selbige dennoch so beschaffen, dass sie in keiner Weise eine Lösung zuzulassen scheint."}\\

Was Euler hier $N$ nennt, ist gerade das $n$--te Legendre--Polynom $P_n(t)$. Dass die Differenzen--Differentialgleichung\footnote{Euler benutzt diesen Begriff, um Differentialgleichungen zweiter Ordnung zu beschreiben.}, welcher die Legendre'schen Polynome Genüge leisten, keine triviale Lösung aufweist, ist heute ebenfalls einsichtig. Sie werden mancherorts ja gerade umgekehrt als Lösung der Legendre'schen Differentialgleichung \textit{definiert}, was Euler fern lag, wie oben (Abschnitt \ref{subsubsec: Die Rolle der Definition in der Mathematik}) diskutiert worden ist. \\

Unabhängig davon hat Euler mit (\ref{eq: Laplace Legendre}) die ersehnte Lösung in \cite{E606} mit Integralen nachgereicht, was ihm allerdings ebenfalls entgangen zu sein scheint; er macht jedenfalls keine Erwähnung davon, was gleichermaßen das Auftreten der Differenzengleichung der Legendre--Polynome in den drei Arbeiten \cite{E551}, \cite{E606} und \cite{E673} betrifft. Allerdings ist der Zusammenhang zwischen den Koeffizienten einer erzeugenden Funktion (\cite{E551}), Kettenbrüchen (\cite{E606}) und Untersuchungen zu speziellen Integralen (\cite{E673}) überaus schwierig zu erkennen und wurde auch nicht in den Vorworten zu Eulers \textit{Opera Omnia} von den jeweiligen Editoren angemerkt.

\subsubsection{Die Darstellung betreffend -- Die hypergeometrische Reihe}
\label{subsubsec: Die Darstellung betreffend -- Die hypergeometrische Reihe}

\epigraph{The real voyage of discovery consists not in seeking new landscapes, but in having new eyes.}{Marcel Proust}

Den Abschluss der Gegenstände, welche bei Euler zwar auftreten, jedoch nach seinen Nachfolgern benannt werden, soll in dieser Ausarbeitung die hypergeometrische Reihe bilden:

\begin{equation}
\label{eq: Hypergeometric Series}
    _2F_1(\alpha, \beta, \gamma;z) = 1+\frac{\alpha \beta}{1\cdot \gamma} \dfrac{z}{1!}+\frac{\alpha (\alpha +1)\beta (\beta+1)}{\gamma (\gamma +1)}\dfrac{z^2}{2!}+\text{etc.},
\end{equation}
 also die Funktion, welche heute mit der linken Seite bezeichnet und seit seiner wegweisenden Arbeit \cite{Ga13} auch \textit{Gauß}'sche hypergeometrische Funktion genannt wird. Im Vorwort zu Band 19 der ersten Serie der \textit{Opera Omnia} (\cite{OO19}, 1936) fasst Faber (1877--1966) die Situation  wie folgt zusammen: \\

\begin{center}
    \textit{``\textsc{Euler} scheint die Größe des Schatzes, den er aufgeschürft hatte, nicht voll ermessen zu haben; er hat manches davon, was er selbst hätte heben können, seinen Nachfahren, allen voran \textsc{Gauß}, überlassen."}
\end{center}
Es soll der Nachweis erbracht werden, dass  dieses Zitat mit gewissen Einschränkungen zu versehen ist, um noch der Wahrheit zu entsprechen. Schon Klein schreibt auf Seite 269 seines Buchs \cite{Kl56}:\\

\textit{``Die Gaußsche Reihe -- der Name ist wiederum höchst unhistorisch, denn schon Euler hat sie gekannt und die merkwürdigsten Eigenschaften gefunden, auch wenn er Gauß die Konvergenzfragen erledigt hat -- [...]"}\\

Es wird nämlich nachfolgend argumentiert werden, dass sich die Euler'schen von den Gauß'schen Entdeckungen bezüglich der hypergeometrischen mehr in ihrer Darstellung als in ihrem Umfang unterscheiden. 

\begin{figure}
    \centering
  \includegraphics[scale=0.9]{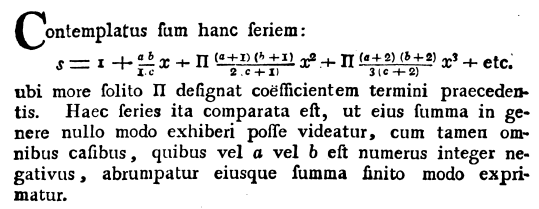}
    \caption{ Euler führt die hypergeometrische Reihe in seiner Arbeit \cite{E710} ein. Euler verwendet dabei das Symbol $\Pi$, um den gesamten vorausgehenden Zahlenkoeffizienten der vorausgehenden Potenz anzuzeigen.}
    \label{fig:E710Hyp}
\end{figure}

\paragraph{Themenübersicht}
\label{para: Themenübersicht}

Die Gauß'sche Arbeit \cite{Ga13}  gibt gleichsam die zu besprechenden Themen vor\footnote{Gauß hat selbst zu Lebzeiten nur die Paragraphen 1 bis 37 veröffentlicht, die restlichen §§ 38 -- 55 sind aus seinem Nachlass entnommen. Da besagter Nachlass jedoch inhaltlich nahtlos an den ersten Teil anknüpft, wird im Folgenden mit \cite{Ga13} auch der Nachlass gekennzeichnet.}. Es wird förderlich sein, eine entsprechende Übersicht mit Angabe, entsprechende Ergebnis bei Euler findet,  der Untersuchung vorausgeschickt zu haben. 

\begin{itemize}
    \item[] §3: Konvergenzbetrachtungen der hypergeometrischen Reihe in der Variable $z$ -- bei Euler in \cite{E43}.
    \item[] §6: Die Entwicklung der Reihe
    \begin{equation*}
    (aa+bb-2ab\cos \varphi)^{-n}= A +2A'\cos \varphi +2A'' \cos 2\varphi +2A'''\cos 3\varphi +\text{etc.}.
\end{equation*} 
 Bei Euler an verschiedenen Stellen, etwa \cite{E464}.
\item[]  § 7: Die Nachbarschaftsrelationen-- bei Euler  in \cite{E123}. 
\item[]  § 12: Die Kettenbruchentwicklung des Quotienten zweier hypergeometrischer Funktionen -- bei Euler in äquivalenter Form in \cite{E123}. 
\item[]  § 15: Die Gauß`sche Summationsformel für die hypergeometrische Funktion, Ermittelung des Wertes: $_2F_1{}(a,b,c;1)$ -- bei Euler in äquivalenter Form in \cite{E575}, \cite{E663}.
\item[]  § 18: Die Produktformel (\ref{eq: Euler Produkt Gamma}) für die Gamma-Funktion $\Gamma$-Funktion -- bei Euler in \cite{E19},\cite{E613}, \cite{E652}. 
\item[]  § 25: Die Reflexionsformel für die $\Gamma$-Funktion -- bei Euler in \cite{E352}, \cite{E421}.
\item[]  § 26: Die Multiplikationsformel für die $\Gamma$-Funktion, also $\Gamma(nx)$ geschrieben als Produkt von $\Gamma$-Funktionen kleineren Arguments -- bei Euler in äquivalenter Form in \cite{E421}. 
\item[]  § 27: \textit{Elastica}--Integrale, bei Euler in \cite{E122}, \cite{E321}.
\item[]  § 28: Die typische Integraldarstellung für die $\Gamma$--Funktion mit Herleitung -- bei  Euler in \cite{E421}, \cite{E19}, \cite{E594}. 
\item[] § 31: Die  $\psi$-Funktion -- bei  Euler in \cite{E212}, \cite{E613} 
\item[] § 35: Eine Integraldarstellung für die $\psi$--Funktion -- nicht von Euler angegeben, aber leicht herzuleiten.
\item[] § 36: Spezielle Frullian'sche Integrale -- bei Euler in \cite{E464}.
\item[] § 38: Die Differentialgleichung für die hypergeometrische Funktion $_2F_{1}(a,b,c;z)$ -- bei Euler in \cite{E710}. 
\item[] §§ 49--55: Quadratische Transformationen für das letzte Argument -- nicht von Euler betrachtet und auch nicht aus seinen Überlegungen abzuleiten. 
\item[] § 55: Ermittelung des Wertes $_2F_1{}(a,b,c; \frac{1}{2})$ -- nicht von Euler ermittelt und aus seinen Formeln heraus nicht möglich.
\end{itemize}

Bis auf die beiden letztgenannten Themen hat Euler demnach in verschiedenen Arbeiten  die Gauß'schen Resultate ebenfalls zutage gefördert  oder ist diesen aus auszuführenden Gründen dermaßen nahe gekommen, dass sie mit wenigen Strichen aus seinen Mitteilungen abgeleitet werden können. Die einzelnen Punkte werden nachstehend gemäß des Umfanges der beizufügenden Erklärungen geordnet abgehandelt.

\paragraph{Konvergenzbetrachtungen}
\label{para: Konvergenzbetrachtungen}

Gauß beginnt seine Untersuchung mit der Frage der Konvergenz der hypergeometrischen Reihe  (\ref{eq: Hypergeometric Series}) im Parameter $z$. Das Ergebnis ist, wie heute wohlbekannt ist, dass $|z|<1$ zu setzen ist, sofern man die anderen Parameter $\alpha$, $\beta$, $\gamma$ keinen weiteren nicht notwendigen Einschränkungen unterwirft. Das Gauß'sche Argument stützt sich dabei auf das Konzept, was später als Quotientenkritierum\footnote{Für absolute Konvergenz der Reihe $S:=\sum_{n=0}^{\infty}a_n$ mit $a_n \in \mathbb{C}$ ist die Bedingung $\lim_{n\rightarrow \infty}\frac{|a_n|}{|a_{n+1}|}\leq x<1$ hinreichend, Gauß wendet dieses Theorem auf den allgemeinen Term der hypergeometrischen Reihe an, und zeigt so, dass $|z|<1$ zur Gewährleistung der Konvergenz sein muss.} von Reihen in die Mathematik eingegangen ist. Gauß schreibt in § 2 seiner Arbeit, wo er $x$ statt $z$ als vierte Variable in (\ref{eq: Hypergeometric Series}) setzt:\\

\textit{``Die Koeffizienten der Potenzen $x^m$, $x^{m+1}$ in unserer Reihe verhalten sich wie}

\begin{equation*}
    1+\dfrac{\gamma+1}{m}+\dfrac{\gamma}{mm}: 1+\dfrac{1+\beta}{m}+\dfrac{\alpha \beta}{mm}
\end{equation*}
\textit{und nähern sich daher dem Verhältnis der Gleichheit umso mehr, umso größer $m$ angenommen wird. Wenn also auch dem vierten Element $x$ ein bestimmter Wert zugeteilt wird, wird von dessen Gestalt die Konvergenz oder Divergenz abhängen. Sooft freilich $x$ ein reeller Wert kleiner als die Einheit zugeteilt wird, ob positiv oder negativ, wird die Reihe, wenn nicht schon gleich zu Beginn, dennoch ab einem gewissen Punkt konvergent sein, und zu einer daraus gänzlich bestimmten endlichen Summe konvergieren. Dasselbe wird durch einen imaginären Wert von $x$ der Form $a+b\sqrt{-1}$ geschehen, sofern $aa+bb<1$ ist. Andererseits wird für einen reellen Wert von $x$ größer als die Einheit oder einen imaginären der Form $a+b\sqrt{-1}$, mit $aa+bb>1$, die Reihe, wenn nicht schon zu Beginn, dennoch nach einem gewissen Intervall divergent sein, sodass von ihrer Summe keine Rede sein kann. Schließlich wird für den Wert $x=1$ (oder allgemeiner einen der Form $a+b\sqrt{-1}$ mit $aa+bb=1$) die Konvergenz von der Beschaffenheit von $\alpha$, $\beta$, $\gamma$ abhängen, worüber wir, und im Speziellen von der Summe der Reihe für $x=1$, im dritten Abschnitt sprechen werden."} \\

Trotz seines Konzeptes von divergenten Reihen\footnote{Die Euler'sche Behandlung divergenter Reihen  wird in dieser Arbeit noch in Abschnitt (\ref{subsubsec: Der Begriff der Summe einer Reihe}) besprochen werden.}, welches Euler etwa in seiner Arbeit \textit{``De seriebus divergentibus"} (\cite{E247}, 1760, ges. 1746) (E247: ``Über divergente Reihen") und seinen \textit{Calculi Differentialis} \cite{E212} vorstellt,  besitzt er durchaus eine Auffassung der Konvergenz einer Reihe; diese bringt er in § 2 seiner Arbeit \cite{E43} sprachlich zum Ausdruck:\\

\textit{``Aber eine Reihe, welche ins Unendliche fortgesetzt eine endliche Summe hat, wird, auch wenn sie doppelt so weit fortgesetzt wird, keine Vermehrung erfahren, sondern das, was nach dem Unendlichen gedanklich hinzugefügt wird, wird in Wahrheit unendlich klein sein. Denn wenn es sich nicht so verhielte, wäre die Summe der unendlichen Reihe, auch wenn sie bis ins Unendliche fortgeführt wird, nicht bestimmt und deswegen auch nicht endlich. Daraus folgt, dass, wenn das, was aus der Fortsetzung über den infinitesimalen Term hinaus entspringt, von endlicher Größe ist, die Summe der Reihe notwendig unendlich sein muss."}\\

Möchte man die Euler'schen Erläuterungen in Formeln überführen, könnte man schreiben:

\begin{equation*}
    \sum_{n=1}^{\infty} a_n =s \quad \Rightarrow \quad \sum_{n=N}^{\infty}a_n =\tilde{s},
\end{equation*}
wobei $s$,$\tilde{s}$ und $N$ endliche Zahlen sind, sodass seine Definition  modern anmutet. Man kann in seiner Formulierung die notwendige Bedingung für die Konvergenz  oder äquivalent eine hinreichende Bedingung für Divergenz einer Reihe erkennen, allerdings nur für \textit{reelle} Reihenglieder. Explizit verwendet scheint er diese Definition indes selten zu haben, sofern man von \cite{E43} absieht. Jedoch ist Euler von der Idee eines Konvergenz\textit{radius} weit entfernt, welches Konzept unter anderem zeigt, dass die betrachtete Funktion als Funktion über den komplexen Zahlen zu verstehen ist\footnote{Euler beschränkt dahingegen, aus Gründen, die  unten in Abschnitt (\ref{subsubsec: Komplexe Analysis}) noch erläutert werden, seine Betrachtungen meist auf die reelle Achse.}, einen Punkt den Gauß vehement unterstrichen hat. Somit erfährt das Klein'sche Resümee,  Gauß habe die Konvergenzfragen der hypergeometrischen Reihe erledigt, seine Bestätigung.

\paragraph{Entwicklung einer Funktion in eine Reihe}
\label{para: Entwicklung einer Funktion in eine Reihe}

In § 6 von \cite{Ga13}  betrachtet  Gauß die Entwicklung der Funktion
\begin{equation*}
    \dfrac{1}{(a^2-2ab \cos \varphi +b^2)^n}
\end{equation*}
für eine beliebige natürliche Zahl $n$ in eine Reihe von Kosinus: 

\begin{equation}
\label{eq: Expansion}
    \dfrac{1}{(a^2-2ab \cos \varphi +b^2)^n} = A_0 +2A_1\cos(\varphi)+2A_2 \cos(2\varphi)+2A_3 \cos(3\varphi)+\text{etc.}, 
\end{equation}
Sein Ergebnis für die Koeffizienten lautet:

\begin{equation*}
    A_k= \dfrac{n(n+1)(n+2)\cdots (n+k-1)}{1\cdot 2 \cdot 3 \cdots k}a^{-2n-k}b^k{}_2F_1{}\left(n,n+k,k+1,\dfrac{b^2}{a^2}\right).
\end{equation*}
Faktorisiert man $a^2-2ab\cos \varphi +b^2$ in lineare Faktoren, entwickelt beide Faktoren in eine Reihe und multipliziert sie, ergibt sich aus einem Koeffizientenvergleich mit (\ref{eq: Expansion}), gelangt man wie Gauß zur alternativen Darstellung:

\begin{equation*}
    A_k=\dfrac{n(n+1)(n+2)\cdots (n+k-1)}{1\cdot 2 \cdot 3 \cdots k} (a^2+b^2)^{n-k}a^kb^k \times
\end{equation*}
\begin{equation*}
     {}_2F_1{}\left(\dfrac{1}{2}n+\dfrac{3}{2},\dfrac{1}{2}n+2, k+1, \dfrac{4a^2b^2}{(a^2+b^2)^2}\right).
\end{equation*}
Demnach implizieren die beiden Ausdrücke die folgende Transformationsformel für die hypergeometrische Reihe
\begin{equation}
    \label{eq: Transformation Gauss quadratic}
    _2F_1{} (a,b,1-a-b;z)= (1+z)^{-a}{}_2F_1{} \left(\dfrac{a}{2},\dfrac{a+1}{2},1-a-b;\dfrac{4z}{(1+z)^2}\right),
\end{equation}
welche üblicherweise E. Kummer (1810--1893) zugeschrieben wird, welcher sie in der Arbeit \textit{``Über die hypergeometrische Reihe $1+\frac{\alpha \cdot \beta}{1\cdot \gamma}x+\frac{\alpha (\alpha+1)  \beta (\beta +1)}{1\cdot2 \gamma(\gamma +1)}x^2+\frac{\alpha (\alpha+1)(\alpha+2)  \beta (\beta +1)(\beta +2)}{1\cdot 2 \cdot 3 \gamma(\gamma +1)(\gamma +2)}x^3+\text{etc.}$"} (\cite{Ku36}, 1836) auf andere Weise aus einer allgemeineren Betrachtung heraus bewiesen hat.\\

Euler scheint solche Transformationen nie betrachtet zu haben, in welchen das letzte Argument einer quadratische Transformation untergeht, worauf gegen Ende dieses Abschnittes noch detaillierter eingegangen werden wird. Allerdings stößt Euler bei mehreren Begebenheiten auf die Reihe (\ref{eq: Expansion});  exemplarisch sei  seine Arbeit \textit{``Nova methodus quantitates integrales determinandi"} (\cite{E464}, 1775, ges. 1774) (E464: ``Eine neue Methode Integrale zu bestimmen"), insbesondere  §§ 32 -- 39 erwähnt. Hier beginnt Euler von der Reihe

\begin{equation*}
    Q(z)= z \sin (u) +z^2 \sin (2u)+z^3 \sin (3u)+z^4\sin(4u)+\text{etc.}
\end{equation*}
und zeigt, dass sie diesem Ausdruck gleichwertig ist

\begin{equation*}
    Q(z)= \dfrac{z \sin (u)}{1-2z \cos(u)+z^2}.
\end{equation*}
Über iterierte Integration und Differentiation lassen sich mehr Reihen aus der geschlossenen Formel für $Q(z)$ ableiten, wie es Euler in besagter Arbeit auch tut; (\ref{eq: Expansion}) ist ebenfalls darunter.\\

Ein anderer Weg, weder von Gauß noch von Euler bestritten, beginnt mit der geometrischen Reihe 

\begin{equation*}
    \dfrac{1}{1-z}=1+z+z^2+z^3+z^4+\text{etc.}.
\end{equation*}
Setzt man nun $z=re^{i\varphi}$, sodass 

\begin{equation*}
    \dfrac{1}{1-re^{i\varphi}}= 1+re^{i \varphi}+r^2e^{2i \varphi}+r^3e^{3i \varphi}+r^4e^{4i \varphi}+\text{etc.}
\end{equation*}
\begin{equation*}
    =1+r\cos(\varphi)+r^2 \cos(2\varphi)+ \cdots +i(r\sin(\varphi)+r^2\sin(2\varphi)+\cdots)
\end{equation*}
Die linke Seite kann beschrieben werden als

\begin{equation*}
    \dfrac{1}{1-re^{i\varphi}}=\dfrac{1}{1-re^{i\varphi}}\cdot \dfrac{1-re^{-i\varphi}}{1-re^{-i\varphi}}= \dfrac{1-re^{-i \varphi}}{1-2re^{-i\varphi}-2re^{+i\varphi}+r^2}
\end{equation*}
\begin{equation*}
    = \dfrac{1-re^{-i\varphi}}{1-2\cos(\varphi)r+r^2}= \dfrac{1-r\cos(\varphi)}{1-2\cos(\varphi)r+r^2}+ \dfrac{ir\sin(\varphi)}{1-2\cos(\varphi)r+r^2}.
\end{equation*}
Also gibt ein Vergleich der Imaginärteile der beiden Seiten direkt die Euler'sche Formel. 
Entsprechend gelangt man durch  $(n-1)$-fache Differentiation zu (\ref{eq: Expansion}). Die Reihe taucht demnach unter anderem im Kontext anderer Integrale in Eulers Opus auf. Sie ließe sich ebenfalls in den schon diskutierten \cite{E672}, \cite{E673}, \cite{E674} entdecken (Abschnitt \ref{subsubsec: Den Kontext betreffend: Die Legendre Polynome}), zumal sie dort den Integranden des von Euler studierten Integrals bildet.

\paragraph{Die Differentialgleichung für die hypergeometrische Funktion}
\label{para: Die Differentialgleichung für die hypergeometrische Funktion}

Wohingegen Gauß seine Investigationen mit der Differentialgleichung für die hypergeometrische Funktion beginnt\footnote{Dies bezieht sich auf den zweiten Teil von \cite{Ga13} aus dem Nachlass, ab § 38. Gauß erkennt hier insbesondere, dass es zweckmäßig ist, das Studium der hypergeometrischen Funktion von ihrer Differentialgleichung aus zu beginnen. Bei Euler bleibt sie ein Nebenprodukt.}, haben Eulers Studien -- diametral zu Gauß --  in \cite{E710} ihren Ursprung in der Reihe (\ref{eq: Hypergeometric Series}), aus welcher er dann die Differentialgleichung deduziert. Seine Rechnungen sind sehr gradlinig und führen ihn zum bekannten  Ergebnis, dass, wenn die hypergeometrische Reihe mit $F$ abbreviiert wird, gilt:

\begin{equation}
\label{eq: DGL Hyp}
    0=x(1-x)\dfrac{d^2F}{dx^2}+[c-(a+b+1)x]\dfrac{dF}{dx}-abF.
\end{equation}
Die genauen Schritte können an dieser Stelle ausgespart werden.

\begin{figure}
    \centering
     \includegraphics[scale=0.9]{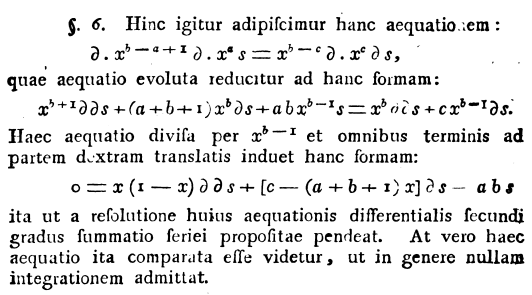}
    \caption{ Euler gelangt in seiner Arbeit \cite{E710} zur Differentialgleichung für die hypergeometrische Funktion.}
    \label{fig:E710DGLHypergeometrischeReihe}
\end{figure}

\paragraph{Integrale zur Curva elastica}

Die Elastica-Integrale, welche so genannt werden, weil sie bei Betrachtung der sogenannten \textit{Curva elastica} auftreten, haben Euler an verschiedenen Stellen seiner Laufbahn eingenommen. Er behandelt sie unter anderem im zweiten Anhang der \textit{Methodus} \cite{E65}; die Arbeit \textit{``De miris proprietatibus curvae elasticae sub aequatione $y=\int \frac{xx}{\sqrt{1-x^4}}$ contentae"} (\cite{E605}, 1786, ges. 1775) (E605: ``Über die wunderbaren Eigenschaften der Curva Elastica, welche mit der Gleichung $y=\int \frac{xx}{\sqrt{1-x^4}}$ ausgedrückt wird") ist in Gänze dieser Kurve gewidmet. Insbesondere die spezielle Formel

\begin{equation}
\label{eq: Elastica}
    \int\limits_{0}^1 \dfrac{x^2dx}{\sqrt{1-x^4}} \cdot \int\limits_{0}^1 \dfrac{dx}{\sqrt{1-x^4}}= \dfrac{\pi}{4}
\end{equation}
hebt Euler wegen ihrer Symmetrie besonders hervor. Gauß tut es ihm in dieser Hinsicht gleich und bemerkt bei dieser Gelegenheit außerdem, dass Euler diese Formel und ähnliche andere ``mit großer Arbeit"{} erlangt hat; ``aus seinen Formeln folge sie hingegen leicht"{}. Es ist nicht ganz ersichtlich, auf welche Arbeit von Euler sich Gauß hier bezieht. Er könnte zum einen  die Arbeit \textit{``De productis ex infinitis factoribus ortis"} (\cite{E122}, 1750, ges, 1739) (E122: ``Über Produkte, die aus unendlich vielen Faktoren entspringen") meinen, wo Gleichung (\ref{eq: Elastica}) sich in § 17 neben vielen anderen ähnlicher Gestalt findet. Jedoch ist besagte Abhandlung vornehmlich Produkten gewidmet, welche  aus einem Quotienten bestehen, dessen Zähler und Nenner jeweils  unendliche Produkte sind, deren Terme in einer arithmetischen Progression fortschreiten. Euler bemerkt, dass sich selbige als Quotient zweier Beta--Integrale darstellen lassen. Diese Integrale  untersucht Euler dann gesondert in der Arbeit \textit{``Observationes circa integralia formularum $\int x^{p-1}dx (1-x^n)^{\frac{q}{n}-1}$ posito post integrationem $x = 1$"} (\cite{E321}, 1766, ges. 1765) (E321: ``Beobachtungen zu den Integralen der Form $\int x^{p-1}dx (1-x^n)^{\frac{q}{n}-1}$, wenn nach der Integration $x=1$ gesetzt wird"), wo sich auch alle von Gauß in \cite{Ga13} zu diesem Thema behandelten Integrale finden. Es ist  wahrscheinlich, dass Gauß die letzte genannte Euler'sche Arbeit nicht bekannt war, weil die erwähnte Eigenschaft der Elastica ein leichtes Korollar der dortigen Untersuchungen ist. So wird Gauß sich aller Wahrscheinlichkeit nach auf die Arbeit \cite{E605} beziehen, wo Euler in seinen Studien zur Curva elastica erst in § 21 zur oben erwähnten Eigenschaft gelangt. Da Euler in dieser Abhandlungen alles aus den Eigenschaften der Kurve selbst heraus ableitet, ist der Beweis  naturgemäß umfassender als der  Gauß'sche Zugang (und identische Euler'sche in \cite{E321}) über das Beta--Integral. \\

Seit Legendre  wird (\ref{eq: Elastica}) sie indes als Spezialfall der nach ihm benannten Relation

\begin{equation*}
    K(\varepsilon)E(\sqrt{1-\varepsilon^2})+E(\varepsilon)K(\sqrt{1-\varepsilon^2})+K(\varepsilon)K(\sqrt{1-\varepsilon^2})=\dfrac{\pi}{2}
\end{equation*}
mit

\begin{equation*}
    K(\varepsilon)= \int\limits_{0}^1 \dfrac{dx}{\sqrt{(1-x^2)(1-\varepsilon^2 x^2)}} \quad \text{und} \quad E(\varepsilon)=\int\limits_{0}^1 \dfrac{\sqrt{1-\varepsilon^2x^2}}{\sqrt{1-x^2}}dx
\end{equation*}
gesehen. Diese Formel findet sich bei Legendre in seinem Buch \cite{Le11} (S. 61).

\begin{figure}
    \centering
   \includegraphics[scale=1.4]{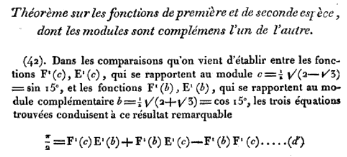}
    \caption{Legendre gibt in seinem Buch \cite{Le11} die nach ihm benannte Relation für die elliptischen Integrale für die speziellen Werte $b=\frac{1}{2}\sqrt{2+\sqrt{3}}$ und $c=\frac{1}{2}\sqrt{2-\sqrt{3}}$  an.}
    \label{fig:Legendrerelation}
\end{figure}

\paragraph{Die Integraldarstellung der Gamma-Funktion mitsamt Herleitung}

In § 28 gibt Gauß eine Herleitung der Integraldarstellung der Fakultät an, genauer:

\begin{equation*}
    n!= \int\limits_{0}^{\infty} y^{n-1}e^{-y}dy.
\end{equation*}
Gauß' Idee ist dabei, das Integral

\begin{equation*}
    \int\limits_{0}^m y^{n-1}\left(1-\dfrac{y}{m}\right)^mdy
\end{equation*}
zunächst für endliches $m$ zu  berechnen und anschließend $m$ unendlich groß werden zu lassen. Gauß nimmt also an, dass

\begin{equation*}
  \lim_{m \rightarrow \infty} \int\limits_{0}^m y^{n-1}\left(1-\dfrac{y}{m}\right)^mdy 
\end{equation*}
\begin{equation*}
    = \int\limits_{0}^{\lim_{m \rightarrow \infty}m}   y^{n-1} \lim_{m \rightarrow \infty} \left(1-\dfrac{y}{m}\right)^mdy   = \int\limits_{0}^{\infty} y^{n-1}e^{-y}dy.
\end{equation*}
Das letzte Integral ist gerade die $\Gamma$-Funktion. Natürlich ist die Gauß'sche Herleitung aus moderner Sicht wegen der zwei gleichzeitig ins Unendliche wachsenden Ausdrücke nicht vollkommen streng. Ein rigoroser Beweis nach dem Gauß'schen Vorbild findet man beispielsweise in der Monografie \textit{``Die Gammafunktion"} von Nielsen (1865--1931) (\cite{Ni05}, 1905).\\

Gauß leitet  die $\Gamma$-Funktion aus einem Integral über algebraische Funktionen, genauer Beta--Integrale, her, was dem Euler'schen Ansatz aus \cite{E19} ähnelt\footnote{Euler wiederholt die auf diesem Gedankengang basierende Herleitung in der deutlich später verfassten Arbeit \cite{E421}. Eulers Argumentationskette ist an vielen Stellen auseinander gesetzt, etwa im Buch \cite{Du99} und umfassender in \cite{Ay21a}.}.  Euler gelangt mit seinen Überlegungen im Unterschied zu Gauß zum Ausdruck

\begin{equation*}
    n! = \int\limits_{0}^1 \left(\log\left(\dfrac{1}{x}\right)\right)^ndx,
\end{equation*}
welcher obigem vermöge der Substitution $x=e^y$ gleichwertig ist.  Weiterhin illustriert Euler durch die beigefügten verbalen Erläuterungen ein Beispiel für die Methodus inveniendi, welches hier nicht vorenthalten werden soll. In §§ 3--7 des besagten Papiers \cite{E19} schreibt er: \\

\textit{``Ich hatte zuvor geglaubt, dass der allgemeine Term der Reihe $1$, $2$, $6$, $24$ etc., wenn nicht algebraisch, dennoch als Exponentialfunktion gegeben ist. Aber nachdem ich herausgefunden hatte, dass gewisse Terme von der Quadratur des Kreises abhängen}\footnote{Euler spricht hier die Tatsache an, dass modern ausgedrückt, etwa $\Gamma \left(\frac{1}{2}\right)=\sqrt{\pi}$ ist. Zu diesem Ergebnis ist er über das Wallis'sche Produkt für $\pi$ gelangt, also der Formel
\begin{equation*}
    \frac{\pi}{2}= \frac{2 \cdot 2}{1\cdot 3}\cdot  \frac{4 \cdot 4}{3\cdot 5} \cdot  \frac{6 \cdot 6}{5\cdot 7} \cdots
\end{equation*}
Man vergleiche seine Ausführungen in § 2 von \cite{E19}.}\textit{, bin ich zur Erkenntnis gelangt, dass weder algebraische noch exponentielle Größen hinreichen, um sie auszudrücken. [...] Nachdem ich beobachtet hatte, dass unter den Differentialen Formeln solcher Art existieren, welche eine Integration in gewissen Fällen zulassen und dann auch algebraische Größen geben, aber in anderen keine Integration zulassen und dann Quadraturen darbieten, welche von der Quadratur des Kreises abhängen, kam es mir in den Sinn, dass solche Formeln womöglich geeignet sind, den allgemeinen Term der erwähnten und anderen Progressionen auszudrücken. [...] Aber die Differentialform muss eine variable Größe beinhalten. [...] Der Klarheit wegen sage ich, dass $\int pdx$ der allgemeine Term der aufzufindenden Progression ist, welche man wie folgt findet; man lasse $p$ eine Funktion von $x$ und Konstanten bedeuten, unter welchen $n$ enthalten sein muss. Man stelle sich $pdx$ integriert vor und eine solche Konstante addiert, dass für $x=0$ das ganze Integral verschwindet; dann setze man $x$ einer bekannten Größe gleich. Danach, wenn  im gefundenen Integral nur Größen, welche sich  auf die Progression beziehen, aufgefunden werden, wird es den Term ausdrücken, dessen Index $n$ ist. Anders ausgedrückt, wird das in dieser Weise ausgedrückte Integral der allgemeine Term sein. [...] Daher habe ich viele Differentialformeln betrachtet, welche nur eine Integration zulassen, wenn man $n$ als eine ganze positive Zahl nimmt, sodass die Hauptterme algebraisch werden, und habe daraus Progressionen gebildet."}\\

Eulers grundlegende Einsicht besteht im Postulat der Existenz einer Funktion $p(x,n)$ solcher Art, dass

\begin{equation*}
    n ! = \int\limits_{a}^{b}p(x,n)dx \quad \text{für} \quad n \in \mathbb{N}.
    \end{equation*}
In den folgenden Paragraphen probiert Euler systematisch mannigfache Funktionen\footnote{Alle von Euler betrachteten Funktionen lassen sich leicht auf das Beta--Integral reduzieren.} und gelangt schließlich zur Integraldarstellung  der $\Gamma$-Funktion. Euler zeigt nun zuerst

\begin{equation*}
    \dfrac{1 \cdot 2 \cdot 3 \cdots n}{(f+g)(f+2g)\cdots (f+ng)}= \dfrac{f+(n+1)g}{g^{n+1}}\int\limits_{0}^{1} x^{\frac{f}{g}}dx(1-x)^n,
\end{equation*}
was leicht mit Induktion gezeigt wird\footnote{Euler tut dies allerdings in seiner Arbeit \cite{E19} nicht, ist sich aber angesichts seiner vorherigen Ausführung der Gültigkeit dieses Ausdrucks ebenfalls sicher.}. Für $g=0$ reduziert sich der Ausdruck linker Hand nun auf $1\cdot 2 \cdots n =n!$, also

\begin{equation*}
    \dfrac{n!}{f^n} = f\lim_{g \rightarrow 0} \int\limits_{0}^{1}\dfrac{x^{\frac{f}{g}}(1-x)^n}{g^{n+1}}dx.
\end{equation*}
Um sich nun $g$  im Nenner zu entledigen, setzt Euler  $x=y^{\frac{g}{f+g}}$ und gelangt zu

\begin{equation*}
   \dfrac{n!}{f^n} = f\lim_{g \rightarrow 0} \int\limits_{0}^{1}\dfrac{g}{f+g}\dfrac{\left(1-x^{\frac{g}{f+g}}\right)^n}{g^{n+1}}dx.
\end{equation*}
Der Fall $f=1$ ist nun von besonderem Interesse, denn:

\begin{equation*}
    n! =\lim_{g \rightarrow 0} \int\limits_{0}^{1}\dfrac{g}{1+g}\dfrac{\left(1-x^{\frac{g}{1+g}}\right)^n}{g^{n+1}}dx.
\end{equation*}
Dies ist aber dasselbe wie:

\begin{equation*}
    n! = \lim_{g \rightarrow 0} \int\limits_{0}^{1}\dfrac{\left(1-x^{\frac{g}{1+g}}\right)^n}{g^{n}}dx.
\end{equation*}
 Den Grenzwert lässt sich in das Integral hineinziehen
 
\begin{equation*}
    n! =  \int\limits_{0}^{1}\lim_{g \rightarrow 0}\dfrac{\left(1-x^{g}\right)^n}{g^{n}}dx.
\end{equation*}
Der Grenzwert ist in letzterer Formel im Nenner der Potenz $x$ bereits gezogen, und Umschreiben besagter Formel gibt

\begin{equation*}
    n! =  \int\limits_{0}^{1}\left(\lim_{g \rightarrow 0}\dfrac{(1-x^{g})}{g}\right)^ndx,
\end{equation*}
was die Stetigkeit der Funktion $x^n$ zu tun gestattet. \\

\begin{figure}
    \centering
    \includegraphics[scale=0.7]{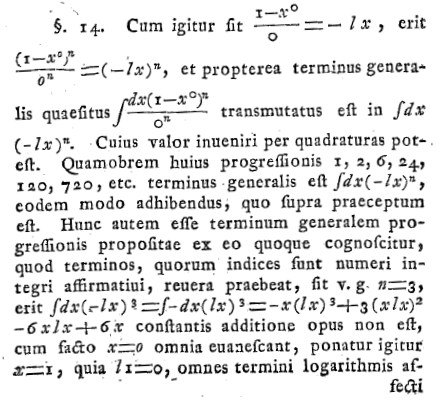}
    \caption{Euler gelangt in seiner Arbeit \cite{E19} durch Grenzwertbildung zur Integraldarstellung der $\Gamma$-Funktion. }
    \label{fig:E19GammaInt}
\end{figure}

Für natürliches, aber auch nicht natürliches, $n$ kann dieser Grenzwert mit der l'Hospital'schen Regel ermittelt werden. Man findet:

\begin{equation*}
    n ! = \int\limits_{0}^{1} \left(\log \dfrac{1}{x}\right)^ndx.
\end{equation*}
Und dies ist natürlich die Integraldarstellung für $\Gamma(n+1)$. Daher erschließt sich die Euler'sche Präferenz dieser Formel gegenüber der heute üblichen

\begin{equation*}
    \Gamma(n+1)= \int\limits_{0}^{\infty} t^n e^{-t}dt,
\end{equation*}
 denn Euler wurde  auf natürlichem Wege zu seiner Version geführt.

\paragraph{Die Frullian'schen Integrale}

Gauß betrachtet im letzten Paragraphen seiner Publikation \cite{Ga13} (§ 37) das Integral

\begin{equation}
\label{eq: Log Frulliani}
    \int\limits_{0}^1 \dfrac{u^{\alpha-1}-u^{\beta -1}}{\log u}du= \log \left(\dfrac{\alpha}{\beta}\right), 
\end{equation}
welches heute zu den Frullian'schen Integralen gerechnet wird.  Dies sind seit ihrer Einführung in der Arbeit  \textit{``Sopra Gli Integrali Definiti"} (\cite{Fr28}, 1828) (``Über gewisse bestimmte Integrale") uneigentliche Integrale der Form

\begin{equation*}
    \int\limits_{0}^{\infty} \dfrac{f(ax)-f(bx)}{x}dx.
\end{equation*}
Setzt man etwa für eine stetige Funktion $f$ zum einen $f(\infty)=B$ und $f(0)=A$, entspringt der allgemeine Ausdruck:

\begin{equation*}
    \int\limits_{0}^{\infty} \dfrac{f(ax)-f(bx)}{x}dx=(B-A)\log \dfrac{a}{b}.
\end{equation*}
Von dieser Formel findet sich bei Euler in § 27 von \cite{E464} (in leicht anderer Form vorgestellt) die Behandlung für die spezielle Wahl $f(x)=e^{-x}$\footnote{Integrale der Form

\begin{equation*}
    \int\limits_{0}^1 \dfrac{R(x)}{\log x}dx
\end{equation*}
mit einer Funktion $R$, sodass $R(0)=R(1)=0$, haben wiederholt Eulers Untersuchungen eingenommen. Stellvertretend sei Arbeit \textit{``De valore formulae integralis $\int \frac{x^{a-1}dx}{\log x}\cdot \frac{(1-x^b)(1-x^c)}{(1-x^n)}$  a termino x = 0 usque ad x = 1 extensae"} (\cite{E500}, 1780, ges. 1776) (E500: ``Über den Wert des Integrals $\int \frac{x^{a-1}dx}{\log x}\cdot \frac{(1-x^b)(1-x^c)}{(1-x^n)}$ erstreckt von der Grenze $x=0$ bis hin zur Grenze $x=1$") erwähnt.}.\\

Euler betrachtet explizit ebenfalls obiges Integral (\ref{eq: Log Frulliani}) um seiner selbst willen. Seine erste systematische Untersuchung findet sich in seiner Arbeit \cite{E464}. Obiges, auch von Gauß behandeltes, Integral findet sich dabei in § 6 dieses Papiers. Euler stellt eine Herleitung über Differenzieren unter dem Integralzeichen nach einem Parameter vor, eine bis dahin nicht geläufige Methode, eine \textit{``Methoda nova"}, wie Euler im Titel schreibt. Sein Argument lautet wie folgt: Bekanntermaßen gilt
\begin{equation*}
    \int\limits_{0}^1 x^{\alpha-1}dx= \dfrac{1}{\alpha} \quad \text{für} \quad \alpha>0.
\end{equation*}
Integriert man diese Gleichheit nun über $\alpha$ von $x=a$ bis $x=b$, erhält man

\begin{equation*}
    \int\limits_{a}^b  \int\limits_{0}^1 x^{\alpha-1}dx d\alpha = \int\limits_{0}^1 \dfrac{x^{b-1}-x^{a-1}}{\log x}dx =\left. \log \alpha \right|_{a}^{b}= \log \left(\dfrac{b}{a}\right),
\end{equation*}
was  gerade (\ref{eq: Log Frulliani}), lediglich mit anderen Buchstaben, ist. \\

\begin{figure}
    \centering
      \includegraphics[scale=1.1]{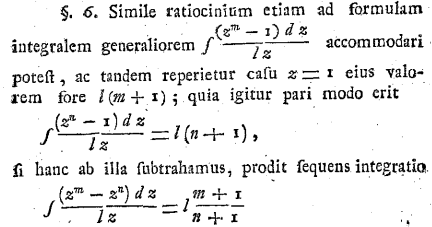}
    \caption{Euler gibt in seiner Arbeit \cite{E464} ein Frulliani'sches Integral an. Die Grenzen der Integration sind $z=0$, was Euler nicht explizit erwähnt, und $z=1$ für die obere Grenze.}
    \label{fig:E464Log(ab)}
\end{figure}

In seiner  Arbeit \cite{Ga13} sagt Gauß aber  explizit, dass Euler ebenfalls -- allerdings mit einer anderen Methode -- zu diesem Integral gelangt ist. Er selbst deduziert (\ref{eq: Log Frulliani}) als speziellen Fall der $\psi$-Funktion, welche eigens noch gleich besprochen wird.\\

Mit den ``anderen"{} Methoden bezieht sich Gauß allerdings wahrscheinlich  auf die just vorgestellte Methode des Differenzierens unter dem Integralzeichen als auf die Euler'sche Abhandlung \textit{``Observationes in aliquot theoremata illustrissimi de la Grange"} (\cite{E587}, 1785, ges. 1775) (E587: ``Bemerkungen zu einigen Themen des höchst illustren Lagrange"). Denn hier, neben der in \cite{E464} präferierten Methode des Differenzierens unter dem Integralzeichen, bemerkt Euler in §§ 15--18, dass

\begin{equation*}
    x^{\alpha}= e^{\alpha\log  x}= 1+ \dfrac{\alpha\log (x)}{1}+\dfrac{\alpha^2 \log^2 (x)}{1\cdot 2}+\dfrac{\alpha^3 \log^3 (x)}{1\cdot 2 \cdot 3}+\text{etc.}.
\end{equation*}
Weiterhin

\begin{equation*}
    x^{\alpha}-x^{\beta}= 1+ \dfrac{(\alpha-\beta)\log (x)}{1}+\dfrac{(\alpha^2-\beta^2) \log^2 (x)}{1\cdot 2}+\dfrac{(\alpha^3-\beta^3) \log^3 (x)}{1\cdot 2 \cdot 3}+\text{etc.}
\end{equation*}
Zieht man nun die folgende Darstellung für die Fakultät heran

\begin{equation*}
  \int \limits_{0}^{1}\log^{n}(x)dx =(-1)^n n!,
\end{equation*}
entspringt durch durch gliedweises Integrieren die Formel

\begin{equation*}
    \int\limits_{0}^1 \dfrac{x^{\alpha}-x^{\beta}}{\log x}dx= \dfrac{\alpha-\beta}{1} -\dfrac{\alpha^2-\beta^2}{2}+\dfrac{\alpha^3-\beta^3}{3}-\dfrac{\alpha^4-\beta^4}{4}+\text{etc.}.
\end{equation*}
Demnach zeigt die bekannte Potenzreihe $\log(1+x)=x-\frac{x^2}{2}+\frac{x^3}{3}-\frac{x^4}{4}+\cdots$, dass gerade gilt:

\begin{equation*}
   \int\limits_{0}^1 \dfrac{x^{\alpha}-x^{\beta}}{\log x}dx = \log (\alpha+1)-\log (\beta +1)=\log \dfrac{\alpha+1}{\beta +1},
\end{equation*}
welches (\ref{eq: Log Frulliani}) mit leicht anderen Buchstaben ist. Inhaltlich geht also Gauß auch hier nicht über Euler hinaus. Insbesondere der letzte Euler'sche Beweis ist gar unmittelbarer als der Gauß'sche, welcher die Formel auch hier  als speziellen Wert einer Funktion auffasst. \\

\paragraph{Die $\psi$-Funktion}
\label{para: Die Digamma-Funktion}

Euler kommt in seinen Arbeiten der logarithmischen Ableitung der $\Gamma$-Funktion, welche oft mit $\psi$ notiert wird, sehr nahe, wovon weiter oben (Abschnitt \ref{subsubsec: Ein anderes Vorhaben: Das Weierstraß-Produkt}) schon Zeugnis abgelegt worden ist. In seiner Arbeit \cite{Ga13} definiert Gauß selbige als 

\begin{equation}
\label{eq: Gauss psi}
    \Psi(z)= \Psi(0) +1+\dfrac{1}{2}+\dfrac{1}{3}+\cdots +\dfrac{1}{z}.
\end{equation}
Demnach ist $\Psi$ eine Interpolation der harmonischen Summe  $H_n:=\sum_{k=1}^n \frac{1}{k}$. Die Gauß'sche Funktion $\Psi(z)$ ist mit der heute gebräuchlicheren $\psi$-Funktion über die simple Gleichung $\Psi(z-1)=\psi(z)$ verknüpft. Euler behandelt Interpolationsprobleme dieser Art an verschiedenen Stellen mit immer wieder anderen Methoden. Die harmonische Summe tritt etwa in \cite{E20} auf, genauer betrachtet Euler dort in § 4 das Integral

\begin{equation}
    \label{eq: Interpolation Harmonic Series}
    h(n)= \int\limits_{0}^1 \dfrac{1-x^n}{1-x}dx. 
\end{equation}
Bemerkt man nun

\begin{equation*}
    1+x+x^2+x^3+\cdots +x^{n-1}= \dfrac{1-x^n}{1-x},
\end{equation*}
gilt

\begin{equation*}
    h(n)= \sum_{k=1}^n \dfrac{1}{k},
\end{equation*}
sodass $h(n)$ ebenfalls die harmonische Reihe bedeutet. In besagter Abhandlung gibt Euler außerdem  die exakten Summationen für die Werte $n=\frac{1}{2}, \frac{3}{2}, \frac{5}{2}$ etc. an, deren erster

\begin{equation*}
    h\left(\dfrac{1}{2}\right)= 2 -2\log 2
\end{equation*}
ist.\\

Nach dem Vorwort zu Band 19 der Serie 1 der \textit{Opera Omnia} (\cite{OO19}, 1936) (S. LXV) findet man die Formel (Gleichung [77] in Gauß' Arbeit \cite{Ga13}) 

\begin{equation}
    \label{eq: Gauss psi int}
    \psi(x)= \int\limits_{0}^1 \left(- \dfrac{1}{\log z}-\dfrac{z^{x-1}}{1-z}\right)dz
\end{equation}
nicht in den Euler'schen Arbeiten. Lediglich der Spezialfall $x=1$ taucht in Eulers Arbeiten auf; etwa  in seinen Arbeiten \textit{``De numero memorabili in summatione progressionis harmonicae naturalis occurrente"} (\cite{E583}, 1785, ges. 1776) (E583: ``Über die bemerkenswerte Zahl, die in der Summation der natürlichen harmonischen Progression auftritt"), besonders in § 18, und \textit{``Evolutio formulae integralis $\int dx\left(\frac{1}{1-x}+\frac{1}{\log x}\right)$ a termino $x = 0$ ad $x = 1$ extensae"} (\cite{E629}, 1789, ges. 1776) (E629: ``Entwicklung der Integralformel $\int dx\left(\frac{1}{1-x}+\frac{1}{\log x}\right)$  erstreckt von der unteren Grenze $x=0$ bis hin zu $x=1$"). Eulers Aufmerksamkeit richtet sich dabei  auf die Euler--Mascheroni--Konstante, heute mit $\gamma$ bezeichnet\footnote{Die Ursache, warum diese Konstante heute mit $\gamma$ und nicht etwa mit $C$, wie Euler sie bei seiner Erstentdeckung in \cite{E43} genannt hat, notiert wird, eruiert Sandifer in seinem Artikel \textit{``Gamma the constant"} (\cite{Sa07oct}, 2007).}, definiert als

\begin{equation*}
    \gamma := \lim_{n \rightarrow \infty} \left(\sum_{k=1}^n \dfrac{1}{k}- \log (n)\right).
\end{equation*}
Im Gegensatz zu Gauß betrachtet er das obige titelgebende  Integral aus \cite{E629} nicht als einen Spezialfall einer allgemeineren Formel. Jedoch, Eulers Formel aus dem zuletzt genannten Werk gebrauchend
\begin{equation*}
    \gamma = \int\limits_{0}^1 \left(\dfrac{dx}{1-x}+\dfrac{1}{\log x}\right)dx
\end{equation*}
und sie von seiner Interpolationsformel für die harmonische Reihe subtrahierend (\ref{eq: Interpolation Harmonic Series}), offenbart unmittelbar die Gauß'sche Formel (\ref{eq: Gauss psi int}). \\

Demnach hätte Euler ohne große Mühe auch diese Identität eruieren können. Zumal für Gauß  (\ref{eq: Gauss psi}) den Ausgangspunkt bildet, weiß er um die Wichtigkeit des Wertes  $\Psi(0)(=\psi(1))=-\gamma$  für seine späteren Untersuchungen, wohingegen Euler aller Wahrscheinlichkeit nach diese Integraldarstellung nicht nieder geschrieben hat, weil andere Fragen sein Interesse auf sich zogen. Diese Unterschiedlichkeit der Ausgangsfrage bzw. der betrachteten Gegenstände kann überdies als ein zentrale Unterschied zwischen der Mathematikphilosophie Gauß' und Eulers angesehen werden.

\paragraph{Die Nachbarschaftsrelationen und Kettenbrüche}
\label{para: Die Nachbarschaftrelationen und Kettenbrüche}

Die sogenannten Nachbarschaftsrelation werden von Gauß in §§ 7--11 von \cite{Ga13} eingeführt und anschließend einem Studium unterworfen. Diese können als homogene Differenzengleichungen der ersten drei Argumente der hypergeometrischen Funktion $_2F_1{}(\alpha,\beta,\gamma;z)$ gedeutet werden. In § 7 listet Gauß  alle $15$ möglichen Gleichungen dieser Art auf. Die Verbindung zum Euler'schen Werk wird erkannt, sofern man die Nachbarschaftsrelationen betrachtet, welche sich nur in einem Argument auf einmal verändern. Als Beispiel sei Gleichung [1] aus Gauß' Arbeit herangezogen:

\begin{equation}
\label{eq: Gauss contiguous}
    0=(\gamma -2\alpha -(\beta -\alpha)z){}_2F_1(\alpha,\beta, \gamma; z)
\end{equation}
\begin{equation*}
    +\alpha (1-x){}_2F_1(\alpha +1,\beta, \gamma;z)-(\gamma -\alpha){}_2F_1(\alpha-1,\beta,\gamma;z).
\end{equation*}
Diese ist also eine Differenzengleichung in $\alpha$, von welchem Parameter die Koeffizienten lineare Funktionen sind, Gleichungen welcher Gattung oben (Abschnitt \ref{subsubsec: Die Mellin--Transformierte bei Euler}) auseinander gesetzt worden sind. Die Verbindung besagter Untersuchungen zur hypergeometrischen Funktion verlangt einige vorauszuschickende Erklärungen: Genauer bedarf es Eulers Arbeit \textit{``Constructio aequationis differentio-differentialis $Ay du^2 + (B+Cu)dudy + (D+Eu+Fuu)ddy = 0$"} (\cite{E274}, 1763, ges. 1760) (E274: ``Konstruktion der Differentialgleichung $Aydu^2+(B+Cu)dudy+(D+Eu+Fuu)ddy=0$"). Diese Abhandlung, wie der Titel bereits andeutet, handelt von der Konstruktion einer Lösung zur allgemeinen Differentialgleichung 

\begin{equation}
\label{eq: DEQ E274}
    A\dfrac{d^2y}{du^2}+(B+Cu)\dfrac{dy}{du}+(D+Eu+Fu^2)y=0.
\end{equation}
Dafür nimmt Euler in § 3 von besagter Abhandlung an, dass eine Lösung in dieser Form gegeben ist

\begin{equation*}
    y=\int Pdx(u+x)^n,
\end{equation*}
wo $P$ eine zu diesem Zeitpunkt noch unbekannte Funktion von $x$ und $n$ ein von diesem $x$ unabhängiger Parameter ist. Die Grenzen des Integrals sind  an die jeweiligen Bedingungen anzupassen. \\

\begin{figure}
    \centering
    \includegraphics[scale=1.0]{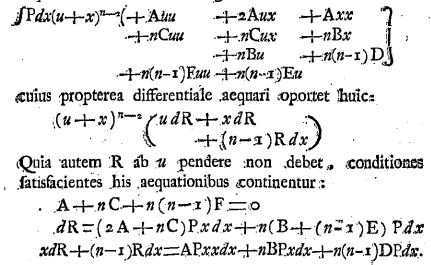}
    \caption{Eulers Erläuterung  der Wahl der Funktion $P$ zur Erfüllung einer vorgelegten Differentialgleichung mit dem Ausdruck $\int Pdx(u+x)^n$ aus seiner Arbeit \cite{E274}.}
    \label{fig:E274DGL}
\end{figure}

Der Vergleich von  (\ref{eq: DEQ E274})  mit der Differentialgleichung der hypergeometrischen Reihe  (\ref{eq: DGL Hyp}) enthüllt letztere als Spezialfall der ersten. Die Arbeit \cite{E273} bricht unvermittelt nach ersten Beispielen ab\footnote{ In der \textit{Opera Omnia} Version (\cite{OO22}, 1936) der Arbeit \cite{E274} erfährt der Leser von A. Speiser in einer Fußnote, dass der Rest der originalen Arbeit im editorialen Prozess verloren gegangen ist.}. So bleibt es im Unklaren, wie nahe  Euler der Integraldarstellung  von (\ref{eq: Hypergeometric Series})  

\begin{equation}
\label{eq: Integral Rep}
  B(\beta, \gamma- \beta)   _2F_1(\alpha, \beta, \gamma; z)= \int\limits_{0}^{1} t^{\beta-1}(1-t)^{\gamma -\beta -1}(1-z t)^{-a}dt
\end{equation}
wirklich kommt. Hier ist $B(x,y)$ wieder das Beta--Integral definiert als

\begin{equation*}
    \label{Def: Beta}
    B(x,y):= \int\limits_{0}^1 t^{x-1}(1-t)^{y-1}dt \quad \text{für} \quad \operatorname{Re}(x)>0,\operatorname{Re}(y)>0.
\end{equation*}

Die Darstellung (\ref{eq: Integral Rep}), die oft als  Integral des Euler'schen Typs bezeichnet wird\footnote{Man konsultiere beispielsweise das Buch \textit{``Theory of Hypergeometric Functions"} (\cite{Ao11}, 2011).}, obwohl Euler selbst diese Formel nie selbst niedergeschrieben zu haben scheint\footnote{Verschiedene Autoren beziehen wie etwa \cite{We24} sich dabei auf Baileys Buch 
 \textit{``Generalized hypergeometric Series"} (\cite{Ba35}, 1935) als das erstzitierende Buch mit der Angabe, wo Euler diese Formel aufgeschrieben haben soll. Aber weder der entsprechende Abschnitt noch die Bibliographie von Baileys Monografie beinhalten eine Referenz zu Eulers Arbeiten.}. Am nächsten kommt der Integraldarstellung noch ein Ausdruck aus dem zweiten Teil seiner \textit{Calculi Integralis} \cite{E366} Hier gibt Euler in § 1033, modern ausgedrückt

 \begin{equation*}
     y(u):= \int\limits_{0}^{c} x^{n-1}(u^2+x^2)^{\mu}(c^2-x^2)^{\nu}dx
 \end{equation*}
 als Lösung der Differentialgleichung

 \begin{equation*}
     u(c^2+u^2)y''-((n+2\mu-1)(c^2+u^2)+2(\mu+\nu)u^2)y'+2\mu(n+2\mu+2\nu)uy=0
 \end{equation*}
 an. Die Wesensnähe dieser Gleichung zur Differentialgleichung der hypergeometrischen Reihe (\ref{eq: DGL Hyp}) und der Integraldarstellung für $y(u)$ zu (\ref{eq: Integral Rep}) ist offenkundig.\\

 \begin{figure}
     \centering
   \includegraphics[scale=1.0]{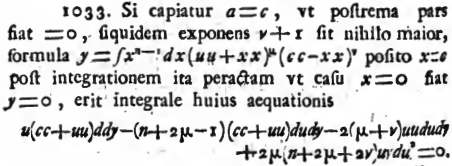}
     \caption{Euler gelangt in den Ausführungen seines Lehrbuchs zur Integralrechnung \cite{E366} zu einer Integraldarstellung einer Funktion, welche derjenigen der hypergeometrischen Reihe (\ref{eq: Integral Rep}) sehr ähnlich ist.}
     \label{fig:E366EulerHypInt - Kopie}
 \end{figure}

Mit diesen Ausführungen kann die Verbindung von Eulers Werk zu (\ref{eq: Gauss contiguous}) geschlagen werden. Genauer sind die Ergebnisse und Ideen der Arbeiten  \cite{E123},  \textit{``De formatione fractionum continuarum"} (\cite{E522}, 1782, ges. 1775) (E522: ``Über die Konstruktion von Kettenbrüchen"),  \cite{E594} vonnöten: Wie Gauß die Kettenbrüche aus den Nachbarschaftsrelationen ableitet, so tut Euler es ihm in den erwähnten Arbeiten gleich, indem er homogene Differenzengleichungen mit linearen Koeffizienten der Form wie folgt umschreibt (siehe auch Abschnitt (\ref{subsubsec: Die Mellin--Transformierte bei Euler}) für den Bezug zur Mellin--Transformierten.):

\begin{equation}
\label{eq: Euler difference}
    (a+\alpha n)f(n)= (b+\beta n)f(n+1)+(c+ \gamma n)f(n+2),
\end{equation}
sodass die Beziehung zum Gauß'schen Ergebnis aufgedeckt wird. Der Hauptunterschied der beiden Autoren besteht darin, dass Euler  primär an Kettenbrüchen interessiert ist und bemerkt, dass er sie aus der Differenzengleichung heraus konstruieren kann, wohingegen Gauß  von Differenzengleichung aus beginnt und die Kettenbrüche als Nebenprodukt erhält. Dieses Diametrale zwischen der Gau\ss {}schen und Euler'schen Herangehensweise  spiegelt zugleich dem Unterschied der Arbeitsweise der beiden Mathematiker insgesamt wider.\\

Bündig soll betrachtet werden, wie Gauß die Kettenbrüche für den Quotienten zweier hypergeometrischer Reihen abgeleitet hat. Neben vielen anderen nennt er den folgenden Kettenbruch (Gleichung [25] in \cite{Ga13}):

\begin{equation}
\label{eq: Gauss Cfrac}
    \dfrac{_2F_1(\alpha,\beta +1,\gamma +1;x)}{_2F_1(\alpha,\beta, \gamma; x)}= \cfrac{1}{1-\dfrac{ax}{1-\cfrac{bx}{1-\cfrac{cx}{1-\cfrac{dx}{1-\text{etc.}}}}}}
\end{equation}
mit

\begin{equation*}
    \renewcommand{\arraystretch}{2,0}
    \setlength{\arraycolsep}{0.0mm}
    \begin{array}{rclrcl}
         a &~=~& \dfrac{\alpha (\gamma- \beta)}{\gamma (\gamma +1)} \qquad \qquad & b &~=~& \dfrac{(\beta+1)(\gamma+1-\alpha)}{(\gamma+1)(\gamma +2)} \\
          c &~=~& \dfrac{(\alpha +1) (\gamma+1- \beta)}{(\gamma +2) (\gamma +3)} \qquad \qquad & d &~=~& \dfrac{(\beta+2)(\gamma+2-\alpha)}{(\gamma+3)(\gamma +4)} \\
           e &~=~& \dfrac{(\alpha +2) (\gamma+2- \beta)}{(\gamma +4) (\gamma +5)} \qquad \qquad & f &~=~& \dfrac{(\beta+3)(\gamma+3-\alpha)}{(\gamma+5)(\gamma +6)} \\
             &  &\text{etc.}  &  & &\text{etc.}
    \end{array}
\end{equation*}
Gauß gibt weiter Kettenbrüche für andere Quotienten (siehe die Gleichungen [26]-[30] in \cite{Ga13}), von welchen  

\begin{equation*}
    \dfrac{_2F_1(\alpha,\beta+1,\gamma;x)}{_2F_1(\alpha,\beta,\gamma;x)}= \cfrac{1}{1-\cfrac{\alpha x}{\gamma}\cdot \dfrac{_2F_1(\alpha,\beta+2,\gamma;x)}{_2F_1(\alpha,\beta,\gamma;x)}}
\end{equation*}
erwähnt sei.\\

 Auch Euler schreibt die Differenzengleichung  (\ref{eq: Euler difference}) um, nämlich wie nachstehend:
\begin{equation*}
    \dfrac{f(n)}{f(n+1)}= \dfrac{b+\beta n}{a+\alpha n}+\dfrac{c+\gamma n}{a+\alpha n} \dfrac{f(n+2)}{f(n+1)}.
\end{equation*}
Auf diese Weise, indem man die letzte Gleichung wiederholt anwendet, gelangt man zu der Gleichung (siehe § 50 in \cite{E123}):

\begin{equation*}
    \dfrac{f(1)}{f(0)}=\cfrac{a}{b+\cfrac{(a+\alpha)c}{b+\beta+\cfrac{(a+2\alpha)(c+\gamma)}{b+2\beta+\cfrac{(a+3\alpha)(c+2\gamma)}{b+3\beta+\cfrac{(a+4\alpha)(c+3\gamma)}{b+4\beta +\text{etc.}}}}}}.
\end{equation*}
Nach vielen Beispielen, setzt Euler  schließlich in § 63 von \cite{E123}

\begin{figure}
    \centering
    \includegraphics[scale=0.9]{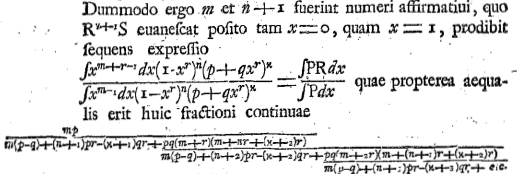}
    \caption{Euler stellt in seiner Arbeit \cite{E123} den Quotienten von zwei hypergeometrischen Funktionen, hier ausgedrückt über Integrale, als Kettenbruch dar.}
    \label{fig:E123KettenbruchHyper}
\end{figure}

\begin{equation}
    P=x^{m-1}(1-x^r)^n(p+qx^r)^{\kappa} \quad \text{und} \quad R=x^r,
\end{equation}
was die folgende Differenzengleichung produziert

\begin{equation*}
    (a+ \nu \alpha)\int\limits_{0}^{1} PR^{\nu}dx=(b+\nu \beta) \int\limits_{0}^1 PR^{\nu +1}dx+(c+ \nu \gamma)\int\limits_{0}^1 P R^{\nu +2}dx,
\end{equation*}
sodass in der Notation von (\ref{eq: Euler difference}) dann $f(n)=\int\limits_{0}^1 PR^{n}dx$ ist. Diese Gleichung führt  zum Kettenbruch für den Quotienten
\begin{equation}
\label{eq: Quotient Euler hyp Integrals}
    \dfrac{\int\limits_{0}^1 PRdx}{\int\limits_{0}^1 Pdx}= \dfrac{\int\limits_{0}^1 x^{m+r-1}dx(1-x^r)^n(p+qx^r)^{\kappa}}{\int\limits_{0}^1 x^{m-1}dx(1-x^r)^n(p+qx^r)^{\kappa}}.
\end{equation}
Bemerkt man weiter, dass

\begin{equation*}
    \int\limits_{0}^1 x^{m+r-1}dx(1-x^r)^n(p+qx^r)^{\kappa} 
\end{equation*}
vermöge (\ref{eq: Integral Rep}) ausgedrückt werden kann, weist dies (\ref{eq: Quotient Euler hyp Integrals}) als Quotienten zweier hypergeometrischer Reihen nach\footnote{Diesen findet man auch  im Buch \textit{``Geschichte der Mathematik 1700–1900"} \cite{Di85} (p. 38), jedoch ohne Bezugnahme auf die hypergeometrische Reihe.}. Demnach hat also auch Euler schon die Kettenbrüche für die benachbarten hypergeometrischen Funktionen zutage gefördert, obwohl er keine Erwähnung davon macht, nicht zuletzt weil er wie erwähnt die Integraldarstellung (\ref{eq: Integral Rep}) nicht explizit niedergeschrieben zu haben scheint.\\

Den Abschluss dieses Themenblockes bildet die Erwähnung, dass der von Gauß ([36] in \cite{Ga13}) gegebene Kettenbruch

\begin{equation*}
    _2F_1 \left(\dfrac{m}{n},1,\gamma,- \gamma nt\right) = 1-mt+m(m+n)tt-m(m+n)(m+2n)t^3+\text{etc.}
\end{equation*}
\begin{equation*}
    =\cfrac{1}{1+\dfrac{mt}{1+\cfrac{nt}{1+\cfrac{(m+n)t}{1+\cfrac{2nt}{1+\cfrac{(m+2nt)t}{1+3nt+\text{etc.}}}}}}}
\end{equation*}
von Euer in der Abhandlung \textit{``De transformatione seriei divergentis $1 - mx + m(m+n)x^2 - m(m+n)(m+2n)x^3 + \text{etc.}$ in fractionem continuam"} (\cite{E616}, 1788, ges. 1776) (E616: ``Über die Transformation der divergenten Reihe $1 - mx + m(m+n)x^2 - m(m+n)(m+2n)x^3 + \text{etc.}$ in einen Kettenbruch") einer Untersuchung unterworfen wird. Auch dieses Exempel zeigt einen zentralen Unterschied zwischen Euler und Gauß: Euler widmet seine Arbeiten sehr oft einer einzigen sehr speziellen Frage, einem einzelnen Objekt oder einer Funktion, Technik oder Methode, Gauß  hingegen stellt erst die allgemeine Theorie vor und geht dann zu Beispielen über, womit die Hankel'sche Analyse des Euler'schen Arbeitshabitus eine weitere Bestätigung erfährt (siehe  den entsprechenden Paragraphen aus Abschnitt \ref{para: Ausführliche Darstellung des Gedankengangs und Präsentation von Rechnungen}). 

\paragraph{Die Gauß'sche Summationsformel für die hypergeometrische Reihe}
\label{para: Die Gauß'sche Summationsformel für die hypergeometrische Reihe}

Hier wird nachgewiesen, dass Euler die sogenannte Gauß'sche Summationsformel für die hypergeometrische Reihe bereits gekannt hat.  Grundlage des schon in der Arbeit \textit{``Euler and the Gaussian Summation Formula for the Hypergeometric Series"} (\cite{Ay24b}, 2024) vorgetragenen Arguments konstituiert die Euler'sche Ausarbeitung \textit{``Plenior expositio serierum illarum memorabilium, quae ex unciis potestatum binomii formantur"} (\cite{E663}, 1794, ges. 1776) (E663: ``Eine umfassendere Erklärung jener merkwürdigen Reihen, welche aus den Binomialkoeffizienten gebildet werden"). Es soll

\begin{equation}
\label{eq: Gauss Summation}
    _2F_1(\alpha, \beta, \gamma; 1) = \dfrac{\Gamma(\gamma -\alpha -\beta)\Gamma(\gamma)}{\Gamma(\gamma -\alpha)\Gamma(\gamma -\beta)} \quad \text{mit} \quad \operatorname{Re}(\gamma)>\operatorname{Re(\alpha+\beta)}
\end{equation}
gezeigt werden. Die Darstellung (\ref{eq: Integral Rep}) darf  angesichts der vorherigen Betrachtungen als Euler ebenfalls bekannt vorausgesetzt werden. Natürlich kann letztgenannte leicht gezeigt werden, indem man $(1-zt)^{-a}$ in eine Potenzreihe entwickelt und termweise integriert, von welchem Grundgedanken die gesamte Euler'sche Abhandlung \cite{E663} durchzogen ebenfalls durchzogen ist. So möchte Euler in § 19  die Reihe

\begin{equation*}
    S= 1+ \mathfrak{A}A+\mathfrak{B}B +\mathfrak{C}C+\mathfrak{D}D+\text{etc.}
\end{equation*}
summieren, wo er die Zahlen  $A$, $B$, $C$, $D$ etc.  und $\mathfrak{A}$, $\mathfrak{B}$, $\mathfrak{C}$, $\mathfrak{D}$ etc. wie folgt definiert: Er setzt

\begin{equation*}
    (1-x^b)^{-\frac{a}{b}}=1+Ax^b+Bx^{2b}+Cx^{3b}+\text{etc.} 
\end{equation*}
und

\begin{equation*}
    (1-x^b)^{-\frac{\alpha}{b}}=1+\mathfrak{A}x^b+ \mathfrak{B}x^{2b}+\mathfrak{C}x^{3b}+\text{etc.}
\end{equation*}
In moderner Notation ausgedrückt, nimmt Euler demnach die Summation der folgenden Reihe in Angriff:

\begin{equation}
\label{eq: S}
   S= 1 + \binom{\frac{a}{b}}{1}\binom{\frac{\alpha}{b}}{1}+\binom{\frac{a}{b}}{2}\binom{\frac{\alpha}{b}}{2}+\binom{\frac{a}{b}}{2}\binom{\frac{\alpha}{b}}{2}+\text{etc.}= \sum_{k=0}^{\infty} \binom{\frac{a}{b}}{k}\binom{\frac{\alpha}{b}}{k}.
\end{equation}
Sein Endresultat, immer noch im selben Paragraphen, ist

\begin{equation}
    \label{eq: Euler Sum}
    S= \dfrac{1}{\Delta} \int\limits_{0}^{1} \dfrac{x^{a-1}dx}{\sqrt[b]{(1-x^b)^{a+\alpha}}} \quad \text{mit} \quad \Delta = \dfrac{\pi}{b \sin \frac{a \pi}{b}}.
\end{equation}
Vergleicht man das Integral in Eulers Formel mit (\ref{eq: Integral Rep}) möchte man auf Grund der drei Parameter annehmen, dass die Euler'sche Formel schon zu  (\ref{eq: Gauss Summation}) äquivalent ist. Allerdings zeigt eine schnelle Rechnung:

\begin{equation*}
    S=  \sum_{k=0}^{\infty} \binom{\frac{a}{b}}{k}\binom{\frac{\alpha}{b}}{k} = {}_2F_1\left(\dfrac{a}{b}, \dfrac{\alpha}{b},1;1\right).
\end{equation*}
Dass (\ref{eq: Euler Sum}) tatsächlich nicht drei \textit{unabhängige} Parameter beinhaltet bzw. nur von den zwei Verhältnissen $a:b$ und $\alpha:b$ abhängig ist,  wird unter der Verwendung der Substitution $y=x^b$ und anschließender Vereinfachung schnell eingesehen. \\

Um die Summationsformel (\ref{eq: Gauss Summation}) aus den Euler'schen Formeln heraus zu zeigen, bedarf es  des Ergebnisses aus § 22 von \cite{E663}. In der Beschreibung des dort vorgestellten Problems möchte er die Summe der folgenden Reihe finden:

\begin{equation*}
    \binom{\frac{a}{b}}{n}\binom{\frac{\alpha}{b}}{0}+ \binom{\frac{a}{b}}{n+1}\binom{\frac{\alpha}{b}}{1}+\binom{\frac{a}{b}}{n+1}\binom{\frac{\alpha}{b}}{2}+\text{etc.}= \sum_{k=0}^{\infty} \binom{\frac{a}{b}}{n+k}\binom{\frac{\alpha}{b}}{k},
\end{equation*}
was gerade die ``verschobene"{} Version der Reihe $S$ in (\ref{eq: S}) ist; Euler notiert diese neue Reihe entsprechend $S^{(n)}$. Am Ende desselben Paragraphen gelangt er zur allgemeinen Formel
\begin{equation}
    \label{eq: Euler general summation formula}
   S^{(n)}= \dfrac{1}{\Delta}\int\limits_{0}^{1} \dfrac{x^{a+nb-1}dx}{\sqrt[b]{(1-x^b)^{\alpha +a}}},   \quad \text{mit} \quad \Delta = \dfrac{\pi}{b\sin \frac{a\pi}{b}}.
\end{equation}
Das just zuvor über  den Buchstaben  $b$ Bemerkte, behält auch hier seine Gültigkeit: Die Allgemeinheit der Formel erfährt keine Einschränkung der Allgemeinheit durch die Wahl  $b=1$. Demnach beinhaltet (\ref{eq: Euler general summation formula}) tatsächlich drei unabhängige Parameter und  (\ref{eq: Gauss Summation}) lässt sich gewiss daraus ableiten, zumal  Eulers Formel  natürlich ein spezieller Fall von (\ref{eq: Integral Rep}) ist -- nämlich der für $z=1$ .\\

\begin{figure}
    \centering
    \includegraphics[scale=0.8]{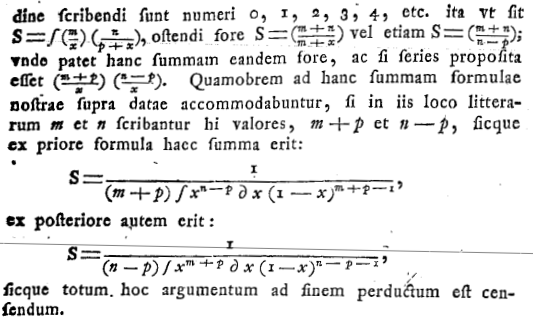}
    \caption{ Eulers Formel, welche der Gauß'schen Summationsformel (\ref{eq: Gauss Summation}) gleichwertig ist, aus seiner Arbeit \cite{E663}. Das Symbol $\left(\frac{a}{b}\right)$ bezeichnet bei Euler hier den Binomalkoeffizienten $\binom{a}{b}$. Weiter benutzt Euler hier das Integralzeichen $\int$ statt $\sum$, um eine Summe anzuzeigen.}
    \label{fig:E663F(a,b,c,1)}
\end{figure}

Um aber die Gauß'sche Summationsformel herzuleiten,  verwende man die von Euler in § 42 von \cite{E663} angegebene Identität:

\begin{equation}
\label{eq: Euler Sum with Binomial}
  \binom{m}{0}\binom{n}{p}+ \binom{m}{1}\binom{n}{p+1} +\text{etc.}=\sum_{k=0}^{\infty} \binom{m}{k}\binom{n}{p+k}  = \binom{m+n}{m+p}.
\end{equation}
Euler\footnote{Die Umbenennung der Buchstaben in dieser letzten Gleichung im Vergleich zu  (\ref{eq: S}) und (\ref{eq: Euler general summation formula}) wird auch von Euler ohne weitere Begründung vorgenommen.} drückt somit eine Summe von Binomialkoeffizienten als einen einzelnen Binomialkoeffizienten aus. Die Buchstaben $n$, $m$ und $p$ sind dabei nicht notwendig ganze Zahlen. Der Kürze wegen sei die linke Seite der letzten Gleichung $S_g$ genannt, was $S_g$ über $_2F_1$ wie folgt ausgedrückt:

\begin{equation}
    \label{eq: Summation Formula Euler modern}
    S_g = \binom{n}{p}{}_2F_1(-m,p-n, p+1;1).
\end{equation}
Um die Brücke zu (\ref{eq: Gauss Summation}) zu schlagen, hat man  $-m=\alpha$, $\beta=p-n$ und $\gamma=p+1$ zu setzen, dass  $m=-\alpha$, $p=\gamma-1$ und $n=\gamma -\beta -1$ ist. Daher, indem man  (\ref{eq: Euler Sum with Binomial}) und (\ref{eq: Summation Formula Euler modern}) kombiniert, entspringt die Gleichung

\begin{equation*}
    \binom{\gamma -\beta -1}{\gamma -1}{}_2F_1(\alpha, \beta, \gamma; 1) = \binom{-\alpha +\gamma-\beta -1}{-\alpha +\gamma-1}.
\end{equation*}
Auflösen nach $_2F_1$ auf und Schreiben der Binomialkoeffizienten mithilfe von $\Gamma$-Funktionen liefert:

\begin{equation*}
    _2F(\alpha, \beta, \gamma;1) =  \binom{-\alpha +\gamma-\beta -1}{-\alpha +\gamma-1}{ \binom{\gamma -\beta -1}{\gamma -1}}^{-1}
\end{equation*}
\begin{equation*}
    = \dfrac{\Gamma(\gamma-\beta -\alpha)}{\Gamma(\gamma-\alpha)\Gamma(1-\beta)} \cdot \dfrac{\Gamma(\gamma)\Gamma(1-\beta)}{\Gamma(\gamma-\beta)}=  \dfrac{\Gamma(\gamma -\alpha -\beta)\Gamma(\gamma)}{\Gamma(\gamma -\alpha)\Gamma(\gamma -\beta)}.
\end{equation*}
Das ist gerade die Gauß'sche Summationsformel (\ref{eq: Gauss Summation}), womit der Nachweis erbracht ist, dass Gauß' Formel aus den Euler'schen Formeln   (\ref{eq: Euler general summation formula}) und (\ref{eq: Euler Sum with Binomial}) mit leichter Rechnung folgt.

\paragraph{Gauß eigene Formeln}
\label{para: Gauß eingene Formeln}

Es sollen nun noch die Gauß allein zuzusprechenden Formeln der hypergeometrischen Reihe besprochen werden. Dies betrifft insbesondere die quadratische Transformation im letzten Argument der hypergeometrischen Reihe.  Diesen Gegenstand am  nächsten kommt Euler  in der Arbeit
\textit{``Consideratio aequationis differentio-differentialis $(a+bx)ddz+(c+ex)\frac{dxdz}{x}+(f+gx)\frac{zdx^2}{xx}=0$"} (\cite{E431}, 1773, ges. 1772) (E431: ``Betrachtung der Differenzen--Differentialgleichung $(a+bx)ddz+(c+ex)\frac{dxdz}{x}+(f+gx)\frac{zdx^2}{xx}=0$"), wo er die titelgebende Differentialgleichung betrachtet. Dass sich selbige für entsprechende Auswahl der Koeffizienten  auf die Differenzialgleichung für die hypergeometrische Reihe (\ref{eq: DGL Hyp}) reduziert, ist leicht einzushen. In seiner Untersuchung stellt Euler nun die Frage, eine Relation welcher Art  zwischen den Koeffizienten $a,b,c,e,f,g$ zu bestehen hat, sodass ein Potenzreihenansatz zur Lösung der titelgebenden Differentialgleichung eine polynomiale ist. Findet man zwei solcher Lösungen, kann man indes die vollständige Lösung konstruieren. In der Sprache der hypergeometrischen Funktion $_2F_1(\alpha,\beta,\gamma;z)$ sucht Euler also eine Relation zwischen den ersten drei Argumenten in $\alpha$,$\beta$ und $\gamma$, sodass man eine abbrechende Reihe erhält. Ähnliche Ideen sollten sich in der später verfassten Arbeit \cite{E710} als fruchtbar erweisen und ihn zur Gleichung

\begin{equation*}
    _2F_1(\alpha, \beta, \gamma,z)= (1-x)^{c-\alpha-\beta}{} _2F_1(\gamma -\alpha, \gamma-\beta, \gamma;z)
\end{equation*}
führen, welche entsprechend als Euler'sche Transformation der hypergeometrischen Reihe bezeichnet wird, und den Legendre--Polynomen leiten\footnote{Genauer führte ihn dies zu Erkenntnissen über Integraldarstellung (\ref{eq: Euler Family}), welche in Abschnitt (\ref{subsubsec: Den Kontext betreffend: Die Legendre Polynome}) auseinander gesetzt worden ist.}. Selbiger ließe sich noch die sog. Pfaff'sche Transformation zur Seite stellen:

\begin{equation*}
    _2F_1(\alpha, \beta, \gamma,z)= (1-x)^{c-\alpha}{} _2F_1\left(\gamma, \gamma-\beta, \gamma;\dfrac{z}{z-1}\right),
\end{equation*}
welche sich zwar nicht direkt bei Euler findet, jedoch unmittelbar aus den allgemeineren Transformationsformeln für Reihen fließt, welche er in § 7 vom zweiten Teil seiner \textit{Calculi Differentialis} \cite{E212} vorstellt. \\

Indes scheint  Euler nicht die Frage gestellt zu haben, ob sich in der hypergeometrische Reihen nicht gar das vierte Argument $z$ verändern ließe, ohne dass selbiges auf Differentialgleichung zutrifft, wenn freilich sie in der Variable $u=f(z)$ betrachtet wird\footnote{Diese Frage leitet   Kummer in seiner Arbeit  \cite{Ku36}  und lässt ihn überdies die Theorie der hypergeometrischen Reihe weiter ausbauen.}.  Gauß kann  aus seinen quadratischen Transformationen unter anderem seine zweite Summationsformel nachweisen:

\begin{equation}
    \label{eq: Gauss zweite Summationsformel}
    _2F_1 \left(a,b,\frac{1}{2}(a+b+1),\frac{1}{2}\right)= \dfrac{\Gamma\left(\frac{1}{2}\right)\Gamma \left(\frac{1}{2}(a+b+1)\right)}{\Gamma \left(\frac{1}{2}(a+1)\right)\Gamma \left(\frac{1}{2}(b+1)\right)},
\end{equation}
welche ein direktes Korollar aus (\ref{eq: Gauss Summation}) und (\ref{eq: Transformation Gauss quadratic}) darstellt, sich aber bei Euler nicht findet und auch aus seinen Resultaten nicht unmittelbar abgeleitet werden kann.

\paragraph{Euler und die konfluente hypergeometrische Funktion}
\label{para: Euler und die konfluente hypergeometrische Funktion}

Den Abschluss Euler'schen Beitrages zur hypergeometrischen Funktion soll ein Gegenstand einnehmen, welcher sich bei Euler findet, jedoch von Gauß in \cite{Ga13} nicht behandelt wird.  Geht man von (\ref{eq: Integral Rep}) aus und setzt zunächst $\frac{t}{\alpha}$ für $t$ und betrachtet im Anschluss $\alpha \rightarrow \infty$, findet man einen Ausdruck der Form

\begin{equation*}
     B(\beta, \gamma- \beta)  \int\limits_{0}^{1} t^{\beta-1}(1-t)^{\gamma -\beta -1}e^{\alpha t}dt,
\end{equation*}
wenn man den bekannten Grenzwert für $e^x$,

\begin{equation*}
    e^x=\lim_{n\rightarrow \infty} \left(1+\dfrac{x}{n}\right)^n,
\end{equation*}
verwendet. Analog kann man eine Funktion proportional zum Ausdruck

\begin{equation*}
    \int\limits_{0}^{\infty} t^{\beta -1}(1-t)^{\gamma-\beta -1}e^{\alpha t}dt
\end{equation*}
ableiten. Diese beiden letzten sind heute in der Theorie der konfluenten hypergeometrischen Funktion relevant. Eingeführt wurden sie auf die beschriebene Weise von Kummer in seiner Arbeit \textit{``De integralibus quibusdam definitis et seriebus infinitis"} (\cite{Ku37}, 1837) (``Über gewisse bestimmte Integrale und unendliche Reihen"). Euler indes scheint diese Reihen nie explizit in Bezug auf die hypergeometrische Funktion betrachtet zu haben. Man findet sie allerdings in gänzlich anderem Zusammenhang in seiner Arbeit \cite{E70} mitgeteilt.\\

\begin{figure}
    \centering
    \includegraphics[scale=0.8]{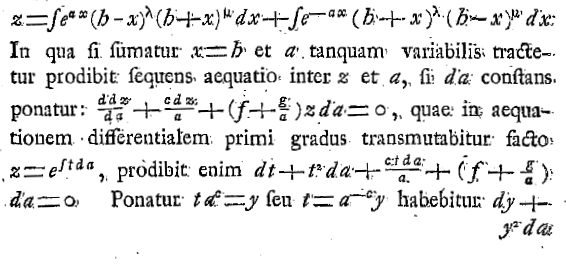}
    \caption{Euler konstruiert in seiner Arbeit \cite{E70} eine Differentialgleichung für einen Spezialfall der konfluenten hypergeometrischen Funktion.}
    \label{fig:E70KonfluenteHypergeometrische}
\end{figure}

 Während Kummer   von der Reihendarstellung der konfluenten hypergeometrischen Reihe aus beginnt, wählt Euler  ein Parameterintegral als Ausgangspunk und konstruiert aus selbigem die entsprechende Differentialgleichung. Sein allgemeinstes Ergebnis formuliert er in § 11 und 12 von \cite{E70}: Der allgemeinen Differentialgleichung (§ 11)

\begin{equation*}
    \dfrac{d^2z}{ds^2}+\left(\dfrac{\eta}{\varepsilon}+\dfrac{\theta}{\zeta}+\dfrac{\lambda+\mu-2}{a}\right)\dfrac{dz}{da}+\left(\dfrac{\eta(\mu+1)}{\varepsilon a}+\dfrac{\theta(\lambda+1)}{\zeta a}+\dfrac{\eta \theta}{\varepsilon \zeta}\right)z=A
\end{equation*}
mit beliebigen Konstanten $\eta,\varepsilon, \theta, \mu, \lambda$ und $A$ wird von einem Ausdruck der Form 

\begin{equation}
\label{eq: Euler Ansatz konfluent}
    z=E\int e^{ax}(\eta +\varepsilon x)^{\lambda}(\theta +\zeta x)^{\mu}dx+F\int e^{-ax}(\eta -\varepsilon x)^{\lambda}(\theta -\zeta x)^{\mu}dx.
\end{equation}
Genüge geleistet. Euler behauptet sogar noch vielmehr, dass sich für beliebige Wahl der Konstanten $E$,$F$  in der obigen Differentialgleichung sogar $A=0$ setzen lässt, sofern man eine gewisse Relation für die anderen freien Parameter voraussetzt, wovon Euler auch ein explizites Beispiel gibt.\\

Ein Vergleich mit der Differentialgleichung der konfluenten hypergeometrischen Reihe
\begin{equation*}
    z \dfrac{d^2w}{dz^2}+(b-z)\dfrac{dw}{dz}-aw=0,
\end{equation*}
welche neben der klassischen von Kummer gegebenen Reihe

\begin{equation}
\label{eq: Konfluente Reihe}
    M(a,b,z):=1+\frac{a}{b}\dfrac{z}{1!}+\dfrac{a(a+1)}{b(b+1)}\dfrac{z^2}{2!}+\dfrac{a(a+1)(a+2)}{b(b+1)(b+2)}\dfrac{z^3}{3!}+\cdots
\end{equation}
auch vom Integral:

\begin{equation*}
    M(a,b,z)= \dfrac{\Gamma(b)}{\Gamma(a)\Gamma(b-a)} \int\limits_{0}^1 e^{az}u^{a-1}(1-u)^{b-a-1}du \quad \text{für} \quad \operatorname{Re}(b)>\operatorname{Re}(a)>0
\end{equation*}
gelöst wird, hilft den Euler'schen Ansatz vom modernen Standpunkt aus zu verstehen. Weiter ist

\begin{equation*}
    U(a,b,z):=\dfrac{1}{\Gamma(a)}\int\limits_{0}^{\infty}e^{-zu}u^{a-1}(1+u)^{b-a-1}du \quad \text{für} \quad  \operatorname{Re}(a)>0
\end{equation*}
mit $M(a,b,z)$ bekanntermaßen über die Beziehung

\begin{equation*}
    U(a,b,z)= \dfrac{\Gamma(1-b)}{\Gamma(a+1-b)}M(a,b,z)+\dfrac{\Gamma(b-1)}{\Gamma(a)}z^{1-b}M(a+1-b,2-b,z)
\end{equation*}
verknüpft. Dies gibt also Anhaltspunkte, wie Euler seine Konstanten in (\ref{eq: Euler Ansatz konfluent}) und die Integrationsgrenzen auch hätte wählen können. Euler scheint die Verbindung zur hypergeometrischen Reihe nicht hergestellt zu haben, so findet man auch die Reihendarstellung (\ref{eq: Konfluente Reihe}) im Euler'schen Opus nicht explizit nieder geschrieben. Jedoch findet sich in der Arbeit \cite{E522} in den Paragraphen 34 und 35 die Kettenbruchentwicklung für den Quotienten zweier konsekutiver Ausdrücke der Form

\begin{equation*}
    f(n)= \int\limits_{0}^1 x^ne^{\alpha x}(1-x)^{\lambda}dx.
\end{equation*}
Selbiger lässt sich ohne Mühe als konfluente hypergeometrische Funktion deuten.

\begin{figure}
    \centering
     \includegraphics[scale=0.9]{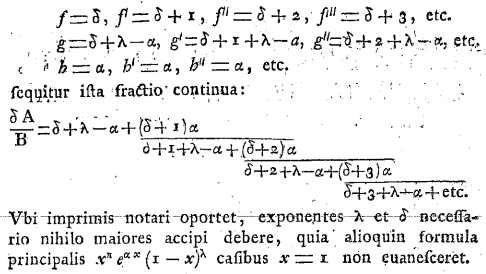}
    \caption{Euler gibt in seiner Arbeit \cite{E522} den Kettenbruch für zwei konsekutive konfluente konfluente hypergeometrische Reihe an. Hier ist in moderner Notation $A= \int\limits_{0}^1 x^{\delta -1}dx e^{\alpha x}(1-x)^{\lambda-1}$ und $B=\int\limits_{0}^1 x^{\delta}dx e^{\alpha x}(1-x)^{\lambda-1}$.}
    \label{fig:E522Kettenbruchkonfluente}
\end{figure}

\newpage

\section{Von Euler nicht bewiesene Entdeckungen}
\label{sec: Von Euler Entdecktes, aber nicht Bewiesenes}

\epigraph{An author never does more damage to his readers than when he hides a difficulty.}{Evariste Galois}

Während die Legendre--Polynome (Abschnitt \ref{subsubsec: Den Kontext betreffend: Die Legendre Polynome}) einige von Eulers Untersuchungen miteinander verknüpft hätten, dies aber für ihn selbst nicht zuletzt wegen des fehlenden Anreizes, diese Verbindungen zu knüpfen, nahezu unmöglich zu sehen war, verhält sich dies bei den folgenden Themen anders. Zum einen betrifft dies Theoreme, die Euler selbst hätte beweisen können (Abschnitt \ref{subsec: Von Euler selbst Beweisbares}), da sich das dafür Notwendige bereits in seinen Arbeiten befindet. Zum anderen betrifft  dies Entdeckungen, die er jedoch nicht mit einem Beweis zum Rang eines Theorems zu erheben vermochte (Abschnitt \ref{subsec: Von Euler Entdecktes, aber nicht Beweisbares}). 

\subsection{Von Euler selbst beweisbare Lehrsätze}
\label{subsec: Von Euler selbst Beweisbares}

\epigraph{Sometimes the situation is only a problem because it is looked at in a certain way. Looked at in another way, the right course of action may be so obvious that the problem no longer exists.}{Edward de Bono}

Als explizite Beispiele solcher Theoreme, welche Euler trotz bestehender Möglichkeit nicht bewiesen hat, werden eine Frage nach der einer Asymptotik (Abschnitt \ref{subsubsec: Durch Anwendung einer Methode: Seine Konstante A}) und die Funktionalgleichung der Riemnann'schen $\zeta$--Funktion angeführt (Abschnitt \ref{subsubsec: Durch Kombinieren von Ergebnissen: Die zeta-Funktion}).  Auch Fragestellungen, bei denen die Schaffung der nötigen Mittel nur einen Schritt entfernt lagen, welchen Euler nicht gegangen ist, finden hier am Exempel der $\vartheta$--Funktionen (Abschnitt \ref{subsubsec: Durch einen neuen Gedanken: Die Theta-Funktion}) ihre Erwähnung.

\subsubsection{Durch Anwendung einer Methode: Seine Konstante $A$}
\label{subsubsec: Durch Anwendung einer Methode: Seine Konstante A}

\epigraph{[W]hen solving a problem, we should always profit from previously solved problems, using their result or method, or experience aquicered in solving them.}{George Poyla}

Wie in der Abhandlung \textit{``Answer to a question concerning Euler’s paper "Variae considerationes circa series hypergeometricas”"} (\cite{Ay22}, 2022) kann die Euler'sche  Methode zur Auflösung homogener Differenzengleichungen mit linearen Koeffizienten (Siehe Abschnitt \ref{subsubsec: Die Mellin--Transformierte bei Euler} für die Erläuterung der Methode.) auch angewandt werden, um das spezifische von Faber  im Vorwort von Band 16,2 von Serie 1 der \textit{Opera Omnia} (\cite{OO162}, 1935) formulierte Problem zu lösen:\\

\textit{``Es würde sich vielleicht lohnen, die Euler'schen Bemühungen um die Ermittelung dieser Zahl $A$ wieder aufzunehmen."}\\

Dabei bezieht sich Faber auf Eulers Abhandlung \textit{``Variae considerationes circa series hypergeometricas"} (\cite{E661}, 1793, ges. 1776) (E661:``Verschiedene Betrachtungen über hypergeometrische Reihen"). Die Ätiologie dieser Konstante ist diese: In \cite{E661} betrachtet Euler die Funktion

\begin{equation}\label{eq: Euler'sche Darstellung}
     \Gamma_E(x) = a\cdot (a+b) \cdot (a+2b) \cdots (a+(x-1)b) \quad \text{mit} \quad \operatorname{Re(a)},\operatorname{Re}(b)>0,
\end{equation}
für welche er vermöge der Euler--Maclaurin'schen Summenformel nachstehende Asymptotik\footnote{Euler selbst benutzt in seiner Arbeit den Majuskel $\Gamma$. Um die Verwechslung mit der $\Gamma$--Funktion zu vermeiden, aber gleichzeitig möglichst nahe an der Euler'schen Notation zu bleiben, wird in dieser Arbeit $\Gamma_E$ für die in (\ref{eq: Euler'sche Darstellung}) eingeführte Funktion geschrieben.} ableitet

\begin{equation}\label{eq: Asymptotiken}
     \Gamma_E(x) \sim A \cdot e^{-x} (a-(b-1)x)^{\frac{a}{b}+x-\frac{1}{2}},
\end{equation}
welche Asymptotik für den Grenzwert $x \rightarrow \infty$ gilt. Nun soll die Konstante $A$ gefunden werden, jedoch schreibt Euler in § 17: \\

\textit{``Indes ist daraus dennoch nicht ersichtlich, wie die Konstante $A$ als absoluter Wert bestimmt werden kann, [...]."}

\begin{figure}
    \centering
   \includegraphics[scale=0.8]{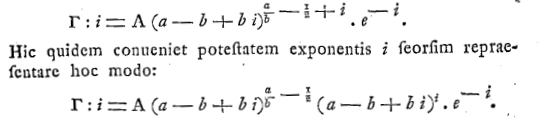}
    \caption{Euler leitet in seiner Arbeit \cite{E661} die Asymptotik seiner verallgemeinerten Fakultät her.}
    \label{fig:E661GeneralGamma}
\end{figure}

\paragraph{Ermittelung der Konstanten aus Eulers Formeln heraus}

Bemerkt man nun, dass (\ref{eq: Euler'sche Darstellung}) die Gleichung

\begin{equation} \label{eq: Funktionalgleichung}
    \Gamma_E(x+1)=(a+bx)\Gamma_E(x)
\end{equation}
erfüllt, lässt sich Eulers Ansatz zur Lösung von homogenen Differenzengleichungen mit linearen Koeffizienten auf diese anwenden (Abschnitt \ref{subsubsec: Die Mellin--Transformierte bei Euler}). Man findet auf diesem Wege

\begin{equation}
\label{eq: GammaE und Gamma}
    \Gamma_E(x)= \dfrac{b^x}{\Gamma\left(\frac{a}{b}\right)}\cdot\Gamma \left(x+\dfrac{a}{b}\right),
\end{equation}
wobei hier $\Gamma$ wieder die bekannte $\Gamma$--Funktion bedeutet.
Zur Ermittlung der Konstante $A$ in (\ref{eq: Asymptotiken}), gebrauche man die Stirling'sche Formel\footnote{Selbige hat Euler wie an entsprechender Stelle erwähnt (Abschnitt \ref{subsubsec: Stirling'sche Formel nach Euler}) in seinen \textit{Calculi Differentialis} \cite{E212} korrekt angegeben, weshalb sie in dieser Beweisführung ihre Verwendung finden darf.}

\begin{equation*}
    \Gamma(x+1) \sim \sqrt{2\pi x}\cdot x^x \cdot e^{-x} \quad \text{für} \quad x \rightarrow \infty.
\end{equation*}
oder auch, falls man $x-1$ statt $x$ schreibt,

\begin{equation*}
    \Gamma(x) \sim  \sqrt{2\pi}(x-1)^{\frac{1}{2}}\cdot (x-1)^{x-1}e^{-(x-1)} \quad \text{für} \quad x \rightarrow \infty.
\end{equation*}
Dementsprechend entspringt aus (\ref{eq: GammaE und Gamma})  die Asymptotik

\begin{equation*}
   \Gamma_E\left(x\right)\sim \dfrac{b^x}{\Gamma\left(\frac{a}{b}\right)}\cdot \sqrt{2\pi}\left(x-1+\dfrac{a}{b}\right)^{\frac{1}{2}}\left(x+\dfrac{a}{b}-1\right)^{x+\frac{a}{b}-1}e^{-\left(x+\frac{a}{b}-1\right)}.
\end{equation*}
was durch Vereinfachen zu 

\begin{equation*}
    \Gamma_E(x) \sim \dfrac{\sqrt{2\pi}}{\Gamma\left(\frac{a}{b}\right)}\cdot b^{-\frac{1}{2}-\frac{a}{b}+1}\left(a+bx-b\right)^{-\frac{1}{2}+x+\frac{a}{b}}e^{-x} \cdot e^{-\frac{a}{b}+1}
\end{equation*}
wird, sodass 

\begin{equation}
\label{eq: Constant A}
    A= \dfrac{\sqrt{2\pi}}{\Gamma \left(\frac{a}{b}\right)}\cdot e^{-\frac{a}{b}+1}\cdot b^{\frac{1}{2}-\frac{a}{b}}.
\end{equation}
Somit hätte Euler das sich selbst gesetzte Problem ebenfalls auflösen können.

\paragraph{Ein Korollar seiner Untersuchungen: Die Legendre'sche Verdopplungsformel}

Überdies, wie auch in der Arbeit \textit{``Euler and the Duplication Formula for the Gamma-Function"} (\cite{Ay23c}, 2023) gezeigt, lässt sich aus den Euler'schen Überlegungen und nach Ermittlung der Konstante $A$ die Legendre'sche Verdopplungsformel nachweisen. Dies bedeutet die Formel:

\begin{equation}
\label{eq: Legendre Doppel}
    \Gamma(x)= \dfrac{2^{x-1}}{\sqrt{\pi}}\Gamma \left(\dfrac{x}{2}\right)\Gamma\left(\dfrac{x}{2}+1\right).
\end{equation}
Dies soll hier der Vollständigkeit wegen bündig demonstriert werden. Abgesehen von der Funktion $\Gamma_E$, führt Euler in seiner Arbeit \cite{E661} zwei eng mit dieser verwandte Funktionen ein:

\begin{equation}
\label{eq: Euler's other Functions}
        \renewcommand{\arraystretch}{1,5}
\setlength{\arraycolsep}{0.0mm}
\begin{array}{llllll}
    \Delta(x) &~=~&~ a \cdot (a+2b)\cdot (a+4b)\cdot (a+6b)\cdot \cdots \cdot (a+(2x-2)b), \\
    \Theta(x) &~=~&~ (a+b)\cdot (a+3b) \cdot (a+5b) \cdot \cdots \cdot (a+(2x-1)b).
\end{array}
\end{equation}
In gänzlich gleicher Weise wie bei $\Gamma_E$ findet Euler die asymptotischen Entwicklungen für seine Funktionen $\Delta$ und $\Theta$:

\begin{equation}
\label{eq: Delta}
        \renewcommand{\arraystretch}{1,5}
\setlength{\arraycolsep}{0.0mm}
\begin{array}{llllll}
     \Delta(x) &~\sim ~&~ B\cdot e^{-x} \cdot (a-2b+2bx)^{\frac{a}{2b}+x-\frac{1}{2}} \\
     \Theta(x) &~\sim ~&~ C \cdot e^{-x} \cdot (a-b+2bx)^{\frac{a}{2b}+x},
\end{array}
\end{equation}
wo $B$ und $C$ die Rolle der Konstante gleich der von $A$ bei $\Gamma_E$ einnehmen, weswegen die Asymptotiken gleichermaßen für den Fall $x \rightarrow \infty$ zu verstehen sind.  Wie bereits angesprochen, sieht Euler sich nicht in der Lage auch nur eine der Konstanten $A$, $B$ und $C$ zu bestimmen, es gelingt ihm aber aus den Relationen zwischen $\Gamma_E$, $\Delta$ und $\Theta$ entsprechende zwischen den Konstanten abzuleiten:

\begin{equation}
    \label{eq: Relation among the Constants A}
   A= \dfrac{B \cdot C}{\sqrt{e}} 
\end{equation}
und

\begin{equation}
     \label{eq: Relation among the Constants B}
     B = C\cdot k \cdot \sqrt{e}
\end{equation}
mit $ k=\Delta\left(\frac{1}{2}\right)$. Hinter diesen Relationen verbirgt sich bereits die Legendre'sche Verdopplungsformel (\ref{eq: Legendre Doppel}), wozu es aber zunächst der Konstante $k$ bedarf. \\

Um selbige ausfindig zu machen, bemerke man zuerst, dass die Ersetzung $b \mapsto 2b$ in $\Gamma_E$ den Ausdruck für $\Delta$ aus (\ref{eq: Delta}) entspringen lässt. Gleiches in der expliziten Formel (\ref{eq: GammaE und Gamma}) durchgeführt produziert:

\begin{equation*}
    \Delta(x) = \dfrac{(2b)^x}{\Gamma\left(\frac{a}{2b}\right)}\cdot \Gamma \left(x+\dfrac{a}{2b}\right).
\end{equation*}
Demnach für $x=\frac{1}{2}$

\begin{equation}
\label{eq: Value k}
    k = \Delta \left(\dfrac{1}{2}\right)=  \dfrac{(2b)^{\frac{1}{2}}}{\Gamma\left(\frac{a}{2b}\right)}\cdot \Gamma \left(\dfrac{1}{2}+\dfrac{a}{2b}\right).
\end{equation}
Mit diesem Wert von $k$ ausgestattet nutze man (\ref{eq: Relation among the Constants A}) sowie (\ref{eq: Relation among the Constants B}). Setzt man den Wert von $C$ aus zweiter Gleichung in die erste ein, gibt dies:

\begin{equation}
\label{eq: Relation AB}
    A= \dfrac{B^2}{\Delta \left(\frac{1}{2}\right)}e^{-1}.
\end{equation}
Da sich $\Delta(x)$ aus der Substitution $b \mapsto 2b$ aus $\Gamma_{E}(x)$ ergibt, gilt selbiges für den Wert $B$ aus $A$:

\begin{equation}
    \label{eq: Constant B}
    B = \dfrac{\sqrt{2\pi}}{\Gamma \left(\frac{a}{2b}\right)}\cdot (2b)^{\frac{1}{2}-\frac{a}{2b}}\cdot e^{1-\frac{a}{2b}}.
\end{equation}
Einsetzen der Werte von $A$ aus (\ref{eq: Constant A}), $B$ aus (\ref{eq: Constant B}) und $k$ aus (\ref{eq: Value k}) wandelt die Relation (\ref{eq: Relation AB}) in:

\begin{equation*}
    \dfrac{\sqrt{2\pi}}{\Gamma\left(\frac{a}{b}\right)}\cdot e^{1-\frac{a}{b}}\cdot b^{\frac{1}{2}-\frac{a}{b}} = \dfrac{\left(\frac{\sqrt{2\pi}}{\Gamma \left(\frac{a}{2b}\right)}\cdot (2b)^{\frac{1}{2}-\frac{a}{2b}}\cdot e^{1-\frac{a}{2b}}\right)^2}{\frac{(2b)^{\frac{1}{2}}}{\Gamma\left(\frac{a}{2b}\right)}\cdot \Gamma \left(\frac{1}{2}+\frac{a}{2b}\right)} \cdot e^{-1}.
\end{equation*}
Viele Terme streichen sich, sodass verbleibt

\begin{equation*}
    \dfrac{1}{\Gamma \left(\frac{a}{b}\right)}= \dfrac{\sqrt{2\pi}\cdot 2^{\frac{1}{2}-\frac{a}{b}}}{\Gamma \left(\frac{a}{2b}\right)\cdot \Gamma \left(\frac{1}{2}+\frac{a}{2b}\right)}.
\end{equation*}
Schreibt man schließlich $x$ anstatt $\frac{a}{b}$ und löst nach $\Gamma(x)$, steht nach einiger Vereinfachung:

\begin{equation*}
    \Gamma(x) = \dfrac{2^{x-1}}{\sqrt{\pi}}\cdot \Gamma \left(\dfrac{x}{2}\right)\cdot \Gamma \left(\dfrac{x+1}{2}\right),
\end{equation*}
was gerade die Legendre'sche Verdopplungsformel (\ref{eq: Legendre Doppel}) ist.

\subsubsection{Durch Kombinieren von Ergebnissen: Die $\zeta$-Funktion}
\label{subsubsec: Durch Kombinieren von Ergebnissen: Die zeta-Funktion}

\epigraph{If you search everywhere, yet cannot find what you are seeking, it is because what you seek is already in your possession.}{Laotse}

Dass Euler in seiner Arbeit \cite{E352} als erster die Funktionalgleichung für die Riemann'sche $\zeta$-Funktion,

 \begin{equation}
     \label{eq: Def zeta}
     \zeta(s):= \sum_{n=1}^{\infty} \dfrac{1}{n^s}, \quad \text{für} \quad \operatorname{Re}(s)>1,
 \end{equation}
nieder geschrieben hat, ist heute sehr bekannt. Man siehe etwa \cite{Ha49}\footnote{Auf Seite 22 schreibt Hardy  Edmund Landau (1877--1938) und Cahen (1865--1941) das Verdienst zu, als erster erwähnt zu haben, dass Euler die Funktionalgleichung der Riemann'schen $\zeta$-Funktion in \cite{E352} angegeben hat. Hardy referiert dabei wohl auf Landaus Arbeit \textit{``Euler und die Funktionalgleichung der Riemannschen Zetafunktion"} (\cite{La06}, 1906). Dieser schreibt allerdings direkt zu Beginn, dass Cahen in seiner Arbeit \textit{``Sur la fonction $\zeta(s)$ de Riemann et sur les fonctions analogues"} (\cite{Ca94}, 1894) (``Über die Funktion $\zeta(s)$ von Riemann und über analoge Funktion") bereits vor ihm die Euler'sche Erstentdeckerschaft der Funktionalgleichung bemerkt hat.} oder konsultiere den Übersichtsartikel \cite{Ay74}.\\

Ergänzend sei vorausgeschickt, dass Euler die Funktionalgleichung für die Dirichlet'sche $\eta$--Funktion

\begin{equation}
\label{eq: Def eta}
    \eta(s):= \sum_{n=1}^{\infty} \dfrac{(-1)^{n+1}}{n^s} \quad \text{für} \quad \operatorname{Re}(s)>0.
\end{equation}
und nicht die für die $\zeta$--Funktion angibt. Sie lautet in moderner Darstellung

\begin{equation}
\label{eq: Functional Equation Eta}
    \eta(1-s)=  \dfrac{2^s-1}{1-2^{s-1}}\pi^{-s}\cos \left(\dfrac{\pi s}{2}\right)\Gamma(s)\eta(s).
\end{equation}
Euler sagt explizit, dass er diese Gleichung zwar nicht streng beweisen kann, jedoch seine Ausführungen ihn von ihrer Richtigkeit überzeugen\footnote{Diese Gleichung ist bereits oben (Abschnitt \ref{subsubsec: Eulers Auffassung eines Beweises}) bei der Diskussion der Euler'schen Auffassung eines Beweises aufgetreten.}. Entdeckt hat er sie mit routiniertem Umgang mit divergenten Reihen, die Euler'sche Auffassung zu welchen weiter unten (Abschnitt \ref{subsubsec: Der Begriff der Summe einer Reihe}) umfassender auseinander gesetzt werden wird.

\paragraph{Eulers Anwendung einer divergenten Reihe -- Die Zeta-Funktion}
\label{para: Eulers Anwendung einer divergenten Reihe}

Es sollen nur die Grundbausteine seiner Gedanken aus \cite{E352} an dieser Stelle vorgestellt werden, für Details ist auf \cite{Ay74} oder \cite{Va06} verwiesen. Die fundamentale Idee  zur Erlangung der Funktionalgleichung für die $\eta$--Funktion besteht in der Verwendung  zweier die Bernoulli--Zahlen involvierender Identitäten\footnote{Auch Ramanujan hat überdies die Funktionalgleichung der $\zeta$--Funktion auf diese Weise abgeleitet. Man siehe etwa das Buch \textit{``Ramanujan's Notebooks -- Part I} (\cite{Be85}, 1985).}. Bereits aus  \cite{E130} weiß Euler, dass gilt:

\begin{equation*}
    \eta(2n)= (-1)^{n+1} \dfrac{2^{2n}-1}{(2n)!}B_{2n}\pi^{2n},
\end{equation*}
welche er später  interpolieren möchte, weswegen er sie  wie folgt ausdrückt:

\begin{equation*}
    \eta(2n) = -\dfrac{1}{\cos \left(n\pi\right)} \dfrac{2^{2n}-1}{\Gamma(2n+1)}B_{2n}\pi^{2n},
\end{equation*}
somit kann Euler auch schreiben:

\begin{equation}
    \label{eq: Eta(n) Euler}
    \eta(n) = -\dfrac{1}{\cos \left(\dfrac{n\pi}{2}\right)} \dfrac{2^{n}-1}{\Gamma(n+1)}B_{n}\pi^{n}.
\end{equation}
Mithilfe der Methode der Abel'schen Summation\footnote{Bei dieser wird eine Summe wie folgt als Grenzwert definiert:
\begin{equation*}
    s:= \lim_{x \rightarrow 1^{-}} \sum_{n=1}^{\infty}a_n x^n,
\end{equation*}
wobei statt einer Potenzreihe prinzipiell auch andere Funktionenreihen stehen können, wie etwa Fourier-- oder Dirichletreihen usw.}, aber auch der Euler-Maclaurin'schen Summationsformel, gelangt er in \cite{E352} zur Formel für $\eta(1-n)$:

\begin{equation*}
    \eta(1-n)= \dfrac{2^n-1}{n}B_n.
\end{equation*}

\begin{figure}
    \centering
   \includegraphics[scale=0.7]{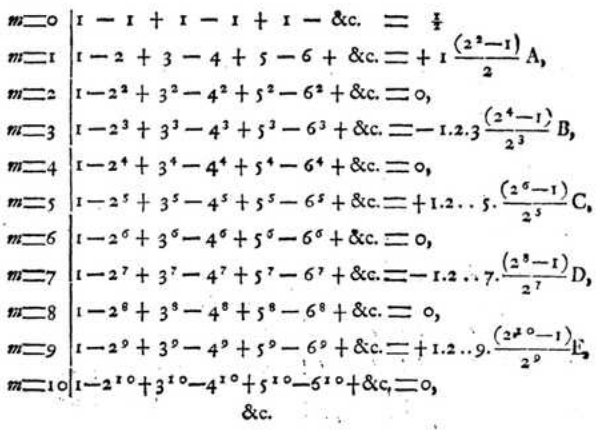}
    \caption{Euler gibt in seiner Arbeit \cite{E352} die Werte einiger divergenten Reihen an. Die Zahlen $A$, $B$, $C$, $D$, $E$ etc. sind die Bernoulli-Zahlen.}
    \label{fig:E352Divergent}
\end{figure}

Also gibt eine Division der beiden letzten Formeln, gerade die Gleichung:

\begin{equation*}
    \dfrac{\eta(1-n)}{\eta(n)}= -\dfrac{\Gamma(n)}{\pi^n} \cdot \dfrac{2^n-1}{2^{n-1}-1}\cdot \cos \left(\dfrac{n\pi}{2}\right), 
\end{equation*}
was gerade die Version von (\ref{eq: Functional Equation Eta}) ist, welche schon in Abschnitt (\ref{subsubsec: Eulers Auffassung eines Beweises}) erwähnt worden ist. 

\paragraph{Beweis der Funktionalgleichung aus Eulers Formeln}
\label{para: Beweis der Funktionalgleichung aus Eulers Formeln}

Der Beweis von (\ref{eq: Functional Equation Eta}) benötigt zwei Formeln von Euler aus anderer Quelle. Die erste stammt aus § 12 von Eulers Arbeit \textit{``De resolutione fractionum transcendentium in infinitas fractiones simplices"} (\cite{E592}, 1785, ges. 1775) (E592: ``Über die Auflösung von transzendenten Brüchen in unendliche viele einfache Brüche"), sie lautet

\begin{equation}
\label{eq: frac sin}
    \dfrac{\lambda \pi}{2\sin \lambda \pi}= \sum_{n=1}^{\infty} \dfrac{(-1)^{n+1}n^2}{n^2 -\lambda^2}.
\end{equation}
Über die Gleichung

\begin{equation*}
    \sin x= \dfrac{e^{ix}-e^{-ix}}{2i}
\end{equation*}
findet man daraus leicht 

\begin{equation}
\label{eq: Partial e}
    \dfrac{u}{e^u-e^{-u}}= \sum_{n=1}^{\infty} (-1)^{n+1} \dfrac{n^2}{n^2+\frac{u^2}{\pi^2}}= \sum_{n=1}^{\infty} \dfrac{(-1)^{n+1}}{1+\left(\frac{u}{\pi n}\right)^2}.
\end{equation}
Die zweite grundlegende Formel entnehme man aus § 8 (aus der \textit{Opera Omnia} Version\footnote{Die Originalarbeit weist eine andere Nummerierung der Paragraphen auf. Hier befindet sich die Formel in § 131.}) der Arbeit \textit{``De valoribus integralium a termino variabilis $x = 0$ usque ad $x = \infty$ extensorum"} (\cite{E675}, 1794, ges. 1778) (E675: ``Über die Werte von Integralen, die von $x=0$ bis zu $x=\infty$ erstreckt worden sind"):

\begin{equation*}
    \int y^{n-1}\partial y e^{-ky}=\dfrac{\triangle}{k^n}, \quad \text{mit} \quad \triangle=1 \cdot 2 \cdot 3 \cdot 4 \cdots (n-1).
\end{equation*}
Hier ist das Integral von $y=0$ bis hin zu $y=\infty$ zu nehmen, $n$ sieht Euler als beliebige natürliche Zahl größer als $1$ und $k$ soll sogar eine beliebige komplexe Zahl sein\footnote{$k$ wäre noch der Einschränkung $\operatorname{Re}(k)>0$ zu versehen, um die Konvergenz des Integrals zu gewährleisten, was Euler indes nicht bemerkt zu haben scheint.}.\\

\begin{figure}
    \centering
    \includegraphics[scale=0.8]{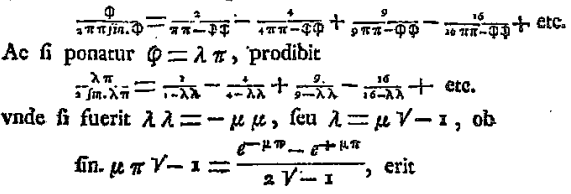}
    \caption{Euler gibt in der Arbeit \cite{E592} die Partialbruchzerlegung der Funktion $\frac{\varphi}{\sin (\varphi)}$ an.}
    \label{fig:E592Fig1original}
\end{figure}

Man kann Eulers Formel  wie folgt schreiben

\begin{equation}
\label{eq: Gamma Start}
    \dfrac{\Gamma(x)}{k^x} = \int\limits_{0}^{\infty} e^{-ku}u^{x-1}du,
\end{equation}
wobei aus Gründen besserer Übersichtlichkeit $n$ zu $x$ geändert ist;  $x$ darf freilich auch eine reelle Zahl sein. Aus dieser Formel ergibt sich nun direkt

\begin{equation*}
      \dfrac{\Gamma(x)}{(2n+1)^x} = \int\limits_{0}^{\infty} e^{-(2n+1)u}u^{x-1}du.
\end{equation*}
Summieren über $n$ liefert

\begin{equation*}
    \sum_{n=0}^{\infty} \dfrac{\Gamma(x)}{(2n+1)^x}= \Gamma(x)\sum_{n=0}^{\infty} \dfrac{1}{(2n+1)^x}=:\Gamma(x)\lambda(x),
\end{equation*}
wobei

\begin{equation}
    \label{eq: Def Lambda}
    \lambda(x):=\sum_{n=0}^{\infty}\frac{1}{(2n+1)^x},
\end{equation}
eingeführt wurde. Nun gilt aber auch

\begin{equation*}
 \sum_{n=0}^{\infty}   \int\limits_{0}^{\infty} e^{-(2n+1)u}u^{x-1}du = \int\limits_{0}^{\infty} \sum_{n=0}^{\infty}\left(  e^{-(2n+1)u}\right)u^{x-1}du 
\end{equation*}
\begin{equation*}
    = \int\limits_{0}^{\infty} \dfrac{ u^{x-1}e^{-u}du}{1-e^{-2u}}= \int\limits_{0}^{\infty}  \dfrac{u^{x-1}du}{e^{u}-e^{-u}}.
\end{equation*}
Die Vertauschung von Summation und Integration ist leicht zu rechtfertigen; im zweiten Schritt wurde die geometrische Reihe angewandt. Insgesamt hat man:

\begin{equation}
\label{eq: Intermediate}
    \lambda(x)\Gamma(x)= \int\limits_{0}^{\infty} u^{x-2} \cdot\dfrac{udu}{e^u-e^{-u}},
\end{equation}
sodass  nun die Partialbruchzerlegung (\ref{eq: Partial e}) ihre Anwendung finden kann, was zu

\begin{equation*}
     \lambda(x)\Gamma(x)=\int\limits_{0}^{\infty} u^{x-2}du \sum_{n=1}^{\infty} \dfrac{(-1)^{n+1}}{1+\left(\frac{u}{\pi n}\right)^2}.
\end{equation*}
führt. Substituiert man nun $\frac{u}{\pi n}=v$, gelangt man zu

\begin{equation*}
    \lambda(x)\Gamma(x)= \int\limits_{0}^{\infty} \sum_{n=1}^{\infty} (-1)^{n+1}\pi^{x-1}n^{x-1} \cdot \dfrac{v^{x-2}dv}{1+v^2}  =\pi^{x-1}\eta(1-x)\int\limits_{0}^{\infty} \dfrac{v^{x-2}dv}{1+v^2},
\end{equation*}
wobei  im letzten Schritt die Definition der $\eta$-Funktion (\ref{eq: Def eta}) genutzt wurde. Das verbleibende Integral berechnet Euler ebenfalls an entsprechender Stelle; es tritt etwa in der Arbeit \textit{``Investigatio formulae integralis $\int \frac{x^{m-1} dx}{(1+x^k)^n}$ casu, quo post intagrationem statuitur $x = \infty$"} (\cite{E588}, 1785, ges. 1775) (E588: ``Untersuchung des Integrals $\int \frac{x^{m-1} dx}{(1+x^k)^n}$ im Fall, im welchem man nach der Integration $x=\infty$ setzt") auf. Daraus fließt die Formel:

\begin{equation*}
    \int\limits_{0}^{\infty} \dfrac{v^{x-2}dv}{1+v^2}= -\dfrac{1}{2} \dfrac{\pi}{\cos \left(\frac{\pi x}{2}\right)},
\end{equation*}
welche in der letzten Gleichung eingesetzt gibt:

\begin{equation*}
    \Gamma(x)\lambda(x)= -\dfrac{1}{2} \dfrac{\pi}{\cos \left(\frac{\pi x}{2}\right)} \pi^{x-1}\eta(1-x).
\end{equation*}
Unter Verwendung der elementaren Identität

\begin{equation*}
    \lambda(x)= \dfrac{2^x-1}{2^x-2}\eta(x)
\end{equation*}
gelangt man schlussendlich zu

\begin{equation*}
    \eta(1-x)= \dfrac{2^x-1}{1-2^{x-1}}\pi^{-x}\cos \left(\dfrac{\pi x}{2}\right)\Gamma(x) \eta(x),
\end{equation*}
was  gerade (\ref{eq: Functional Equation Eta}) ist. Somit ist ein Beweis geführt, welcher sich lediglich von  Euler selbst eruierter Theoreme bedient.

\paragraph{Bezug zur $\zeta$-Funktion}
\label{para: Bezug zur zeta-Funktion}

Die explizite Betonung, Euler habe die $\eta$--Funktion und nicht die $\zeta$--Funktion betrachtet, bedarf angesichts der sie verknüpfenden und Euler auch bekannten (siehe § 170 seiner \textit{Introductio} \cite{E101}) Formel 
\begin{equation}
 \label{eq: etazeta}
     \eta(s)= (1-2^{1-s})\zeta(s)
 \end{equation}
weiterer Erläuterung. Freilich impliziert die Funktionalgleichung (\ref{eq: Functional Equation Eta}) unmittelbar die der $\zeta$-Funktion: 

  \begin{equation}
     \label{eq: Func zeta}
     \zeta(s)= 2^s \pi^{s-1} \sin \left(\dfrac{\pi s}{2}\right)\Gamma(1-s)\zeta(1-s).
 \end{equation}
 Einen direkten Beweis dieser Gleichung analog zu dem der $\eta$-Funktion misslingt allerdings und insgesamt scheint ein direkter Beweis aus den Euler'schen Formeln nicht möglich zu sein\footnote{Zum Zeitpunkt der Abgabe der hiesigen Ausführungen befindet sich eine Arbeit des Verfassers mit einem solchen Beweisversuch unter Begutachtung.}. Auch seine Definition der Summe eine divergenten Reihe (siehe Abschnitt \ref{subsubsec: Der Begriff der Summe einer Reihe}) hätte ihm diesbezüglich nicht zur Hilfe gereicht\footnote{Seine Auffassung einer divergenten Reihen zwingt Euler, $\zeta(-k)$ für $k \in \mathbb{N}$ alle zu $\infty$ zu summieren, was natürlich unrichtig ist.}. Vielmehr ließe sich argumentieren, dass Euler zwar  (\ref{eq: Functional Equation Eta}) in just aufgezeigter Manier von einer Vermutung in die Form eines Theorems hätte überführen können, aber selbiges für (\ref{eq: Func zeta}) nicht hätte leisten können, ohne sich mit Widersprüchen seiner Theorie divergenter Reihen konfrontiert zu sehen. Nichtsdestotrotz darf konstatiert werden, dass Euler alles über die $\zeta$--Funktion zutage gefördert hat, was möglich ist, wenn man sie als Funktion einer \textit{reellen} Variablen betrachtet. Die Auflösung der noch bestehenden Schwierigkeiten bedarf die Kenntnisse der komplexen Analysis, welche Euler jedoch aus unten (Abschnitt \ref{subsubsec: Komplexe Analysis}) auszuführenden Gründen verschlossen bleiben musste.

 \paragraph{Bemerkung zu Riemann und der Zeta-Funktion}

Auch der von Riemann in seiner Arbeit (\textit{``Über die Anzahl der Primzahlen unter einer gegebenen Größe"} \cite{Ri60}, 1860, ges. 1859)  mitgeteilte Beweis von (\ref{eq: Func zeta}) greift zu Mitteln, welche Euler nicht zur Verfügung standen. Riemann gibt in seiner Arbeit zwei Beweise, von denen ersterer auf Integration in der komplexen Ebene fußt, der zweite verwendet hingegen eine Transformationsformel für die Jacobi'sche Thetafunktion. Dies  stellt in diesem Zusammenhang die Identität

\begin{equation*}
    \vartheta(-\dfrac{1}{\tau}) = \sqrt{\dfrac{\tau}{i}}\vartheta(\tau)
\end{equation*}
dar, wobei 

\begin{equation}
\label{eq: Theta}
    \vartheta(\tau)= \sum_{n=-\infty}^{\infty} e^{i\pi n^2 \tau}
\end{equation}
ist. Diese Transformationsformel ergibt sich aus den Ergebnissen von Jacobis Arbeit \textit{``Fundamenta nova theoriae functionum ellipticarum"} (\cite{Ja29}, 1829) (``Neue Grundlagen der Theorie der elliptischen Funktionen"). Sie findet sich aber, wie auch Jacobi selbst in seiner Notiz \textit{``Suite des Notices sur les Functions Ellitiques"} (\cite{Ja28}, 1828) (``Weitere Notizen zu elliptischen Funktionen") zeitlich vor eben erwähnter verfassten Arbeit anmerkt, schon in S. Poissons (1781--1840) Arbeit \textit{``Suite du Mémoire sur les Intégrales définies et sur la sommation des Séries"} (\cite{Po23}, 1823) (``Fortsetzung zur Abhandlung über bestimmte Integrale und die Summierung von Reihen"). Dort steht sie auf Seite 284.\\

\begin{figure}
    \centering
     \includegraphics[scale=0.8]{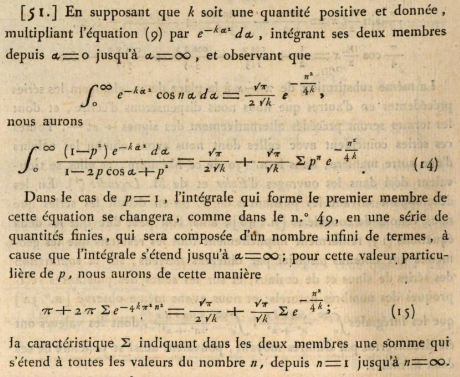}
    \caption{Poission gelangt in seiner Arbeit \cite{Po23} zur Funktionalgleichung für die Thetafunktion (\ref{eq: Theta}). Es ist Gleichung (15) in seiner Arbeit.}
    \label{fig:PoissonTheta}
\end{figure}

Man mag sich an dieser Stelle wundern, ob Riemann Kenntnis von der Euler'schen Arbeit \cite{E352} oder zumindest \cite{E432} hatte. Er macht jedenfalls  keine direkte Erwähnung davon, lediglich das Euler--Produkt (\ref{eq: Euler Produkt}) nennt er explizit, was also wohl auf das Riemann'sche Studium der Arbeit \cite{E72} oder -- unlängst wahrscheinlicher --   der \textit{Introductio} \cite{E101} hindeutet.\\

Es darf jedoch berechtigt davon ausgegangen werden, dass Riemann von einer externen Quelle auf die Funktionalgleichung hingewiesen worden ist. Denn unmittelbar auf den von Riemann angezeigten Integrationsweg zum Beweis der Funktionalgleichung zu stoßen, wäre ein übergroßer Zufall. Auch das Umschreiben auf eine Form, welche die $\vartheta$--Funktion verwenden lässt, setzt eine Intention mit dem Endergebnis im Blick voraus. Riemanns Anregungen mögen tatsächlich auch die erstgenannten Euler'schen Arbeiten gewesen sein, oder aber auch die von J. Malmstén (1814--1886) \textit{``De integralibus quibusdam definitis seriebusque infinitis"} (\cite{Ma46}, 1846) (``Über gewisse bestimmte Integrale und unendliche Reihen"), welche allerdings nur die Funktionalgleichung für die sogenannte Dirichlet'sche $\beta$-Funktion, definiert als

 \begin{equation}
 \label{eq: Def Beta}
     \beta(s):= \sum_{n=0}^{\infty} \dfrac{(-1)^n}{(2n+1)^s} \quad \text{für} \quad \operatorname{Re}(s)>0
 \end{equation}
 enthält. Die entsprechende Gleichung findet sich auch bei Euler in \cite{E352}. \\

\begin{figure}
    \centering
       \includegraphics[scale=0.8]{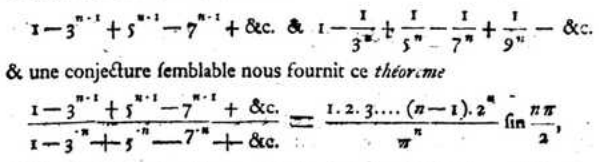}
    \caption{Euler formuliert in seiner Arbeit \cite{E352} die Funktionalgleichung für die heute sogenannte Dirichelt'sche $\beta$-Funktion, definiert über  (\ref{eq: Def Beta}), als Vermutung und gleichzeitig als Theorem.}
    \label{FuncBetaE352Original}
\end{figure}

 Modern ausgedrückt liest sich die Euler'sche Formulierung wie folgt:

 \begin{equation}
     \label{eq: Functional Equation Beta}
     \beta(1-s)= \left(\dfrac{\pi}{2}\right)^{-s}\sin \left(\dfrac{\pi s}{2}\right)\beta(s) \Gamma(s).
 \end{equation}
 Der Beweis kann analog zu dem für die $\eta$--Funktion geführt werden, die nötige Partialbruchzerlegung findet sich  in § 29 von \cite{E592}. Malmstén indes sagt in seiner zitierten Arbeit, dass er die Funktionalgleichung ``irgendwo bei Euler"{} gesehen hat. Er referiert wohl \cite{E352}, da in \cite{E432} die $\beta$--Funktion nicht genannt wird. Der Vollständigkeit wegen  sei noch die Arbeit \textit{``On Eisenstein’s Copy of the Disquisitiones"} (\cite{We89}, 1989) von André Weil (1906--1990) erwähnt; selbige argumentiert überzeugend, dass Eisenstein (1829--1852) bereits im Besitz eines Beweises von letztgenannter Funktionalgleichung (\ref{eq: Functional Equation Beta}) war und so Riemanns Untersuchung ihren Ursprung gegeben haben könnte. Das Eisenstein'sche Argument, datiert auf das Jahr 1849, ist ebenfalls, wie das entsprechende Riemnann'sche, auf die Funktionalgleichung der Jacobi'schen $\vartheta$--Funktion gestützt. \\

 Oft, wie etwa von Speiser im Vorwort von Band 9 der 1. Serie der \textit{Opera Omnia} (\cite{OO9}, 1945), wird Oskar Schlömilch (1823--1901) als Entdecker von (\ref{eq: Func zeta}) genannt. Indes scheint Schlömilch diese Gleichung nicht gezeigt zu haben. Die Referenzen Speisers führen zur Arbeit (\cite{Sc58}, 1858); in dieser Arbeit betrachtet der Autor jedoch ebenfalls die $\beta$-Funktion (\ref{eq: Def Beta}), liefert aber einen neuen Beweis über Lemmata aus der Fourier'schen Analysis. Schlömilch selbst verweist auch noch auf seine Arbeit \textit{``Über das Integral $\int \limits_{0}^{\infty} \dfrac{x^{\mu}dx}{r^2+2rx\cos u +x^2}$"} (\cite{Sc49}, 1849), welche das Integral

 \begin{equation*}
     \int\limits_{0}^{\infty} \dfrac{y^{\mu}dy}{1+2y\cos u+y^2},
 \end{equation*}
 für welches die für $-1 \leq \mu \leq 1$ gültige Gleichheit

 \begin{equation*}
      \int\limits_{0}^{\infty} \dfrac{y^{\mu}dy}{1+2y\cos u+y^2}=  \int\limits_{0}^{\infty} \dfrac{y^{-\mu}dy}{1+2y\cos u+y^2}
 \end{equation*}
 abgeleitet wird, untersucht. Die letzte Gleichung besitzt, wie in der Arbeit \textit{``Euler and a Proof of the Functional Equation for the Riemann Zeta-Function He Could Have Given"} (\cite{Ay24c}, 2024) genauer dargestellt, eine Schlüsselrolle im Beweis der Funktionalgleichung für $\eta(s)$ (\ref{eq: Functional Equation Eta}), sofern sie um ein entsprechendes Grenzwertargument ergänzt wird, welches Schlömilch indes nicht explizit angibt. Einen elementaren Beweis der Funktionalgleichung der $\eta$--Funktion gibt allerdings Hardy in seiner Arbeit \textit{``A new proof of the functional equation for the Zeta--function"} (\cite{Har22}, 1922). Der Beweis ist insofern elementar, als dass er lediglich die Formel

 \begin{equation*}
     (-1)^k \dfrac{\pi}{4}= \sum_{m=0}^{\infty} \dfrac{\sin (2m+1)x}{2m+1}
 \end{equation*}
 mit einer ganzen Zahl $k$ und $k\pi <x< (k+1)\pi$ voraussetzt, welche Euler in entsprechender Form ebenfalls bekannt war. Man findet etwa in § 16 von \cite{E555} die Formel:

 \begin{equation*}
     \dfrac{\pi-\omega}{2}= \sin \omega +\dfrac{\sin 2\omega}{2}+\dfrac{\sin 3\omega}{3}+\dfrac{\sin 4 \omega}{4}+\text{etc.},
 \end{equation*}
 aus welcher sich die von Hardy ohne Mühe ergibt.

\subsubsection{Durch einen neuen Gedanken: Die $\vartheta$-Funktion}
\label{subsubsec: Durch einen neuen Gedanken: Die Theta-Funktion}

\epigraph{The answers you get depend upon the questions you ask.}{Thomas Kuhn}

Rückblickend lässt sich bei einigen Begebenheiten konstatieren, dass Euler lediglich einen grundlegenden Gedankengangen von der Bergung eines neuen mathematischen Schatzes entfernt war. 
Um dies an der Jacobi'schen Thetafunktionen nachvollziehen zu können, ist allem voran die Euler'sche Methode Einführung neuer Funktionen in die Analysis zu erläutern.
Diesbezüglich findet man in seinen Werken zwei Wege; beide scheint er jedoch explizit nicht umfassender diskutiert zu haben, vielmehr erwähnt er sie im Verlauf einer anderen Untersuchungen beiläufig.

\paragraph{Einführung neuer Funktionen über bekannte Potenzreihen}
\label{para: Einführung neuer Funktionen über bekannte Potenzreihen}

Insbesondere seiner \textit{Introductio} \cite{E101} ist die Einführung der elementaren Funktionen, wie sie  heute noch geläufig sind\footnote{Dies bedeutet Polynome, gebrochen rationale Funktionen, Wurzelausdrücke derselben, Winkel-- und ihre Umkehrfunktionen, sowie Exponentialfunktion und Logarithmus.}, zu verdanken.  In seinem Lehrbuch geht er dabei von der Annahme aus, dass sich jede Funktion in der Form

\begin{equation*}
    f(x)= Ax^{\alpha}+Bx^{\beta}+Cx^{\gamma}+Dx^{\delta}+\cdots,
\end{equation*}
mit beliebigen Koeffizienten $A$, $B$, $C$, $D$ etc. und Potenzen $\alpha$, $\beta$, $\gamma$, $\delta$ etc., darstellen lässt. Dies setzt er in Kapitel 4 von \cite{E101} auseinander, wo  er  in § 59 schreibt:\\

\textit{``Es ist aber ersichtlich, dass keine nicht ganze Funktion von $z$ mit Termen der Art $A+Bz+Cz^2+\text{etc.}$ dargestellt werden kann; denn ansonsten wären sie ja ganz; indes kann sie aber durch eine unendliche Reihe dargestellt werden, wenn jemand daran zweifelt, wird dieser Zweifel durch die Entwicklung der entsprechenden Funktion zerstreut werden. Damit sich diese Ausführungen aber weiter erstrecken, müssen neben den ganzen positiven Potenzen von $z$ beliebige zugelassen werden. So wird freilich kein Zweifel bestehen, dass eine jede Funktion von $z$ in einen unendlichen Ausdruck von dieser Art überführt werden kann:}\\

\begin{equation*}
    Az^{\alpha}+Bz^{\beta}+Cz^{\gamma}+Dz^{\delta}+\text{etc.},
\end{equation*}
\textit{während die Exponenten $\alpha$, $\beta$, $\gamma$, $\delta$ etc. irgendwelche Zahlen bezeichnen."}\\

Nachdem die ganz rationalen und die rationalen Funktionen zuvor bereits abgehandelt worden sind, möchte Euler die irrationalen und elementaren transzendenten Funktionen einführen. Grundlage ist dabei  (§ 71) der binomische Lehrsatz. Die Wichtigkeit dieses Satzes für Euler kann kaum überschätzt werden. In der Arbeit  \textit{``Demonstratio theorematis Neutoniani de evolutione potestatum binomii pro casibus, quibus exponentes non sunt numeri integri"} (\cite{E465}, 1775, ges. 1773) (E465: ``Beweis des Newton'schen Theorems über die Entwicklung der Potenzen des Binoms für die Fälle, in denen die Exponenten keine ganzen Zahlen sind"), in welcher er ihn zum ersten Mal beweist\footnote{Tatsächlich sind auch moderne Beweise dem Euler'schen Beweis nachempfunden, siehe etwa das Lehrbuch von K. Knopp (1882--1957) \textit{``Theorie und Anwendung der unendlichen Reihen"} (\cite{Kn96}, 1996).}, schreibt er direkt im ersten Satz:\\

\textit{``Dieser Lehrsatz, für gewöhnlich so dargestellt}
\begin{equation*}
    (a+b)^n=a^n+\frac{n}{1}\cdot a^{n-1}b+\dfrac{n}{1}\cdot \dfrac{n-1}{2}a^{n-2}b^2+\dfrac{n}{1}\cdot \frac{n-1}{2} \cdot \dfrac{n-2}{3}a^{n-3}b^3+\text{etc.},
\end{equation*}
\textit{sofern verstanden wird, dass er sich weitest möglich erstreckt und in dem Exponenten $n$ alle denkbaren Zahlen erfasst werden, bildet das gesamte Fundament der höheren Analysis;[...]"}\\

Heute  wird er mithilfe der Binomialkoeffizienten wie folgt präsentiert:
\begin{equation}
    \label{eq: Binomischer Lehrsatz}
    (1+x)^{\alpha}=1+\binom{\alpha}{1}x+\binom{\alpha}{2}x^2+\binom{\alpha}{3}x^3+\cdots,
\end{equation}
in welcher Form Euler ihn auch ab § 72 verwendet.\\

\begin{figure}
    \centering
        \includegraphics[scale=1.15]{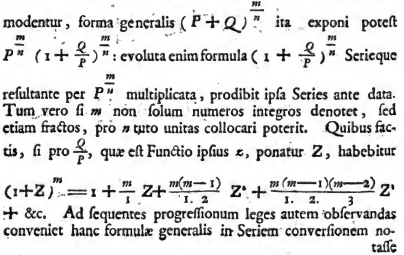}
    \caption{Euler stellt in seiner \textit{Introductio} \cite{E101} den binomischen Lehrsatz und die Binomalreihe  vor.}
    \label{fig:E101Binomi}
\end{figure}
 Damit sind automatisch auch alle irrationalen Funktionen als Reihen dargestellt bzw. lassen sich mit seiner Hilfe angeben. In Kapitel 5 wendet er sich nun auch noch den elementaren transzendenten Funktion $e^x$ und $\log(1+x)$ zu. Um erstere als Reihe darzustellen, bemerkt Euler

\begin{equation*}
    \left(1+\dfrac{x}{n}\right)^n =e^x
\end{equation*}
für eine unendlich große Zahl $n$. Für die linke Seite findet er vermöge (\ref{eq: Binomischer Lehrsatz})

\begin{equation*}
   \left(1+\dfrac{x}{n}\right)^n =1 + \dfrac{n}{1}\cdot \dfrac{x}{n}+\dfrac{n(n-1)}{1\cdot 2}\dfrac{x^2}{n^2}+\dfrac{n(n-1)(n-2)}{1\cdot 2 \cdot 3}\dfrac{x^3}{n^3}+\text{etc.} 
\end{equation*}
Da $n$ unendlich groß ist, so sagt Euler, ist

\begin{equation*}
    e^x=1+\dfrac{x}{1}+\dfrac{x^2}{1\cdot 2}+\dfrac{x^3}{1\cdot 2 \cdot 3}+\cdots,
\end{equation*}
was ja gerade die bekannte Potenzreihe von $e^x$ ist. Auch wenn aus moderner Sicht betrachtet  Eulers Herleitung wegen der unzulässigen Grenzwertvertauschung nicht vollständig ist, zeigt sie dennoch zeigt  seine Findigkeit, da er in seiner \textit{Introductio} \cite{E101} vollkommen auf die Differentialrechnung verzichten möchte und somit den klassischen Taylor'schen Satz nicht für den Beweis heranziehen kann. Die Potenzreihe für  $\log(1+x)$ leitet er ähnlich aus dem Grenzwert

\begin{equation*}
    \log(1+x)= \lim_{n \rightarrow \infty} n((1+x)^{\frac{1}{n}}-1)
\end{equation*}
her.  Hieraus zeigt sich, dass Euler neue Funktionen mit Vorliebe aus den \textit{bekannten} Reihenentwicklungen ableitet, was in der Folge zur Kenntnis neuer Funktionen immer neuere Methoden verlangt, immer kompliziertere Potenzreihen geschlossen aufsummieren zu können. Dies bietet weiter eine Erklärung, weshalb  das \textit{Finden} von Summen in den Euler'schen Arbeiten einen solch großen Raum einnimmt\footnote{Allein 4 Bände seiner \textit{Opera Omnia} (Band 14, 15, 16,1 und 16,2 der ersten Serie) enthalten einzig Arbeiten, die sich der Reihenlehre widmen.}.

\paragraph{Einführung über Integration von algebraischen Funktionen}
\label{para: Einführung über Integration von algebraischen Funktionen}

Einen weiteren Weg, neue Funktionen in die Analysis einzuführen, hat Euler in seiner Arbeit \textit{``De plurimis quantitatibus transcendentibus, quas nullo modo per formulas integrales exprimere licet"} (\cite{E565}, 1784, ges. 1775) (E565: ``Über die vielen transzendenten Größen, welche sich in keiner Weise durch Integralformeln ausdrücken lassen") präsentiert, wo sich die Methode gleich in § 1 zeigt: Neue transzendente Funktionen können über die Integration algebraischer Differentialformen eingeführt werden.  Als elementare Beispiele nennt Euler hier

\begin{equation*}
    \log x= \int\limits_{1}^x \dfrac{dx}{x}, \quad \arctan (x)= \int\limits_{0}^x \dfrac{dt}{1+tt}
\end{equation*}
und das nicht mehr triviale Beispiel

\begin{equation*}
    \int \sqrt{\dfrac{f+bxx}{g+kxx}}dx,
\end{equation*}
welches  zu elliptischen Integralen führt und -- wie Euler auch an vielen Stellen in seinem Opus bemerkt  -- sich nicht mit Logarithmen und Kreisbögen\footnote{Diese würden heute eher als inverse trigonometrische Funktionen oder Arkusfunktionen bezeichnet. Euler nennt sie zumeist Kreisbögen.}  ausdrücken lässt.\\

Die erwähnte Arbeit  zieht jedoch ihr Interesse auf sich, weil, wie auch der Titel bereits suggeriert, viele Reihen existieren, die sich nicht auf die beschriebene Weise gewinnen lassen, wovon in § 5  Euler als Beispiel die Reihe

\begin{equation}
\label{eq: Quadratsumme}
    1+x^1+x^4+x^9+x^{16}+\text{etc.}, 
\end{equation}
nennt, welche man als heute $\vartheta$--Funktion ausmacht, da

\begin{equation*}
    1+x^1+x^4+x^9+x^{16}+\text{etc.}=\dfrac{1}{2}\left(\vartheta_{3}(x)+1\right), \quad \text{da} \quad \vartheta_3(x)= \sum_{k=-\infty}^{\infty} x^{k^2}.
\end{equation*}
 Euler ist öfters zu der Reihe (\ref{eq: Quadratsumme}) geführt worden, da sie bzw. ihre Potenzen Aussagen zur Darstellbarkeit von Zahlen als Quadratsummen erlaubt, sofern man denn einen leicht potenzierbaren geschlossenen Ausdruck für sie finden könnte. Dies mag  sein Ausgangspunkt gewesen sein, überhaupt (\ref{eq: Quadratsumme}) zu erwähnen. Den Bezug zur Zahlentheorie macht er dabei explizit in § 6. \\

Aus den Erklärungen  mag man bereits erahnen, dass es Euler zeitlebens nicht gelingen sollte, einen geschlossenen Ausdruck für (\ref{eq: Quadratsumme}) zu finden und ähnliche Reihen, deren Potenzen gemäß einer quadratischen Progression fortschreiten, aufzusummieren. In § 10 von \cite{E475} schreibt er diesem Zusammenhang:\\

\textit{``Weil weiter unter Verwendung einer unendlichen Reihe gilt:}

\begin{equation*}
    \dfrac{\pi}{4}= \arctan \left(\dfrac{1}{2}\right)+\arctan \left(\dfrac{1}{8}\right)+ \arctan \left(\dfrac{1}{18}\right)+\arctan \left(\dfrac{1}{32}\right)
\end{equation*}
\begin{equation*}
    +\arctan \left(\dfrac{1}{50}\right)+\text{etc.},
\end{equation*}
\textit{der allgemeine Terme welcher Reihe $\arctan\left(\frac{1}{2nn}\right)$ ist, werden wir diese bemerkenswerte Integration haben:}

\begin{equation*}
    \int \dfrac{dx}{x\log x}(x^2+x^8+x^{18}+x^{32}+\text{etc.})\sin \log x = \dfrac{\pi}{4},
\end{equation*}
\textit{welche umso bemerkenswerter ist, weil die unendliche Reihe $x^2+x^8+x^{18}+x^{32}+\text{etc.}$ in keiner Weise auf eine endliche Summe reduziert werden kann.``}\footnote{Das Integral, so erklärt Euler vorher, ist von $x=0$ bis $x=1$ zu nehmen.} \\

\begin{figure}
    \centering
 \includegraphics[scale=1.0]{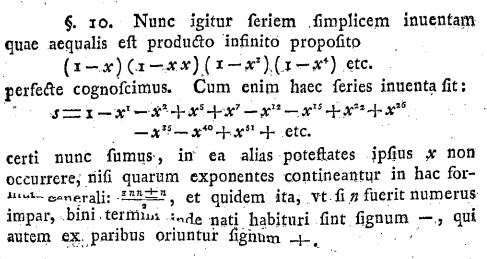}
    \caption{Euler spricht in seiner Arbeit \cite{E541} den Pentagonalzahlensatz aus.}
    \label{fig:E541Pentagonalzahlensatz}
\end{figure}

Jacobi, welcher Summen dieser Gestalt nachdrücklich in die Analysis eingebracht hat,  ist solche Probleme umgegangen, indem er solche Reihen als summiert \textit{festlegt}, wenn man sie auf elliptische Funktionen wie eben die just genannte $\vartheta_3$ zurückgeführt hat. Man vergleiche  seine  Arbeit \textit{``Suite des Notices sur les Functions Ellitiques"} (\cite{Ja28}, 1828) (``Weitere Bemerkungen über elliptische Funktionen").\\ 

Eine Summierung im eigentlichen Sinne gelingt Euler ebenfalls nicht bei seinem berühmten Pentagonalzahlensatz:
\begin{equation*}
    \prod_{k=1}^{\infty}(1-q^k) = \sum_{n=-\infty}^{\infty} (-1)^n q^{\frac{3n^2+n}{2}},
\end{equation*}
welchen er bereits zum ersten Mal in Kapitel 16, § 323 von \cite{E101} angibt und durch Tabellen zu erhärten versucht. Einen Beweis reicht er sehr viel später in der Arbeit \textit{``Evolutio producti infiniti $(1-x)(1-xx)(1-x^3)(1-x^4)(1-x^5)\cdot \text{etc.}$ in seriem simplicem"} (\cite{E541}, 1783, ges. 1775) (E541: ``Entwicklung des unendlichen Produktes $(1-x)(1-xx)(1-x^3)(1-x^4)(1-x^5)\cdot \text{etc.}$ in eine unendliche Reihe") nach\footnote{Der Euler'sche Beweis des Pentagonalzahlensatzes wird auch Polyas Buch \textit{``Mathematics and Plausible Reasoning"} (\cite{Po14}, 2014) ausführlich als Beispiel für Induktion in Eulers Werk besprochen.}. Euler findet demnach aus dem Produkt die Summe, nicht aus der Summe das Produkt. Hierfür bedurfte es erst der Idee von Abel \textit{``Recherches sur les fonctions elliptiques"} (\cite{Ab27}, 1827) (``Untersuchungen zu elliptischen Funktionen"), die elliptischen Funktionen als Umkehrfunktionen zu den elliptischen Integralen einzuführen und sie in  Analogie zu den aus der Trigonometrie bekannten Arkusfunktionen und Sinusfunktionen zu behandeln.  Der Abel'sche Ansatz liegt insofern diametral zu den Euler'schen Herangehensweisen, da er die Frage aufwirft, bis zu welcher Grenze eine Quadratur mit vorgegebenen Integranden erstreckt werden muss\footnote{In Abels erwähnter Arbeit \cite{Ab27} waren dies die elliptischen Integrale, ein Spezialfall der nach Abel benannten Abel'schen Integrale.}, sodass selbige Quadratur einem beliebigen Wert gleich wird, wohingegen Euler Quadraturen mit variablen Grenzen überhaupt nicht eingehender untersucht zu haben scheint. Sie treten dann als Lösungen von Differentialgleichungen auf, werden aber dann nicht weiter diskutiert. Selbst an den Beispielen der $\Gamma$-Funktion etwa in \cite{E19} oder der hypergeometrischen Funktion in \cite{E274} findet sich bei Euler die variable Größe stets im Integranden und nicht in der Integralgrenze selbst.  

\subsection{Von Euler nicht beweisbare Entdeckungen}
\label{subsec: Von Euler Entdecktes, aber nicht Beweisbares}

\epigraph{Today's accomplishments were yesterday's impossibilities.}{Robert H. Schuller}

Nach der Besprechung von Instanzen, in welchen Euler die Fähigkeit hatte, einen Beweis für seine Entdeckung zu liefern, es jedoch nicht tat, werden nun von Euler nicht beweisbare Entdeckungen in den Fokus rücken. Die Gründe für das Ausbleiben eines Nachweises reichen dabei von einer unzureichenden Formulierung wie beim Primzahlsatz (Abschnitt \ref{subsubsec: Weg fehlender Formulierung: Der Primzahlsatz}) über das schlichte Fehlen von Mitteln wie beim quadratischen Reziprozitätsgesetz (Abschnitt \ref{subsubsec: Wegen fehlender Mittel: Das Reziprozitätsgesetz})  bis hin zu einer für Euler nicht erkennbaren Unvollständigkeit des Beweis wie beim großen Satz von Fermat für den Fall $n=3$ (Abschnitt \ref{subsubsec: Wegen übersehener Unvollständigkeit: Der große Satz von Fermat für n=3}).

\subsubsection{Wegen fehlender Formulierung: Der Primzahlsatz}
\label{subsubsec: Weg fehlender Formulierung: Der Primzahlsatz}


\epigraph{God may not play dice with the universe, but something strange is going on with the prime numbers.}{Carl Pomerance}

Der Primzahlsatz ist eines der großen Resultate der Zahlentheorie, sodass es verwundern mag, dass Euler nicht zumindest auf ihn gestoßen sein soll\footnote{In der Literatur wird oft Gauß als der Erstentdecker genannt, der ihn basierend auf Primzahltabellen vermutet hat. Es existiert ein Brief aus seinem Nachlass, in dem Gauß die Entdeckung des Primzahlsatzes auf das Jahr 1796 datiert; man vergleiche \cite{Ga49}.}. Denn, wie wenig verwundern dürfte,  die Primzahlen auch Eulers Interesse auf sich, obwohl sie ihn bisweilen an seine Grenzen gebracht haben. Dies bringt er eines Beispiels wegen direkt in § 1 seiner Arbeit \textit{``Découverte d'une loi tout extraordinaire des nombres par rapport à la somme de leurs diviseurs"} (\cite{E175}, 1751, ges. 1747) (E175: ``Entdeckung eines völlig außergewöhnlichen Gesetzes bezüglich der Summe ihrer Teiler") folgendermaßen zum Ausdruck:\\

\textit{``Die Mathematiker haben sich bis zum heutigen Tage vergebens bemüht, irgendeine Ordnung in der Folge der Primzahlen zu entdecken, und man ist geneigt zu glauben, dies sei ein Geheimnis, das der menschliche Geist niemals zu durchdringen wissen wird."}\\

Jedoch liegt in den Euler'schen Arbeiten ein heuristisches Argument vergraben, mit welchem man zumindest zur Formulierung des Primzahlsatzes gelangen kann. Grundlage ist dabei das Euler--Produkt

\begin{equation}
    \label{eq: Euler Produkt}
    \sum_{n=1}^{\infty} \dfrac{1}{n^s} = \prod_{p \in \mathbb{P}} \dfrac{1}{1-\frac{1}{p^s}},
\end{equation}
wobei $\mathbb{P}$ die Menge der Primzahlen bedeutet, welches zu ersten Mal von Euler in Theorem 8 in der Arbeit \textit{``Variae observationes circa Series infinitas"} (\cite{E72}, 1744, ges. 1737) (E72: ``Verschiedene Beobachtungen über unendliche Reihen")  abgeleitet wird. Eine geordnete Darstellung gibt Euler  in § 173 seiner \textit{Introductio} \cite{E101}.\\

\begin{figure}
    \centering
  \includegraphics[scale=1.0]{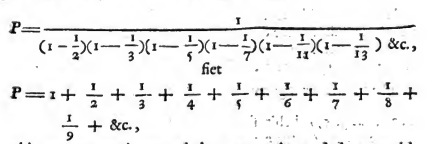}
    \caption{Euler formuliert das Euler--Produkt für die harmonische Reihe. Hier ist die Version von \cite{E101} zu sehen.}
    \label{fig:E101EulerProd}
\end{figure}

Hier soll aber die Aufmerksamkeit auf Theorem 7 und die Korollare aus selbigem in \cite{E72} gerichtet werden. In Korollar 3 zu diesem Theorem 7 schreibt Euler:\\

\textit{``Aber daraus sieht man auch ein, dass es unendlich mal weniger Primzahlen als ganze Zahlen gibt, denn dieser Ausdruck $\frac{2\cdot 3 \cdot 4 \cdot 5 \cdot 6 \cdot 7 \cdot \text{etc.}}{1\cdot 2\cdot 3 \cdot 4 \cdot 5 \cdot 6  \cdot \text{etc.}}$ hat einen absolut unendlichen Wert, während derselbe, aber lediglich aus Primzahlen entstanden, der Logarithmus von jenem Unendlichen ist."}\\

Modern mag man Eulers Aussage vielleicht so ausdrücken:
\begin{equation*}
   \log\left(\prod_{k=1}^{\infty} \dfrac{k+1}{k}\right)=  \prod_{p \in \mathbb{P}} \dfrac{p}{p-1}=\prod_{p \in \mathbb{P}} \dfrac{1}{1-\frac{1}{p}},
\end{equation*}
wobei  das $=$ natürlich formal zu verstehen ist. Dass  im Ausdruck rechter Hand die Primzahlen im Zähler und im Nenner dieselben weniger $1$ zu schreiben sind, lehrt Euler  im Beweis zu seinem Theorem 7. Anders ausgedrückt formuliert Euler  hier, dass die harmonische Reihe wie der Logarithmus wächst:

\begin{equation*}
    \prod_{p \in \mathbb{P}} \dfrac{1}{1-\frac{1}{p}}= \sum_{n=1}^{\infty} \dfrac{1}{n}= \log(\infty),
\end{equation*}
wie Euler  auch  in Korollar 1 zu selbigem Theorem 7 sagt\footnote{Mertens (1840--1927) hat später in seiner Arbeit \textit{``Ein Beitrag zur analytischen Zahlentheorie"} (\cite{Me74}, 1874) folgende Formel nachgewiesen:

\begin{equation*}
    \lim_{n\rightarrow \infty} \dfrac{1}{\log(n)}\prod_{p \in \mathbb{P}}^{n} \dfrac{1}{1-\frac{1}{p}}=e^{\gamma},
\end{equation*}
wobei  $\gamma$ die Euler--Mascheroni--Konstante bedeutet.}. Damit ist unter anderem die Unendlichkeit der Menge der Primzahlen nachgewiesen, wie auch Sandifer in seiner Notiz \textit{``Infinitely many primes"} (\cite{Sa06mar}, 2006) hervorhebt. Aus heutiger Sicht mag Eulers Ausdrucksweise mehr zur Verwirrung als zur Erklärung beitragen und hat womöglich auch deswegen Sandifer in seinem Buch \textit{``The Early Mathematics of Leonhard Euler"} (\cite{Sa07a}, 2007) in der Besprechung von Eulers Arbeit \cite{E72} bewogen, die oben zitierte Aussage wie folgt in Formeln ausgedrückt zu sehen

\begin{equation*}
    \lim_{n\rightarrow \infty} \dfrac{\pi(n)}{\log(n)}=1,
\end{equation*}
mit $\pi(n)$ als der Anzahl der Primzahlen unter $n$. Sandifer merkt aber auch an, dass er sich bei dieser Interpretation wegen Eulers Formulierungsweise nicht völlig sicher ist.  Trotz alledem ist es   möglich, aus den Euler'schen Überlegungen heraus eine Verbindung zu Primzahlsatz herzustellen. \\

Auf der formalen Ebene sagt Euler nämlich auch

\begin{equation*}
    \log \zeta(1)= \log (\log(\infty)),
\end{equation*}
wobei $\zeta(1)$ hier der Kürze wegen für die harmonische Reihe steht, sodass

\begin{equation*}
    \zeta(1)= 1+\dfrac{1}{2}+\dfrac{1}{3}+\dfrac{1}{4}+\cdots
\end{equation*}
Nach dem Euler--Produkt (\ref{eq: Euler Produkt}) gilt jedoch überdies:
\begin{equation*}
    \log \zeta(1)=- \sum_{p\in \mathbb{P}} \log \left(1-\dfrac{1}{p}\right)=-\sum_{n=2}^{\infty}(\pi(n)-\pi(n-1))\log \left(1-\dfrac{1}{n}\right).
\end{equation*}
Ein Indexshift gibt:

\begin{equation*}
    \log \zeta(1)= -\sum_{n=2}^{\infty}\pi(n) \left(\log\left(1-\dfrac{1}{n}\right)-\log\left(1-\dfrac{1}{n+1}\right)\right).
\end{equation*}
Weil 

\begin{equation*}
    \int \dfrac{dx}{x(x-1)}=\log(1-x)-\log(x)
\end{equation*}
gilt, hat man folglich

\begin{equation*}
    \log \zeta(1)= \sum_{n=2}^{\infty} \pi(n)\int\limits_{n}^{n+1} \dfrac{dx}{x(x-1)}=\int\limits_{2}^{\infty}\dfrac{\pi(x)dx}{x(x-1)}.
\end{equation*}
Euler hätte  dies  wohl zunächst so geschrieben:

\begin{equation*}
    \log(\log(\infty))=\int\limits_{2}^{\infty}\dfrac{\pi(x)dx}{x(x-1)}.
\end{equation*}
 Die letzte Gleichung\footnote{Dies ist ein spezieller Fall von Formel (1.1.3) aus dem Buch \textit{``The Theory Of The Riemann Zeta-Function"}  (\cite{Ti87}, 1987) von Titchmarsh (1899--1963)):  

\begin{equation*}
    \log(\zeta(s))= \int\limits_{2}^{\infty} \dfrac{s\pi(x)dx}{x(x^s-1)},
\end{equation*} 
wobei $\zeta(s)$ die Riemann'sche $\zeta$--Funktion (\ref{eq: Def zeta}) bezeichnet.} drückt aus, dass das Integral auf der rechten Seite divergiert wie $\log(\log(x))$ für $x \rightarrow \infty$, was die inverse Frage impliziert, wie  die Funktion $\pi(x)$ zu wählen ist, dass sich das Integral zu $\log(\log(x))$ integriert. Da

\begin{equation*}
    \dfrac{d}{dx}\log(\log(x))= \dfrac{1}{\log(x)}\cdot \dfrac{1}{x}
\end{equation*}
gilt, müsste also 

\begin{equation*}
    \dfrac{1}{x\cdot \log(x)}= \dfrac{\pi(x)}{x(x-1)}
\end{equation*}
gelten. Dies ergibt 

\begin{equation*}
    \pi(x)= \dfrac{x-1}{\log(x)}.
\end{equation*}
Für unendliches $x$ demnach insbesondere:

\begin{equation*}
    \pi(x) \sim \dfrac{x}{\log(x)} \quad \text{für} \quad x \rightarrow \infty,
\end{equation*}
was gerade der Primzahlsatz ist, wie in Gauß in seinen Tagebüchern zuerst formuliert hat\footnote{Man vergleiche dazu den Brief aus \cite{Ga49} aus seinem Nachlass.}. Das vorgestellte Argument ist selbstredend rein heuristisch. Es zeigt allerdings, wie Euler zumindest auch aus seinen Überlegungen zum Primzahlsatz hätte gelangen können, sich wegen umständlicher Formulierungen diesen Weg jedoch  versperrte.

\subsubsection{Wegen fehlender Mittel: Das Reziprozitätsgesetz}
\label{subsubsec: Wegen fehlender Mittel: Das Reziprozitätsgesetz}

\epigraph{The Difficult is that which can be done immediately; the Impossible that which takes a little longer.}{George Santayana}

Entsprechende algebraische Konzepte werden bereits schmerzlich bei kleineren Sätzen der Zahlentheorie vermisst, wie etwa dem Satz von Wilson (1741--1793), also den Satz

\begin{equation*}
    p| (p-1)! +1 \quad \Longleftrightarrow \quad p ~ \in \mathbb{P}.
\end{equation*}
Während Euler diesen im Vergleich zu Lagrange\footnote{Lagrange liefert den ersten Beweis dieses Resultats in seiner Arbeit \textit{``Démonstration d'un théorème nouveuau concernant les nombres premiers"} (\cite{La73}, 1773, ges. 1771) (``Beweis eines neuen Lehrsatzes  über Primzahlen").} noch elegant zu beweisen vermag -- Eulers Beweis aus seiner Arbeit \textit{``Miscellanea analytica"} (\cite{E560}, 1783, ges. 1773) (E560: ``Verschiedenes Analytisches") ist der heute übliche über die Primitivwurzel, welches Konzept er freilich nicht in der heutigen Form kennt, man aber aus seinen Erläuterungen leicht herauslesen kann -- sind spätestens beim quadratischen Reziprozitätsgesetz Lehrsätze aus der höheren Algebra für einen ersten Beweis wohl unumgänglich. Freilich existieren heute analytische Beweise, sie setzen aber stets bei der Auswertung einer Gauß'schen Summe an, welche Euler ebenfalls unbekannt war. Wieder andere Beweise ziehen die Jacobi'schen $\vartheta$--Funktionen heran, die Euler wie oben (Abschnitt \ref{subsubsec: Durch einen neuen Gedanken: Die Theta-Funktion}) diskutiert nur oberflächlich berührt. Am ehesten, was die analytischen Mittel betrifft, wäre Euler wohl zum dem  Eisenstein'schen Beweis  des Reziprozitätsgesetzes aus der Arbeit \textit{``Applications de l'Algébre à l'Arithmètique transcendente"} (\cite{Ei45}, 1845) (``Anwendungen der Algebra auf die transzendente Zahlentheorie") fähig gewesen. Selbiger nutzt neben den dem bereits erwähnten Gauß'schen Lemma überdies das Legendre--Symbol\footnote{Das Legendre--Symbol $\left(\frac{a}{p}\right)$ für eine ganze Zahl $a$ und eine Primzahl $p$ kann bekanntermaßen nur die drei Werte $-1,0$ und $1$ annehmen. Es ist $1$, sofern $a$ ein quadratischer Rest modulo $p$ ist, $-1$, falls $a$ kein quadratischer Rest modulo $p$ ist und $0$, dass $p$ bereits $a$ teilt.}. Es gilt nämlich zunächst:

\begin{Thm}[Gauß'sches Lemma]
    Es gilt

    \begin{equation*}
        \left(\dfrac{a}{p}\right)= (-1)^{\sum_{j=1}^{\frac{p-1}{2}}\left[ \frac{2aj}{p} \right]}
    \end{equation*}
    mit dem Legendre--Symbol $\left(\frac{a}{p}\right)$.
\end{Thm}
Eisensteins Beweis nutzt die Identität:

\begin{equation*}
    \dfrac{\sin (mx)}{\sin (x)}=  (-4)^{\frac{m-1}{2}}\prod_{j=1}^{\frac{m-1}{2}} \left(\sin^2 (x)-\sin^2 \left(\dfrac{2\pi j}{m}\right)\right)
\end{equation*}
für eine ungerade Zahl $m$. Nennt man nun die Menge der Reste modulo einer ungeraden Primzahl $p$, die $\leq \frac{p-1}{2}$ sind, der Kürze wegen $R_p$, so gilt mit dem Lemma von Gauß und der eben genannten Identität:

\begin{equation*}
   \left(\frac{q}{p}\right)= \prod_{j \in R_p} \dfrac{\sin \left(\frac{2\pi q j}{p}\right)}{\sin \left(\frac{2 \pi j}{p}\right)}= \prod_{j \in R_p} (-4)^{\frac{q-1}{2}} \prod_{k \in R_q} \left(\sin^2 \left(\dfrac{2\pi j}{q}\right)-\sin^2 \left(\dfrac{2 \pi k}{p}\right)\right)
\end{equation*}
oder unter Vertauschung der Summanden in der letzten Formel:

\begin{equation}
\label{eq: Eisenstein--Schlüssel}
    \left(\frac{q}{p}\right) = (-4)^{\frac{p-1}{2}\frac{q-1}{2}} \prod_{j \in R_p} \prod_{j \in R_q} \left(\sin^2 \left(\dfrac{2\pi k}{p}\right)-\sin^2 \left(\dfrac{2 \pi j}{q}\right)\right)
\end{equation}
Also gilt analog für vertauschte Rollen von $p$ und $q$:

\begin{equation*}
     \left(\frac{p}{q}\right)=  (-4)^{\frac{p-1}{2}\frac{q-1}{2}} \prod_{j \in R_p} \prod_{j \in R_q} \left(\sin^2 \left(\dfrac{2\pi k}{q}\right)-\sin^2 \left(\dfrac{2 \pi j}{p}\right)\right).
\end{equation*}
Da das Produkt über $R_p$ nun $\frac{p-1}{2}$ und  das über $R_q$ entsprechend $\frac{q-1}{2}$ Faktoren beinhaltet, gilt auch:

\begin{equation*}
     \left(\frac{p}{q}\right)= (-4)^{\frac{p-1}{2}\frac{q-1}{2}} (-1)^{\frac{p-1}{2}\frac{q-1}{2}}\prod_{j \in R_p} \prod_{j \in R_q} \left(\sin^2 \left(\dfrac{2\pi k}{p}\right)-\sin^2 \left(\dfrac{2 \pi j}{q}\right)\right),
\end{equation*}
was jedoch nach (\ref{eq: Eisenstein--Schlüssel}) gerade wieder mit $\left(\frac{q}{p}\right)$ ausgedrückt werden kann, sodass insgesamt:

\begin{equation*}
     \left(\frac{p}{q}\right) = (-1)^{\frac{p-1}{2}\frac{q-1}{2}} \left(\frac{q}{p}\right).
\end{equation*}
Dies ist genau die Formulierung des quadratischen Reziprozitätsgesetzes, wie sie Legendre in seinem Lehrbuch \textit{``Essai sur la théorie des nombres"} (\cite{Le97}, 1797) (``Essay zur Zahlentheorie") gegeben hat. \\

Euler selbst formuliert das Reziprozitätsgesetz  in seiner Arbeit \textit{``Observationes circa divisionem quadratorum per numeros primos"} (\cite{E552}, 1783, ges. 1772) (E552: ``Beobachtungen zur Teilung von Quadratzahlen durch Primzahlen") wie folgt, sofern man seine Aussage gedrungen in moderner Terminologie darstellt:\\

\begin{figure}
    \centering
    \includegraphics[scale=0.8]{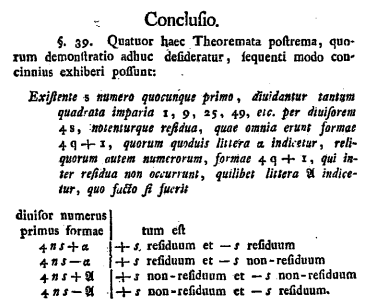}
    \caption{Euler formuliert im letzten Paragraphen seiner Arbeit \cite{E552} seine Version des quadratischen Reziprozitätsgesetzes und räumt ein, dass er es nicht beweisen kann.}
    \label{fig:E552Rez}
\end{figure}

\textit{``Für eine ungerade Primzahl $p$ mit $(a,q)=1$, mit einer ganzen Zahl $a$, gelte $p \equiv \pm q ~\operatorname{mod~} (4a)$ für eine weitere Primzahl $q$, dann gilt $\left(\frac{a}{q}\right)= \left(\frac{a}{p}\right)$, wobei hier die Legendre--Symbole gemeint sind."} \\

Zum Vergleich und zur Untermalung des durch die Legendre'sche Formulierung gewonnenen Vorteils sei auch noch Eulers Beschreibung aus besagter Arbeit beigefügt. Er schreibt: \\

\textit{``Während $s$ eine beliebige Primzahl ist, teile man die ungeraden Quadrate $1$, $9$, $25$ etc. durch den Teiler $4s$, dann werden Reste entstehen, die alle von der Form $4q+1$ sein werden, von welchen man jeden mit dem Buchstaben $\alpha$ notiere, aber von den übrigen Zahlen der Form $4q+1$, welche nicht unter den Resten auftreten, zeige man jedwede mit dem Buchstaben $\mathfrak{A}$ an, wonach es sich wie nachstehend verhalten wird:}

\begin{equation*}
    \renewcommand{\arraystretch}{1,5}
 \setlength{\arraycolsep}{3.5mm}
\begin{array}{l|l}
   \textit{Für Primteiler der Form}  & \textit{ist dann} \\ 
   \qquad \qquad 4ns+ \alpha & \textit{$+s$ ein Rest und $-s$ ein Rest} \\ 
    \qquad \qquad 4ns - \alpha & \textit{$+s$ ein Rest und $-s$ ein Nicht--Rest} \\ 
   \qquad \qquad 4ns+ \mathfrak{A} & \textit{$+s$ ein Nicht--Rest und $-s$ ein Nicht--Rest} \\
    \qquad  \qquad 4ns- \mathfrak{A} & \textit{$+s$ ein Nicht--Rest und $-s$ ein Rest}"  
\end{array}
\end{equation*}

Im Falle des quadratischen Reziprozitätsgesetzes ist demnach eine echte Grenze Eulers erreicht, angesichts des fehlenden algebraischen Rüstzeugs sowie einer entsprechenden Formulierung lässt sich zur Behauptung gelangen, dass sich ein Beweis des Reziprozitätsgesetzes aus dem Euler'schen Opus heraus nicht geben lässt, sondern eines darüber hinausgehenden Konzepts bedarf. \\

\begin{figure}
    \centering
  \includegraphics[scale=1.2]{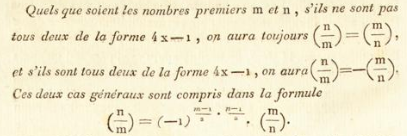}
    \caption{Legendre formuliert in seinem Buch \cite{Le97} das quadratische Reziprozitätsgesetz.}
    \label{fig:LegendreReci}
\end{figure}
All dem zum Trotze sei an dieser Stelle noch ein Gedankengang mitgeteilt, wie sich Euler, sofern er von der Bedeutung des Legendre--Symbols gewusst hätte, sich von der Gültigkeit des Reziprozitätsgesetzes aus seinen Überlegungen heraus überzeugt haben könnte. Bereits in seiner Arbeit \textit{``Theoremata circa divisores numerorum in hac forma $paa \pm qbb$ contentorum"} (\cite{E164}, 1751, ges. 1747) (E164: ``Theoreme über die in der Form $paa \pm qbb$ enthaltenen Zahlen") findet sich, wenn auch verbal formuliert, in Annotatio 16 eine (von Euler nicht bewiesene) Aussage, welche sich wie folgt in moderne Terminologie überführen lässt:\\

\begin{Thm}[Eulers erste Version des Reziprozitätsgesetzes]
\label{Theorem: Euler Rezi}
    Für ungerade Primzahlen $p$ und $q$ hat man folgende Äquivalenz
    \begin{equation*}
        \left(\dfrac{q}{p}\right)=1 \quad \Longleftrightarrow \quad p \equiv \pm \beta^2 ~ \operatorname{mod} (4q),
    \end{equation*}
    wobei $\left(\frac{q}{p}\right)$ das Legendre--Symbol meint und $\beta$ eine ungerade Zahl ist.
\end{Thm}
Nun wäre für Euler wohl naheliegend gewesen, den nachstehenden Ausdruck einer Untersuchung zu unterwerfen:

\begin{equation}
    \label{eq: f(p,q)}
    f(p,q):=  \left(\dfrac{q}{p}\right) \left(\dfrac{p}{q}\right)
\end{equation}
und diesen zunächst weitestmöglich zu vereinfachen. Aus Eulers Vermutung gilt $p= \pm\beta +4ql $ mit einer ganzen Zahl $l$. Damit hat man

\begin{equation*}
     f(p,q)=  \left(\dfrac{q}{p}\right) \left(\dfrac{\pm\beta^2 +4ql}{q}\right).
\end{equation*}
Aus den Eigenschaften des Legendre--Symbols folgt weiter

\begin{equation*}
      f(p,q)= \left(\dfrac{q}{p}\right) \left(\dfrac{\pm\beta^2}{q}\right)= \left(\dfrac{q}{p}\right) \left(\dfrac{\beta^2}{q}\right)\left(\dfrac{\pm 1}{q}\right)=  \left(\dfrac{q}{p}\right) \left(\dfrac{\pm 1}{q}\right).
\end{equation*}
Nun sind Fälle zu unterscheiden, je nachdem, ob das obere oder das untere Vorzeichen gilt. Es gelte zunächst das obere. Dann findet man folgendes vor:

\begin{equation*}
    f(p,q) = \left(\dfrac{q}{p}\right) \left(\dfrac{1}{q}\right)= \left(\dfrac{q}{p}\right)= 1,
\end{equation*}
wobei im letzten Schritt die Euler'sche Annahme, dass $\left(\frac{q}{p}\right)=1$ ist, aus dem Theorem genutzt wurde. \\
Jetzt gelte das untere Vorzeichen, woraus sich ergibt:

\begin{equation*}
    f(p,q) =  \left(\dfrac{q}{p}\right) \left(\dfrac{-1}{q}\right)=  \left(\dfrac{q}{p}\right) (-1)^{\frac{q-1}{2}}=(-1)^{\frac{q-1}{2}}.
\end{equation*}
Im letzten Schritt wurde wieder die Voraussetzung $\left(\frac{q}{p}\right)=1$ verwendet. Dass $\left(\frac{-1}{q}\right)=(-1)^{\frac{q-1}{2}}$ gilt, was gerade der erste Ergänzungssatz zum quadratischen Reziprozitätsgesetz ist, hat Euler gleichsam durch den Nachweis des Zwei--Quadratesatzes aus Theorem (\ref{Theorem: Zwei--Quadrate--Satz}) gezeigt\footnote{Dieser wird unten in Abschnitt (\ref{subsubsec: Praxis über Abstraktion}) noch ausführlich diskutiert werden.}. \\

Weiter lässt sich $f(p,q)$ zunächst einmal nicht einschränken. Man sieht jedoch leicht ein, dass die Auswahl

\begin{equation*}
    f(p,q)= (-1)^{\frac{p-1}{2}\frac{q-1}{2}}
\end{equation*}
den herausgearbeiteten Bedingungen Genüge leistet, zumal man zunächst einmal findet:

\begin{equation*}
     (-1)^{\frac{p-1}{2}\frac{q-1}{2}}= (-1)^{\frac{\pm \beta^2+4ql -1}{2}\frac{q-1}{2}},
\end{equation*}
weshalb wieder nach den Vorzeichen unterschieden werden muss. Es gelte das obere und es sei $\beta=2\alpha +1$, da $\beta$ als ungerade Zahl vorausgesetzt wird, sodass

\begin{equation*}
    (-1)^{\frac{\beta^2+4ql -1}{2}\frac{q-1}{2}}= (-1)^{\frac{ 4\alpha^2+4\alpha+1+4ql -1}{2}\frac{q-1}{2}}=  (-1)^{(2\alpha^2+2\alpha+2ql)\frac{q-1}{2}} =1,
\end{equation*}
da $(2\alpha^2+2\alpha+2ql)$ eine gerade Zahl ist. \\
Für den anderen Fall gilt, sofern wieder $\beta =2\alpha +1$ gesetzt wird, 

\begin{equation*}
        (-1)^{\frac{+ \beta^2+4ql -1}{2}\frac{q-1}{2}}  =  (-1)^{ (-2\alpha^2-2\alpha +2ql -1)\frac{q-1}{2}}=(-1)^{\frac{q-1}{2}},
\end{equation*}
da $-2\alpha^2-2\alpha +2ql -1$ eine ungerade Zahl ist, womit auch diese aus den  Euler'schen Annahme herstammende Bedingung erfüllt wird, sodass Euler sich zweifelsohne bewogen sähe, zu behaupten,  (\ref{eq: f(p,q)}) stelle sich wie folgt dar:

\begin{equation*}
    f(p,q):=  \left(\dfrac{q}{p}\right) \left(\dfrac{p}{q}\right)= (-1)^{\frac{p-1}{2}\frac{q-1}{2}},
\end{equation*}
was natürlich das Reziprozitätsgesetz in Legendre'scher Formulierung ist. Dass neben der Wahl von $f(p,q)$ als  $(-1)^{\frac{p-1}{2}\frac{q-1}{2}}$ noch weitere mit seinen Annahmen verträgliche existieren, wäre Euler vermutlich auch bewusst gewesen. Jedoch unterscheiden sie sich lediglich im Exponenten  um ein gerades Vielfaches vom Ausdruck $\frac{p-1}{2}\frac{q-1}{2}$, sodass Euler hier mit guten Recht sein Argument der Einfachheit als Maßstab für die Richtigkeit (wie oben im Abschnitt (\ref{subsubsec: Eulers Auffassung eines Beweises}) am Beispiel der Funktionalgleichung $\zeta$--Funktion gesehen) seiner Auswahl geltend gemacht hätte. Überdies lässt sich mit den vorgeführten ähnlichen Rechnungen zeigen, dass das Reziprozitätsgesetz auch Eulers Vermutung aus Theorem (\ref{Theorem: Euler Rezi}) impliziert.

\subsubsection{Wegen nicht zu sehender Unvollständigkeit: Der große Satz von Fermat für n=3}
\label{subsubsec: Wegen übersehener Unvollständigkeit: Der große Satz von Fermat für n=3}

\epigraph{I have discovered a truly marvelous proof of this, which however the margin is not large enough to contain.}{Pierre de Fermat}

Wie bereits an verschiedenen Stellen angeklungen,  sollte das Euler'sche Interesse an zahlentheoretischen Fragen, insbesondere an den Entdeckungen von Fermat,  nie zum Erliegen kommen. Einen ausführlichen Bericht mit mathematischen Details finden man im Buch von  Weil   \textit{``Number theory. An approach through history. From Hammurapi to Legendre"} (\cite{We84}, 1984) im Kapitel über Euler. Hier soll es um einen Spezialfall des großen Fermat'schen Satzes gehen, welchen Euler in § 243 des zweiten Teils seines Buchs \textit{``Vollständige Anleitung zur Algebra"}, bestehend aus den zwei Teilen  (\cite{E387}, 1770, ges. 1767) und (\cite{E388}, 1770, ges. 1776), bewiesen hat. Er intendiert zu demonstrieren, dass

\begin{equation*}
    x^3+y^3=z^3
\end{equation*}
abgesehen von der offenkundigen Lösung $x=y=z=0$ keine andere ganzzahlige Lösung zulässt.  Der Schlüssel zum gesuchten Nachweis besteht in folgender Erkenntnis: Ist $ax^2+cy^2$, wobei $-ac$ kein Quadrat sein soll\footnote{Euler drückt dies aus, indem er fordert, dass die Form $ax^2+by^2$ sich nicht in lineare Faktoren zerlegen lassen soll.}, und überdies ein Kubus, wo $x$ und $y$ teilerfremd sein sollen, so kann man stets ganze Zahlen $p$ und $q$ finden, sodass

\begin{equation*}
    x=ap^3-3cpq^2 \quad \text{und}\quad y=3ap^2q-cq^3
\end{equation*}
gilt. Dass diese Wahl von $x$ und $y$ die Form $ax^2+by^2$ zu einem Kubus macht, ist durch direkte Rechnung schnell überprüft.\\

\begin{figure}
    \centering
    \includegraphics[width=0.8\linewidth]{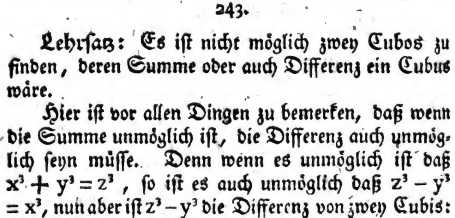}
    \caption{Euler formuliert in seiner \textit{Algebra} \cite{E387} den großen Satz von Fermat für $n=3$.}
    \label{fig:SatzFermat3}
\end{figure}

Der Euler'sche Gedankengang lässt sich aus dem Ansatz

\begin{equation*}
    x\sqrt{a}\pm y\sqrt{-c}=(p\sqrt{a}\pm q\sqrt{-c})^3
\end{equation*}
durch Ausmultiplizieren der rechten Seite und anschließendem Koeffizientenvergleich unmittelbar nachvollziehen. Konkret braucht Euler hier den Fall, in welchem $a=1$ und $c=3$ ist, sodass man im Erweiterungskörper $\mathbb{Q}[\sqrt{-3}]$ arbeitet, die Eigenschaften welcher Zahlen Euler in seiner Abhandlung \textit{``Solutio generalis quorundam problematum Diophanteorum, quae vulgo nonnisi solutiones speciales admittere videntur"} (\cite{E255}, 1761, ges. 1754) (E255: ``Die allgemeine Lösung gewisser Diophant'scher Probleme, die für gewöhnlich nur spezielle Lösungen zuzulassen scheinen") herausgearbeitet hat. Eulers Beweis ist auf einem unendlichen Abstieg begründet und die Methode lässt sich auch als mathematisch valide nachweisen, da sich der Euklid'sche Algorithmus  $\mathbb{Q}[\sqrt{-3}]$ durchführen lässt. Das Problem des Euler'schen Arguments ist demnach folgendes: Freilich lässt sich, wie Euler es in § 243 seiner \textit{Algebra} \cite{E387} tut,

\begin{equation*}
    p+q\sqrt{-3}= (t+u\sqrt{-3})^3 \quad \text{bzw.} \quad p-q\sqrt{-3}= (t-u\sqrt{-3})^3
\end{equation*}
und damit auch

\begin{equation*}
    p^2+3q^2=(t^2+3u^2)^3
\end{equation*}
setzen, jedoch gilt die Umkehrung, aus dem Kubus auf die linke Seite zu schließen, nur wegen der Euklid'schen Struktur von  $\mathbb{Q}(\sqrt{-3})$. Euler nimmt es allerdings ohne Beweis an, dass $p^2+3q^2$ \textit{nur} die Zerlegung $(p+\sqrt{-3}q)(p-\sqrt{-3}q)$ zulässt. Er hat hier freilich Recht, jedoch ist dies für andere Ringe ähnlicher Gestalt falsch\footnote{Die Wichtigkeit einer bis auf Einheiten und Reihenfolge eindeutigen Primzahlzerlegung in Untersuchungen dieser Art hat wohl als erster Gauß in seiner Arbeit \textit{``Theoria residuorum biquadraticorum. Commentatio secunda."} (\cite{Ga31}, 1831) (``Theorie der biquadratischen Reste -- zweite Mitteilung") im Kontext der nach ihm benannten Gauß'schen Zahlen $\mathbb{Z}[\sqrt{-1}]$ explizit erwähnt. Euler hat bei seinem Nachweis gleichsam die Eindeutigkeit der Eisensteinzahlen $\mathbb{Z}[\omega]$ mit $\omega^2+\omega+1=0$ genutzt. Diese wurden von Eisenstein in seiner Arbeit \textit{``Beweis des Reciprocitätssatzes für die cubischen Reste in der Theorie der aus den dritten Wurzeln der Einheit zusammengesetzten Zahlen"} (\cite{Ei44}, 1844) entsprechend eingeführt und zum Beweis des kubischen Reziprozitätsgesetzes herangezogen.}.
\\

Bei diesem Euler'schen Beweis folgender interessanter Fall vor: Einerseits muss, das heutige Stringenzmaßstäbe anlegend, Eulers Argument als unvollständig angesehen werden --  er weist die Legitimität des unendlichen Abstiegs nicht nach\footnote{Dazu benötigt man die Eindeutigkeit der Primfaktorzerlegung im Ring $\mathbb{Z}[\sqrt{-n}]$, welche für beispielsweise für $n=5$ nicht vorliegt, zumal zum einen $6=(1+\sqrt{-5})(1-\sqrt{-5})=2 \cdot 3$ gilt, jedoch zum anderen all die Zahlen $2$, $3$, $1+\sqrt{-5}$, $1-\sqrt{-5}$ Einheiten in diesem Ring sind. Man vergleiche in diesem Zusammenhang auch die Ausführung von Weyl in seinem Buch \textit{``Number theory. An approach through history. From Hammurapi to Legendre"} (\cite{We84}, 1984), welcher den Euler'schen Beweis des großen Satzes von Fermat für den Fall $n=3$, trotz des erwähnten Mangels, als Meilenstein bei der Entstehung der algebraischen Zahlentheorie einordnet.}, andererseits wurden solche Existenzfragen\footnote{Überdies scheint sich bei Euler kein Beispiel einer Zahl aus $\mathbb{Z}[\sqrt{-n}]$ zu finden, woran sich eine nichteindeutige Zerlegung für ihn hätte erkennen lassen, wenn man alleinig seine Veröffentlichungen heranzieht. Im Buch \textit{``4000 Jahre Zahlentheorie"} (\cite{Le23}, 2023) wird indes auf einen Tagebucheintrag \textit{``Fragmenta arithmetica ex Adversariis mathematicis deprompta"} (\cite{E806}, 1862, ges. ?) (E806: ``Auszüge aus den mathematischen Tagebüchern zum Gegenstand der Zahlentheorie") verwiesen, wo die Gleichung $181^2+7=32^3$ die Euler'sche Schlussweise zum Nachweis des Fermat'schen Satzes für den Fall $n=3$ als unvollständig ausweist, zumal hier nicht $181 \pm \sqrt{-7}=(5\pm \sqrt{-7})^3$ gilt. Es mag von Interesse sein, dass das von Euler gefundene Beispiel aus den elliptischen Kurven der Gestalt $Y^2=X^3-N$ mit $N \in \mathbb{N}$ das für $N=7$ ist, demnach das für das kleinste $N$, sodass man auf der elliptischen Kurve mehr als eine Lösung über den natürlichen Zahlen findet -- das andere ist $Y=1$ und $X=2$. In seiner \textit{``Algebra"} \cite{E387} findet man in dieser Hinsicht nur den Fall $\mathbb{Z}[\sqrt{10}]$ bzw. Zahlen der Form $a+b\sqrt{10}$ diskutiert.} zu Eulers Lebzeiten  selten bis gar nicht gestellt, sodass -- aus dem Euler'schen Paradigma heraus argumentiert -- der Beweis als vollständig angesehen werden muss. Diese Betrachtung kann demnach als stellvertretendes Exempel für den stattgefundenen Paradigmenwechsel in der Mathematik herangezogen werden.

\newpage

\section{Mathematische Grenzen Eulers}
\label{sec: Grenzen Eulers gezogen durch das Paradigma}

\epigraph{The discovery of truth is prevented more effectively, not by the false appearance things present and which mislead into error, not directly by weakness of the reasoning powers, but by preconceived opinion, by prejudice.}{Arthur Schopenhauer}

Am Beispiel des quadratischen Reziprozitätsgesetzes (Abschnitt \ref{subsubsec: Wegen fehlender Mittel: Das Reziprozitätsgesetz}) ist bereits auf von Euler nicht zu überwindende Grenzen verwiesen worden. War dies in diesem Fall wegen einer unhandlichen Formulierung und fehlender algebraischer Mittel mehreren Gründen geschuldet, manifestieren sich bei anderen Begebenheiten überdies  Grenzen gänzlich anderer, nämlich psychologischer, Natur: Hier sind vor allem die Begriffe selbst zu nennen, deren Verständnis nicht erwartete Restriktionen schaffen oder gar zu unauflösbaren Widersprüchen führen  können, sofern sie nicht erweitert werden. Solche restringierenden Konzeptbildungen Eulers werden in Abschnitt (\ref{subsec: Grenzen wegen Begriffsbildung}) angeführt. Jedoch verhindern selbst korrekte Definitionen und Begriffe nicht das Stellen einer irreleitenden Frage, welchem Umstand Euler auch anheim gefallen ist (Abschnitt \ref{subsec: Durch eine falsche Frage}). Während eine solche  Frage explizit formuliert werden muss, existieren gleichermaßen analoge implizite Fragestellungen, die ihren Ursprung im vorherrschenden Paradigma selbst haben, wie unterem Kuhn in seinem Buch \textit{``The Structure of Scientific Revolutions"} (\cite{Ku12}, 2012)\footnote{Es handelt sich bei der zitierten Quelle um  50th Anniversary Edition.} herausgearbeitetet hat. Bedenkt man die führende Rolle Eulers bei der Konstruktion desselben in der Mathematik im 18. Jahrhundert, erlaubt dies zugleich einen Blick auf Eulers Gedankengebäude selbst.  Verwandt damit ist die Untersuchung nach durch den eigenen Arbeitsethos gleichsam selbst unwillentlich gesetzten Schranken, welche in Abschnitt (\ref{subsec: Grenzen durch den eigenen Arbeitsethos}) anhand von Beispielen illustriert werden.

\subsection{Aus der Begriffsbildung resultierende Grenzen}
\label{subsec: Grenzen wegen Begriffsbildung}

\epigraph{To the reader of today much in the conception and mode of expression of that time appears strange and unusual. Between us and the mathematicians of the late seventeenth century stands Leonhard Euler [...] He is the real founder of our modern conception.}{Josef Ehrenfried Hofman}

Es ist bereits an diversen Stellen angeklungen, dass teilweise allein das Fehlen eines modernen Begriffs Euler am Weiterkommen gehindert hat, welche Behauptung nun exemplarisch untermauert werden soll. Zu diesem Zweck wird auf den zentralen Begriff der Funktion (Abschnitt \ref{subsubsec: Der Begriff der Funktion}) und damit verknüpfter Konzepte eingegangen werden, auch Eulers Auffassung von Grenzwerten wird besprochen werden (Abschnitt \ref{subsubsec: Der Begriff des Grenzwerts}), ein Begriff, welchen Euler noch nicht im modernen Sinne begreift. Am Beispiel seines Konzepts einer Summe einer Reihe wird anschließend (Abschnitt \ref{subsubsec: Der Begriff der Summe einer Reihe}) demonstriert werden, wie viel Euler mit einer aus moderner Sicht zum Scheitern verurteilten Definition zu leisten vermochte, bevor schlussendlich die Konsequenzen für ein bei Euler gänzlich fehlendes Konzept beleuchtet werden. Es wird sich zeigen (Abschnitt \ref{subsubsec: Mehrwertige Funktionen: Das fehlende Konzept der Riemann'schen Fläche}), wie mehrwertige Funktionen Euler wegen des Fehlens der Idee der Riemann'schen Fläche zu -- auch nach seinem eigenen Verständnis  -- falschen Aussagen diesbezüglich führen musste.

\subsubsection{Der Begriff der Funktion und damit eng verwandte}
\label{subsubsec: Der Begriff der Funktion}

\epigraph{The difficulty lies, not in the new ideas, but in escaping from the old ones.}{John Maynard Keynes}

Euler hat das Konzept der Funktion nachdrücklich in die Analysis eingebracht; heute ist sich die Mathematik nicht mehr ohne diesen Begriff zu denken. Darum wird es förderlich sein, die Bedeutung des Begriffs bei Euler selbst zu beleuchten und Unterschiede zur modernen Auffassung zu illustrieren. Dies wird zeigen, dass Euler neben dem Funktionsbegriff selbst überdies abweichende Anschauungen von ihren Eigenschaftsbegriffen wie  etwa dem der Stetigkeit hatte. 

\paragraph{Entwicklung des Funktionsbegriffs bei Euler}

 In §4  seiner \textit{Introductio} \cite{E101}  gibt Euler die folgende Definition des Begriffs \textit{Funktion}: \\

     \textit{``Eine Funktion einer variablen Größe ist ein analytischer Ausdruck, beliebig zusammengesetzt aus dieser variablen Größe und Zahlen oder konstanten Größen."} \\

 Später schreibt er im Vorwort seiner \textit{Calculi Differentialis}  \cite{E212}: \\

    \textit{``Wenn [...] $x$ eine variable Größe bezeichnet, werden alle Größen, die von ihr in beliebiger Weise abhängen, Funktionen von ihr genannt."} \\

Schließlich findet man in seinem Übersichtsartikel \cite{E322} in § 1  die Definition wiederholt, dass  \textit{eine Funktion eine Größe ist, beliebig definiert aus bestimmten Größen.}\\

\begin{figure}
    \centering
    \includegraphics[width=0.8\linewidth]{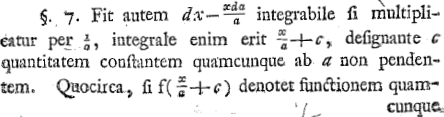}
    \caption{Euler nutzt in \cite{E45} erstmalig die Notation $f(x)$ für eine Funktion von $x$.}
    \label{fig:E45Funktion}
\end{figure}

Euler ändert seine Vorstellung des Funktionsbegriffs über die Jahre hinweg demnach wenig. Vergleicht man dies nun mit der modernen Auffassung des Begriffs der Funktion, also als Abbildungsvorschrift zwischen Mengen, die jedem Element der Ausgangsmenge ein Element aus der Zielmenge zuordnet, so fällt das Fehlen des Konzepts der Menge auf und somit auch all dessen, was sich daraus ableitet. Genauer findet im Buch \textit{``3000 Jahre Analysis: Geschichte - Kulturen - Menschen"}  (\cite{So16}, 2016) (SS. 512--513) die Bourbaki'sche Definition aus dem Buch \textit{``Theory of Sets"}  (\cite{Bo68}, 1968) von Sonar (1958--) wie folgt zitiert:\\

\textit{``Seien $E$ und $F$ zwei Mengen, die nicht notwendig verschieden sein müssen: Eine Relation zwischen einem veränderlichem Element $x$ von $E$ und einem veränderlichem Element $y$ von $F$ heißt funktionale Relation in $y$, wenn für alle $x \in E$ ein eindeutig bestimmtes $y \in F$ existiert, das in der gegebenen Relation mit $x$ steht. Wir vergeben den Namen Funktion für die Operation, die in dieser Art mit jedem Element $x \in E$ das Element $y \in F$ assoziiert, das in der gegebenen Relation $x$ steht; $y$ heißt der Wert der Funktion an dem Element $x$ und die Funktion heißt bestimmt durch die gegebene funktionale Relation. Zwei äquivalente funktionale Relationen bestimmen dieselbe Funktion."}\\

Das Fehlen der eindeutigen Zuordnung der Elemente in der Euler'schen Definition lässt Euler, eines Beispiels wegen, den Ausdruck $f: x \mapsto \sqrt{x}$ als eine Funktion sehen, die   für jeden reellen (und auch komplexen Wert) \textit{zwei} Werte hat.  Euler spricht in diesem Zusammenhang von mehrwertigen Funktionen, welche den Inhalt von Abschnitt (\ref{subsubsec: Mehrwertige Funktionen: Das fehlende Konzept der Riemann'schen Fläche}) bilden werden.  Die moderne Auffassung verbietet hingegen solche Auslegungen wie die Euler'sche und schreibt im Fall der obigen Funktion die Wahl eines Zweiges der Wurzel vor\footnote{Der Euler'schen Auffassung bezüglich mehrdeutiger Funktionen kommt heute das Konzept der \textit{Korrespondenz} noch am nächsten, selbiges verlangt aber auch die Kenntnis von Mengen, welche bei Euler noch fehlen.}.

\paragraph{Stetigkeit}

Neben Unterschieden im Begriff der Funktion selbst, haben einige Euler'sche  Eigenschaftsbegriffe derselben in der von ihm vorgeschlagenen Form keinen Eingang in moderne Mathematik gefunden, was hier am Beispiel des Konzepts der Stetigkeit gezeigt werden soll. Eine präzise Definition von Stetigkeit gibt Euler in seiner Arbeit \cite{E322} ein. Hier schreibt er nämlich in § 2:\\

\textit{``Nun ist indes allbekannt, dass in der höheren Geometrie lediglich Kurven betrachtet zu werden pflegen, deren Natur, mit einer gewissen Relation zwischen den Koordinaten ausgedrückt, über eine gewisse Gleichung bestimmt ist, sodass all ihre Punkte durch diese Gleichung gleichsam wie durch ein Gesetz festgelegt werden. Weil dieses Gesetz verstanden wird, das Prinzip der Stetigkeit in sich zu umfassen, nach welchem freilich alle Teile der Kurve mit einem so starken Band miteinander zusammenhängen, dass in jenen eine Veränderung bei vorliegender Stetigkeit nicht auftreten kann, werden dieses Grund wegen diese gekrümmten Linien stetig genannt, [...], solange wir nur verstehen, dass eine gewisse Gleichung gegeben ist, mit welcher die Natur dieser Kurven zum Ausdruck gebracht wird."} \\

So fallen für Euler etwa, trotz ihrer beiden Äste, die Hyperbeln unter die stetigen Kurven, so wie  alle anderen Kegelschnitte. Gleichermaßen wie alle anderen Kurven, welche durch eine Gleichung der Form $F(x,y)=0$ gegeben sind. Zu unstetigen Funktionen schreibt Euler indes in § 3 derselben Arbeit:\\

\textit{``Nachdem also das Kriterium für Stetigkeit festgelegt worden ist, ist unmittelbar ersichtlich, was eine unstetige  oder eine dem Kontinuitätsgesetz nicht unterworfene Funktion ist: Denn alle durch keine feste Gleichung bestimmten gekrümmten Linien, von welcher Art etwa freihändig gezeichnete bezeichnet zu werden pflegen, geben solche unstetigen Funktionen an die Hand, weil ja in  ihnen die Werte der Ordinaten nach keiner festen Vorschrift aus den Abszissen bestimmt werden können."} \\

Kontrastiert man die Euler'schen Ausführungen etwa mit der modernen Definition, tritt die Euler'sche Emphase bezüglich der notwendig vorhandenen Gleichung hervor,  was  in modernen Definitionen von Stetigkeit keine Erwähnung findet. Ebenso ist Eulers Klassifikation der freihändig gezeichneten Funktionen zu den unstetigen Funktion mit der modernen Auffassung unverträglich. Zwecks Illustration sei an dieser Stelle die Dirichet'sche Definition des Funktionen-- und Stetigkeitsbegriffs zitiert. In seiner Arbeit \textit{``Über die Darstellung ganz willkürlicher Functionen durch Sinus-- und Cosinusreihen"} (\cite{Di37}, 1837) findet sich seine Erläuterung zu diesem Gegenstand: \\

\textit{``Entspricht nun jedem $x$ ein einziges, endliches $y$,  und zwar so, dass, während $x$ das Intervall von $a$ bis $b$ stetig durchläuft, $y=f(x)$ sich ebenfalls allmählich verändert, so heisst $y$ eine stetige oder continuierliche Function von $x$ für dieses Intervall. Es ist dabei gar nicht nöthig, dass $y$ in diesem ganzen Intervalle nach demselben Gesetze von $x$ abhängig sei, ja man braucht nicht einmal an eine durch mathematische Operationen ausdrückbare Abhängigkeit zu denken."}\\

Dirichlets Konzept der Stetigkeit weist somit insbesondere bereits das Euler'sche  zurück. Was Euler zu seiner Auffassung von Stetigkeit bewogen hat, wird wohl Spekulation bleiben. Jedoch wird dabei die praktische Anwendbarkeit die Inklusion einer Gleichung in der Definition begünstigt haben, welche -- wie Dirichlet erkannt hat -- zur Ableitung allgemeiner Eigenschaften nicht bekannt zu sein braucht. Überdies sind natürlich heute unstetige Funktionen bekannt, die sich mit einer Gleichung explizit angeben lassen\footnote{Euler scheint in seinen Arbeiten nicht das formal definiert zu haben, was man heute eine abschnittsweise definierte Funktion nennt. Er beschreibt die Funktion gegebenenfalls mit Worten, wie etwa in \cite{E322}, aber belässt es dann dabei.}. Man denke etwa an das bekannte Beispiel der Dirichlet'schen Sprungfunktion.

\begin{equation}
\label{eq: Dirichlet}
    D(x):=
    \renewcommand{\arraystretch}{1,0} 
    \setlength{\arraycolsep}{0.0mm}
    \left\lbrace 
    \begin{array}{l}
       1, \quad \text{für $x$ rational}  \\
       0, \quad \text{für $x$ irrational}  \\
    \end{array}
    \right.
\end{equation}

Hier könnte man, in Eulers Sinne, entgegnen, dass solche Funktionen nicht betrachtet zu werden brauchen. Begründet wird dies von Euler von § 4 in \cite{E322} über ihre Unnatürlichkeit. Er schreibt:\\

\textit{``Dass nun allen unstetigen Linien und Funktionen von dieser Art in der geometrischen Analysis kein Raum gelassen wird, ist per se offenkundig, weil diese ganze Betrachtung vom Ausfindigmachen Eigenschaften der Linien, welche untersucht werden, eingenommen ist, welche Aufgabe in keiner Weise in Angriff genommen werden könnte, wenn die Natur dieser Linien nicht in einem gewissen Gesetz oder einer Gleichung enthalten wäre."}\\

Euler plädiert jedoch später in seiner Arbeit dafür,  auch solche Funktionen in der Analysis zuzulassen -- man ist gar dazu gezwungen -- weil andernfalls das Problem der schwingenden Saite\footnote{Dies ist unter anderem seiner Untersuchung in \cite{E119} zur schwingenden Saite geschuldet.} sich als unlösbar herausstellen würde. Genauer schreibt er diesbezüglich in  § 6 von \cite{E322}:\\

\textit{``Um diesen Streit beizulegen, bemerke ich, dass weder in der gemeinen Algebra noch im dem Teil der Analysis des Unendlichen, der bisher hauptsächlich behandelt worden ist, unstetige Funktionen zugelassen werden können. Aber die Analysis des Unendlichen ist natürlich zu verstehen, sich um vieles weiter zu erstrecken und auch solche Zweige zu besitzen, welche von unstetigen Funktionen nicht nur nicht zurückweichen, sondern sie gar nach ihrer Natur  so beinhalten, dass kein sich darauf beziehendes Problem als vollständig gelöst  anzusehen ist, wenn nicht vollkommen beliebige, und daher auch unstetige,  Funktionen in die Lösung eingeführt worden sind."}\\

Zusammenfassend gesteht Euler also -- notgedrungen -- die Wichtigkeit unstetiger Funktionen in der Analysis ein\footnote{Diese Ansicht wird auch von Lützen in seiner Arbeit \textit{``Euler’s Vision of a General Differential Calculus for a Generalized Kind of Function"} (\cite{Lu83}, 1983) geteilt, welcher die Euler'sche Auffassung des Funktionsbegriffs zum Gegenstand hat.}, jedoch nur im Falle mehrerer Variablen,  gleichsam naturalistisch begründet\footnote{Einmal mehr mag man aus diesen Aussagen den Euler'schen Gedanken herauslesen, die Mathematik sei so eng mit den Naturgesetzen verwoben, dass sie gleichsam deckungsgleich sind. Was nicht durch physikalische Prinzipien erläutert werden kann, ist zugleich auch unnütz für die Mathematik.}. in der Mathematik  unstetige Funktionen zuzulassen. Demnach wäre es Euler wohl auch nicht in den Sinn gekommen, Funktionen zu betrachten, die zwar stetig sind, jedoch nicht differenzierbar. Hat  man nämlich eine Gleichung vorgegeben, so hat man nach Euler zum Bilden der Differentiale statt $x$ entsprechend $x+dx$ und statt $y$ entsprechend $y+dy$ zu setzen, um anschließend unter Verwendung der Ausgangsgleichung das Verhältnis $\frac{dy}{dx}$ zu ermitteln, was bei Vorlage einer die Funktion definierenden Gleichung, für Euler, stets möglich ist. \\

\paragraph{Integrale}

Neben dem Konzept der Stetigkeit, welcher einer Funktion als Eigenschaft zukommt, soll mit dem Integral ein Begriff betrachtet werden, welcher auf Funktionen angewendet wird. Das Integral definiert er konkret in seinem Buch \textit{``Institutionum calculi integralis volumen primum"} \cite{E342} (``Grundlagen des Integralkalküls, erster Band"). Hier schreibt er in § 1:\\

\textit{``Das Integralkalkül ist die Methode aus einem gegebenen Verhältnis der Differentiale die Beziehung der Größen selbst zu finden: Und die Operation, mit welcher dies geleistet wird, pflegt Integration genannt zu werden."}\\

Obwohl Euler an anderer Stelle, die geometrische Interpretation der Integrale als Flächeninhaltsbilanz zwischen einem Graphen und der $x$--Achse verwendet, etwa in seinen Arbeiten \textit{``Methodus universalis serierum convergentium summas quam proxime inveniendi"} (\cite{E46}, 1741, ges. 1735) (E46: ``Eine allgemeine Methode, die Summen konvergenter Reihen näherungsweise zu ermitteln"), die einen wichtigen Vorstoß zur Euler--Maclaurin'schen Summenformel bedeutet, sowie   \cite{E587}, wo die elementaren Eigenschaften des Integrals wie Linearität, Additivität und Verhalten und Grenzvertauschung  aus der geometrischen Anschauung heraus erklärt werden (dort §§ 1--7) --, leitet Euler  zumeist die Auffassung des Integrals als Inverse Operation des Differenzierens. So schreibt Euler im Zusammenhang mit Differentialgleichungen nahezu nie, dass diese \textit{gelöst} werden müssen, sondern dass diese eben zu \textit{integrieren} sind. Dass Euler auch wegen des fehlenden Begriffs des Grenzwerts die moderne Zugänge von Riemann vermöge Unter-- und Obersummen oder gar Lebesgue (1875 -- 1941) über Treppenfunktion verwehrt bleiben mussten, braucht an dieser Stelle entsprechend nicht weiter ausgeführt zu werden. Euler wäre auf immense Probleme in den Beispielen gestoßen, welche sich mit der einen Methode evaluieren lassen, mit der anderen hingegen nicht. Ein klassisches Beispiel ist:

\begin{equation*}
    \int\limits_{0}^{\infty} \dfrac{\sin x}{x}dx= \dfrac{\pi}{2},
\end{equation*}
welches sich bei Euler auch schon in \cite{E675} als Spezialfall einer allgemeineren Formel findet und mit der Riemann'schen Idee ebenfalls ermittelt werden kann. Mit der Lebesgue'schen Integrationsvorschrift kann man diesem Integral hingegen keinen Wert zuordnen. Diese erlaubt dahingegen die Berechnung des Integrals der Dirchlet'schen Funktion (\ref{eq: Dirichlet}), Riemanns Vorgehen kann den Wert  nicht ausgeben. \\

\begin{figure}
    \centering
   \includegraphics[scale=0.9]{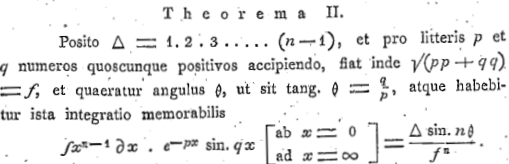}
    \caption{Euler gibt in seiner Arbeit \cite{E675} eine allgemeine Formel zur Berechnung vieler Integrale an.}
    \label{fig:E675GammaInt}
\end{figure}

Euler versteht ein Integral gleichermaßen als natürliche Verallgemeinerung des Begriffs der Summe. Der Begriff wird von den natürlichen Zahlen auf die reellen Zahlen ausgeweitet. Dies zeigt sich insbesondere in seinem Buch zur Variationsrechnung, der \textit{Methodus}, \cite{E65} und auch in physikalischen Rechnungen war diese Herangehensweise für ihn überaus hilfreich.  Euler scheint mit dieser Auffassung sogar nie auf konzeptuelle Probleme gestoßen zu sein\footnote{Divergente Integrale, welche jedoch physikalisch motiviert endliche Werte ausgeben müssten -- man denke hier etwa an die Selbstenergie eines Elektrons oder gar die divergenten Integrale in der Quantenfeldtheorie -- scheinen bei Euler noch nicht aufgetreten zu sein.}.   Insgesamt drängte Euler nichts zu einer Modifikation seiner Vorstellung des Konzepts der Integration.


\subsubsection{Der Begriff des Grenzwerts und damit verknüpfte Konzepte}
\label{subsubsec: Der Begriff des Grenzwerts}

\epigraph{A new word is like a fresh seed sown on the ground of the discussion.}{Ludwig Wittgenstein}

Grenzwerte im heutigen Sinne waren Euler noch fremd, selbige wurden erst  später ca. 40 Jahre nach Eulers Tod formal eingeführt;  der moderne Begriff stammt von Cauchy aus seinem Lehrbuch \textit{``Cours d'analyse de l'École royale polytechnique"} (\cite{Ca21}, 1821) (``Kursus der Analysis der École Royale Polytechnique"). Die Notwendigkeit eines präzisen Grenzwertbegriffs zeigt sich beim Auftreten der ersten  Widersprüche, zu denen man  geleitet wird, was bei Euler bei Differentialen und seiner Idee von infinitesimalen Größen geschah.

\paragraph{Bei Differentialen}
\label{para: Bei Differentialen}

Euler hat sehr detailliert beschrieben, was man unter Differentialen zu verstehen hat. Eine erste Erklärung findet man im Vorwort von seiner \textit{Calculi Differentialis} \cite{E212}, woraus jedoch zunächst einmal seine Definition der Differentialkalküls zitiert sei. Euler schreibt auf Seite VII:\\

\textit{``Das Differentialkalkül, das ist die Methode, das Verhältnis verschwindender Inkremente zueinander zu bestimmen, welche beliebige Funktionen erhalten, während der variablen Größe, von welcher sie Funktionen sind, verschwindende Inkremente zugeteilt werden."}\\

Dieser Ausführung sei folgende aus \cite{E322} an die Seite gestellt, wo Euler in § 9 schreibt:\\

\textit{``Wenn also zuerst $y$ irgendeine Funktion einer einzigen Variable $x$ war, pflegen die Zuwächse oder Abnahmen jener Funktion $y$ betrachtet zu werden, während die Größe $x$ um einen beliebigen Wert wächst. Dann wird dieses Inkrement aufgefasst, immer weiter verringert zu werden, bis es schließlich völlig verschwindet, in welchem Fall freilich auch das Inkrement jener Funktion $y$ verschwindet; weil diese verschwindenden Inkremente Differentiale genannt zu werden pflegen, ist es ersichtlich, dass sie überhaupt keine Größe haben und daher gleich Null sind, sodass bezüglich ihrer Größe überhaupt keine Frage aufkommen kann. Aber dennoch besteht das Differentialkalkül nicht aus dem Finden der Größe der Differentiale selbst, welche ja nicht vorhanden ist, sondern in der Bestimmung ihres Verhältnisses zueinander."}\\

 Die Euler'sche Begriffserklärung kann etwaige Schwierigkeiten, da wegen $dx=0$ auch $2dx=dx$, woraus man $2=1$ ableiten müsste, obwohl $\frac{dx}{dx}=1$ sein soll, nicht umgehen, wenn man auch, wie Euler es in seinen Arbeiten tut und überdies in seinen \textit{Calculi Differentialis} \cite{E212} auch erläutert, mit den Differentialen als absolute Größen rechnen möchte. Im Gegensatz dazu umgeht die moderne Definition über Grenzwerte

\begin{equation*}
   y'(x):=\lim_{h \rightarrow 0} \dfrac{y(x+h)-y(x)}{x+h-x}
\end{equation*}
das Problem der Euler'schen Definition von Differentialen, indem sie gar nicht erst genannt werden\footnote{Zudem versteht man unter einem Differential heute viel allgemeinere Objekte als Euler.}. An seinen Ausführungen lässt Euler erneut seine charakteristische Art durchscheinen,  jedem Objekt in irgendeiner Form eine Größe zuzuweisen, um mit selbiger greifbar und auch absolut, nicht nur in Verhältnissen\footnote{Sonar bezeichnet daher in seinem Buch \textit{``3000 Jahre Analysis: Geschichte - Kulturen - Menschen"} (\cite{So16}, 2016) die Euler'schen Ausführungen zu seinem Differentialkalkül als ``Nullenrechnung"{} und würdigt auch die Leistungen Eulers in dieser Hinsicht entsprechend.}, rechnen zu können\footnote{Das Euler'sche Kalkül hat mehr als 150 Jahre nach Eulers Tod durch die Nicht-Standardanalysis, erläutert  im Buch \textit{``Nonstandard Analysis"} (\cite{Ro66}, 1966), auch eine strenge Rechtfertigung erhalten.}. \\

 Die Grundlage der Differentialrechnung besteht also für Euler darin, dass Differentiale eigentlich Differenzen sind\footnote{Eine umfassende Diskussion der Euler'schen Ideen findet sich in Kapitel 3  des Buchs \textit{``Euler as Physicist"} (\cite{Su10}, 2010) von Suisky.}; aber nicht nur unendlich kleine, sondern sogar null. Da aber der Begriff des Grenzwerts bei Euler noch fehlt, hätte er etwa nicht zu erklären vermocht, weshalb die partiellen Ableitungen einer Funktion von zwei oder mehr Veränderlichen im Allgemeinen nicht vertauscht werden dürfen\footnote{Man vergleiche auch Sandifers Ausführungen in seiner Arbeit \textit{``Mixed Partial Derivatives"} (\cite{Sa04mai}, 2004) zu diesem Gegenstand.}, obgleich er die Vertauschbarkeit wie folgt in § 7 von \cite{E44} beweisen will:\\

 \textit{``Obwohl aber die Gültigkeit dieses Lehrsatzes demjenigen, der ihn prüfen will, leicht ersichtlich wird, möchte ich dennoch den folgenden aus der Natur der Differentiale abgeleiteten Beweis hinzufügen. Während $A$  eine Funktion von $t$ und $u$ ist, gehe $A$ in $B$ über, wenn $t+dt$ anstelle von $t$ geschrieben wird; aber nach Setzen von $u+du$ anstelle von $u$ gehe $A$ in $C$ über. Nachdem aber zugleich $t+dt$ anstelle von $t$ und $u+du$ anstelle von $u$ geschrieben worden ist,  ändere sich $A$ zu $D$. Aus diesen Festlegungen ist ersichtlich, dass, wenn  in $B$ dann $u+du$ anstelle von $u$ geschrieben wird,  $D$ hervorgeht und wenn in $C$ entsprechend $t+dt$ anstelle von $t$ gesetzt wird, dass gleichermaßen $D$ hervorgeht. Nachdem all dies vorausgeschickt worden ist, wird, wenn nun $A$ für konstantes $t$ differenziert wird, $C-A$ hervorgehen, denn nach Setzen von $u+du$ anstelle von $u$ geht $A$ in $C$ über, aber das Differential ist $C-A$. Wenn in $C-A$ weiter $t+dt$ anstelle von $t$ gesetzt wird, wird $D-B$ hervorgehen, woraus das Differential $D-B-C+A$ sein wird. Nachdem aber in umgekehrter Reihenfolge $t+dt$ anstelle von $t$ in $A$ gesetzt worden ist, wird man $B$ haben, und das Differential von $A$ für allein variables $t$ wird $B-A$ sein. Dieses Differential geht für $u+du$ anstelle von $u$ gesetzt in $D-C$ über, woraus das Differential davon wiederum $D-B-C+A$ sein wird, was mit dem Differential übereinstimmt, das mit der ersten Operation gefunden worden ist."}\\
 
 Eulers Argumentation wird durch das von Peano (1858--1932) angegebene\footnote{Es wird unter anderem in Hobsons (1856--1933) Buch \textit{``The theory of functions of a real variable and the theory of Fourier's series. Vol. I"} (\cite{Ho21}, 1921) Peano zugeschrieben.} und auch noch heute im Kontext des Satzes von Schwarz (1843--1921) vorgestellte Beispiel

\begin{equation*}
    f(x,y)= \dfrac{xy^3-x^3y}{x^2+y^2} \quad \text{für} \quad (x,y) \neq (0,0)
\end{equation*}
als fehlerhaft entlarvt\footnote{Der Fehler in Eulers Beweis besteht in der Annahme der Möglichkeit, man könne \textit{gleichzeitig} in beiden Variablen die Differentiale einsetzen. Dies entspräche modern aufgefasst der Annahme, dass beim Differenzieren nach den Variablen $x$ und $y$ keine Reihenfolge festzulegen ist.}, für welches eben nicht gilt:

\begin{equation*}
    \dfrac{\partial}{\partial x}\dfrac{\partial}{\partial y}f= \dfrac{\partial}{\partial y}\dfrac{\partial}{\partial x}f,
\end{equation*}
wenn man den Ursprung betrachtet. Daher findet sich im Satz von Schwarz die Bedingung der Stetigkeit der zweiten partiellen Ableitungen, welche  notwendig ist, will man die Vertauschbarkeit der zweiten Ableitungen garantieren. Selbst wenn man die eben vorgestellte Euler'sche  Definition von Stetigkeit aus \cite{E322} anwendet, hätte Euler den Schwarz'schen Satz nicht rechtfertigen können, da er auch hier zu einem Widerspruch in seinen Ansichten geführt worden wäre. Man sieht das wie folgt ein:\\
 Durch seine Untersuchungen zu Funktionen im Allgemeinen  sah Euler zur Annahme veranlasst (siehe die entsprechenden Ausführungen in Abschnitt \ref{subsubsec: Durch einen neuen Gedanken: Die Theta-Funktion}.), mit Potenzreihen viele, wenn nicht gar alle, Funktionen darstellen zu können, sofern man beliebige Potenzen zulässt. Insbesondere findet man  auch die Ableitungen der Funktionen im Ursprung unmittelbar aus der Potenzreihenentwicklung, wie  der Satz von Taylor lehrt\footnote{Auf diesem Wege \textit{definiert} Lagrange gar die Ableitung einer Funktion in seinem Buch \textit{``Théorie des Fonctions analytiques"} (\cite{La97}, 1797) (``Theorie der analytischen Funktionen"). Er umgeht damit ebenfalls die unendlichen kleinen Größen von Euler.}. Bezogen auf das Peano'sche Exempel führt dies zur Ableitung $ \frac{\partial}{\partial x}\frac{\partial}{\partial y}f(0,0)$ sowie die von $ \frac{\partial}{\partial y}\frac{\partial}{\partial x}f(0,0)$ aus dem Koeffizienten von $xy$ bzw. $yx$. Entwickelt man aber nun 

\begin{equation*}
    \dfrac{1}{x^2+y^2}= \dfrac{1}{y^2}\left(1-\dfrac{x^2}{y^2}+\cdots\right) \quad \text{bzw.} \quad = \dfrac{1}{x^2}\left(1-\dfrac{y^2}{x^2}+\cdots\right)
\end{equation*}
findet man entsprechend aus den Koeffizienten $xy$ von $f$ hieraus einmal $+1$ und einmal $-1$, somit auch

\begin{equation*}
\dfrac{\partial}{\partial x}\dfrac{\partial}{\partial y}f(0,0)=-1 \quad \text{und} \quad \dfrac{\partial}{\partial y}\dfrac{\partial}{\partial x}f(0,0)=+1,
\end{equation*}
was Euler vermutlich aufgestoßen wäre, hätte er davon eine Erwähnung gemacht. Andererseits zeigt eine allgemeine Rechnung gemäß des Euler'schen ``Beweises"{} die Nicht--Vertauschbarkeit nicht auf. Dann findet man natürlich zunächst:

\begin{equation*}
    \dfrac{\partial^2}{\partial x \partial y} f=  \dfrac{\partial^2}{\partial y \partial x} f = \dfrac{x^6+9x^4y^2-9x^2y^4-y^6}{(x^2+y^2)^3}.
\end{equation*}
Lediglich die Reihenfolge, in welcher man $x$ und $y$ gegen $0$ streben lässt,  kann den Widerspruch aufzeigen. Es ist eine interessante Diskussion, wie Euler dieses Problem wohl gedeutet hätte, wäre es ihm nur aufgefallen. Bezugnehmend auf seinen Beweisversuch des Vertauschungssatzes, ist anzunehmen, dass Euler $x$ und $y$ als gleichberechtigt behandelt hätte,  zunächst $x=y=a$ gesetzt hätte und dann $a$ gegen $0$ laufen lassen hätte oder, in seiner Sprache, $a$ unendlich klein hätte werden lassen. Man findet auf diesem Wege

\begin{equation*}
     \dfrac{a^6+9a^4a^2-9a^2a^4-a^6}{(a^2+a^2)^3}= 0,
\end{equation*}
was ja gerade der Mittelwert der obigen Werte $+1$ und $-1$ ist\footnote{Diesem Ergebnis wäre Euler zweifelsohne sehr zugetan gewesen, zumal sich die Situation ähnlich verhält wie bei der Tatsache, dass 

\begin{equation*}
    1-1+1-1+1-1+\text{etc.}=\dfrac{1}{2}
\end{equation*}
als Mittelwert von $0$ und $1$ aufgefasst werden kann.}. Diesen Fehlschluss entlarvt man an solchen, wenn man $x=r\cos(\alpha)$ und $y=r\sin(\alpha)$ werden lässt, sodass unter Verwendung trigonometrischer Identitäten gilt:

\begin{equation*}
   \dfrac{x^6+9x^4y^2-9x^2y^4-y^6}{(x^2+y^2)^3}= \dfrac{1}{2}\left(2\cos(2\alpha)-\cos(6\alpha)\right),
\end{equation*}
was wegen des Herausfallens des Radius $r$ die Wegabhängigkeit des Grenzwerts beim Zulaufen auf den Ursprung zeigt. Jedoch erfährt die vermutliche Euler'sche Ansicht durch die Betrachtung des ``Durchschnitts"{} aller Grenzwerte

\begin{equation*}
  \left\langle \dfrac{1}{2}\left(2\cos(2\alpha)-\cos(6\alpha)\right) \right\rangle :=  \dfrac{1}{2\pi} \int\limits_{0}^{2\pi} \dfrac{1}{2}\left(2\cos(2\alpha)-\cos(6\alpha)\right)dx=0
\end{equation*}
eine Rechtfertigung. \\

Neben dem Beispiel der Differentiation stellten die Differentiale Euler auch bei der Integration vor Schwierigkeiten -- wieder sind die Funktionen dabei von zwei Variablen abhängig. Das tiefer liegende Problem ist freilich, dass Euler in seinen Betrachtungen  zwar $(du)^2=0$ bzw. $(dt)^2=0$ benutzt, dies aber damit begründet, dass die Quadrate lediglich in Bezug auf $du$ und $dt$  verworfen werden müssen.  Heute ist dies eine Konsequenz der Antikommutativität der Differentiale bzw. der Gleichung $du \wedge dt=-dt \wedge du$ mit dem Symbol für das Dachprodukt $\wedge$. Führt man den oben zitierten Euler'schen Beweis der Vertauschbarkeit durch, ist man jedoch gezwungen $dudt=dtdu$ zu verwenden, und die höheren Potenzen der Differentiale fortzulassen. Des antikommutativen Charakters von Differentialen wird Euler erstmalig  explizit bei seiner Behandlung von Doppelintegralen in \textit{``De formulis integralibus duplicatis"} (\cite{E391}, 1770, ges. 1768) (E391: ``Über Doppelintegrale") gewahr. Dort bemerkt er, dass man bei einer Substitution in beiden Integrationsvariablen, welche  $u(a,b)$ und $t(a,b)$ seien, nicht schlicht $dudt=(Ada+Bdb)(Cda+Ddb)$ mit entsprechenden Funktionen $A,B,C,D$ schreiben darf, sondern mit einem Faktor zu multiplizieren hat, der heute als Jacobi--Determinante bezeichnet wird. Diese findet Euler (siehe §§ 23--27 der erwähnten Arbeit)  für den Fall zweier Variablen, der allgemeine Fall von einer beliebigen Anzahl an Variablen wurde erst von Jacobi behandelt. Jacobi untersucht die nach ihm benannten Gebilde in den Arbeiten \textit{``De Determinantium functionalibus"} (\cite{Ja41a}, 1841) (``Über Funktionaldeterminanten") sowie \textit{``De formatione et proprietatibus Determinatium"} (\cite{Ja41b}, 1841) (``Über die Bildung und die Eigenschaften von Determinanten"). Für die Einordnung des Euler'schen Beitrags zur Theorie der Multiplikation von Differentialformen konsultiere man den Artikel \textit{``Euler and Differentias"} (\cite{Fe94}, 1994), wo die Geschichte des Differentials von Euler bis hin zu E. Cartan (1869--1951) nachgezeichnet wird.

\paragraph{Bei der Infinitesimalrechnung}
\label{para: Infinitesimalrechnung}

Euler verwendet in seinen Arbeiten häufig den Begriff von \textit{infinitesimalen Größen} oder spricht alternativ  von unendlich kleinen oder unendlich großen Größen, was durch das Fehlen des Grenzwertbegriffs notwendig war. Trotz der generellen Invalidität der Euler'schen Argumente zu Grenzwertprozessen haben   sie ihn  überzufällig häufig  zu richtigen Ergebnissen geleitet, sodass es aus seiner Sicht vermutlich keines neuen bzw. alternativen Begriffs bedurft hat. Seine \textit{Introductio} \cite{E101} präsentiert viele Resultate mittels Infinitesimalrechnung, welche sich heute nur mithilfe von Grenzwerten deuten lassen. Neben dem schon erwähnten Beispiel des Sinus--Produktes (\ref{eq: sine-Product}) (siehe Abschnitt \ref{subsubsec: Ermittlung des exakten Wertes}.), welches sich auch in seinem Lehrbuch findet, sei exemplarisch die Euler'sche Herleitung der Potenzreihe des Logarithmus aus dem binomischen Lehrsatz heraus demonstriert.  Euler geht nun vom Grenzwert 

\begin{equation}
\label{eq: Grenzwert Log}
      \log(1+x)= \lim_{n \rightarrow \infty} n((1+x)^{\frac{1}{n}}-1)
\end{equation}
aus, wobei er statt dem heute gebräuchlichen $\lim$  von unendlich großem $n$ spricht. Kompakt zusammengefasst begründet Euler wie folgt:

\begin{equation*}
    \lim_{n \rightarrow \infty} n((1+x)^{\frac{1}{n}}-1)=  \lim_{n \rightarrow \infty} \left(n \sum_{k=0}^{\infty} \binom{\frac{1}{n}}{k}x^k -1\right)= \sum_{k=1}^{\infty} \lim_{n \rightarrow \infty} n \binom{\frac{1}{n}}{k}x^k
\end{equation*}
\begin{equation*}
    = \sum_{k=1}^{\infty} (-1)^{k+1}\dfrac{x^k}{k},
\end{equation*}
was die Reihenentwicklung für $\log(1+x)$ ist.  \\

Natürlich ist die Vertauschung der Grenzwertprozesse gesondert zu rechtfertigen, was  in diesem Fall auch leicht zu leisten wäre.  Für Euler war es angesichts fehlender Gegenbeispiele nahezu unmöglich, die allgemeine Invalidität seines Vorgehens zu erkennen, zumal er aus anderer Quelle  -- etwa dem Satz von Taylor  oder Integration der Reihe für $\frac{1}{1+x}$ --  auf alternativem Weg\footnote{Ergänzend sei angemerkt, dass sich aus dem Grenzwert (\ref{eq: Grenzwert Log}) auch die Verbindung des Logarithmus zum Integral über die Hyperbel herstellen lässt.  Denn bekanntermaßen gilt für alle $s\neq 0$:

\begin{equation*}
    \int\limits_{0}^x (1+t)^{s-1}dt= \dfrac{(1+x)^{s}-1}{s}.
\end{equation*}
Betrachtet man nun den Grenzwert $s\rightarrow 0$, hat man für die rechte Seite:

\begin{equation*}
    \lim_{s\rightarrow 0} \dfrac{(1+x)^{s}-1}{s} = \lim_{s\rightarrow \infty } s((1+x)^s-1)=\log(1+x),
\end{equation*}
wobei  im letzten Schritt (\ref{eq: Grenzwert Log}) benutzt wurde. Für die linke Seite ergibt sich:

\begin{equation*}
     \lim_{s\rightarrow 0}  \int\limits_{0}^x (1+t)^{s-1}dt = \int\limits_{0}^x (1+t)^{-1}dt,
\end{equation*}
wo die Grenzwertvertauschung hier natürlich rechtens ist, weil das Integral rechter Hand existiert. Man gelangt  zur bekannten Tatsache, dass

\begin{equation*}
    \log(1+x)=  \int\limits_{0}^x (1+t)^{-1}dt
\end{equation*}
Diesen Zugang scheint Euler selbst allerdings nicht in seinen Werken erläutert zu haben.} zur Reihendarstellung des Logarithmus zu gelangen wusste\footnote{Dies illustriert ein weiteres Mal das von Euler präferierte Prinzip der konvergenten Evidenz.}. Selbiges gilt für die weiteren bekannten Funktionen wie $e^x, \sin x, \cos x$; sie lassen sich ebenfalls streng begründet nach dem Euler'schen Vorbild aus der Binomialreihe herleiten. Die Verwendung des Infinitesimalen scheint Euler überdies bei keiner Begebenheit zu einem falschen Resultat geführt zu haben.\\

\paragraph{Exakte bzw. alternative Konvergenzbegriffe} 
\label{para: Exakte Konvergenzbegriffe}

Die bisher erläuterten Beispiele stehen auch stellvertretend für das häufige Auftauchen von Grenzwertvertauschungen bei Euler, worunter auch die zuvor erwähnten die Frulliani'schen Integrale (\ref{eq: Log Frulliani}),  als Vertauschung von Ableitung und Integral (Abschnitt \ref{subsubsec: Die Darstellung betreffend -- Die hypergeometrische Reihe}), fallen. Dass dieser Prozess bei letztgenannten gestattet ist, zeigt man heute mit wenig Mühe. Ebenso verhält es sich bei der gliedweisen Integration und Differentiation von Potenzreihen innerhalb ihres Konvergenzradius, einer präferierten Methode Eulers. Weniger offenkundig ist dies bei seinen Grenzwertvertauschungen, welche er  bei trigonometrischen Reihen durchführt. Dazu betrachte man eines Beispiels wegen  § 40 von  \cite{E464}, in welcher er folgende Reihe untersucht:

\begin{equation*}
    \dfrac{\sin u}{1}+\dfrac{\sin 2u}{2}+\dfrac{\sin 3u}{3}+\dfrac{\sin 4u}{4}+\cdots =A-\dfrac{u}{2},
\end{equation*}
wo er die Konstante $A$ aus einem speziellen Fall zu bestimmen sucht. Aus dem Fall $u=0$ gelangt er wegen $\sin u=u$ für ``unendlich kleine $u$"{} zu der Gleichung

\begin{equation*}
    (u+u+u+\cdots)=A,
\end{equation*}
eine Gleichung der Gestalt $\infty \cdot 0=A$, welchen jeden beliebigen Wert annehmen kann. Daher betrachtet Euler nun den Wert $u=\pi$, genauer den Fall $u=\pi-\omega$, für unendlich kleines $\omega$, was diese Reihe zur Folge hat

\begin{equation*}
    -\omega +\omega -\omega +\omega -\omega+ \cdots. 
\end{equation*}
Diese  Summe ist für Euler endlich, er hätte sie wohl unter Anwendung seiner Definition von divergenten Reihen (siehe unten in Abschnitt \ref{subsubsec: Der Begriff der Summe einer Reihe}) zu $-\frac{1}{2}\omega$ berechnet. Da nun aber $\omega$ unendlich klein ist, ist diese Summe nun in der Tat $=0$. Dies erlaubt Euler nun $A=\frac{\pi}{2}$ zu bestimmen. Tatsächlich ist

\begin{equation}
\label{eq: Log Fourier}
    \dfrac{\sin u}{1}+\dfrac{\sin 2u}{2}+\dfrac{\sin 3u}{3}+\dfrac{\sin 4u}{4}+\cdots =\dfrac{\pi-u}{2}
\end{equation}
 für $u \in \left(-\pi,\pi\right)$ richtig. Gegen das Euler'sche Argument lässt sich der Einwand erheben, dass man statt $u=\pi-\omega$ auch den Wert $u=\pi+\omega$ einsetzen könnte, wodurch man wieder auf den Fall $0\cdot \infty$ geführt werden würde, also nichts gewonnen hätte. Der Widerspruch lässt sich erst durch den eigens Fourierreihen zukommenden Konvergenzbegriff auflösen. Euler selbst hat vermutlich durch folgende Betrachtung Bestätigung für die Richtigkeit der letzten Gleichung erhalten: Die linke Seite ergibt sich gerade als imaginärer Teil der Reihe

\begin{equation*}
    -\log (1-e^{iu})= \dfrac{e^{iu}}{1}+ \dfrac{e^{2iu}}{2}+ \dfrac{e^{3iu}}{3}+\cdots. 
\end{equation*}
Die Berechnung von Imaginärteilen von Logarithmen hat Euler an mehreren Stellen auseinander gesetzt, sie reduziert sich auf Gleichung (\ref{eq: log complex}). Nun hat man jedoch für $ 1-e^{iu}=1-\cos (u)-i \sin (u)$ nach Eulers Formel für das Argument einer komplexen Zahl:

\begin{equation*}
  -\arcsin \left( \dfrac{-\sin (u)}{\sqrt{(1-\cos(u))^2+\sin^2(u)}}\right)=\arcsin \left( \dfrac{\sin (u)}{\sqrt{2-2\cos(u)}}\right),
\end{equation*}
was mit den Identitäten $\sin(u)=2 \sin \left(\frac{u}{2}\right)\cos \left(\frac{u}{2}\right)$ und $\sin \left(\frac{x}{2}\right)= \sqrt{\frac{1-\cos(x)}{2}}$ zu

\begin{equation*}
    \arcsin \left(\cos \left(\frac{u}{2}\right)\right)= \arcsin \left(\sin\left(\dfrac{\pi}{2}-\dfrac{u}{2}\right)\right)= \dfrac{\pi-u}{2}
\end{equation*}
wird und somit (\ref{eq: Log Fourier}) bestätigt. Diese Reihe stellt wieder ein Beispiel für ein mehrfach von Euler nachgewiesenes Ergebnis dar, so findet man sie -- wie auch in Abschnitt (\ref{subsubsec: Durch Kombinieren von Ergebnissen: Die zeta-Funktion}) erwähnt -- in der Arbeit \cite{E555} aus der Lagrange--Interpolation heraus abgeleitet. In § 92 des zweiten Teils seiner \textit{Calculi Differentialis} \cite{E212} stellt Euler eine Herleitung aus der Potenzreihe des Arkustangens vor. Auf dem just vorgestellten Weg gelangt Euler in seiner Arbeit \textit{``Summatio progressionum $\sin^{\lambda}(\phi)+\sin^{\lambda}(2\phi)+\sin^{\lambda}(3\phi)+\cdots+\sin^{\lambda}(n\phi)$, $\cos^{\lambda}(\phi)+\cos^{\lambda}(2\phi)+\cos^{\lambda}(3\phi)+\cdots+\cos^{\lambda}(n\phi)$"} (\cite{E447}, 1774, ges. 1773) (E447: ``Summation der Progressionen $\sin^{\lambda}(\phi)+\sin^{\lambda}(2\phi)+\sin^{\lambda}(3\phi)+\cdots+\sin^{\lambda}(n\phi)$, $\cos^{\lambda}(\phi)+\cos^{\lambda}(2\phi)+\cos^{\lambda}(3\phi)+\cdots+\cos^{\lambda}(n\phi)$") zu dieser Reihe, wo sie sich im letzten Beispiel (Beispiel 2) findet.

\subsubsection{Der Begriff der Summe einer Reihe}
\label{subsubsec: Der Begriff der Summe einer Reihe}

\epigraph{Discovery consists of looking at the same thing as everyone else and thinking something different.}{Albert Szent--Gjorgyi}


Wie in der Besprechung des Euler'schen Beitrags zur Riemann'schen $\zeta$--Funktion (Abschnitt  \ref{subsubsec: Durch Kombinieren von Ergebnissen: Die zeta-Funktion}) schon angedeutet worden ist, ist Euler über divergente Reihen zur  Funktionalgleichung (\ref{eq: Functional Equation Eta}) für die Riemann'sche $\zeta$--Funktion gelangt. Dies ist zweifelsohne eine von Eulers erstaunlichsten Entdeckungen, was Abel indes nicht abhielt, in einem Brief  (Brief XX. aus seinen gesammelten Werken \cite{Ab12b}) zu schreiben:\\

\textit{``Divergente Reihen sind im Allgemeinen etwas sehr Verhängnisvolles, und es ist eine Schande es zu wagen, auf sie irgendeinen Beweis zu stützen. Man kann alles Mögliche mit ihrer Hilfe beweisen. Und sie haben so viel Unheil verursacht, so viele Paradoxa hervorgebracht."}\\

Trotz der (berechtigten) Abel'schen Kritik hat sich  Euler an verschiedenen Stellen mit divergenten Reihen beschäftigt und bemerkenswerte Ergebnisse zutage gefördert. Es soll seine Definition der divergenten Reihe selbst vorgestellt werden, anschließend die Anwendung auf die $\zeta$--Funktion erläutert werden, aber auch die Widersprüche aus der Definition sollen nicht unerwähnt bleiben.

\paragraph{Eulers Definition von divergenten Reihen}
\label{para: Eulers Definition  von divergenten Reihen}

Euler stellt seine Theorie der divergenten Reihen ausführlich in seiner Abhandlung \textit{``De seriebus divergentibus"} (\cite{E247}, 1760, ges. 1746) (E247: ``Über divergente Reihen") vor, aus welcher zunächst die Fundamente seiner Auffassung  zitiert werden sollen. Direkt in § 1 schreibt er:\\

\textit{``Während konvergente Reihen so definiert werden, dass sie aus stetig schrumpfenden Termen bestehen, die schließlich, wenn die Reihe ins Unendliche fortgesetzt wird, verschwinden}\footnote{Dass dieser Satz nicht gänzlich wörtlich zu verstehen ist, wird Euler bewusst gewesen sein, zumal er in der vorher verfassten Arbeit \cite{E43} Reihen wie 

\begin{equation*}
    \log(3)=1+\frac{1}{2}-\frac{2}{3}+\frac{1}{4}+\frac{1}{5}-\frac{2}{6}+\frac{1}{7}+\frac{1}{8}-\frac{1}{9}+\frac{1}{10}+\frac{1}{11}-\frac{2}{12}+\text{etc.}
\end{equation*}
betrachtet, welche zwar, absolute Werte betrachtend, aus im Unendlichen verschwindenden, aber nicht stetig kleiner werdenden Termen besteht. Diese Reihe und ähnliche andere studiert Euler in § 6 der erwähnten Arbeit.},\textit{ sieht man leicht ein, dass die Reihen, deren infinitesimale Terme  nicht verschwinden,  sondern entweder endlich bleiben oder ins Unendliche anwachsen, weil sie nicht konvergent sind, zur Klasse der divergenten Reihen gerechnet werden müssen. Je nachdem, ob die letzten Terme der Reihe, zu welchem beim Fortsetzen der Progression bis ins Unendliche gelangt wird, von endlicher oder unendlicher Größe waren, wird man zwei Gattungen von divergenten Reihen haben, welche beide weiter in zwei Klassen unterteilt werden, je nachdem ob alle Terme mit demselben Vorzeichen behaftet sind oder die Vorzeichen $+$ und $-$ alternieren. Insgesamt werden wir also vier Klassen von divergenten Reihen haben, aus welchen ich zwecks der Verständlichkeit einige Beispiele beifügen möchte."}\\

Die von Euler angedeuteten Exempel sind in nachstehender Tabelle zu sehen:

\begin{equation*}
    \renewcommand{\arraystretch}{2,0}
\begin{array}{llllllllllllllllll}
     \text{I.}) \quad &  1 & + & 1 &+ & 1 &+ &1 &+ & 1 &+ & 1 &+ &\text{etc.} \\
     \quad &  \frac{1}{2} & + & \frac{2}{3} & + & \frac{3}{4} &+ & \frac{4}{5} &+ & \frac{5}{6} &+ & \frac{6}{7}  &+ &\text{etc.} \\ \hline
     \text{II.}) \quad &  1 & - & 1 &+ & 1 & - &1 &+ & 1 & - & 1 &+ &\text{etc.} \\
     \quad &  \frac{1}{2} & - & \frac{2}{3} & + & \frac{3}{4} & - & \frac{4}{5} &+ & \frac{5}{6} & - & \frac{6}{7}  &+ &\text{etc.} \\ \hline
     \text{III.}) \quad &  1 & + & 2 &+ & 3 &+ & 4 &+ & 5 &+ & 6 &+ &\text{etc.} \\
     \quad &  1 & + & 2 & + & 4 &+ & 8 &+ & 16 &+ & 32  &+ &\text{etc.} \\ \hline
     \text{IV.}) \quad &  1 & + & 2 &+ & 3 &+ & 4 &+ & 5 &+ & 6 &+ &\text{etc.} \\
     \quad &  1 & - & 2 & + & 4 &- & 8 &+ & 16 &- & 32  &+ &\text{etc.} \\ \hline
\end{array}
\end{equation*}

Es wird förderlich sein, die grundlegenden Gedanken Eulers zu jeder ersten drei Klassen mitgeteilt zu haben, bevor die vierte gesondert betrachtet wird. Zur ersten Klasse schreibt er in § 2:\\

\textit{``Und zuerst ist freilich ersichtlich, dass die Summen der Reihen, welche ich der ersten Klasse zugeordnet habe, in Wirklichkeit unendlich groß sind}\footnote{Dies ist stark abhängig von der Summationsmethode, welcher man sich bedient, so ließe sich etwa das erste Beispiel von Euler mithilfe der Riemann'schen $\zeta$--Funktion als $\zeta(0)=-\frac{1}{2}$ angeben. Auf die Widersprüche in der Euler'schen Definition von Summen divergenten Reihen wird auch weiter unten in diesem Abschnitt noch eingegangen werden.}, \textit{weil durch tatsächliches Sammeln  der Terme zu einer Summe größer als jede gegebene Zahl gelangt wird. Daher besteht freilich kein Zweifel, dass die Summen dieser Reihen durch Ausdrücke der Gestalt $\frac{a}{0}$ dargeboten werden können."}\\

Seine Ausführungen zur zweiten Klasse findet man in § 3. Insbesondere das von Leibniz betrachtete Beispiel ist Kern seiner Diskussion. Euler schreibt:\\

\textit{``Aus der zweiten Klasse hat als erster Leibniz diese Reihe betrachtet}

\begin{equation*}
    1-1+1-1+1-1+1-1+\text{etc.},
\end{equation*}
\textit{deren Summe er  festgelegt hatte $=\frac{1}{2}$ zu sein, wobei er sich auf diese hinreichend soliden Gründe berief: Zuerst geht diese Reihe hervor, wenn dieser Bruch $\frac{1}{1+a} $ durch unendlich viele Male wiederholte Division auf gewohnte Weise in diese Reihe $1-a+a^2-a^3+a^4+\text{etc.}$ aufgelöst wird, und anschließend der Buchstabe $a$ der Einheit gleich gewählt wird. Dann [...] gibt er folgendes Argument: Wenn die Reihe irgendwo abgebrochen wird, und die Anzahl der Terme gerade war, wird ihr Wert $=0$ sein, wenn aber die Anzahl der Terme ungerade war, wird der Wert der Reihe $=1$ sein. Wenn also die Reihe ins Unendliche fortschreitet und die Anzahl der Terme weder als gerade noch als ungerade angesehen werden kann, kann er daraus schließen, dass die Summe weder $=0$ noch $=1$ ist,  sondern einen gewissen Zwischenwert, welcher von den beiden gleich weit entfernt ist, annehmen muss, welcher $=\frac{1}{2}$ ist."}\\

Bereits in diesen Erläuterungen kristallisiert sich die Euler'sche Auffassung zum Umgang mit divergenten Reihen heraus, welche in seiner Besprechung der Reihen der dritten Klasse noch intensiver durchscheint. Er  schreibt in § 6:\\

\textit{``Obwohl nämlich die Terme dieser Reihen [der dritten Klasse] ununterbrochen wachsen, und daher durch tatsächliches Sammeln der Terme zu einer Summe größer als jede angebbare Zahl gelangt werden kann, was die Definition des Unendlichen ist, sind dennoch die Fürsprecher von Reihen bei dieser Gattung gezwungen Reihen von dieser Art zuzulassen, deren Summen endlich sind, freilich sogar negativ oder kleiner als $0$. Denn weil der Bruch $\frac{1}{1-a}$ durch Division entwickelt $1+a+a^2+a^3+a^4+\text{etc.}$ gibt, müsste auch gelten:}

\begin{equation*}
     \renewcommand{\arraystretch}{1,5}
\begin{array}{rclllllllllllll}
     -1 & = & 1 & + & 2 & + & 4 & + & 8 & + &16 &+ & \text{etc.} \\ 
      -\frac{1}{2} & = & 1 & + & 3 & + & 9 & + & 27 & + & 81 &+ & \text{etc.}"
\end{array}
\end{equation*}

Nach einigen Erläuterungen allgemeinerer Natur präsentiert Euler in § 12 seine Definition der Summe einer divergenten Reihe. Hier schreibt er:\\

\textit{``Wenn wir also die übliche Auffassung einer Summe lediglich so abwandeln, dass wir sagen, dass die Summe einer Reihe der endliche Ausdruck ist, aus dessen Entwicklung jene Reihe selbst entsteht, werden alle Schwierigkeiten [...] verschwinden."}\\

Da  heute verschiedene Methoden zur Summation von divergenten Reihen bekannt sind, welche unterschiedliche Werte für ein und dieselbe divergente Reihe geben können, ist Eulers Aussage in der Allgemeinheit zwar zurückzuweisen, was für Euler selbst auch zu erkennen gewesen wäre, wie weiter unten in diesem Abschnitt deutlich werden wird. In seiner Arbeit \cite{E247} behandelt er jedoch unabhängig davon aus seiner vierten Klasse von Reihen das  berühmte Beispiel der Fakultätenreihe

\begin{equation}
\label{eq: Fakulatätenreihe}
    s(1)=0!-1!+2!-3!+4!-\text{etc.}
\end{equation}
mit

\begin{equation}
\label{eq: PowerFak}
    s(x)=0!-1!x+2!x^2-3!x^3+4!x^4-\text{etc.},
\end{equation}
 zur Summierung welcher er 4 verschiedene Methoden vorstellt, die ihn alle zum selben Wert leiten\footnote{Die erste bedient sich der Differenzenrechnung (§ 13), die zweite führt die Summierung auf ein Interpolationsproblem zurück (§ 17), in § 19 wird dann die zur obigen Reihe assoziierte Potenzreihe zunächst auf eine Differentialgleichung reduziert, welche anschließend gelöst wird. Schließlich wird dieselbe Potenzreihe auch in einen Kettenbruch umgewandelt, welcher konvergiert. Besagter Kettenbruch ist ein Spezialfall desjenigen aus  \cite{E616}.}. Er findet mit all seinen Methoden annähernd denselben numerischen Wert, sodass das Prinzip der konvergenten Evidenz ihn (in diesem Fall berechtigt) von der Richtigkeit nachstehender Gleichheit überzeugt:

 \begin{equation*}
    0!-1!+2!-3!+4!-\text{etc.} \approx  0.5963473621372
 \end{equation*}

 \begin{figure}
     \centering
    \includegraphics[scale=0.7]{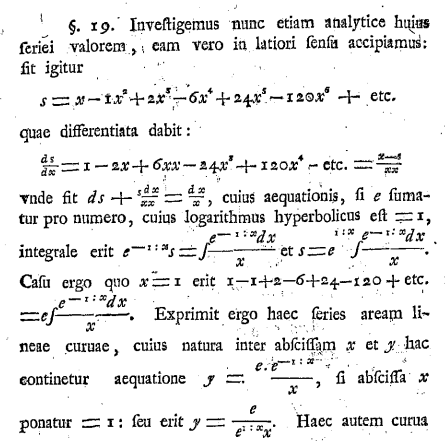}
     \caption{Euler leitet in seiner Arbeit \cite{E247} eine Differentialgleichung (und ihre Lösung) zwecks Summierung der divergenten Fakultätenreihe her.}
     \label{fig:E247Potenzreihe}
 \end{figure}

All diese Methoden sind genauer in de Artikeln \textit{``Euler Subdues a Very Obsteperous Series"} (\cite{Ba07}, 2007) und \textit{``Divergent series"} (\cite{Sa06jun}, 2006) aus Sandifers Kolumne vorgestellt. Allgemein werden die Methoden in den Büchern \textit{``Divergent Series"} (\cite{Ha49}, 1949) und  \textit{`` Euler Through Time: A New Look at Old Themes"} (\cite{Va06}, 2006) besprochen, das Buch von Balser \textit{``From Divergent Power Series to Analytic Functions: Theory and Application of Multisummable Power Series"} (\cite{Ba09}, 2009) nimmt den Ausgangspunkt bei dem Ansatz über Potenzreihen -- dem dritten Euler'schen Ansatz zur Summation der Fakultätenreihe.  Aus moderner Sicht ist die Äquivalenz seiner Methoden in diesem Fall dem Alternieren der Reihe zuzuschreiben\footnote{In der Tat hat Euler oft allein alternierende Reihen betrachtet -- die Arbeit \textit{``De summatione serierum, in quibus terminorum signa alternantur"} (\cite{E617}, 1788, ges. 1776) (E617: ``Über die Summation der Reihen, in welchen die Vorzeichen der Terme alternieren") ist ein Zeugnis dessen. Die Gründe sind wohl in der leichteren Handhabbarkeit der teils auch numerischen Rechnungen zu suchen, welche sich über Ideen aus der komplexen Analysis auch vom mathematischen Standpunkt aus  nachvollziehen lässt. Man vergleiche dazu auch die Ausführungen in \cite{Ha49} für den Nachweis der Äquivalenz der Methoden.}. An dieser Stelle sei lediglich kleine Ergänzung zu seiner Methode über die Interpolation  eingebracht.  In § 17 bemerkt Euler, dass sich die gesuchte Reihe (\ref{eq: Fakulatätenreihe}) aus der Funktion

\begin{equation*}
    P(n):= 1+(n-1)+(n-1)(n-2)+(n-1)(n-2)(n-3)+\text{etc.}
\end{equation*}
für den speziellen Wert $n=0$ ergibt. Euler findet diesen numerischen Zahlenwert, welchen er $A$ nennt, über Differenzenrechnung am Ende desselben Paragraphen. Jedoch bemerkt er bei der Konstruktion seiner Tabelle die Gleichung

\begin{equation*}
    P(n+1)=nP(n)+1,
\end{equation*}
welche man alternativ als homogene Differenzengleichung mit linearen Koeffizienten darstellen könnte. Es gilt gleichermaßen:

\begin{equation*}
    P(n)=nP(n-1)-(n-2)P(n-2).
\end{equation*}
Dies erlaubt nun die Applikation der oben (Abschnitt \ref{subsubsec: Die Mellin--Transformierte bei Euler}) vorgestellten Methode, welche zur Gleichung

\begin{equation*}
    P(n)=e \cdot \int\limits_{1}^{\infty}e^{-t}t^{n-1}dt
\end{equation*}
führt. Demnach ergibt sich für den Wert von (\ref{eq: Fakulatätenreihe}) auf diesem Wege

\begin{equation*}
    P(0)= e\int\limits_{1}^{\infty} \dfrac{e^{-t}}{t}dt, 
\end{equation*}
welches Integral numerisch ausgewertet ebenfalls den von Euler gefundenen Wert gibt. Euler gelangt überdies vermöge seines Ansatzes über die Lösung der Differentialgleichung, welche von (\ref{eq: PowerFak}) erfüllt wird, in § 19 zu dem gleichwertigen Integral

\begin{equation*}
   e\cdot \int\limits_{0}^1 \dfrac{e^{-\frac{1}{u}}}{u}du
\end{equation*}
für die Fakultätenreihe (\ref{eq: Fakulatätenreihe}). Es geht aus dem vorherigen durch die Substitution $x=\frac{1}{u}$ hervor. In ähnlicher Manier findet man zur Reihe

\begin{equation*}
    z=1-1+1\cdot 3 -1\cdot 3 \cdot 5 +1\cdot 3 \cdot 5 \cdot 7- \cdots
\end{equation*}
über die assoziierte Differenzengleichung

\begin{equation*}
    P_z(n+2)=(n+1)P_z(n)+1
\end{equation*}
die Summe

\begin{equation*}
    z=P_z(0)= \sqrt{\dfrac{e\pi}{2}}\operatorname{erf}\left(\dfrac{1}{\sqrt{2}}\right) \approx 0,65568 \quad \text{mit} \quad \operatorname{erf}(x):= \int\limits_x^{\infty} e^{-\frac{t^2}{2}}dt.
\end{equation*}
oder gar die Summe dieser allgemeineren Reihe

\begin{equation*}
    f(p,q):= 1-p+p(p+q)-p(p+q)(p+2q)+p(p+q)(p+2q)(p+3q)-\cdots
\end{equation*}
aus der assoziierten Differenzengleichung

\begin{equation*}
    P_f(n+q)=(n-p+q)P_f(n)+1
\end{equation*}
als

\begin{equation*}
    P_f(0)= \sqrt[q]{e} \int\limits_{1}^{\infty} u^{q-p-1}e^{-\frac{u^q}{q}}du,
\end{equation*}
welche sich entsprechend auch bei Euler finden\footnote{Bei Euler steht in § 29
\begin{equation*}
    1x-1x^3+1\cdot 3x^5 -1\cdot 3 \cdot 5x^7 +1\cdot 3 \cdot 5 \cdot 7x^9- \cdots = e^{\frac{1}{2x^2}}\int\limits_{0}^{1} \dfrac{e^{-\frac{-1}{2x^2}}}{xx}dx
\end{equation*}
und bereits vorher in § 28 von \cite{E247}:

\begin{equation*}
    x^m-px^{m+q}+p(p+q)x^{m+2q}-p(p+q)(p+2q)x^{m+3q}+\cdots= e^{-\frac{1}{x^q}}x^{m-p} \int\limits_{0}^{1} e^{-\frac{1}{x^q}}x^{p-q-1}dx.
\end{equation*}
}.\\

Alternativ mag man die Potenzreihe (\ref{eq: PowerFak}) auch wie folgt umformen: Man schreibe

\begin{equation*}
    s(x)= \sum_{n=0}^{\infty} (-1)^{n}n!x^n.
\end{equation*}
Mit der $\Gamma$--Funktion findet man auch:

\begin{equation*}
     s(x)= \sum_{n=0}^{\infty} (-1)^{n}\cdot \left(\int\limits_{0}^{\infty}e^{-t}t^ndt \right)  x^ndt.
\end{equation*}
Eine Vertauschung von Summe und Integral und anschließende Verwendung der geometrischen Reihe gibt demnach:

\begin{equation*}
     s(x) = \int\limits_{0}^{\infty}e^{-t} \sum_{n=0}^{\infty} (-xt)^ndt = \int\limits_{0}^{\infty} \dfrac{e^{-t}}{1+xt}dt,
\end{equation*}
Integrale welcher Gestalt für Borel (1871--1956) in seinen \textit{``Le\c{c}ons sur les séries divergentes"} (\cite{Bo01}, 1901) (``Vorlesungen über divergente Reihen") den Ausgangspunkt seiner Untersuchungen gebildet haben.

\paragraph{Widersprüche aus der Euler'schen Definition}
\label{para: Widersprüche aus der Euler'schen Definition}

Wie erwähnt, führt Eulers Definition zu nicht auflösbaren Widersprüchen. Denn sie impliziert insbesondere die (in Wirklichkeit nicht vorhandene) Eindeutigkeit der Summe einer divergenten Reihe, welche Auffassung  Eulers vermutlich aus der Betrachtung von Potenzreihen herrührt\footnote{In diesem Kontext lässt sich anmerken, dass Euler selbst, wenn auch lediglich implizit, in seiner Arbeit \cite{E47} über die Euler--Maclaurin'sche Summenformel eine Summationsmethode mitgeteilt hat, welche von denen in der gerade besprochenen \cite{E247} verschieden ist. Diese hat er, entsprechend interpretiert, unter anderem zur Summierung der harmonischen Reihe $1+\frac{1}{2}+\frac{1}{3}+\frac{1}{4}+\cdots$ angewandt und ihren Wert als $\gamma$, die Euler--Mascheroni--Konstante, berechnet. Ramanujan hingegen hat dieselbe Methode explizit benutzt und divergente Reihen mithilfe der Euler--Maclaurin'schen Summenformel evaluiert. Man vergleiche dazu insbesondere Kapitel 6 des Buches \textit{``Ramanujan’s Notebooks: Part I"} (\cite{Be85}, 1985).}. Gleichwohl die Anwendung seiner Ideen gar in der Funktionalgleichung der Riemann'schen $\zeta$--Funktion bzw. genauer der $\eta$--Funktion (\ref{eq: Functional Equation Eta}) mündeten, erzeugen sie am Beispiel dieser beiden Funktionen bereits Paradoxa. Beide Funktionen sind in elementarer Weise durch die Formel 

\begin{equation}
\label{eq: zeta-eta elementar}
    1-\dfrac{1}{2^n}+\dfrac{1}{3^n}-\cdots =\left(1-\dfrac{1}{2^{n-1}}\right)\left( 1+\dfrac{1}{2^n}+\dfrac{1}{3^n}-\cdots\right)
\end{equation}
verbunden, welche Gleichung formal für beliebiges $n$ richtig ist\footnote{Euler gibt sie in § 170 seiner \textit{Introductio} an.}, auch wenn nur für $\operatorname{Re}(n)>1$ auf beiden Seiten zugleich konvergente Reihen stehen. Wie oben in Abschnitt (\ref{subsubsec: Durch Kombinieren von Ergebnissen: Die zeta-Funktion}) angedeutet, kann Euler die divergenten Reihen linker Hand evaluieren, sofern $n$ negative ganze Zahlen sind. In \cite{E352} nennt er neben anderen explizit

\begin{equation*}
    1-1+1-1-1+\cdots = \dfrac{1}{2}, \quad 1-2+3-4+\cdots =\dfrac{1}{4}, \quad 1-2^2+3^2-4^2+\cdots=0 \quad \text{etc.}
\end{equation*}
Dabei bedient er sich entsprechender Potenzreihen wie etwa

\begin{equation*}
    \dfrac{1}{1+x}=1-x+x^2-x^3+\cdots, \quad \dfrac{1}{(1+x)^2}=1-2x+3x^2-\text{etc.}, 
\end{equation*}
\begin{equation*}
    \dfrac{1-x}{(1+x)^3}=1-2^2x+3^2x^2-\text{etc.}, 
\end{equation*}
aus welchen für $x=1$ die obigen Reihen folgen\footnote{Diesen Vorgang der Zuschreibung einer Summe zu einer divergenten Reihe nennt man Abel'sche Summation. Der Name rührt von der Analogie zum Abel'schen Grenzwertsatz her, welcher die hier Grenzwertvertauschung unter entsprechenden Voraussetzungen gestattet. Abel hat den nach ihm benannten Satz in seiner Arbeit \textit{``Untersuchungen über die Reihe $1+mx+\frac{m\cdot (m-1)}{2\cdot 1}x^2+\frac{m\cdot (m-1) \cdot (m-2)}{3 \cdot 2\cdot 1}x^3+\cdots$."} (\cite{Ab26}, 1826) bewiesen.}. Weiter findet man aus (\ref{eq: zeta-eta elementar})  formal:

\begin{equation*}
    1-1+1-1+\text{etc.}= \left(1-2\right)\left(1+1+1+1+\text{etc.}\right),
\end{equation*}
woraus man mit dem bekannten Wert für die linke Seite auch

\begin{equation*}
    1+1+1+1+\text{etc.}=-\dfrac{1}{2}
\end{equation*}
erhält, was dem Wert $\zeta(0)$ entspricht. Alternativ hätte Euler diesen Wert aus der formalen Identität

\begin{equation*}
    \frac{1}{2}= \sum_{n=0}^{\infty} \cos(nx)
\end{equation*}
herleiten können, welche er etwa in § 11  seiner Arbeit \cite{E447} angibt\footnote{In dieser Arbeit begründet er die Richtigkeit dieser Formel ebenfalls mit seiner just vorgestellten Konzeption einer Reihe.}. Man hat lediglich $x=0$ zu setzen, Euler macht jedoch keine Erwähnung davon\footnote{Gleichermaßen ergeben sich aus der formal gültigen Gleichheit $  \frac{1}{2}= \sum_{n=0}^{\infty} \cos(nx)$ durch wiederholte Differentiation unmittelbar die Werte $\zeta(-2k)$ für alle natürlichen $k$ unmittelbar als $=0$.}. \\

Entsprechend fände man auch 

\begin{equation*}
    1+2+3+4+\text{etc.}=-\dfrac{1}{12}, \quad  1+2^2+3^2+4^2+\text{etc.}=0 \quad \text{etc.}
\end{equation*}
und so weiter. \\

Will man jedoch direkt mit der Euler'schen Idee beginnen und etwa die Summe

\begin{equation*}
    1+2+3+4+5+\cdots 
\end{equation*}
berechnen, würde man zunächst die Reihe 

\begin{equation*}
    \dfrac{1}{1-x}=1+x+x^2-\text{etc.}
\end{equation*}
verwenden wollen. Setzt man hier allerdings $x=1$, so gelangt man zu der Aussage

\begin{equation*}
     1+2+3+4+5+\cdots= \infty,
\end{equation*}
was dem zur mitgeteilten Ergebnis $1+2+3+4+\cdots= \zeta(-1)=-\frac{1}{12}$ widerspricht. Eulers Vorgaben leiten demnach zu konträren Ergebnissen\footnote{Diese Reihen sind in der oben vorgestellten Taxonomie von divergenten Reihen der ersten Klasse, welche für Euler tatsächlich den Wert $\infty$ haben.}. \\

Nichtsdestotrotz stellt  Eulers Definition von divergenten Reihen ein sehr illustratives Beispiel für den Prozess der mathematischen Innovation dar. Zum einen sieht er sich durch viele Beispiele und sorgfältige Diskussion gleichsam gezwungen, diesen Begriff auf dem beschriebenen Wege zu definieren. Er hat ihn also gewissermaßen durch  Induktion geleitet eingeführt, dann jedoch weiter  streng deduktiv aus ihm heraus die Funktionalgleichung (\ref{eq: Functional Equation Eta}) abgeleitet, muss aber -- obschon er es in seinen Arbeiten nicht explizit erwähnt zu haben scheint  -- gleichzeitig der Inkonsistenz seiner Auffassung von divergenten Reihen gewahr gewesen sein. Das Euler'sche Beispiel der Behandlung divergenter Reihen kann  exemplarisch für die Illustration wissenschaftlichen Arbeitens herangezogen werden: Euler geht bis an die Grenze des für ihn Leistbaren und scheut nicht, auch neue Gebiete zu eröffnen, welches dann auch seine Nachfolger mit glücklichem Erfolg anzuwenden wussten. Es seien konkret Fourier und Poincaré (1854--1912) als Beispiele angeführt. Ersterer hat in seinem Hauptwerk \textit{``Théorie analytique de la chaleur"} (\cite{Fo22}, 1822) (``Die analytische Theorie der Wärme") Teile seiner Resultate über das Divergente erlangt\footnote{Man vergleiche auch die entsprechenden Ausführungen aus \cite{Ha49}.}. Die Rolle der divergenten Reihen in Poincarés Arbeit \textit{``Sur le problème des trois corps et les équations de la dynamique"} (\cite{Po90}, 1890) (``Über das Dreikörperproblem und die Gleichungen des Dynamik") zum Dreikörperproblem beim Nachweis seiner generellen Unlösbarkeit ist hinlänglich bekannt. Schließen soll diesen Paragraphen  ein Zitat von Andreas Speiser aus Vorwort zu Band 9 von Serie 1 der \textit{Opera Omnia} (\cite{OO9}, 1945), wo er auf den Seiten IX und X zu Eulers Handhabung der divergenten Reihen schreibt:\\

\textit{``Euler hat in seinen unermeßlichen Expeditionen im Reiche der Reihen die überwältigende Entdeckung gemacht, dass gerade die divergenten Reihen das kräftigste Mittel zur Auffindung unerwarteter Tatsachen bilden. Der Durchgang durch das Divergente ist  noch ertragreicher als der Durchgang durch das Komplexe in der Funktionentheorie. So fand er auf diesem Weg die Produktentwicklung der Gammafunktion, ferner die Funktionalgleichung derselben und der Zetafunktion. Letztere Entdeckung, die zum Merkwürdigsten  gehört, was in der Mathematik gefunden wurde, konnte er nicht beweisen."}\\

Der letzte Satz von Speiser ist indes durch das oben (Abschnitt \ref{subsubsec: Durch Kombinieren von Ergebnissen: Die zeta-Funktion}) Gesagte mit entsprechenden Erläuterungen bezüglich des Wortes ``konnte"{} zu ergänzen.

\subsubsection{Mehrwertige Funktionen: Das fehlende Konzept der Riemann'schen Fläche}
\label{subsubsec: Mehrwertige Funktionen: Das fehlende Konzept der Riemann'schen Fläche}

\epigraph{It is tempting, if the only tool you have is a hammer, to treat everything as if it were a nail.}{Abraham Maslow}

Die Euler'sche Meinung,  Paradoxa werden sich durch Einführung eines synthesestiftenden übergeordneten Begriffs auflösen, ist oben (siehe den entsprechenden Paragraphen in Abschnitt \ref{para: Dialektischer Ansatz}) am Beispiel des komplexen Logarithmus illustriert worden. Die Auflösung aller Schwierigkeiten bestand dort darin, den Logarithmus im Bereich der komplexen Zahlen als unendlichwertige Funktion zu erkennen. Dass jedoch mehrwertige Funktionen selbst  eines neuen widersprüchebeseitigenden Begriffs bedürfen, scheint bei Euler keine Erwähnung gefunden zu haben. Das Fehlen des Konzeptes des Riemann'schen Fläche, im Rahmen welcher verzweigte komplexe Funktionen heute behandelt werden, führte Euler bei verschiedenen Begebenheiten zu falschen Aussagen. Dies soll hier aufgezeigt werden. Überdies soll plausibel gemacht werden, dass es für Euler schwierig gewesen sein muss, überhaupt die Frage nach einem neuen Konzept für die mehrwertigen Funktionen aufkommen zu lassen.

\paragraph{Zum Logarithmus}

Es ist dem Verständnis zuträglich, die Erläuterungen von  Eulers Theorie des Logarithmus aus zu beginnen.  In den beiden an entsprechender Stelle erwähnten Abhandlungen \cite{E168} und \cite{E170} nimmt Euler implizit die Gültigkeit der Funktionalgleichung $\log  (x\cdot y)= \log  (x)+\log (y)$ für alle komplexen Zahlen an -- ein Irrtum, wie heute bekannt ist\footnote{Die Funktionalgleichungen behält nur ihre Gültigkeit, sofern man bei der Addition der beiden Logarithmen nicht das entsprechende Blatt der Riemann'schen Fläche wechseln muss.}. Für den Beweis aus der ersten Abhandlung (dort findet sich die Annahme in Problem 1) ist das auch zwingend nötig, der Beweis über die Differentialformen aus der zweiten benötigt diese Annahme hingegen nicht. An dieser Stelle betritt Euler allerdings das Gebiet der inexakten Differentiale, woraus man die Mehrdeutigkeit des $\log $ ebenfalls ableiten kann\footnote{In der Funktionentheorie hingegen ist das Integral

\begin{equation*}
    f(z)=\frac{1}{z}
\end{equation*}
bekannt, keine Stammfunktion zu haben, was auch in der Vieldeutigkeit des $\log $ begründet liegt.}.\\

Alternativ ließe sich für Euler aus dem Grenzwert (\ref{eq: Grenzwert ln})  direkt die Gültigkeit der Funktionalgleichung

\begin{equation*}
    \log  (x \cdot y)= \log (x)+\log (y),
\end{equation*}
nachweisen, zumal 

\begin{equation*}
    \log  (x \cdot y) = \lim_{n\rightarrow \infty} n \left((xy)^{\frac{1}{n}}-1\right)= \lim_{n\rightarrow \infty} n \left((xy)^{\frac{1}{n}}-x^{\frac{1}{n}}+x^{\frac{1}{n}}-1\right).
\end{equation*}
Eine Aufteilung in zwei Grenzwerte gibt

\begin{equation*}
    \renewcommand{\arraystretch}{2,0}
\begin{array}{rcL}
    \log  (x \cdot y)  & = &  \lim_{n\rightarrow \infty} n \left((xy)^{\frac{1}{n}}-x^{\frac{1}{n}}\right)+\lim_{n\rightarrow \infty} n \left(x^{\frac{1}{n}}-1\right)\\
     & = & \lim_{n\rightarrow \infty} x^{\frac{1}{n}} \cdot n \left(y^{\frac{1}{n}}-1\right)+\lim_{n\rightarrow \infty} n \left(x^{\frac{1}{n}}-1\right) \\
     & = & \lim_{n\rightarrow \infty} x^{\frac{1}{n}}\cdot\lim_{n\rightarrow \infty}   n \left(y^{\frac{1}{n}}-1\right)+\lim_{n\rightarrow \infty} n \left(x^{\frac{1}{n}}-1\right).
\end{array}
\end{equation*}
Unter Gebrauch des  Grenzwertes\footnote{Selbst unter der Mehrdeutigkeit des Ausdrucks $x^{\frac{1}{n}}=|x|^{\frac{1}{n}}e^{\frac{i\varphi+2k\pi i}{n}}$ wird dieser Grenzwert für unendliches $n$ zu $1$.}

\begin{equation*}
    \lim_{n\rightarrow \infty} x^{\frac{1}{n}}=1,
\end{equation*}
 wo $\varphi$ das Argument der komplexen Zahl $x$ und $k$ alle natürlichen Zahlen bis hin zu $n$ bedeutet, ist der Grenzwert $=1$, und der Definition (\ref{eq: Grenzwert ln}), findet man schließlich

\begin{equation*}
    \log (x\cdot y) = \log  (x)+\log  (y),
\end{equation*}
wie gewünscht. Einen Beweis, welchen Euler wohl -- in seiner Wortwahl versteht sich -- anerkannt haben würde.\\

Jedoch gelangt man damit zu solch einer bizarren Gleichung
\begin{equation*}
  0= \log (1)=\log  ((-1) \cdot (-1))=  \log  (-1)+\log (-1)=2 \log  (-1).
\end{equation*}
Dies hätte $\log(-1)=0$ zur Folge, was natürlich nicht der Wahrheit entspricht, da ja bereits aus der Identität

\begin{equation*}
    e^{i \pi }+1=0
\end{equation*}
und der zuvor erwiesenen Mehrdeutigkeit des Logarithmus die Gleichung

\begin{equation*}
    \log  (-1)= i  \pi +2k\pi i \quad \text{für} \quad k \in \mathbb{Z}
\end{equation*}
folgt. Dies ist für die Euler'sche Theorie ein Problem, welches er auch bemerkt hat. In der Lösung zu Problem 2 in seiner Arbeit \cite{E168} sieht sich  Euler daher veranlasst   zu behaupten, dass zwar $2\log (-a)=2 \log(+a)$, aber nicht $\log(-a)=\log(+a)$ ist. Er versucht nun, diese Unannehmlichkeit zu umgehen, indem er richtig anmerkt, dass ja eigentlich

\begin{equation*}
    \log (-1)+\log (-1)= (2k+1)\pi \cdot i + (2l+1)\pi \cdot i
\end{equation*}
für ganze Zahlen $k$ und $l$ ist. Und letztere lassen sich eben stets so auswählen, dass tatsächlich

\begin{equation*}
    \log (-1)+\log (-1)=0
\end{equation*}
gilt. Dies ist ebenfalls korrekt, suggeriert aber keine Erklärung, warum jeweils verschiedene Werte zu wählen sind, geschweige denn, welche. Denn die spezielle Wahl $l=-k$ führt  zur richtigen bzw. gewollten Gleichung.  Unter diesen unendlich vielen Möglichkeiten sind aber anderweitig keine besonders ausgezeichnet. Dieser Umstand ist Euler ebenfalls bewusst, jedoch begnügt er sich mit just nachgezeichneten Ausführungen.\\

In anderem Zusammenhang erläutert er die paradoxe Rechnung
\begin{equation*}
    1=\sqrt{1\cdot 1}=\sqrt{(-1) \cdot (-1)}=\sqrt{-1}\cdot \sqrt{-1}=-1
\end{equation*}
nicht weiter  als zu erwähnen, dass eine solche Rechnung unrichtig sei. Ursächlich verhält es dabei  wie beim Logarithmus: Der Ausdruck $\sqrt{-1}$ kann   verschiedene Bedeutungen annehmen, modern gesprochen lässt sich für ihn $i$ und $-i$ setzen, von welchen jedoch meist ersterer gewählt wird. In obigen Beispiel tritt eben der Fall ein, in welchem $\sqrt{-1}$ einmal $=i$ und einmal $-i$ gesetzt werden muss, um die Richtigkeit der allgemeinen Gleichung

\begin{equation*}
    \sqrt{x\cdot y}=\sqrt{x}\cdot \sqrt{y}
\end{equation*}
auch für negative Zahlen $x$ und $y$ zu erhalten. Deswegen kann man auch nicht schlicht $\sqrt{-1}$ als die \textit{einzige} Zahl definieren, die mit sich selbst multipliziert $-1$ ergibt, wie Euler es in § 145 seiner \textit{Algebra} \cite{E387} tut. Er schreibt in seiner altertümlichen Ausdrucksweise:\\

\textit{``Ungeacht aber diese Zahlen als z. E. $\sqrt{-4}$, ihrer Natur nach ganz und gar ohnmöglich sind, so haben wir davon doch einen hinlänglichen Begriff, indem wir wißen, daß dadurch eine solche Zahl angedeutet werde, welche mit sich selbsten multipliziert zum Product $-4$ hervorbringe; und dieser Begriff ist zureichend um diese Zahlen in der Rechnung gehörig zu behandeln."} \\

Um dieses Problem zu umgehen, kann  die Quadratwurzel mithilfe von Logarithmen $\sqrt{x}:=e^{\frac{1}{2}\log x}$ definiert werden, sofern man auf der komplexen Ebene arbeiten möchte. Es wird demnach gleichsam das Euler'sche Argument, den Logarithmus aus einer Wurzel wie in Gleichung (\ref{eq: Grenzwert ln}) zu definieren, umgekehrt, und die Wurzel aus dem Logarithmus heraus definiert. Jedoch ist diese Erklärung \textit{a posteriori}, zumal die Bedeutungen entsprechend des gewünschten Endergebnisses ausgewählt werden. Dies wäre für Euler angesichts seiner  physiko--teleologischen Auffassung von der Mathematik und Natur höchst unbefriedigend gewesen. Jedoch fehlte ihm wie bereits erwähnt, das Konzept der Riemann'schen Fläche -- und gar die geometrische Interpretation der komplexen Zahlen, was unten (Abschnitt \ref{subsubsec: Komplexe Analysis}) noch eingehender erläutert werden wird --, welche  eine entsprechende Auswahl \textit{a priori} gestattet.\\

\paragraph{Falsche Aussagen zu mehrwertigen Funktion}
\label{para: Falsche Aussagen}

Die Schwierigkeiten, welche die mehrdeutigen Funktionen Euler bereitet haben, zeigen sich auch an den falschen Aussagen, welche er über sie tätigt. Dies betrifft, noch einmal auf den Logarithmus bezogen, solche Formulierungen wie diejenige, es ergäben sich die unendlich vielen Werte des Logarithmus daraus, dass die Gleichung 

\begin{equation*}
    \left(1+\dfrac{\log (x)}{n}\right)^n -x =0
\end{equation*}
ja $n$ Wurzeln haben muss, etwa wegen der Gültigkeit des Fundamentalsatzes der Algebra\footnote{Auch wenn Eulers Beweis in \cite{E170} erst etwas später erschien als seine Theorie des komplexen Logarithmus in \cite{E168}, wird er dennoch von der Gültigkeit des Fundamentalsatzes überzeugt gewesen sein.}. Dass $n$ hier eine unendlich große Zahl ist, ändert dabei für ihn nichts. Auch wenn Euler in diesem Fall Recht behalten sollte, ist dieses Argument im Allgemeinen nicht schlüssig. Denn das Polynom

\begin{equation*}
    P_n(x)=1+x+x^2+\cdots +x^n =\dfrac{1-x^{n+1}}{1-x}
\end{equation*}
besitzt gleichermaßen $n$ Nullstellen -- sie sind  $e^{\frac{2k \pi i}{n}}$, wobei $k$ alle natürlichen Zahlen bis einschließlich $n$ durchläuft --, aber im Grenzwert $n \rightarrow \infty$
hat der entsprechende Ausdruck keine endliche Nullstelle mehr. Es ist zu vermuten, dass Euler es auch im Falle des Logarithmus nicht als zwingend ansieht, es jedoch erwähnt, um seine Leser durch die Analogie zum endlichen Fall von der Unendlichwertigkeit des Logarithmus überzeugen. Unabhängig davon ist sein Beweis im Wesentlichen auf dieselbe Idee gestützt wie der des Sinusproduktes (\ref{eq: sine-Product}) -- es wird die Zerlegung von $a^n \pm b^n$ in trinomische Faktoren gesucht und der Grenzübergang $n \rightarrow \infty$ vollzogen. Grundlegend waren also die Euler'sche Formel

\begin{equation*}
    e^{ix}= \cos (x)+i \sin (x)
\end{equation*}
und die Auffassung des Logarithmus als Umkehrfunktion zur Exponentialfunktion. \\

Kann man im Fall des Logarithmus die Euler'schen Ausführungen noch mit einer entsprechenden Auslegung auch im heutigen Verständnis als korrekt ansehen, ist dies beim nachstehenden Fall nicht mehr möglich. Hierfür ziehe man eine Aussage aus  der Arbeit  \textit{``Methodus facilis inveniendi integrale huius formulae $\int \frac{\partial x}{x}\cdot \frac{x^{n+p}-2x^n\cos \zeta +x^{n-p}}{x^{2n}-2x^n \cos \theta +1}$  casu quo post integrationem ponitur vel $x = 1$ vel $x = \infty$"} (\cite{E620}, 1788, ges. 1776) (E620: ``Eine leichte Methode, die Integrale dieser Formel $\int \frac{\partial x}{x}\cdot \frac{x^{n+p}-2x^n\cos \zeta +x^{n-p}}{x^{2n}-2x^n \cos \theta +1}$ im Fall zu finden, in welchem nach der Integration entweder $x=1$ oder $x=\infty$ gesetzt wird") heran. Dort findet er in § IX zunächst die Formel 

\begin{equation*}
    \int\limits_{0}^{1} \dfrac{x^p+x^{-p}}{x^n+\left(f+\frac{1}{f}\right)+x^{-n}}\cdot \dfrac{dx}{x}= \dfrac{\pi \left(f^{\frac{p}{n}}-f^{-\frac{p}{n}}\right)}{n(f-f^{-1})\sin \frac{p\pi}{n}}.
\end{equation*}
Ist die obere Grenze $x=\infty$, dann ist der Wert doppelt so groß, wie Euler ebenfalls mitteilt. Gleichwohl diese Formel aus dem Fall abgleitet worden ist, in welchem $n$ eine natürlich Zahl ist, und $p$ so genommen werden muss, dass das Integral konvergiert, ist es zunächst ersichtlich -- was Euler auch bemerkt --, dass die Formel auch für negatives ganzzahliges $n$ gültig ist. Ebenso behält sie ihre Gültigkeit für rationale und sogar beliebige reelle $n$. Aber in § XII[a]\footnote{In der ursprünglichen Arbeit wird § XII irrtümlich wiederholt, weshalb in der \textit{Opera Omnia} Version der wiederholte Paragraph mit § XII[a] bezeichnet wird, was  hier übernommen ist.} gelangt er zur Formel

\begin{equation}
\label{eq: Euler False}
    \int\limits_{0}^1 \dfrac{x^p+x^{-p}}{2\cos (m\log x)-2\cos \theta}\cdot \dfrac{dx}{x}= \dfrac{p}{m\sqrt{-1}}\cdot \dfrac{e^{\frac{p}{m}(\pi-\theta)}-e^{-\frac{p}{m}(\pi-\theta)}}{\sin \theta \left(e^{\frac{\pi p}{m}}-e^{-\frac{\pi p}{m}}\right)},
\end{equation}
welche aus der vorherigen für $f+\frac{1}{f}=2\cos(\pi-\theta)$ und $n=\sqrt{-1}\cdot m$  nach einer leichten Vereinfachung hervorgeht. Euler behauptet demnach, dass dieses Integral einen rein imaginären Wert hat, obgleich man über eine rein reelle Funktion integriert, da nach Eulers Festlegungen alle involvierten Variablen $m$, $p$ und $\theta$ reell sind. Euler argumentiert diesbezüglich folgendermaßen: Da 

\begin{equation*}
    \int\limits_{0}^1 \dfrac{d x}{x \cos (m\log x)}= -\int\limits_{0}^{\infty} \dfrac{ dz}{\cos (mz)}
\end{equation*}
sowie

\begin{equation*}
    \int \dfrac{d \varphi}{\sin \varphi}= \log \tan \dfrac{1}{2}\varphi
\end{equation*}
gilt, wie man schnell prüft, wird man zu einem Logarithmus geführt. Dieser erlaubt nun wiederum aus seiner Eigenschaft, für reelle Argumente auch komplexwertig werden zu können, dass (\ref{eq: Euler False})  imaginär wird. Dabei belässt Euler seine Erklärungen. \\

Diese könnte man wie folgt fortführen, sofern man für $m>0$ die Formel

\begin{equation*}
    f(m):= \int\limits_{0}^{\infty} \dfrac{dx}{\cosh(mx)}= \left.\dfrac{2}{m}\cdot\arctan \left(\tanh \left(\dfrac{mx}{2}\right)\right)\right|_{0}^{\infty}= \dfrac{2}{m}\cdot\dfrac{\pi}{4}= \dfrac{\pi}{2m}.
\end{equation*}
bemerkt. Nun setze $im$ statt $m$, sodass

\begin{equation*}
    f(im)=\int\limits_{0}^{\infty} \dfrac{dx}{\cos(mx)}= \dfrac{\pi}{2im},
\end{equation*}
was  mit Eulers Ergebnis in Einklang steht. Ursächlich für imaginären Werte der Integrale ist in diesem Zusammenhang  ihre Divergenz. Diesen Integralen lässt sich, ebenso wie es Euler für Reihen in \cite{E247} getan hat (Abschnitt \ref{subsubsec: Der Begriff der Summe einer Reihe}), ein  endlicher Wert zuordnen. Für Integrale scheint Euler eine analoge Untersuchung nicht unternommen zu haben, was wiederum seiner Auffassung von Integralen als spezielle Summen geschuldet sein mag.\\

In diesem Zusammenhang sei bemerkt, dass nach Euler gewisse divergente Reihen einen rein imaginären Wert haben, sofern ihnen gemäß seiner Auffassung (siehe Abschnitt \ref{subsubsec: Der Begriff der Summe einer Reihe}) ein endlicher Wert zugeordnet wird. Er schreibt zum Beispiel in seiner Arbeit \textit{``De serie Lambertina plurimisque eius insignibus proprietatibus"} (\cite{E532}, 1783, ges. 1779) (E532: ``Über die Lambert'sche Reihe und ihre wunderbaren Eigenschaften")  in § 11, dass die stetig wachsenden Terme einer Reihe ein Zeichen für die Evaluation ihrer Summe zu einem imaginären Wert sind, genauer:\\

\textit{``[...], dass eine divergente Reihe resultieren würde, weil in der allgemeinen Form $\frac{n^n}{1\cdot 2 \cdot 3 \cdots n}$ der Zähler den Nenner immer mehr übersteigen würde und daher alle Terme sogar unendlich groß werden würden, was ein Anzeichen einer imaginären Summe ist."}\\

Während sich also das Integral (\ref{eq: Euler False}) noch in gewisser Weise konsistent deuten lässt, ist dies bei seinen Erläuterungen in § V von \cite{E620} nicht mehr möglich. Hier findet sich zunächst die Formel

\begin{equation*}
    \int\limits_{0}^{\infty} \dfrac{x^{\pm p}}{x^n-2\cos \theta +x^{-n}}\cdot \dfrac{dx}{x}=\dfrac{\pi \sin \frac{p}{n}(\pi-\theta)}{n\sin \theta \sin \frac{p\pi }{n}},
\end{equation*}
welche, entsprechenden Restriktionen unterworfen, korrekt ist. Jedoch sagt Euler dann in § VI, die Formel behielte ihre Richtigkeit auch nach der Translation $\theta \mapsto \theta +2 \pi k$ für beliebiges ganzzahliges $k$ bei, was nicht richtig ist. Obige Formel gilt nur für $0<\theta<2\pi$. Euler fühlt sich dabei zu seiner Behauptung bewogen, \textit{obwohl} das Integral bzw. der Integrand unter dieser Umbenennung keine Veränderung erfährt, anders als der Ausdruck auf der rechten Seite. Dies stellt für Euler jedoch kein Problem dar, sondern muss sich vielmehr so verhalten. Denn, so sagt Euler, das Integral linker Hand ist eigentlich mehrdeutig, was sich  im analogen Beispiel $\int \frac{d x}{1+xx}$ oder ausintegriert $\arctan x$ einfach nur unmittelbarer zeigt, da  $\arctan(x)$ aufgefasst als Logarithmus $=\frac{1}{2i}\log \frac{1+ix}{1-ix}$ als mehrwertige Funktion zu erkennen ist\footnote{In der Tat ist der $\arctan(x)$ bzw. das ihn darstellende Integral in Eulers Formel erfasst und sie gibt gar, mit der Interpretation über $\log$ als unendlichwertige Funktion, die richtigen Werte für alle Werte $\theta= \frac{\pi}{2} +2k\pi$.}. Auch hier schlösse sich, wie beim Auswahlproblem des Logarithmus,  die Frage an, welchen der möglichen unendlich vielen Werte man in speziellen Betrachtungen -- wie etwa die Addition mehrerer solcher Integrale -- zu wählen hat. Jedoch bleibt Euler dem Leser diese Antwort schuldig.


\subsection{Durch eine irreleitende Frage}
\label{subsec: Durch eine falsche Frage}

\epigraph{When things get too complicated, it sometimes makes sense to stop and wonder: Have I asked the right question?} {Enrico Bombieri}%

Manchmal trägt es sich zu, dass einem eine Einsicht verwehrt bleibt, weil man nicht die richtige Frage gestellt hat, wovor selbst Euler nicht gefeit war. Obwohl er alle Bausteine beisammen hatte und auch die entsprechenden Fähigkeiten mitbrachte, haben ihn bisweilen seine Fragestellungen an einem weiteren Vorankommen gehindert. Dies wird im Folgenden am Beispiel der Reduktion von elliptischen Integralen auf Normalform (Abschnitt \ref{subsubsec: Die Normalform von elliptischen Integralen}) und der Auflösbarkeit von polynomialen Gleichungen durch Radikale (Abschnitt \ref{subsubsec: Wegen Unbeweisbarkeit: Wurzeln von Polynomen}) verdeutlicht.

\subsubsection{Die Normalform von elliptischen Integralen}
\label{subsubsec: Die Normalform von elliptischen Integralen}

\epigraph{The mind when it has an old experience will add that data into its current experience, and it keeps coming up with wrong answers.}{L. Ron Hubbard}

Euler ist unbestreitbar der Begründer der Theorie der elliptischen Integrale. Sein Additionstheorem für die elliptischen Integrale gehört zu seinen schönsten Entdeckungen. Andererseits konnte er trotz seiner umfassenden Untersuchungen zu diesem Gegenstand seine Ergebnisse nicht zu einer prägnanten Theorie bündeln\footnote{Einer Tatsache, welcher er sich trotz seiner gesamten Erfolge in diesem Gebiet auch bewusst war.}, anders als  es ihm etwa bei der Variationsrechnung in seiner \textit{Methodus} \cite{E65} gelang. Die Gründe hierfür sollen im Folgenden diskutiert werden. Seine Herleitung des Additionstheorems für elliptische Integrale soll dabei nicht ausgespart werden und bildet den Anfang der Besprechung;  anschließend erfolgt der Übergang zu seinen Versuchen einer Reduktion auf Normalform. 

\paragraph{Erste Untersuchungen -- Fagnanos Verdopplungsformel}
\label{para: Herleitung seiner Additionstheoreme}

Elliptische Integrale, wie man heute unter anderem die Integrale nennt, die aus der Rektifizierung von Ellipsenbogen entspringen, treten vielerorts bei Euler auf, bevor er ab der Arbeit  \textit{``De integratione aequationis differentialis $\frac{mdx}{\sqrt{1-x^4}}=\frac{ndy}{\sqrt{1-y^4}}$"} (\cite{E251}, 1761, ges. 1751) (E251: ``Über die Integration der Differentialgleichung $\frac{mdx}{\sqrt{1-x^4}}=\frac{ndy}{\sqrt{1-y^4}}$") die Angelegenheit systematisch in Angriff nimmt.
Das vehemente Interesse ist dabei einer Entdeckung von G. Fagnano (1682--1766) zu verdanken, dass man für

\begin{equation*}
    F(x)= \int\limits_{0}^{x} \dfrac{dt}{\sqrt{1-t^4}}
\end{equation*}
die Verdopplungsformel

\begin{equation}
    \label{eq: Verdopplungsformel Fagnano}
    2F(x)=F \left(\dfrac{2x\sqrt{1-x^4}}{1+x^4}\right)
\end{equation}
nachweisen kann. Genauer beginnt die Euler'sche Entdeckungsreise in diesem Gebiet mit der zeitlich leicht vor \cite{E251} liegenden Arbeit \textit{``Observationes de comparatione arcuum curvarum irrectificibilium"} (\cite{E252}, 1761, ges. 1753) (E252: ``Beobachtungen über den Vergleich von nicht rektifizierbaren Kurven"), wo Euler in Theorem 5 bzw. § 30 auch die Gleichung  (\ref{eq: Verdopplungsformel Fagnano}) beweist. Euler schreibt dies in leicht anderer jedoch äquivalenter Form, indem er sagt, dass die Gleichung $u=\frac{2z\sqrt{1-z^4}}{1+z^4}$ der Differentialgleichung

\begin{equation*}
    \dfrac{du}{\sqrt{1-u^4}}=\dfrac{dz}{\sqrt{1-z^4}}
\end{equation*}
Genüge leistet. Diese Arbeit ist neben  dem ersten erkennbaren Schritt zum Euler'schen Additionstheorem für elliptische Integrale wegen der findigen Versuche  in den ersten Paragraphen, die vorgelegten Differentialgleichungen durch geschicktes Raten zu lösen, von Interesse. Charakteristisch für Euler werden dem Leser auch alle diversen, insbesondere die fehlgeschlagenen, Anläufe bis zum Ende vollends berechnet präsentiert. So sieht Euler sich in § 4 beispielsweise gezwungen, die Lösung

\begin{equation*}
    \int du\sqrt{\dfrac{1-nuu}{1-uu}}+\int dx \sqrt{\dfrac{1-nxx}{1-xx}}=xu\sqrt{n}+C,
\end{equation*}
wobei $C$ eine Konstante ist und

\begin{equation*}
    u=\sqrt{\dfrac{1-nxx}{n-nxx}}
\end{equation*}
erfüllt, zu verwerfen, weil diese  für erlaubte Werte von $n$ zu einem komplexwertigen Bogen\footnote{Da Euler dieses Problem als die Addition zweier \textit{reeller} Ellipsenbogen auffasst, muss für ihn $n<1$ sein.}, jedoch nicht nicht zu Ellipsenbogen führt\footnote{Diese und ähnliche Vermeidungen des Komplexen wird Euler am Vordringen in das Gebiet der Funktionentheorie hindern, wie  unten (Abschnitt \ref{subsubsec: Komplexe Analysis}) berichtet werden wird. Sie werden sich auch für das Scheitern bei der Reduktion der elliptischen Integrale als ursächlich herausstellen.}. Die Lösung

\begin{equation*}
    \int du\sqrt{\dfrac{1-nuu}{1-uu}}+\int dx \sqrt{\dfrac{1-nxx}{1-xx}}+nxu = C,
\end{equation*}
wobei

\begin{equation*}
    \dfrac{1-nxx}{1-xx}=\dfrac{1}{u}
\end{equation*}
gilt, genügt indes seinen Ansprüchen. 

\paragraph{Die kanonische Gleichung und das Additionstheorem}
\label{para: Die kanonische Gleichung und das Additionstheorem}

In § 7 der Arbeit \textit{``Specimen alterum methodi novae quantitates transcendentes inter se comparandi: De comparatione arcuum ellipsis"} (\cite{E261}, 1761, ges. 1755) (E261: ``Ein anderes Beispiel der neuen Methode transzendente Größen miteinander zu vergleichen; über den Vergleich von Ellipsenbogen") legt Euler die -- wie er sie nennt -- kanonische Gleichung vor
\begin{equation}
    \label{eq: Euler ellip. kanonisch}
    0= \alpha +\gamma (xx+yy)+2\delta xy+ \zeta xxyy,
\end{equation}
in welcher man aus heutiger Sicht leicht eine elliptische Kurve erkennt. Euler nutzt diese nun, um einmal $x$ als Funktion von $y$ und umgekehrt $y$ als Funktion von $x$ zu schreiben. Sein Ergebnis lautet:

\begin{equation*}
 \renewcommand{\arraystretch}{2,5}
 \setlength{\arraycolsep}{0.0mm}   
 \begin{array}{rcl}
     y &~=~& \dfrac{-\delta x+\sqrt{\delta \delta xx-(\alpha +\gamma xx)(\gamma +\zeta xx)}}{\gamma +\zeta xx}, \\
      x &~=~& \dfrac{-\delta y-\sqrt{\delta \delta yy-(\alpha +\gamma yy)(\gamma +\zeta yy)}}{\gamma +\zeta yy},
 \end{array}
\end{equation*}
wo er später abkürzend

\begin{equation*}
    X= \sqrt{\delta \delta xx-(\alpha +\gamma xx)(\gamma +\zeta xx)} \quad \text{und} \quad Y=\sqrt{\delta \delta yy-(\alpha +\gamma yy)(\gamma +\zeta yy)}
\end{equation*}
schreibt. Damit berechnet man

\begin{equation}
\label{eq: X,Y als f(x,y)}
    X=\gamma y+\delta x +\zeta xx y ~~ \text{und}  ~~ -Y=\gamma x +\delta y +\zeta xyy.
\end{equation}
Differenziert man nun (\ref{eq: Euler ellip. kanonisch}), findet man

\begin{equation*}
    0=dx(\gamma x+\delta y+\zeta xyy)+dy(\gamma y+\delta x+\zeta xxy),
\end{equation*}
oder mit (\ref{eq: X,Y als f(x,y)}) auch

\begin{equation*}
    \dfrac{dy}{Y}-\dfrac{dx}{X}=0.
\end{equation*}
Das ist  gerade eine Differentialgleichung mit Polynomen vom Grad $4$ unter der Wurzel, womit die kanonische Gleichung (\ref{eq: Euler ellip. kanonisch}) ihr Integral ist. \\

Dass man statt der allgemeinen Gleichung

\begin{equation}
\label{eq: Elliptisches Integral}
    \dfrac{dx}{\sqrt{A+Bx+Cx^2+Dx^3+Ex^4}}=\dfrac{dy}{\sqrt{A+By+Cy^2+Dy^3+Ey^4}}
\end{equation}
für beliebige Koeffizienten $A$, $B$, $C$, $D$, $E$ ohne Einschränkung der Allgemeinheit auch diese reduzierte

\begin{equation*}
    \dfrac{dx}{\sqrt{A+Cx^2+Ex^4}}=\dfrac{dy}{\sqrt{A+Cy^2+Ey^4}}
\end{equation*}
betrachten kann, hat Euler später unter Verwendung einer gebrochen linearen Substitution bewiesen. Er schreibt in § 2 seiner Abhandlung  \textit{``Integratio aequationis $\frac{dx}{\sqrt{A+Bx+Cx^2+Dx^3+Ex^4}}=\frac{dy}{\sqrt{A+By+Cy^2+Dy^3+Ey^4}}$"} (\cite{E345}, 1768, ges. 1765) (E345: ``Integration der Gleichung $\frac{dx}{\sqrt{A+Bx+Cx^2+Dx^3+Ex^4}}=\frac{dy}{\sqrt{A+By+Cy^2+Dy^3+Ey^4}}$"):\\

\textit{``Und zuerst bemerke ich freilich, dass die vorgelegte Gleichung immer in eine Form solcher Art überführt werden kann, in welcher die Koeffizienten $B$ und $D$ verschwinden, was freilich für je einen der beiden aus den Elementen hinreichend bekannt ist. Damit aber beide zugleich verschwinden können, bedarf es einer dafür eigenen Form: Für $x=\frac{mz+a}{nz+b}$ geht nämlich die erste Form, welcher die andere ja ähnlich ist, zunächst in diese über [...]"}\\

Euler präsentiert die Rechnungen sehr ausführlich und gelangt in § 4 schließlich zu letztgenannter Gleichung. Auch von dieser ist (\ref{eq: Euler ellip. kanonisch}) wieder die Integralgleichung  der letzten Gleichung, wenn die Koeffizienten $\alpha$, $\gamma$, $\delta$, $\zeta$ entsprechend angepasst werden\footnote{Da die Anzahl um eins größer ist als die von $A$,$B$, $C$, ist eine der ersten beliebig und nimmt den Platz der Integrationskonstante ein.}.  Euler selbst ist mit dieser  vom Gesuchten ausgehenden Lösung unzufrieden gewesen, weil eben die Lösung bereits bekannt sein muss\footnote{Dies läuft wie schon bei anderen Begebenheiten seinem physiko--teleologischen Weltbild zuwider.}. Eine Methode, die gemäß seiner Ansprüche hinreichend direkt zur Integration von (\ref{eq: Elliptisches Integral}) führt,  stellt Euler später in seiner Arbeit \textit{``Dilucidationes super methodo elegantissima, qua illustris de la Grange usus est in integranda aequatione differentiali $\frac{dx}{\sqrt{X}}=\frac{dy}{\sqrt{Y}}$"} (\cite{E506}, 1781, ges. 1778) (E506: ``Erläuterungen zur höchst eleganten Methode, welche der illustre Lagrange  für die Integration der Differentialgleichung $\frac{dx}{\sqrt{X}}=\frac{dy}{\sqrt{Y}}$ verwendet hat") vor. {Zwar präsentiert er bereits in  der schon genannten Arbeit \cite{E345} eine Methode, welche sich ebenfalls als direkte Methode sehen ließe, ihm selbst ist das dortige Vorgehen jedoch wegen der vielen Substitutionen zu unnatürlich. Genauer schreibt er in § 16:\\

\textit{``Zwar habe ich hier das Integral der vorgelegten Differentialgleichung mit einer direkten Methode erhalten, ich kann aber dennoch nicht in Abrede stellen, dass dies über viele Umwege geleistet worden ist, so dass kaum zu erwarten ist, dass jemandem diese Operationen in den Sinn hätten kommen können."}\footnote{Die in \cite{E506} erwähnte Operation, wie bereits aus dem Titel ersichtlich, ist aber Lagrange zuzuschreiben, wie Euler auch gleich zu Beginn der Arbeit bemerkt. Sie beruht gleichermaßen beruht  auf einem eher künstlichen Ansatz, welcher sich nicht auf höhere Grade von Polynomen verallgemeinern lässt.} \\

Zusammenfassend lässt sich feststellen, dass Euler selbst keine Methode abzuleiten vermochte, die Differentialgleichung (\ref{eq: Elliptisches Integral}) direkt und gleichzeitig a priori zu integrieren, hatte er doch bereits in § 1 von \cite{E345} zu seiner Methode über die kanonische Gleichung (\ref{eq: Euler ellip. kanonisch}) resümiert:\\

\textit{``Mit einer völlig einzigartigen und zugleich merkwürdigen Methode war ich einst zur Integration dieser Gleichung gelangt, deren Integral, und gar das vollständige, ich entdeckt habe, in einer algebraischen Gleichung enthalten zu sein. Das scheint umso wundersamer, zumal das Integral jeder der beiden Formeln einzeln nicht nur nicht algebraisch, sondern nicht einmal  durch die Quadratur des Kreises oder der Hyperbel ausgedrückt werden konnte. Dann trug  sich indes dieses besonders bemerkenswerte Phänomen zu, dass kein direkter Weg offen stand, dieses algebraische Integral zu finden. Aber keine Gelegenheit scheint mehr geeignet, die Grenzen der Analysis zu erweitern, als wenn wir uns bemühen, dasselbe, was wir vermöge einer ungradlinigen Methode ermittelt haben, mit einer direkten Methode ausfindig zu machen."}\\

Aus heutiger Sicht fehlten Euler für dieses Unterfangen die entsprechenden nicht elementaren Funktionen. Die Lösung von (\ref{eq: Elliptisches Integral}) verlangt nämlich die bewusste Einführung von neuen transzendenten Funktionen, die eben genau über diese  Eigenschaft definiert werden und davon ausgehend untersucht werden\footnote{Euler selbst scheint in seinem Opus bei keiner Gelegenheit eine Funktion über ihre Differentialgleichung definiert zu haben. Differentialgleichungen waren für ihn immer gleichsam ein Problem, das es aufzulösen galt.}. Dass die üblichen elementaren Funktionen -- also Logarithmen und Kreisfunktionen -- hierfür nicht hinreichen, ist Euler vermutlich seit seiner Abhandlung \textit{``Theorematum quorundam arithmeticorum demonstrationes"} (\cite{E98}, 1747, ges. 1738) (E98: ``Beweise gewisser zahlentheoretischer Theoreme") bewusst gewesen. Hier beweist er in Problem 2 die Unmöglichkeit der Auflösbarkeit der Gleichung

\begin{equation*}
    x^4-y^4=z^2
\end{equation*}
über den natürlichen Zahlen, woraus sich wiederum ergibt, dass bereits das Integral

\begin{equation*}
    \dfrac{du}{\sqrt{1-u^4}}
\end{equation*}
durch keine bekannte Substitution rational gemacht werden kann\footnote{Der Zusammenhang der Existenz einer solchen Substitution mit dem Geschlecht der Lemniskate und generell dem Geschlecht einer Kurve war Euler noch unbekannt. Der Begriff des Geschlechts einer algebraischen Kurve kam erst wesentlich später mit Clebsch (1833--1872), wie man etwa Kleins Buch \textit{``Vorlesungen über die Entwicklung der Mathematik im 19. Jahrhundert"} (\cite{Kl56}, 1956) entnehmen kann.}. Somit implizieren Integrale von Formen wie das über letztgenannte Differentialformeln die Existenz neuer Transzendenten.  Diese sind  seit Legendres drei Werken \textit{``Exercices du calcul intégral -- Tome premier"} (\cite{Le11}, 1811) (``Übungen zur Integralrechnung -- Erstes Buch"), \textit{``Exercices du calcul intégral -- Tome second"} (\cite{Le17}, 1817) (``Übungen zur Integralrechnung -- Zweites Buch"), \textit{ ``Exercices du calcul intégral -- Tome troisième"} (\cite{Le16}, 1816) (``Übungen zur Integralrechnung -- Drittes Buch") die elliptischen Integrale. Genau diese bewusste Untersuchung von Funktionen aus ihren Eigenschaften heraus -- eben ohne Zuhilfenahme von geometrischen Überlegungen -- findet sich Euler bei Euler noch nicht. Er ist vielmehr an Fragen über die Eigenschaften der bekannten Funktionen selbst interessiert.\\

Auf die Legendre'sche Klassifikation von elliptischen Integralen Bezug nehmend, hat  Euler jedoch -- natürlich ohne Kenntnis von Legendres Einteilung  gehabt zu haben -- die Additionstheoreme für die elliptischen Integrale erster und zweiter Art in \cite{E345} vollständig abgehandelt, dasjenige für die dritter Art findet sich in seinen Formeln der Abhandlung \textit{``Plenior explicatio circa comparationem quantitatum in formula integrali $\int \frac{Zdz}{\sqrt{1+mzz+nz^4}}$ contentarum denotante $Z$ functionem quamcunque rationalem ipsius $zz$""} (\cite{E581}, 1785, ges. 1775) (E581: ``Eine umfassendere Erklärung zum Vergleich der in der Integralform $\int \frac{Zdz}{\sqrt{1+mzz+nz^4}}$ enthaltenen Größen, wobei $Z$ irgendeine rationale Funktion von $zz$ ist").

\paragraph{Versuch einer Reduktion auf Normalformen}
\label{para: Versuch einer Reduktion}

 Euler ist  ebenfalls von der Entdeckung einer tiefer liegenden Struktur der elliptischen Integrale eingenommen gewesen, was man heute am ehesten als die Reduktion auf die just erwähnten Normalformen bezeichnen würde. Insbesondere die Arbeiten \textit{``Consideratio formularum, quarum integratio per arcus sectionum conicarum absolvi potest"} (\cite{E273}, 1763, ges. 1760) (E273: ``Betrachtung von Formeln, deren Integration mit den Bogen von Kegelschnitten durchgeführt werden kann") sowie \textit{``De reductione formularum integralium ad rectificationem ellipsis ac hyperbolae"} (\cite{E295}, 1766, ges. 1764) (E295: ``Über die Reduktion von Integralformeln auf die Rektifizierung der Ellipse und der Hyperbel").\\

\begin{figure}
    \centering
     \includegraphics[scale=0.8]{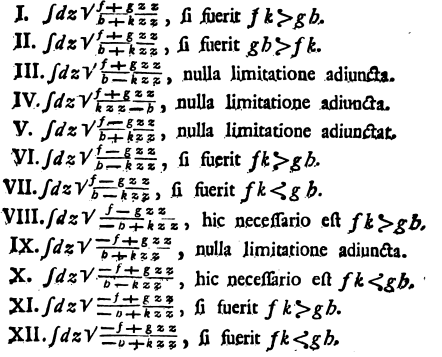}
    \caption{Eulers Einteilung der elliptischen Integrale aus § 38 seiner Arbeit \cite{E295} in $12$ Klassen.}
    \label{fig:E295EulersListeElliptischeIntegral}
\end{figure}

 Hier lässt sich Euler allerdings von einem geometrischen Bild bei der Einteilung leiten, was -- wie seit Legendres Untersuchungen bekannt ist -- nicht zielführend  sein konnte, wenn man die Legendre'sche Einteilung als Maßstab heranzieht. In \cite{E273} gelangte er dabei insgesamt zu 12 Formeln, die sich auf die 3 von Legendre reduzieren ließen, sofern man das geometrische Bild fallen lässt und die Integrale auch über der komplexen Ebene betrachtet. Die Situation wird von Adolf Krazer (1858--1926) in dem Vorwort zu Band 20 von Serie 1 der \textit{Opera Omnia} (\cite{OO20}, 1912)  treffend zusammengefasst: \\

\textit{``Fragen wir aber, warum das in dem anderen Hauptteil der Abhandlungen (insbesondere 295, 273) behandelte Problem der Reduktion der Integrale auf feste Normalformeln und der Zurückführung des allgemeinen Integrals auf diese, trotzdem es für Euler, der rechnen konnte und wollte wie kein anderer, wie geschaffen war, keine so glückliche Lösung gefunden hat, so müssen wir, wie schon oben erwähnt, dem Nichtloskommen von geometrischen Vorstellungen die Schuld geben. Ein entscheidender Fortschritt in dieser Richtung konnte erst geschehen, wenn die geometrische Grundlage, der allerdings die Lehre von den elliptischen Integralen bisher fast alles verdankte, zurücktrat und einer Behandlung der Integrale um ihrer selbst Willen Platz machte; [...]"}\\

\paragraph{Die Variationsrechnung als weiteres analoges Beispiel}

Bei dem Versuch der Reduktion auf Normalformen erweist sich  insbesondere die Geometrie, die Euler bis zu diesem Punkt hatte vordringen lassen, als Hindernis für das Vorankommen bei der Theorie der elliptischen Integrale für Euler, was in merkwürdiger Parallelität auf die Variationsrechnung zutrifft. Eulers \textit{Methodus}  \cite{E65}  liefert hierfür ein treffliches Beispiel. \\

\begin{figure}
    \centering
  \includegraphics[scale=1.0]{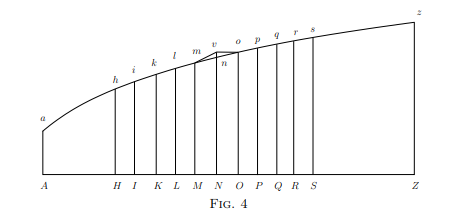}
    \caption{Eulers (mit Latex nachgezeichnete) Figur aus seinem Buch \cite{E65}, welche er zur Herleitung von den Euler-Lagrange'schen Differentialgleichung verwendet.}
    \label{fig:E65Figur}
\end{figure}

Denn die geometrischen Methoden haben Euler zu den fundamentalen Erkenntnissen der Variationsrechnung geführt, so leitet er etwa die Euler--Lagrange--Gleichungen\footnote{Dass die Euler--Lagrange Gleichungen der Variationsrechnung  eine notwendige Bedingung für die Lösung eines Variationsproblems sind, hat Euler spätestens in seiner Arbeit \textit{``De insigni paradoxo, quod in analysi maximorum et minimorum occurrit"} (\cite{E735}, 1811, ges. 1779) (E735: ``Über ein außergewöhnliches Paradoxon, welches bei der Analysis von Maxima und Minima auftritt") bemerkt. Das Paradoxon, welches hier auftritt, besteht darin, dass die Lösung des dort untersuchten Variationsproblems nicht die Euler--Langrange'schen Gleichungen auf dem gesamten untersuchten Intervall erfüllt. Die Lösungsfunktion ist eine lediglich Lipschitz--stetige Funktion.} ab und kann auch Variationsprobleme mit Nebenbedingungen von geometrischen Vorstellungen geleitet entsprechend als Differentialgleichungen formulieren. Zusätzlich ist Euler hier  von der Auffassung geleitet, dass Variationen stets existieren. Im  zweiten Anhang seines Buchs, direkt in § 1, schreibt er:\\

\textit{``Weil ja alle Wirkungen der Natur einem gewissen Gesetz von Maximum oder Minimum folgen, besteht kein Zweifel, dass bei den Kurven, welche Körper beschreiben, wenn sie von irgendwelchen Kräften angegriffen werden, eine gewisse Eigenschaft des Maximums oder Minimums auftritt."}\\

Auf die mögliche Frage, warum alle Wirkungen einem solchen Gesetz folgen mögen, mag man die Euler'schen Ausführungen im ersten Anhang in \cite{E65} (über Elastische Kurven) heranziehen. Hier sagt er, gleichzeitig sein physiko--teleologisches Weltbild bestätigend, in § 1:\\

\textit{``Denn weil die ganze Welt in vollkommener Art errichtet worden ist, und dies vom weisesten Schöpfer, gibt es in der Welt nichts, worin sich nicht das Prinzip der Maxima und Minima zeigt; deswegen besteht auch keinerlei Zweifel, dass jede Wirkung der Welt in endgültigen Ursachen mithilfe der Methode der Maxima und der Minima mit gleichem Erfolg bestimmt werden kann wie aus den sie bewirkenden Gründen selbst."}\\

Damit ist für Euler zugleich die Frage der Existenz von Variationen geklärt\footnote{Gleichzeitig schimmert hier der Glaube von Euler an die Gültigkeit des Prinzips vom hinreichenden Grund hervor.}.  Die Überführung der Variationsrechnung in ein reines Kalkül schreibt Euler indes Lagrange zu, dessen Methode, die konsequente Verwendung eines Symbols zur Anzeige der Variation $\delta$ in vollkommener Analogie zum Differential $d$, Euler in aller Ausführlichkeit in seinen beiden Arbeiten  \textit{``Elementa calculi variationum"} (\cite{E296}, 1766, ges. 1756) (E296: ``Elemente der Variationsrechnung") -- in § 12 tritt die Definition  des Variationssymbols $\delta$  zum ersten Mal auf -- und \textit{``Analytica explicatio methodi maximorum et minimorum"} \cite{E297} (E297: ``Analytische Erklärung der Methode der Maxima und Minima") vorstellt. Die Arbeit \textit{``Methodus nova et facilis calculum variationum tractandi"} (\cite{E420}, 1772, ges. 1771) (E420: ``Eine neue und leichte Methode das Variationsrechnung zu behandeln") mutet indes sehr modern an (insbesondere die Ausführungen in §§ 5--7), weil Eulers Zugang zum Finden der Variation an die Methoden erinnert, wie sie noch gegenwärtig in einführenden Veranstaltungen zur theoretischen Physik gelehrt werden. Die Beifügung der Euler'schen Erläuterungen zum Nachweis des Gesagten hinreichen: \\

\textit{``Natürlich betrachte ich für jenes Inkrement, welches ich Variation genannt habe, die Größe $y$ nicht weiter als eine Funktion der Variable $x$ allein, sondern führe sie als Funktion der zwei Variablen $x$ und $t$ in die Rechnung ein, während auf diese Weise nämlich $dx\left(\frac{dy}{dx}\right)$ das tatsächliche Differential von bezeichnet,  wird diese Formel $dt \left(\frac{dy}{dt}\right)$ dasselbe bedeuten können, was wir zuvor mit dem Symbol $\delta y$ angezeigt haben. Um all das besser zu illustrieren, wollen wir $y$ als die einer Abszisse $x$ entsprechende Ordinate einer gewissen Kurve auffassen und im Variationskalkül wird nun eine andere Relation verlangt, welche alle anderen dieser hinreichend nahe kommenden Kurven zugleich erfasst; aber alle Kurven dieser Art, wenn $X$ jene Funktion bedeutet, welcher $y$ gleich wird, werden offenkundig in einer solchen Gleichung erfasst werden können $y=X+tV$, während $V$ eine beliebige Funktion von $x$ bedeutet. [...] Für das Finden der Variation ist demnach $x$ wie eine Konstante zu behandeln und das Differential von $y$ aus der Veränderlichkeit von $t$ zu entnehmen."} \\

Die Anwendung auf Funktionale der Form $\int ydx$ (Euler selbst benutzt wie in seiner \textit{Methodus} \cite{E65} den Formelbuchstaben $Z$) erfolgt schließlich ab § 8 bis § 15, wo sich in entsprechender Form auch die Euler--Lagrange Gleichungen finden. \\

\begin{figure}
    \centering
   \includegraphics[scale=0.7]{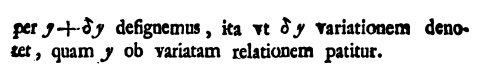}
    \caption{Euler Definition  des Wortes Variation und des Zeichens $\delta$ für sein Kalkül in seiner Arbeit \cite{E296}.}
    \label{fig:E296Variation}
\end{figure}

Die Geometrie bzw. das Arbeiten mit Bildern hat Euler also nie ganz losgelassen, aber teilweise an weiteren Fortschritten gehindert. Dies mutet umso ironischer an, da Euler im Vorwort seiner \textit{Calculi Differentialis} \cite{E212} nicht ohne Stolz verkündet, dass er seine Theorie der Differentialrechnung gänzlich ohne Zuhilfenahme von Bildern präsentieren wird. Lagrange sollte es ihm später in seinem Grundlagenwerk zur Mechanik \textit{``Mécanique Analytique"} (\cite{La88}, 1788)  später gleichtun. Man vergleiche das Zitat zu seinem Buch, welches man im Buch \textit{``Mathematics, from the Points of View of the Mathematician and of the Physicist"}  (\cite{Ho12}, 1912) wie folgt  ins Englische übersetzt findet:\\

\textit{``We have already various treatises on Mechanics, but the plan of this one is entirely new. I intend to reduce the theory of this Science, and the art of solving problems relating to it, to general formulae, the simple development of which provides all the equations necessary for the solution of each problem. I hope that the manner in which I have tried to attain this object will leave nothing to be desired. No diagrams will be found in this work. The methods that I explain require neither geometrical, nor mechanical, constructions or reasoning, but only algebraical operations in accordance with regular and uniform procedure. Those who love Analysis will see with pleasure that Mechanics has become a branch of it, and will be grateful to me for having thus extended its domain."}

\subsubsection{Wegen Unbeweisbarkeit: Auflösbarkeit von polynomialen Gleichungen mit Radikalen}
\label{subsubsec: Wegen Unbeweisbarkeit: Wurzeln von Polynomen}

\epigraph{Niels Abel shocked the mathematical world by showing that no ``solution by radicals"{} was possible for fifth-- or higher--degree equations. [...] Abel's proof [...] stands as a landmark in mathematics history.}{William Dunham}

Endlich soll noch ein Fall betrachtet werden, bei welchem die eigentliche Unbeweisbarkeit den Euler'schen Fortschritt verhindert hat: Die Lösung von polynomialen Gleichungen mithilfe von Radikalen. Stellt man hier wie Euler die Frage, wie die Auflösung solcher Gleichungen mit Wurzelausdrücken gelingt, ist man von Beginn an zu Scheitern verurteilt, wie man seit Abels Arbeit \textit{Mémoire sur les èquations algébriques, ou l'on démontre l'Impossibilité de la Résolution de l'équation générale du cinquième Degré."} (\cite{Ab24}, 1824) (``Abhandlung über algebraische Gleichungen, oder der Nachweis der Unmöglichkeit der Auflösung der allgemeinen Gleichung vom fünften Grad") weiß. Unabhängig davon wird im Folgenden neben den Euler'schen Auflösungsformeln mit Radikalen, welche er für spezielle Fälle zu finden vermochte, die Lambert'sche Reihe vorgestellt. Zudem lieferte Euler einen  Ansatz, welcher die zu lösende algebraische Gleichung auf eine differentielle zurückführt.

\paragraph{Auflösbarkeit mit Radikalen}

 Gegeben sei  die Gleichung:
\begin{equation}
\label{eq: P_n(x)=0}
    a_nx^n+a_{n-1}x^{n-1}+\cdots +a_1x +a_0=0.
\end{equation}
Euler versteht dabei stets die Koeffizienten als reell und ist in seiner Schaffenszeit wiederholt zur Frage nach einer expliziten Lösungsformel für die Wurzeln des Polynoms zurückgekehrt. Stellvertretend seien die Arbeiten \textit{``De formis radicum aequationum cuiusque ordinis coniectatio"} (\cite{E30}, 1738, ges. 1733) (E30: ``Eine Vermutung über die Formen der Wurzeln von Gleichungen einer jeden Ordnung"), \textit{``De resolutione aequationum cuiusvis gradus"} (\cite{E282}, 1764, ges. 1753) (E282: ``Über die Auflösung von Gleichungen jedweden Grades"), \textit{``Innumerae aequationum formae ex omnibus ordinibus, quarum resolutio exhiberi potest"} (\cite{E644}, 1790, ges. 1776) (E644: ``Unzählige Formen von Gleichungen aus allen Ordnungen, deren Auflösung dargeboten werden kann") in diesem Zusammenhang genannt. In den ersten beiden Arbeiten ist er dabei von der Vermutung

\begin{equation}
\label{eq: Euler Vermutung Radikale}
    x= \sqrt[n]{A_1}+\sqrt[n]{A_2}+\sqrt[n]{A_3}+\cdots + \sqrt[n]{A_{n-1}}
\end{equation}
mit konstanten $A_1$, $A_2$, $\cdots$, $A_{n-1}$, die noch zu bestimmen sind, geleitet. Für die Fälle $n=1$, $n=2$, $n=3$, $n=4$ gelangt er so zu den bekannten Lösungsformeln ($pq$-Formel für $n=2$, die Cardano'sche Formel für $n=3$ und die Ferrari'sche Formel\footnote{Eulers Lösungsformel weicht leicht von der von Ferrari (1522--1565) ab. Der Ansatz, die biquadratische Gleichung auf die Lösung einer kubischen zu reduzieren, ist bei beiden Autoren zu finden. Man vergleiche auch Eulers Ausführungen in Kapitel 14 und 15 des zweiten Teils seiner \textit{Algebra} zu diesem Gegenstand.} für $n=4$). Dass er  für $n\geq 5$  seinen Erfolg nicht wiederholen kann, ist erst seit Abels eingangs erwähnter Abhandlung  \cite{Ab24} aus dem Jahr 1824 bekannt. Seine Misserfolge bei der allgemeinen Auflösung haben Euler  zu speziellen Fällen übergehen lassen, in welchen die Auflösbarkeit vermöge (\ref{eq: Euler Vermutung Radikale}) gelingt. Die letzte der oben genannten Abhandlungen betrachtet hingegen spezielle Formen von Gleichungen beliebigen Grades, welche eine Reduktion auf eine Resolvente von einem niederen Grad zulassen. Es seien an dieser Stelle jedoch zunächst die wesentlichen Gedanken zu den Polynomen vom $5.$ Grade aus \cite{E282} mitgeteilt. Eulers Rechnungen beginnen ab § 30. Er setzt gemäß (\ref{eq: Euler Vermutung Radikale})

\begin{equation*}
    x= \mathfrak{A}\sqrt[5]{v}+\mathfrak{B}\sqrt[5]{v^2}+\mathfrak{C}\sqrt[5]{v^3}+\mathfrak{D}\sqrt[5]{v^4}
\end{equation*}
als Lösung der Gleichung

\begin{equation*}
    x^5=Ax^3+Bx^4+Cx+D
\end{equation*}
an. Die Wurzeln dieser Gleichung nennt Euler $\alpha$, $\beta$, $\gamma$, $\delta$, $\varepsilon$ und sagt, sie sind wie folgt gegeben:\\

\begin{equation*}
    \renewcommand{\arraystretch}{1,5}
\setlength{\arraycolsep}{0.5mm}
\begin{array}{rcrcrcrcrcrclclcl}
    \alpha & = & \mathfrak{A} \cdot 1\cdot  \sqrt[5]{v} & + & \mathfrak{B} \cdot 1^2 \cdot  \sqrt[5]{v^2} & + & \mathfrak{C} \cdot 1^3 \cdot  \sqrt[5]{v^3} & + & \mathfrak{D} \cdot 1^4 \cdot  \sqrt[5]{v^4} &  \\
     \beta & = & \mathfrak{A} \cdot \mathfrak{a} \cdot  \sqrt[5]{v} & + & \mathfrak{B} \cdot \mathfrak{a}^2 \cdot  \sqrt[5]{v^2} & + & \mathfrak{C} \cdot \mathfrak{a}^3 \cdot  \sqrt[5]{v^3} & + & \mathfrak{D} \cdot \mathfrak{a}^4 \cdot  \sqrt[5]{v^4}  \\
     \gamma & = & \mathfrak{A} \cdot \mathfrak{b} \cdot  \sqrt[5]{v} & + & \mathfrak{B} \cdot \mathfrak{b}^2 \cdot  \sqrt[5]{v^2} & + & \mathfrak{C} \cdot \mathfrak{b}^3 \cdot  \sqrt[5]{v^3} & + & \mathfrak{D} \cdot \mathfrak{b}^4 \cdot  \sqrt[5]{v^4}  \\
     \delta & = & \mathfrak{A} \cdot \mathfrak{c} \cdot  \sqrt[5]{v} & + & \mathfrak{B} \cdot \mathfrak{c}^2 \cdot  \sqrt[5]{v^2} & + & \mathfrak{C} \cdot \mathfrak{c}^3 \cdot  \sqrt[5]{v^3} & + & \mathfrak{D} \cdot \mathfrak{c}^4 \cdot  \sqrt[5]{v^4}  \\
     \varepsilon & = & \mathfrak{A} \cdot \mathfrak{d} \cdot  \sqrt[5]{v} & + & \mathfrak{B} \cdot \mathfrak{d}^2 \cdot  \sqrt[5]{v^2} & + & \mathfrak{C} \cdot \mathfrak{d}^3 \cdot  \sqrt[5]{v^3} & + & \mathfrak{D} \cdot \mathfrak{d}^4 \cdot  \sqrt[5]{v^4}
\end{array}
\end{equation*}
wobei $1, \mathfrak{a}, \mathfrak{b}, \mathfrak{c}, \mathfrak{d}$ die fünften Einheitswurzeln sind und somit $x^5-1=0$ lösen. Aus dieser Annahme leitet Euler unter Verwendung der Newton--Girad'schen Formeln  Beziehungen zwischen den gegebenen Koeffizienten $A$, $B$, $C$, $D$ und den angenommenen Größen $\mathfrak{a}, \mathfrak{b}, \mathfrak{c}, \mathfrak{d}$ und $\mathfrak{A}, \mathfrak{B}, \mathfrak{C}, \mathfrak{D}$ sowie $v$ her. In § 36 findet man die Relationen:

\begin{equation*}
    \renewcommand{\arraystretch}{1,5}
\setlength{\arraycolsep}{0.5mm}
\begin{array}{rcl}
     A & = & 5(\mathfrak{AD}+\mathfrak{BC})v, \\
     B & = & 5 (\mathfrak{A}^2\mathfrak{C}+ \mathfrak{AB}^2+\mathfrak{BD}^2v+\mathfrak{C}^2\mathfrak{D}v)v, \\
     C & = & 5 (\mathfrak{A}^3\mathfrak{B}+\mathfrak{B}^3\mathfrak{D}v+\mathfrak{AC}^3v+\mathfrak{CD}^3v^2)v-5(\mathfrak{A}^2\mathfrak{D^2}+\mathfrak{B}^2\mathfrak{C}^2)v^2+ 5\mathfrak{ABCD}v^2, \\
     D & = & \mathfrak{A}^5v+\mathfrak{B}^5v^2+\mathfrak{C}^5v^3+\mathfrak{D}^5v^4-5(\mathfrak{A}^3\mathfrak{CD}+\mathfrak{AB}^3\mathfrak{C}+\mathfrak{BC}^3\mathfrak{D}v+\mathfrak{ABD}^3v)v^2 \\
      &  + & 5(\mathfrak{A}^2\mathfrak{B}^2\mathfrak{D}+\mathfrak{A}^2\mathfrak{B}\mathfrak{C}^2+\mathfrak{AC}^2\mathfrak{D}^2v+\mathfrak{B}^2\mathfrak{CD}^2v)v^2.
\end{array}
\end{equation*}
Im folgenden Paragraphen erklärt Euler nun, dass man diese Gleichungen nach $v$ aufzulösen hat, erwähnt aber auch die damit verbundenen Schwierigkeiten. Er schreibt weiter:\\

\textit{``Mit hinreichender Gewissheit lässt sich aber vermuten, wenn diese Elimination entsprechend durchgeführt werden würde, dass schließlich zu einer Gleichung vierten Grades gelangt werden kann, von welcher der Wert von $v$ bestimmt wird. Wenn nämlich eine Gleichung höheren Grades hervorginge, dann würde auch der Wert von $v$ Wurzelzeichen  des selben Grades beinhalten, was absurd erscheint. Weil aber die Menge der Terme diese Aufgabe dermaßen schwierig gestaltet, dass sie nicht einmal mit Erfolg versucht werden kann, wird es nicht unangemessen sein, gewisse weniger allgemeine Fälle zu entwickeln, welche nicht zu dermaßen komplizierten Fällen führen."}\\

In § 44 der Arbeit \cite{E282} gelangt er nach seinen Rechnungen zu einer \textit{hinreichenden} Bedingung, wann eine Gleichung vom Grad $5$ mit Radikalen auflösbar ist. Ist wie oben

\begin{equation*}
    x^5=Ax^3+Bx^2+Cx+D
\end{equation*}
vorgelegt, gelingt eine Auflösung für die folgende Wahl der Koeffizienten:

\begin{equation*}
     \renewcommand{\arraystretch}{2,0}
\setlength{\arraycolsep}{0.5mm}
\begin{array}{rcl}
     A & = & \dfrac{5}{gk}(g^3+k^3), \\
     B & = & \dfrac{5}{mnr}((m+n)(m^2g^3-n^2k^3)-(m-n)r^2), \\
     C & = & \dfrac{5}{mng^2k^2r^2}(g^3(m^2g^3-n^2k^3)^2\\
       & -&\left(m(m+n)g^6-(m^2+mn-n^2)g^3k^3+n(m-n)k^6)r^2-k^3r^4\right)), \\
     D & = & \dfrac{g^2}{m^2nk^4r^3}\left(m^2g^3-n^2k^3)^3-(m^2g^3-n^2k^3)(m^2g^3+2n^2k^3)r^2-n^2k^3r^4\right) \\
       &+ & \dfrac{k^2}{mn^2g^4r}\left(m^2g^3(m^2g^3-n^2k^3)-(2m^2g^3+n^2k^3)r^2+r^4\right) \\
       &+ & \dfrac{5(m-n)(g^3-k^3)(m^2g^3-n^2k^3)}{mngkr}-\dfrac{5(m+n)(g^3-k^3)r}{mngk}.
\end{array}
\end{equation*}
Setzt man weiter
\begin{scriptsize}
\begin{equation*}
     \renewcommand{\arraystretch}{2,5}
\setlength{\arraycolsep}{0.5mm}
\begin{array}{rcl}
     T & = & (m^2g^3-n^2k^3)^2-2(m^2g^3+n^2k^3)r^2+r^4 \\
     P & = & \dfrac{(m^2g^3-n^2k^3)^3-(m^2g^3-n^2k^3)(m^2g^3+2n^2k^3)r^2-n^2k^3r^4 +((m^2g^3-n^2k^3)^2-n^2k^3r^2)\sqrt{T}}{2m^2nr^3} \\
     Q & = & \dfrac{(m^2g^3-n^2k^3)^3-(m^2g^3-n^2k^3)(m^2g^3+2n^2k^3)r^2-n^2k^3r^4 -((m^2g^3-n^2k^3)^2-n^2k^3r^2)\sqrt{T}}{2m^2nr^3} \\
     R & = & \dfrac{(m^2g^3-n^2k^3)m^2g^3-(2m^2g^3+n^2k^3)r^2+r^4+(m^2g^3-r^2)\sqrt{T}}{2mn^2r} \\
     S & = & \dfrac{(m^2g^3-n^2k^3)m^2g^3-(2m^2g^3+n^2k^3)r^2+r^4-(m^2g^3-r^2)\sqrt{T}}{2mn^2r},
\end{array}
\end{equation*}
\end{scriptsize}so ist 

\begin{equation*}
    x= \mathfrak{a}\sqrt[5]{\dfrac{g^2}{k^4}P}+\mathfrak{a}^2\sqrt[5]{\dfrac{k^2}{g^4}R}+\mathfrak{a}^3\sqrt[5]{\dfrac{k^2}{g^4}S}+\mathfrak{a}^4\sqrt[5]{\dfrac{g^2}{k^4}Q},
\end{equation*}
wo $\mathfrak{a}$ die erste von $1$ verschiedene $5.$ Einheitswurzel ist. \\

Beispiele dieser allgemeineren Lösungsformel gibt Euler in § 46; im ersten findet er zur vorgelegten Gleichung

\begin{equation*}
    x^5=40x^3+70xx-50x-98
\end{equation*}
die Lösung

\begin{equation*}
    x=\sqrt[5]{-31+3\sqrt{-7}}+\sqrt[5]{-18+10\sqrt{-7}}+\sqrt[5]{-18+-10\sqrt{-7}}
\end{equation*}
\begin{equation*}
    +\sqrt[5]{-31-3\sqrt{-7}}.
\end{equation*}
Im zweiten Exempel gibt Euler für

\begin{equation*}
    x^5=2625x+16600
\end{equation*}
 die Lösung

\begin{equation*}
    x= \sqrt[5]{75(5+4\sqrt{10})}+\sqrt[5]{225(35+11\sqrt{10})}+\sqrt[5]{225(35-11\sqrt{10})}
\end{equation*}
\begin{equation*}
    +\sqrt[5]{75(5-4\sqrt{10})}
\end{equation*}
an.\\

Welche Gleichungen Euler hingegen für beliebiges $n$ lösen konnte, sind die Kreisteilungsgleichungen:

\begin{equation*}
    x^n-1=0,
\end{equation*}
welche bekanntermaßen die Einheitswurzeln 

\begin{equation*}
    x_k=\cos\left(\dfrac{2k\pi}{n}\right)+i\sin \left(\dfrac{2k\pi}{n}\right)
\end{equation*}
als Lösungen besitzen. Es gelang ihm in seiner Arbeit \textit{``De extractione radicum ex quantitatibus irrationalibus"} (\cite{E157}, 1751, ges. 1740) (E157: ``Über das Ziehen von Wurzeln aus irrationalen Größen")  darüber hinaus, die Einheitswurzeln für die Fälle $n=\lbrace 1,2, \cdots, 10 \rbrace$ mithilfe von Radikalen entsprechend der Form (\ref{eq: Euler Vermutung Radikale}) auszudrücken. Bei den $11$--ten Einheitswurzeln muss er hingegen kapitulieren und schreibt im letzten Paragraphen dieser Arbeit: \\

\textit{``Aber die elf Wurzeln der Gleichung $x^{11}-1=0$ können hingegen nur mithilfe einer Gleichung von fünf Dimensionen dargeboten werden, deren Auflösung noch im Verborgenen liegt, weshalb wir an dieser Stelle gezwungen sind aufzuhören."}\\

Bereits an dieser Stelle schimmert Eulers Überzeugung der  Möglichkeit der Auflösung der Gleichung vom $5.$ Grade hervor\footnote{Euler wäre mit der Methode in seiner Arbeit \cite{E157} aus der Gleichung $x^{11}-1=0$ durch die Ersetzung $x+\frac{1}{x}=y$ in dem Ausdruck $1+x+x^2+\cdots +x^{10}=0$ zur einer Gleichung vom Grad $5$ gelangt, führt dies aber nicht aus. Man erhielte die Gleichung $0=1-3y+3y^2-4y^3+y^4+y^5$. Für den Fall $n=7$ zeigt er alle Schritte bis in Detail. Die wesentliche Erkenntnis ist dabei, dass die Gleichung für die Ersetzung von $x$ mit $\frac{1}{x}$ unverändert bleibt. Euler benutzt diese Erkenntnis indirekt, indem er eine Zerlegung $1+x+x^2+x^3+x^4+x^5+x^6=(x^2+px+1)(x^2+rx+1)(x^2+qx+1)$ sucht.}. Sein gesamter Beitrag zu diesem Bereich wird detailliert im Artikel \textit{``Cyclotomy: From Euler through Vandermonde to Gauss"} (\cite{Ne07}, 2007) dargestellt, aber auch der entsprechende Abschnitt im Buch \textit{``Geschichte der Mathematik 1700–1900"}  (\cite{Di85}, 1985) ordnet die Kreisteilungsgleichungen in den historischen Kontext ein, weswegen die entsprechende Diskussion an dieser Stelle ausgespart wird.\\

Gleichwohl Euler noch weit von den grundlegenden Gedanken der Galois--Theorie entfernt war, drängt sich  die Frage auf, zu welchen Kriterien der Möglichkeit der Auflösbarkeit von (\ref{eq: P_n(x)=0}) mittels Radikalen  Euler hätte vordringen können.  Hätte er es etwa vermocht, zu folgendem Satz zu gelangen? 

\begin{Thm}[Auflösbarkeit von Gleichungen 5. Grades mit Radikalen]
    Die Gleichung 

\begin{equation*}
    x^5+ax+b=0
\end{equation*}
ist genau dann mit Wurzeln auflösbar, sofern sie diese Form aufweist:

\begin{equation}
\label{eq: Resolvierbar}
    x^5+\dfrac{5\mu^4(4\nu +3)}{\nu^2+1}x+\dfrac{4\mu^5(2v+1)(4\nu+3)}{\nu^2+1}=0,
\end{equation}
wobei $\mu$ und $\nu$ rationale Zahlen sind.
\end{Thm}

Entsprechende Auflösungsversuche und den erwähnten Lehrsatz findet man respektive in den Arbeiten \textit{``Solution of Solvable Irreducible Quintic Equations, Without the Aid of a Resolvent Sextic"}\footnote{Der scheinbare Widerspruch des Titels von Youngs Arbeit zum Gesagten löst sich auf, wenn man erwähnt, dass Young  explizit nur die (algebraisch) auflösbaren Fälle betrachtet. Für diese, nicht jedoch für den allgemeinen Fall, gelingt es ihm, eine Resolvente vom Grad $4$ zu finden.} (\cite{Yo85}, 1885) und \textit{``Über die auflösbaren Gleichungen von der Form $x^5 + ux + v = 0$"} (\cite{Ru85}, 1885) gezeigt wird. Eine Durchsicht der erwähnten Abhandlungen, insbesondere der zweitgenannten, zeigt die notwendigen gruppentheoretischen Überlegungen, welche Euler natürlich noch fremd waren. Der fundamentale Gedanke  ist, über die symmetrischen Funktionen nachzuweisen, dass die Auflösbarkeit der Gleichung vom 5. Grade mit Wurzelausdrücken zu einer Resolvente vom 6. Grade führt. Man kann folgenden Satz zeigen:

\begin{Thm}[Cayley--Resolvente für Polynome von Grad 5]
    Das Polynom
\begin{equation*}
    y^5+py^3+qy^2+ry+s=0
\end{equation*}
ist genau dann mit Radikalen auflösbar, wenn der Ausdruck (auch Cayley--Resolvente genannt)

\begin{equation*}
    P^2-1024 z \Delta
\end{equation*}
eine rationale Wurzel in $z$ hat. Dabei ist

\begin{small}
\begin{equation*}
   \renewcommand{\arraystretch}{1,5}
\begin{array}{rcl}
     P & = & z^3-z^2(20r+3p^2)-z(8p^2r-16pq^2-240r^2+400sq-3p^4) \\
       & - & p^6 +28p^4r-16p^3q^2-176p^2r^2-80p^2sq \\
       &+ & 224prq^2-64q^4+4000ps^2+320r^3-1600rsq
\end{array}
\end{equation*}
\end{small}
und der Diskriminante

\begin{small}
\begin{equation*}
    \renewcommand{\arraystretch}{1,5}
\begin{array}{rcl}
    \Delta & = & -128p^2r^4+3125s^4-72p^4pqrs+560p^2qr^2s+16p^4r^3+256r^5 +108p^5s^2  \\
       & - & 1600qr^3s+144pq^2r^3-900p^3rs^2+2000pr^2s^2-3750pqs^3+825p^2q^2s^2 \\
       & + & 2250 q^2rs^2 +108q^5s -27q^4r^2-630 pq^3rs+16p^3q^3s -4p^3q^2r^2.
\end{array}
\end{equation*}
\end{small}
\end{Thm}
In diesem Zusammenhang sei des Weiteren auf die Arbeit \textit{``Observatiunculae ad Theoriam Aequationum pertinentes"} (\cite{Ja34}, 1834) (``Kleinere Bemerkungen zur Theorie von Gleichungen") von Jacobi und die Abhandlung \textit{``On a new auxiliary equation in the theory of equations of the fifth order"} (\cite{Ca61}, 1861) Cayley (1821--1895) verwiesen. Während Jacobi noch keine explizite Formel für die Koeffizienten der Resolvente $6.$ Grades angibt\footnote{Jacobi gibt einen der insgesamt drei Koeffizienten an und sagt, für die weiteren benötige man eine ``ein klein wenig längere Rechnung", was eine sehr euphemistische Formulierung ist.}, gelingt Cayley dies\footnote{Für weitere Ausführungen vergleiche man auch den Artikel \textit{``Resolvent sextics of quintic Equations"} (\cite{Di25}, 1925) oder auch den moderneren \textit{``Characterization of Solvable Quintics $x^5 + ax + b$"} (\cite{Sp94}, 1994). Der Artikel \textit{``The Quintic, the Icosahedron, and Elliptic Curves"} (\cite{Ba24}, 2024) fasst die Geschichte der quintischen Gleichung zusammen und betont auch die Verbindungen zur Ikosaedergruppe und elliptischen Kurve.}.\\

Euler auf der anderen Seite schreibt in diesem Zusammenhang in § 14 der Arbeit \cite{E644}: \\

\textit{``Wenn wir die Formen, welche wir für die Wurzeln dieser Gleichungen angegeben haben, genauer betrachten, werden sie alle entdeckt, hervorragend mit jener Vermutung übereinzustimmen, welche ich einst zu äußern gewagt habe, während ich für die Auflösung einer Gleichung beliebigen Grades, in welcher der zweite Term fehlt, wie etwa}

\begin{equation*}
    x^n=px^{n-2}+qx^{n-3}+rx^{n-4}+\text{etc.},
\end{equation*}
\textit{behauptete, dass immer eine auflösende Gleichung von einem Grad weniger gegeben ist:}

\begin{equation*}
    y^{n-1}-Ay^{n-2}+By^{n-3}-Cy^{n-4}+Dy^{n-5}- \text{etc.}=0,
\end{equation*}
\textit{wenn deren insgesamt $n-1$ Wurzeln  $\alpha$, $\beta$, $\gamma$, $\delta$, $\varepsilon$ etc. waren, dass dann stets gilt:}

\begin{equation*}
    x= \sqrt[n]{\alpha}+\sqrt[n]{\beta}+\sqrt[n]{\gamma}+\sqrt[n]{\delta}+\text{etc.}"
\end{equation*}

Die zwei Arbeiten \cite{Yo85} und \cite{Ru85} sowie die anderen genannten \cite{Ja34} und \cite{Ca61} weisen die Euler'sche Vermutung demnach als falsch aus. Es ist indes ein Kuriosum, aus welchem Grunde Euler für die Grade $1$, $2$, $3$ und $4$ Recht behält, jedoch ab dem $5.$ Grad der Resolvente den der Ausgangsgleichung zu übersteigen beginnt. In seiner Arbeit \cite{E644} gibt Euler zum Schluss (§ 15) jedoch bis hin zum $6.$ Grad Bespiele an, die seine Vermutung unterstützen.\\

\begin{figure}
    \centering
      \includegraphics[scale=0.7]{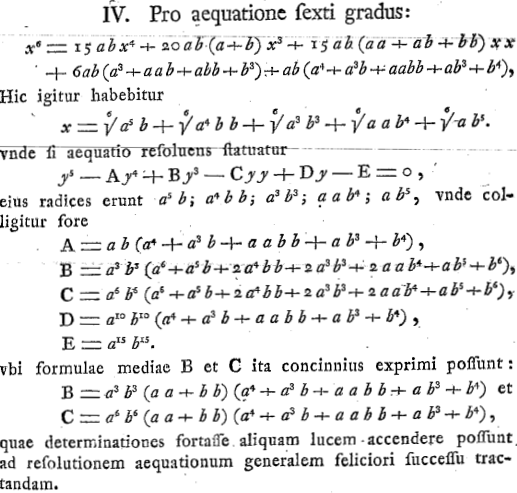}
    \caption{Euler gibt in seiner Arbeit \cite{E644} ein allgemeines Beispiel einer Gleichung vom Grad $6$ an, die eine Resolvente vom Grad $5$ besitzt. Beim ersten Term in der zweiten Zeile sollte wohl $+6ab(a^3+aab+abb+b^3)x$ statt nur $+6ab(a^3+aab+abb+b^3)$ stehen.}
    \label{fig:E644ResolventeEuler}
\end{figure}
Zu bemerken ist hierbei aber, dass es sich sowohl  bei der Euler'schen Vermutung  als auch bei der Aussage bezüglich (\ref{eq: Resolvierbar}) um eine bloße Existenzaussage handelt, sodass dieses Ergebnis Euler seinem Ziele nur einen Schritt näher gebracht hätte. Man sieht sich demnach insgesamt zu dem Schluss veranlasst, dass Euler ein allgemeines Kriterium wegen der fehlenden gruppentheoretischen Lehrsätze verborgen bleiben musste, obschon er in \cite{E282} und \cite{E644}  zu speziellen Fällen des obigen Lehrsatzes gelangt. Zu den Überlegungen von Jacobi aus \cite{Ja34} wäre er zweifelsohne in der Lage gewesen, scheint aber im Gegensatz zu seinem Nachfolger die Möglichkeit er höhergradigen Resolvente nicht zugelassen zu haben, sodass er sich den Weg zu dieser Erkenntnis wegen der falschen Ausgangsfrage bereits versperrt hatte.

\paragraph{Ein Zugang über Differentialgleichungen}
\label{para: Differentialgleichung}

Es ist typisch für Euler, die Auflösung von (\ref{eq: P_n(x)=0}) nicht nur auf die ``klassische"{} Weise mit Radikalen zu versuchen, sondern für dasselbe Unterfangen auch die Methoden der Analysis zu gebrauchen. In § 24 der Arbeit \textit{``Observationes circa radices aequationum"} (\cite{E406}, 1771) (E406: ``Beobachtungen über die Wurzeln von Gleichungen") stellt er die Lösung der kubischen Gleichung 

\begin{equation*}
    y^3+py+q=0
\end{equation*}
vermöge einer Differentialgleichung vor. Dazu möchte Euler sie auf die Differentialgleichung

\begin{equation*}
    ddy+Qdy+Ry=0
\end{equation*}
mit noch zu bestimmenden Funktionen $Q$ und $R$ in den Variablen $p$ und $q$ überführen.  Er findet nach zweimaliger Differentiation der Ausgangsgleichung und anschließender Vereinfachung:

\begin{equation*}
    Q= \dfrac{18p^2qdp^2-2(8p^3-27q^2)dpdq-54pqdq^2}{(3qdp-2pdq)(4p^3+27^2)}+\dfrac{2pddq-3qddp}{3qdp-2pdq},
\end{equation*}
\begin{equation*}
    R=\dfrac{6p(dq^3+pdp^2dq-qdp^3)}{(3qdp-2pdq)(4p^3+27qq)}+\dfrac{dqddp-dpddq}{3qdp-2pdq}.
\end{equation*}
Setzt man $q^2=\frac{4p^3x}{27}$, wird die angenommene Differentialgleichung

\begin{equation*}
    ddy-dy\left(\dfrac{ddx}{dx}+\dfrac{dp}{p}-\dfrac{dx}{2x}-\dfrac{dx}{2(1+x)}\right)
\end{equation*}
\begin{equation*}
    +y\left(\dfrac{dpddx}{2pdx}-\dfrac{ddp}{2p}+\dfrac{3dp^2}{4pp}-\dfrac{dpdx}{4px}-\dfrac{dpdx}{4p(1+x)}-\dfrac{dx^2}{36x(1+x)}\right)=0.
\end{equation*}
Um sie integrierbar zu machen, multipliziert Euler diese mit

\begin{equation*}
    \dfrac{x(1+x)}{p^2dx^2}(2pdy-ydp)
\end{equation*}
und integriert sie zu

\begin{equation*}
    \dfrac{x(1+x)}{pdx^2}\left(dy-\dfrac{ydp}{2p}\right)^2=\dfrac{C}{36}+\dfrac{y^2}{36p}
\end{equation*}
mit einer beliebigen Konstante $C$. Für  $y=z\sqrt{p}$ wird die letzte Gleichung zu

\begin{equation*}
    \dfrac{6dz}{\sqrt{C+z^2}}=\dfrac{dx}{\sqrt{x(1+x)}}.
\end{equation*}
Beide Seiten können uneingeschränkt mithilfe von Logarithmen integriert werden. Euler gibt direkt die Lösung

\begin{equation*}
    (z+\sqrt{C+z^2})^6= D\left(\dfrac{1}{2}+x+\sqrt{x(1+x)}\right)=\dfrac{1}{2}D\left(\sqrt{x}+\sqrt{1+x}\right)^2
\end{equation*}
mit einer konstanten Größe $D$. Daher erhält er für $z$

\begin{equation*}
    z=\dfrac{y}{\sqrt{p}}=A\left(\sqrt{x}+\sqrt{1+x}\right)^{\frac{1}{3}}+B \left(\sqrt{x}-\sqrt{1+x}\right)^{\frac{1}{3}},
\end{equation*}
wobei $A$ und $B$ beliebige Größen sind. Der Kubus gibt schließlich

\begin{equation*}
    z^3=-3ABz+(A^3+B^3)\sqrt{x}+(A^3-B^3)\sqrt{1+x}.
\end{equation*}
Mit diesem Ergebnis, gerade die Cardano'sche Formel zur Lösung von kubischen Gleichungen, sofern $p$ und $q$ wieder eingesetzt werden, schließt Euler seine Abhandlung \cite{E406}. \\

Die Euler'sche Vorgehensweise ist höchst bemerkenswert und verdient gewiss eine weitere Untersuchung, welche jedoch an dieser Stelle nicht unternommen werden kann. Einschränkend ist zu festzuhalten, dass Euler im Wissen der korrekten Lösungsformel geleitet die entsprechenden Substitutionen zur Integration der Differentialgleichungen durchführt.  Es wäre gewiss ebenfalls lohnenswert, die Euler'schen Einsichten -- unter Kenntnis der Lösungen  polynomialer Gleichungen  höherer Grade -- weiter auszuarbeiten. Man würde gewiss zum oben  erwähnten Satz, welcher die Bedingungen zur Auflösbarkeit mit Radikalen beschreibt, geführt werden und vielleicht zum Satz von Hermite (1822--1901), wie man diese Gleichungen mit elliptischen Funktionen aufzulösen vermag\footnote{Hermite beweist diesen Satz in seiner Arbeit \textit{``Sur la résolution de l’Équation du cinquiéme degré "} (\cite{He58}, 1858) (``Über die Auflösung der Gleichung vom 5. Grade").}. 

\paragraph{Lambert'sche Reihen}
\label{para: Lambert'sche Reihen}

Ein weiterer Auflösungsversuch Eulers von (\ref{eq: P_n(x)=0}) besteht in der Lambert'schen Reihe, wie Euler sie etwa in \cite{E532} nennt, gleichwohl er schon bei der Betrachtung  in der Arbeit \cite{E406} zu ihr gelangt war. Ihren Ursprung hat die Reihe in der Gleichung

\begin{equation*}
    1=\dfrac{A}{x}+\dfrac{B}{x^2}+\dfrac{C}{x^3}+\dfrac{D}{x^4}+\text{etc.},
\end{equation*}
wobei sich die Koeffizienten  beliebig wählen lassen. In der Arbeit \textit{``Analysis facilis et plana ad eas series maxime abstrusas perducens, quibus omnium aequationum algebraicarum non solum radices ipsae sed etiam quaevis earum potestates exprimi possunt"} (\cite{E631}, 1789, ges. 1776) (E631: ``Eine leichte und verständliche Analysis, welche zu den höchst bizarren Reihen führt, mit welchen nicht nur die Wurzeln der Gleichungen selbst, sondern auch jede Potenz derselben ausgedrückt werden können") betrachtet er in § 26 die verwandte Gleichung

\begin{equation*}
    1-\dfrac{1}{Z^{\alpha}}=\dfrac{B}{Z^{\beta}}+\dfrac{C}{Z^{\gamma}}+\dfrac{D}{Z^{\delta}}+\dfrac{E}{Z^{\varepsilon}}+\text{etc.}.
\end{equation*}
 Für $Z^n$ leitet er in Abhängigkeit der völlig beliebigen Potenzen $\alpha$, $\beta$, $\cdots$ und der Koeffizienten $A$, $B$, $\cdots$ und natürlich $n$ eine Reihenentwicklung her, für deren allgemeinen Term er das Bildungsgesetz in Worten beschreibt\footnote{Das von Euler beschriebene Bildungsgesetz ließe sich über äußerst verwickelte Rechnungen mit der von Lagrange in seiner Arbeit \textit{``Nouvelle méthode pour résoudre les équations littérales par le moyen des séries"} (\cite{La68}, 1770, ges. 1768) (``Eine neue Methode zur Auflösung von Gleichungen mithilfe von Reihen") heute nach ihm benannten Inversionsformel deuten. Langrange beweist, dass eine Gleichung $z=f(x)$ mit $f'(a)\neq 0$  die formale Inverse 
\begin{equation*}
    x=g(z)=a+\sum_{n=1}^{\infty} g_n \dfrac{(z-f(a))^n}{n!}
\end{equation*}
mit

\begin{equation*}
    g_n = \lim_{z \rightarrow a} \dfrac{d^{n-1}}{dz^{n-1}} \left(\dfrac{z-a}{f(z)-f(a)}\right)^n
\end{equation*}
besitzt.}. Die vollkommen allgemeine Formel wendet Euler indes nicht über drei Terme hinaus an. Es sei an dieser Stelle Eulers zwecks einer Illustration due allgemeine Formel (sein Theorema generale) aus § 11 von \cite{E631} erwähnt. Ist die Gleichung

\begin{equation*}
    1= \dfrac{A}{x^{\alpha}}+\dfrac{B}{x^{\beta}}
\end{equation*}
vorgelegt, so gilt

\begin{equation*}
    \renewcommand{\arraystretch}{2,0}
     \setlength{\arraycolsep}{0.5mm}
     \begin{array}{rcl}
          x^n & = & A^{\frac{n}{\alpha}}+ A^{\frac{n-\beta}{\alpha}}B \cdot \dfrac{n}{\alpha}+ A^{\frac{n-2\beta}{\alpha}}B^2  \cdot \dfrac{n}{\alpha} \cdot \dfrac{n+\alpha -2\beta}{2\alpha} \\
          &+ & A^{\frac{n-3\beta}{\alpha}}B^3 \cdot \dfrac{n}{\alpha} \cdot \dfrac{n+\alpha -3\beta}{2\alpha} \cdot \dfrac{n+2\alpha -3\beta}{3\alpha} \\
          &+ & A^{\frac{n-4\beta}{\alpha}}B^4 \cdot \dfrac{n}{\alpha} \cdot \dfrac{n+\alpha -4\beta}{2\alpha} \cdot \dfrac{n+2\alpha -4\beta}{3\alpha} \cdot \dfrac{n+3\alpha -4\beta}{4\alpha} \\
          &+ & A^{\frac{n-5\beta}{\alpha}}B^5 \cdot \dfrac{n}{\alpha} \cdot \dfrac{n+\alpha -5\beta}{2\alpha} \cdot \dfrac{n+2\alpha -5\beta}{3\alpha} \cdot \dfrac{n+3\alpha -5\beta}{4\alpha} \cdot \dfrac{n+4\alpha -5\beta}{5\alpha} \\
          &+ &\text{etc.}
     \end{array}
\end{equation*}

 Eulers Beispiel 4 aus § 15 von \cite{E631} ist die Gleichung

\begin{equation*}
    1=\dfrac{A}{xx}+\dfrac{B}{x},
\end{equation*}
welche man natürlich direkt mit

\begin{equation*}
    x=\dfrac{B+\sqrt{BB+4A}}{2}
\end{equation*}
auflösen kann. Vermöge der angesprochenen Reihenentwicklung  findet Euler indes zunächst, dass  für die Potenz $x^n$ gelten müsste:

\begin{equation*}
    (b+\sqrt{b^2+a^2})^n= a^n+\dfrac{n}{2}a^{n-1}\cdot 2b+\dfrac{n}{2}\cdot \dfrac{n}{4}{n-2}\cdot 4bb
\end{equation*}
\begin{equation*}
    +\dfrac{n}{2}\cdot \dfrac{n-1}{4}\cdot \dfrac{n+1}{6}a^{n-3}\cdot 8b^3+\dfrac{n}{2}\cdot \dfrac{n-2}{4}\cdot \dfrac{n}{6}\cdot \dfrac{n+2}{8}a^{n-4}\cdot 16b^4+\text{etc.}
\end{equation*}
mit $A=a^2$ und $2b=B$; und daher für $n=1$:

\begin{equation*}
    \sqrt{b^2+a^2}=a+\dfrac{1}{2}\cdot \dfrac{b^2}{a}-\dfrac{1\cdot 1}{2\cdot 4}\cdot \dfrac{b^4}{a^3}+\dfrac{1\cdot 1\cdot 3}{2\cdot 4 \cdot 6}\cdot \dfrac{b^6}{a^5}-\text{etc.}.
\end{equation*}
Dies bestätigt die Richtigkeit der mit der Reihe gefundenen Lösung, wenn man den binomischen Lehrsatz (\ref{eq: Binomischer Lehrsatz}) zum Vergleich heranzieht.


\subsection{Grenzen durch den eigenen Arbeitsethos}
\label{subsec: Grenzen durch den eigenen Arbeitsethos}

\epigraph{Mathematical discoveries, like springtime violets in the woods, have their season which no human can hasten or retard.}{Carl Friedrich Gauß}

Die Arbeitsweise eines jeden Mathematikers ist von einer ihm eigenen Philosophie geleitet, welche selbigen eine bestimmte Sicht auf jeweilige Problemstellungen entwickeln lässt\footnote{Dies ist ein in der Psychologie umfassend untersuchtes Phänomen und unter dem Überbegriff \textit{Einstellungseffekt} bekannt. Dieser Effekt besagt, dass man von einer Lösungsmethode mit einer Wahrscheinlichkeit umgekehrt proportional zur mit ihr erfolgreich gelösten Probleme wieder abweichen wird, \textit{obwohl} die präferierte zuvor gewinnbringende Methode nicht zwingend das geeignetste Verfahren für eine andere, gegenwärtig zu bewältigende, ähnliche Aufgabe darstellt. Man konsultiere etwa die Arbeit \textit{``Mechanization in problem solving: The Effect of Einstellung"} (\cite{Lu42}, 1942) und ähnliche für konkrete Beispiele zu diesem Gegenstand.}. Es liegt in der Natur der Dinge, dass die subjektive Arbeitsweise einen bei den einen Problemstellungen unerwartet weit vorzudringen gestattet, bei anderen hingegen zur einer vollkommenen Stagnation des Fortschritts führen kann. Während der Großteil der vorliegenden Abhandlung Beispiele für den ersten Fall im Kontext von Eulers Arbeiten besprochen hat, sollen  nun solche Begebenheiten den Diskussionsgegenstand bilden, in welchen ein tieferes Verständnis allein durch den Arbeitsethos unmöglich gemacht worden ist. So wird gezeigt werden (Abschnitt \ref{subsubsec: Komplexe Analysis}), dass Euler der Weg zur Theorie von Funktionen einer komplexen Variable trotz Entdeckungen einiger ihrer Grundbausteine -- wie den Riemann--Cauchy--Differentialgleichungen -- verwehrt bleiben musste. Weiter werden Beispiele für das ``Missglücken"{} der von Euler präferierten Methodus inveniendi angeführt (Abschnitt \ref{subsubsec: Methodus Inveniendi über Methodus Demonstrandi}), bevor schließlich zu möglichen Erkenntnissen übergegangen wird, die praktischen Anwendungen weichen mussten (\ref{subsubsec: Praxis über Abstraktion}).

\subsubsection{Eulers Nichtentdeckung der komplexen Analysis}
\label{subsubsec: Komplexe Analysis}

\epigraph{The shortest path between two truths in the real domain passes through the complex domain.}{Jacques Hadamard}

Euler war zumeist bemüht, Fragestellungen aus anderen Bereichen auf die reelle Analysis zu reduzieren, wo die Betonung auf dem Wort \textit{reell} liegt. Es wird sich nämlich zeigen, dass durch diesen Ansatz Euler die Theorie einer Funktion einer \textit{komplexen} Variable nicht zur vollen Blüte entwickeln konnte. Obwohl, insbesondere aus moderner Sicht betrachtet, Euler wesentliche Beiträge zu diesem Teilbereich geleistet hat\footnote{Das Buch \textit{``Ostwalds Klassiker der exakten Wissenschaften Band 261: Zur Theorie komplexer Funktionen von Leonhard Euler"} (\cite{OK07}, 2007) ist eigens Beiträgen dieser Art gewidmet.}, sind dennoch Hürden, die Euler nicht überwinden konnte, gemeinsam mit den Ursachen ihres Auftretens deutlich erkennbar\footnote{So mag es begründet sein, warum das Buch über die Historie der komplexen Analysis \textit{``The Real and the Complex: A History of Analysis in the 19th Century"} (\cite{Gr15}, 2015) von Gray (1947--) nicht mit Euler beginnt, sondern mit Lagrange. Auch das in diesem Thema noch umfassender gewidmete Buch \textit{``The Higher Calculus: A History of Real and Complex Analysis from Euler to Weierstrass"} (\cite{Bo86}, 1986) von Bottazini (1947--) behandelt den Euler'schen Beitrag vergleichsweise kurz.}. 

\paragraph{Euler und die Riemann-Cauchy'schen Differentialgleichungen}
\label{paragraph: Die Riemann-Cauchy'schen Differentialgleichungen}

Einer der grundlegenden Errungenschaften der Funktionentheorie sind die Riemann--Cauchy--Differentialgleichungen. Selbige sind Euler nicht entgangen, wenn er auch deren Tragweite nicht erfasst hat\footnote{Die Erstentdeckerschaft gebührt allerdings d'Alembert (1717--1783), der sie in seinem Buch \textit{``Essai d’une nouvelle théorie de la résistance des fluides"} (\cite{DA52}, 1752) (``Essay über eine neue Theorie des Widerstands von Fluiden") niedergeschrieben hat.}. Dies war seinen  Nachfolgern, vor allem {Riemann}\footnote{Riemann stellt die Theorie einer Funktion einer komplexen Variable in seiner Dissertation \textit{``Grundlagen für eine allgemeine Theorie der Funktionen einer veränderlichen komplexen Grösse"} (\cite{Ri51}, 1851) vor.}  und {Cauchy}\footnote{Cauchy gibt die nach ihm benannten Differentialgleichungen zum ersten Mal in seiner Arbeit \textit{``Mémoire sur les intégrales définies, prises entre des limites imaginaires"} (\cite{Ca14}, 1825, ges. 1814) (``Abhandlung über bestimmte Integrale, erstreckt zwischen imaginären Grenzen") an.}, vorbehalten. Man kann sie wie folgt formulieren:  Hat man eine differenzierbare Funktion $f$ der variablen Größe $z$ und schreibt diese als $z=x+iy$ und setzt $f(x,y)=u(x,y)+iv(x,y)$, dann gelten die folgenden Differentialgleichungen:

\begin{equation}
\label{eq: Riemann-Cauchy-DGL}
    \dfrac{\partial u}{\partial x}= \dfrac{\partial v}{\partial y} \quad \text{und} \quad  \dfrac{\partial u}{\partial y}=-\dfrac{\partial v}{\partial x}.
\end{equation}
Jede über den komplexen Zahlen differenzierbare Funktion erfüllt \textit{notwendig} die obigen Differentialgleichungen, was ihr Aufkommen in Eulers Untersuchungen weniger überraschend macht; man findet sie in der Abhandlung \textit{``De integrationibus maxime memorabilibus ex calculo imaginariorum oriundis"} (\cite{E656}, 1793, ges. 1777) (E656: ``Über höchst merkwürdige Integrationen, die aus dem Kalkül mit imaginären Größen entspringen"), direkt im ersten Paragraphen. Euler gelangt dabei, ebenso wie später in § 1--2 der Arbeit \textit{``Ulterior disquitio de formulis integralibus imaginariis"} (\cite{E694}, 1797, ges. 1777) (E694: ``Eine weitere Untersuchung über imaginäre Differentialformen") noch einmal, wie folgt zu (\ref{eq: Riemann-Cauchy-DGL}): \\

\begin{figure}
    \centering
 \includegraphics[scale=0.9]{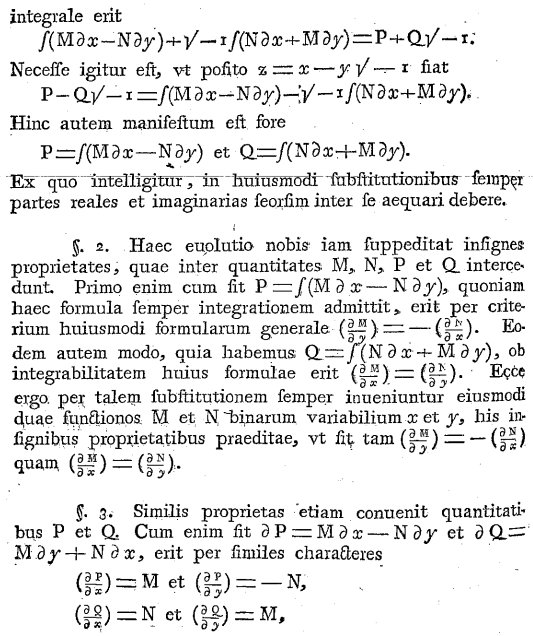}
    \caption{Eulers Herleitung der Riemann--Cauchy Differentialgleichung aus seiner Arbeit \cite{E694} über die Theorie von Funktionen zweier reeller Variablen.}
    \label{fig:E694Riemann-Cauchy}
\end{figure}

Er schreibt in erstgenannter Arbeit zunächst:

\begin{equation*}
    \Delta:(x+\sqrt{-1}\cdot y)= A:(x+\sqrt{-1}\cdot y)+\sqrt{-1}\cdot B:(x+\sqrt{-1}\cdot y),
\end{equation*}
wobei Euler statt dem heute geläufigem $i$ hier wieder $\sqrt{-1}$ verwendet\footnote{Das Symbol $i$ für $\sqrt{-1}$ hat in der Tat ebenfalls Euler eingeführt und zwar in § 2 seiner Abhandlung  \textit{ ``De formulis differentialibus angularibus maxime irrationalibus, quas tamen per logarithmos et arcus circulares integrare licet"} (\cite{E671}, 1794, ges. 1777) (E671: ``Über höchst irrationale angulare Differentialformeln, welche sich dennoch durch Logarithmen und Kreisbogen integrieren lassen"). Jedoch verwendet Euler in seinen nachfolgenden Werken den Buchstaben $i$ nicht, um $\sqrt{-1}$ anzuzeigen.}. Euler betrachtet somit die Funktion nicht als $\Delta: \mathbb{C}\mapsto \mathbb{C}$, sondern als $\Delta: \mathbb{R}^2 \mapsto \mathbb{R}^2$, wodurch $\sqrt{-1}$ auf den Rang einer weiteren Konstante ohne eine tiefere Bedeutung herabgestuft wird. Auf die Funktionen, welche Euler $M$ und $N$ nennt, wendet Euler dann nun einen Satz aus der Integrationstheorie an. Er nimmt an

\begin{equation*}
    \Delta: z = \int Zdz,
\end{equation*}
dass also $Zdz$ ein vollständiges Differential ist. Nun schreibt er:

\begin{equation*}
    Zdz = (M+N\sqrt{-1})(dx+\sqrt{-1}dy)=(Mdx-Ndy)+\sqrt{-1}(Mdy+Ndx).
\end{equation*}
Nun, so argumentiert Euler, ist $Zdz$ ein vollständiges Differential, somit muss auch die rechte Seite integrierbar sein. Da $\int Zdz$ ebenfalls wieder in Real- und Imaginärteil zerlegt werden kann, so müssen auch der Real- und Imaginärteil auf der rechten Seite  integrierbar sein. Daraus ergibt sich aus der Integrationstheorie die Implikation

\begin{equation*}
    \dfrac{\partial M}{\partial y}=-\dfrac{\partial N}{\partial x}
\end{equation*}
aus dem Realteil, analog folgt aus dem Imaginärteil die Gleichung

\begin{equation*}
    \dfrac{\partial M}{\partial x}=\dfrac{\partial N}{\partial y},
\end{equation*}
was  gerade die Riemann-Cauchy'schen Differentialgleichungen (\ref{eq: Riemann-Cauchy-DGL}) sind. Sie folgen also, wie bereits angedeutet, notwendig aus der Annahme der Differenzierbarkeit einer Funktion einer komplexen Variablen. \\

Man findet die Riemann--Cauchy Gleichungen in anderem Zusammenhang bereits in der Arbeit \textit{``Considerationes de traiectoriis orthogonalibus"} (\cite{E390}, 1770, ges. 1768) (E390: ``Betrachtungen zu orthogonalen Trajektorien"), welche der Untersuchung von orthogonalen Trajektorien, also einem geometrischen Problem, gewidmet ist.  Hier gelangt Euler über die Frage, welche Gleichungen zwei Kurven $t(p,q)=t$, $u(p,q)=u$, die sich orthogonal zu schneiden haben, erfüllen müssen, zu der Bedingung:

\begin{equation*}
    \dfrac{\partial t}{\partial p}\dfrac{\partial t}{\partial q}+\dfrac{\partial u}{\partial p}\dfrac{\partial u}{\partial q}=0.
\end{equation*}
Schnell kann man nachprüfen, dass diese äquivalent zu den Riemann-Cauchy--Differentialgleichungen sind. Obwohl diese Tatsache den Weg zur geometrischen Interpretation der komplexen Zahlen ebnet, scheint sie von Euler nicht erkannt worden zu sein\footnote{Diese Ansicht wird von Wußing (1927--2011) und Juschkewitsch (1906--1993),  den Autoren von \textit{``Ostwalds Klassiker der exakten Wissenschaften Band 261: Zur Theorie komplexer Funktionen von Leonhard Euler"} (\cite{OK07}, 2007),  geteilt. Man konsultiere dazu das einleitende Vorwort.}.

\paragraph{Eulers Verhältnis zu den komplexen Zahlen}

Die bahnbrechende Leistung Cauchys und Riemanns besteht demnach weniger  in der Derivation der Gleichungen, sondern vielmehr im Nachweis, dass umgekehrt (\ref{eq: Riemann-Cauchy-DGL}) auch bereits die Differenzierbarkeit von komplexen Funktionen impliziert, also insgesamt Differenzierbarkeit in der Variable $z=x+iy$ und die Riemann--Cauchy--Differentialgleichungen äquivalent sind,  ein Faktum, welches Euler entgangen ist\footnote{Angesichts seiner Meinung (Abschnitt \ref{subsubsec: Durch einen neuen Gedanken: Die Theta-Funktion}), man könne jede Funktion als verallgemeinerte Potenzreihe auffassen, kann man die Ansicht vertreten, Euler würde diese Frage sogar gar nicht gestellt haben.}.\\

Den Hauptgrund für Eulers Übersehen der Implikation von (\ref{eq: Riemann-Cauchy-DGL}) entnimmt man indirekt aus § 14 seines Papiers \textit{``Integratio formulae differentialis maxime irrationalis, quam tamen per logarithmos et arcus circulares expedire licet"} (\cite{E689}, 1795, ges. 1777) (E689: ``Integration einer höchst irrationalen Integralformel, welche sich dennoch mit Logarithmen und Kreisbogen erledigen lässt"). In dieser Arbeit schreibt er seiner Lösung des gestellten Problems eine Sonderstellung zu, weil sie das Rechnen mit imaginären Größen \textit{verlangt}\footnote{Er sieht sich jedenfalls außer Stande, einen direkten Weg, womit er wohl einen über die reelle Analysis im Sinn hat, anzugeben.}. Dementsprechend sind für Euler die komplexen Zahlen, wie auch schon oben angeklungen ist, ein reines Hilfsmittel für die reelle Analysis und bedürfen keiner eigenen Untersuchung. Von dieser Auffassung zeugt auch die Arbeit \textit{``Nova methodus integrandi formulas differentiales rationales sine subsidio quantitatum imaginariarum"} (\cite{E572}, 1784, ges. 1775) (E572: ``Eine neue Methode rationale Differentialformeln ohne Hilfe von imaginären Größen zu integrieren"). Komplexe Größen werden also, wenn möglich, sogar umgangen, ein Phänomen, was auch schon oben (Abschnitt \ref{subsubsec: Differentialgleichungen unendlicher Ordnung: Der homogene Fall}) bei Differentialgleichungen aufgetreten ist. Sie mögen zum Auffinden neuer Sachverhalte sehr von Nutzen sein, aber endgültig können alle Sachverhalte über reelle Funktionen auch hinreichend leicht und bequem formuliert werden, ohne das Reich der reellen Funktionen verlassen zu müssen. Weiter lässt sich seine Auffassung der Rolle komplexer Zahlen aus seinem Lehrbuch zur Algebra \cite{E387} entnehmen; hier schreibt er in altertümlicher Sprache in  § 145 zu den Wurzeln negativer Zahlen:\\

\textit{``Dahero bedeuten alle diese Ausdrücke $\sqrt{-1}$, $\sqrt{-2}$, $\sqrt{-3}$, $\sqrt{-4}$, etc. solche ohnmögliche oder Imaginäre Zahlen, weil dadurch Quadratwurzeln von negativen Zahlen angezeigt werden. Von diesen behauptet man also mit allem Recht daß sie weder größer noch kleiner sind als nichts; und auch nicht einmal nichts selbsten, als aus welchem Grund sie folglich für ohnmöglich gehalten werden müßen.}\\

Eine weitere Bestätigung dieser Einstellung gegenüber den imaginären Zahlen  gibt Aussage am Ende von § 3 von \cite{E807}.  Hier schreibt Euler nämlich bezüglich der Leibniz'schen Auffassung, die Logarithmen können  für negative Zahlen nicht definiert werden:\\

\textit{``Ich muss indes gestehen, dass diese Antwort, wenn sie korrekt wäre, das Fundament der ganzen Analysis erschüttern würde, welches insbesondere in der Allgemeinheit der Regeln und der Operationen besteht, welche als wahr eingeschätzt worden sind, von welcher Natur auch immer man die Größen zu sein annimmt, auf welche sie angewendet werden."}\\

Übertragen auf die komplexen Zahlen ließe sich das etwa so deuten, dass eine Funktion betrachtet als $f:\mathbb{C} \rightarrow \mathbb{C}$ keine anderen Eigenschaften haben kann als eine betrachtet als $f:\mathbb{R}^2 \rightarrow \mathbb{R}^2$. Dass diese Auffassung unrichtig ist, ist heute wohl bekannt. So ist eine Funktion, sofern sie nach $z$ differenzierbar ist, auch nach $x$ und $y$ mit $z=x+iy$ differenzierbar. Die Umkehrung dieser Aussage ist indes nicht korrekt, wie das Beispiel der Funktion $f(z)= \overline{z}= x-iy$ durch das Nichterfüllen von (\ref{eq: Riemann-Cauchy-DGL}) zeigt. \\ 

So überrascht es auch nicht, dass weitere Grundpfeiler der Theorie  von Funktionen komplexer Variablen Euler verborgen bleiben mussten. Vor allem anderen ist hier der Residuensatz zu nennen, welchen man wohl überhaupt nur zu formulieren gedenkt, sofern man zahlreiche Integrationen in der komplexen Ebene entlang von Wegen betrachtet hat.  Euler scheint nicht einmal explizit Wegintegrale in der \textit{reellen} Ebene  $\mathbb{R}^2$ betrachtet zu haben. Doppelintegrale scheint er um ihre selbst willen lediglich in der einen einzigen Abhandlung \cite{E391} einer eingehenden Untersuchung unterworfen zu haben.

\subsubsection{Methodus Inveniendi über Methodus Demonstrandi}
\label{subsubsec: Methodus Inveniendi über Methodus Demonstrandi}

\epigraph{The art of being wise is the art of knowing what to overlook.}{William James}

Die Frage nach dem Finden einer Lösung setzt implizit ihre Existenz  voraus. Umgekehrt klammert dies die Frage nach der Existenz aus. Dies zieht, insbesondere vor dem Hintergrund eines modernen Verständnisses von mathematischer Strenge, unter Umständen etwaige Probleme nach sich, von welchen nachstehend einige angeführt werden sollen.

\paragraph{Existenzfragen bei Euler im Allgemeinen}
\label{para: Existenzfragen}

Allen voran ist zu erwähnen, dass das vollständige Umgehen von Existenzfragen in der Mathematik es Euler überhaupt erst ermöglicht hat, gewisse Schätze zu heben. Seine \textit{Methodus} \cite{E65}, welche die Variationsrechnung zu einem eigenständigen Bereich der Mathematik erhoben hat, ist eines der eindrücklichsten Beispiele dessen. Es sei noch einmal das Zitat aus  § 1 des ersten Anhangs zu diesem Buch wiederholt:\\

\textit{``Weil ja alle Wirkungen der Natur einem gewissen Gesetz von Maximum oder Minimum folgen, besteht kein Zweifel, dass bei den Kurven, welche Körper beschreiben, wenn sie von irgendwelchen Kräften angegriffen werden, eine gewisse Eigenschaft des Maximums oder Minimums auftritt."}\\

Es es freilich schwerlich vorstellbar, wie  sich ohne die in obiger Aussage enthaltene Ansicht der Existenz einer Variation überhaupt der Grundstein ihres Kalküls legen ließe.  Dieser Behauptung verleiht die Bemerkung von Speiser im Vorwort von Band 29 der Serie 1 der \textit{Opera Omnia} (\cite{OO29}, 1956) Nachdruck. Er schreibt dort bezüglich des Prozesses zum Auffinden einer Variation (pp. IX--X):\\

\textit{``Der Übergang vom Variationsproblem zur Differentialgleichung gehört der Methodus inveniendi, der produktiven Phantasie, an. Wenn sich der Beweis daran versucht, so verstrickt er sich in Tautologien. So steht in einem Lehrbuch ein Satz von folgender Art: Es gibt eine Funktion, welche die Variation zum Verschwinden bringt. In einer Anmerkung wird aber gesagt: Dieser Satz kann nur bewiesen werden, wenn man für die Funktion die Existenz einer Ableitung voraussetzt. Das heißt offenbar: Um die Existenz der Funktion zu beweisen, muß man ihre Existenz voraussetzen."}\\

Trotz der Fruchtbarkeit dieser Herangehensweise, welche schon in der \textit{Methodus} \cite{E65} in den Euler--Lagrange--Gleichungen der Variationsrechnung und dem in der Physik seither nicht mehr wegzudenkenden Wirkungsprinzip\footnote{Hamilton (1805--1865) hat in seinen Arbeiten in expliziter Analogie in der Optik das nach ihm benannte Hamilton'sche Extremalprinzip formuliert. Man konsultiere etwa seine Arbeit \textit{``On the application to dynamics of a general mathematical method previously applied to optics"} (\cite{Ha34}, 1834). Sein Prinzip soll später Schrödinger (1887--1961) in der Formulierung der Quantenmechanik leiten. Er erwähnt Hamiltons Arbeit in seiner eigenen \textit{``An undulatory theory of the mechanics of atoms and molecules"} (\cite{Sc26}, 1926). In der Quantenfeldtheorie ließe sich das Feynman'sche Pfadintegral von Feynman (1918--1988), aber noch vielmehr das Quantenwirkungsprinzip von Schwinger (1918--1994) nennen. Letzterer führt sein Variationsprinzip in der Arbeit \textit{``On Theory of quantized fields I"} (\cite{Sc51}, 1951) ein. Feynman erklärt das nach ihm benannte Integral in seiner Arbeit \textit{``Space-Time Approach to Non-Relativistic Quantum Mechanics"} (\cite{Fe48}, 1948). Sowohl Feynman als auch Schwinger haben dabei von der Dirac'schen Arbeit \textit{``The Lagrangian in Quantum Mechanics"} (\cite{Di33}, 1933) ausgehend ihre Untersuchungen begonnen.} ihre Höhepunkte gefunden haben, sollen auch die aufgrund dessen von Euler nicht tangierten Aspekte des Variationskalüls nicht verschwiegen werden. Sind Euler selbst die Grenzen des Anwendungsbereiches der Euler--Lagrange--Differentialgleichungen in seiner Arbeit \textit{``De insigni paradoxo, quod in analysi maximorum et minimorum occurrit"} (\cite{E735}, 1811, ges. 1779) (E735: ``Über ein außergewöhnliches Paradoxon, welches bei der Analysis von Maxima und Minima auftritt") noch vor Augen getreten -- das Paradoxon besteht dabei darin, dass die Lösungsfunktion zu dem Variationsproblem von Interesse nicht von den Euler--Lagrange--Gleichungen erfasst wird\footnote{Euler betrachtet hier explizit das aus physikalischen Überlegungen herstammende Funktional $\int\limits_{a}^{b}\sqrt{x}\sqrt{1+(y'(x))^2}dx$. Die zugehörigen Euler--Langrange--Gleichung -- sie lautet in diesem Fall $0=\frac{d}{dx}\sqrt{x}\frac{y'}{\sqrt{1+(y')^2}}$ -- werden von $y(x)=2c\sqrt{x-c^2}+d$ mit beliebigen Konstanten $c,d$ gelöst. Jedoch, wie Euler auch in § 7 bemerkt, gibt es unzählige Funktionen, die das Funktional kleiner werden lassen als diese Lösung, wovon Euler auch explizite Beispiele vorrechnet. Die lokale Natur des durch die Euler--Lagrange--Gleichungen eruierten Minimums erläutert er im letzten Paragraphen der Arbeit. Euler streift hier somit den Unterschied zwischen den heute sogenannten \textit{schwachen} und \textit{starken Extrema} bei Funktionalen, hatte aber selbstredend noch kein Konzept davon gebildet. Die Euler--Lagrange'schen Differentialgleichungen erlauben, bei vorausgesetzter zweimaliger stetiger Differenzierbarkeit, lediglich ein schwaches Extrema zu finden.} -- und vermochte er durch explizites Einsetzen von Werten das Konzept der zweiten Variation\footnote{Erste entscheidende Anstöße zu diesem Konzept hat Legendre in seiner Arbeit \textit{``Mémoire sur la manière de distinguer les Maxima des Minima dans le Calcul des Variation"} (\cite{Le88}, 1788, ges. 1786) (``Abhandlung über die Art und Weise Maxima und Minima im Variationskalkül zu unterscheiden") gegeben. Die von Legendre mitgeteilte Bedingung ist jedoch wie die Euler--Lagrange'schen Gleichungen noch eine notwendige Bedingung für ein schwaches Extremum. Erst Weierstraß konnte in seinen \textit{``Vorlesungen ueber Variationsrechnung"} (\cite{We79}, 1879) eine umfassende Darstellung des Gegenstandes geben und auch notwendige Bedingungen für das Vorliegen eines starken Extremums formulieren. Die explizite Ausnutzung der zweiten Variation als quadratische Form in einem Hilbert--Raum geschieht erst viel später durch Morse (1892--1974) in seinem Werk \textit{``Calculus of variations in the large"} (\cite{Mo34}, 1934).} zur Überprüfung, ob das Extremum ein Minimum oder Maximum ist, noch zu umschiffen, musste ihm die etwaige globale Natur von Variationsproblemen, wie sie Jacobi durch sein Konzept der konjugierten Punkte in seinem Brief \textit{``Zur Theorie der Variations--Rechnung und der Differential--Gleichungen"} (\cite{Ja36}, 1836) völlig klar herausarbeitet, verborgen bleiben.\\

Ähnlich trägt es sich bei der Euler--Maclaurin'schen Summenformel zu. Geleitet vom Willen, die allgemeine Summe

\begin{equation*}
    f(x):= \sum_{k=1}^{x} g(k)
\end{equation*}
mit einer beliebigen Funktion $g$, auszuwerten zusammen mit der Erkenntnis, dass selbige der Gleichung

\begin{equation*}
    f(x)-f(x-1)=g(x)
\end{equation*}
genügt,  lässt einen wohl nur die Überzeugung der Äquivalenz der formal über eine Differentialgleichung unendlicher Ordnung gefundenen Lösung mit der obigen Summe erst zu dem Ausdruck gelangen, welcher heute -- freilich in anderer Form präsentiert -- als Euler--Maclaurin'sche Summenformel bekannt ist. In der Form, in welcher Euler die Formel niedergeschrieben hat (siehe Abschnitt \ref{subsubsec: Ein Fehlschluss von Euler}), gelangt man nämlich zu einer asymptotischen Entwicklung. Die Existenz zweier Lösungen zu derselben Gleichung impliziert für Euler bereits die Gleichheit, den Aspekt der lediglich asymptotischen Gleichheit schob Euler meist beiseite, obschon er von der generellen Divergenz der durch die Summenformel ausgedrückten Reihe wusste. Eine Ansicht, welche sich auch bei der Euler'schen Untersuchung zur Fakultät bzw. der $\Gamma$--Funktion widerspiegelt.  In seiner Arbeit \textit{``De curva hypergeometrica hac aequatione expressa $y = 1\cdot 2 \cdot 3 \cdots x$"} (\cite{E368}, 1769, ges. 1765) (E368: ``Über die mit der Gleichung $y=1\cdot  2 \cdot 3 \cdots x$ ausgedrückte hypergeometrische Kurve") stellt er alle verschiedenen  Ausdrücke, welche  im Verlauf seiner Schaffenszeit als Interpolation der Fakultät vorgestellt hatte, zusammen. Sie finden sich in §§ 11--12 des besagten Papiers. Die Herleitungen und Diskussionen der entsprechenden Ausdrücke finden sich in der Arbeit \cite{Ay21a} des Verfassers der hiesigen Ausarbeitung umfassend diskutiert. Euler behauptet, wenn auch indirekt, etwa, dass die Stirling'sche Formel (\ref{eq: Stirling}) (siehe Abschnitt \ref{subsubsec: Ein Fehlschluss von Euler}), die bekannte Integraldarstellung (\ref{eq: Solution Gamma}) (siehe Abschnitt \ref{subsubsec: Die Mellin--Transformierte bei Euler}) und anderem die Produktdarstellung (\ref{eq: Euler Produkt Gamma}) (siehe Abschnitt \ref{subsubsec: Produktdarstellung für die Gamma-Funktion}) identisch gleich sind, obwohl erstgenannte eine asymptotische Entwicklung ist, zweites Integral  Realteile größer als $0$ verlangt, das Produkt jedoch für alle nicht--ganzen Zahlen sinnvoll ist\footnote{Dass die Produktformel und das Integral vom funktionentheoretischen Konzept der analytischen Fortsetzung aus betrachtet gleich sind, wusste Euler noch nicht.}.  \\

\begin{figure}
    \centering
    \includegraphics[scale=0.7]{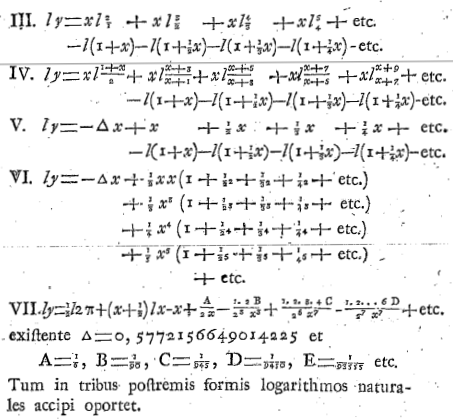}
    \caption{Euler gibt in seiner Arbeit \cite{E368} einige Ausdrücke für den Logarithmus der Fakultät an.}
    \label{fig:E368}
\end{figure}

Tatsächlich stellen alle von Euler mitgeteilte Formeln die auf die eine oder andere $\Gamma$--Funktion dar, obschon Euler auch hier die Gleichheit der Ausdrücke nie explizit gezeigt hat; er hat hingegen auf jeweils unterschiedliche Art nachgewiesen, dass seine vorgestellten Ausdrücke der Funktionalgleichung $f(x+1)=xf(x)$ zusammen mit Bedingung $f(1)=1$ Genüge leisten. Es ist anzunehmen, dass Euler bewusst gewesen ist, dass diese beiden Eigenschaften die Fakultät nicht eindeutig bestimmen, zumal ein Interpolationsproblem nie eindeutig lösbar ist, sofern nicht gewisse Einschränkungen hinzutreten\footnote{Euler erwähnt die Nicht--Eindeutigkeit des Interpolationproblems in mehreren Arbeiten explizit, unter anderem in der schon umfassend diskutierten Arbeit \cite{E189}. Aber auch in der Abhandlung  \cite{E555}, welche eigens dem Thema Interpolation gewidmet ist, wo dies gleich in der Einleitung mitgeteilt wird.}. \\

Es war Euler nicht möglich, bei seinen oben erwähnten Ausdrücken eine Diskrepanz auszumachen, was ihn von der Gleichheit all jener überzeugt haben wird. Heute wird die $\Gamma$--Funktion eindeutig charakterisiert, indem zu den oben genannten beiden Bedingungen noch mindestens eine weitere hinzutritt. Vielerorts ist dies eine Eigenschaft, welche die komplexe Analysis für ihr Verständnis verlangt, und somit auch für Euler nicht zugänglich war, wie eben   auseinandergesetzt worden ist (Abschnitt \ref{subsubsec: Komplexe Analysis}). Am ehesten wäre Euler wohl der zusätzlichen Eigenschaft der logarithmischen Konvexität zugetan gewesen, wie sie im Satz von Bohr--Mollerup\footnote{Die Namensgeber dieses Lehrsatzes haben ihn in ihrem Buch \textit{``Lærebog i Kompleks Analyse vol. III"} (\cite{Bo22}, 1922) (``Lehrbuch der komplexen Analysis, Band 3") bewiesen.} zur eindeutigen Charakterisierung von $\Gamma(x)$ gefordert wird. \\

\paragraph{Partialbruchzerlegungen transzendenter Funktionen}

Wohingegen bei den obigen Beispielen das Ignorieren von Existenzfragen der Richtigkeit der Ergebnisse noch keinen Abbruch tut, verhält sich dies bei nachfolgenden Sachverhalten anders. Wie schon oben (Abschnitt \ref{subsubsec: Ein Fehlschluss von Euler}) angedeutet, soll das Beispiel der Partialbruchzerlegung diskutiert werden. Seine Ideen zur Partialbruchzerlegungen von transzendenten Funktionen, welche Euler in den Arbeiten \textit{``Nova methodus fractiones quascumque rationales in fractiones simplices resolvendi"} (\cite{E540}, 1783, ges. 1775) (E540: ``Eine neue Methode, beliebige rationalen Funktionen in einfache Brüche aufzulösen") und  \textit{``De resolutione fractionum transcendentium in infinitas fractiones simplices"} \cite{E592} (E592: ``Über die Auflösung von transzendenten Brüchen in unendlich viele einfache Brüche") darlegt, nehmen die Existenz einer solchen Zerlegung in völliger Analogie zum endlichen Fall an. Wie im Fall von gebrochen rationalen Funktionen\footnote{Die allgemeine Integration von gebrochen rationalen Funktionen mithilfe von Partialbruchzerlegungen hat Euler in seinen Arbeiten \textit{``Methodus integrandi formulas differentiales rationales unicam variabilem involventes"} (\cite{E162}, 1751, ges. 1748) (E162: ``Eine Methode rationale Differentialformen einer Variable zu integrieren") und \textit{``Methodus facilior atque expeditior integrandi formulas differentiales rationales"} (\cite{E163}, 1751, ges. 1748) (E163: ``Eine leichtere und zügigere Methode rationale Differentialformen zu integrieren") behandelt, im letzten Kapitel seiner  \textit{Calculi Differentials} \cite{E212} hat er eine weitere Methode für dieselbe Aufgabe angegeben.}  versucht Euler, die gesuchte Entwicklung aus den Nullstellen des Nenners allein heraus zu konstruieren. Für den endlichen Fall führt dies zu keinerlei Problemen, im Gegensatz zum Fall unendlich vieler Nullstellen. \\

Eine kurze Illustration seiner Methode am Beispiel der Funktion $\cot x= \frac{\cos x}{\sin x}$ wird dies bereits aufzeigen. Euler würde hier argumentieren, dass die Nullstellen des Nenners gegeben sind als $x=0,-1\pi,1\pi,-2\pi,2\pi,\cdots$, was er wie oben gesehen (\ref{subsubsec: Ein strengerer Beweis}) auch nachgewiesen hatte. Sei also allgemein $x=k\pi$ eine Nullstelle des Nenners. Euler setzt dann, weil die Nullstelle eine einfache ist, an:

\begin{equation*}
    \cot x = \dfrac{A_k}{x- k\pi}+R(x),
\end{equation*}
wobei $R(x)$ eine Funktion ist, die den Faktor $x-k \pi$ nicht enthält. Basierend darauf hat man 

\begin{equation*}
    (x-k \pi)\cot x=A_k +(x-k \pi)R(x)
\end{equation*}
oder

\begin{equation*}
    A_k= (x-k \pi)R(x) +\dfrac{(x-k \pi)\cos x}{\sin x}.
\end{equation*}
Um nun $A_k$ zu berechnen, betrachtet Euler den Grenzwert $x \rightarrow k \pi$:

\begin{equation*}
    A_k= \lim_{x\rightarrow k\pi} \left((x-k \pi)R(x) +\dfrac{(x-k \pi)\cos x}{\sin x}\right).
\end{equation*}
Da $R(x)$ per Annahme den Faktor $x-k\pi$ nicht enthält, kann der erste Term der rechten Seite ignoriert werden, sodass

\begin{equation*}
    A_k= \lim_{x\rightarrow k\pi} \dfrac{(x-k \pi)\cos x}{\sin x}=\cos(k \pi) \lim_{x\rightarrow k\pi} \dfrac{x-k \pi}{\sin x}= \dfrac{\cos(k\pi)}{\cos(k \pi)}=1,
\end{equation*}
wobei  im vorletzten letzten Schritt die Formel $\sin(x+k\pi)=\sin (x)\cos(k\pi)$ und der bekannte Grenzwert $\lim_{x\rightarrow 0}\frac{\sin x}{x}=1$ benutzt wurde. Alternativ ließe sich  die l'Hospital'sche Regel anwenden. Davon ausgehend, so argumentiert Euler, hat man demnach die Partialbruchzerlegung

\begin{equation*}
   \cot x=  \sum_{k \in \mathbb{Z}} \dfrac{1}{x-k\pi},
\end{equation*}
was sich bei Euler ausgeschrieben wie folgt findet

\begin{equation*}
    \cot x = \dfrac{1}{x}+\dfrac{1}{x-\pi}+\dfrac{1}{x+\pi}+\dfrac{1}{x-2\pi}+\dfrac{1}{x+2\pi}+\cdots,
\end{equation*}
was natürlich die bekannte Zerlegung des $\cot x$ ist, wie sie schon oben (Abschnitt \ref{subsubsec: Ein strengerer Beweis}) angegeben worden ist, dort ist es Gleichung (\ref{eq: Partialbruchzerlegung cot}). Die anderen in \cite{E592} behandelten Beispiele sind ebenfalls gebrochen rationale Funktionen der trigonometrischen Funktionen, sodass Euler deren Richtigkeit schon auf andere Weise hatte nachweisen können\footnote{Man findet ein Teil der in \cite{E592} gegebenen Beispiele schon in seiner \textit{Introductio} \cite{E101}, aber auch in seinen \textit{Calculi Differentialis} \cite{E212} mit jeweils anderen Mitteln hergeleitet.}. Obgleich alle von Euler in der erwähnten mitgeteilten Zerlegungen richtig sind, würde sein Vorgehen bei der Funktion

\begin{equation*}
    \dfrac{1}{e^z-1}
\end{equation*}
scheitern, da sein Ansatz die Partialbruchzerlegung allein aus den Nullstellen des Nenners konstruieren will. Oben (Abschnitt \ref{subsubsec: Ein Fehlschluss von Euler}) wurde die richtige Partialbruchzerlegung in (\ref{eq: Partial Fraction Decomposition ex-1}) angegeben und der in ihr auftretende Term  $-\frac{1}{2}$ als ursächlich für falsche Stirling'sche Formel ausgemacht.  \\

Das Erläuterte offenbart bereits, welche größtenteils impliziten Fragestellungen Euler in seinen Forschungen geleitet haben: Ausgehend von einer vorgegebenen Funktion unternimmt Euler etwa den Versuch,  ihre Eigenschaften  wie Nullstellen, Extremstellen oder Polstellen zu ermitteln, woraus er in der Folge  etwaige Alternativdarstellungen ableitet. Das Sinus--Produkt (\ref{eq: sine-Product}) liefert wie die just besprochene Partialbruchzerlegung von $\cot x$  ein Beispiel. Angesichts der damit einhergehenden Schwierigkeiten bezüglich der Validität wird es wenig überraschen, dass der moderne Ansatz gleichsam diametral zum Euler'schen ist: So gibt eines Beispiels wegen Weierstraß bei seinem Produktsatz  die Nullstellen einer Funktion mitsamt weiterer allgemeiner Eigenschaften vor und geht erst anschließend  dazu über, die Existenz einer solchen Funktion nachzuweisen. Man vergleiche die Ausführung zum Weierstraß'schen Produktsatz in Abschnitt (\ref{subsubsec: Ein anderes Vorhaben: Das Weierstraß-Produkt}). Ähnliches gilt natürlich für die Vorgabe von Polstellen beim Mittag-Leffler'schen Satz für Partialbruchzerlegungen. Diese modernen Ansätze geben indes keine direkte Möglichkeit, eine speziell vorgelegte Funktion auf ihre Eigenschaften hin zu untersuchen, sie werden ja im Gegenteil in der Formulierung der Theoreme mit vorgeschrieben. Man kann so etwa mithilfe des Weierstraß'schen Satzes zeigen, dass 

\begin{equation*}
    f(x):= x \prod_{k=1}^{\infty} \left(1-\dfrac{x^2}{k^2\pi^2}\right)
\end{equation*}
eine holomorphe Funktion darstellt. Ihre Gleichheit zu $\sin x$ erfordert jedoch einen gesonderten Beweis\footnote{Euler stellt sich der Aufgabe, von der rechten Seite der letzten Gleichung ausgehend nachzuweisen, dass die linke tatsächlich $\sin x$ ist, in seiner Arbeit \textit{``Exercitatio analytica"} (\cite{E664}, 1794, ges. 1776) (E664: ``Eine analytische Übung"). Er sieht dies allerdings eher als eine kleine mathematische Fingerübung und hat diese ``Übung"{} für andere Funktionen nicht wiederholt.}.  Dieses Studium von ganzen Klassen von Funktionen im Gegensatz zur Untersuchung einzelner Funktionen um ihrer selbst Willen ist später von Gauß,  illustriert am Beispiel der hypergeometrischen Funktion (Abschnitt \ref{subsubsec: Die Darstellung betreffend -- Die hypergeometrische Reihe}), und Riemann -- man denke hier etwa an Theoreme wie den Satz von Riemann--Roch -- weit umfassender betrieben worden und bestimmten seitdem den mathematischen Duktus.

\paragraph{Eulers Beweisversuch des Fundamentalsatzes der Algebra}
\label{para: Eulers Beweisversuch des Fundamentalsatzes der Algebra}

Bei kaum einem anderen Sachverhalt ist die Vernachlässigung  der Existenzfrage dermaßen Inhalt der Kritik anderer Mathematiker gewesen wie beim Euler'schen Beweisversuch des Fundamentalsatzes der Algebra. Dies mag insbesondere darin begründet sein, dass der Euler'sche Beweis, trotz Richtigkeit des zu beweisenden Resultats, leicht nachweisbar unvollständig ist. Die Kritiken richten sich auf seine Argumente in der Arbeit  \textit{``Recherches sur les racines imaginaires des equations"} (\cite{E170}, 1751, ges. 1746) (E170: ``Untersuchueng über die imaginären Wurzeln von Gleichungen"), in welcher Euler also demonstrieren möchte, dass die Gleichung

\begin{equation*}
    a_nx^n+a_{n-1}x^{n-1}+\cdots +a_0=0
\end{equation*}
mit reellen\footnote{Euler nimmt die Koeffizienten  als reell an und diese Annahme ist auch für seinen Beweisversuch von fundamentaler Wichtigkeit.}Koeffizienten stets eine Lösung über den komplexen Zahlen zulässt. Für die Grade $1$,$2$,$3$,$4$ gibt Euler sogar konkrete Formeln, was er wie oben besprochen (Abschnitt \ref{subsubsec: Wegen Unbeweisbarkeit: Wurzeln von Polynomen}) für höhere Grade nicht wiederholen kann.  Seine grundlegende Idee ist die folgende: Da er bereits weiß, wie man ein Polynom vom Grad $2$ behandelt, setzt er für ein Polynom vom Grad $4$ an (§ 27 in \cite{E170}):

\begin{equation}
\label{eq: Ansatz Grad 4}
    x^4+Bx^2+Cx+D=(x^2+ux+\alpha)(x^2-ux+\beta).
\end{equation}
Dieser Ansatz ist allgemein, weil sich die Potenz $x^{n-1}$ durch eine generische Transformation in einem Polynom vom Grad $n$ stets wegschaffen lässt\footnote{Speziell für ein Polynom vom vierten Grad findet sich der Ansatz, selbiges in zwei entsprechende vom Grad 2 zu zerlegen, bereits bei Descartes (1596--1650) in seiner \textit{``La Géométrie"} von 1637. In der Übersetzung \textit{``The Geometry of Rene Descartes"} (\cite{Sm54}, 1954) findet man ein erstes Beispiel auf Seite 184. Aus Gründen der Vollständigkeit sei weiter angemerkt, dass in ähnlicher Manier, der Fall vom Grad $n=6$ mit glücklichem Erfolg in der Arbeit \textit{``Johannis Huddenii Epistula prima de Reductione Aequationum"} (\cite{Hu37}, 1637) (``Erster Brief von Johann Hudde: Über die Reduktion von Gleichungen") von Hudde (1628--1704) abgehandelt worden. Dort wird  das Polynom vom Grad $6$ in eines vom vierten und zweiten Grad zerlegt, was zu einem Polynom vom Grad $15= \binom{6}{2}$ für die zu eruierenden Koeffizienten führt. Da dieses Polynom als eines vom ungeraden Grad stets eine Nullstelle besitzt, beweist Hudde somit zugleich die Gültigkeit des Fundamentalsatzes für Polynome vom Grad $6$. Euler betrachtet diesen Fall in § 50 von \cite{E170}.} . Die Koeffizienten $B$, $C$, $D$ sind dabei vorgegeben, die Zahlen $u$, $\alpha$, $\beta$ sind ausfindig zu machen. Multipliziert man die rechte Seite aus und führt einen Koeffizientenvergleich durch, gelangt man schließlich zur folgenden Gleichung für $u$:

\begin{equation*}
    u^6-2Bu^4-(B^2-4D)u^2-C^2=0.
\end{equation*}
Dies ist ein Polynom vom Grad $6=\binom{4}{2}$ in $u$ und hat wegen des negativen absoluten Terms $-C^2$ stets eine Lösung. Demnach kann man den Ansatz (\ref{eq: Ansatz Grad 4}) stets verwirklichen. Da aber zusätzlich der Fall vom Grad $2$ abgehandelt ist, ist es auch der Fall $n=4$. \\

Damit bewaffnet geht Euler den Fall vom Grad $8$ an und setzt an (§ 34 in \cite{E170}):

\begin{equation}
    \label{eq: Ansatz Grad 8}
    x^8+Bx^6+Cx^5+Dx^4+Ex^3+Fx^2+Gx+H
\end{equation}
\begin{equation*}
    =(x^4+ u x^3+\alpha x^2+ \beta x+\gamma)(x^4-ux^3+\delta x^2+\varepsilon x+ \eta )
\end{equation*}
Die Idee ist dieselbe wie zuvor: Ein Polynom in $u$ zu finden, welches sicher eine reell Nullstelle besitzt, sodass der Ansatz tatsächlich verwirklicht werden kann. Wie Euler auf kombinatorische Begründungen zurückgreifend erklärt, hat dieses Polynom einen Grad von $70=\binom{8}{4}$, was die explizite Berechnung unmöglich werden lässt. Es genügt aber, wie im vorherigen Fall, nachzuweisen, dass der absolute Term negativ ist. Dies begründet Euler entsprechend -- unter Annahme einer Zerlegung des Polynoms in lineare Faktoren -- und kann somit nachweisen, dass (\ref{eq: Ansatz Grad 8}) tatsächlich realisierbar ist. Da aber der Fall vom Grad $4$ bereits abgehandelt ist, ist es somit auch der vom Grad $8$. \\

Nun ist das allgemeine Argument (§ 45 in \cite{E170}),  den Fall vom Grad $2^n$ auf den Fall vom Grad $2^{n-1}$ zu reduzieren, was Euler mit (\ref{eq: Ansatz Grad 4}) und (\ref{eq: Ansatz Grad 8}) entsprechenden Ansätzen erreicht.  Damit ist der Euler'sche Beweis des Fundamentalsatzes der Algebra beendet. Und er selbst muss mit selbigem überaus zufrieden gewesen sein. \\

\begin{figure}
    \centering
    \includegraphics[width=0.8\linewidth]{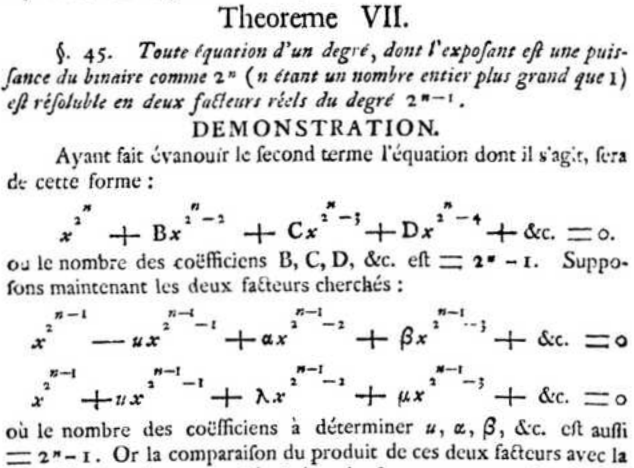}
    \caption{Eulers Ansatz aus \cite{E170} zur Zerlegung eines Polynoms vom Grad $2^n$ in zwei vom Grad $2^{n-1}$.}
    \label{fig:E170Ansatz}
\end{figure}

Allerdings gab es wie erwähnt, wenn auch nicht mehr zu seinen Lebzeiten  kritische Gegenstimmen. Allen voran die von Gauß, welcher den ersten allgemein akzeptierten Beweis in der Abhandlung \textit{``Demonstratio nova theorematis omnem functionem algebraicam rationalem integram unius variabilis in factores reales primi vel secundi gradus resolvi posse"} (\cite{Ga99}, 1799) (``Ein neuer Beweis des Lehrsatzes, dass je ganz rationale Funktion von einer Variable in reelle Faktoren ersten oder zweiten Grades aufgelöst werden kann") gibt. Der von Gauß gewählte Titel zeigt zwar einerseits indirekt an, dass er zuvor gegebene Beweise anerkennt, andererseits schreibt er in § 8, nach vorausgehender Reflexion  des just vorgestellten Euler'schen Arguments:\\

\textit{``Die Regel, nach welcher Euler schließt, aus $2m-1$ Gleichungen die $2m-2$ Unbekannten $\alpha, \beta ~~ \text{etc.}, \mu, \nu ~~ \text{etc.}$ alle rational bestimmen zu können, ist keineswegs allgemein gültig, sondern unterliegt oftmals einer Ausnahme."}\\

Damit ist bereits ein kritischer Punkt in Eulers Beweis getroffen. Letzterer nimmt nämlich an, dass sich mit hinreichender Rechenkraft wie bei (\ref{eq: Ansatz Grad 4}) vorgehen ließe und man das entsprechende Polynom in $u$ gewiss finden könnte. Ein Irrtum, wie Gauß aufzeigt. Unabhängig davon trifft aber folgender Einwand von Gauß noch härter:\\

\textit{``Euler nimmt stillschweigend an, dass die Gleichung $X=0$ $2m$ Wurzeln hat, und setzt deren Summe $=0$, weil der zweite Term in $X$ fehlt. Wie ich über diese Annahme (welche alle Autoren über diesen Gegenstand gebrauchen) richte, habe ich schon in § 3 erläutert."}\\

An nämlicher Stelle (§ 3) schreibt Gauß, dass man zuallererst die \textit{Existenz} der Wurzeln nachzuweisen hat. Hier kristallisiert sich bereits ein zentraler Unterschied zwischen Euler  und Gauß heraus: Nämlich dass Euler, im  Gegensatz zu seinem Nachfolger, die Existenz der Nullstellen überhaupt nicht in Frage stellt. Das \textit{Finden} der Wurzeln bestimmt Eulers Handeln, wohingegen Gauß noch einen Schritt weiter zurücktritt und zunächst die Existenz überhaupt erst einmal gesichert haben will. Euler begeht also hier einen Modalitätsfehlschluss und verwechselt die notwendige mit der hinreichenden Bedingung. Dieses Vorgehen, aus einer notwendigen Bedingung entsprechende Gleichungen abzuleiten, welche  dann still schweigend auch als hinreichend annimmt, ist schon bei der Besprechung seiner Behandlung von Differentialgleichungen mit konstanten Koeffizienten zutage getreten (Siehe Abschnitt \ref{subsubsec: Der inhomogene Fall}).\\

Über den vorgestellten Euler'schen Beweis ist vielerorts in der Literatur diskutiert worden, mit unterschiedlicher Meinung bezüglich der Korrektheit. So nimmt etwa Dunham in seinem Artikel \textit{``Euler and the Fundamental Theorem of Algebra"} (\cite{Du91}, 1991) eine ähnliche Position ein wie Gauß in  \cite{Ga99}. Dahingegen vertritt Speiser im Vorwort von Band 29 der Serie 1 der \textit{Opera Omnia}, genauer Seite X von (\cite{OO29}, 1956), eine fast gegenteilige Position und untermauert dort, wie weit diese Arbeitsweise Euler auch bei anderen Gelegenheiten, wie etwa der Variationsrechnung, geführt hat. Überdies lässt Speiser  F. Frobenius (1849--1917)  im Vorwort von Band 6 der Serie 1 der \textit{Opera Omnia}  (\cite{OO6}, 1921) mit folgendem Zitat zu diesem Gegenstand zu Wort kommen: \\

\textit{``Für die Existenz der Wurzeln einer Gleichung führt er [Euler] jenen am meisten algebraischen Beweis, der darauf fußt, dass jede reelle Gleichung unpaaren Grades eine reelle Wurzel besitzt. Ich halte es für unrecht, diesen Beweis ausschließlich Gauss zuzuschreiben, der doch nur die letzte Feile daran gelegt hat."}\\

Dem Frobenius--Zitat ist im Hinblick auf die moderneren Zugänge über algebraische Konstruktionen der Lösungen mithilfe von Körpererweiterungen umso mehr beizupflichten. In diesem Ansatz werden die Lösungen in Euler'scher Manier ebenfalls vorausgesetzt und deren Existenz durch Konstruktion der Erweiterungskörper auf entsprechende Weise gesichert\footnote{Die Grundidee reicht bis auf den Vollständigkeitssatz von Ostrowski (1893--1986) (über vollständige archimedisch bewertete Körper) aus seiner Arbeit \textit{``Über einige Lösungen der Funktionalgleichung $\varphi(x)\cdot \varphi(y)= \varphi(xy)$"} (\cite{Os18}, 1918) zurück, welchen man mit den Ideen von Artin (1898--1962) aus seiner Abhandlung \textit{``Über die Bewertungen algebraischer Zahlkörper"} (\cite{Ar31}, 1931) oder seinem Buch \textit{``Algebraic Numbers and Algebraic Functions"} (\cite{Ar68}, 1968) insbesondere Theorem 3 aus Kapitel 1 zu einem ``algebraischen"{} Beweis des Fundamentalsatzes der Algebra überführen kann, welcher in seinen Grundzügen Euler'sch ist.}.

\subsubsection{Praxis über Abstraktion}
\label{subsubsec: Praxis über Abstraktion}

\epigraph{Die Wissenschaft soll die Freundin der Praxis sein, aber nicht ihre Sklavin.}{Carl Friedrich Gauß}



Zuletzt soll noch zur Sprache kommen, wie der von Euler angestrebte Praxisbezug sein Vorankommen bei gewissen Begebenheiten obstruiert hat. Der jeweilige Vergleich, mit den Ergebnissen seiner Nachfolger wird dem Gesagten dabei Nachdruck verleihen. Dabei bilden die Gebiete der quadratischen Formen in der Zahlentheorie und die Differentialgeometrie von Flächen den Diskussionsrahmen.

\paragraph{Euler und quadratische Formen in der Zahlentheorie}
\label{para: Quadratische Formen}

Der starke Praxisbezug Eulers kommt eindrücklich zum Vorschein, wenn man seine Forschungen zu quadratischen Formen in der Zahlentheorie betrachtet, welche ihn über seine Laufbahn hinweg immer wieder eingenommen haben. Es wird zwecks der Illustration hinreichen, das Augenmerk auf den Zwei-- und Vier--Quadrate--Satz sowie die Zahlen der Form $nx^2+y^2$  gerichtet zu haben. Eine Vorbereitung auf das nachstehend Präsentierte bietet die folgende Charakterisierung Eulers von Fueter (1880--1950) im vierten Band der ersten Serie der \textit{Opera Omnia} (\cite{OO4}, 1941). Er schreibt auf Seite XVII:\\

\textit{``[Es] ist ein typischer Zug seiner [Eulers] Forschungsmethode, sich ein riesiges Zahlenmaterial zu verschaffen, aus dem er vermöge seiner genialen Divinationsgabe die weitesgehenden Schlüsse zu ziehen vermochte. Heute verfährt man zuweilen umgekehrt, indem man große Theorien entwickelt und kaum oder erst später nach den Möglichkeiten der Realisierung fragt."}\\

Nachdem diese Beschreibung vorausgeschickt worden ist,  sei also mit dem Zwei--Quadrate Satz begonnen, sprich der Aussage, dass eine Primzahl $p$ der Form $x^2+y^2$ eine der Form $p=4n+1$ ist, aber eben auch die Umkehrung dieser Tatsache richtig ist\footnote{Die Aussage, dass $\left(\frac{-1}{p}\right)=(-1)^{\frac{p-1}{2}}$ mit dem Legendre--Symbol $\left(\frac{p}{q}\right)$ ist, ist der von Euler bewiesenen gleichwertig, womit Euler demnach den ersten Ergänzungssatz zum quadratischen Reziprozitätsgesetz bewiesen hat.}. Dieser Satz und sein Nachweis hat  Euler einiges abverlangt. Den ersten großen Schritt zu seinem Beweis ist er dabei in der Arbeit \textit{``De numeris, qui sunt aggregata duorum quadratorum"} (\cite{E228}, 1758, ges. 1749) (E228: ``Über die Zahlen, die  Aggregate zweier Quadrate sind") gegangen. Hier demonstriert er folgende auf das Ziel vorbereitende Sätze. Zunächst in § 5: Für zwei Primzahlen $p=a^2+b^2$ und $q=c^2+d^2$ gilt:

\begin{equation*}
    p \cdot q = (a^2+b^2)(c^2+d^2)= (ac+bd)^2+(ad-bc)^2,
\end{equation*}
von welchem Satz er in § 8 (Proposition 1) eine Art Umkehrung vorstellt. Gilt $pq=e^2+f^2$ und ist $p=a^2+b^2$ eine Primzahl, so ist $q$ zwingend die Summe zweier Quadrate. In § 14 (Proposition 2) zeigt er den verwandten Satz, dass falls $pq=e^2+f^2$, aber nicht $p=a^2+b^2$ gilt, die Zahl $q$ ebenfalls nicht die Summe zweier Quadrate ist. In § 19 (Proposition 4) präsentiert er dann einen zentralen Satz:

\begin{Thm}[Teiler der Summe zweier Quadrate]
\label{Theorem: Teiler Quadrate}
    Gilt für eine Primzahl $p$, dass $p|(a^2+b^2)$, dann existiert eine Summe von zwei Quadraten $c^2+d^2$, sodass $p|(c^2+d^2)$.
\end{Thm}
Dieser Satz wird dann ergänzt von Proposition 4 aus § 22, dass die Summe $a^2+b^2$ überhaupt nur von Primzahlen $p$ geteilt werden kann, die selbst die Summe zweiter Quadrate sind. In § 28[a]\footnote{In der Originalarbeit wird irrtümlich § 28 wiederholt, weshalb hier die Notation der \textit{Opera Omnia} Version übernommen wird.} formuliert Euler schließlich den Satz, den er zeigen möchte:

\begin{Thm}[Zwei--Quadrate--Satz]
\label{Theorem: Zwei--Quadrate--Satz}
    Ist $p$ eine Primzahl der Form $4n+1$, so ist sie auch die Summe zweier Quadrate $p=a^2+b^2$. Umgekehrt ist jede Primzahl der Form $p=a^2+b^2$ gleichzeitig eine der Form $4n+1$.
\end{Thm}

Diesen Satz vermochte Euler in \cite{E228} noch nicht zu beweisen, und weist seine Anstrengungen diesbezüglich entsprechend als ``Tentamen demonstrationis"{} aus. In seiner Abhandlung \textit{``Demonstratio theorematis Fermatiani omnem numerum primum formae $4n+1$ esse summam duorum quadratum"} (\cite{E241}, 1760, ges. 1750) (E241: ``Beweis des Fermat'schen Satzes, dass jede Primzahl der Form $4n+1$ die Summe zweier Quadrate ist") reicht Euler den ersehnten Nachweis nach. Die Schwierigkeit besteht für Euler darin, ein Argument über die Methode des unendlichen Abstiegs zu finden, was ihm dann im letztgenannten Werk schließlich gelingt. Seine Idee sei hier in Kürze wiedergegeben. Für $p=4n+1$ Euler geht von der allgemeinen Form $a^{4n}-b^{4n}$ aus. Da $p$ eine Primzahl ist, weiß er aus dem kleinen Satz von Fermat\footnote{Wie oben in Abschnitt  (\ref{subsubsec: Beweistechnik der Induktion bei Euler}) beschrieben, hat Euler diesen Satz bereits in \cite{E54} bewiesen.}, dass $p|(a^{4n}-b^{4n})$ gilt. Wegen $a^{4n}-b^{4n}=(a^{2n}-b^{2n})(a^{2n}+b^{2n})$ gilt aber auch $p|(a^{2n}-b^{2n})$ oder $p|(a^{2n}+b^{2n})$. Euler möchte nun zeigen, dass $p$ den zweiten Faktor teilt, wozu er nachweist, dass der erste $p$ nicht als Teiler enthält. Euler argumentiert nun per Widerspruch und nimmt gegenteilig an, dass doch $p|(a^{2n}-b^{2n})$ gilt. Dementsprechend müssen insbesondere die Differenzen:

\begin{equation*}
    2^{2n}-1^{2n}, ~~  3^{2n}-2^{2n}, ~~  4^{2n}-3^{2n}, ~~ \cdots, ~~  (4n)^{2n}-(4n-1)^{2n}
\end{equation*}
alle durch $p$ teilbar sein. Genauso  müssen es die Differenzen dieser 1. Differenzen sein, weiter deren Differenzen, und erneut die Differenzen von letzteren usw. Euler nutzt nun aus, dass alle Differenzen der Ordnung $2n$ einander gleich sind und kann sie auch zu $(2n)!$ berechnen. Jedoch gilt offensichtlich nicht $p|(2n)!$, sodass in der Konsequenz auch nicht $p|(a^{2n}-b^{2n})$ gelten kann, und damit das Gewünschte bewiesen ist: $p=4n+1$ teilt stets die Summe zweiter Quadrate, sodass verbunden mit der Eigenschaft aus Theorem (\ref{Theorem: Teiler Quadrate})  der Zwei--Quadrate--Satz (\ref{Theorem: Zwei--Quadrate--Satz})  folgt.\\

Betrachtet man den Euler'schen Beweis\footnote{In gleicher Manier, also über Betrachtung der Differenzen, hat Euler überdies beweisen, dass für eine Primzahl $p$ gilt: $p=8n+1 \Longrightarrow p= x^2+2y^2$ und $p=3$ oder $p=6n+1 \Longleftrightarrow p=x^2+3y^2$, wobei die Darstellungen als Summen von Quadraten jeweils eindeutig sind. Den Beweis der ersten Aussage findet man in \cite{E256}, die zweite wird in der Arbeit \textit{``Supplementum quorundam theorematum arithmeticorum, quae in nonnullis demonstrationibus supponuntur"} (\cite{E272}, 1763, ges. 1759) (E272: ``Supplement gewisser zahlentheoretischer Theoreme, welche bei einigen Beweisen vorausgesetzt werden") demonstriert. Den Nachweis, dass $x^2+2y^2$ auch die Zahlen der Form $8n+3$ gleichermaßen wie die der Form $8n+1$ erfasst, konnte Euler, wie er auch in \cite{E256} einräumt, zunächst nicht beweisen. Es sollte ihm in der Arbeit \textit{``Demonstrationes circa residua ex divisione potestatum per numeros primos resultantia"} (\cite{E449}, 1774, ges. 1772) (E449: ``Beweis über die Reste, welche aus der Teilung von Potenzen durch Primzahlen entstehen") gelingen.  Dieses Verdienst wird oft Lagrange zugeschrieben (siehe etwa das Buch \textit{``Elliptic Functions} (\cite{Ch85}, 1985), wo dies auf S. 153 erwähnt wird), welcher dies aus seinen Untersuchungen in seinen Arbeiten \textit{``Recherches d'arithmétique -- Première Partie"} (\cite{La75}, 1775) (``Untersuchungen zur Zahlentheorie -- Erster Teil") sowie \textit{``Recherches d'arithmétique -- Seconde Partie"} (\cite{La77}, 1777) (``Zahlentheoretische Untersuchungen -- Zweiter Teil") zu binären quadratischen Formen bewiesen hat. Lagrange und Euler haben damit gleichermaßen den zweiten Ergänzungssatz zum quadratischen Reziprozitätsgesetz entdeckt, dass $\left(\frac{2}{p}\right)=(-1)^{\frac{p^2-1}{8}}$ mit dem Legendre--Symbol $\left(\frac{p}{q}\right)$ gilt.}, mag man sich vielleicht vor der Euler'schen Kreativität verneigen, weiß aber indes auch den Wert der modernen Konzepte wie das Gauß'schen Zahlen, des Euklidischen Ringes und der Norm umso mehr zu schätzen\footnote{Man bemerke etwa, dass Euler Formeln wie $(a^2+b^2)(c^2+d^2)=(ac-bd)^2+(ad-bc)^2$ ohne die Eigenschaften von Normfunktionen finden musste. Er wird in diesem Fall vermutlich über die komplexen Zahlen vorgegangen sein.}. Mit seinen Methoden und seiner Auffassung von Mathematik ist es für Euler meist schwerer gewesen, Beweise in der Zahlentheorie zu führen, da eine Analytisierung Euler bei dieser und ähnlichen Fragestellung meistens, wenn überhaupt, lediglich partiell gelungen ist. Als Beispiel hierfür sei die Arbeit \cite{E256} herangezogen, welche oben (Abschnitt \ref{subsubsec: Beweistechnik der Induktion bei Euler}) im Kontext der Induktion schon einmal genannt wurde. Dies Abhandlung hat einen eher didaktischen bzw. erläuternden Charakter und nicht zuletzt aus diesem Grunde sind die in dieser Arbeit über alte Sätze hinausgehenden Theoreme in ihrer Anzahl klein. Genauer beschränken sie sich auf den Satz, dass eine Zahl der Form $2a^2+b^2$ nur Primteiler derselben Form zulässt. Dies ist Beobachtung 6 in Eulers Arbeit, bewiesen wird sie dann in Theorem 9 (§ 42).  Der Beweis wird  mithilfe eines unendlichen Abstiegs geführt. In § 46 versucht Euler, den Beweis  auf die Form $ma^2+b^2$ zu übertragen und gelangt zur Bedingung $\frac{m+1}{4}\leq 1$, was aber nur den Fall $a^2+b^2$ abdeckt, womit die angestrebte Verallgemeinerung misslingt. Diese Untersuchung und ohnehin die ganze Arbeit steht als pars pro toto für Eulers Forschungen zu den quadratischen Formen: Durch eine ihresgleichen suchende Beobachtungsgabe findet Euler sehr  tiefe Wahrheiten in der Arithmetik, indes gelingen allgemeine Beweisversuche nicht, wodurch er  gezwungen ist, jedes Problem  mit einem neuen Geniestreich zu überwinden. Die nötigen Verallgemeinerungen konnten dann erst nach den entsprechenden Begriffsbildungen erfolgen\footnote{Aus dem Gesagten erhellt sich auch, warum Euler seinen Traktat zur Zahlentheorie \textit{``Tractatus de numerorum doctrina capita sedecim, quae supersunt"} (\cite{E792}, begonnen ca. 1745)  (E792: ``Traktat über die Zahlentheorie in sechzehn Kapiteln") nicht fertig gestellt hat. In selbigem findet man gar Ergänzungssätze des kubischen und biquadratischen Reziprozitätsgesetzes, welche Euler wie viele andere Beobachtungen aber zu seinen Lebzeiten nicht zu beweisen vermochte. Man konsultiere auch das Vorwort zu Band 5 der ersten Serie der \textit{Opera Omnia} (\cite{OO5}, 1944).}. \\

Einen Teilbeweis des Vier--Quadrate--Satzes, also der Aussage, dass jede natürlich Zahl die Summe vierer Quadrate ist, legt Euler  in der Abhandlung \textit{``Demonstratio theorematis Fermatiani omnem numerum sive integrum sive fractum esse summam quatuor pauciorumve quadratorum"} (\cite{E242}, 1760, ges. 1751) (E242: ``Beweis des Fermat'schen Satzes, dass jede ganze oder rationale Zahl die Summe von vier oder weniger Quadraten ist") vor. Er beweist die Richtigkeit im Falle rationaler Zahlen, muss allerdings  Lagrange den Vortritt beim Fall natürlicher Zahlen lassen\footnote{Euler selbst reicht basierend auf der Lagrange'schen Arbeit einen Beweis für die ganzen Zahlen in der Abhandlung \textit{``Novae demonstrationes circa resolutionem numerorum in quadrata"} (\cite{E445}, 1773, ges. 1772) (E445: ``Neue Beweise zur Auflösbarkeit von Zahlen in Quadrate") nach.}. Zentral für den Lagrange'schen Beweis in seiner Arbeit \textit{``Démonstration d'un théorème d'arithmétique"} (\cite{La70}, 1770) (``Beweis eines Satzes in der Zahlentheorie") war die Identität, die  Euler in § 93 von \cite{E242} bewiesen hat. Modern würde man diese wohl mithilfe der Quaternionen\footnote{Die Quaternionen werden oft Hamilton zugeschrieben, der sie in 1843 entdeckt hat. Jedoch finden sich alle Eigenschaften bereits bei Rodrigues (1795--1851). Man konsultiere dazu die Arbeit \textit{``Hamilton, Rodrigues and the quaternion scandal"} (\cite{Al89}, 1989).} ableiten, Euler muss besagte Formel:

\begin{equation*}
    (a^2+b^2+c^2+d^2)(p^2+q^2+r^2+s^2)=x^2+y^2+z^2+v^2
\end{equation*}
mit

\begin{equation*}
    x=ap+bq+cr+ds,~~y=aq-bq\pm cs\mp dr,~~ z=ar \mp bs-cp\pm dq \quad \text{sowie}
\end{equation*}
\begin{equation*}
    ~~v=as\pm br \mp cq -dp
\end{equation*}
dahingegen auf anderem, leider nicht zweifelsfrei rekonstruierbarem, Wege gefunden haben\footnote{Eine mögliche Erklärung findet sich in der Arbeit \textit{``Leonhard Eulers handschriftlicher Nachlass zur Zahlentheorie"} aus dem Buch (\textit{``Leonhard Euler 1707–1783 – Beiträge zu Leben und Werk"} (\cite{Bu83}, 1983) angedeutet.}.\\

\begin{figure}
    \centering
    \includegraphics[scale=1.1]{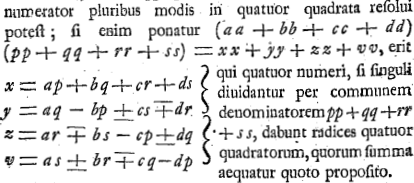}
    \caption{Eulers berühmte Identität für das Produkt von zwei Summanden aus je vier Quadraten aus seiner Arbeit \cite{E242}.}
    \label{fig:E242vierQuadrate}
\end{figure}

 Der Vier--Quadrate--Satz stellt gleichsam den Punkt dar, ab welchem Lagrange Euler in Bezug auf die Zahlentheorie überkommen hat. Noch eindrücklicher sollte dies dann schließlich bei den Numeri idonei werden, wozu auch auf das Buch \textit{``Primes of the Form $x^2+ny^2$"} (\cite{Co13}, 2013), insbesondere auf das erste Kapitel, verwiesen sei. Bei nämlichen benötigte Euler explizit ein Ergebnis von Lagrange über quadratische Formen, um auch selbst noch Fortschritte zu verzeichnen\footnote{Lagrange hatte seinen ``Vorteil"{}  aus seiner Präferenz des algebraischen Zugangs gegenüber demjenigen über die Analysis, welchen Euler bevorzugte.}. Die Numeri idonei, welche Euler insbesondere in den drei Arbeiten  \textit{``De insigni promotione scientiae numerorum"} (\cite{E598}, 1785, ges. 1775) (E598: ``Über einen außergewöhnlichen Fortschritt der Zahlentheorie"), \textit{``Novae demonstrationes circa divisores numerorum formae $xx + nyy$"} (\cite{E610}, 1787, ges, 1775) (E610: ``Neue Beweise über die Teiler der Form $xx+nyy$"), \textit{``De formulis speciei $mxx + nyy$ ad numeros primos explorandos idoneis earumque mirabilibus proprietatibus"} (\cite{E708}, 1801, ges. 1778) (E708: ``Über die Formeln der Gattung $mxx+nyy$, welche zur Ermittlung von Primzahlen geeignet sind, und deren wundersame Eigenschaften") untersucht. Das Prädikat ``ideonei"{} (was sich mit geeignet oder tauglich übersetzen lässt) rührt aus ihrer Verwendbarkeit zum Testen auf Primalität her\footnote{Euler hat der Erstellung von Primzahltabellen gar einzelne Arbeiten gewidmet. So ist etwa seine Arbeit \textit{``Utrum hic numerus $1000009$ sit primus necne inquiritur"} (\cite{E699}, 1797, ges. 1778) (E699: ``Ob die Zahl $1000009$ eine Primzahl ist oder nicht, wird untersucht") allein der Frage verschrieben, ob $1000009$ eine Primzahl ist oder nicht. Dies kann Euler vermöge der Darstellung $1000009=1000^2+3^3=972^2+235^2$ verneinen.}. Letzteres musste Euler fast  zu den Numeri ideonei führen. Zum einen war, wie gesehen, Euler für die Summe $p=x^2+y^2$  das Ergebnis schon lange bekannt, dass man $p$ als Primzahl nachgewiesen hat, wenn es keine weitere Darstellung $a^2+b^2$ gibt, sodass $x^2+y^2=a^2+b^2$, aber die alle natürlichen Zahlen $x$, $y$, $a$, $b$ voneinander verschieden sind. Die Euler'sche Frage, für welche $N$ in $x^2+Ny^2$ dieses Kriterium analog Geltung hat, ist demnach nur konsequent. Explizit formuliert er die Frage dabei in § 6 von\cite{E708}, wo er schreibt: \\

 \begin{figure}
     \centering
     \includegraphics[width=0.8\linewidth]{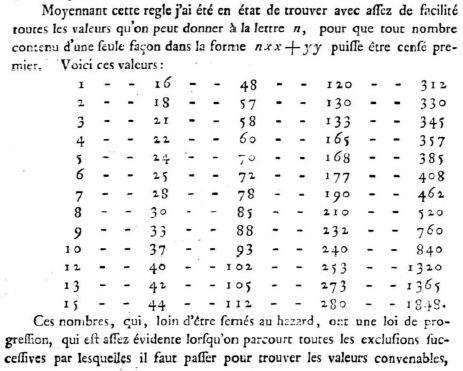}
     \caption{Euler listet seine $65$ ``Numeri idonei"{} in der Arbeit \cite{E498} auf.}
     \label{fig:E498.PNG}
 \end{figure}

 \textit{``Weil also diese Proposition höchst gewiss ist, dass alle Zahlen, die auf mehr als eine einzige Weise in derselben Formel $mxx+nyy$ enthalten sind, nicht prim sondern zusammengesetzt sind, und daher eine Primzahl nur auf eine einzige Weise in einer solchen Form enthalten sein kann, wollen wir die Umkehrung dieser Proposition betrachten, welche man wie folgt formulieren kann: Alle zusammengesetzten Zahlen, die in der Form $mxx+nyy$ enthalten sind, sind auch auf mehr als auf eine Weise  in derselben Form enthalten, oder alle auf nur eine einzige Weise in dieser Form enthaltenen Zahlen sind gewiss prim."}\\

 Leicht sieht man ein, dass Eulers Proposition aus Sicht der modernen Zahlentheorie zu der Frage, welche quadratisch--imaginären Zahlkörper lediglich ambige Klassen besitzen, äquivalent ist.
  Gauß hat bekanntermaßen hingegen  in seinen \textit{Disquistiones} (\cite{Ga01}, 1801) einen vollkommen anderen Zugang zu diesen Zahlen gefunden und verfolgt, nämlich den über das Geschlecht\footnote{Euler spricht in seinen oben erwähnten Arbeiten zwar ebenfalls von \textit{``Genera"}, was sich als Geschlechter übersetzten ließe, bezeichnet dabei aber nicht das, was  Gauß in seinen \textit{Disquistiones} darunter verstanden hat und seither darunter verstanden wird.}. Man weiß heute, dass die Anzahl der Klassen  gleich der Anzahl der Geschlechter im Ring $r(m)$ ist, welcher von $(1, \sqrt{-m})$ erzeugt wird, weswegen man auch -- wie Gauß es getan hat -- untersuchen kann, in welchen dieser Ringe jedes Geschlecht genau eine Klasse enthält\footnote{Gauß war sich der Äquivalenz seiner Fragestellung mit der von Eulers bewusst; er zitiert Eulers Beitrag entsprechend. Gauß formuliert diese Äquivalenz zwar ebenfalls nicht explizit, aber sie ist implizit in seinen vorgestellten Faktorisierungsmethoden zu finden. Man betrachte insbesondere §§ 329--334 seiner \textit{Disquisitiones} \cite{Ga01}.}. Dies ist ein rein algebraischer und überdies nicht anwendungsorientierter Zugang, welcher es Gauß ermöglichte, die Natur der Numeri idonei{} tiefer zu verstehen als Euler. Rückblickend mag man konstatieren, dass die Vermeidung des Komplexen und fehlende Begriffe wie den der Diskriminante einer quadratischen Form Euler seine Numeri idonei nicht in einem für ihn über einen gewissen Punkt hinaus erhellendem Licht erscheinen ließen. Dennoch hat er durch explizites Probieren\footnote{Das Probieren basiert dabei auf einem Verfahren, dass Euler in \cite{E708} genauer vorstellt.}, wie nach ihm Gauß, auf diesem Wege $65$ entsprechende Zahlen ausfindig gemacht, von denen $1848$ die größte ist\footnote{Euler hat, wie er in seiner Arbeit \textit{``Illustratio paradoxi circa progressionem numerorum idoneorum sive congruorum"} (\cite{E725}, 1806, ges. 1778) (E725: ``Beleuchtung des Paradoxons über die Progression der geeigneten oder kongruenten Zahlen") schreibt, seine Rechnungen bis hin zu $10000$ ausgedehnt und keine weitere Zahl gefunden. Die scheinbare Endlichkeit der Menge der Numeri ideonei bezeichnet er gar als paradox. Die Endlichkeit wurde erst 1973 in der Arbeit \textit{``Exponents of the class groups of complex quadratic fields"} (\cite{We73}, 1973) bewiesen. Ob die von Euler und Gauß gefundenen $65$ Zahlen tatsächlich die einzigen sind, ist derzeit eine offene Frage.}. Die vollständige Liste findet man bereits in Eulers Brief \textit{``Extrait d'un lettre de M. Euler à M. Beguelin en Mai 1778"} (\cite{E498}, 1779, ges. 1778) (E498: ``Auszug eines Briefs an Hr. Beguelin im Mai 1778"). Euler schreibt bereits hier: \\

  \textit{``Zu diesem Zwecke habe ich folgende Regel gefunden und bewiesen: Wenn alle in der Form $n+yy$ enthaltenen Zahlen, die zugleich kleiner sind als $4n$, (sofern man freilich für $n$ zu $y$ prime Zahlen einsetzt), entweder  Primzahlen $p$ oder das Doppelte von Primzahlen $2p$ oder die Quadrate von Primzahlen $pp$ oder schließlich irgendeine Potenz von $2$ sind, dann wird der Wert von $n$, der diesen Bedingungen Genüge leistet, als geeignet für die Untersuchung einer solchen beliebigen vorgelegten Zahl, zugelassen werden können. [...] Vermöge dieser Regel war ich in der Lage, mit hinreichender Leichtigkeit alle Werte zu finden, welche man dem Buchstaben $n$ zuteilen kann, sodass jede in auf nur eine einzige Weise in der Form $nxx+yy$ enthaltene Zahl als Primzahl angesehen werden kann."}\\

  Nach diesen Ausführungen zu zahlentheoretischen Beiträgen Eulers, soll der Übergang zur Differentialgeometrie statt haben, in welchem Gebiet Euler wieder -- und dies ist noch höherem Maße -- von seinem Nachfolger Gauß überkommen wurde.

\paragraph{Differentialgeometrie}
\label{para: Differentialgeometrie}

Der zweite in diesem Abschnitt zu besprechende Zweig der Mathematik ist die Differentialgeometrie von Flächen, in welchem Euler höchst bemerkenswerte Ergebnisse zutage gefördert hat. Dennoch hat die Geschichte gezeigt, dass seine Nachfolger -- allen voran Gauß und Riemann -- seiner Ergebnisse nicht unbedingt bedurft haben und zugleich sehr weit über sie hinausgegangen sind, welche Aussage auch in der Arbeit \textit{``Euler’s Contribution to Differential Geometry and its Reception"} (\cite{Re07}, 2007) weiter elaboriert wird. \\

Ein Großteil von Eulers Beiträgen zu diesem Themenfeld haben ihren Ursprung im Problem der Erstellung von Landkarten. Neben den beiden explizit der Kartenerstellung gewidmeten Arbeiten 
\textit{``De proiectione geographica superficiei sphaericae"} (\cite{E491}, 1778, ges. 1775), (E491: ``Über die stereographische Projektion auf die Kugeloberfläche") und der unmittelbar folgenden \textit{``De proiectione geographica Deslisliana in mappa generali imperii russici usitata"} (\cite{E492}, 1778, ges. 1775), (E492: ``Über die steographische Projektion, welche Deslisle für die Karte des russischen Reichs verwendet hat") sei als stellvertretendes Exempel die Arbeit  \textit{``De solidis quorum superficiem in planum explicare licet"} (\cite{E419}, 1772, ges. 1770), (E419: ``Über Festkörper, deren Oberfläche sich in die Ebene ausbreiten lassen") genannt. Letztgenannte handelt aus heutiger Sicht über die Art von Flächen, die als abwickelbar bezeichnet werden. Geleitet hat Euler dabei vermutlich ursprünglich die Frage nach der Abwickelbarkeit der Kugelfläche auf die Ebene, welche natürlich für die Kartenprojektion vonnöten ist. Einen ausführlichen Beweis  der Unmöglichkeit einer längen-- und winkeltreuen Abbildung von der Sphäre auf die Ebene hat Euler  in der Arbeit \textit{``De repraesentatione superficiei sphaericae super plano"} (\cite{E490}, 1778, ges. 1775), (E490: ``Über die Darstellung der Kugeloberfläche in der Ebene") (§§ 1-8) gegeben\footnote{Dieser Beweis wird von Speiser im Vorwort zu Band 28 von Serie 1 der \textit{Opera Omnia} (\cite{OO28}, 1955) ebenfalls nachvollzogen. Die Kugelgeometrie hat Euler auch unabhängig von der Kartenerstellung um ihrer selbst willen beschäftigt. Die allgemeinen Formeln, welche Euler in der Arbeit \textit{``Trigonometria sphaerica universa, ex primis principiis breviter et dilucide derivata"} (\cite{E514}, 1782, ges. 1775), (E514: ``Die ganze Kugeltrigonometrie, aus ersten Grundlagen zügig und klar hergeleitet") vorstellt, zeugen davon. Auch die Arbeit \textit{``Variae speculationes super area triangulorum sphaericorum"} (\cite{E698}, 1797, ges. 1775), (E698: ``Verschiedene Betrachtungen über die Fläche von Kugeldreiecken") kann in diesem Zusammenhang genannt werden.}. Eulers Beweis benutzt nur das Linienelement und reduziert sich auf die Gleichung

\begin{equation*}
    dx^2+dy^2= du^2+\cos^2 udv
\end{equation*}
mit den rechtwinkligen Koordinaten $x$, $y$ sowie den Längen-- und Breitengraden $u$, $v$ auf der Kugel. Für eine kongruente Abbildung von der Kugel auf die Ebene muss die obige Gleichung die Identität werden. Setzt man

\begin{equation*}
    dx=pdu+qdv \quad \text{und} \quad dy=rdu+sdv,
\end{equation*}
erzwingt dies die Gleichungen:

\begin{equation*}
    p^2+r^2=1, \quad pq+rs=0, \quad q^2+s^2=\cos^2 u, 
\end{equation*}
welche Euler mit 

\begin{equation*}
    p= \cos \varphi, \quad q=- \sin \varphi \cos u, \quad r=\sin \varphi, \quad s= \cos \varphi \cos u
\end{equation*}
löst. Damit gelangt Euler zu den Darstellungen:

\begin{equation*}
    dx=\cos \varphi du- \sin \varphi \cos u dv \quad \text{und} \quad dy= \sin \varphi du +\cos \varphi \cos u dv.
\end{equation*}
Die notwendige Integrierbarkeit dieser Differentialformen induziert dieses System von Gleichungen:

\begin{equation*}
    \renewcommand{\arraystretch}{2,0}
\setlength{\arraycolsep}{0.5mm}
\begin{array}{rclcl}
   - \sin \varphi \dfrac{\partial \varphi}{\partial v}&=& - \cos u \cos \varphi \dfrac{\partial \varphi}{\partial u}&+& \sin \varphi \sin u, \\
   \cos \varphi \dfrac{\partial \varphi}{\partial v}&=& - \cos u \sin \varphi \dfrac{\partial \varphi}{\partial u}&-& \cos \varphi \sin u,
\end{array}
\end{equation*}
welche auf die zwei sich widersprechenden Gleichungen

\begin{equation*}
    -\cos u \dfrac{\partial \varphi}{\partial u}=0 \quad \text{und} \quad \dfrac{\partial \varphi}{\partial v}=-\sin u
\end{equation*}
geführt werden. Die zuerst erwähnte Arbeit, eigens über abwickelbare Flächen, \cite{E419} erstreckt sich insofern weiter, da in ihr allgemeine Gleichungen für abwickelbare Flächen abgeleitet werden. Das Problem formuliert Euler dabei in § 1:\\

\textit{``Eine allgemeine Gleichung für alle Festkörper zu finden, deren Oberfläche sich in die Ebene ausbreiten lassen."}\\

Die analytische Lösung\footnote{Euler gibt weiterhin eine geometrische Lösung (§§ 8--37) und eine aus der Projektionstheorie (§§ 38--54) an.} gibt er schließlich in § 7 als folgenden Satz von  Gleichungen:\\

\textit{``Nach Vorlage der zwei Variablen $t$ und $u$ sechs so beschaffene Funktionen $l$, $m$, $n$ und $\lambda$, $\mu$, $\nu$ derer zu finden, dass den folgenden sechs Bedingungen Genüge geleistet wird:}

\begin{equation*}
\renewcommand{\arraystretch}{2,0}
\setlength{\arraycolsep}{0.5mm}
\begin{array}{ll}
     \text{I.} ~~ & \dfrac{\partial l}{\partial u}= \dfrac{\partial \lambda}{\partial t}, \quad \text{II.} ~~ \dfrac{\partial m}{\partial \mu}= \dfrac{\partial \lambda}{\partial t}, \quad \text{III.} ~~ \dfrac{\partial n}{\partial \mu}= \dfrac{\nu \lambda}{\partial t}, \\
     \text{IV.} & l^2+m^2+n^2=1, \quad \text{V.} ~~ \lambda^2+\mu^2+\nu^2=1, \\
     \text{VI.} & l\lambda +m \mu +n \nu =0."
\end{array}
\end{equation*}
Ein Anwendungsbeispiel findet man  im letzten Paragraphen der Abhandlung. Euler behauptet von der Fläche definiert über

\begin{equation*}
    4y^3x+72y^2xxz-yy-18yxz+27xxzz+2z=0,
\end{equation*}
eine abwickelbare Fläche zu sein. Dies ist wegen eines Rechenfehlers nicht der Fall, das Euler'sche Beispiel führt indes zur abwickelbaren Fläche

\begin{equation*}
    -4xy^3-y^2+18xyz+27x^2z^2+4z=0,
\end{equation*}
wie Speiser in der \textit{Opera Omnia} Version der Arbeit auch anmerkt. \\

Obwohl Euler in dieser Arbeit \cite{E419} sogar die Objekte benutzt, welche man heute Gauß'sche Koordinaten in der Flächentheorie nennt,  ist der Euler'sche Beitrag von Gauß in seiner berühmten Arbeit \textit{``Disquisitiones generales circa superficies curvas"} (\cite{Ga28}, 1828, ges. 1827) (``Allgemeine Untersuchungen über gekrümmte Flächen") nicht erwähnt. Unabhängig davon geht Gauß in dieser für die gesamte Differentialgeometrie wegweisenden Arbeit dermaßen über Eulers Erkenntnisse heraus, dass letztere auch von anderen Mathematikern nach Gauß kaum noch aufgegriffen worden sind.  Der wesentliche Unterschied zwischen Gauß und Euler betrifft das Anwendungsgebiet: Sind die Flächen, die Euler untersucht, stets Grenzflächen von Festkörpern\footnote{Gauß selbst resümiert die Situation in § 12 von \cite{Ga28} -- unmittelbar nach Formulierung des Theorema egregium -- dahingehend, dass der Spezialfall der abwickelbaren Flächen der ist, auf den die Mathematiker bis dahin ihre Untersuchungen beschränkt haben.}, sind die Flächen um ihrer selbst willen bei Gauß Gegenstand der Untersuchung\footnote{Natürlich hatten die Gauß'schen Untersuchungen auch in praktischen Fragestellungen des Autors ihren Ursprung, wie Klein in seinem Buch \cite{Kl56} auch entsprechend herausarbeitet, jedoch hat Gauß sich in dieser Arbeit von selbigen völlig gelöst gehabt.}, wodurch er zu einer zentralen Invariante der Flächentheorie gelangt ist: Der Gauß'schen Krümmung. Euler konnte durch seine Gedankengänge hingegen nur  extrinsische Eigenschaften von Flächen behandeln. Den größten Grad an Allgemeinheit erreicht er dabei in seiner Arbeit \textit{``Recherches sur la courbure des surfaces"} (\cite{E333}, 1767, ges. 1763), (E333: ``Untersuchungen über die Krümmung von Flächen"). Das Problem formuliert er in § 1:\\

\textit{``Eine Fläche, deren Natur bekannt ist, sei von einer beliebigen Ebene geschnitten worden, die Krümmung des Schnittes zu ermitteln, welcher daraus entstanden ist."}\\

Die zentrale Formel teilt Euler im letzten Paragraphen mit, sie lautet:

\begin{equation*}
    r= \dfrac{2fg}{f+g-(f-g)\cos (2\varphi)}
\end{equation*}
wobei $r$ der Krümmungsradius des Normalschnittes ist, entsprechend dem Winkel $\varphi$ durchgeführt, $f$ und $g$ sind die Krümmungsradien der Hauptschnitte. Euler beschreibt seine Formel selbst wie folgt:\\

\textit{``So reduziert sich die Beurteilung über die Krümmung von Flächen, wie kompliziert sie auch immer anfangs erschienen mögen, für jedes Element auf die Kenntnis von zwei Krümmungsradien, von welchen der eine der größte und der andere der kleinste in diesem Element ist; diese zwei Dinge bestimmen gänzlich die Natur der Krümmung , indem sie uns die Krümmung von allen möglichen Schnitten enthüllt, welche senkrecht zum vorgelegten Element sind."}\\

Vergleicht man dies mit dem ebenfalls Euler zugeschriebenen Satz:

\begin{equation}
\label{eq: Satz Euler Diffgeo}
    k_n = k_1 \cos^2 \alpha +k_2 \sin^2 \alpha
\end{equation}
mit der Normalkrümmung $k_n$ und den Hauptkrümmungen $k_1$, $k_2$ sowie dem Winkel $\alpha$ zwischen gegebener Tangentialrichtung und zu $k_1$ gehöriger Tangentialrichtung, erkennt man, dass 
Euler auch in diesem Fall die Bestandteile des allgemeineren Konzepts der Gaußkrümmung

\begin{equation}
\label{eq: Def Gaußkrümmung}
    K_{G}:=k_1 \cdot k_2
\end{equation}
zwar zutage zu fördern vermochte, sie aber nicht wie Gauß verwoben. In § 8 (IV.) seiner Arbeit \cite{Ga28} fasst Gauß das Gesagte, bezugnehmend auf Gleichung (\ref{eq: Satz Euler Diffgeo}) folgendermaßen zusammen:\\

\textit{``Diese Schlussfolgerungen beinhalten nahezu all jenes, was der illustre Herr Euler über die Krümmung gekrümmter Flächen als erster gelehrt hat."}\\

Ebenfalls in § 8 (V.) gibt er  die Definition (\ref{eq: Def Gaußkrümmung}) und reicht folgende Erklärung in Form eines Theorems nach:\\

\textit{``Das Krümmungsmaß in jedwedem Punkt einer Fläche ist gleich dem Bruch, dessen Zähler die Einheit ist, der Nenner hingegen das Produkt der zwei extremalen Krümmungsradien in den Schnitten durch eine orthogonale Ebene ist."}\\

Zu dem heute als ``Theorema egregium"{} bezeichneten Satz gelangt er nach langwieriger Rechnung in § 12 und formuliert es in Worten wie folgt:\\

\textit{``Wenn eine gekrümmte Fläche in eine beliebige andere Fläche ausgebreitet wird, bleibt das Krümmungsmaß in den einzelnen Punkten unverändert."}

\newpage


\section{Herleitungen aus Eulers Formeln und Ideen}
\label{sec: Herleitungen aus Eulers Formeln und Ideen}

\epigraph{Noch viele Schätze sind in Eulers Werk zu heben, und wer Prioritäten jagen will, findet kein dankbareres Gefilde.}{Andreas Speiser}

Dieser den Schlusspunkt der vorliegenden Arbeit bildende Abschnitt ist der Anwendung der vorgetragenen Euler'schen Ideen gewidmet. Vorgestellt werden Korollare aus den Hauptergebnissen Eulers (Abschnitt \ref{subsec: Unmittelbare Korollare aus den vorgestellten Euler'schen Formeln}). Desweiteren soll, indem aus Eulers Gedanken bekanntere Formeln von Ramanujan abgeleitet werden, eine Verbindung zwischen diesen beiden Mathematikern hergestellt werden (Abschnitt \ref{subsec: Mit Eulers Ideen zu Formeln von Ramanujan}).

\subsection{Unmittelbare Korollare aus den vorgestellten Euler'schen Formeln}
\label{subsec: Unmittelbare Korollare aus den vorgestellten Euler'schen Formeln}

\epigraph{It appears to me that if one wants to make progress in mathematics, one should study the masters and not the pupils.}{Niels Henrik Abel}

Euler hat im Laufe seiner Karriere dermaßen zahlreiche Entdeckungen gemacht, dass es für ihn unmöglich gewesen ist, alle unmittelbaren Korollare aus selbigen ebenfalls mitzuteilen. Dies soll hier für ausgewählte im oberen Teil präsentierte Gegenstände nachgeholt werden. So wird aus seinem Ergebnis zur Lösung der einfachen Differenzengleichung (Abschnitt \ref{subsubsec: Eulers Lösung der einfachen Differenzengleichung}) ein weiterer Nachweis der Formel für die Summe der reziproken Potenzsummen abgeleitet (Abschnitt \ref{subsubsec: Herleitung Formel für die Potenzsummen der Reziproken}), weiter wird seine Theorie zur Auflösung von homogenen Differenzengleichungen mit linearen Koeffizienten (siehe Abschnitt \ref{subsubsec: Die Mellin--Transformierte bei Euler}) auf andere bekannte Orthogonalpolynome angewendet werden (\ref{subsubsec: Explizite Formeln aus Eulers Theorie zu Differenzengleichungen}), bevor endlich sehr bekannte sowie weniger bekannte Eigenschaften der Legendre--Polynome (Abschnitt \ref{subsubsec: Weitere Untersuchungen zu den Legendre-Polynomen}) aus den Euler'schen Formeln bewiesen werden.

\subsubsection{Herleitung Formel für die Potenzsummen der Reziproken}
\label{subsubsec: Herleitung Formel für die Potenzsummen der Reziproken}

\epigraph{Once a day [...] call yourselves to an account what new ideas, what new proposition or truth you have gained, what further confirmation of known truths, and what advances you have made in any part of knowledge.}{Isaac Watts}

Es ist förderlich, mit einem Beweis der Formeln für die Summen über die Potenzen der Reziproken der natürlichen Zahlen zu beginnen. Grundlage bildet die Euler'sche Formel zur Auflösung der einfachen Differenzengleichung (siehe Abschnitt \ref{subsubsec: Eulers Lösung der einfachen Differenzengleichung})

\begin{equation*}
    f(x+1)-f(x)=g(x).
\end{equation*}

Eulers Lösung lautet in korrigierter Fassung 

\begin{equation}
    \label{eq: Solution Simple Difference Equation}
    f(x) = \int g(x)dx - \dfrac{1}{2}g(x) + \sum_{k \in \mathbb{Z}\setminus \lbrace 0 \rbrace} e^{2k \pi ix} \int e^{-2 k \pi i x}g(x)dx,
\end{equation}
mithilfe welcher die Summen 

\begin{equation}
\label{eq: zeta(2n)}
    \zeta(2n):=\sum_{k=1}^{\infty} \dfrac{1}{k^{2n}}
\end{equation}
berechnen werden sollen\footnote{Dies wurde auch in der Arbeit \textit{``On Euler's Solution of the Simple Difference Equation"} (\cite{Ay23b}, 2023) aufgezeigt.}, was Euler vielerorts geleistet hat, insbesondere in \cite{E130}. Vermöge der mit elementaren Mitteln nachweisbaren Formel

\begin{equation*}
    e^{a x} \int e^{-ax}x^n dx = -\dfrac{1}{a^{n+1}}\sum_{k=0}^{n} \dfrac{n!}{k!}a^kx^k
\end{equation*}
vereinfacht sich (\ref{eq: Solution Simple Difference Equation}) für die Wahl $g(x)=x^n$ zu

\begin{equation}
    \label{eq: Solution x^n}
    f(x)= \dfrac{x^{n+1}}{n+1}-\dfrac{x^n}{2}- \sum_{k \in \mathbb{Z}\setminus \lbrace 0 \rbrace} \dfrac{1}{(2k \pi i)^{n+1}}\cdot \sum_{j=0}^n \dfrac{n!}{j!}(2k\pi i)^j x^j+h(x),
\end{equation}
wobei $h$ die Gleichung $h(x+1)=h(x)$ erfüllt. Andererseits wird dieselbe Differenzengleichung gleichfalls von

\begin{equation}
\label{eq: Faulhaber}
    f(x) = \sum_{k=1}^{x-1}k^n= \sum_{k=1}^{x}k^n - x^n = \dfrac{x^{n+1}}{n+1}+\dfrac{x^n}{2}+\dfrac{1}{n+1}\sum_{j=2}^{n} \binom{n+1}{j}B_j x^{n+1-j}-x^n
\end{equation}
\begin{equation*}
    = \dfrac{x^{n+1}}{n+1}-\dfrac{x^n}{2}+\dfrac{1}{n+1}\sum_{j=2}^{n} \binom{n+1}{j}B_j x^{n+1-j}
\end{equation*}
erfüllt, was die Faulhaber'sche Formel für die Potenzsummen darstellt. Ignoriert man nun die periodische Funktion $h$ in (\ref{eq: Solution x^n}), gelangt man nach einer kurzen Rechnung zu:

\begin{equation}
\label{eq: Ordered x^n}
    f(x)= \dfrac{x^{n+1}}{n+1}-\dfrac{x^n}{2}-  \sum_{j=0}^n \dfrac{n!}{j!} \sum_{k \in \mathbb{Z}\setminus \lbrace 0 \rbrace}(2k\pi i)^{j-(n+1)} x^j.
\end{equation}
Da sich diese Formel und (\ref{eq: Faulhaber}) lediglich um eine periodische Funktion unterscheiden, ist ein Koeffizientenvergleich der Potenzen von $x$ möglich. Mit

\begin{equation}
    \label{eq: Definition Coefficients B}
    B(n,j):= \dfrac{1}{n+1}\binom{n+1}{j}B_j \quad \text{für} \quad j \geq 2,
\end{equation}
wobei alle anderen Werte von $j$ den Ausdruck verschwinden lassen, sowie 

\begin{equation}
    \label{eq: Definition Coefficients A}
A(n,j):= - \dfrac{n!}{j!}\sum_{k \in \mathbb{Z}\setminus \lbrace 0 \rbrace}(2k\pi i)^{j-(n+1)}
\end{equation}
führt dieser Koeffizientenvergleich zu 

\begin{equation*}
     A(n,n+1-j)=B(n,j).
\end{equation*}
Mit den jeweiligen Definitionen wird dies zu

\begin{equation*}
     - \dfrac{n!}{(n+1-j)!}\sum_{k \in \mathbb{Z}\setminus \lbrace 0 \rbrace}(2k\pi i)^{-j}= \dfrac{1}{n+1}\binom{n+1}{j}B_j.
\end{equation*}
Umstellen nach der Summe gibt 

\begin{equation*}
    \sum_{k \in \mathbb{Z}\setminus \lbrace 0 \rbrace}(2k\pi i)^{-j}=- \dfrac{(n+1-j)!}{n!}\cdot \dfrac{1}{n+1}\cdot \dfrac{(n+1)!}{j!(n+1-j)!}B_j=-\dfrac{B_j}{j!}
\end{equation*}
oder auch 

\begin{equation*}
     \sum_{k \in \mathbb{Z}\setminus \lbrace 0 \rbrace}k^{-j} = -(2\pi i)^{j}\dfrac{B_j}{j!}.
\end{equation*}
Weil nun die ungeraden Potenzen von $j$ jeweils zu sich aufhebenden Termen führen, findet sich

\begin{equation}
\label{eq: Formula zeta(2n)}
    \zeta(2j)=\sum_{k=1}^{\infty} \dfrac{1}{k^{2j}}=\dfrac{(-1)^{j-1}(2\pi)^{2j} B_{2j}}{2(2j)!},
\end{equation}
was den Beweis abschließt.

\subsubsection{Explizite Formeln aus Eulers Theorie zu Differenzengleichungen}
\label{subsubsec: Explizite Formeln aus Eulers Theorie zu Differenzengleichungen}

\epigraph{Knowing is not enough. We must apply. Willing is not enough. We must do.}{Bruce Lee}

Hier werden  weitere konkrete Applikationen von Eulers Theorie zu Differenzengleichungen vorgestellt. Den Anfang bildet eine von Euler angegebene Gleichung, welche er selbst auf anderem als dem vorzustellendem Wege aufgelöst hat. Gefolgt wird dies von der Ableitung von alternativen Darstellungen für spezielle Orthogonalpolynome und Funktionen.

\paragraph{Eine Anwendung auf eine Formel von Euler}
\label{para: Eine Anwendung auf eine Formel von Euler}

In seiner Abhandlung \textit{``Solutio quaestionis curiosae ex doctrina combinationum"} (\cite{E738}, 1811, ges. 1779) (E738: ``Die Lösung einer interessanten Frage aus der Kombinatorik") gelangt Euler in der Diskussion eines Problems aus der Wahrscheinlichkeitsrechnung -- genauer die Anzahl der Fehlstände einer Menge von Zahlen -- zu einer Funktion, welche folgende Differenzengleichung erfüllt (siehe § VII der erwähnten Arbeit):

\begin{equation}
    \label{eq: Euler DG Fehlstände}
    \Pi(n)=(n-1)\left[\Pi(n-1)+\Pi(n-2)\right].
\end{equation}
Selbige ist eine homogene Differenzengleichung mit linearen Koeffizienten, welche sich demnach mit der in Abschnitt (\ref{subsubsec: Die Mellin--Transformierte bei Euler}) vorgestellten Methode auflösen lässt. Gemäß selbiger nehme man an, dass eine Lösung zu selbiger in der Form

\begin{equation*}
    \Pi(n)=\int\limits_{a}^b t^np(t)dt
\end{equation*}
 mit einer unbekannten Funktion $p(t)$ gegeben ist. Man hat  die Hilfsgleichung

\begin{equation*}
    \int t^n p(t)dt= (n-1)\int t^{n-1}p(t)dt+(n-1)\int t^{n-2}p(t)+t^nq(t)
\end{equation*}
mit einer weiteren unbekannten Funktion $q(t)$ zu betrachten, welche differenziert gibt:

\begin{equation*}
    t^{n-2} p(t)= (n-1)t^{n-1}p(t)+(n-1)t^{n-2}p(t)+nt^{n-1}q(t)+t^nq'(t).
\end{equation*}
Man teile diese durch $t^n$:

\begin{equation*}
    t^2 p(t)= (n-1)t^{1}p(t)+(n-1)p(t)+ntq(t)+t^2q'(t).
\end{equation*}
Aus dem Koeffizientenvergleich leitet man folgendes System von gekoppelten Differentialgleichungen ab:

\begin{equation*}
\text{Aus}~~n^0: ~~ t^2p(t)=-tp(t)-p(t)+t^2q'(t) \quad \text{und aus}~~ n^1:~~ 0=tp(t)+p(t)+tq(t).
\end{equation*}
Dieses System wird gelöst von:

\begin{equation*}
    p(t)=C \cdot e^{-t} \quad \text{und} \quad q(t)=C\cdot e^{-t}\dfrac{t+1}{t}
\end{equation*}
mit $C\neq 0$, womit

\begin{equation*}
    t^n\cdot q(t)= C\cdot t^{n-1}e^{-t}\cdot(t+1)=0
\end{equation*}
zur Auffindung der Grenzen zu lösen ist. Man findet

\begin{equation*}
    1) \quad t=0, \quad 2) \quad t=-1 \quad \text{und} \quad t=\infty.
\end{equation*}
Dementsprechend kommen als mögliche Lösungen von (\ref{eq: Euler DG Fehlstände}) die folgenden in Frage:

\begin{equation*}
    1.) ~~ f_1(n)= C\int\limits_{-1}^{0} t^n e^{-t}dt, \quad 2.) ~~ f_2(n)= C\int\limits_{-1}^{\infty} t^n e^{-t}dt 
\end{equation*}
und
\begin{equation*}
    3.) ~~ f_3(n)= C\int\limits_{0}^{\infty} t^n e^{-t}dt.
\end{equation*}
Nun bleibt an dieser Stelle nichts anderes übrig, als die Lösungen nacheinander auf ihre Verträglichkeit mit der Problemstellung hin zu prüfen. Hierzu benötigt man die Anfangsbedingungen $\Pi(1)=0$ und $\Pi(2)=1$, welche Euler in seiner Arbeit \cite{E738} ebenfalls angibt. $f_3$ kann also nicht als Lösung in Betracht kommen, da $f_3(1)=C\cdot 1$ und somit müsste $C=0$ sein. Nimmt man nun $f_1$, ist freilich zunächst

\begin{equation*}
    f_1(1)=C\int\limits_{-1}^0 t^1e^{-t}dt =C \cdot (-1).
\end{equation*}
Allerdings müsste hier wieder $C=0$ sein, um die Anfangsbedingung $f_1(1)=0$ zu gewährleisten, sodass diese Lösung ebenfalls nicht in Frage kommt.\\

Probiert man nun $f_2(n)$, führt dies zu

\begin{equation*}
    f_2(1)=C\int\limits_{-1}^{\infty} t^1 e^{-t}dt = C\cdot 0,
\end{equation*}
sodass die erste Bedingung automatisch erfüllt ist. Unglücklicherweise ist die Konstante $C$ hieraus nicht zu bestimmen, weshalb man sich der zweiten Gleichung bediene:

\begin{equation*}
    f_2(2)= C\int\limits_{-1}^{\infty} t^2 e^{-t}dt =C\cdot e.
\end{equation*}
Wegen $f_2(2)=1$ muss $C=\frac{1}{e}$ sein, was zu nachstehendem Ausdruck für (\ref{eq: Euler DG Fehlstände}) führt:

\begin{equation}
\label{eq: Fehlstände explizit}
    \Pi(n)= f_2(n)= \dfrac{1}{e}\cdot \int\limits_{-1}^{\infty} t^n e^{-t}dt.
\end{equation}
Diese Formel scheint in der Tat richtig zu sein, denn man findet durch direktes Rechnen:

\begin{equation*}
    f_2(3)= \dfrac{2e}{e}=2, ~~ f_2(3)=\dfrac{9e}{e} =9, ~~ f_2(5)=\dfrac{44e}{e} =44, ~~ f_2(6)=\dfrac{265e}{e} =265,
\end{equation*}
was gerade die von Euler in \cite{E738} mitgeteilten Werte sind. \\

\begin{figure}
    \centering
  \includegraphics[scale=0.7]{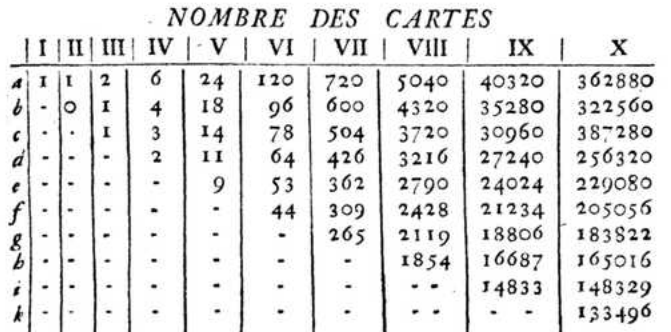}
    \caption{ Eulers Tabelle aus \cite{E201} zu den Fehlständen beim Spiel ``Rencontre".}
    \label{fig:E201Karten}
\end{figure}

In seiner Abhandlung \textit{``Calcul de la probabilite dans le jeu de rencontre"} (\cite{E201}, 1753) (E201: ``Die Wahrscheinlichkeitsrechnung im Spiel Rencontre") ist Euler beim Betrachten des Kartenspiels ``Rencontre"{} zum selben Problem der Fehlstände gelangt. In der späteren Arbeit \cite{E738} hat er diese Verbindung allerdings nicht erwähnt. Unabhängig davon gelangt er in § XLIV des zuerst genannten Papiers zur expliziten Formel:

\begin{equation}
\label{eq: Euler Fehlstände explizit}
    \Pi(n) = n! \sum_{k=0}^n (-1)^k \dfrac{1}{k!}.
\end{equation}
Die Gleichheit dieser Formel mit (\ref{eq: Fehlstände explizit}) lässt sich wie folgt nachweisen. Zunächst bemerke man die für $k>0$ gültige Formel: 

\begin{equation*}
    \int\limits_{-1}^{\infty} e^{-kt}dt = \left[-\dfrac{e^{-kt}}{k}\right]_{-1}^{\infty}= \dfrac{e^k}{k}.
\end{equation*}
$n$-faches Differenzieren gibt für die linke Seite

\begin{equation*}
    \dfrac{d^n}{dk^n}  \int\limits_{-1}^{\infty} e^{-kt}dt 
\end{equation*}
\begin{equation*}
    =  \int\limits_{-1}^{\infty} \dfrac{d^n}{dk^n} e^{-kt}dt =  \int\limits_{-1}^{\infty} (-t)^n e^{-kt}dt 
\end{equation*}
und für die rechte Seite findet man unter Verwendung der Leibniz--Regel zur Ableitung von Produkten

\begin{equation*}
    \dfrac{d^n}{dk^n}e^{-k}\cdot (k^{-1})= \sum_{l=0}^n \binom{n}{l} \dfrac{d^l}{dk^l}e^{k}\cdot \dfrac{d^{n-l}}{dk^{n-l}}(k^{-1}) 
\end{equation*}
\begin{equation*}
    = \sum_{l=0}^n \binom{n}{l} e^k \cdot (-1)^{n-l}(n-l)!k^{-(n+l+1)}.
\end{equation*}
Demnach zusammengefasst:

\begin{equation*}
    \int\limits_{-1}^{\infty} (-t)^n e^{-kt}dt =  \sum_{l=0}^n \binom{n}{l} e^k \cdot (-1)^{n-l}(n-l)!k^{-(n+l+1)}.
\end{equation*}
Somit hat man für $k=1$ und nach Division beider Seiten durch $(-1)^n$

\begin{equation*}
    \int\limits_{-1}^{\infty} t^n e^{-t}dt = e\cdot \sum_{l=0}^n \binom{n}{l}  \cdot (-1)^{l}(n-l)!
\end{equation*}
oder alternativ

\begin{equation*}
    \dfrac{1}{e}\cdot  \int\limits_{-1}^{\infty} t^n e^{-t}dt= \sum_{l=0}^n \binom{n}{l}  \cdot (-1)^{l}(n-l)!= \sum_{l=0}^n \dfrac{(-1)^l\cdot n!}{l!(n-l)!}\cdot (n-l)!= n!\sum_{l=0}^n\dfrac{(-1)^l}{l!}.
\end{equation*}
Die linke Seite ist die eben gezeigte Lösung (\ref{eq: Fehlstände explizit}) und die rechte Seite ist Eulers Formel (\ref{eq: Euler Fehlstände explizit}), womit  die Gleichheit der beiden Formeln nachgewiesen ist. In \cite{Du99} findet man einen alternativen Beweis von (\ref{eq: Euler Fehlstände explizit}) mittels Induktion.

\paragraph{Orthogonale Polynome}

Die Theorie der Orthogonalpolynome wird heute angesehen, mit der wegweisenden Arbeit von Stieltjes \textit{``Recherches sur les fractiones continues"} (\cite{St94}, 1894) (``Untersuchungen über Kettenbrüche") entscheiden vorangebracht worden zu sein. Stieltjes führte in dieser Arbeit das nach ihm benannte Stieltjes--Integral ein und verband  das Momentenproblem mit Kettenbrüchen und auch den Orthogonalpolynomen. Obschon Stieltjes seinen Vorgänger Euler weit in seinen Erkenntnissen zu diesen Gegenständen übersteigt, lassen sich in den Euler'schen Arbeiten auch explizite Formeln für Orthogonal--Polynome herauslesen, wovon oben (Abschnitt \ref{subsubsec: Den Kontext betreffend: Die Legendre Polynome}) mit den Legendre--Polynomen schon ein Beispiel umfassender diskutiert worden ist\footnote{Man beachte in diesem Kontext den Lehrsatz von Favard (1902--1965) aus seiner Arbeit \textit{``Sur les polynomes de Tchebicheff"} (\cite{Fa35}, 1935) (``Über die Tschebycheff'schen Polynome"), dass für Polynome $p_n(x)$, welche nachstehender Rekurrenzrelation genüge leisten:

\begin{equation*}
    p_n(x)=(A_n x+B_n)p_{n-1}(x)+C_n p_{n-2}(x)
\end{equation*}
gewiss ist, dass sie einen Satz von Orthogonalpolynomen bilden. Überdies ist die Umkehrung dieser Aussage ebenfalls wahr. Demnach erkennt man in der für Euler in seinen Untersuchungen zu Kettenbrüchen führenden Gleichung (\ref{eq: Linear Difference Equation}) auch eine für entsprechende Orthogonalpolynome. Überdies lassen sich seine Beiträge zur Mellin--Transformierten (\ref{subsubsec: Die Mellin--Transformierte bei Euler}) als Beiträge zum Momenten--Problem umdeuten. Stieltjes hat  demnach in seiner Arbeit \cite{St94} die drei scheinbar unverwandten Gebiete von Kettenbrüchen, Momentproblemen -- ein Begriff, der zu Eulers Lebenszeiten noch gar nicht existierte -- und Orthogonalpolynomen synthetisiert.}. \\

An dieser Stelle wird mit den Tschebycheff--Polynomen ein Beispiel für eine Differenzengleichung mit konstanten Koeffizienten präsentiert, sodass die Ergebnisse aus Abschnitt (\ref{subsubsec: Eulers Lösung der einfachen Differenzengleichung}) hier zum Tragen kommen. Die Hermite--Polynome erfüllen hingegen eine Differenzengleichung mit linearen Koeffizienten, womit sie von dem in Abschnitt (\ref{subsubsec: Die Mellin--Transformierte bei Euler}) Diskutierten erfasst werden. Gefolgt  wird dies von den Laguerre--Polynomen -- ein Beispiel, in welchem die Euler'sche Methode an seine Grenzen stößt. Schließlich soll mit den Omega--Funktionen,  noch ein Beispiel einer etwas allgemeineren Funktionenklasse mit Eulers Ansatz behandelt werden.

\subparagraph{Die Tschebycheff-Polynome}

Die Tschebycheff--Polynome erfüllen die folgende Differenzengleichung:

\begin{equation}
    \label{eq: Diff. Tschebycheff}
    T_{n+1}(x)= 2x T_n(x)-T_{n-1}(x),
\end{equation}
also  eine Differenzengleichung mit konstanten Koeffizienten in der Variable $n$.  Eulers Methode befolgend, ist diese Gleichung zu betrachten:

\begin{equation*}
    0=e^z-2x+e^{-z}
\end{equation*}
und nach $z$ aufzulösen. Man findet diese unendlichen vielen Lösungen:

\begin{equation*}
    z_k= \log (x\pm \sqrt{x^2-1}) + 2k\pi i \quad \text{mit} \quad k \in \mathbb{Z}.
\end{equation*}
Eulers Ansatz führt also zur folgenden Lösung von (\ref{eq: Diff. Tschebycheff}):

\begin{equation*}
    T_n(x) = \sum_{k=-\infty}^{\infty} c_k e^{n(\log(x+\sqrt{x^2-1})+2k\pi i)}+b_k e^{n(\log(x-\sqrt{x^2-1})+2k\pi i)}
\end{equation*}
mit beliebigen Konstanten $c_k$ und $b_k$. Dies vereinfacht sich zu:

\begin{equation*}
    (x+\sqrt{x^2-1})^n \sum_{k=-\infty}^{\infty} c_k  e^{2k\pi in} +(x-\sqrt{x^2-1})^n \sum_{k=-\infty}^{\infty} b_k  e^{2k\pi in},
\end{equation*}
sodass man zur folgenden Lösung von (\ref{eq: Diff. Tschebycheff}) gelangt

\begin{equation*}
  T_{n,E}(x)=   f(n)(x+\sqrt{x^2-1})^n+g(n) (x-\sqrt{x^2-1})^n,  
\end{equation*}
wenn $f(n)=f(n+1)$ und $g(n)=g(n+1)$ gilt. Der zusätzliche Index $E$ zeigt an, dass es eine aus der Euler'schen Vorgehensweise derivierte Lösung ist.\\

Vergleicht man dieses Ergebnis mit der bekannten Darstellung:

\begin{equation}
\label{eq: Tscheby explizit}
    T_n(x)= \dfrac{(x+\sqrt{x^2-1})^n+(x-\sqrt{x^2-1})^n}{2},
\end{equation}
ergibt sich $T_n(x)$ aus $T_{n,E}(x)$ für die spezielle Wahl $f(n)=g(n)=\frac{1}{2}$. Man kann dies auch aus $T_{n,E}(x)$ selbst ableiten, sofern man das Polynom $T_1(x)=x$ betrachtet. Dann müsste gelten

\begin{equation*}
    x= T_{1,E}=f(1)(x+\sqrt{x^2-1})+g(1)(x-\sqrt{x^2-1})
\end{equation*}
\begin{equation*}
    =(f(0)+g(0))x+(f(0)-g(0))\sqrt{x^2-1},
\end{equation*}
wobei $f(n+1)=f(n)$ und $g(n)=g(n+1)$  benutzt wurde. Ein Koeffizientenvergleich führt nun zu den Gleichungen:

\begin{equation*}
    1.)~~ f(0)+g(0)=1 \quad \text{und} \quad 2.)~~ f(0)-g(0)=0.
\end{equation*}
Die Lösungen dieses Gleichungssystems sind $f(0)=g(0)=\frac{1}{2}$. Somit reduziert sich für natürliche und ganzzahlige $n$ die Euler'sche Lösung $T_{n,E}(x)$ genau auf die explizite Formel für die Tschebyscheff Polynome (\ref{eq: Tscheby explizit}).

\subparagraph{Die Hermite-Polynome}

Die Hermite-Polynome erfüllen die Differenzengleichung

\begin{equation}
    \label{Diff. Hermite}
    H_{n+1}(x)=2x H_n(x)-2nH_{n-1}(x)
\end{equation}
mit $H_0(x)=1$ und $H_1(x)=2x$. Nach der Euler'schen Methode setze man

\begin{equation*}
    \int t^{n+1}p(t)dt= 2x\int t^np(t)dt-2n\int t^{n-1}p(t)dt+t^nq(t),
\end{equation*}
wobei  die Funktionen $p(t)$ und $q(t)$ zu bestimmen sind. Differenzieren der Gleichungen und anschließendes Vereinfachen ergibt

\begin{equation*}
    t^2p(t)=2xtp(t)-2np(t)+nq(t)+tq'(t).
\end{equation*}
Der Vergleich der Koeffizienten der jeweiligen Potenzen von $n$ impliziert:

\begin{equation*}
    1)~~ t^2p(t)=2xtp(t)+tq'(t) \quad \text{und} \quad 2)~~  0=-2p(t)+q(t).
\end{equation*}
Dieses System von gekoppelten Differenzialgleichungen besitzt die Lösungen:

\begin{equation*}
    p(t)=\dfrac{1}{2} C(x)e^{\left(\frac{t^2}{4}-xt\right)} \quad \text{sowie} \quad q(t)=C(x)e^{\left(\frac{t^2}{4}-xt\right)},
\end{equation*}
wobei $C(x)$ eine Funktion $\neq 0$ ist. Um nun die Grenzen zu bestimmen, hat man nach Eulers Vorgabe die Gleichung $t^n\cdot q(t)=0$ zu betrachten, was $t_{1,2}=\pm i \infty$ bedeutet, sodass zunächst folgender Ansatz gemacht werden kann:

\begin{equation*}
    H_{n}(x)= 2C(x) \int\limits_{-i\infty}^{i\infty} t^n e^{\left(\frac{t^2}{4}-xt\right)}dt.
\end{equation*}
Um sich der imaginären Grenzen zu entledigen,  setze man $t=iy$; damit wird 

\begin{equation*}
     H_{n}(x)= 2C(x) i^{n+1} \int\limits_{-\infty}^{\infty} y^n e^{-\frac{y^2}{4}-iyx}dy.
\end{equation*}
Aus $H_0(x)=1$ muss demnach gelten:

\begin{equation*}
    H_0(x)=1=2C(x) i \int\limits_{-\infty}^{\infty} y^n e^{-\frac{y^2}{4}-iyx}dy= 2C(x)i\cdot 2\sqrt{\pi}e^{-x^2},
\end{equation*}
sodass

\begin{equation*}
    C(x)=\dfrac{e^{x^2}}{4i\sqrt{\pi}}. 
\end{equation*}
Einsetzen dieses Wertes von $C(x)$  in den Ansatz gibt als finalen Ausdruck:

\begin{equation*}
    H_n(x)= \dfrac{e^{x^2}}{2\sqrt{\pi}}\int\limits_{-\infty}^{\infty} \left(iy\right)^n e^{-\frac{y^2}{4}-iyx}dy.
\end{equation*}

\subparagraph{Die Laguerre-Polynome}

Anhand der Laguerre--Polynome sei zuletzt ein Exempel angegeben,  beim welchem Eulers Ansatz an seine Grenzen stößt. Die Laguerre--Polynome erfüllen die Differenzengleichung

\begin{equation}
    \label{eq: Diff. Laguerre}
    (n+1)L_{n+1}(x)= (2n+1-x)L_n(x)-nL_{n-1}(x).
\end{equation}
Vermöge des Ansatzes 

\begin{equation*}
    (n+1)\int t^{n+1}p(t)dt= (2n+1-x)\int t^n p(t)dt-n \int t^{n-1}p(t)dt+t^nq(t)
\end{equation*}
gelangt man zum folgenden System von Differentialgleichungen für die zu findenden Funktion $p(t)$ und $q(t)$:

\begin{equation*}
    1.) ~~ t^2p(t)=2tp(t)-p(t)+q(t) \quad \text{und} \quad 2.)~~ tp(t)=(1-x)p(t)+q'(t).
\end{equation*}
Die Lösung lautet

\begin{equation*}
    p(t)= \dfrac{C(x)}{t-1}e^{-\frac{x}{t-1}} \quad \text{und} \quad q(t)= C(x)(t-1)e^{-\frac{x}{t-1}}
\end{equation*}
mit einer Funktion $C(x)\neq 0$. \\

Will man nun die Grenzen der Integration ermitteln, hat man die Gleichung

\begin{equation*}
    0=q(t)=C(x)(t-1)e^{-\frac{x}{t-1}}
\end{equation*}
nach $t$ aufzulösen. Jedoch findet man hier lediglich die eine Lösung $t=1$, obwohl man zweier zur Anwendung von Eulers Methode bedürfte. Mithilfe der komplexen Analysis findet man

\begin{equation*}
    L_n(x) = \dfrac{1}{2\pi i}\oint_{C} \dfrac{e^{-\frac{xt}{t-1}}}{(1-t)t^{n+1}}dt
\end{equation*}
wobei $C$ ein geschlossener Weg um den Ursprung ist, welcher $t=1$ nicht mit einschließt. \\

Um diese Vorstellung der Euler'schen Ideen nicht auf einem Missklang enden zu lassen, sei die Bemerkung gestattet, dass die Euler'sche Herangehensweise an Differenzengleichungen, insbesondere der Ansatz aus (\ref{eq: Aux}), dem modernen Ansatz im Studium der Feynman'schen Integrale in der Physik entspricht.  Die Arbeit \textit{``Decomposition of Feynman Integrals on the Maximal Cut by Intersection Numbers"} (\cite{Fr19}, 2019) war eine der ersten, welche sich Methoden aus der Mathematik bedient, welche ihre Wurzel in der Euler'schen Methode haben. Das Buch \textit{``Feynman Integrals: A Comprehensive Treatment for Students and Researchers"} (\cite{We22}, 2022) stellt die Theorie und Methodik umfassend dar. Man konsultiere auch das entsprechende Kapitel in \cite{Ay21a} für einen Bezug der Überlegungen in der Physik zu den Euler'schen Formeln.

\paragraph{Omega--Funktionen}

Die Beispiele der Orthogonalpolynome haben gezeigt, dass die Euler'sche Methode zur Auflösungsmethode einiges zu leisten vermag, jedoch auch an ihre Grenzen stößt, wenn die zu behandelnde Differenzengleichung entsprechende unliebsame Eingeschaften aufweist. Daher soll abschließend noch ein Beispiel diskutiert werden, wo Eulers Ansatz die gewünschten Ergebnisse zutage fördert. Die zum Verfassungszeitpunkt dieser Arbeit noch nicht publizierte Arbeit \textit{``Difference Equations and Omega Functions"} (\cite{Pe25}, vorauss. 2025)\footnote{Die Arbeit von Marco--Perez wurde dem Verfasser der gegenwärtigen Ausarbeitung mit einem Gesuch nach Revision vom \textit{American Journal of Mathematics} zugesandt und anschließend auch mit Perez selbst diskutiert.} hat die sogenannten $\Omega$--Funktionen zum Gegenstand, welche mit funktionentheoretischen Mitteln untersucht werden. Sie werden in erwähnter Abhandlung über ihre Differenzengleichung eingeführt. Sie erfüllen diese:

\begin{equation*}
    s\Omega(s)= \sum_{k=1}^{d} \alpha_k \Omega(s+d),
\end{equation*}
sodass sie sich mit der Euler'schen Methode behandeln lassen. Mit dem Ansatz

\begin{equation}
    \label{eq: Ansatz Omega}
    \Omega(s)= \int\limits_{a}^{b} t^{s-1}p(t)dt
\end{equation}
und der Hilfsgleichung

\begin{equation*}
    s\int t^{s-1}p(t)dt = \sum_{k=1}^d \alpha_k \int t^{s-1+k}p(t)dt +t^{s}q(t)
\end{equation*}
mit zu bestimmenden Funktionen $p(t)$ und $q(t)$ und zu eruierenden Integrationsgrenzen $a$ und $b$ wird man analog zu den vorherigen Beispielen verfahrend zu folgenden System von Gleichungen für $p(t)$ und $q(t)$ geleitet:

\begin{equation*}
    1.~~ p(t)=q(t) \quad \text{und} \quad 0= \sum_{k=1}^{d} \alpha_k t^k p(t) +tq'(t).
\end{equation*}
Dieses System besitzt die Lösung

\begin{equation*}
    p(t)=q(t)= C \cdot e^{\left(-\sum_{k=1}^d \alpha_k \frac{t^k}{k}\right)},
\end{equation*}
wobei die Integrationskonstante $C$ nicht verschwindend gewählt werden kann. Für die Grenzen ist demnach die Gleichung

\begin{equation*}
    t^sq(t)= t^s \cdot e^{\left(-\sum_{k=1}^d \alpha_k \frac{t^k}{k}\right)}
\end{equation*}
zu betrachten. Die Lösung $t=0$ ist leicht ersichtlich. Für eine weitere $t_0$ kann man an dieser Stelle nur sagen, dass sie so zu wählen ist, dass das Polynom $\sum_{k=1}^d \alpha_k \frac{t^k}{k}$ unendlich wird. Dies hängt insbesondere von der Natur des Koeffizienten $\alpha_k$ ab. Aus der Funktiontheorie ist jedoch bekannt, dass ein Polynom als unbeschränkte holomorphe Funktion auch unendlich werden kann, wenn $t=\infty$ als Wert zugelassen wird. Die Existenz eines solchen $t_0$ kann demnach garantiert werden, sodass man mit dem Ansatz (\ref{eq: Ansatz Omega}) in jedem Fall eine Darstellung der Form

\begin{equation*}
    \Omega(s)= C \int\limits_{0}^{t_0} t^{s-1} \cdot e^{\left(-\sum_{k=1}^d \alpha_k \frac{t^k}{k}\right)}dt
\end{equation*}
findet. Eine weitere Untersuchung soll aber hier nicht erfolgen, zumal sie in der bald erscheinenden Arbeit \cite{Pe25} nachgereicht werden würden. Ebenso wird der Vergleich der Euler'schen Methoden zur Auflösung von Differenzengleichungen mit denen von N\o rlund aus seinem Buch \textit{``Vorlesungen über Differenzenrechnung"} (\cite{No24}, 1924) auf einen zukünftigen Zeitpunkt verlegt. Stattdessen sollen noch die Bessel'schen Funktionen eine Behandlung mit dieser Methode finden.

\paragraph{Besselfunktionen}

Bekanntermaßen, man konsultiere etwa das Buch \textit{``A Treatise on the Theory of Bessel Functions"} (\cite{Wa95}, 1995), erfüllen die Bessel'schen Funktionen, definiert über die Reihe

\begin{equation}
    \label{eq: Def Bessel}
    J_{\alpha}(x):= \sum_{n=0}^{\infty}\dfrac{(-1)^n}{n!\Gamma(n+\alpha+1)}\left(\dfrac{x}{2}\right)^{2n+\alpha},
\end{equation}
folgende Rekursionsvorschrift im Index $\alpha$:

\begin{equation}
\label{eq: Bessel Rek}
    2\alpha \cdot J_{\alpha}= x(J_{\alpha-1}+J_{\alpha+1}),
\end{equation}
wobei das Argument der Kürze wegen fortgelassen wurde. Nun mache man zur Lösung dieser entsprechend Euler'scher Vorgabe hier den Ansatz

\begin{equation}
\label{eq: Ansatz Bessel}
    J_{\alpha} = \int\limits_{a}^{b}e^{-\alpha t}p(t)dt,
\end{equation}
sodass im Gegensatz zu oben (Abschnitt \ref{subsubsec: Die Mellin--Transformierte bei Euler}) statt der Mellin--Transformierten eine Laplace--Transformierte als Lösungsform gewählt wird\footnote{Dies bedeutet für die Methode an sich keinen Unterschied, zumal die beiden Transformierten vermöge einer Substitution ineinander überführt werden können. Jedoch wird bei praktischen Rechnungen bisweilen der Aufwand für entsprechend anders gewähltem Ansatz merklich reduziert. Die Bessel'schen Funktionen statuieren ein Exempel dessen. Euler selbst scheint in seinen Untersuchungen stets den Ansatz des Mellin'schen Typus' gewählt zu haben.}. Die entsprechende Hilfsgleichung lautet entsprechend:

\begin{equation*}
    2\alpha \int e^{-\alpha t}p(t)dt =x \left(\int e^{-(\alpha-1) t}p(t)dt+\int e^{-(\alpha+1) t}p(t)dt\right)+e^{-\alpha t}q(t)
\end{equation*}
mit zu bestimmenden Funktionen $p(t)$ und $q(t)$, für welche man, wie oben verfahrend, folgende Gleichungen ableitet:

\begin{equation*}
    1. ~~ 2p(t)=-q(t) \quad \text{sowie} \quad  0= x(e^t+e^{-t})p(t)+q'(t).
\end{equation*}
Eine kurze Rechnung ergibt für $q(t)$ und $p(t)$:

\begin{equation*}
    p(t)= Ce^{x\sinh(t)} \quad \text{und} \quad q(t)=-2Ce^{x\sinh(t)}, 
\end{equation*}
wobei $C$ eine nicht verschwindende Konstante ist. Demnach sind die Grenzen der Integration aus der Gleichung

\begin{equation*}
    0= e^{-\alpha t+x \sinh(t)}
\end{equation*}
zu ermitteln. Demnach muss das Argument der Exponentialfunktion $-\infty$ werden. Für die weitere Betrachtung seien $\alpha, x>0$. Dann ist eine erste Lösung mit $t=-\infty$ gegeben, zumal das Wachstum von $\sinh(t)$ das der linearen Funktion $-\alpha t$ dominiert. Zur Ermittlung einer weiteren bemerke man die elementare Identität

\begin{equation*}
    \sinh(u \pm i \pi) = - \sinh(u)= \sinh(-u),
\end{equation*}
sodass $t=\infty \pm i \pi$ ebenfalls Lösungen sind. Insgesamt führt der Ansatz (\ref{eq: Ansatz Bessel}) also zunächst zu

\begin{equation}
\label{eq: Eigenlich Hankel}
    J_{\alpha}^{1,2} = C_{1,2} \int\limits_{-\infty}^{\infty \pm i \pi}e^{-\alpha t +x \sinh(t)}dt,
\end{equation}
woraus allerdings die Konstanten $C_{1,2}$ schwer zu ermitteln sind. Darum behelfe man sich durch Betrachtung der Summe der beiden Lösungen, welche natürlich immer noch (\ref{eq: Bessel Rek}) Genüge leistet. Zudem sei $C_1=-C_2 =C$. Dementsprechend ist zunächst

\begin{equation*}
    J_{\alpha}^{1}+J_{\alpha}^{2}= C \int\limits_{-\infty}^{\infty + i \pi}e^{-\alpha t +x \sinh(t)}dt- C \int\limits_{-\infty}^{\infty - i \pi}e^{-\alpha t +x \sinh(t)}dt,
\end{equation*}
welche sich mit elementaren Regeln der Integralrechnung zu

\begin{equation*}
     C \int\limits_{\infty- i\pi}^{\infty + i \pi}e^{-\alpha t +x \sinh(t)}dt
\end{equation*}
kontrahieren lässt. Eine leichte Anwendung des Residuensatzes und Hinzunahme der Periodizität von $\sinh(t)$ gestattet\footnote{Man integriere entgegen des Uhrzeigersinnes über das Rechteck mit den Eckpunkten $0-i\pi$, $\infty-i \pi$, $\infty+i \pi$, $0+i \pi$. Der Residuensatz gibt dieses Integral als $=0$ aus. Aus der Periodizität von $\sinh(t)$ und $e^{i\alpha t}$ folgert man indes 
\begin{equation*}
    \int\limits_{i\pi}^{\infty+i\pi}e^{x \sinh (t)-\alpha t}dt= \int\limits_{-i\pi}^{\infty-i\pi}e^{x \sinh (t)-\alpha t}dt,
\end{equation*}
sodass die unendlichen Teilwege sich beim angezeigten Rundwege gegenseitig aufheben.}, für ganzzahliges $\alpha$, die Reduktion des letzten auf das folgende Integral:

\begin{equation*}
    C \int\limits_{-i\pi}^{+i\pi} e^{x\sinh(t)-\alpha t}dt.
\end{equation*}
Hier lasse man schließlich $t=i\omega$ werden, sodass man zu nachstehendem Integral gelangt:

\begin{equation*}
    F_{\alpha}:= Ci \int\limits_{-\pi}^{\pi}e^{ix \sin (\omega)-i \omega \alpha}d\omega.
\end{equation*}
Die einfachere Gestalt erlaubt nun die Ermittelung der Konstanten $C$ durch den Vergleich mit der Bessel--Funktion in (\ref{eq: Def Bessel}). Man nehme $\alpha=0$ an und betrachte zunächst $F_{\alpha}$ für diesen Wert, sprich, das Integral

\begin{equation*}
    Ci  \int\limits_{-\pi}^{\pi}e^{ix \sin (\omega)}d\omega,
\end{equation*}
welches keine Integration in elementaren Funktionen zulässt. Daher entwickle man den Integranden in eine Potenzreihe, wodurch man hat

\begin{equation*}
    F_{0}= Ci\int\limits_{-\pi}^{\pi}e^{ix \sin (\omega)}d\omega = \int\limits_{-\pi}^{\pi} \sum_{n=0}^{\infty}\dfrac{(ix\sin(\omega))^n}{n!} d\omega,
\end{equation*}
sodass sich nach Vertauschung von Summe und Integral folgende Potenzreihe ergibt:

\begin{equation}
\label{eq: Power F_0}
    F_{0}=Ci \sum_{n=0}^{\infty} \left(\int\limits_{-\pi}^{\pi}\sin^n(\omega)d\omega\right) \dfrac{(ix)^n}{n!}.
\end{equation}
Es verbleibt demnach die Berechnung des Integrals, welches gleichsam der Koeffizient der Reihe ist. Für ungerades $n$ verschwindet besagtes als ein Integral über die Periode einer punktsymmetrischen Funktion. Für gerades $n$ findet man vermöge der Formel $\sin (\omega)=\frac{e^{i\omega}-e^{-i\omega}}{2i}$ und dem binomischen Lehrsatz in leichter Rechnung

\begin{equation*}
    \int\limits_{-\pi}^{\pi} \sin^{2n}(\omega)d\omega = 2\pi \cdot \dfrac{1}{2^{2n}} \binom{2n}{n},
\end{equation*}
womit (\ref{eq: Power F_0}) nachstehende Form annimmt:

\begin{equation*}
    F_{0}=2\pi Ci \sum_{n=0}^{\infty} \dfrac{1}{2^{2n}} \binom{2n}{n} \cdot \dfrac{(-1)^nx^{2n}}{(2n)!}= 2\pi Ci \sum_{n=0}^{\infty}  \dfrac{(-1)^n}{(n!)^2}\left(\dfrac{x}{2}\right)^{2n}.
\end{equation*}
Ein Vergleich mit $J_0(x)$ aus (\ref{eq: Def Bessel}) weist die Konstante $C$ als $=\frac{1}{2\pi i}$ aus, sodass man aus dem Ansatz (\ref{eq: Ansatz Bessel}) schließlich

\begin{equation*}
    J_{\alpha}(x)= \dfrac{1}{2\pi} \int\limits_{-\pi}^{\pi} e^{ix \sin (\omega)-i \alpha \omega}d\omega
\end{equation*}
erlangt, sofern $\alpha$ ganzzahlig ist. In dieser Form hat bereits Bessel (1784--1846) selbst in seiner Arbeit \textit{``Untersuchung des Theils der planetarischen Störungen, welcher aus der Bewegung der Sonne entsteht"} (\cite{Be24}, 1824)  die heute sogenannten Bessel'schen Funktionen erster Art aus (\ref{eq: Def Bessel}) angegeben und verwendet. Nach Ermittlung der Konstante $C$ lassen sich die beiden Zwischenlösungen aus (\ref{eq: Eigenlich Hankel}) als

\begin{equation*}
   H^{(1)}_{\alpha}(x)=\dfrac{1}{i\pi} \int\limits_{-\infty}^{\infty + i \pi}e^{-\alpha t +x \sinh(t)}dt \quad \text{und} \quad H^{(2)}_{\alpha}(x)= -\dfrac{1}{i\pi }\int\limits_{-\infty}^{\infty - i \pi}e^{-\alpha t +x \sinh(t)}dt
\end{equation*}
identifizieren, welche heute als Hankel'sche Funktionen bezeichnet werden. Man  konsultiere diesbezüglich wieder  Watsons  (1886--1965) Buch \cite{Wa95}, insbesondere Abschnitt 6.21., wo man die Verbindung

\begin{equation*}
    H^{(1)}_{\alpha}(x)+H^{(2)}_{\alpha}(x)=2 J_{\alpha}(x)
\end{equation*}
herauslesen zu den Bessel'schen Funktionen erster Art aus (\ref{eq: Def Bessel}) kann\footnote{Hankel hat die nach ihm benannten Funktionen -- in leicht anderer Form -- unter Verwendung von Kurvenintegralen in der komplexen Ebene in seiner Arbeit \textit{``Die Zylinderfunktionen erster und zweiter Art"} (\cite{Ha69a}, 1869) eingeführt und studiert.}.

\subsubsection{Weitere Untersuchungen zu den Legendre--Polynomen}
\label{subsubsec: Weitere Untersuchungen zu den Legendre-Polynomen}

\epigraph{Things are recalibrated according to new perspectives and perceptions. It's fascinating to me.}{Tim Gunn}

Hier werden die Euler'schen Erkenntnisse bezüglich der Legendre--Polynome ergänzt, indem von (\ref{eq: New Legendre Integral}) aus Abschnitt  (\ref{subsubsec: Den Kontext betreffend: Die Legendre Polynome}) ausgehend, weitere Formeln und Einsichten abgeleitet werden.

\paragraph{Reduktion auf eine bekannte Darstellung}
\label{para: Untersuchung dieser neuen Formel}

Da Euler, an anderen Dingen interessiert, (\ref{eq: New Legendre Integral}) nicht nieder geschrieben hat, soll hier zunächst eine kleine Untersuchung der Formel nachgereicht werden. Als erstes kann sie auf die Laplace'sche Darstellung (\ref{eq: Laplace Legendre}) reduziert werden, was man wie folgt einsieht: Man hat mit der Festlegung

\begin{equation*}
    t+\sqrt{t^2-1}\cos \varphi = y  
\end{equation*}
in dieser Darstellung für die Differentiale
\begin{equation*}
  \sqrt{t^2-1}\sin \varphi d \varphi = dy \quad \text{oder auch} \quad   d\varphi = \dfrac{dy}{-\sqrt{t^2-1}\sin \varphi}.
\end{equation*}
Weil vermöge der Substitution $\cos \varphi = \frac{y-t}{\sqrt{t^2-1}}$ gilt, hat man

\begin{equation*}
    \sin \varphi = \sqrt{1-\dfrac{(y-t)^2}{t^2-1}}.
\end{equation*}
All dies gibt in (\ref{eq: Laplace Legendre}) eingesetzt

\begin{equation*}
    P_{n}(t) =- \dfrac{1}{\pi} \int\limits_{t+\sqrt{t^2-1}}^{t-\sqrt{t^2-1}} \dfrac{y^ndy}{\sqrt{t^2-1}\cdot \sqrt{1-\frac{y-t}{(t^2-1)}}}.
\end{equation*}
Das ``$-$"{} kann durch Vertauschen der Integrationsgrenzen verarbeitet werden, sodass

\begin{equation*}
    P_n(t)= \dfrac{1}{\pi} \int\limits_{t-\sqrt{t^2-1}}^{t+\sqrt{t^2-1}} \dfrac{y^ndy}{\sqrt{(t^2-1)-(y^2-2ty+t^2)}}= \dfrac{1}{\pi} \int\limits_{t-\sqrt{t^2-1}}^{t+\sqrt{t^2-1}} \dfrac{y^ndy}{\sqrt{-1+2ty-y^2}}.
\end{equation*}
Schreibt man $\sqrt{-1}=i$ und zieht $-1$ aus der Wurzel heraus, gelangt man zu

\begin{equation*}
    P_n(t)= \dfrac{1}{i \pi}  \int\limits_{t-\sqrt{t^2-1}}^{t+\sqrt{t^2-1}} \dfrac{y^ndy}{\sqrt{+1-2ty+y^2}}.
\end{equation*}
Indem man also (\ref{eq: New Legendre Integral}) mit dieser letzten Formel vergleicht, findet man auch hier über einen Umweg  $\log(-1)=i \pi$ basierend auf der Wahl $\sqrt{-1}=i$. Für die andere mögliche Wahl $\sqrt{-1}=-i$ erhielte man $\log(-1)=-i\pi$. Die erste Wahl erscheint indes ``etwas natürlicher", da  $\log(-1)=i\pi$ resultiert, wenn  der Hauptzweig des komplexen Logarithmus gewählt wird. Wie an entsprechender Stelle erwähnt, erfährt die natürliche Wahl durch die Arbeit \textit{``Sur les polynômes de Legendre"} (\cite{St90}, 1890) (``Über die Legendre--Polynome") ihre Bestätigung.

\paragraph{Elementare Eigenschaften aus dieser Formel abgeleitet}
\label{para: Parität}

(\ref{eq: New Legendre Integral}) erlaubt indes eine bequemere Ableitung gewisser bekannter Eigenschaften der Legendre--Polynome. Dies sei am Beispiel der Parität aufgezeigt. Dies ist die Eigenschaft:

\begin{equation}
\label{eq: Parity}
    P_n(-t)=(-1)^n P_n(t).
\end{equation}
Aus (\ref{eq: New Legendre Integral}) ergibt sich unmittelbar:

\begin{equation*}
    i \pi P_n(t) = \int\limits_{t-\sqrt{t^2-1}}^{t+\sqrt{t^2-1}} \dfrac{x^ndx}{\sqrt{1-2xt +x^2}}
\end{equation*}
und das Schreiben von $-t$ statt $t$ gibt

\begin{equation*}
    i \pi P_n(-t) = \int\limits_{-t-\sqrt{t^2-1}}^{-t+\sqrt{t^2-1}} \dfrac{x^ndx}{\sqrt{1+2xt +x^2}}.
\end{equation*}
Weiter erhält man für $x=-y$ und daher $dx=-dy$

\begin{equation*}
     i \pi P_n(-t) = -\int\limits_{-(-t-\sqrt{t^2-1})}^{-(-t+\sqrt{t^2-1})} \dfrac{(-y)^ndy}{\sqrt{1-2ty+y^2}}. 
\end{equation*}
Vereinfacht man nun und vertauscht die Integrationsgrenzen, entspringt
\begin{equation*}
       i \pi P_n(t) = (-1)^n \int\limits_{t-\sqrt{t^2-1}}^{t+\sqrt{t^2-1}} \dfrac{y^ndy}{\sqrt{1-2yt+y^2}} =(-1)^n i \pi P_n(t),
\end{equation*}
was  gerade (\ref{eq: Parity}) ist. \\

Auch die Funktionalgleichung (\ref{eq: Functional Equation Legendre}), welche Euler in den  Arbeiten \cite{E672}, \cite{E673}, \cite{E674} in großem Ausmaß eingenommen hat, lässt sich aus (\ref{eq: New Legendre Integral}) beweisen. In moderner Sprache bedeutet dies die Gleichung

\begin{equation}
\label{eq: Functional Equation Legendre}
    P_n(t)=P_{-(n+1)}(t) \quad \text{für} \quad n \in \mathbb{N}.
\end{equation}
Zu diesem Zwecke präsentiere man (\ref{eq: New Legendre Integral}) wie folgt

\begin{equation*}
    i \pi P_n(t)= \int\limits_{a(t)}^{b(t)} \dfrac{x^ndx}{\sqrt{1-2xt+x^2}}
\end{equation*}
mit $a(t)=t-\sqrt{t^2-1}$ und $b(t)=t+\sqrt{t^2-1}$. Dann ist mit $x=\frac{1}{y}$, und daher $dx=-\frac{dy}{y^2}$:

\begin{equation*}
    i \pi P_n(t) = -\int\limits_{\frac{1}{a(t)}}^{\frac{1}{b(t)}} \dfrac{y^{-n}}{\sqrt{1-2t y^{-1}+y^{-2}}} \dfrac{dy}{y^2},
\end{equation*}
welche nach Vertauschen der Integrationsgrenzen diese Form annimmt

\begin{equation*}
    i \pi  P_n(t) =\int\limits_{\frac{1}{b(t)}}^{\frac{1}{a(t)}} \dfrac{y^{-n-1}dy}{\sqrt{1-2t y^{1}+y^{2}}}.
\end{equation*}
Mit der leicht nachweisbaren Identität

\begin{equation}
\label{eq: a*b=1}
    \dfrac{1}{a(t)}= \dfrac{1}{t-\sqrt{t^2-1}}= \dfrac{t+\sqrt{t^2-1}}{t^2-t^2+1}= t+\sqrt{t^2-1}= b(t),
\end{equation}
lässt sich die letzte Gleichung wie folgt darstellen:

\begin{equation*}
    i \pi P_n(t) = i \pi P_{-(n+1)}(t),
\end{equation*}
welche natürlich mit (\ref{eq: Functional Equation Legendre}) gleichwertig ist. \\

Aber auch weniger bekannte, wenn nicht gar neue Formeln lassen sich aus (\ref{eq: New Legendre Integral}) herleiten, wie der Grenzwert:

\begin{equation*}
    i \pi P_n(-1) =i\pi \cdot(-1)^n  = \lim_{t \rightarrow -1} \int\limits_{t-\sqrt{t^2-1}}^{t+\sqrt{t^2-1}} \dfrac{x^ndx}{\sqrt{1-2xt+x^2}} 
\end{equation*}
\begin{equation*}
    = \dfrac{\sqrt{\pi}\Gamma\left(n+1\right)}{\Gamma\left(n+\frac{3}{2}\right)} \lim_{t \rightarrow -1} {}_2F_1\left(\dfrac{1}{2},n+1,n+\frac{3}{2}, (b(t))^2\right) \cdot \left((b(t))^{n+1}-(a(t))^{n+1}\right).
\end{equation*}

\paragraph{Integral--und Reihenidentitäten}

Von größerem Interesse mögen allerdings gewisse Integral-- und Reihenidentitäten sein, womit  auch die erzeugende Funktion der Legendre'schen Polynome  (\ref{eq: Definition Legendre polynomials}) Anwendung findet, sie lautete:
\begin{equation*}
    \dfrac{1}{\sqrt{1-2xt+x^2}}= \sum_{k=0}^{\infty} P_k(t)x^k.
\end{equation*}
Da  der Ausdruck der auf der linken Seite in  (\ref{eq: New Legendre Integral}) auftritt, kann  der entsprechende Anteil  durch den Integranden ersetzen werden. Lässt man auch hier $a(t)=t-\sqrt{t^2-1}$ und $b(t)=t+\sqrt{t^2-1}$ sein und setzt  (\ref{eq: Definition Legendre polynomials}) ein, gelangt man zu:

\begin{equation*}
    i \pi P_n(t)= \int\limits_{a(t)}^{b(t)} x^n \sum_{k=0}^{\infty} P_k(t)x^k = \sum_{k=0}^{\infty} P_k(t) \int\limits_{a(t)}^{b(t)} x^{n+k}dx 
\end{equation*}
\begin{equation*}
    =  \sum_{k=0}^{\infty} P_k(t) \dfrac{(b(t))^{n+k+1}-(a(t))^{n+k+1}}{n+k+1}.
\end{equation*}
Der Kürze wegen setze man $t=\cos(\varphi)$, sodass

\begin{equation*}
    b(t)= \cos (\varphi) +\sqrt{\cos^2(\varphi)-1}= \cos (\varphi) + i \sin (\varphi) =e^{i\varphi}
\end{equation*}
und gleichermaßen

\begin{equation*}
    a(t)= \cos (\varphi)- i \sin (\varphi)=e^{-i \varphi}.
\end{equation*}
Demnach ist

\begin{equation*}
    i \pi P_n(\cos (\varphi))= \sum_{k=0}^{\infty} P_k(\cos (\varphi)) \dfrac{e^{i(n+k+1)\varphi}-e^{-i(n+k+1)\varphi}}{n+k+1}
\end{equation*}
\begin{equation*}
    =  2i \cdot \sum_{k=0}^{\infty} P_k(\cos (\varphi)) \dfrac{\sin ((n+k+1)\varphi)}{n+k+1}.
\end{equation*}
Schließlich hat man

\begin{equation}
\label{eq: Marvelous}
    \dfrac{\pi}{2}P_n(\cos (\varphi)) =  \sum_{k=0}^{\infty} P_k(\cos (\varphi)) \dfrac{\sin ((n+k+1)\varphi)}{n+k+1},
\end{equation}
welche Formel ein einzelnes Legendre--Polynom als Summe aller anderen ausdrückt. Diesem Ausdruck lassen sich mit ähnlichen Überlegungen viele weitere derselben Gestalt an die Seite stellen.\\

Jedoch soll auch die Herleitung von Integralidentitäten dargestellt werden. Aus (\ref{eq: Definition Legendre polynomials}) leitet sich folgende Formel ab:

\begin{equation*}
    \sum_{k=0}^{\infty}P_k(t)= \dfrac{1}{\sqrt{2-2t}}.
\end{equation*}
Für $|t|<1$, ist diese Formel natürlich richtig, sodass vermöge (\ref{eq: New Legendre Integral}) entspringt:

\begin{equation*}
    \dfrac{1}{\sqrt{2-2t}}  = \sum_{n=0}^{\infty}   \dfrac{1}{i\pi} \cdot  \int\limits_{t-\sqrt{t^2-1}}^{t+\sqrt{t^2-1}}\dfrac{x^n dx}{\sqrt{1-2xt+x^2}}= \dfrac{1}{i \pi}  \int\limits_{t-\sqrt{t^2-1}}^{t+\sqrt{t^2-1}} \dfrac{dx}{(1-x)\sqrt{1-2tx+x^2}},
\end{equation*}
wo  die Reihenfolge der Summation und Integration vertauscht und im letzten Schritt die geometrische Reihe verwendet wurde. Auch hier setze man $t=\cos (\varphi)$. Mit der elementaren Identität  $2-2\cos(2x)=4\sin^2 (x)$ gibt dies:

\begin{equation*}
\label{eq: sine Integral}
    \dfrac{1}{2\sin \frac{\varphi}{2}}= \dfrac{1}{i \pi} \int\limits_{e^{-i\varphi}}^{e^{i\varphi}} \dfrac{dx}{(1-x)\sqrt{1-2\cos(\varphi)x+x^2}},
\end{equation*}
wobei, wie schon zuvor, $t\pm \sqrt{t^2-1}=e^{\pm i \varphi}$, genutzt wurde. Statt einer Diskussion der Spezialfälle dieser Formel sei die folgende mitgeteilt

\begin{equation*}
    \dfrac{1}{2\cos \left(\frac{\varphi}{2}\right)}= \dfrac{1}{i \pi} \int\limits_{e^{-i\varphi}}^{e^{i\varphi}} \dfrac{dx}{(1+x)\sqrt{1-2\cos(\varphi)x+x^2}},
\end{equation*}
welche sich durch die Ersetzung $\varphi \mapsto \varphi+ \pi$ aus der vorherigen ergibt.

\paragraph{Interpolation}
\label{para: Interpolation}

Euler war in seiner Arbeit \cite{E606} freilich nicht am Fall $n \notin \mathbb{N}$ interessiert. Jedoch erlaubt  (\ref{eq: New Legendre Integral}) auch eine Darstellung über die hypergeometrische Reihe (\ref{eq: Hypergeometric Series}), welche natürlich eine unmittelbare Interpolation auf die nicht--ganzzahligen Fälle erlaubt. Faktorisiert man in

\begin{equation*}
    i\pi \cdot P_n(t) =  \int\limits_{t-\sqrt{t^2-1}}^{t+\sqrt{t^2-1}} \dfrac{x^ndx}{\sqrt{1-2xt +x^2}}
\end{equation*}
 den Integranden, gelangt man nach eine etwas längeren, aber gradlinigen Rechnung unter Zuhilfenahme der Integraldarstellung der hypergeometrischen Reihe (siehe auch (\ref{eq: Integral Rep}))

\begin{equation*}
    \dfrac{\Gamma(\gamma -\beta)\Gamma(\beta)}{\Gamma(\gamma)}{} _2F_1(\alpha, \beta, \gamma,z)= \int\limits_{0}^1 t^{\beta -1}(1-x)^{\gamma -\beta -1}(1-z x)^{-\alpha}dx,
\end{equation*}
zu

\begin{equation}
\label{eq: Legendre Hyper}
    i \pi \left(\dfrac{\sqrt{\pi}\Gamma\left(n+1\right)}{\Gamma\left(n+\frac{3}{2}\right)}\right)^{-1} P_n(t)=
\end{equation}
\begin{equation*}
    \left((b(t))^{n+1} {}_2F_1\left(\dfrac{1}{2},n+1,n+\frac{3}{2}, (a(t))^2\right)-(a(t))^{n+1} {}_2F_1\left(\dfrac{1}{2},n+1,n+\frac{3}{2}, (b(t))^2\right)\right).
\end{equation*}
wobei die Formel $\Gamma\left(\frac{1}{2}\right)=\sqrt{\pi}$ genutzt wurde. \\

Für den speziellen Fall $(a(t))^2=(b(t))^2$ vereinfacht sich diese Formel zu

\begin{equation*}
  i \pi \left( \dfrac{\sqrt{\pi}\Gamma\left(n+1\right)}{\Gamma\left(n+\frac{3}{2}\right)}\right)^{-1}  P_n(t)  
\end{equation*}
\begin{equation}
\label{eq: Contracted}
 = {}_2F_1\left(\dfrac{1}{2},n+1,n+\frac{3}{2}, (b(t))^2\right) \cdot \left((b(t))^{n+1}-(a(t))^{n+1}\right).
\end{equation}
Man bemerke, dass trotz $(a(t))^2=(b(t))^2$ nicht notwendig $a(t)=b(t)$ gilt\footnote{Es sei eines Beispiels wegen $t=0$. Dann hat man $b(0)=\sqrt{-1}=i$ sowie $a(0)=-\sqrt{-1}=-i$ und $(b(0))^2=(a(0))^2=-1$.}. Um nun die resultierende hypergeometrische Reihe zu evaluieren, bemühe man den Kummer'schen Satz aus \cite{Ku37} (siehe auch (\ref{eq: Transformation Gauss quadratic}))

\begin{equation}
\label{eq: Kummer}
    {}_2F_1(a,b,1-a-b,-1) = \dfrac{\Gamma(1-a-b)\Gamma\left(1+\frac{1}{2}a\right)}{\Gamma(1+a)\Gamma\left(1+\frac{1}{2}a-b\right)}.
\end{equation}
Setzt man hier also $b=n+1$ und $a=\frac{1}{2}$ in der Formel, liefert (\ref{eq: Legendre Hyper}) für diesen Fall:

\begin{equation*}
    i \pi P_n(0) = \dfrac{\sqrt{\pi}\Gamma\left(n+1\right)}{\Gamma\left(n+\frac{3}{2}\right)} \cdot \dfrac{\Gamma\left(n+\frac{3}{2}\right)\Gamma \left(\frac{n}{2}+\frac{3}{2}\right)}{\Gamma(n+2)\Gamma\left(\frac{n}{2}+1\right)}\cdot \left( (i)^{n+1}-(-i)^{n+1}\right).
\end{equation*}
Mit $i=e^{\frac{i\pi}{2}}$ und $-i=e^{-\frac{i\pi}{2}}$ und unter Verwendung der Formel $2i\sin x =e^{ix}-e^{-ix}$ hat man

\begin{equation*}
    i \pi P_n(0)= \dfrac{\sqrt{\pi}\Gamma\left(n+1\right)}{\Gamma\left(n+\frac{3}{2}\right)} \cdot \dfrac{\Gamma\left(n+\frac{3}{2}\right)\Gamma \left(\frac{n}{2}+\frac{3}{2}\right)}{\Gamma(n+2)\Gamma\left(\frac{n}{2}+1\right)}\cdot 2i \sin \left(\dfrac{n+1}{2}\pi\right).
\end{equation*}
Vereinfachen, Verwendung der Relation $\Gamma(x+1)=x\Gamma(x)$  und erneute Vereinfachung, reduziert dies auf

\begin{equation*}
  i \pi P_n(0) =\dfrac{\pi}{2} \cdot \dfrac{\Gamma\left(\frac{1}{2}+\frac{n}{2}\right)}{\Gamma \left(\frac{n}{2}+1\right)} \cdot 2i \sin \left(\dfrac{n+1}{2}\pi\right).
\end{equation*}
Mit der Reflektionsformel 

\begin{equation*}
    \Gamma(x)\Gamma(1-x) = \dfrac{\pi}{\sin (\pi x)}
\end{equation*}
gibt dies schließlich

\begin{equation*}
    i \pi P_n(0)= i \sqrt{\pi} \dfrac{\Gamma\left(\frac{n}{2}+\frac{1}{2}\right)}{\Gamma\left(\frac{n}{2}+1\right)}\cdot  \dfrac{\pi}{\Gamma\left(\frac{n+1}{2}\right)\Gamma \left(\frac{1-n}{2}\right)}
\end{equation*}
oder

\begin{equation*}
        P_n(0)= \dfrac{\sqrt{\pi}}{\Gamma \left(1+\frac{n}{2}\right)\Gamma \left(\frac{1-n}{2}\right)},
\end{equation*}
welche Formel sich auch aus der Rekursionsvorschrift (\ref{eq: Difference Equation for Legendre}) und den Werten $P_0(0)=1$ und $P_1(0)=0$ mit den Eigenschaften der $\Gamma$--Funktion unmittelbar nachweisen lässt.

\subsection{Mit Eulers Ideen zu Formeln von Ramanujan}
\label{subsec: Mit Eulers Ideen zu Formeln von Ramanujan}

\epigraph{Euler and Ramanujan are mathematicians of the greatest importance in the history of constants (and of course in the history of Mathematics [...])}{Edgar William Middlemast}

Da Euler und Ramanujan bezüglich ihrer Fähigkeiten in der Herleitung neuer Formeln, ihrer Virtuosität zur Manipulation arithmetischer Identitäten und mathematischer Intuition vielerorts bis zum heutigen  Tage miteinander verglichen werden, scheint ein diese Verbindung illustrierender Abschnitt den passenden Abschluss der vorliegenden Ausarbeitung zu bilden. Exemplarisch werden drei Resultate von Ramanujan in Euler'scher Manier behandelt. Zunächst wird ein spezielles bestimmtes Integral abgeleitet (Abschnitt \ref{subsubsec: Ein bestimmtes Integral von Ramanujan}), gefolgt vom Ramanujan'schen Mastertheorem (Abschnitt \ref{subsubsec: Ramanujans Mastertheorem}). Den Abschluss bilden die Formeln zur Kreisquadratur (Abschnitt \ref{subsubsec: Ramanujans Formeln zur Kreisquadratur}). Bei alledem wird auf nichts zurückgegriffen, was  Eulers Kenntnisstand wesentlich übersteigt.

\subsubsection{Ein bestimmtes Integral von Ramanujan}
\label{subsubsec: Ein bestimmtes Integral von Ramanujan}

\epigraph{An equation means nothing to me unless it expresses a thought of God.}{Srinivasa Ramanujan}

In 1913 erhielt Godfrey Harold Hardy  einen heute berühmten Brief aus Inden von Srinivasa Ramanujan. Einen Nachdruck dieses Briefes findet man unter anderem in Kapitel $2$ des Buches \textit{``Ramanujan: Letters and Commentary"} (\cite{Ra95}, 1995). In seinem Schreiben listet Ramanujan mehrere Formeln auf. Unter diesen befindet sich diese

\begin{equation*}
    \int\limits_{0}^{\infty} \dfrac{\left(1+\frac{x^2}{(b+1)^2}\right)}{\left(1+\frac{x^2}{a^2}\right)} \cdot  \dfrac{\left(1+\frac{x^2}{(b+2)^2}\right)}{\left(1+\frac{x^2}{(a+1)^2}\right)}\cdot \cdots dx= \dfrac{1}{2}\sqrt{\pi} \cdot \dfrac{\Gamma\left(a+\frac{1}{2}\right)\Gamma(b+1)\Gamma\left(b-a+\frac{1}{2}\right)}{\Gamma(a)\Gamma \left(b+\frac{1}{2}\right)\Gamma(b-a+1)}
\end{equation*}
mit der Einschränkung  $ 0<a<b+\frac{1}{2}$. $\Gamma(x)$ ist auch hier  die $\Gamma$-Funktion, definiert über

\begin{equation}
   \Gamma(x):= \int\limits_{0}^{\infty} t^{x-1}e^{-t}dt \quad \text{für} \quad \operatorname{Re}(x)>0.
\end{equation}
In seiner Arbeit \textit{``Some Definite Integrals"} (\cite{Ra14}, 1915) hat Ramanujan die obige Formel hingegen wie folgt vorgestellt

\begin{equation}
\label{eq: Ramanujan Main}
    \int\limits_{0}^{\infty} \dfrac{\left(1+\frac{x^2}{b^2}\right)}{\left(1+\frac{x^2}{a^2}\right)} \cdot  \dfrac{\left(1+\frac{x^2}{(b+1)^2}\right)}{\left(1+\frac{x^2}{(a+1)^2}\right)}\cdot \cdots dx= \dfrac{1}{2}\sqrt{\pi}\cdot \dfrac{\Gamma\left(a+\frac{1}{2}\right)\Gamma(b)\Gamma\left(b-a-\frac{1}{2}\right)}{\Gamma(a)\Gamma \left(b-\frac{1}{2}\right)\Gamma(b-a)}.
\end{equation}
  Euler hat Integrale ähnlicher Gestalt in seiner Arbeit \textit{``De summo usu calculi imaginariorum in analysi"} (\cite{E621}, 1788, ges. 1776) (E621:`` Über den gewaltigen Nutzen des der imaginären Größen in der Analysis") betrachtet, wenn man statt den Produkten in Ramanujans Formel sich $\sin$ und $\cos$--Funktionen bei Euler denkt. Eine kurze Darstellung der für unsere Zwecke nötigen Ergebnisse  wird der Herleitung der Ramanujan'schen Formel (\ref{eq: Ramanujan Main}) vorangestellt.

\paragraph{Herleitung von Ramanujans Integral aus Eulers Ideen}
\label{para: Herleitung von Ramanujans Integral aus Eulers Ideen}

Ramanujans Formel (\ref{eq: Ramanujan Main}) scheint kein Spezialfall einer von Euler entwickelten Formel aus \cite{E621} zu sein. Nichtsdestoweniger lässt sie sich durch Kombination einiger Ergebnisse von Euler ableiten, die nachstehend noch einmal genannt werden.

\begin{Thm}[Legendre'sche Verdopplungsformel]
\label{Theorem: Legendre'sche Verdopplungsformel}
 Die $\Gamma$-Funktion erfüllt die Identität:

 \begin{equation*}
     \Gamma(2\alpha)= \dfrac{2^{2\alpha -1}}{\sqrt{\pi}}\cdot \Gamma(\alpha)\cdot \Gamma \left(\alpha+\dfrac{1}{2}\right).
 \end{equation*}
\end{Thm}
Obwohl diese Formel nach Legendre benannt ist, ist in den Arbeiten \cite{Ay21} sowie \cite{Ay23c} erläutert und auch in Abschnitt (\ref{subsubsec: Durch Anwendung einer Methode: Seine Konstante A}) dieser Arbeit diskutiert worden, dass sie sich auf andere Weise ausgedrückt bereits in Eulers Arbeiten findet.

\begin{Thm}[Gauß'sche Summationsformel]
\label{Theorem: Gauß'sche Summationsformel}
    Definiert man
    \begin{equation*}
        _2{}F_1(\alpha,\beta,\gamma;x):=
    \end{equation*}
\begin{equation*}
    1+\dfrac{\alpha \beta}{\gamma}\dfrac{x}{1!}+\dfrac{\alpha (\alpha+1)\beta (\beta +1)}{\gamma (\gamma+1)}\dfrac{x^2}{2!}+\dfrac{\alpha (\alpha+1)(\alpha+2)\beta (\beta +1)(\beta+2)}{\gamma (\gamma+1)(\gamma+2)}\dfrac{x^3}{3!}\cdots,
\end{equation*}
gilt die folgende Formel:

\begin{equation*}
    _2{}F_1(\alpha,\beta,\gamma;1)= \dfrac{\Gamma(\gamma)\Gamma(\gamma-\beta-\alpha)}{\Gamma(\gamma-\alpha)\Gamma(\gamma-\beta)}.
\end{equation*}
\end{Thm}
Wie in \cite{Ay24b} argumentiert und oben (Abschnitt \ref{subsubsec: Die Darstellung betreffend -- Die hypergeometrische Reihe}) auch besprochen, findet sich dieser Lehrsatz bereits in anderer Form in der Euler'schen Arbeit \cite{E663}. Gauß hat sie im Jahr 1813 in der Form, wie sie im Theorem zu sehen ist, in seiner Arbeit  \cite{Ga13} bewiesen.

\begin{Thm}[Eine Produktformel]
\label{Theorem: Eine Produktformel}
    Die folgende Formel ist gültig
    \begin{equation*}
        \left(1+\left(\dfrac{x}{a}\right)^2\right) \left(1+\left(\dfrac{x}{a+1}\right)^2\right)\cdot \cdots = \dfrac{\Gamma^2(a)}{\Gamma(a+ix)\Gamma(a-ix)}.
    \end{equation*}
\end{Thm}

\begin{proof}
Diese Formel teilt Ramanujan in \cite{Ra15} ohne Beweis mit. Um sie aus den Euler'schen Formeln heraus zu beweisen, bedarf es eines Ausdrucks aus § 1 seiner Arbeit \cite{E19}. In moderner Schreibweise, findet sich hier die Formel:

    \begin{equation*}
        \Gamma(1+n)= \dfrac{1^{1-n}\cdot 2^n}{1+n}\cdot \dfrac{2^{1-n}\cdot 3^n}{2+n}\cdot \dfrac{3^{1-n}\cdot 4^n}{3+n}\cdot \cdots.
    \end{equation*}
   Also gilt auch

    \begin{equation*}
         \Gamma(1+in+\alpha)= \dfrac{1^{1-(in+\alpha)}\cdot 2^{in+\alpha}}{1+ni+\alpha}\cdot \dfrac{2^{1-(in+\alpha)}\cdot 3^{in+\alpha}}{2+ni+\alpha}\cdot \dfrac{3^{1-(in+\alpha)}\cdot 4^{in+\alpha}}{3+ni+\alpha}\cdot \cdots.
    \end{equation*}
    Analog hat man

\begin{equation*}
         \Gamma(1-in+\alpha)= \dfrac{1^{1-(-in+\alpha)}\cdot 2^{-in+\alpha}}{1-ni+\alpha}\cdot \dfrac{2^{1-(-in+\alpha)}\cdot 3^{-in+\alpha}}{2-in+\alpha}\cdot \dfrac{3^{1-(-in+\alpha)}\cdot 4^{-in+\alpha}}{3-ni+\alpha}\cdot \cdots.
    \end{equation*}
    Daher gibt das Produkt der beiden vorherigen Ausdrücke

\begin{equation*}
    \dfrac{1^{-2\alpha}\cdot 2^{2\alpha}}{(1+\alpha)^2+n^2}\cdot \dfrac{2^{-2\alpha}\cdot 3^{2\alpha}}{(2+\alpha)^2+n^2}\cdot \dfrac{3^{-2\alpha}\cdot 4^{2\alpha}}{(3+\alpha)^2+n^2}\cdot \cdots
\end{equation*}
  Dies vereinfacht sich zu

   \begin{equation*}
 \Gamma(1+in+\alpha)   \Gamma(1-in+\alpha)=   
   \end{equation*}
   \begin{equation*}
       =\dfrac{1}{(1+\alpha)^2}\cdot \dfrac{1}{(2+\alpha)^2}\cdot \dfrac{1}{(3+\alpha)^2}\cdots \times \dfrac{1}{1+\left(\frac{n}{1+\alpha}\right)^2}\cdot \dfrac{1}{1+\left(\frac{n}{2+\alpha}\right)^2}\cdot \dfrac{1}{1+\left(\frac{n}{3+\alpha}\right)^2}\cdot \cdots
   \end{equation*}
   Der Spezialfall $n=0$ reduziert sich auf

   \begin{equation*}
       \Gamma^2(1+\alpha)=\dfrac{1}{(1+\alpha)^2}\cdot \dfrac{1}{(2+\alpha)^2}\cdot \dfrac{1}{(3+\alpha)^2}\cdots
   \end{equation*}
Demnach

\begin{equation*}
    \Gamma(1+in+\alpha)   \Gamma(1-in+\alpha)= \Gamma^2(1+\alpha)\cdot \dfrac{1}{1+\left(\frac{n}{1+\alpha}\right)^2}\cdot \dfrac{1}{1+\left(\frac{n}{2+\alpha}\right)^2}\cdot \dfrac{1}{1+\left(\frac{n}{3+\alpha}\right)^2}\cdot \cdots
\end{equation*}
Somit, indem man $a$ anstelle von von $1+\alpha$ und $x$ anstatt  $n$ schreibt, folgt das Theorem.
\end{proof}

\begin{Thm}[Eine erste Partialbruchzerlegung]
\label{Theorem: Eine erste Partialbruchzerlegung}
    Das Produkt

    \begin{equation*}
        \dfrac{1}{1+\frac{x^2}{a^2}}\cdot \dfrac{1}{1+\frac{x^2}{(a+1)^2}}\cdot \dfrac{1}{1+\frac{x^2}{(a+2)^2}}\cdot \cdots
    \end{equation*}
    ist der folgenden unendlichen Partialbruchzerlegung gleichwertig

    \begin{equation*}
        A_0 \cdot \dfrac{a}{a^2+x^2}+A_1 \cdot \dfrac{a+1}{(a+1)^2+x^2}+A_2 \cdot \dfrac{a+2}{(a+2)^2+x^2}+\cdots
    \end{equation*}
    mit 

    \begin{equation*}
        A_0=\dfrac{2\Gamma(2a)}{\Gamma^2(a)}, \quad A_1= \dfrac{2\Gamma(2a)}{\Gamma^2(a)}\cdot \left(-\frac{2a}{1!}\right), \quad A_2= \dfrac{2\Gamma(2a)}{\Gamma^2(a)}\cdot \dfrac{2a(2a+1)}{2!} \quad \text{etc.}
    \end{equation*}
\end{Thm}

\begin{proof}
Für die Berechnung der Koeffizienten $A_0$, $A_1$, $A_2$ etc. kann man sich der Methode bedienen, welche Euler in seiner Arbeit \cite{E592} vorgestellt hat und welche in Abschnitt (\ref{subsubsec: Methodus Inveniendi über Methodus Demonstrandi}) zur am Beispiel der Partialbruchzerlegung von $\cot(x)$ Sprache kam. Es wird für das Verständnis hinreichen, stellvertretend für den allgemeinen Fall den ersten Koeffizienten $A_0$ berechnet zu haben. Nach Eulers Vorgehen nehme man an:

    \begin{equation*}
          \dfrac{1}{1+\frac{x^2}{a^2}}\cdot \dfrac{1}{1+\frac{x^2}{(a+1)^2}}\cdot \dfrac{1}{1+\frac{x^2}{(a+2)^2}}\cdot \cdots= A_0 \dfrac{a}{a^2+x^2}+R_0,
    \end{equation*}
    wo $R_0$ eine Funktion ist, die den Faktor $a^2+x^2$ nicht beinhaltet. Gemäß Theorem (\ref{Theorem: Eine Produktformel}) liest sich diese Gleichung als

    \begin{equation*}
        \dfrac{\Gamma(a+ix)\Gamma(a-ix)}{\Gamma^2(a)}=A_0 \cdot \dfrac{a}{a^2+x^2}+R_0,
    \end{equation*}
   welche gleichwertig ist zu

    \begin{equation*}
        1= A_0 \cdot \dfrac{a}{a^2+x^2}\dfrac{\Gamma^2(a)}{\Gamma(a+ix)\Gamma(a-ix)}+R_0 \cdot \dfrac{\Gamma^2(a)}{\Gamma(a+ix)\Gamma(a-ix)}.
    \end{equation*}
 Man nehme den Grenzwert $x \rightarrow ia$. Dies gibt

 \begin{equation*}
     1=A_0 \cdot a \lim_{x\rightarrow ia} \dfrac{\Gamma^2(a)}{\Gamma(a+ix)\Gamma(a-ix)}\dfrac{1}{a^2+x^2},
 \end{equation*}
 weil der zweite Term auf der rechten Seite für diesen Grenzwert verschwindet. Unter Bemerken der Gleichheit $a^2+x^2=(a+ix)(a-ix)$ und Verwendung der Relation $\Gamma(z+1)=z\Gamma(z)$ erhält man:

 \begin{equation*}
     1=A_0 \cdot a \lim_{x \rightarrow ia} \dfrac{\Gamma^2(a)}{\Gamma(a+ix+1)\Gamma(a-ix+1)}.
 \end{equation*}
 Der Grenzwert kann evaluiert werden und führt zur Gleichung

 \begin{equation*}
     1= A_0 \cdot a \cdot \dfrac{\Gamma^2(a)}{\Gamma(1)\Gamma(2a+1)}=  A_0 \cdot a \cdot \dfrac{\Gamma^2(a)}{2a\Gamma(2a)},
 \end{equation*}
 wo $\Gamma(z+1)=z\Gamma(z)$ und $\Gamma(1)=1$ im letzten Schritt gebraucht wurden. Durch Auflösen nach $A_0$ gelangt man zu

 \begin{equation*}
     A_0 = \dfrac{2\Gamma(2a)}{\Gamma^2(a)}.
 \end{equation*}
 Die Berechnung der weiteren Koeffizienten $A_1$, $A_2$ verläuft analog. Das Muster der Werte wird schnell sichtbar, sodass eine detaillierte Untersuchung hier ausgespart wird. \end{proof}

 \begin{Thm}[Eine zweite Partialbruchzerlegung]
 \label{Theorem: Eine zweite Partialbruchzerlegung}
     Das Produkt

     \begin{equation*}
         \dfrac{\Gamma^2(b)\Gamma(a+ix)\Gamma(a-ix)}{\Gamma^2(a)\Gamma(b+ix)\Gamma(b-ix)}= \dfrac{1+\frac{x^2}{b^2}}{1+\frac{x^2}{a^2}}\cdot \dfrac{1+\frac{x^2}{(b+1)^2}}{1+\frac{x^2}{(a+1)^2}}\cdot \dfrac{1+\frac{x^2}{(b+1)^2}}{1+\frac{x^2}{(a+1)^2}} \cdot \cdots
     \end{equation*}
     kann als die folgende Partialbruchzerlegung dargestellt werden

     \begin{equation*}
          \dfrac{\Gamma^2(b)\Gamma(a+ix)\Gamma(a-ix)}{\Gamma^2(a)\Gamma(b+ix)\Gamma(b-ix)}= B_0 \cdot \dfrac{a}{a^2+x^2}+B_1 \cdot \dfrac{a+1}{(a+1)^2+x^2}+B_2 \cdot \dfrac{a+2}{(a+2)^2+x^2}\cdots
     \end{equation*}
     mit Koeffizienten

     \begin{equation*}
         B_0= A_0 \cdot \dfrac{\Gamma^2(b)}{\Gamma(b-a)\Gamma(b+a)}, \quad B_1= \dfrac{\Gamma^2(b)}{\Gamma(b-a)\Gamma(b+a)}\cdot A_1 \cdot \dfrac{b-a-1}{b+a},
     \end{equation*}
     \begin{equation*}
         B_2= \dfrac{\Gamma^2(b)}{\Gamma(b-a)\Gamma(b+a)}\cdot A_2 \cdot \dfrac{b-a-1}{b+a} \cdot \dfrac{b-a-2}{b+a+1} \quad \text{etc.}
     \end{equation*}
     wo $A_0$, $A_1$, $A_2$ etc. die Koeffizienten von Theorem (\ref{Theorem: Eine erste Partialbruchzerlegung}) andeutet, also

     \begin{equation*}
          A_0=\dfrac{2\Gamma(2a)}{\Gamma^2(a)}, \quad A_1= \dfrac{2\Gamma(2a)}{\Gamma^2(a)}\cdot \left(-\frac{2a}{1!}\right), \quad A_2= \dfrac{2\Gamma(2a)}{\Gamma^2(a)}\cdot \dfrac{2a(2a+1)}{2!} \quad \text{etc.}
     \end{equation*}
 \end{Thm}

 \begin{proof}
     Der Beweis verläuft analog zu dem von Theorem (\ref{Theorem: Eine erste Partialbruchzerlegung}). Zwecks Illustration wird der Koeffizient $B_0$ hier explizit ausgerechnet. Wieder dem Euler'schen Vorgehen aus  \cite{E592} folgend, betrachte man den Ansatz

     \begin{equation*}
         \dfrac{\Gamma^2(b)\Gamma(a+ix)\Gamma(a-ix)}{\Gamma^2(a)\Gamma(b+ix)\Gamma(b-ix)}= B_0 \cdot \dfrac{a}{a^2+x^2}+R_0,
     \end{equation*}
     wo $R_0$ eine Funktion ist, die den Faktor $a^2+x^2$ nicht aufweist. Daraus inferiert man den Grenzwert

     \begin{equation*}
         \dfrac{\Gamma^2(b)}{\Gamma(b-a)\Gamma(b+a)}=B_0 \cdot a \cdot \lim_{x\rightarrow ia} \dfrac{1}{x^2+a^2}\dfrac{1}{\Gamma(a+ix)\Gamma(a-ix)}.
     \end{equation*}
     Der noch vorhandene Grenzwert ist aus dem vorherigen Theorem bekannt. Unter Verwendung des Koeffizienten $A_0$ aus nämlichen Theorem (\ref{Theorem: Eine erste Partialbruchzerlegung}), gibt dies insgesamt die Gleichung

     \begin{equation*}
         B_0=A_0 \cdot \dfrac{\Gamma^2(b)}{\Gamma(b-a)\Gamma(b+a)}.
     \end{equation*}
  Gleichermaßen ist

   \begin{equation*}
       B_1= \dfrac{\Gamma^2(b)}{(b+a)\Gamma(b+a)}\cdot \dfrac{b-a-1}{\Gamma(b-a)}\cdot A_1=- \dfrac{\Gamma^2(b)}{\Gamma(b+a)\Gamma(b-a)} \cdot A_1 \cdot \dfrac{a-b+ 1}{b+a}
   \end{equation*}
   und

   \begin{equation*}
       B_2= \dfrac{\Gamma^2(b)}{\Gamma(b-a)\Gamma(b+a)}\cdot A_2 \cdot \dfrac{a-b+1}{b+a} \cdot\dfrac{a-b+2}{b+a+1}.
   \end{equation*}
  Analoges gilt für die weiteren Koeffizienten.
 \end{proof}

 \paragraph{Herleitung des Ramanujan'schen Integrals}
 \label{para: Herleitung des Ramanujan'schen Integrals}

 Nach diesen Vorbereitungen kann nun (\ref{eq: Ramanujan Main}) hergeleitet werden. Aus Theorem (\ref{Theorem: Eine zweite Partialbruchzerlegung}) lässt sich das Integral wie folgt umschreiben 

 \begin{equation*}
     \int\limits_{0}^{\infty} \left(B_0 \cdot \dfrac{a}{a^2+x^2}+B_1 \cdot \dfrac{a+1}{(a+1)^2+x^2}+B_2 \cdot \dfrac{a+2}{(a+2)^2+x^2}\cdots\right)dx.
 \end{equation*}
Im Allgemeinen gilt

 \begin{equation*}
     \int\limits_{0}^{\infty} \dfrac{kdx}{x^2+k^2}= \dfrac{\pi}{2}
 \end{equation*}
 für eine reelle Zahl $k$. Demnach erhält man durch Einsetzen der Koeffizienten $B_0$, $B_1$, $B_2$, $\cdots$ aus Lehrsatz (\ref{Theorem: Eine zweite Partialbruchzerlegung}) und separates Integrieren jedes Terms:

 \begin{equation*}
     \int\limits_{0}^{\infty}  \dfrac{\Gamma^2(b)\Gamma(a+ix)\Gamma(a-ix)}{\Gamma^2(a)\Gamma(b+ix)\Gamma(b-ix)}dx= \dfrac{\pi}{2}\cdot\dfrac{\Gamma^2(b)}{\Gamma(b-a)\Gamma(b+a)} \cdot \dfrac{2 \Gamma^2(2a)}{\Gamma^2(a)}\times
 \end{equation*}
 \begin{equation}
 \label{eq: Main Intermediate}
     \left(1+ \dfrac{2a}{1}\cdot \dfrac{a-b+1}{a+b}+\dfrac{2a(2a+1)}{2!}\cdot \dfrac{a-b+1}{a+b}\cdot \dfrac{a-b+2}{a+b+1}+\cdots\right).
 \end{equation}
 Die Summe, welche den zweiten Faktor bildet, kann vermöge Theorem  (\ref{Theorem: Gauß'sche Summationsformel}) gefunden werden. Die Gauß'sche Formel liefert:

 \begin{equation*}
     _2F_1(2a,a-b+1,a+b;1)= \dfrac{\Gamma(a+b)\Gamma(2(b-a)-1)}{\Gamma(b-a)\Gamma(2b-1)}.
 \end{equation*}
 Unter Verwendung von $\Gamma(x+1)=x\Gamma(x)$ und anschließend der Legendre'schen Verdopplungsformel (Theorem \ref{Theorem: Legendre'sche Verdopplungsformel}) ergibt sich die rechte Seite zu

 \begin{equation*}
     \dfrac{\Gamma(a+b)}{\Gamma(b-a)} \cdot \dfrac{2^{2(b-a)-1}}{2^{2b-1}} \cdot \dfrac{\Gamma(b-a)\Gamma\left(b-a+\frac{1}{2}\right)}{\Gamma(b)\Gamma \left(b+\frac{1}{2}\right)} \cdot \dfrac{2\left(b-\frac{1}{2}\right)}{2\left(b-1-\frac{1}{2}\right)},
 \end{equation*}
 sodass sich viele Terme aufheben. Schlussendlich gelangt man zum Ausdruck:

 \begin{equation*}
     \dfrac{\Gamma(a+b)}{\Gamma(b)}\cdot 2^{-2a} \cdot \dfrac{\Gamma \left(b-a-\frac{1}{2}\right)}{\Gamma \left(b-\frac{1}{2}\right)}.
 \end{equation*}
Setzt man dies in  (\ref{eq: Main Intermediate}) ein, ergibt sich

 \begin{equation*}
     \int\limits_{0}^{\infty}  \dfrac{\Gamma^2(b)\Gamma(a+ix)\Gamma(a-ix)}{\Gamma^2(a)\Gamma(b+ix)\Gamma(b-ix)}dx= 
 \end{equation*}
 \begin{equation*}
     \dfrac{\pi}{2}\cdot\dfrac{\Gamma^2(b)}{\Gamma(b-a)\Gamma(b+a)}\dfrac{2 \Gamma^2(2a)}{\Gamma^2(a)} \times \dfrac{\Gamma(a+b)}{\Gamma(b)}\cdot 2^{-2a} \cdot \dfrac{\Gamma \left(b-a-\frac{1}{2}\right)}{\Gamma \left(b-\frac{1}{2}\right)}.
 \end{equation*}
Mit Theorem (\ref{Theorem: Legendre'sche Verdopplungsformel}) ersetze man $\Gamma(2a)$ und vereinfache anschließend; die rechte Seite wird dann zu

\begin{equation*}
    \dfrac{1}{2}\sqrt{\pi}\cdot \dfrac{\Gamma \left(a+\frac{1}{2}\right)\Gamma(b)\Gamma\left(b-a-\frac{1}{2}\right)}{\Gamma(a)\Gamma\left(b-\frac{1}{2}\right)\Gamma(b-a)}.
\end{equation*}
Vermöge der Produktformel (\ref{Theorem: Eine Produktformel}) kann der Integrand als Produkt ausgedrückt werden, womit  (\ref{eq: Ramanujan Main}) bewiesen ist -- allein unter Verwendung Euler'scher Formeln und Ideen.

\paragraph{Gründe warum Euler Ramanujans Formel nicht selbst entdeckt hat}
\label{para: Gründe warum Euler Ramanujans Formel nicht selbst entdeckt hat}


Obwohl Euler alle Bausteine und Fähigkeiten besaß, um (\ref{eq: Ramanujan Main}) selbst zu entdecken, findet sich die Ramanujan'sche Formel nicht unter den seinen. Euler scheint überdies die Produkte in der Form von Theorem (\ref{Theorem: Eine Produktformel}) nur für den Fall betrachtet  $a=1$ zu haben\footnote{In diesem Fall wird das Produkt $\frac{\sinh(\pi x)}{\pi x}$.}. Unabhängig davon hat Euler in seinen Untersuchungen über die $\Gamma$--Funktion seine Forschungen nie \textit{explizit} auf komplexe Zahlen ausgedehnt. Vielmehr präferierte Euler generell die reelle Gerade in seinen Ausführungen, wie oben  (Abschnitt \ref{subsubsec: Komplexe Analysis}) dargelegt wurde. Daher lässt sich die Nichtentdeckung von (\ref{eq: Ramanujan Main}) durch Euler wohl zum großen Teil auf Eulers Art der Behandlung von Funktionen einer komplexen Variable zurückführen.

\subsubsection{Ramanujans Mastertheorem}
\label{subsubsec: Ramanujans Mastertheorem}

\epigraph{If you have to prove a theorem, do not rush. First of all, understand fully what the theorem says, try to see clearly what it means. [...] When you have satisfied yourself that the theorem is true, you can start proving it.}{George Polya}

Ramanujan hat in seiner Arbeit \cite{Ra15} weiterhin die allgemeine Formel  

\begin{equation}
´\label{eq: Ramanujans Master}
    \int\limits_{0}^{\infty} x^{s-1}f(x)dx= \varphi(-s)\Gamma(s) \quad \text{und} \quad f(x)= \sum_{k=0}^{\infty} \dfrac{\varphi(k)}{k!}(-x)^k
\end{equation}
angegeben. Während Ramanujan in erwähnter Abhandlung keinen Beweis mitteilt, soll hier ein Weg aufgezeigt werden, wie ein solcher sich aus Eulers Formeln finden lässt. Dafür bedarf es zweier Ausdrücke aus Eulers Arbeiten. 

\begin{Thm}[Euler'sche Transformationsformel]
    \label{Theorem: Euler'sche Transformationsformel}
    Die Potenzreihe
    \begin{equation*}
        B(x):=b_1 x-b_2x^2+b_3x^3-b_4x^4+\cdots 
    \end{equation*}
    kann auch wie folgt geschrieben werden:

\begin{equation*}
    B(x)= b_1 \dfrac{x}{1+x}- \Delta^1 b_1 \left(\dfrac{x}{1+x}\right)^2+\Delta^2 b_1 \left(\dfrac{x}{1+x}\right)^3- \Delta^3 b_1 \left(\dfrac{x}{1+x}\right)^4-\cdots
\end{equation*}
Dabei sind $\Delta^1 b_1$, $\Delta^2 b_1$, $\Delta^3 b_1$ etc. die Differenzen erster, zweiter, dritter etc. Ordnung.
\end{Thm}
Euler weist diese Identität unter anderem in § 8 des zweiten Teils seiner \textit{Calculi Differentialis} \cite{E212} nach\footnote{Diese Formel ließe sich auch in der Untersuchung über von Euler mitgeteilte Kuriositäten (Abschnitt \ref{subsubsec: Eulers Auffassung eines Beweises}) herauslesen.}. Weiter benötigt man den folgenden Lehrsatz, welchen Euler in § 4 seiner Arbeit \cite{E613} beweist. 

\begin{Thm}[Euler-Newton-Interpolation]
\label{Theorem: Euler-Newton-Interpolation}
    Ist eine Funktion $f$ auf den natürlichen Zahlen gegeben, kann sie wie folgt geschrieben werden:

\begin{small}
    \begin{equation*}
        f(n)= f(0)+\dfrac{n-1}{1}\Delta^1 f(0) +\dfrac{n-1}{1}\dfrac{n-2}{2}\Delta^2 f(0) +\dfrac{n-1}{1}\dfrac{n-2}{2}\dfrac{n-3}{3} \Delta^3 f(0)+\text{etc.}
    \end{equation*}
   \end{small} 
    mit den Differenzenoperatoren $\Delta^1$, $\Delta^2$, $\Delta^3$ etc. erster, zweiter, dritter etc. Ordnung, sodass die rechte Seite auch für nicht--natürliche $n$ sinnvoll ist.
\end{Thm}
Diese Interpolationsformel ist auch als Newton--Interpolation bekannt\footnote{Euler hat sie auch  in § 17 seiner viel früher verfassten Arbeit \cite{E125} aus dem Jahr 1739 angegeben und beschreibt sie dort als ``bekannt". Damit könnte Euler auf Newtons Werk referieren, wo man sie im dritten Buch seiner \textit{``Philosophiae naturalis principia mathematica"} (\cite{Ne87}, 1687) findet, mit welchem Euler  gut vertraut war.}. Mit ihr und dem ersten Satz kann nun Ramanujans Formel (\ref{eq: Ramanujans Master}) bewiesen werden. Wie in Theorem (\ref{Theorem: Euler'sche Transformationsformel}) setze man 

\begin{equation*}
    B(x) = b_1 x-b_2x^2+b_3x^3-b_4x^4+\cdots = \sum_{n=1}^{\infty} (-1)^{n+1} b_n x^n
\end{equation*}
und betrachte das Integral

\begin{equation*}
    F(s) := \int\limits_{0}^{\infty} x^{s-1}B(x)dx.
\end{equation*}
Nach der Euler'schen Transformationsformel (\ref{Theorem: Euler'sche Transformationsformel}) gilt:

\begin{equation*}
    F(s) = \int\limits_{0}^{\infty} x^{s-1}\sum_{n=0}^{\infty} \Delta^n b_1 \left(\dfrac{x}{1+x}\right)^{n+1}dx.
\end{equation*}
Unter Annahme entsprechender Annahmen über die Koeffizienten, lassen sich die Integration und Summation vertauschen, sodass gilt:

\begin{equation*}
    F(s) = \sum_{n=0}^{\infty} (-1)^n \Delta^n b_1 \int\limits_{0}^{\infty} \dfrac{x^{s+n}dx}{(1+x)^{n+1}}.
\end{equation*}
Mit den Darstellungen des Beta--Integrals

\begin{equation*}
    B(x,y)= \int\limits_{0}^1 \dfrac{t^{x-1}dt}{(1+t)^{x+y}}= \frac{\Gamma(x)\Gamma(y)}{\Gamma(x+y)}
\end{equation*}
vereinfacht sich die letzte Formel zu:

\begin{equation*}
    F(s)= \sum_{n=0}^{\infty} (-1)^n \Delta^n b_1  \dfrac{\Gamma(-s)\Gamma(s+n+1)}{\Gamma(n+1)}.
\end{equation*}
Unter wiederholter Anwendung der Funktionalgleichung $\Gamma(x+1)=x\Gamma(x)$ ergibt sich

\begin{equation*}
    F(s)= \sum_{n=0}^{\infty} (-1)^n \Delta^n b_1  \Gamma(s)\Gamma(1-s) \dfrac{s+1}{1}\dfrac{s+2}{2}\cdots \dfrac{s+n}{n},
\end{equation*}
welche Formel sich unter Berücksichtigung der Formel

\begin{equation*}
    \Gamma(s)\Gamma(1-s)=\dfrac{\pi}{\sin(\pi s)}
\end{equation*}
als

\begin{equation*}
    F(s)= \dfrac{\pi}{\sin(-s\pi)} \sum_{n=0}^{\infty} (-1)^n \Delta^n b_1  \dfrac{s+1}{1}\dfrac{s+2}{2}\cdots \dfrac{s+n}{n}
\end{equation*}
schreiben lässt. Hier lasse man nun $s$ in $-s$ übergehen: 

\begin{equation*}
    F(-s)= \dfrac{\pi}{\sin (s\pi)} \sum_{n=0}^{\infty}  (-1)^n \Delta^n b_1  \dfrac{-s+1}{1}\dfrac{-s+2}{2}\cdots \dfrac{-s+n}{n}
\end{equation*}
oder

\begin{equation*}
    F(-s) = \dfrac{\pi}{\sin (s\pi)} \sum_{n=0}^{\infty} \Delta^n b_1 \dfrac{s-1}{1}\dfrac{s-2}{2}\cdots \dfrac{s+n}{n}.
\end{equation*}
Die Summe kann nun mit der Newton'schen--Interpolationformel aus Theorem (\ref{Theorem: Euler-Newton-Interpolation}) ausgewertet werden, sodass insgesamt gilt:

\begin{equation*}
    F(-s) =\dfrac{\pi}{\sin(\pi s)}b(s)
\end{equation*}
oder auch

\begin{equation}
\label{eq: Euler Master}
    F(s)= \int\limits_{0}^{\infty} x^{s-1}B(x)dx= -\dfrac{\pi}{\sin(\pi s)}b(-s).
\end{equation}
Weil 

\begin{equation*}
    -B(x)= \sum_{n=1}^{\infty} (-1)^{n} b_n x^n,
\end{equation*}
interpoliert die soeben bewiesene Formel die Koeffizienten $b_n$. Um nun zur Ramanujan'schen Formel (\ref{eq: Ramanujans Master}) zu gelangen, ist in den Formeln $b_n$ in $\frac{\varphi(n)}{n!}$ übergehen zu lassen, woraus die Formel unmittelbar folgt.\\

Die Ramanujan'sche Interpolationsformel (\ref{eq: Ramanujans Master}) kann also ebenfalls aus Euler'schen Formeln hergeleitet werden. Die Verbindung zur Euler'schen Transformationsformel aus Theorem (\ref{Theorem: Euler-Newton-Interpolation}) ist auch von Hardy auf Seite 192 seines Buches \textit{``Ramanujan"} (\cite{Ha02}, 2002)  bemerkt worden. In selbigem Buch findet sich auch ein Beweis von Ramanujans Theorem und ähnlichen derselben Gestalt mit den Mitteln komplexer Analysis.

\subsubsection{Ramanujans Formeln zur Kreisquadratur}
\label{subsubsec: Ramanujans Formeln zur Kreisquadratur}

\epigraph{The mathematician's patterns, like the painter's or the poet's must be beautiful; the ideas, like the colours or the words must fit together in a harmonious way. Beauty is the first test: there is no permanent place in the world for ugly mathematics.}{Godfrey Harold Hardy}

Euler hat seine Definition  der Summe von divergenten Reihen (siehe Abschnitt (\ref{subsubsec: Der Begriff der Summe einer Reihe}) für die Definition) an vielen anderen Stellen verwendet, um gewisse Hindernisse zu überwinden -- das angesprochene Beispiel der $\zeta$--Funktion (Abschnitt \ref{subsubsec: Durch Kombinieren von Ergebnissen: Die zeta-Funktion}) kann als pars pro toto dessen angesehen werden.  Statt einer Besprechung weiterer Beispiele von Euler sich kann sich auf modernere Arbeiten wie  \textit{``Proof of some conjectured formulas for $\frac{1}{\pi}$ by Z.W.Sun"} (\cite{Al11}, 2011)  und \textit{``On proving some of Ramanujan's formulas for $\frac{1}{\pi}$ with an elementary method"} (\cite{Ay13}, 2013) berufen werden, welche im Zusammenhang von Reihen der Form

\begin{equation}
\label{eq: 1/pi}
    \sum_{n=0}^{\infty}\dfrac{\left(\frac{1}{2}\right)_n \left(s\right)_n\left(1-s\right)_n}{(1)_n^3}(An+B)(x_0)^n =\dfrac{C}{\pi},
\end{equation}
wo $s=\frac{1}{2},\frac{1}{3}$ oder $\frac{1}{6}$ und $A,B$ rationale sowie $C$ und $x_0$ komplexe Zahlen sind und 

\begin{equation*}
    (x)_n:= x(x+1)\cdots (x+n)
\end{equation*}
das Pochhammer--Symbol bedeutet,  ebenfalls in Euler'scher Manier mit divergenten Reihen hantieren, um gewisse Resultate nachzuweisen.  Die Betrachtung und Behandlung dieser Reihen stellt eine schöne Parallele zur Euler'schen Lösung des Basel-Problems (Abschnitt \ref{subsec: Die Lösung des Baseler Problems}) dar. Von den ursprünglichen Ramanujan'schen Arbeiten, wie \textit{``Modular equations and approximations to $\pi$"} (\cite{Ra14}, 1914), motiviert sind in der Arbeit \textit{``Some new series for $\frac{1}{\pi}$ and related congruences"}  (\cite{Su11}, 2011) viele Reihen der Form (\ref{eq: 1/pi}) mit Hilfe computerunterstützter Methoden gefunden worden, verlangten aber nach einem expliziten Beweis, was für einen Teil der dort mitgeteilten Identitäten  in \cite{Al11}  mit elementaren Methoden geschehen ist, welche in ihrer Essenz nicht über die aus den Euler'schen Arbeiten und \cite{Ga13} hinausgehen. Die grundlegende Identität ist die Gauß'sche Formel (\ref{eq: Gauss zweite Summationsformel}), welche dann mit den zahlreichen Transformationsformeln für die hypergeometrische Funktion, welche Euler und Gauß noch nicht bekannt waren, und etwa im Buch \textit{``Special Functions"} (\cite{An11}, 2011) zu finden sind, in die in \cite{Su11} entdeckten Reihen überführt worden sind. Diese Idee wird in der Übersichtsarbeit \cite{Ay13}  dargestellt. Hier sollen die oben genannten Arbeiten am Beweis der sogenannten Bauer'schen Reihe nachvollzogen werden. Dies ist die Summierung

\begin{equation}
    \label{eq: Bauer}
    \dfrac{2}{\pi}=  1- 5 \left(\dfrac{1}{2}\right)^3+9\left(\dfrac{1 \cdot 3}{2 \cdot 4}\right)^3- 13\left(\dfrac{1 \cdot 3 \cdot 5}{2 \cdot 4 \cdot 6}\right)^3+\cdots
\end{equation}
Diese Reihe wurde ebenfalls von Ramanujan in seinem oben (Abschnitt \ref{subsubsec: Ein bestimmtes Integral von Ramanujan}) erwähnten Brief an Hardy  mitgeteilt. Daneben wurde ihr unlängst die Arbeit \textit{``The Bauer–Ramanujan formula: historical analyses and perspectives"} (\cite{Ca24}, 2024) gewidmet, welche verschiedene Beweise dieser Formel gesammelt darstellt. 

\paragraph{Die grundlegende Formel}
\label{para: Grundformel}

Für den Beweis der Bauer'schen Reihe bedarf es einer Hilfsformel aus \cite{Al11}.

\begin{Thm}[Hilfsformel]
    \label{Theorem: Hilfsformel 1/pi}
    Es gilt
    \begin{equation*}
         \sum_{n=0}^{\infty} n \left(\dfrac{1}{2}\right)^n \sum_{k=0}^n \left(\dfrac{(s)_k(1-s)_k}{(k!)}\right)^2  \left(\dfrac{(s)_{n-k}(1-s)_{n-k}}{(n-k)!}\right)^2
    \end{equation*}
    mit dem Pochhammer--Symbol 

    \begin{equation*}
        (s)_n = \dfrac{\Gamma(s+n)}{\Gamma(s)}= (s+n-1)(s+n-2)\cdots s.
    \end{equation*}
\end{Thm}

\begin{proof}
    Man betrachte das Quadrat der hypergeometrischen Funktion $_2F_1(s,1-s,1;x)$, denn mit dem Cauchy--Produkt für Potenzreihen gilt:

    \begin{equation*}
        ({}_2F_1(s,1-s,1;x))^2= \sum_{n=0}^{\infty} x^n\sum_{k=0}^n \left(\dfrac{(s)_k(1-s)_k}{(k!)}\right)^2  \left(\dfrac{(s)_{n-k}(1-s)_{n-k}}{(n-k)!}\right)^2.
    \end{equation*}
    Nun werde der Operator $x\frac{d}{dx}$ auf beide Seiten appliziert:

    \begin{equation*}
        \sum_{n=0}^{\infty} nx^n \sum_{k=0}^n \left(\dfrac{(s)_k(1-s)_k}{(k!)}\right)^2  \left(\dfrac{(s)_{n-k}(1-s)_{n-k}}{(n-k)!}\right)^2= 
    \end{equation*}
    \begin{equation*}
        2s(1-s)x{}_2F_1(s,1-s,1;x){}_2F_1(1+s,2-s,2;x).
    \end{equation*}
    Die linke Seite hat bereits die gewünschte Form. Zur Vereinfachung der rechten Seite greife man auf die zweite Gauß'sche Summationsformel zurück,  welche auch oben in (\ref{eq: Gauss zweite Summationsformel}) schon auftrat. Sie lautet 

    \begin{equation*}
        _2F_1 \left(a,b,\dfrac{a+b+1}{2};\dfrac{1}{2}\right)= \sqrt{\pi} \dfrac{\Gamma \left(\frac{a+b+1}{2}\right)}{\Gamma \left(\frac{a+1}{2}\right)\Gamma \left(\frac{b+1}{2}\right)},
    \end{equation*}
    mit welcher die rechte Seite übergeht in:

    \begin{equation*}
        2 \cdot \dfrac{1}{2} s(1-s) \pi \cdot \dfrac{\Gamma(1)}{\Gamma\left(\frac{1}{2}+\frac{s}{2}\right)\Gamma \left(1-\frac{s}{2}\right)} \cdot \frac{\Gamma(2)}{\Gamma \left(1+\frac{s}{2}\right)\Gamma \left(\frac{3}{2}-s\right)}.
    \end{equation*}
    Mit der Formel $\Gamma(x+1)=x\Gamma(x)$ folgt dann:

        \begin{equation*}
        2 \cdot \dfrac{1}{2} s(1-s) \pi \cdot \dfrac{\Gamma(1)}{\Gamma\left(\frac{s}{2}\right)\frac{s}{2}\Gamma \left(1-\frac{s}{2}\right)} \cdot \frac{\Gamma(2)}{\Gamma \left(1+\frac{s}{2}\right)\left(\frac{1}{2}-\frac{s}{2}\right)\Gamma \left(\frac{1}{2}-\frac{s}{2}\right)}.
    \end{equation*}
    Dies lässt sich vermöge der Reflexionsformel $\frac{\pi}{\sin (\pi s)}=\Gamma(s)\Gamma(1-s)$ vereinfachen zu:

    \begin{equation*}
        4\pi \cdot \dfrac{\sin \left(\frac{\pi s}{2}\right)}{\pi} \cdot \dfrac{\sin \left( \pi \left(\frac{1}{2}-\frac{s}{2}\right)\right)}{\pi}= \dfrac{4}{\pi} \sin \left(\dfrac{\pi s}{2}\right) \cos \left(\dfrac{\pi s}{2}\right)= \dfrac{2}{\pi}\sin (\pi s),
    \end{equation*}
    wobei  im letzten Schritt die Verdopplungsformel für den Sinus benutzt wurde.
\end{proof}

\paragraph{Das Beispiel der Bauer'schen Reihe}
\label{para: Beispiel der Bauer'schen Reihe}

Nun erfolgt der Beweis von (\ref{eq: Bauer}), welche mithilfe des Pochhammer--Symbols auch wie folgt angeschrieben werden kann: 

\begin{equation*}
    \sum_{n=0}^{\infty} \dfrac{\left(\frac{1}{2}\right)_n \left(\frac{1}{2}\right)_n \left(\frac{1}{2}\right)_n}{(1)_n)(1)_n)(1)_n)}(4n+1)(-1)^n = \dfrac{2}{\pi}.
\end{equation*}

Denn die Verwendung einer weiteren Transformationsformel für hypergeometrische  Reihe, die Johann Pfaff (1765--1825) zugeschrieben wird, sich aber wie erwähnt (Abschnitt \ref{subsubsec: Die Darstellung betreffend -- Die hypergeometrische Reihe}) sich aus der Euler'schen Formel aus Theorem (\ref{Theorem: Euler'sche Transformationsformel}) ergibt, nämlich:

\begin{equation}
    \label{eq: Pfaff-Trafo}
    _2F_1(\alpha, \beta, \gamma; x) =(1-x)^{-\beta} {}_2F_1 \left(\gamma -\alpha, \beta, \gamma; \dfrac{x}{x-1}\right),
\end{equation}
für den Fall $\alpha=\frac{3}{4}$, $\beta=\frac{1}{4}$, $\gamma=1$ und anschließender Quadratur ergibt:

\begin{equation*}
    \left(_2F_1 \left(\dfrac{1}{4}, \dfrac{3}{4}, 1;x\right)\right)^2 =(1-x)^{-\frac{1}{2}}\left(_2F_1 \left(\dfrac{1}{4}, \dfrac{1}{4}, 1; \dfrac{x}{x-1}\right)\right)^2
\end{equation*}
Mithilfe der folgenden Clausen'schen Formel aus der Arbeit (\textit{``Ueber die Fälle, wenn die Reihe von der Form $y=1+\frac{\alpha}{1}\cdot \frac{\beta}{\gamma}x+\frac{\alpha \cdot \alpha+1}{1\cdot 2}\cdot \frac{\beta \cdot \beta +1}{\gamma \cdot \gamma +1}x^2+\text{etc.}$ ein Quadrat der Form $z= 1+\frac{\alpha'}{1}\cdot \frac{\beta'}{\gamma'}\cdot \frac{\delta'}{\varepsilon'}x+\frac{\alpha' \cdot \alpha'+1}{1\cdot 2}\cdot \frac{\beta' \cdot \beta'+1}{\gamma' \cdot \gamma' +1}\cdot \frac{\delta' \cdot \delta' +1}{\varepsilon' \cdot \varepsilon' +1}x^2+\text{etc.}$ hat"}) (\cite{Cl28}, 1828)

\begin{equation}
\label{eq: Clausen Quadrat}
    _2F_1 \left(a,b,a+b+\frac{1}{2}; z\right)^2= {}_3F_2 \left(2a,2b,a+b; 2a+2b, a+b+\frac{1}{2}\right)
\end{equation}
mit

\begin{equation}
\label{eq: 3F2}
    _3F_2(\alpha,\beta,\gamma; a,b;z):= \sum_{n=0}^{\infty} \dfrac{(\alpha)_n (\beta)_n (\gamma)_n}{(a)_n (b)_n} \dfrac{z^n}{n!},
\end{equation}
hat man demnach insgesamt:

\begin{equation*}
     \left(_2F_1 \left(\dfrac{1}{4}, \dfrac{3}{4}, 1;x\right)\right)^2= (1-x)^{-\frac{1}{2}}{}_3F_2 \left(\dfrac{1}{2},\dfrac{1}{2},\dfrac{1}{2}; 1,1; \dfrac{x}{x-1}\right).
\end{equation*}
Nun wende man auf beide Seiten dieser Gleichung den Operator $x\frac{d}{dx}$ und setze anschließend $x=\frac{1}{2}$. Die linke Seite lässt sich dann mit dem Hilfssatz (\ref{Theorem: Hilfsformel 1/pi}) berechnen. Es bleibt die rechte Seite zu betrachten. Mit der Festlegung

\begin{equation*}
    _3F_2 \left(\dfrac{1}{2},\dfrac{1}{2},\dfrac{1}{2}; 1,1; x\right)= \sum_{n=0}^{\infty}a_n x^n,
\end{equation*}
gilt:

\begin{equation*}
    x\dfrac{d}{dx} (1-x)^{-\frac{1}{2}} {}_3F_2 \left(\dfrac{1}{2},\dfrac{1}{2},\dfrac{1}{2}; 1,1; \dfrac{x}{x-1}\right) = x\dfrac{d}{dx}\sum_{n=0}^{\infty} a_n \left(\dfrac{x}{x-1}\right)^n
\end{equation*}
\begin{equation*}
    = \dfrac{x(1-x)^{-\frac{3}{2}}}{2}\sum_{n=0}^{\infty}a_n \left(\dfrac{x}{x-1}\right)^n +(1-x)^{-\frac{3}{2}}\sum_{n=0}^{\infty} na_n \left(\dfrac{x}{x-1}\right)^n.
\end{equation*}
Setzt man in der Hilfsformel (\ref{Theorem: Hilfsformel 1/pi}) nun $s=\frac{1}{4}$, findet man schließlich:

\begin{equation*}
    \dfrac{\sin \frac{\pi}{4}}{\pi}= \dfrac{\sqrt{2}}{2\pi}= \dfrac{\sqrt{2}}{2}\sum_{n=0}^{\infty} \dfrac{\left(\frac{1}{2}\right)_n^3}{(1)_n^3}(4n+1)(-1)^n,
\end{equation*}
was gerade die Bauer'sche Reihe (\ref{eq: Bauer}) ist. Wie auch in \cite{Ha02} erwähnt, ist unklar, wie Ramanujan selbst die Bauer'sche Reihe gezeigt hat. Obwohl man Reihen ähnlicher Natur, die auch in seinen Brief \cite{Ra95} zu finden sind, mithilfe der Clausen'schen Formel (\ref{eq: Clausen Quadrat}) und Summationstheoremen für die hypergeometrische Reihe zeigen kann, ist zu vermuten, dass Ramanujan sie aus einer allgemeineren Summationsformel für die verallgemeinerte hypergeometrische Reihe abgeleitet hat. Während also der Ramanujan'sche Zugang eher ein \textit{``top--down"} Vorgang ist, ist der eben präsentierte ein eher \textit{``bottom-up"} Vorgang.

\paragraph{Divergente Reihen derselben Form}
\label{para: Divergnete Reihen}

Mit derselben Methode, mit welcher die Bauer'sche Formel (\ref{eq: Bauer}) nachgewiesen worden ist, lassen sich weitere Reihen der Form (\ref{eq: 1/pi}) zeigen. Man konsultiere die schon oben genannte Arbeit \cite{Al11} für mehrere Beispiele. Darunter finden sich unter anderen die divergenten Reihen:

\begin{equation*}
    \sum_{n=0}^{\infty}\dfrac{\left(\frac{1}{2}\right)^3_n}{(1)_n^3}(3n+1)(-8)^n =\dfrac{1}{\pi}
\end{equation*}
sowie

\begin{equation*}
    \sum_{n=0}^{\infty}\dfrac{\left(\frac{1}{2}\right)^3_n}{(1)_n^3}(3n+1)(4)^n =-\dfrac{2i}{\pi}.
\end{equation*}
 Die erste alternierende Reihe lässt sich mit den von Euler in \cite{E247} vorgestellten Ideen  ebenfalls zum selben Wert berechnen. Die zweite Reihe hingegen nicht, allerdings erinnert sie an die oben (Abschnitt \ref{subsubsec: Mehrwertige Funktionen: Das fehlende Konzept der Riemann'schen Fläche}) zitierte Euler'sche Äußerung, die Addition von lauter positiven Termen führe zu einer imaginären Summe, sofern man natürlich Summen einer divergenten Reihe im Euler'schen Sinne zulässt. Mit dem Beispiel des Ramanujan'schen Typs erfährt demnach auch diese Euler'sche Behauptung eine weitere Bestätigung und in der Tat ließen sich beide angeführten Reihen verwenden, um weitere ebenfalls richtige Reihen der Form (\ref{eq: 1/pi}) nachzuweisen.  Streng bewiesen und eingeordnet wurden diese divergenten Reihen dann von vielen weiteren Autoren. Eine zusammenfassende Darstellung findet man etwa im entsprechenden Kapitel des Buchs \textit{``Ramanujan's Theta Functions"} (\cite{Co17}, 2017).

\paragraph{Ein Beispiel ähnlicher Form bei Euler}
\label{para: Ein Beispiel ähnlicher Form bei Euler}

Abschließend soll noch eine Reihe der Gestalt (\ref{eq: 1/pi}) aus dem Euler'schen Opus angegeben werden. In seiner Arbeit \textit{``De eximio usu methodi interpolationum in serierum doctrina"}  (\cite{E555}, 1783, ges. 1773) (E555: ``Über den riesigen Nutzen der Interpolationsmethoden bei der Reihenlehre")  gelangt Euler in seinen Untersuchungen zu Interpolationen  in § 54 zu der allgemeinen Formel:

\begin{equation*}
    \dfrac{\sin (\omega \pi)}{\omega\pi}=1 + \dfrac{\omega^2 (\omega^2-1)}{1\cdot 2 \cdot 1^2}+\dfrac{\omega^2(\omega^2-1)^2(\omega^2-4)}{1\cdot 2 \cdot 3 \cdot 4 \cdot 1^2 \cdot 2^2}+\cdots
\end{equation*}
\begin{equation*}
    -\dfrac{\omega^2}{1}-\dfrac{\omega^2(\omega^2-1)^2}{1\cdot 2 \cdot 3 \cdot 1^2 \cdot 2}-\dfrac{\omega^2(\omega^2-1)^2(\omega^2-4)^2}{1\cdot 2 \cdot 3 \cdot 4 \cdot 5 \cdot 1^2 \cdot 2^2 \cdot 3}- \cdots
\end{equation*}
Euler betrachtet hiervon den Spezialfall $\omega=\frac{1}{2}$, welche Reihe sich wie nachstehend schreiben lässt:

\begin{equation*}
    1- \sum_{n=1}^{\infty} \left(\dfrac{(2n-1)!!}{n!}\right)^3 \dfrac{2n+1}{2n-1}\dfrac{1}{2^{5n}}-\dfrac{1}{4}\left(1+\sum_{n=1}^{\infty} \left(\dfrac{(2n-1)!!}{n!}\right)^3 \dfrac{2n+1}{n+1}\dfrac{1}{2^{5n}}\right).
\end{equation*}
Mithilfe der verallgemeinerten hypergeometrischen Funktion $_3F_2$ aus (\ref{eq: 3F2}) lässt sich dies wie folgt ausdrücken:

\begin{equation*}
    \dfrac{2}{\pi}={}_3F_2\left(-\dfrac{1}{2},\dfrac{1}{2},\dfrac{1}{2};1,1;\dfrac{1}{4}\right)-\dfrac{1}{16}\cdot {}_3F_2\left(\dfrac{1}{2},\dfrac{3}{2},\dfrac{3}{2};2,2;\dfrac{1}{4}\right)
\end{equation*}
\begin{equation*}
    -\dfrac{1}{4}\left({}_3F_2\left(\dfrac{1}{2},\dfrac{1}{2},\dfrac{1}{2};1,2;\dfrac{1}{4}\right)+\dfrac{1}{32}{}\cdot {}_3F_2\left(\dfrac{3}{2},\dfrac{3}{2},\dfrac{3}{2};2,3;\dfrac{1}{4}\right)\right).
\end{equation*}
Man kann sie also nicht ganz in die Ramanujan'sche Form bringen, jedoch ist die Ähnlichkeit vorhanden. Denn setzt man

\begin{equation*}
    _3F_2 \left(-\dfrac{1}{2},\dfrac{1}{2},\dfrac{1}{2};1,1;x\right)= \sum_{n=0}^{\infty} a_n x^n \quad \text{und} \quad _3F_2 \left(\dfrac{1}{2},\dfrac{1}{2},\dfrac{1}{2};1,1;x\right)= \sum_{n=0}^{\infty} b_n x^n,
\end{equation*}
kann man auch schreiben:

\begin{equation*}
    \dfrac{2}{\pi}= \sum_{n=0}^{\infty} (2n+1)a_n \left(\dfrac{1}{4}\right)^n - \dfrac{1}{4}  \sum_{n=0}^{\infty} (2n+1)b_n \left(\dfrac{1}{4}\right)^n
\end{equation*}
oder 

\begin{equation*}
    \dfrac{8}{\pi}= \sum_{n=0}^{\infty} (4a_n-b_n)(2n+1) \left(\dfrac{1}{4}\right)^n
\end{equation*}
oder gar nach Vereinfachung der Differenz $4a_n-b_n$, indem man nur $a_n$ übrig behält:

\begin{equation*}
    \dfrac{8}{\pi}= \sum_{n=0}^{\infty} 3\cdot\dfrac{(2n+1)^2}{n+1} a_n\left(\dfrac{1}{4}\right)^n,
\end{equation*}
welche Reihe ``fast"{} den Ramanujan'schen Typ (\ref{eq: 1/pi}) hat, zumal der Term $\frac{(2n+1)^2}{n+1}$ das Polynom $4n$ als Asymptote besitzt.


\newpage

\section{Zusammenfassung}
\label{sec: Zusammenfassung}

\epigraph{To do easily what is difficult for others is the mark of talent. To do what is impossible for talent is the mark of genius.}{Henri-Frédéric Amiel}


Die vorliegende Arbeit soll mit einer knappen Zusammenfassung der vorgetragenen Euler'schen Arbeitsexzerpte ihren Schlusspunkt finden.

\subsection{Grundlegendes}
\label{subsec: Grundlegendes}

\epigraph{In a different time, in a different place, it is always some other side of our common human nature that has been developing itself. The actual truth is the sum of all these.}{Thomas Carlyle}

Zu Eulers Zeiten, dem 18. Jahrhundert, ist Mathematik gänzlich anders betrieben und präsentiert worden als heutzutage. Die Kernaspekte der Euler'schen Philosophie und insbesondere die Unterschiede zum modernen Gegenstück sind anhand kurzer Beispiele illustriert worden (Abschnitt \ref{sec: Eulers Ansicht zur Mathematik}). Die wesentlichen Differenzen betreffen unter anderem die Auffassung eines Beweises (\ref{subsubsec: Eulers Auffassung eines Beweises}), die Rolle der Definition (\ref{subsubsec: Die Rolle der Definition in der Mathematik}) und auch die bewusste Darstellung kurioser Sachverhalte (\ref{para: Auch Falsches bzw. Bizarres wird mitgeteilt}). Während bei Euler die Beweise bisweilen eine heuristische Essenz haben und Euler den Leser durch induktive Argumente zu überzeugen versucht, wird heute wieder der synthetisch--deduktive Stil nach Euklid'schem Vorbild gewählt.\\

Neben der Vorstellung der Euler'schen Schaffensweise anhand bekannter Beispiele wie dem Baseler Problem (Abschnitt \ref{subsec: Die Lösung des Baseler Problems}) sind außerdem weniger in der Literatur diskutierte Exempel wie seine Auflösung von Differenzengleichungen (Abschnitt \ref{subsec: Ein Prototyp -- Differenzengleichungen}) vorgestellt worden. Insbesondere seine Behandlung des inhomogenen Falls mit konstanten Koeffizienten enkapsuliert die gesamte Bandbreite seiner Schaffensweise in einer einzigen Arbeit: Sie lässt die umfangreichen Euler'schen verbalen Erklärungen erkennen, illustriert Eulers genauen Kenntnisstand über von seinen Vorgängern erzielte Resultate, zeugt von seiner eigenen Kreativität und demonstriert überdies Eulers Fähigkeiten, sein gesamtes Wissen ertragreich zu synthetisieren. Überdies darf der Leser miterleben, wie ihn mathematisch waghalsige Manöver Bekanntes auf neue Art demonstrieren lassen, ihm erlauben die Grenzen der Mathematik auszudehnen, aber ihn auch auf Irrwege führen.  \\

Überdies wurden  einige Beispiele von modernen Konzepten genannt, welche Euler in seinen Arbeiten vorweg nimmt, wovon vielleicht die Mellin--Transformation (Abschnitt \ref{subsubsec: Die Mellin--Transformierte bei Euler}) und die Legendre--Polynome (Abschnitt \ref{subsubsec: Den Kontext betreffend: Die Legendre Polynome}) die am wenigsten erwarteten Erwähnungen sein dürften.

\subsection{Grenzen Eulers}
\label{subsec: Grenzen}

\epigraph{You can do anything, but not everything.}{David Allen}

Jedoch sind auch die Grenzen, an welche Euler gestoßen ist, Teil der Diskussionen gewesen (Abschnitt \ref{sec: Grenzen Eulers gezogen durch das Paradigma}), welche am eindrücklichsten in der Zahlentheorie ersichtlich geworden sind. So ist es nicht dem Zufall zuzuschreiben, dass  in diesem Zweig der Mathematik das Verhältnis zwar gemachter, jedoch nicht bewiesener Euler'scher Entdeckungen mit weitem Abstand das größte ist.   Das Reziprozitätsgesetz (Abschnitt \ref{subsubsec: Wegen fehlender Mittel: Das Reziprozitätsgesetz}) und (unter wohlwollender Auslegung seiner Ausführungen) gar den Primzahlsatz (Abschnitt \ref{subsubsec: Weg fehlender Formulierung: Der Primzahlsatz}) mag Euler zwar entdeckt haben, konnte selbige aber nicht mit einem Beweis sichern (Abschnitt \ref{subsec: Von Euler Entdecktes, aber nicht Beweisbares}). \\

An anderer Stelle waren schlicht die fehlenden modernen Auslegungen  von ihm selbst eingebrachter Konzepte, wofür dasjenige der Funktion (\ref{subsubsec: Der Begriff der Funktion}) und ihrer Eigenschaftsbegriffe (wie die Notion der Stetigkeit) sowie auf selbige angewandte Begriffe (wie das Integral) stellvertretend stehen, einem Vorankommen abträglich und verhinderten gar die Entwicklung neuer Begrifflichkeiten. Während eine gewisse Stagnation durch Stellen einer irreleitenden Frage (Abschnitt  \ref{subsec: Durch eine falsche Frage}), wie bei der Reduktion von elliptischen Integralen auf Normalformen (Abschnitt \ref{subsubsec: Die Normalform von elliptischen Integralen}) und Auflösungsformeln für die Wurzeln von Polynomen (Abschnitt \ref{subsubsec: Wegen Unbeweisbarkeit: Wurzeln von Polynomen}), leicht nachvollziehbar wird, mag sich das angesichts der vielen Einzelbeiträge beim Gegenstand der komplexen Analysis anders verhalten, welche Euler (Abschnitt \ref{subsubsec: Komplexe Analysis}) bis auf ihre Grundzüge versperrt geblieben ist. Dies war seinem Arbeitsethos geschuldet (Abschnitt \ref{subsec: Grenzen durch den eigenen Arbeitsethos}), welcher sich unter anderem in der, aus heutiger Sicht geurteilt, überbetonten Präferenz des Reellen und der damit einhergehenden unverstandenen Sonderstellung der komplexen Zahlen in der Analysis manifestiert. Daraus erklärt sich wohl auch, jedenfalls zum Teil, warum im Euler'schen Opus wenig zu den Jacobi'schen Theta--Funktionen gefunden wird (Abschnitt \ref{subsubsec: Durch einen neuen Gedanken: Die Theta-Funktion}), die ihre Eigenschaften in Gänze erst als Funktion über der komplexen Ebene zu erkennen geben. Nichtsdestotrotz erlaubte Euler sein Arbeitsduktus neben anderen Dingen die Etablierung der Variationsrechnung als eigenen Teilbereich der Mathematik, die Entdeckung der Euler--Maclaurin'schen Summenformel und die Derivation diverser Summen-- und Integralidentitäten, welche anderweitig schwer zu finden sind. Nicht zuletzt dieses Grundes wegen, ist die Euler'sche Priorität an vielerlei Stellen erst nach seinen Lebzeiten (wie bei der $\zeta$--Funktion (Abschnitt \ref{subsubsec: Durch Kombinieren von Ergebnissen: Die zeta-Funktion})) von nachfolgenden Mathematikergenerationen bemerkt worden. Bei Euler noch verstreut zu findende Ergebnisse wie beim Beispiel der hypergeometrischen Funktion (Abschnitt \ref{subsubsec: Die Darstellung betreffend -- Die hypergeometrische Reihe}) sind zwar zur Kenntnis genommen worden, jedoch scheint entsprechendes bei anderen Kontexten (insbesondere bei den Legendre--Polynomen (Abschnitt \ref{subsubsec: Den Kontext betreffend: Die Legendre Polynome}) und der Mellin--Transformation (Anschnitt \ref{subsubsec: Die Mellin--Transformierte bei Euler})) erst vor kurzem  erstmalig geschehen zu sein.   Die bei diesen und ähnlichen Fragen von Euler vernachlässigte Frage nach der Existenz einer Lösung hatte zwar bei seinem Beweis des Fundamentalsatzes der Algebra einen unglücklichen Ausgang (Abschnitt \ref{subsubsec: Methodus Inveniendi über Methodus Demonstrandi}), was allerdings eine der wenigen Ausnahmen dieser Art bleiben sollte. Überdies hat, wie ebenfalls ausgeführt, der Euler'sche Ansatz durch die Entwicklung der Algebra über 150 Jahre nach seinem Tod eine Rechtfertigung erfahren.\\

Schlussendlich enthielt der letzte Abschnitt (\ref{sec: Herleitungen aus Eulers Formeln und Ideen}) Beispiele dafür, wie weit sich mit den Euler'schen Kenntnissen allein noch  gelangen ließe. So wurden etwa über 150 Jahre später von Ramanujan bewiesene Formeln  lediglich unter Zuhilfenahme Euler'scher Resultate abgeleitet. Es wird ist zweifelsohne eine lohnenswerte Aufgabe sein, zu untersuchen, wie weit zum einen Euler selbst bei anderen Fragestellungen noch hätte vordringen können -- der für ihn selbst mögliche Nachweis der Funktionalgleichung der Riemann'schen $\zeta$--Funktion (Abschnitt (\ref{subsubsec: Durch Kombinieren von Ergebnissen: Die zeta-Funktion})) ist hierfür ein Beispiel -- und zum anderen, wie vieles  sich unter alleiniger Verwendung Euler'scher Resultate derivieren ließe\footnote{Zum gegenwärtigen Zeitpunkt planen die Betreiber des \textit{Euleriana} Journals eine Kolumne über diesen Gegenstand.}.

\subsection{Diskussion der Euler'schen Arbeitsweise}
\label{subsubsec: Diskussion über Arbeitsweise}

\epigraph{Euler fehlt nur eine Eigenschaft zu einem vollkommenen Genie: nämlich unverständlich zu sein.}{Georg Ferdinand Frobenius}

Obschon die Euler'schen Argumentationen, nach heutigem Verständnis, mancherorts rein formaler Natur sind und daher einen heute praktizierenden Mathematiker in dieser Hinsicht kaum überzeugen dürfen,  entbehren sie selten eines gewissen Reizes. Denn selbst in den Fällen eines Misserfolgs sind die Gründe dessen zumeist wenig offenkundig, insbesondere sofern man den tiefer liegenden Grund ausfindig zu machen sucht (Siehe eines Bespiele wegen Abschnitt (\ref{subsubsec: Ein Fehlschluss von Euler}), wo ein Fehler Eulers zum Beweis der Stirling'schen Formel für die Fakultät diskutiert worden ist.). \\

Abgesehen davon sind, nicht zuletzt  wegen ihrer heutigen Seltenheit, die Euler'schen Exempel  der \textit{Methodus inveniendi} äußerst illustrativ für die Arbeitsweise eines Mathematikers abseits der Publikationen und bringen einen bei eingehendem Studium womöglich der Frage näher, wie  sich in der Mathematik (und  anderen Wissenschaften) überhaupt zu Einsichten gelangen lässt. Denn trotz vieler Erklärungsversuche aus der Philosophie wie etwa \textit{``The Logic of Scientific Discovery"} (\cite{Po02}, 2002) oder auch \textit{``Erkenntnistheorie"} (\cite{Ba14}, 2014), der Psychologie zur Kreativitätsforschung \textit{``FLOW und Kreativität: Wie Sie Ihre
Grenzen überwinden und das Unmögliche schaffen"} (\cite{Cs14}, 2014), aber auch aus den Naturwissenschaften \textit{``The Art of Scientific Investigation"} (\cite{Be04}, 2004) und der Mathematik selbst wie etwa \textit{``The Mathematician’s Mind: The Psychology of
Invention in the Mathematical Field"} (\cite{Ha96}, 1996) bleibt Eulers Werk eine Bereicherung in dieser Hinsicht. Hat nämlich Euler selbst die wohl zuträglichste Methode zum \textit{Erlernen} der Mathematik erläutert, muss  seine Methode, Mathematik zu betreiben und  zu schaffen, aus seinen Arbeiten dahingegen in indirekter Weise entnommen werden. Dies ist jedoch eine Untersuchung für eine andere Abhandlung. In seiner auch in \cite{Fe95} (S. 11--12) zitierten Autobiographie schreibt Euler über seinen Unterricht bei Johann Bernoulli: \\

\textit{``Privat Lectionen schlug er mir zwar wegen seiner Geschäfte gänzlich ab: er gab mir aber einen weit heilsameren Rath, welcher darin bestund, dass ich selbsten einige mathematische Bücher vor mich nehmen, und mit allem Fleiss durchgehen sollte, und wo ich einigen Anstoss oder Schwierigkeiten finden möchte, gab er mir alle Sonnabend Nachmittag einen freyen  Zutritt bey sich, und hatte die Güte mir die gesammelete Schwierigkeiten  zu erläutern, welches mit so erwünschtem Vortheile geschah, dass wann er mir einen Anstoss gegeben hatte, dadurch zehn andere auf einmahl verschwanden, welches gewiss die beste Methode ist, um in den mathematischen Wissenschaften glückliche Progressen zu machen."}\\

Diese Euler'schen Ausführungen decken sich mit entsprechenden modernen Psychologie--Forschungen im Bereich Gedächtnis und Lernen. Genauer findet man in Eulers Zitat unter anderem den sogenannten ``Processing Difficulty Effect" \footnote{Dieser Effekt wird zum Beispiel in Arbeit \textit{``When comprehension difficulty improves memory for text"} (\cite{OB85}, 1985) in Zusammenhang des Memorierens von Textinhalten diskutiert, lässt sich aber analog auf andere Gebiete übertragen. Der Effekt besagt, dass man sich Inhalte umso besser einprägt, je schwieriger sie für einen selbst zu Beginn zu durchdringen waren.}. Hingegen gibt es bezüglich des kreativen Schaffens von Mathematik noch keine universell konsensuellen Theorien. Am trefflichsten fasst dies wohl der Weierstraß'sche Ausspruch zusammen (siehe \cite{Fe95}, Seite 124): \\

\textit{``Dass dem Forscher, solange er sucht, jeder Weg gestattet ist, versteht sich von selbst."}\\

Es scheint jedoch indes eine Notwendigkeit, sich das Wissen in einer solchen Weise angeeignet zu haben, dass man  es konsistent miteinander verweben kann. Den Schluss  bildet nachstehender in \textit{``Euler and the Calculus of Variations"} (\cite{Th07}, 2007) zitierter Ausspruch von Kneser (1898--1973):  \\

\textit{``Why do we rummage in rubble for some antiques? To enrich the ars inveniendi, to explain the methods by excellent examples, and last but not least to appreciate the intellectual company."}


\end{document}